\pdfoutput=1  %for arXiv rendering of xymatrix to come out right

\documentclass[11pt]{article}
\usepackage[margin=.8in,left=.8in]{geometry}
\usepackage{amsmath}
\usepackage{amsfonts}
\usepackage{amssymb}
\usepackage{amsthm}
\usepackage{mathtools}
\usepackage{subcaption}
\usepackage{graphicx}
\usepackage{color}
\usepackage{xcolor}
\usepackage{stmaryrd}
\usepackage{hyperref}
\usepackage[all]{xy}

\usepackage{tikz}
\usetikzlibrary{decorations.pathmorphing}
\usetikzlibrary{cd}

\usepackage[titletoc]{appendix}

\usepackage{wasysym}

\usepackage[safe]{tipa} % for \esh

\definecolor{darkblue}{rgb}{0.05,0.25,0.65}
\definecolor{greenii}{RGB}{20,140,10}
\definecolor{lightgray}{rgb}{0.9,0.9,0.9}
\definecolor{orangeii}{RGB}{200,100,5}

\usepackage{multirow}

\usepackage{stmaryrd}    % for \sslash

\usepackage{enumerate} %Roman enumerate

\usepackage{amssymb,amsmath}
\newcommand{\squig}{{\scriptstyle\sim\mkern-3.9mu}}

\newcommand{\rsquigend}{{\scriptstyle\rule{.1ex}{0ex}\rhd}}
\newcounter{sqindex}

\usepackage{mathptmx}
\usepackage{amsmath}
\usepackage{graphicx}
\DeclareRobustCommand{\coprod}{\mathop{\text{\fakecoprod}}}
\newcommand{\fakecoprod}{%
  \sbox0{$\prod$}%
  \smash{\raisebox{\dimexpr.9625\depth-\dp0}{\scalebox{1}[-1]{$\prod$}}}%
  \vphantom{$\prod$}%
}

\newcommand{\Sets}{
  \mathrm{Sets}
}

\newcommand{\Groups}{
  \mathrm{Grps}
}

\newcommand{\GroupOfUnits}[1]{
  \mathrm{GL}\left(1,{#1}\right)
}

\newcommand{\Abelian}{
  \mathrm{Ab}
}

\newcommand{\AbelianGroups}{
  \Abelian\Groups
}

\newcommand{\Groupoids}{
  \mathrm{Grpds}
}

\newcommand{\Categories}{
  \mathrm{Cats}
}

\newcommand{\Derived}{
  \mathbb{D}
}

\newcommand{\LeftDerived}{
  \Derived
}

\newcommand{\RightDerived}{
  \Derived
}

\newcommand{\Algebras}{
  \mathrm{Algs}
}

\newcommand{\dgcAlgebras}[1]{
  \mathrm{dgc}\Algebras
    _{\scalebox{.5}{${#1}$}}
    ^{\scalebox{.5}{$\geq 0$}}
}

\newcommand{\Bexp}{
  B\,\mathrm{exp}
}

\newcommand{\flatBexp}{
  \flat \mathbf{\mathbf{B}}\,\mathrm{exp}
}

\newcommand{\SullivanModels}{
  \mathrm{SullModels}
    _{\scalebox{.5}{$\mathbb{R}$}}
}

\newcommand{\SullivanModelsConnected}{
  \mathrm{SullModels}
    _{\scalebox{.5}{$\mathbb{R}$}}
    ^{\scalebox{.5}{$\geq 1$}}
}

\newcommand{\SullivanModelsConnectedOp}{
  \big(
    \SullivanModelsConnected
  \big)^{\mathrm{op}}
}

\newcommand{\TopologicalSpaces}{
  \mathrm{TopSp}
}

\newcommand{\Simplicial}{
  \Delta
}

\newcommand{\SimplicialSets}{
  \Simplicial\Sets
}

\newcommand{\SimplicialAbelianGroups}{
  \Simplicial\Abelian\Groups
}

\newcommand{\HomotopyTypes}{
  \mathrm{Ho}
  \big(
    \SimplicialSets_{\mathrm{Qu}}
  \big)
}

\newcommand{\PointedHomotopyTypes}{
  \mathrm{Ho}
  \big(
    \SimplicialSets^{\ast/}_{\mathrm{Qu}}
  \big)
}

\newcommand{\StableHomotopyTypes}{
  \mathrm{Ho}
  \big(
    \mathrm{SequentialSpectra}_{\mathrm{BF}}
  \big)
}

\newcommand{\TwistedHomotopyTypes}{
  \mathrm{Ho}
  \Big(
      \SimplicialSets_{\mathrm{Qu}}^{/B G}
  \Big)
}

\newcommand{\RationallyTwistedHomotopyTypes}{
  \mathrm{Ho}
  \Big(
    \big(
      \SimplicialSets_{\mathrm{Qu}}
    \big)^{/L_{\mathbb{R}}B G}
  \Big)
}

\newcommand{\Manifolds}{
  \mathrm{Mfds}
}

\newcommand{\Smooth}{
  \mathrm{Smth}
}

\newcommand{\SmoothManifolds}{
  \Smooth\Manifolds
}

\newcommand{\CartesianSpaces}{
  \mathrm{CartSp}
}

\newcommand{\Stacks}{
  \mathrm{Ho}
  (
    \mathrm{SmthStacks}_\infty
  )
}

\newcommand{\StacksTwisted}{
  \mathrm{Ho}
  \big(
    \mathrm{SmthStacks}^{/BG_{\mathrm{diff}}}_\infty
  \big)
}

\newcommand{\StacksOverX}{
  \mathrm{Ho}
  \big(
    \mathrm{SmthStacks}^{/\mathcal{X}}_\infty
  \big)
}

\newcommand{\StacksEtaleOverX}{
  \mathrm{Ho}
  \big(
    \mbox{\bf\'Et}_{\mathcal{X}}
  \big)
}

\newcommand{\LieAlgebras}{
   \mathrm{LieAlgebras}
     _{\scalebox{.5}{$\mathbb{R}, \mathrm{fin}$}}
}

\newcommand{\LInfinityAlgebrasGeneral}{
   L_\infty\Algebras
     _{\scalebox{.5}{$\mathbb{R}$}}
}

\newcommand{\LInfinityAlgebras}{
   L_\infty\Algebras
     ^{\scalebox{.5}{$\geq 0$}}
     _{\scalebox{.5}{$\mathbb{R}, \mathrm{fin}$}}
}

\newcommand{\LInfinityAlgebrasNil}{
   L_\infty\Algebras
     ^{\scalebox{.5}{$\geq 0, \mathrm{nil}$}}
     _{\scalebox{.5}{$\mathbb{R}, \mathrm{fin}$}}
}

\newcommand{\QFiniteHomotopyTypes}{
    \mathrm{Ho}
    \big(
      \SimplicialSets_{\mathrm{Qu}}
    \big)
    ^{\mathrm{fin}_{\mathbb{Q}}}
}

\newcommand{\NilpotentConnectedQFiniteHomotopyTypes}{
    \mathrm{Ho}
    \big(
      \SimplicialSets_{\mathrm{Qu}}
    \big)
    _{\geq 1, \mathrm{nil}}^{\mathrm{fin}_{\mathbb{Q}}}
}

\newcommand{\DifferentialGradedLieAlgebras}{
   \mathrm{DiffGradedLieAlgebras}
     ^{\scalebox{.5}{$\geq 0$}}
     _{\scalebox{.5}{$\mathbb{R}, \mathrm{fin}$}}
}

\newcommand{\dgcAlgebrasOp}[1]{
    \big(
      \dgcAlgebras{#1}
    \big)^{\mathrm{op}}
}

\newcommand{\dgcAlgebrasOpProj}[1]{
    \big(
      \dgcAlgebras{#1}
    \big)^{\mathrm{op}}_{\mathrm{proj}}
}

\newcommand{\dgcAlgebrasProj}[1]{
    \big(
      \dgcAlgebras{#1}
    \big)_{\mathrm{proj}}
}

\newcommand{\GradedAlgebras}{
  \mathrm{gcAlgs}
    ^{\scalebox{.5}{$\geq 0$}}
    _{\scalebox{.5}{$\mathbb{R}$}}
}

\newcommand{\GradedAlgebrasFin}{
  \mathrm{gcAlgs}
    ^{\scalebox{.5}{$\geq 0, \mathrm{fin}$}}
    _{\scalebox{.5}{$\mathbb{R}$}}
}

\newcommand{\CochainComplexes}{
    \mathrm{CochainComplexes}
      ^{\scalebox{.5}{$\geq 0$}}
      _{\scalebox{.5}{$\mathbb{R}$}}
}

\newcommand{\CochainComplexesInj}{
  \big(
    \CochainComplexes
  \big)_{\mathrm{inj}}
}

\newcommand{\ZChainComplexes}{
    \mathrm{ChainComplexes}
      ^{\scalebox{.5}{$\geq 0$}}
      _{\scalebox{.5}{$\mathbb{Z}$}}
}

\newcommand{\ZChainComplexesProj}{
  \big(
    \mathrm{ChainComplexes}
      ^{\scalebox{.5}{$\geq 0$}}
      _{\scalebox{.5}{$\mathbb{Z}$}}
  \big)_{\mathrm{proj}}
}

\newcommand{\ZChainComplexesUnbounded}{
    \mathrm{ChainComplexes}
      _{\scalebox{.5}{$\mathbb{Z}$}}
}

\newcommand{\ZChainComplexesUnboundedProj}{
  \big(
    \mathrm{ChainComplexes}
      _{\scalebox{.5}{$\mathbb{Z}$}}
  \big)_{\mathrm{proj}}
}

\newcommand{\VectorSpaces}{
  \mathrm{VectSp}
}

\newcommand{\Graded}{
  \mathrm{Grd}
}

\newcommand{\GradedVectorSpaces}{
    \Graded\VectorSpaces
      ^{\scalebox{.5}{$\geq 0$}}
      _{\scalebox{.5}{$\mathbb{R}$}}
}

\newcommand{\GradedVectorSpacesFin}{
    \Graded\VectorSpaces
      ^{\scalebox{.5}{$\geq 0, \mathrm{fin}$}}
      _{\scalebox{.5}{$\mathbb{R}$}}
}

\newcommand{\mapsup}{\rotatebox[origin=c]{90}{$\mapsto$}}
\newcommand{\mapsdown}{\rotatebox[origin=c]{-90}{$\mapsto$}}

\usepackage[new]{old-arrows}   %For \longhookrightarrow

\newdir{> }{{}*!/10pt/@{>}}

\usepackage{amssymb}

\DeclareRobustCommand{\rchi}{{\mathpalette\irchi\relax}}
\newcommand{\irchi}[2]{\raisebox{\depth}{$#1\chi$}} % inner command, used by \rchi

\makeatletter
\newif\if@sup
\newtoks\@sups
\def\append@sup#1{\edef\act{\noexpand\@sups={\the\@sups #1}}\act}%
\def\reset@sup{\@supfalse\@sups={}}%
\def\mk@scripts#1#2{\if #2/ \if@sup ^{\the\@sups}\fi \else%
  \ifx #1_ \if@sup ^{\the\@sups}\reset@sup \fi {}_{#2}%
  \else \append@sup#2 \@suptrue \fi%
  \expandafter\mk@scripts\fi}
\def\tensor#1#2{\reset@sup#1\mk@scripts#2_/}
\def\multiscripts#1#2#3{\reset@sup{}\mk@scripts#1_/#2%
  \reset@sup\mk@scripts#3_/}
\makeatother

\makeatletter
\newbox\slashbox \setbox\slashbox=\hbox{$/$}
\def\itex@pslash#1{\setbox\@tempboxa=\hbox{$#1$}
  \@tempdima=0.5\wd\slashbox \advance\@tempdima 0.5\wd\@tempboxa
  \copy\slashbox \kern-\@tempdima \box\@tempboxa}
\def\slash{\protect\itex@pslash}
\makeatother

\def\clap#1{\hbox to 0pt{\hss#1\hss}}
\def\mathllap{\mathpalette\mathllapinternal}
\def\mathrlap{\mathpalette\mathrlapinternal}
\def\mathclap{\mathpalette\mathclapinternal}
\def\mathllapinternal#1#2{\llap{$\mathsurround=0pt#1{#2}$}}
\def\mathrlapinternal#1#2{\rlap{$\mathsurround=0pt#1{#2}$}}
\def\mathclapinternal#1#2{\clap{$\mathsurround=0pt#1{#2}$}}

\let\oldroot\root
\def\root#1#2{\oldroot #1 \of{#2}}
\renewcommand{\sqrt}[2][]{\oldroot #1 \of{#2}}

\DeclareSymbolFont{symbolsC}{U}{txsyc}{m}{n}
\SetSymbolFont{symbolsC}{bold}{U}{txsyc}{bx}{n}
\DeclareFontSubstitution{U}{txsyc}{m}{n}

\DeclareSymbolFont{stmry}{U}{stmry}{m}{n}
\SetSymbolFont{stmry}{bold}{U}{stmry}{b}{n}

\DeclareFontFamily{OMX}{MnSymbolE}{}
\DeclareSymbolFont{mnomx}{OMX}{MnSymbolE}{m}{n}
\SetSymbolFont{mnomx}{bold}{OMX}{MnSymbolE}{b}{n}
\DeclareFontShape{OMX}{MnSymbolE}{m}{n}{
    <-6>  MnSymbolE5
   <6-7>  MnSymbolE6
   <7-8>  MnSymbolE7
   <8-9>  MnSymbolE8
   <9-10> MnSymbolE9
  <10-12> MnSymbolE10
  <12->   MnSymbolE12}{}

\usepackage{cleveref}

\crefformat{section}{\S#2#1#3} % see manual of cleveref, section 8.2.1
\crefformat{subsection}{\S#2#1#3}
\crefformat{subsubsection}{\S#2#1#3}

\theoremstyle{italics}
\newtheorem{theorem}{Theorem}[section]
\newtheorem{lemma}[theorem]{Lemma}
\newtheorem{prop}[theorem]{Proposition}

\theoremstyle{definition}
\newtheorem{defn}[theorem]{Definition}
\newtheorem{notation}[theorem]{Notation}
\newtheorem{example}[theorem]{Example}

\newtheorem{remark}[theorem]{Remark}

\usepackage{amsfonts}
\usepackage{colortbl}

\renewcommand{\emph}{\textit}

\begin{document}

\title{The character map in (twisted differential) non-abelian cohomology}

\author{Domenico Fiorenza, \; Hisham Sati, \; Urs Schreiber}

\maketitle

\begin{abstract}
  The Chern character on K-theory has a natural extension to
  arbitrary generalized cohomology theories, known
  as the Chern-Dold character. Here we further extend this
  to (twisted, differential) non-abelian
  cohomology theories, where its target is a
  non-abelian de Rham cohomology of
  twisted $L_\infty$-algebra valued differential forms.
  The construction amounts to leveraging the fundamental theorem
  of dg-algebraic rational homotopy theory
  to a twisted non-abelian generalization of the de Rham theorem.
  We show that the non-abelian
  character reproduces, besides the Chern-Dold character,
  also the Chern-Weil homomorphism
  as well as its secondary Cheeger-Simons homomorphism
  on (differential) non-abelian cohomology in degree 1,
  represented by principal bundles (with connection);
  and thus generalizes all these to
  higher (twisted, differential) non-abelian cohomology,
  represented by higher bundles/higher gerbes (with higher connections).
  As a fundamental example, we discuss the
  twisted non-abelian character map
  on twistorial Cohomotopy theory over 8-manifolds,
  which may be viewed as a twisted non-abelian enhancement of
  topological modular forms ($\mathrm{tmf}$) in degree 4.
  This turns out to exhibit a list of subtle topological relations
  that in high energy physics are thought to govern
  the charge quantization of fluxes in M-theory.

\end{abstract}

\medskip

\tableofcontents

\newpage

%%%%%%%%%%%%%%%%%%%%%%%%%%%%%
\section{Introduction}
%%%%%%%%%%%%%%%%%%%%%%%%%%%%%

Generalized cohomology theories \cite{Whitehead62}\cite{Adams74}
-- such as K-theory,  elliptic cohomology, stable Cobordism
and stable Cohomotopy -- are rich.
This makes them fascinating but also intricate to deal with.
In algebraic topology it has become commonplace to apply
filtrations by iterative \emph{localizations}
\cite{Bousfield79} (review in \cite[\S V]{EKMM97}\cite{Bauer14})
that allow generalized cohomology
to be approximated in consecutive stages;
a famous example of current interest is the chromatic filtration on
complex oriented cohomology theories
(\cite{MahowaldRavenel87}, review in \cite{Ravenel86}\cite{Lurie10}).

\medskip

\noindent {\bf The Chern-Dold character.}
The primary approximation stage of generalized cohomology theories is
their \emph{rationalization}
(e.g., \cite{Hilton71}\cite{Bauer14})
%[Ex. 1.7 (4)]{Bauer14})
to ordinary cohomology (e.g., singular cohomology) with rational coefficients
or real coefficients (see Remark \ref{RationalHomotopyTheoryOverTheRealNumbers}).
This goes back to \cite{Dold65}; and since
on topological K-theory (Ex. \ref{ChernCharacterInKTheory})
it reduces to the Chern character map
\cite[Thm. 5.8]{Hilton71},
this has been called the \emph{Chern-Dold character}
\cite{Buchstaber70}:

\vspace{-.4cm}
\begin{equation}
  \label{ChernDoldCharacterInIntroduction}
  \begin{tikzcd}[row sep = -8pt]
    &&
    E^n_{\color{blue}\mathbb{Q}}(X)
        \ar[
      rr,
      "\mbox{
        \tiny
        \color{greenii}
        \bf
        Dold's equivalence
      }"{below},
      "\sim"{above, yshift=-1pt}
    ]
    \ar[
      dd,
      "\mbox{
        \tiny
        \color{greenii}
        \bf
        \def\arraystretch{.9}
        \begin{tabular}{c}
          extensions
          \\
          of scalars
        \end{tabular}
      }"{left, xshift=7pt}
    ]
    &[+10pt]
    &
    \underset{k}{\bigoplus}
    \,
    H^{n+k}
    \big(
      X;\,
      \pi_k(E) \otimes_{{}_{\mathbb{Z}}} \color{blue}\mathbb{Q}
    \big)
    \mathrlap{
      \mbox{
       \tiny
       \color{darkblue}
       \bf
       \def\arraystretch{.9}
       \begin{tabular}{c}
         rational
         \\
         cohomology
       \end{tabular}
      }
    }
    \ar[
      dd,
      start anchor={[yshift=6pt]},
      end anchor={[yshift=-1pt]}
    ]
    \\
    \mathllap{
      \mbox{
        \tiny
        \color{darkblue}
        \bf
       \def\arraystretch{.9}
        \begin{tabular}{c}
          generalized
          \\
          cohomology
        \end{tabular}
      }
      \!\!
    }
    E^n(X)
    \ar[
      urr,
      "\mbox{
        \tiny
        \color{greenii}
        \bf
        rationalization
      }"{sloped, above}
    ]
    \ar[
      urrrr,
      rounded corners,
      to path={
           -- ([yshift=+29pt]\tikztostart.north)
           --node[above]{
             {
               \tiny
               \color{greenii}
               \bf
               Chern-Dold character
             }
             \mbox{
               \tiny
               $\mathrm{ch}^n_E$
             }
             } ([yshift=+8pt]\tikztotarget.north)
           -- (\tikztotarget.north)},
    ]
    \ar[
      drrrr,
      rounded corners,
      to path={
           -- ([yshift=-23pt]\tikztostart.south)
           --node[above]{
               \tiny
               \color{greenii}
               \bf
               differential-geometric Chern-Dold character
             } ([yshift=-8pt]\tikztotarget.south)
           -- (\tikztotarget.south)}
    ]
    \ar[drr]
    &&
    &&
    \\
    &&
    E^n_{\color{blue}\mathbb{R}}(X)
    \ar[
      rr,
      "\mbox{
        \tiny
        \color{greenii}
        \bf
        de Rham theorem
      }"{above},
      "\sim"{below, yshift=+1pt}
    ]
    &&
    \mathrm{Hom}_{\scalebox{.7}{\color{blue}$\mathbb{R}$}}
    \big(
      [\pi_\bullet(E),\mathbb{R}],
      \,
      H^{\bullet+n}_{\mathrm{dR}}(X)
    \big)
    \mathrlap{
      \mbox{
        \tiny
        \color{darkblue}
        \bf
        \def\arraystretch{.9}
        \begin{tabular}{c}
          de Rham
          \\
          cohomology
        \end{tabular}
      }
    }
  \end{tikzcd}
\end{equation}
\vspace{-.3cm}

\noindent That the left map in \eqref{ChernDoldCharacterInIntroduction}
is indeed the rationalization approximation on
coefficient spectra is left somewhat implicit in \cite{Buchstaber70}
(rationalization was properly formulated only in \cite{BousfieldKan72});
a fully explicit statement is in \cite[\S 2.1]{LindSatiWesterland16}.
The equivalence on the top of \eqref{ChernDoldCharacterInIntroduction}
serves to make explicit how the result of that rationalization operation
indeed lands in ordinary cohomology, and this was Dold's original
observation \cite[Cor. 4]{Dold65} (see Prop. \ref{DoldEquivalenceViaNonAbelianRealCohomology}).

\medskip

\noindent {\bf At the heart of differential cohomology.}
While rationalization is the coarsest of the localization approximations,
it stands out in that it connects, via the de Rham theorem,
to \emph{differential geometric} data
-- when the base space $X$ has the structure of a smooth manifold,
and the coefficients are taken to be $\mathbb{R}$ instead of
$\mathbb{Q}$.
Indeed, this ``differential-geometric Chern-Dold character''
shown on the bottom of \eqref{ChernDoldCharacterInIntroduction},
underlies (usually without attribution to either Dold or Buchstaber)
the pullback-construction of
differential generalized cohomology theories \cite[\S 4.8]{HopkinsSinger05}
(see \cite[p. 17]{BunkeNikolaus14}\cite[Def. 7]{GS-AHSS}\cite[Def. 17]{GS-KO}\cite[Def. 1]{GS-tAHSS}, recalled as
Def. \ref{DifferentialNonAbelianCohomology} and
Example \ref{DifferentialGeneralizedCohomology} below).

\medskip
\noindent {\bf At the heart of non-perturbative field theory.}
It is in this differential-geometric form
that the Chern-Dold character plays a pivotal role
in high energy physics. Here closed differential forms encode
\emph{flux densities} $F_p \in \Omega^p_{\mathrm{dR}}(X)$ of
generalized electromagnetic fields on spacetime manifolds $X$;
and the condition that these lift through (i.e., are in the image
of) the differential-geometric Chern-Dold character
\eqref{ChernDoldCharacterInIntroduction} for $E$-cohomology theory
encodes a \emph{charge quantization} condition in $E$-theory
(see \cite{Freed00}\cite{tcu}\cite{GS-RR}),
generalizing Dirac's charge quantization of the ordinary
electromagnetic field in ordinary cohomology
\cite{Dirac31} (review in \cite[\S 2]{Alvarez85}\cite[16.4e]{Frankel97}):

\vspace{-.2cm}
\begin{equation}
\label{DiracChargeQuantization}
\hspace{-1cm}
  \xymatrix@R=-12pt@C=4em{
    E^n(X)
    \ar[rr]_-{
        \mathrm{ch}^n_E
      }^-{
        \raisebox{3pt}{
          \tiny
          \color{greenii}
          \bf
          \def\arraystretch{.9}
          \begin{tabular}{c}
            differential-geometric
            \\
            Chern-Dold character
          \end{tabular}
        }
        }
      &&
    \mathrm{Hom}_{\mathbb{R}}
    \Big(
      \big[
        \pi_\bullet(E),\,\mathbb{R}
      \big]
      \,,\,
      H^{n+\bullet}_{\mathrm{dR}}(X)
    \Big)
    \\
    \underset{
      \mathclap{
      \raisebox{-3pt}{
        \tiny
        \color{darkblue}
        \bf
        \def\arraystretch{.9}
        \begin{tabular}{c}
          class in
          \\
          $E$-cohomology
        \end{tabular}
      }
      }
    }{
      [c]
    }
    \ar@{|->}[rr]
    &&
    \Big[
    \underset{
      \mathclap{
      \raisebox{+5pt}{
        \hspace{-2cm}
        \tiny
        \color{darkblue}
        \bf
        \def\arraystretch{.9}
        \begin{tabular}{c}
          charge-quantized
          flux densities
        \end{tabular}
      }
      }
    }{
      \big\{
        F^{\scalebox{.55}{$(a)$}}_{ r_a }
        \in
        \Omega^{r_a}_{\mathrm{dR}}(X)
      \big\}
      _{
        {
          1 \leq a \leq \mathrm{dim}[\pi_\bullet(E),\mathbb{R}]
        }
      }
    }
    \,\big\vert\,
     d\, F^{\scalebox{.55}{$(a)$}}_{r_a}
     \;=\;
     0
    \Big]
  }
\end{equation}
\vspace{-.5cm}

\noindent
This idea of charge quantization in a
generalized cohomology theory
has become famous for the case of
topological K-theory -- $E = \mathrm{KU}, \mathrm{KO}$ --
where it is argued to capture aspects of the expected
nature of the Ramond-Ramond (RR) fields in type II/I string theory
(see \cite{FreedHopkins00}\cite{Freed00}\cite{Evslin06}\cite{GS-RR}\cite{GS-KO}).

\medskip
\noindent
{\bf Need for non-abelian generalization.}
However, various further topological conditions
(recalled in \cite[Table 1]{FSS19b}\cite[p. 2]{FSS19c}\cite[Table 3]{SS20a}\cite[p. 2]{FSS20a},
see Rem. \ref{SummaryAndHypothesisH} below),
in non-perturbative type IIA string theory (``M-theory'')
are not captured by charge-quantization \eqref{DiracChargeQuantization}
in K-theory, nor in any Whitehead-generalized cohomology theory,
since they involve
{\it non-linear} functions \eqref{QuadraticFunctions} in the fluxes.

In order to systematically discuss the rich but under-appreciated area
of non-abelian charge quantization, we introduce and explore,
in \cref{ChernCharacterInNonabelianCohomology} and \cref{TwistedNonabelianCharacterMap},
the natural non-linear/non-abelian generalization of the character map.
This is based on
classical constructions of dg-algebraic rational homotopy theory
which we recall and develop in \cref{NonAbelianDeRhamCohomologyTheory}.

\medskip

\noindent {\bf The non-abelian character map.}
Indeed, despite their
established name, generalized cohomology theories
in the traditional sense of
Whitehead \cite{Whitehead62}\cite{Adams74} are
not general enough for many purposes:

\noindent
{\bf (i)} Already the time-honored non-abelian cohomology that
classifies
principal bundles (Ex. \ref{TraditionalNonAbelianCohomology} below),
being the domain of the Chern-Weil homomorphism \cite{Chern50}
(recalled as Def. \ref{ChernWeilHomomorphism},
Prop. \ref{FundamentalTheoremOfChernWeilTheory} below),
falls outside the scope of Whitehead-generalized cohomology.
Its {\it flat} sector alone, observed by secondary Cheeger-Simons invariants
(re-derived as Thm. \ref{SecondaryDifferentialNonAbelianCharacterSubsumesCheegerSimonsHomomorphism} below), is controlled by the classical {\it Maurer-Cartan equation}
(e.g. \cite[\S 5.6.4]{Nakahara03}\cite[Prop. 1.4.9]{RudolphSchmidt17})
on a Lie algebra valued differential form $A_{\!1}$:

\vspace{-.3cm}
\begin{equation}
  \label{ClassicalMaurerCartanEquation}
  d\, A^{\scalebox{.55}{$(c)$}}_{\!1}
  \;=\;
  f^c_{a b}
  \,
  A^{\scalebox{.55}{$(b)$}}_{\!1}
  \!
  \wedge
  A^{\scalebox{.55}{$(a)$}}_{\!1}
  \;\;\;
  \in
  \;
  \Omega^2_{\mathrm{dR}}(-)
\end{equation}
\vspace{-.5cm}

\noindent
(for $f^c_{a b}$  the structure constants, recalled as Ex. \ref{FlatLieAlgebraValuedDifferentialForms} below)
whose importance in large areas of mathematics
and mathematical physics is hard to overstate
(the ``master equation'', e.g. \cite[Rem. 3.12]{Markl12}\cite{ChuangLazarev13}),
but whose cohomological content is not captured by Whitehead-generalized
abelian cohomology theory.

\noindent
{\bf (ii)}
Similarly outside the scope of Whitehead-generalized cohomology
is the non-abelian cohomology classifying gerbes \cite{Giraud71}
(see Ex. \ref{NonAbelianGerbes} below). In its flat sector this
serves to adjoin to \eqref{ClassicalMaurerCartanEquation}
the higher-degree condition

\vspace{-.3cm}
\begin{equation}
  \label{StringyMaurerCartanEquation}
  d\, B_2
  \;=\;
  \mu_{a b c}
  \,
  A^{(a)}_{1}
  \wedge
  A^{(b)}_{1}
  \wedge
  A^{(a)}_{1}
  \;\;\;
  \in
  \;
  \Omega^3_{\mathrm{dR}}(-)
\end{equation}
\vspace{-.55cm}

\noindent
(for some differential 2-form, see Ex. \ref{FlatStringLie2AlgebraValuedDifferentialForms})
which has come to be recognized as a deep stringy refinement of the
classical Maurer-Cartan equation \eqref{ClassicalMaurerCartanEquation}
(see \cite[App.]{FSS12a} for pointers).

\vspace{1mm}
However
-- and this is our topic here --
these two items are just the first two stages within a truly general concept
of \emph{higher non-abelian cohomology}
(Def. \ref{NonAbelianCohomology} below),
that classifies higher bundles/higher gerbes
(Ex. \ref{NonAbelianCohomologyInUnboundedDegree} below),
whose non-abelian character map (Def. \ref{NonAbelianChernDoldCharacter} below)
takes values in flat $L_\infty$-algebra valued differential forms
(Def. \ref{FlatLInfinityAlgebraValuedDifferentialForms} below)
satisfying non-linear polynomial differential relations
($L_\infty$-algebraic Maurer-Cartan equations,  e.g. \cite[(31)]{DoubekMarklZima07}\cite[Def. 5.1]{Lazarev13})
which in string-theoretic applications
(see \eqref{FlatGSMechanismInIntroduction} and \cref{CohomotopicalChernCharacter} below)
are identified with
higher {\it Bianchi identities} on flux densities:

\vspace{-.8cm}
\small
\begin{equation}
  \label{NonAbelianCharacterMapInIntroduction}
  \begin{tikzcd}[row sep=0pt]
    \overset{
      \mathclap{
        \rotatebox{0}{
          \hspace{-12pt}
          \tiny
          \color{darkblue}
          \bf
          \def\arraystretch{.9}
          \begin{tabular}{c}
            non-abelian
            \\
            cohomology
          \end{tabular}
        }
      }
    }{
      A(X)
    }
    \ar[
      rr,
      "{
        \mathrm{ch}_A
      }"{below},
      "{
        \mbox{
          \tiny
          \color{greenii}
          \bf
          non-abelian character
        }
      }"
    ]
    &&
    \overset{
      \mathrlap{
        \;\;\;\;\;\;\;\;\;\;\;
        \rotatebox{16}{
          \hspace{-8pt}
          \tiny
          \color{darkblue}
          \bf
          \def\arraystretch{.9}
          \begin{tabular}{c}
            non-abelian
            \\
            de Rham cohomology
          \end{tabular}
        }
      }
    }{
      H_{\mathrm{dR}}
    }
    \big(
      X;
      \,
      \overset{
        \mathrlap{
          \;\;\;\;\;\;\;\;\;\;\;\;\;\;\;
          \rotatebox{16}{
            \tiny
            \color{darkblue}
            \bf
            Whitehead $L_\infty$-algebra
          }
        }
      }{
        \mathfrak{l}A
      }
    \big)
    \\
    \underset{
      \mathclap{
      \raisebox{-3pt}{
        \tiny
        \color{darkblue}
        \bf
        \def\arraystretch{.9}
        \begin{tabular}{c}
          class in
          \\
          $A$-cohomology
        \end{tabular}
      }
      }
    }
    {[c]}
    \ar[
      rr,
      |->
    ]
    &&
    \Big[
      \underset{
        \hspace{-2.2cm}
        \raisebox{+4pt}{
          \tiny
          \color{darkblue}
          \bf
          \def\arraystretch{.9}
          \begin{tabular}{c}
            charge-quantized
            flux densities
          \end{tabular}
        }
      }{
      \left.
      \big\{
        F^{\scalebox{.55}{$(a)$}}_{r_a}
        \,\in\,
        \Omega^{r_a}_{\mathrm{dR}}(X)
      \big\}_{1 \leq a \leq \mathrm{dim}[\pi_\bullet(A), \mathbb{R}] }
      \right\vert
      }
      \,
      \underset{
        \raisebox{-2pt}{
          \hspace{-.9cm}
          \tiny
          \color{darkblue}
          \bf
          higher Bianchi identities
        }
      }{
      d\, F^{\scalebox{.55}{$(a)$}}_{r_a}
        \;=\;
      P_{r_a}
      \big(
        \{F^{\scalebox{.55}{$(b)$}}_{r_b}\}_{b \leq a}
      \big)
      }
    \Big]
  \end{tikzcd}
\end{equation}
\vspace{-.4cm}

\noindent
This generalizes
(by Thm. \ref{NonAbelianChernCharacterSubsumesChernDoldCharacter} below)
the Chern-Dold character \eqref{DiracChargeQuantization}
on Whitehead-generalized cohomology,
which is subsumed as the abelian sector within the non-abelian theory
(Ex. \ref{GeneralizedCohomologyAsNonabelianCohomology}).

\medskip

While the non-abelian character map \eqref{NonAbelianCharacterMapInIntroduction}
is built from mostly classical ingredients
of dg-algebraic rational homotopy theory
(recalled and developed in\cref{RationalHomotopyTheory}),
its re-incarnation within non-abelian cohomology
provides a new unifying perspective on mathematical phenomena
expected to be relevant for non-perturbative physics:

\noindent
{\bf Yang-Mills monopoles via higher non-abelian cohomology.}
In modern formulation, Dirac's charge quantization
(e.g. \cite[\S 2]{Alvarez85})s
of the electromagnetic
field around a magnetic monopole with worldline
$\mathbb{R}^{0,1} \xhookrightarrow{\;} \mathbb{R}^{3,1}$,
is the statement that the topological class of the
field is encoded by
a continuous map from the
surrounding spacetime, which in the classical homotopy category
(Ex. \ref{TheClassicalHomotopyCategory})
is the 2-sphere
$
  \mathbb{R}^{3,1} \setminus \mathbb{R}^{0,1}
    \;\simeq\;
  \mathbb{R}^3 \setminus \{0\}
    \;\simeq\;
  S^2
  \;\;
  \in
  \;
  \mathrm{Ho}
  \big(
    \TopologicalSpaces_{\mathrm{Qu}}
  \big)
  \,,
$
to the classifying space of the circle group
$B \mathrm{U}(1) \,\simeq\, K(\mathbb{Z},2)$ \eqref{EilenbergMacLaneSpaces}:

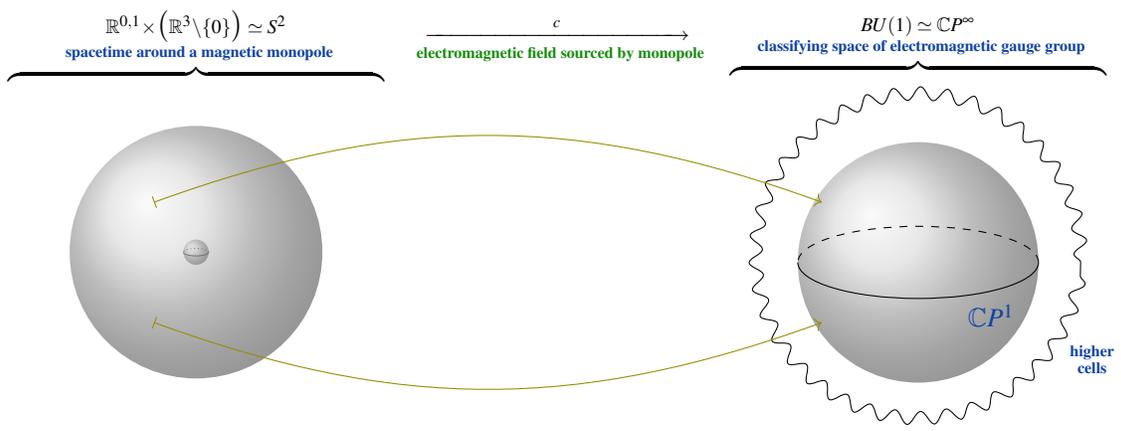
\vspace{-.3cm}

\begin{center}
\begin{tikzpicture}[decoration=snake, scale=.8]

  \begin{scope}[shift={(0,1)}]

  \begin{scope}

  \draw (6,4.1)
    node
    {
      $
        \underset{
               \mbox{
                 \tiny
                 \color{greenii}
                 \bf
                 \def\arraystretch{.9}
                 \begin{tabular}{c}
                   electromagnetic field
                   sourced by monopole
                 \end{tabular}
               }
        }{
        \xrightarrow{
          \mbox{
          \tiny
          $
          \phantom{AAAAAAAAAAAA}
          c
          \phantom{AAAAAAAAAAAA}
          $
          }
        }
        }
      $
     };

  \draw (0,4)
   node
     {
       $
         \overset{
           \underset{
             \mathclap{
             \mbox{
               \tiny
               \color{darkblue}
               \bf
               \def\arraystretch{.9}
               \begin{tabular}{c}
                 spacetime around
                 a magnetic monopole
               \end{tabular}
             }
             }
           }{
           \mathbb{R}^{0,1}
           \times
           \big(
           \mathbb{R}^3
           \setminus \{0\}
           \big)
           \;\simeq\;
           S^2
           }
         }{
           \overbrace{\phantom{\mbox{\hspace{5cm}}}}
         }
      $
     };

  \draw (12,4)
   node
     {
       $
         \overset{
           \underset{
             \mathclap{
             \mbox{
               \tiny
               \color{darkblue}
               \bf
               \def\arraystretch{.9}
               \begin{tabular}{c}
                 classifying space of
                 electromagnetic gauge group
               \end{tabular}
             }
             }
           }{
             B U(1)
             \;\simeq\;
             \mathbb{C}P^\infty
           }
         }{
           \overbrace{\phantom{\mbox{\hspace{5cm}}}}
         }
       $
     };

  \end{scope}

  \begin{scope}[shift={(0,.4)}]

  \begin{scope}[shift={(0,.8)}, scale=(.7)]

  \begin{scope}[scale=1.5]

  \begin{scope}[shift={(0,-.6)}]

  \shade[ball color=gray!40, opacity=.9] (0,0) circle (.2);
  \draw (-.2,0) arc (180:360:.2 and .06);
  \draw[densely dotted] (.2,0) arc (0:180:.2 and .06);

  \shade[ball color=gray!40, opacity=.6] (0,0) circle (2);

  \end{scope}

  \end{scope}

  \end{scope}

  \begin{scope}[shift={(12,0)}]

  \begin{scope}[scale=1]

  \shade[ball color=gray!40, opacity=.6] (0,0) circle (2);
  \draw (-2,0) arc (180:360:2 and 0.6);
  \draw[dashed] (2,0) arc (0:180:2 and 0.6);

  \draw[decorate] (0,0) circle  (2.7);
  \draw (-36:1.5) node
    {

      \color{darkblue}
      $\mathbb{C}P^1$
    };
  \draw (-30:3.3)
    node
    {
      \tiny
      \color{darkblue}
      \bf
      \def\arraystretch{.9}
      \begin{tabular}{c}
        higher
        \\
        cells
      \end{tabular}
    };

  \end{scope}

  \end{scope}

  \draw[|->, olive] (-.7,1) to[bend left=20] (10.4,1);
  \draw[|->, olive] (-.7,-1) to[bend right=20] (10.4,-1);

  \end{scope}

  \end{scope}

\end{tikzpicture}
\end{center}

\noindent
Since this is the classifying space for integral 2-cohomology (Exp. \ref{OrdinaryCohomology}),
one deduces, generally, that magnetic charge of the abelian
$\mathrm{U}(1)$-Yang-Mills field is measured in ordinary abelian cohomology
$H^2(-;\mathbb{Z})$ (e.g. \cite[p. 7]{FreedMooreSegal06}).
But the minimal cell decomposition of this classifying space
is by complex projective $k$-spaces for $k \in \mathbb{N}$:

\vspace{-.2cm}
\begin{equation}
  B \mathrm{U}(1)
  \;\simeq\;
  B^2 \mathbb{Z}
  \;=\;
  K(\mathbb{Z},2)
  \;\simeq\quad
  \overset{
    \mathclap{
    \raisebox{3pt}{
      \tiny
      \color{darkblue}
      \bf
      \def\arraystretch{.9}
      \begin{tabular}{c}
        infinite (stable, abelian)
        \\
        complex projective space
      \end{tabular}
    }
    }
  }{
    \mathbb{C}P^\infty
  }
  \quad \simeq \quad
  \overset{
    \mathclap{
    \raisebox{3pt}{
      \tiny
      \color{greenii}
      \bf
      \def\arraystretch{.9}
      \begin{tabular}{c}
        direct
        \\
        limit
      \end{tabular}
    }
    }
  }{
  \underset{
    \underset{
      k \to \infty
    }{
      \longrightarrow
    }
  }{\mathrm{lim}}
  }
  \;
  \mathbb{C}P^k
  \;\;
  \simeq
  \;\;
  \underset{
    \underset{
      k \to \infty
    }{
      \longrightarrow
    }
  }{\mathrm{lim}}
  \;
  \overset{
    \mathclap{
    \raisebox{3pt}{
      \tiny
      \color{darkblue}
      \bf
      \def\arraystretch{.9}
      \begin{tabular}{c}
        finite (unstable, non-abelian)
        \\
        complex projective
        $k$-space
      \end{tabular}
    }
    }
  }{
    \frac{
      \mathrm{SU}(k+1)
    }{
      \mathrm{U}(k)
    }
  }
  \,.
\end{equation}
\vspace{-.2cm}

\noindent
While none of these finite-dimensional stages $\mathbb{C}P^k$ by themselves classify
an abelian Whiethead-type cohomology theory, each of them classifies a
{\it higher non-abelian cohomology theory} $H^1\big(-; \Omega \mathbb{C}P^k \big)$
(by Ex. \ref{NonAbelianCohomologyInUnboundedDegree} below).

\vspace{1mm}
We observe that this formal mathematical fact
(Prop. \ref{ConnectedHomotopyTypesAreHigherNonAbelianClassifyingSpaces} below)
actually captures fine detail of the motivating physics,
in that this higher non-abelian deformation of abelian cohomology
measures magnetic charge of
{\it non-abelian magnetic monopoles}
in $\mathrm{SU}(k)$-Yang-Mills theory
(review in \cite{AtiyahHitchin88}\cite{Sutcliffe97})
obtained by reduction from
higher dimensional spacetimes $\mathbb{R}^{0,1} \times \mathbb{R}^3 \times X^d$
on a smooth fiber manifold $X^d$:

\vspace{-.6cm}
\begin{equation}
  \label{NonAbelianCPkCohomologyClassifiesMagneticCharges}
  \hspace{-8mm}
  \begin{tikzcd}[column sep=45pt, row sep=2pt]
    \overset{
      \mathclap{
      \raisebox{3pt}{
        \tiny
        \color{darkblue}
        \bf
        \def\arraystretch{.9}
        \begin{tabular}{c}
          gauge-equivalence class of moduli
          \\
          of $N$ magnetic monopoles on $\mathbb{R}^3$
          \\
          of $\mathrm{SU}(k+1)$-Yang-Mills theory
          \\
          minimally broken to $\mathrm{U}(k)$
        \end{tabular}
      }
      }
    }{
      \mathcal{M}
      \big(
        \mathrm{SU}(k+1)
      \big)_N
    }
    \ar[
      r,
      "{
        \mathclap{
        \mbox{
        \tiny
        \color{greenii}
        \bf
        \def\arraystretch{.9}
        \begin{tabular}{c}
          Donaldson-Jarvis'
          \\
          theorem
        \end{tabular}
        }
      }}"{below},
      "{
        \mbox{
          \tiny
          homeomorphism
        }
      }"{above}
    ]
    &[6pt]
    \overset{
      \mathclap{
      \raisebox{3pt}{
        \tiny
        \color{darkblue}
        \bf
        \def\arraystretch{.9}
        \begin{tabular}{c}
          holomorphic maps of algebraic degree $N$
          \\
          from Riemann sphere around monopoles
          \\
          to the complex projective $k$-manifold
        \end{tabular}
      }
      }
    }{
    \mathrm{Maps}_{\mathrm{hol}}
    \big(
      \mathbb{C}P^1\!
      ,
      \,
      \mathbb{C}P^k
    \big)_{\mathrm{deg}=N}
    }
    \ar[
      r,
      "{
        \mbox{
          \tiny
          \color{greenii}
          \bf
          \def\arraystretch{.9}
          \begin{tabular}{c}
            Segal's
            theorem
          \end{tabular}
        }
      }"{below}
    ]
    \ar[
     r,
     hook,
     "{
       \mbox{
         \tiny
         $N(2k-1)$-equivalence
         \eqref{WeakHomotopyEquivalence}
       }
     }"{above}
    ]
    &[16pt]
    \overset{
      \mathclap{
      \raisebox{3pt}{
        \tiny
        \color{darkblue}
        \bf
        \def\arraystretch{.9}
        \begin{tabular}{c}
          continuous maps of topological degree $N$
          \\
          from 2-sphere around monopoles
          \\
          to complex projective $k$-space
        \end{tabular}
      }
      }
    }{
    \mathrm{Maps}
    \big(
      S^2
      \!
      ,
      \,
      \mathbb{C}P^k
    \big)_{\mathrm{deg}=N}
    }
    \\[-3pt]
    \overset{
      \raisebox{3pt}{
        \tiny
        \color{darkblue}
        \bf
        \def\arraystretch{.9}
        \begin{tabular}{c}
          $X^d$-parameterized
          deformation classes
          \\
          of moduli of $N$ magnetic monopoles
        \end{tabular}
      }
    }{
      H
      \Big(
        X^d
        ;
        \,
        \mathcal{M}\big(\mathrm{SU}(k+1)\big)_N
      \Big)
    }
    \ar[
      rr,-,
      shift left=1pt,
      "{
        \mbox{
          \tiny
          \color{greenii}
          \bf
          for $d \,<\, N(2k-1)$
        }
      }"
    ]
    \ar[
      rr,-,
      shift right=1pt
    ]
    &&
    \overset{
      \mathclap{
      \raisebox{3pt}{
        \tiny
        \color{orangeii}
        \bf
        \def\arraystretch{.9}
        \begin{tabular}{c}
          higher non-abelian cohomology
        \end{tabular}
      }
      }
    }{
      H^1
      \big(
        X^d \times S^2
        ;\,
        \Omega\mathbb{C}P^k
      \big)
    }
    \underset{
      \underset{
        \mathclap{
        \raisebox{-1pt}{
          \tiny
          \color{darkblue}
          \bf
          Cohomotopy
        }
        }
      }{
        \pi^2(S^2)
      }
    }{\times}
    \{N\}
  \end{tikzcd}
\end{equation}
\vspace{-.5cm}

\noindent
This is a direct consequence
(using Prop. \ref{HomotopyClassesOfMapsOutOfnManifolds} below)
of classical theorems
shown in the first line of \eqref{NonAbelianCPkCohomologyClassifiesMagneticCharges}:
due to Donaldson (\cite{Donaldson84}, for $k =1$),
Jarvis
(\cite{Jarvis98}\cite{Jarvis00} for general $k$,
originally conjectured by Atiyah \cite[\S 5]{Atiyah84},
review in \cite{IoannidouSutcliffe99})
and Segal (\cite[Prop. 1.2]{Segal79}, see also \cite{CCMM91}\cite{Kamiyama07}).
Notice that the same moduli spaces of
holomorphic maps out of $\mathbb{C}P^1$
(often regarded and referred to as {\it rational maps} out of $\mathbb{C}$),
hence the same non-abelian cohomology sets \eqref{NonAbelianCPkCohomologyClassifiesMagneticCharges},
control numerous other aspects of non-abelian Yang-Mills theory,
notably the topological field configurations known as
{\it Skyrmions}
(an observation due to \cite{HoughtonMantonSutcliffe98}\cite{MantonPiette01}
whose homotopy-theoretic implications through Segal's theorem \eqref{NonAbelianCPkCohomologyClassifiesMagneticCharges}
have been found in \cite{Krusch03}\cite{Krusch06}),
which are  of deep relevance in non-perturbative quantum chromodynamics
(hadrodynamics),
not only theoretically but also experimentally
(review in \cite{RhoZahed16}, see \cite[p. 23]{BattyeMantonSutcliffe10}
for the impact of Segal's theorem \eqref{NonAbelianCPkCohomologyClassifiesMagneticCharges}).
Moreover, the homotopy quotient of these spaces
by the symmetries of $\mathbb{C}P^1$
(after compactification via adjoining of ``stable maps''
on degeneration limits of $\mathbb{C}P^1$)
govern the Gromov-Witten invariants of $\mathbb{C}P^k$
(review in \cite[\S 2]{Bertram02}\cite{Katz06} )
and,
for $k = 3$, the D-instantons of twistor string theory (\cite{Witten03}),
the scattering amplitudes of $\mathcal{N}=4$ super Yang-Mills theory
\cite{RoibanSpradlinVolovich04}
and those of $\mathcal{N} = 8$ supergravity
(\cite{CachazoSkinner12}\cite{Adamo15}),
for mathematical review see \cite[\S 7]{AtiyahDunajskiMason17}.

\medskip

\noindent
{\bf Non-abelian character of Yang-Mills monopoles.}
It turns out (Ex. \ref{NonAbelianCharacterOnCohomotopyTheory})
that the non-abelian character map \eqref{NonAbelianCharacterMapInIntroduction}
on these non-abelian magnetic monopole charges
\eqref{NonAbelianCPkCohomologyClassifiesMagneticCharges}
extracts the underlying abelian magnetic flux density $F_2$
together with a non-linear differential relation:

\vspace{-.7cm}
\begin{equation}
  \label{NonAbelianCharacterOnCPkCohomologyInIntroduction}
  \begin{tikzcd}[row sep=-1pt]
    &&
    \mathbb{C}P^1
    \mathrlap{
      \;
      =
      \,
      \mathrm{SU}(2)/\mathrm{U}(1)
    }
    \ar[
      d,
      hook
    ]
    &&[30pt]
    \def\arraystretch{.9}
    \begin{array}{l}
      \overset{
        \mathclap{
        \raisebox{3pt}{
          \tiny
          \color{darkblue}
          \bf
          \def\arraystretch{.9}
          \begin{tabular}{l}
            non-abelian character:
            \\
            non-linear Bianchi identity
          \end{tabular}
        }
        }
      }{ \scalebox{0.8}{$
        d\, H_3 \;=\; - F_2 \wedge F_2
      $}}
      \\
    \scalebox{0.8}{$  d\, F_2 \;=\; 0$}
    \end{array}
    \\
    &&
    \mathbb{C}P^2
    \mathrlap{
      \;
      =
      \,
      \mathrm{SU}(3)/\mathrm{U}(2)
    }
    \ar[
      d,
      hook
    ]
    &&
    \def\arraystretch{.9}
    \begin{array}{l}
     \scalebox{0.8}{$ d\, H_5 \;=\; - F_2 \wedge F_2 \mathrlap{\wedge F_2} $}
      \\
     \scalebox{0.8}{$ d\, F_2 \;=\; 0 $}
    \end{array}
    \\
    &&
    \mathbb{C}P^3
    \ar[
      d,
      hook,
      dotted
    ]
    \mathrlap{
      \;
      =
      \,
      \mathrm{SU}(4)/\mathrm{U}(3)
    }
    &&
    \def\arraystretch{.9}
    \begin{array}{l}
     \scalebox{0.8}{$ d\, H_7 \;=\; - F_2 \wedge F_2 \mathrlap{\wedge F_2 \wedge F_2} $}
      \\
     \scalebox{0.8}{$ d\, F_2 \;=\; 0 $}
    \end{array}
    \ar[
      d,
      phantom,
      "\vdots"
    ]
    \\[10pt]
    \mathllap{
      \underset{
        \mathrlap{
        \raisebox{-2pt}{
          \tiny
          \color{darkblue}
          \bf
          transversal spacetime
          }
        }
      }{
      X^d
      \times
      \big(
        \mathbb{R}^3 \!\setminus\! \{0\}
      \big)
      }
      \,
      \simeq
      \;\,
    }
    X^d
      \times
    \mathbb{C}P^1
    \ar[
      rr,
      dashed,
      "{
        \mbox{
          \tiny
          \color{greenii}
          \def\arraystretch{.9}
          \begin{tabular}{c}
            underlying abelian
            \\
            magnetic charge
          \end{tabular}
        }
      }"{below}
    ]
    \ar[
      urr,
      dashed,
      "{
        \mbox{
          \tiny
          \color{greenii}
          \bf
          $\mathrm{SU}(4)$-monopole
        }
      }"{above, yshift=-1pt, sloped},
      "{
        \mbox{
          \tiny
          \color{greenii}
          min. symm. break.
        }
      }"{below, yshift=1pt, sloped}
    ]
    \ar[
      uurr,
      dashed,
      bend left=5,
      "{
        \mbox{
          \tiny
          \color{greenii}
          \bf
          $\mathrm{SU}(3)$-monopole
        }
      }"{above, yshift=-1pt, sloped}
    ]
    \ar[
      uuurr,
      dashed,
      bend left=10,
      "{
        \mbox{
          \tiny
          \color{orangeii}
          \scalebox{1.1}{$X^d$}-parameterized
        }
      }"{left, xshift=-7pt, pos=.38},
      "{
        \mbox{
          \tiny
          \color{orangeii}
          \begin{tabular}{c}
            $\mathrm{SU}(2)$-monopole
          \end{tabular}
        }
      }"{sloped}
    ]
    &&
    \underset{
      \mathrlap{
      \;\;\;\;\;
      \raisebox{-3pt}{
        \tiny
        \color{darkblue}
        \bf
        classifying space
      }
      }
    }{
      \mathbb{C}P^\infty
          \mathrlap{
        \, \simeq\,
        B^2 \mathbb{Z}
      }
    }
    &&
    \underset{
      \mathllap{
      \raisebox{-3pt}{
        \tiny
        \color{darkblue}
        \bf
        \def\arraystretch{.9}
        \begin{tabular}{l}
          abelian character:
          \\
          linear Bianchi-identity
        \end{tabular}
      }
      \;\;\;\;\;\;
      }
    }{
     \scalebox{0.8}{$ d\, F_2 \;=\; 0 \phantom{F_2 \wedge F_2} $}
    }
  \end{tikzcd}
\end{equation}

\vspace{-2mm}
\noindent
While the algebraic form of this non-abelian character data
follows readily
--
once the non-abelian character map has been conceived in the first place,
that is, according to our Def. \ref{NonAbelianChernDoldCharacter}
--
from the well-known Sullivan model for
complex projective spaces (Ex. \ref{RationalizationOfComplexProjectiveSpace}),
it is curious to observe \cite{FSS20a}\cite{SS20c}
and seems to have gone unnoticed before\footnote{Recently the string physics community
is picking up some terminology of higher gauge theory
in interpreting the Green-Schwarz mechanism,
following \cite{SatiSchreiberStasheff08}\cite{SSS09}\cite{FSS12a},
identifying
the Green-Schwarz-type Bianchi identity \eqref{FlatGSMechanismInIntroduction},
Ex. \ref{CharacterMapOnJTwistedCohomotopyAndTwistorialCohomotopy}, as reflecting
{\it 2-group symmetry}, e.g. \cite[(1.18)]{CDI00}\cite[(3.3)]{DelZottoOhmori20}.
To justify this
terminology, one has to exhibit the GS-identity as
the higher curvature invariant of a higher gauge bundle,
hence as the non-abelian character of a higher non-abelian cohomology theory,
foundations for which we mean to lay out here. Our results
\cite{FSS20a}\cite{SS20c} (see \cref{CohomotopicalChernCharacter} below) indicate
that to account for the fine structure of string/M-theory a 2-group is not
sufficient, but a full $\infty$-group (Ex. \ref{NonAbelianCohomologyInUnboundedDegree},
such as $\Omega S^4$ (Rem. \ref{NonAbelianInfinityGroups}) equipped with
twisting and equivariance, is required. },
that its non-linear differential relations
\eqref{NonAbelianCharacterOnCPkCohomologyInIntroduction}
on magnetic flux densities are those of
important anomaly cancellation mechanisms in string theory.
Notably the first non-linear relation in the list

\vspace{-.7cm}
\begin{equation}
  \label{FlatGSMechanismInIntroduction}
  \begin{tikzcd}[row sep=0pt, column sep=-4pt]
    \overset{
      \mathclap{
      \raisebox{0pt}{
        \tiny
        \color{darkblue}
        \bf
        2-Cohomotopy
      }
      }
    }{
      \pi^2(-)
    }
    \;\coloneqq\;
    H^1
    \big(
      -;
      \,
      \Omega \mathbb{C}P^1
    \big)
    \ar[
      rr,
      "\mathrm{ch}_{
        \!
        \scalebox{.5}{
          $\mathbb{C}P^{\scalebox{.6}{$1$}}$
        }
      }"{below},
      "{
        \mbox{
          \tiny
          \color{greenii}
          \bf
          non-abelian character map
        }
      }"
    ]
    &&
    H_{\mathrm{dR}}
    \big(
      -
      ;
      \,
      \mathfrak{l}
      \mathbb{C}P^1
    \big)
    \\[-5pt]
    {[c]}
   \qquad  & \qquad \longmapsto& \qquad \qquad \qquad
    \scalebox{.8}{$
     \left[
    \left.
    \begin{array}{l}
      H_3 \,\in\, \Omega^3_{\mathrm{dR}}(-)
      \\
      \,F_2 \,\in\, \Omega^2_{\mathrm{dR}}(-)
    \end{array}
    \right\vert
    \def\arraystretch{1}
    \begin{array}{l}
      d\, H_3 \;=\; - F_2 \wedge F_2
      \\
      d\, \,F_2 \;=\; 0
    \end{array}
    \right]
    $}
    \begin{array}{l}
      \mbox{
        \tiny
        \color{orangeii}
        \bf
        Green-Schwarz-like Bianchi identity
      }
      \mathclap{\phantom{H_3}}
      \\
      \mbox{
        \tiny
        \color{darkblue}
        \bf
        ordinary abelian Bianchi identity
      }
      \mathclap{\phantom{F_2}}
    \end{array}
  \end{tikzcd}
\end{equation}
\vspace{-.4cm}

\noindent
has the non-linear form of
the Bianchi identity governing the {\it Green-Schwarz mechanism}
\cite[(4)-(6)]{GreenSchwarz84}\cite[[p. 49]{CHSW85}
(a mathematical account is in \cite[p. 40]{Freed00})
for anomaly cancellation in heterotic string theory
(the ``first superstring revolution'' \cite{Schwarz07}),
here for the case of {\it heterotic line bundles}
(of phenomenological interest \cite{AGLP12}\cite{AGLP11},
see \cite[\S 4.2]{ADO20a}\cite[\S 2.2]{ADO20b} for the case at hand),
namely heterotic gauge bundles whose gauge group is reduced along the symmetry breaking
$\mathrm{U}(1) \!\xhookrightarrow{\,}\! \mathrm{SU}(2) \!\xhookrightarrow{\,}\! E_8$
of the Yang-Mills monopole \eqref{NonAbelianCharacterOnCPkCohomologyInIntroduction}.
Of course, the full Green-Schwarz mechanism is as in
equation \eqref{FlatGSMechanismInIntroduction} but with a
further contribution from a gravitational flux.
This turns out to arise
through tangential {\it twisting} of the non-abelian character,
which is the main result of \cite{FSS19b}\cite{FSS19c}\cite{FSS20a}\cite{SS20c}
discussed in detail in \cref{CohomotopicalChernCharacter} below,
surveyed in a moment, in \eqref{TwistorialCharacterInIntroduction} below.

\medskip

\noindent {\bf Cohomotopy theory as higher non-abelian cohomology.}
The higher non-abelian cohomology theory
on the left of \eqref{FlatGSMechanismInIntroduction} is an example
of a classical concept in homotopy theory, namely of {\it Cohomotopy} sets
(Ex. \ref{CohomotopyTheory})
of homotopy classes of continuous maps into an $n$-sphere:

\vspace{-.4cm}
\begin{equation}
  \label{UnstableCohomotopyInIntroduction}
  \overset{
    \mathclap{
    \raisebox{1pt}{
      \tiny
      \color{darkblue}
      \bf
      $n$-Cohomotopy
    }
    }
  }{
    \pi^n(-)
  }
  \;\coloneqq\;
  \overset{
    \mathclap{
    \raisebox{-1pt}{
      \tiny
      \color{darkblue}
      \bf
      \def\arraystretch{.9}
      \begin{tabular}{c}
        homotopy classes of
        contin.
        \\
        maps into $n$-sphere
      \end{tabular}
    }
    }
  }{
    \mathrm{Ho}(\TopologicalSpaces_{\mathrm{Qu}})
    \big(
      -,
      \,
      S^n
    \big)
  }
  \; =\;
  \begin{tikzcd}
  \overset{
    \mathclap{
    \raisebox{3pt}{
      \tiny
      \color{darkblue}
      \bf
      \def\arraystretch{.9}
      \begin{tabular}{c}
        higher non-abelian
        \\
        $S^n$-cohomology
      \end{tabular}
    }
    }
  }{
    H^1
    \big(
      -;
      \,
      \Omega S^n
    \big)
  }
  \ar[
    rr,-,
    shift left=1pt
  ]
  \ar[
    rr,-,
    shift right=1pt,
    "{
      \mbox{
        \tiny
        \color{greenii}
        \bf
        Pontrjagin's theorem
      }
    }"{above},
    "{
      \mbox{
        \tiny
        (on clsd. manifolds)
      }
    }"{below}
  ]
  &[10pt]
  &
  \qquad
  \overset{
    \mathclap{
    \raisebox{3pt}{
      \tiny
      \color{darkblue}
      \bf
      \def\arraystretch{.9}
      \begin{tabular}{c}
        cobordism classes of
        submanifolds
        \\
        of codimension $n$
        \\
        and normally framed
      \end{tabular}
    }
    }
  }{
    \mathrm{Cob}^n_{\mathrm{Fr}}(-)\;.
  }
  \end{tikzcd}
\end{equation}
\vspace{-.3cm}

\noindent
The stabilization of \eqref{UnstableCohomotopyInIntroduction}
to {\it stable Cohomotopy} (Ex. \ref{StableCohomotopy})
is a widely recognized Whitehead-generalized cohomology theory,
usually discussed in the context of the stable Pontrjagin-Thom theorem.
But the original Pontrjagin theorem
(\cite{Pontrjagin55}, see \cite{SS21}\cite{SS19a} for review and further pointers)
is decidedly unstable and
as such says that the non-abelian Cohomotopy cohomology theory in \eqref{UnstableCohomotopyInIntroduction} measures
non-abelian charges carried by (normally framed) {\it sub}manifolds 
(``branes''),
which generalize the monopole charges in
\eqref{FlatGSMechanismInIntroduction} to higher codimension.

\medskip

\noindent
{\bf Non-abelian character of Cobordism.}
The non-abelian character \eqref{NonAbelianCharacterMapInIntroduction}
on unstable Cohomotopy/Cobordism \eqref{UnstableCohomotopyInIntroduction}
turns out (Ex. \ref{NonAbelianCharacterOnCohomotopyTheory} below)
to generalize the non-linear Green-Schwarz-type Bianchi identity
\eqref{FlatGSMechanismInIntroduction} to higher even degrees.
In degree 4 this yields the Bianchi identity of the $C_3$-field in
$D=11$ supergravity (due to \cite[\S 2.5]{Sati13}\cite[\S 2]{FSS16a}, review in \cite[\S 7]{FSS19a}),
which merges with the monopole characters \eqref{NonAbelianCharacterOnCPkCohomologyInIntroduction}
to the mixed Bianchi identity\footnote{The physics-inclined reader may want to think
of the broken $\mathrm{SU}(4)$ in \eqref{NonAbelianCharacterOn} as a {\it flavor} symmetry group,
along the lines of \cite{FSS20Flavor}.}
expected in Ho{\v r}ava-Witten's heterotic M-theory
(due to \cite{FSS20a}\cite{SS20c}, see \cref{CohomotopicalChernCharacter} below):

\vspace{-.8cm}
\begin{equation}
  \label{NonAbelianCharacterOn}
  \hspace{-3mm}
  \begin{tikzcd}[row sep=-1pt, column sep=37pt]
    &&
    \mathbb{C}P^1
    \ar[
      d,
      hook
    ]
    &&[-25pt]
    \overset{
      \mathclap{
      \raisebox{-1pt}{
        \hspace{-36pt}
        \tiny
        \color{darkblue}
        \bf
        non-abelian characters
      }
      }
    }{
    \begin{array}{l}
    \scalebox{.8}{$    d\, H_3 \;=\; - F_2 \wedge F_2 \phantom{G_4}$}
      \\
     \scalebox{.8}{$   d\, \, F_2 \;=\; 0 $}
    \end{array}
    }
    &[-40pt]
    \mbox{
      \tiny
      \color{darkblue}
      \bf
      \def\arraystretch{.9}
      \begin{tabular}{l}
        Green-Schwarz-type Bianchi identity
        \\
        on gauge field flux
      \end{tabular}
    }
    \\
    &&
    \mathbb{C}P^3
    \ar[
      d,
      ->>,
      "t_{\mathbb{H}}"{left},
      "
        \mbox{
          \tiny
          \color{greenii}
          \bf
          \def\arraystretch{.9}
          \begin{tabular}{c}
            twistor
            \\
            fibration
          \end{tabular}
        }
      "{right, xshift=-8pt}
    ]
    &&
    \begin{array}{l}
   \scalebox{.8}{$     d\, H_3 \;=\; G_4 - F_2 \wedge F_2 $}
      \\
\scalebox{.8}{$        d\, \, F_2 \;=\; 0 $}
      \\
   \scalebox{.8}{$     d \, G_7 \;=\; - \tfrac{1}{2} G_4 \wedge G_4 $}
      \\
  \scalebox{.8}{$      d \, G_4 \;=\; 0 $}
    \end{array}
    &
    \mbox{
      \tiny
      \color{orangeii}
      \bf
      \def\arraystretch{.9}
      \begin{tabular}{l}
        Ho{\v r}ava-Witten-type Bianchi identity
        \\
        on gauge \& $C_3$-field flux
      \end{tabular}
    }
    \\
    \phantom{
      X^d
        \times
      \big(
        \mathbb{R}^3 \!\setminus\! \{0\}
      \big)
      \simeq
      \;
    }
    \mathllap{
      X^d
        \times
      \big(
        \mathbb{R}^3 \!\setminus\! \{0\}
      \big)
      \simeq
      \;
    }
    X^d
    \times
    \mathbb{C}P^1
    \ar[
      urr,
      dashed,
      "{
        \mbox{
          \tiny
          \color{greenii}
          \bf
          $\mathrm{SU}(4)$-monopole
        }
      }"{above, yshift=-1pt, sloped},
      "{
        \mbox{
          \tiny
          \color{greenii}
          \bf
          breaking to $\mathrm{SU}(3)$
        }
      }"{below, yshift=1pt, sloped}
    ]
    \ar[
      uurr,
      dashed,
      bend left=6,
      "{
        \mbox{
          \tiny
          \color{greenii}
          \bf
          $\mathrm{SU}(2)$-monopole
        }
      }"{yshift=-1pt, sloped},
      "{
        \mbox{
          \tiny
          \color{greenii}
          \bf
          \scalebox{1.1}{$X^d$}-parametrized
        }
      }"{pos=.26, xshift=1pt}
    ]
    \ar[
      rr,
      "{
        \!\!
        \mbox{
          \tiny
          \color{orangeii}
          \bf
          underlying Cohomotopy
        }
      }"{above, yshift=-.5pt},
      "{
        \!\!
        \mbox{
          \tiny
          \color{orangeii}
          \bf
          {\color{black}/}
          underlying Cobordism
        }
      }"{below, yshift=+.5pt}
    ]
    &&
    S^4
    &&
    \begin{array}{l}
   \scalebox{.8}{$     d \, G_7 \;=\; - \tfrac{1}{2} G_4 \wedge G_4 $}
      \\
\scalebox{.8}{$        d \, G_4 \;=\; 0 $}
    \end{array}
    &
    \!\!\!\!\!\!\!\!\!\!\!\!
    \mbox{
      \tiny
      \color{orangeii}
      \bf
      \def\arraystretch{.9}
      \begin{tabular}{l}
        11d SuGra-type Bianchi identity
        \\
        on $C_3$-field flux
      \end{tabular}
    }
  \end{tikzcd}
\end{equation}
\vspace{-.5cm}

\noindent
{\bf Twisted Cohomotopy as twisted non-abelian cohomology.}
Classical constructions in differential topology revolving around
the Poincar{\'e}-Hopf theorem (e.g. \cite[\S 11]{BottTu82})
involve deformation classes of
non-vanishing vector fields on a smooth manifold $X$,
hence of homotopy-classes of sections of the unit sphere bundle
$S(T X)$ in the tangent bundle $T X$.
Generally, for
$\tau$ the class of a real vector bundle of rank $n + 1$ over
a paracompact Hausdoff space $X$,
we may consider the homotopy-classes of sections of
its unit sphere bundle $S(\tau)$ (with respect to any fiberwise metric)
as the $\tau$-{\it twisted} generalization
(Ex. \ref{JTwistedCohomotopyTheory})
of non-abelian $n$-Cohomotopy theory
from \eqref{UnstableCohomotopyInIntroduction}:

\vspace{-.4cm}
\begin{equation}
  \label{TwistedCohomotopyInIntroduction}
  \overset{
    \mathclap{
    \raisebox{3pt}{
      \tiny
      \color{darkblue}
      \bf
      \def\arraystretch{.9}
      \begin{tabular}{c}
        $\tau$-twisted
        \\
        Cohomotopy
      \end{tabular}
    }
    }
  }{
    \pi^\tau
    (-)
  }
  \;\coloneqq\;
  \overset{
    \raisebox{1pt}{
      \tiny
      \color{darkblue}
      \bf
      \def\arraystretch{.9}
      \begin{tabular}{c}
        homotopy classes of cts. sections
      \end{tabular}
    }
  }{
    \mathrm{Ho}
    \left(
      \TopologicalSpaces^{/X}_{\mathrm{Qu}}
    \right)
    \big(
      -
      ,\,
      S(\tau)
    \big)
  }
  \;=\;
  \begin{tikzcd}
    \overset{
      \mathclap{
      \raisebox{3pt}{
        \tiny
        \color{darkblue}
        \bf
        \def\arraystretch{.9}
        \begin{tabular}{c}
          $\tau$-twisted non-abelian
          \\
          $S^n$-cohomology
        \end{tabular}
      }
      }
    }{
      H^{\tau}
      \big(
        -
        ;\,
        S^n
      \big)
    }
    \;\;\;
    \ar[
      rr,-,
      shift left=1pt,
      "{
        \mbox{
          \tiny
          \color{greenii}
          \bf
          \def\arraystretch{.9}
          \begin{tabular}{c}
            twisted
            \\
            Pontrjagin theorem
          \end{tabular}
        }
      }"{yshift=-1pt,above},
      "{
        \mbox{
          \tiny
          \rm
          (for smth. $\tau \,\simeq\, \tau' \oplus 1$ )
        }
      }"{below}
    ]
    \ar[
      rr,-,
      shift right=1pt
    ]
    &[25pt]
    &
    \qquad
    \overset{
      \mathclap{
      \raisebox{3pt}{
        \tiny
        \color{darkblue}
        \bf
        \def\arraystretch{.9}
        \begin{tabular}{c}
          cobordims classes of submanifolds
          \\
          with normal $\tau'$-framing
        \end{tabular}
      }
      }
    }{
      \mathrm{Cob}^{\tau'}_{\mathrm{Fr}}
      (
        -
      )
      \mathrlap{\,.}
    }
  \end{tikzcd}
\end{equation}
\vspace{-.3cm}

\noindent
Indeed, when $\tau$ admits smooth structure and there is any
section of $S(\tau)$ at all, then a
twisted version of Pontrjagin's theorem
still applies (e.g. \cite[Lem. 5.2]{Cruickshank03})
to show that the twisted non-abelian cohomology theory
which we may call
{\it twisted Cohomotopy} \cite[\S 2.1]{FSS19b}
still measures
charges carried by cobordism classes of suitable submanifolds (``branes'').

\newpage

\noindent
{\bf Non-abelian character of twisted Cohomotopy.}
In \cref{TwistedNonabelianCharacterMap} below we construct
the twisted generalization of the
non-abelian character map \eqref{NonAbelianCharacterMapInIntroduction},
serving to extract differential form data
underlying such twisted non-abelian cohomology classes.
For instance, applied to the example of
tangentially twisted Cohomotopy \eqref{TwistedCohomotopyInIntroduction}
on even-dimensional smooth manifolds, this deforms the
Bianchi identity of ordinary odd-degree cohomology
by the Euler form $\rchi$ \eqref{FormulaForEulerForms} of the manifold:

\vspace{-.4cm}
\begin{equation}
  \label{CharacterOnTwistedOddCohomotopyInIntroduction}
  \begin{tikzcd}[row sep=-1pt, column sep=40pt]
    &[+10pt]
    \overset{
      \mathclap{
      \raisebox{3pt}{
        \tiny
        \color{darkblue}
        \bf
        \def\arraystretch{.9}
        \begin{tabular}{c}
          unit tangent
          \\
          bundle
        \end{tabular}
      }
      }
    }{
      S(T X)
    }
    \ar[
      r
    ]
    \ar[
      d
    ]
    &[-10pt]
    \overset{
      \mathclap{
      \raisebox{3pt}{
        \tiny
        \color{darkblue}
        \bf
        \def\arraystretch{.9}
        \begin{tabular}{c}
          universal
          $n$-sphere bundle
        \end{tabular}
      }
      }
    }{
      S^{2k-1} \!\sslash\! \mathrm{O}(2k)
    }
    \ar[
      d
    ]
    &&
    \overset{
      \raisebox{5pt}{
        \hspace{-9pt}
        \tiny
        \color{darkblue}
        \bf
        twisted non-abelian character
      }
    }{
\scalebox{.8}{$      d \, \theta_{2k-1} $}
        \;=\;
      \underset{
        \mathrlap{
        \;\;\;\;\;
        \raisebox{0pt}{
          \tiny
          \color{orangeii}
          \bf
          \def\arraystretch{.9}
          \begin{tabular}{c}
            tangential
            de Rham twist
          \end{tabular}
        }
        }
      }{
      \scalebox{.8}{$  \rchi_{2k}(\nabla) $}
      }
    }
    \\
    X
    \ar[
      r,-,
      shift left=1pt
    ]
    \ar[
      r,-,
      shift right=1pt
    ]
    \ar[
      ur,
      dashed,
      "{
        \mbox{
          \tiny
          \color{greenii}
          \bf
          \def\arraystretch{.9}
          \begin{tabular}{c}
            tangentially twisted
            \\
            Cohomotopy
          \end{tabular}
        }
      }"{sloped}
    ]
    &
    X
    \ar[
      r,
      "{
        \tau
        \;\coloneqq\;
        \vdash T X
      }"{above},
      "{
        \mbox{
          \tiny
          \color{orangeii}
          \bf
          \def\arraystretch{.9}
          \begin{tabular}{c}
            tangential
            twist
          \end{tabular}
        }
      }"{below}
    ]
    &
    B \mathrm{O}(2k)
    &&
    \!\!\!\!\!\!\!\!\!\!\!\!\!\!\!
    \overset{
      \mathrlap{
      \raisebox{+1pt}{
        \hspace{+15pt}
        \tiny
        \color{darkblue}
        \bf
        \def\arraystretch{.9}
        \begin{tabular}{c}
          Chern-Weil character
          of tangential twist
        \end{tabular}
      }
      }
    }{
    {\begin{array}{l}
\scalebox{.8}{$      d\, \rchi_{2k}(\nabla) \;=\; 0 $}
      \\
\scalebox{.8}{$      d\, p_{\bullet}(\nabla) \;=\; 0 $}
    \end{array}}
    \mathrlap{
    \begin{array}{l}
      \mbox{\tiny\color{darkblue}\bf Euler form}
      \\
      \mbox{\tiny\color{darkblue}\bf Pontrjagin form}
    \end{array}
    }
    }
  \end{tikzcd}
\end{equation}
\vspace{-.3cm}

\noindent
Hence the mere existence of the twisted non-abelian character on
odd-degree Cohomotopy
reflects part of the classical Poincar{\'e}-Hopf theorem (e.g. \cite[\S 11]{BottTu82}),
namely the vanishing of the Euler number of a manifold
implied by the existence of a unit vector field.
The further extension of this twisted non-abelian character
to even-degree Cohomotopy yields a tangentially twisted enhancement
of the classical Hopf invariant \cite[\S 4]{FSS19c}.

\medskip

\noindent
{\bf Twisted non-abelian character of Yang-Mills monopoles.}
The exceptional isomorphism
$\mathrm{Sp}(2) \simeq \mathrm{Spin}(5)$ between the quaternion
unitary group (the compact ``symplectic group'') and the spin-group in 5 dimensions,
together with the equivariance of the twistor fibration
\eqref{NonAbelianCharacterOn} under the canonical action of these groups
implies a unification of all the above examples in
a twisted non-abelian cohomology theory
which we will call {\it twistorial Cohomotopy}
(Ex. \ref{TwistorialCohomotopy} below).
The twisted non-abelian character on this theory is interesting in that
it exhibits a variety of aspects expected of non-linear Bianchi identities
in non-perturbative string theory
(due to \cite{FSS19b}\cite{FSS19c}\cite{FSS20a}\cite{SS20c},
surveyed in  Rem. \ref{SummaryAndHypothesisH} below)
which cannot be explained
by traditional twisted Whitehead-generalized cohomology theory
(Ex. \ref{TwistedGeneralizedCohomology}):

\vspace{-.3cm}
\begin{equation}
  \label{TwistorialCharacterInIntroduction}
  \begin{tikzcd}[column sep={between origins, 60pt}]
    &&
    \mathbb{C}P^3 \!\sslash\! \mathrm{Sp}(2)
    \ar[
      d,
      "{
        t_{\mathbb{H}}
          \!\sslash\!
        \mathrm{Sp}(2)
      }"{left},
      "{
        \mbox{
          \tiny
          \color{orangeii}
          \bf
          \def\arraystretch{.9}
          \begin{tabular}{c}
            equivariant
            \\
            twistor fibration
          \end{tabular}
        }
      }"{sloped}
    ]
    &[48pt]
    &
    \begin{array}{l}
      \overset{
        \mathrlap{
        \;\;\;\;\;\;\;\;\;\;\;\;\;\;
        \raisebox{2pt}{
          \tiny
          \color{darkblue}
          \bf
          twisted non-abelian character
        }
        }
      }{
      \;\;\,
  \scalebox{.8}{$       d \, H_3 $}
      }
      \;=\;
      \underset{
        \mathllap{
        \raisebox{-2pt}{
          \tiny
          \color{darkblue}
          \bf
          monopole char.
        }
        }
      }{    \scalebox{.8}{$
        - F_2 \wedge F_2
        $}
      }
       \scalebox{.8}{$    + G_4 $}
        \underset{
          \mathrlap{
          \;\;\;\;\;\;
          \raisebox{-2pt}{
            \tiny
            \color{orangeii}
            \bf
            tangential de Rham twist
          }
          }
        }{\scalebox{.8}{$
          - \tfrac{1}{4}p_1(\nabla)
          $}
        }
      {\phantom{AAAAAAAAAAAA}}
      \\
      \;\;\;  \scalebox{.8}{$    d \, F_2 \;=\; 0 $}
    \end{array}
    \\
    X
    \ar[
      dr,
      %"\tau"{above, pos=.7},
      "{
        \mbox{
          \tiny
          \color{greenii}
          \bf
          tangential twist
        }
      }"{below, sloped}
    ]
    \ar[
      rr,
      dashed,
      "{
        \mbox{
          \tiny
          \color{greenii}
          \bf
          twisted Cohomotopy
        }
      }"{above}
    ]
    \ar[
      urr,
      dashed,
      "{
        \mbox{
          \tiny
          \color{orangeii}
          \bf
          twistorial Cohomotopy
        }
      }"{above, sloped}
    ]
    &&
    S^4 \!\sslash\!  \mathrm{Sp}(2)
    \ar[
      dl
    ]
    &&
    \begin{array}{l}
     \scalebox{.8}{$    d\, 2 G_7 $}
        \;=\;
      \overset{
        \mathrlap{
          \raisebox{2pt}{
            \tiny
            \color{darkblue}
            \bf
            cobordism char.
          }
        }
      }{
      \scalebox{.8}{$
      -
      G_4
      \wedge
      G_4
      $}
      }
      \;
      \overset{
        \;\;\;\;\;\;\;\;\;\;\;
        \mathrlap{
        \raisebox{2pt}{
          \tiny
          \color{orangeii}
          tangential de Rham twist
        }
        }
      }{
     \scalebox{.8}{$    + $}
      }
\scalebox{.8}{$         \tfrac{1}{4}p_1(\nabla)
      \wedge
      \tfrac{1}{4}p_1(\nabla)
      -
      \rchi_8(\nabla)
      $}
      \\
      \;\; \scalebox{.8}{$   d \, G_4 \;=\; 0
      $}
    \end{array}
    \\[-15pt]
    &
    B \mathrm{Sp}(2)
  \end{tikzcd}
\end{equation}
\vspace{-.3cm}

\noindent
The desire to systematicall understand this rich example
(see \cref{CohomotopicalChernCharacter})
as a non-abelian
generalization of the traditional character map on twisted Whitehead-generalized
cohomology
originally motivated us to develop the theory of the
twisted non-abelian character map presented here.

\vspace{1mm}
To fully bring out the unifying picture, we will discuss in detail
how relevant examples of twisted Whitehead-generalized cohomology theories
are subsumed by our twisted non-abelian character map
(the ``inverse Whitehead principle'', Rem. \ref{WhiteheadPrincipleOfNonAbelianCohomology}).

\vspace{1mm}
\noindent
{\bf Non-abelian cohomology via classifying spaces.}
As shown by these motivating examples,
in higher non-abelian cohomology the very conceptualization of
cohomology finds a beautiful culmination, as it is reduced
to the pristine concept of

\vspace{.05cm}
\hspace{-.9cm}
\begin{tabular}{ll}
\begin{minipage}[l]{8cm}
homotopy types of mapping spaces
\eqref{NonabelianCohomologyMappingSpace},
or rather, if geometric (differential, equivariant,...) structures are
incorporated, of higher mapping stacks
(Remark \ref{StructuredNonAbelianCohomology} below).
In particular, the concept of \emph{twisted} non-abelian cohomology is
most natural from this perspective (Def. \ref{NonabelianTwistedCohomology} below)
and naturally subsumes the
traditional concept of twisted generalized
cohomology theories (Prop. \ref{ProofTwistedGeneralizedCohomology} below).
\end{minipage}
&
\hspace{5mm}
$
  \;\;\;\;\;
    \overset{
      \mathclap{
      \raisebox{3pt}{
        \tiny
        \color{darkblue}
        \bf
        \def\arraystretch{.9}
        \begin{tabular}{c}
          twisted
          \\
          non-abelian
          \\
          cohomology
        \end{tabular}
      }
      }
    }{
    H^\tau
    (
      X;
      \,
      A
    )
    }
    \;=\;
    \;
    \left\{\!\!\!\!\!\!\!\!
    \raisebox{13pt}{
    \xymatrix@C=20pt@R=1em{
      \;\;\;X\;\;
      \ar@{-->}[rr]^-{
        \overset{
          \mathclap{
          \raisebox{3pt}{
            \tiny
            \color{greenii}
            \bf
            cocycle
          }
          }
        }{
          c
        }
      }_>>>{\ }="s"
      \ar[dr]_-{
        \mathllap{
          \mbox{
            \tiny
            \color{greenii}
            \bf
            twist
          }
          \;
        }
        \tau
      }^-{\ }="t"
      &&
      \overset{
        \mathllap{
          \tiny
          \color{darkblue}
          \bf
          \def\arraystretch{.9}
          \begin{tabular}{c}
            \bf
            coefficient
            \\
            \bf
            $\infty$-stack
          \end{tabular}
        }
      }{
        A \!\sslash\! G
      }
      \ar[dl]^-{
        \underset{
          \mathrlap{
          \raisebox{-2pt}{
            \tiny
            \color{greenii}
            \bf
            \def\arraystretch{.9}
            \begin{tabular}{c}
              local
              \\
              coefficients
            \end{tabular}
          }
          }
        }{
          \rho
        }
      }
      \\
      &
      \mathbf{B} G
      \ar@{=>}^-{\simeq} "s"; "t"
    }
    }
    \right\}_{
      \!\!
      \big/
      \!\!\!\!\!\!\!
      \mbox{
        \tiny
        \def\arraystretch{.9}
        \begin{tabular}{c}
          homotopy
          \\
          relative $\mathbf{B} G$
        \end{tabular}
      }
    }
$
\end{tabular}

\vspace{.5mm}
\noindent {\bf State of the literature.}
It is fair to say that this
transparent fundamental nature of higher non-abelian cohomology
is not easily recognized in much of the traditional literature on
the topic, which is rife with unwieldy variants of cocycle conditions
presented in combinatorial $n$-category-theoretic language. As
a consequence, the development of non-abelian cohomology theory has
seen little and slow progress,
certainly as compared to the flourishing of
Whitehead-generalized cohomology theory.
In particular, the concepts of higher and of twisted non-abelian
cohomology had tended to remain mysterious (see \cite[p. 1]{Simpson97}).
It is the more recently established homotopy-theoretic formulation
of $\infty$-category theory (see Rem. \ref{InfinityCategoryTheory})
in its guise as  $\infty$-topos theory
({\it $\infty$-stacks},
recalled around Def. \ref{SmoothInfinityStacks} below)
that provides the backdrop on which
twisted higher non-abelian cohomology finds its true nature
\cite{Simpson97}\cite{Simpson99}\cite{Toen02}\cite[\S 7.1]{Lurie09}\cite{SSS09}\cite{NSS12a}\cite{NSS12b}\cite{dcct}\cite{FSS19b}\cite{SS20b};
see \cref{NonabCohomology} for details.

\newpage

\noindent {\bf The non-abelian character map.}
From this homotopy-theoretic perspective,
we observe, in \cref{ChernCharacterInNonabelianCohomology} and \cref{TwistedNonabelianCharacterMap},
that the generalization of the Chern-Dold character \eqref{ChernDoldCharacterInIntroduction}
to twisted non-abelian cohomology
naturally exists 
(Def. \ref{NonAbelianChernDoldCharacter}, Def. \ref{TwistedNonAbelianChernDoldCharacter}),
and that the non-abelian analogue of Dold's equivalence in
\eqref{ChernDoldCharacterInIntroduction} may neatly be understood as being,
up to mild re-conceptualization, the
fundamental theorem of dg-algebraic rational homotopy theory
(recalled as Prop. \ref{FundamentalTheoremOfdgcAlgebraicRationalHomotopyTheory} below).
We highlight that this classical theorem is fruitfully
recast as constituting a \emph{non-abelian de Rham theorem}
(Theorem \ref{NonAbelianDeRhamTheorem} below)
and, more generally, a
\emph{twisted non-abelian de Rham theorem}
(Theorem \ref{TwistedNonAbelianDeRhamTheorem} below).
With this in hand, the notion of the (twisted) non-abelian character map
appears naturally
(Def. \ref{NonAbelianChernDoldCharacter} and Def. \ref{TwistedNonAbelianChernDoldCharacter} below):

\vspace{-.4cm}
\begin{equation}
  \label{NonAbelianCharacterInIntroduction}
    \mathllap{
      \mbox{
        \tiny
        \def\arraystretch{.9}
        \begin{tabular}{c}
          \color{greenii}
          \bf
          twisted
          \\
          \color{greenii}
          \bf
          non-abelian
          \\
          \color{greenii}
          \bf
          character map
          \\
          (Def. \ref{TwistedNonAbelianChernDoldCharacter})
        \end{tabular}
      }
      \;\;
    }
    \mathrm{ch}_\rho
    \;:\;
    \xymatrix@C=3em{
      \overset{
        \raisebox{3pt}{
          \tiny
          \def\arraystretch{.9}
          \begin{tabular}{c}
            \color{darkblue}
            \bf
            twisted
            \\
            \color{darkblue}
            \bf
            non-abelian
            cohomology
            \\
            (Def. \ref{NonabelianTwistedCohomology})
          \end{tabular}
        }
      }{
      H^{\tau}
      (
        X;
        \,
        A
      )
      }
      \ar[rr]^-{
        (\eta^\mathbb{R}_\rho)_\ast
      }_-{
        \mbox{
          \tiny
          \begin{tabular}{c}
            \color{greenii}
            \bf
            $\mathbb{R}$-rationalization
            \\
            (Def. \ref{RationalizationOfCoefficientsInNonabelianCohomology})
          \end{tabular}
        }
      }
      &&
      \overset{
        \mathclap{
        \raisebox{3pt}{
          \tiny
          \begin{tabular}{c}
            \color{darkblue}
            \bf
            twisted
            non-abelian
            \\
            \color{darkblue}
            \bf
            real cohomology
            \\
            (Def. \ref{TwistedNonAbelianRealCohomology})
          \end{tabular}
        }
        }
      }{
      H^{L_{\mathbb{R}}\tau}
      \big(
        X;
        \,
        L_{\mathbb{R}}A
      \big)
      }
      \ar[rr]^-{ \simeq }_-{
        \mathclap{
        \mbox{
          \tiny
          \begin{tabular}{c}
            \color{greenii}
            \bf
            twisted
            non-abelian
            \\
            \color{greenii}
            \bf
            de Rham theorem
            \\
            (Thm. \ref{TwistedNonAbelianDeRhamTheorem})
          \end{tabular}
        }
        }
      }
      &&
      \overset{
        \mathclap{
        \raisebox{3pt}{
          \tiny
          \begin{tabular}{c}
            \color{darkblue}
            \bf
            twisted
            non-abelian
            \\
            \color{darkblue}
            \bf
            de Rham cohomology
            \\
            (Def. \ref{TwistedNonabelianDeRhamCohomology})
          \end{tabular}
        }
        }
      }{
      H^{\tau_{\mathrm{dR}}}_{\mathrm{dR}}
      (
        X;
        \,
        \mathfrak{l}A
      )\;.
      }
    }
\end{equation}

\noindent {\bf Twisted differential non-abelian cohomology.}
Moreover, with the (twisted) non-abelian character in hand,
the notion of (twisted) differential non-abelian cohomology
appears naturally
(Def. \ref{DifferentialNonAbelianCohomology}, Def. \ref{TwistedDifferentialNonAbelianCohomology})
together with the expected natural diagrams of
twisted differential non-abelian cohomology operations:

\vspace{-.9cm}
\begin{equation}
  \label{SystemsOfCohomologyOperationsOnDifferentialCohomology}
  \hspace{-.7cm}
  \raisebox{30pt}{
  \xymatrix@C=41pt@R=2.75em{
    \overset{
      \mathclap{
      \raisebox{3pt}{
        \tiny
        \begin{tabular}{c}
          \color{darkblue}
          \bf
          differential
          \\
          \color{darkblue}
          \bf
          non-abelian cohomology
          \\
          (Def. \ref{DifferentialNonAbelianCohomology})
        \end{tabular}
      }
      }
    }{
    \widehat H
    (
      X;
      \,
      A
    )
    }
    \ar[rr]^-{
      \mbox{
        \tiny
        \begin{tabular}{c}
          \color{greenii}
          \bf
          curvature
          \\
          \eqref{CurvatureOnDifferentialCohomology}
        \end{tabular}
      }
    }
    \ar[drr]|-{
      \mbox{
        \tiny
        \begin{tabular}{c}
          \color{greenii}
          \bf
          $\mathclap{\phantom{\vert^\vert}}$
          differential
          \\
          \color{greenii}
          \bf
          non-abelian character
          \\
          \eqref{DifferentialCharacterOnDifferentialCohomology}
          $\mathclap{\phantom{\vert_\vert}}$
        \end{tabular}
      }
    }
    \ar[d]_-{
      \mbox{
        \tiny
        \begin{tabular}{c}
          \color{greenii}
          \bf
          characteristic
          \\
          \color{greenii}
          \bf
          class
          \\
          \eqref{CharacteristicClassOnDifferentialCohomology}
        \end{tabular}
      }
      \!\!\!\!\!
    }
    &&
    \overset{
      \mathclap{
      \raisebox{3pt}{
        \tiny
        \begin{tabular}{c}
          \color{darkblue}
          \bf
          flat $L_\infty$-algebra valued
          \\
          \color{darkblue}
          \bf
          differential forms
          \\
          (Def. \ref{FlatLInfinityAlgebraValuedDifferentialForms})
        \end{tabular}
      }
      }
    }{
    \Omega_{\mathrm{dR}}
    (
      X;
      \,
      \mathfrak{l}A
    )
    }
    \ar[d]^-{
      \mbox{
        \tiny
        \color{darkblue}
        \bf
        \begin{tabular}{c}
        \end{tabular}
      }
    }
    \\
    \underset{
      \mathclap{
      \raisebox{-3pt}{
        \tiny
        \begin{tabular}{c}
          \color{darkblue}
          \bf
          non-abelian cohomology
          \\
          (Def. \ref{NonAbelianCohomology})
        \end{tabular}
      }
      }
    }{
    H
    (
      X
      ;
      \,
      A
    )
    }
    \ar[rr]_-{
      \mbox{
        \tiny
        \begin{tabular}{c}
          \color{greenii}
          \bf
          non-abelian character
          \\
          (Def. \ref{NonAbelianChernDoldCharacter})
        \end{tabular}
      }
    }^-{
      \mathrm{ch}_{A}
    }
    &&
    \underset{
      \mathclap{
      \raisebox{-3pt}{
        \tiny
        \begin{tabular}{c}
          \color{darkblue}
          \bf
          non-abelian
          \\
          \color{darkblue}
          \bf
          de Rham cohomology
          \\
          (Def. \ref{NonabelianDeRhamCohomology})
        \end{tabular}
      }
      }
    }{
    H_{\mathrm{dR}}
    (
      X
      ;
      \,
      \mathfrak{l}A
    )
    }
  }
  }
  %%%%%%%%%%%%%%%
  \phantom{AAA}
  %%%%%%%%%%%%%%%
  \raisebox{30pt}{
  \xymatrix@C=41pt@R=2.5em{
    \overset{
      \mathclap{
      \raisebox{3pt}{
        \tiny
        \begin{tabular}{c}
          \color{darkblue}
          \bf
          twisted differential
          \\
          \color{darkblue}
          \bf
          non-abelian cohomology
          \\
          (Def. \ref{TwistedDifferentialNonAbelianCohomology})
        \end{tabular}
      }
      }
    }{
    \widehat H^{\tau_{\mathrm{diff}}}
    (
      X;
      \,
      A
    )
    }
    \ar[rr]^-{
      \mbox{
        \tiny
        \begin{tabular}{c}
          \color{greenii}
          \bf
          twisted
          \\
          \color{greenii}
          \bf
          curvature
          \\
          \eqref{CurvatureOnTwistedDifferentialCohomology}
        \end{tabular}
      }
    }
    \ar[drr]|-{
      \mbox{
        \tiny
        \begin{tabular}{c}
          \color{greenii}
          \bf
          $\mathclap{\phantom{\vert^\vert}}$
          twisted differential
          \\
          \color{greenii}
          \bf
          non-abelian character
          \\
          \eqref{DifferentialCharacterOnTwistedDifferentialCohomology}
          $\mathclap{\phantom{\vert_\vert}}$
        \end{tabular}
      }
    }
    \ar[d]_-{
      \mbox{
        \tiny
        \begin{tabular}{c}
          \color{greenii}
          \bf
          twisted
          \\
          \color{greenii}
          \bf
          characteristic
          \\
          \color{greenii}
          \bf
          class
          \\
          \eqref{CharacteristicClassOnTwistedDifferentialCohomology}
        \end{tabular}
      }
      \!\!\!\!\!
    }
    &&
    \overset{
      \mathclap{
      \raisebox{3pt}{
        \tiny
        \begin{tabular}{c}
          \color{darkblue}
          \bf
          twisted flat
          \\
          \color{darkblue}
          \bf
          $L_\infty$-algebra valued
          \\
          \color{darkblue}
          \bf
          differential forms
          \\
          (Def. \ref{FlatTwistedLInfinityAlgebraValuedDifferentialForms})
        \end{tabular}
      }
      }
    }{
    \Omega_{\mathrm{dR}}^{\tau_{\mathrm{dR}}}
    (
      X;
      \,
      \mathfrak{l}A
    )
    }
    \ar[d]^-{
      \mbox{
        \tiny
        \color{darkblue}
        \bf
        \begin{tabular}{c}
        \end{tabular}
      }
    }
    \\
    \underset{
      \mathclap{
      \raisebox{-3pt}{
        \tiny
        \begin{tabular}{c}
          \color{darkblue}
          \bf
          twisted
          \\
          \color{darkblue}
          \bf
          non-abelian cohomology
          \\
          (Def. \ref{NonabelianTwistedCohomology})
        \end{tabular}
      }
      }
    }{
    H^\tau
    (
      X
      ;
      \,
      A
    )
    }
    \ar[rr]_-{
      \mbox{
        \tiny
        \begin{tabular}{c}
          \color{greenii}
          \bf
          twisted
          \\
          \color{greenii}
          \bf
          non-abelian character
          \\
          (Def. \ref{TwistedNonAbelianChernDoldCharacter})
        \end{tabular}
      }
    }^-{
      \mathrm{ch}^{\tau}_{A}
    }
    &&
    \underset{
      \mathclap{
      \raisebox{-3pt}{
        \tiny
        \begin{tabular}{c}
          \color{darkblue}
          \bf
          twisted non-abelian
          \\
          \color{darkblue}
          \bf
          de Rham cohomology
          \\
          (Def. \ref{TwistedNonabelianDeRhamCohomology})
        \end{tabular}
      }
      }
    }{
    H^{\tau_{\mathrm{dR}}}_{\mathrm{dR}}
    (
      X
      ;
      \,
      \mathfrak{l}A
    )
    }
  }
  }
\end{equation}

\noindent {\bf Unifying Chern-Dold, Chern-Weil and Cheeger-Simons.}
In order to show that this generalization
of (twisted) character maps and
(twisted) differential cohomology
to higher non-abelian cohomology
is sound, we proceed to prove that the non-abelian character map
(Def. \ref{NonAbelianChernDoldCharacter})
specializes to:
\begin{center}
\begin{tabular}{ll}
\rowcolor{lightgray}
\begin{tabular}{l}
the Chern-Dold character
\\
\phantom{AA} on generalized cohomology
\end{tabular}
&
(Theorem \ref{NonAbelianChernCharacterSubsumesChernDoldCharacter}),
\\
\begin{tabular}{l}
the Chern-Weil homomorphism
\\
\phantom{AA} on degree-1 non-abelian cohomology
\end{tabular}
&
(Theorem \ref{NonAbelianChernDoldSubsumesChernWeil}),
\\
\rowcolor{lightgray}
\begin{tabular}{l}
the Cheeger-Simons homomorphism
\\
\phantom{AA} on degree-1 differential non-abelian cohomology
\end{tabular}
&
(Theorem \ref{SecondaryDifferentialNonAbelianCharacterSubsumesCheegerSimonsHomomorphism}).
\end{tabular}
\end{center}

\noindent All these classical invariants
are thus seen as different low-degree aspects
of the higher non-abelian character map.

\medskip

\noindent {\bf Examples of twisted higher character maps.}
To illustrate the mechanism, we make explicit
several examples of the (twisted) non-abelian character map
on cohomology theories of relevance in high energy physics:

\medskip
\hspace{2cm}
\def\arraystretch{1.2}
\begin{tabular}{ll}
\rowcolor{lightgray}
the Chern character on complex differential K-theory
&
(Example \ref{ChernCharacterInKTheory}, \ref{CurvatureInDifferentialComplexKTheory}),
\\
the Pontrjagin character on real K-theory
&
(Example \ref{PontrjaginCharacterInKSO}),
\\
\rowcolor{lightgray}
the Chern character on twisted differential K-theory
&
(Example \ref{TwistedChernCharacterInTwistedTopologicalKTheory},
\ref{TwistedDifferentialChernCharacterOnTwistedDifferentialKTheory}),
\\
the MMS-character on cohomotopy-twisted K-theory
&
(Example \ref{ChernCharacterOnHigherCohomotopicallyTwistedKTheory}),
\\
\rowcolor{lightgray}
the LSW-character on twisted higher K-theory
& (Example \ref{TwistedChernCharacterInTwistedHigherKTheory}),
\\
the character on integral Morava K-theory
&
(Example \ref{PontrjaginCharacterInMorava}),
\\
\rowcolor{lightgray}
the character on topological modular forms, $\mathrm{tmf}$
&
(Example \ref{CharacterInTMF}).
\end{tabular}
\vspace{.2cm}

\noindent Once incarnated this way within
the more general context of non-abelian cohomology theory,
we may ask for non-abelian enhancements
(Example \ref{StabilizationCohomologyOperation})
of these abelian character maps:

\medskip

\noindent {\bf Non-abelian enhancement of the
$\mathrm{tmf}$-character -- the cohomotopical character.}
Our culminating example, in \cref{CohomotopicalChernCharacter},
is the character map on twistorial Cohomotopy theory
\cite{FSS19b}\cite{FSS20a}, over 8-manifolds $X^8$ equipped with
tangential $\mathrm{Sp}(2)$-structure $\tau$ \eqref{TangentialSp2Structure}.
This may be understood (Remark \ref{ClarifyingTheRoleOfTmfInStringTheory})
as an enhancement of the $\mathrm{tmf}$-character
(Example \ref{CharacterInTMF})
from traditional generalized cohomology
to twisted differential non-abelian cohomology:

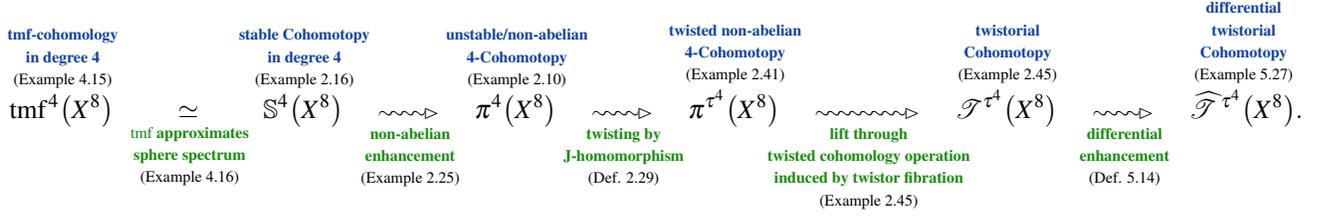
\vspace{-.6cm}
$$
\hspace{3mm}
  \overset{
    \mathclap{
    \raisebox{9pt}{
      \tiny
      \begin{tabular}{c}
        \color{darkblue}
        \bf
        tmf-cohomology
        \\
        \color{darkblue}
        \bf
        in degree 4
        \\
        (Example \ref{CharacterInTMF})
      \end{tabular}
    }
    }
  }{
    \mathrm{tmf}^{\, 4}\big(X^8\big)
  }
  \;\;\;\;\;\;\;\;\;\;
  \underset{
    \mathclap{
    \raisebox{-9pt}{
      \tiny
      \begin{tabular}{c}
        \color{greenii}
        \bf
        $\mathrm{tmf}$ approximates
        \\
        \color{greenii}
        \bf
        sphere spectrum
        \\
        (Example \ref{TheBoardmanHomomorphismIntmf})
      \end{tabular}
    }
    }
  }{
    \simeq
  }
  \;\;\;\;\;\;\;\;\;\;
  \overset{
    \mathclap{
    \raisebox{9pt}{
      \tiny
      \begin{tabular}{c}
        \color{darkblue}
        \bf
        stable Cohomotopy
        \\
        \color{darkblue}
        \bf
        in degree 4
        \\
        (Example \ref{StableCohomotopy})
      \end{tabular}
    }
    }
  }{
    \mathbb{S}^4\big(X^8\big)
  }
  \;\;\;\;\;\;
  \underset{
    \mathclap{
    \raisebox{-9pt}{
      \tiny
      \begin{tabular}{c}
        \color{greenii}
        \bf
        non-abelian
        \\
        \color{greenii}
        \bf
        enhancement
        \\
        (Example \ref{NonabelianEnhancementOfStableCohomotopy})
      \end{tabular}
    }
    }
  }{
    \squig\squig\squig\squig\rsquigend
  }
  \;\;\;\;\;\;
  \overset{
    \mathclap{
    \raisebox{9pt}{
      \tiny
      \begin{tabular}{c}
        \color{darkblue}
        \bf
        unstable/non-abelian
        \\
        \color{darkblue}
        \bf
        4-Cohomotopy
        \\
        (Example \ref{CohomotopyTheory})
      \end{tabular}
    }
    }
  }{
    \pi^4\big(X^8\big)
  }
  \;\;\;\;\;\;
  \underset{
    \mathclap{
    \raisebox{-9pt}{
      \tiny
      \begin{tabular}{c}
        \color{greenii}
        \bf
        twisting by
        \\
        \color{greenii}
        \bf
        J-homomorphism
        \\
        (Def. \ref{NonabelianTwistedCohomology})
      \end{tabular}
    }
    }
  }{
     \squig\squig\squig\squig\rsquigend
  }
  \;\;\;\;\;\;
  \overset{
    \mathclap{
    \raisebox{9pt}{
      \tiny
      \begin{tabular}{c}
        \color{darkblue}
        \bf
        twisted non-abelian
        \\
        \color{darkblue}
        \bf
        4-Cohomotopy
        \\
        (Example \ref{JTwistedCohomotopyTheory})
      \end{tabular}
    }
    }
  }{
    \pi^{\tau^4}\big(X^8\big)
  }
  \;\;\;\;\;\;
  \underset{
    \mathclap{
    \raisebox{-9pt}{
      \tiny
      \begin{tabular}{c}
        \color{greenii}
        \bf
        lift through
         \\
        \color{greenii}
        \bf
        twisted cohomology operation
        \\
        \color{greenii}
        \bf
        induced by twistor fibration
        \\
        (Example \ref{TwistorialCohomotopy})
      \end{tabular}
    }
    }
  }{
      \squig\squig\squig\squig\squig\squig\squig\squig\rsquigend
  }
  \;\;\;\;\;\;
  \overset{
    \mathclap{
    \raisebox{9pt}{
      \tiny
      \begin{tabular}{c}
        \color{darkblue}
        \bf
        twistorial
        \\
        \color{darkblue}
        \bf
        Cohomotopy
        \\
        (Example \ref{TwistorialCohomotopy})
      \end{tabular}
    }
    }
  }{
    \mathcal{T}^{\tau^4}\big(X^8\big)
  }
  \;\;\;\;\;\;
  \underset{
    \mathclap{
    \raisebox{-9pt}{
      \tiny
      \begin{tabular}{c}
        \color{greenii}
        \bf
        differential
        \\
        \color{greenii}
        \bf
        enhancement
        \\
        (Def. \ref{TwistedDifferentialNonAbelianCohomology})
      \end{tabular}
    }
    }
  }{
     \squig\squig\squig\squig\rsquigend
  }
  \;\;\;\;\;\;
  \overset{
    \mathclap{
    \raisebox{9pt}{
      \tiny
      \begin{tabular}{c}
        \color{darkblue}
        \bf
        differential
        \\
        \color{darkblue}
        \bf
        twistorial
        \\
        \color{darkblue}
        \bf
        Cohomotopy
        \\
        (Example \ref{TwistorDifferentialJTwistedCohomotopy})
      \end{tabular}
    }
    }
  }{
    \widehat{\mathcal{T}}^{\; \tau^4}\big(X^8\big)
  }.
$$

\vspace{-1mm}
The non-abelian character map on twistorial Cohomotopy
has the striking property
(Prop. \ref{ChargeQuantizationInJTwistedCohomotopy},
the proof of which is the content of the companion
physics article \cite[Prop. 3.9]{FSS20a})
that the corresponding non-abelian version of
Dirac's charge quantization \eqref{DiracChargeQuantization}
implies
Ho{\v r}ava-Witten's Green-Schwarz mechanism
in heterotic M-theory
for heterotic line bundles $F_2$ (see \cite[\S 1]{FSS20a})
and other subtle effects expected in
non-perturbative high energy
physics; these are discussed in Remark \ref{SummaryAndHypothesisH} below.

\medskip

\noindent {\bf Quadratic character functions from Whitehead brackets in non-abelian coefficient spaces.}
In summary, 
the crucial appearance of \emph{quadratic functions}
in the character map \eqref{QuadraticFunctions}
is brought about by the intrinsic nature of
(twisted) non-abelian cohomology theory, here
specifically of Cohomotopy theory. These non-linearities
originate in non-trivial Whitehead brackets (Rem. \ref{EquivalentLInfinityStructuresOnWhiteheadProducts})
on the non-abelian coefficient spaces,
such as 
$S^4$ (Exp. \ref{RationalizationOfnSpheres})
and $\mathbb{C}P^3$ (Exp. \ref{FlatTwistedDifferentialFormsWithValuesInTwistorSpace}).
Generally, the non-abelian character map
\eqref{NonAbelianCharacterInIntroduction}
involves also higher monomial terms of any order
(cubic, quartic, ...), originating in higher order
Whitehead brackets on the non-abelian coefficient space
(Rem. \ref{EquivalentLInfinityStructuresOnWhiteheadProducts}).

The desire to conceptually grasp character-like
but quadratic functions
appearing in M-theory had been the original motivation
for developing differential generalized cohomology, in \cite{HopkinsSinger05}. 
Here, 
in differential non-abelian
cohomology, they appear intrinsically.

\medskip

\noindent {\bf Acknowledgements.} We thank
John Lind,
Chris Rogers,
Carlos Simpson,
Danny Stevenson,
and
Mathai Varghese
for comments on an earlier version of this text.
Particular thanks go to an anonymous referee for their detailed comments.

\newpage

%%%%%%%%%%%%%%%%%%%%%%%%%%%%%%%%%%%%%%%%%%%%%%%%%%
\section{Non-abelian cohomology}
 \label{NonabCohomology}
%%%%%%%%%%%%%%%%%%%%%%%%%%%%%%%%%%%%%%%%%%%%%%%%%%%

We make explicit the concept of
general non-abelian cohomology (Def. \ref{NonAbelianCohomology} below)
and of twisted non-abelian cohomology (Def. \ref{NonabelianTwistedCohomology} below),
following \cite{Simpson97}\cite{Simpson99}\cite{Toen02}\cite{SSS09}\cite{NSS12a}\cite{NSS12b}\cite{FSS19b}\cite{SS20b};
and we survey how this concept subsumes essentially every notion of
cohomology known.

\medskip

In the following, we make free use of the basic language of category theory
and homotopy theory
(for joint introduction see \cite{Riehl14}\cite{Richter20}).
For $\mathcal{C}$ a category and $X,A \,\in\, \mathcal{C}$
a pair of its objects, we write
\vspace{-1mm}
\begin{equation}
  \label{HomSet}
  \mathcal{C}
  (
    X
    \,,\,
    A
  )
  \;:=\;
  \mathrm{Hom}_{\mathcal{C}}
  (
    X
    \,,\,
    A
  )
  \;\in\;
  \mathrm{Sets}
\end{equation}

\vspace{-2mm}
\noindent for the set of morphisms from $X$ to $A$. These are, of course,
contravariantly and covariantly functorial in their first and second argument,
respectively:
\vspace{-3mm}
\begin{equation}
  \label{HomFunctor}
  \xymatrix{
    \mathcal{C}
    \ar[rr]^-{
      \mathcal{C}
      (
        X
        \,,\,
        -
      )
    }
    &&
    \mathrm{Sets}
  }
  \,,
  \;\;\;\;\;\;
  \xymatrix{
    \mathcal{C}^{\mathrm{op}}
    \ar[rr]^-{
      \mathcal{C}
      (
        -
        \,,\,
        A
      )
    }
    &&
    \mathrm{Sets}
  }.
\end{equation}

\vspace{-2mm}
\noindent
Basic as this is, contravariant hom-functors are of
paramount interest in the case where $\mathcal{C}$
is the \emph{homotopy category}
$\mathrm{Ho}(\mathbf{C})$ (Def. \ref{HomotopyCategory})
of a model category (Def. \ref{ModelCategories}), such as the
classical homotopy category of topological spaces or, equivalently,
of simplicial sets (Example \ref{TheClassicalHomotopyCategory}).

%%%%%%%%%%%%%%%%%%%%%%%%%%%%%%%%%%%%%%%%%%%%%%
\subsection{Non-abelian cohomology theories}
%%%%%%%%%%%%%%%%%%%%%%%%%%%%%%%%%%%%%%%%%%%%%%

\begin{defn}[Non-abelian cohomology]
  \label{NonAbelianCohomology}
  For $X,A \,\in\, \HomotopyTypes$ (Example \ref{TheClassicalHomotopyCategory})
  we say that their
  hom-set \eqref{HomSet} is the
  \emph{non-abelian cohomology} of $X$
  \emph{with coefficients} in $A$, or the
  \emph{non-abelian $A$-cohomology} of $X$, to be denoted:
  \vspace{-3mm}
  \begin{equation}
    \label{NonabelianCohomologyMappingSpace}
    \overset{
      \raisebox{3pt}{
        \tiny
        \color{darkblue}
        \bf
        \begin{tabular}{c}
          non-abelian
          \\
          cohomology
        \end{tabular}
      }
    }{
    H
    (
      X;
      \,
      A
    )
    }
    \;:=\;
    \HomotopyTypes
    (
      X
      \,,\,
      A
    )
    \;\;\;
    =
    \;\;\;
    \left\{
    \xymatrix{
      X
      \ar@/^1.8pc/[rr]^-{
        \overset{
          \raisebox{1pt}{
            \tiny
            \color{darkblue}
            \bf
            map = cocycle
          }
        }{c}
      }_-{\ }="s"
      \ar@/_1.8pc/[rr]_-{
        \underset{
          \raisebox{-1pt}{
            \tiny
            \color{darkblue}
            \bf
            map = cocycle
          }
        }{
          c'
        }
      }^-{\ }="t"
      &&
      A
      \ar@{=>}|-{
        \mbox{
          \tiny
          \color{orangeii}
          \bf
          \begin{tabular}{c}
            homotopy
            =
            \\
            coboundary
          \end{tabular}
        }
      }
        "s"; "t"
    }
    \right\}_{\!\!\big/\mbox{\tiny{homotopy}}}
      \end{equation}

  \vspace{-3mm}
\noindent
  We also call the contravariant hom-functor \eqref{HomFunctor}
  \vspace{-2mm}
\noindent
  \begin{equation}
    \label{NonAbelianCohomologyTheory}
    H
    (
      -;
      \,
      A
    )
    \;:\;
    \xymatrix{
      \HomotopyTypes
      \ar[r]
      &
      \mathrm{Sets}
    }
  \end{equation}

  \vspace{-2mm}
\noindent
  the non-abelian \emph{$A$-cohomology theory}.
\end{defn}

\begin{example}[Ordinary cohomology]
  \label{OrdinaryCohomology}
  For $n \in \mathbb{N}$ and $A$ a discrete abelian group,
  the ordinary cohomology (e.g. singular cohomology)
  in degree $n$ with coefficients in $A$ is equivalently
  (\cite[p. 243]{Eilenberg40}\cite[p. 520-521]{EML54b},
  review in \cite[\S 19]{Steenrod72}\cite[\S 22]{May99}\cite[\S 7.1, Cor. 12.1.20]{AGP02})
  non-abelian cohomology in the sense of Def. \ref{NonAbelianCohomology}
  \vspace{-2mm}
  \begin{equation}
    \label{OrdinaryCohomologyAsNonAbelianCohomology}
    \overset{
      \mathclap{
      \raisebox{3pt}{
        \tiny
        \color{darkblue}
        \bf
        \begin{tabular}{c}
          ordinary
          \\
          cohomology
        \end{tabular}
      }
      }
    }{
    H^n
    (
      -;
      \,
      A
    )
    }
    \;\;\simeq\;\;
    H
    \big(
      -;
      \,
      K(A,n)
    \big)
  \end{equation}

  \vspace{-2mm}
  \noindent
  with coefficients in an {\it Eilenberg-MacLane space}
  \cite{EML53}\cite{EML54a}:
  \begin{equation}
    \label{EilenbergMacLaneSpaces}
    K(A,n) \in \HomotopyTypes
    \;\;\;\;
    \mbox{\rm such that}
    \;\;\;\;
    \pi_k
    \big(
      K(A,n)
    \big)
    \;=\;
    \left\{
    \begin{array}{lcl}
      A & \vert & k = n
      \\
      0 & \vert & k \neq n
      \,.
    \end{array}
    \right.
  \end{equation}
\end{example}

\begin{example}[Traditional non-abelian cohomology]
  \label{TraditionalNonAbelianCohomology}
  For $G$ a well-behaved\footnote{The technical condition is that
   $G$ be {\it well-pointed}, which means that the inclusion
   $\ast \xrightarrow{e} G$ of the neutral element is a
   {\it closed Hurewicz cofibration}, hence that $(G,\{e\})$ is
   an {\it NDR pair}, see \cite{BaezStevenson09} for pointers and
   \cite{RobertsStevenson12} for details. All Lie groups are well-pointed.
  }
  topological group,
  the traditional non-abelian cohomology $H^1(-; G)$
  classifying $G$-principal bundles, is equivalently
  (\cite[\S 19.3]{Steenrod51}\cite[Thm 1.]{RobertsStevenson12}, review in \cite[\S 5]{Addington07})
  non-abelian cohomology in the general sense of Def. \ref{NonAbelianCohomology}
  \vspace{-2mm}
  \begin{equation}
    \label{IsomorphismBetweenPrincipalBundlesAndMapsToBG}
    \overset{
      \mathclap{
      \raisebox{3pt}{
        \tiny
        \color{darkblue}
        \bf
        \begin{tabular}{c}
          classification of
          \\
          principal bundles
        \end{tabular}
      }
      }
    }{
    H^1
    (
      -;
      G
    )
    }
    \;\simeq\;
    H
    (
      -;
      B G
    )
  \end{equation}

  \vspace{-1mm}
  \noindent
  with coefficients in the
  {\it classifying space} $B G$
  (\cite{Milnor56}\cite{Segal68}\cite{Steenrod68}\cite{Stasheff70},
  review in \cite[\S 1.3]{Kochman96}\cite[\S 23.1]{May99} \cite[\S 8.3]{AGP02}\cite[\S 3.7.1]{NSS12b}).
  The latter may
  be given as the homotopy colimit
  (in the classical model structure of  $\TopologicalSpaces_{\mathrm{Qu}}$, Example \ref{ClassicalModelStructureOnTopologicalSpaces})
  over the nerve of the topological group $G$
  (e.g. \cite[Rem. 2.23]{NSS12a}):
  \begin{equation}
    \label{ClassifyingSpace}
    B G
    \;\simeq\;
    \underset{\longrightarrow}{\mathrm{holim}}
    \left(
      \xymatrix{
        \cdots
        \;
        \ar@<+21pt>@{..}[r]
        \ar@<+14pt>@{..}[r]
        \ar@<+7pt>@{..}[r]
        \ar@<+0pt>@{..}[r]
        \ar@<-7pt>@{..}[r]
        \ar@<-14pt>@{..}[r]
        \ar@<-21pt>@{..}[r]
        &
       \;  G \times G \;
        \ar@<+14pt>[rr]
        \ar@{<-}@<+7pt>[rr]
        \ar@<0pt>[rr]|-{\; (-)\cdot (-)\;}
        \ar@{<-}@<-7pt>[rr]
        \ar@<-14pt>[rr]
        &&
       \; G \;
        \ar@<+7pt>[r]
        \ar@{<-}[r]|-{\;\mathrm{e}\;}
        \ar@<-7pt>[r]
        &
     \;   \ast
      }
    \right).
  \end{equation}
\end{example}

\newpage

\begin{example}[Group cohomology and Characteristic classes]
  \label{GroupCohomology}
  Conversely, the ordinary cohomology (Example \ref{OrdinaryCohomology})
  \emph{of} the classifying space $B G$ \eqref{ClassifyingSpace}
  of a Lie group $G$
  with coefficients in a discrete group $A \in \Groups(\Sets)$
  (such as $A = \mathbb{Z}$)
  is, equivalently:\footnote{
    If $G$ is a topological or Lie group, then the
    appropriate
    ({\it continuous} or {\it smooth}, respectively) group cohomology of $G$ is
    (by \cite[Thm. 4.4.36]{dcct})
    in general not that of the classifying space $B G$, but
    of the {\it universal moduli stack} $\mathbf{B}G$
    (Rem. \ref{StructuredNonAbelianCohomology}).
    However, for {\it discrete} coefficients $A$ this reduces
    (by \cite[Prop. 4.4.35]{dcct}) to the
    cohomology of the geometric realization
    of $\mathbf{B}G$,
    which, at least for Lie groups $G$, coincides
    (by \cite[Prop.  4.4.30]{dcct})
    with that of the classifying space $B G$.
  }

\noindent  {\bf (i)} the group cohomology of $G$;

\noindent   {\bf (ii)} the universal characteristic classes of
    $G$-principal bundles:
  \vspace{-2mm}
  $$
    \overset{
      \raisebox{3pt}{
        \tiny
        \color{darkblue}
        \bf
        \begin{tabular}{c}
          group
          \\
          cohomology
        \end{tabular}
      }
    }{
    H
    \big(
      B G;
      \,
      K(A,n)
    \big)
    }
    \;\simeq\;
    H^n
    (
      B G;
      \,
      A
    )
    \;\simeq\;
    H^{n}_{\mathrm{Grp}}
    (
      G;
      \,
      A
    )\;.
  $$
\end{example}

\begin{example}[Non-abelian cohomology in degree 2]
  \label{NonAbelianCohomologyInDegree2}
  For a well-behaved topological 2-group,
  such as the \emph{string 2-group} $\mathrm{String}(G)$
  (of a connected, simply connected semi-simple Lie group $G$)
  \cite{BCSS07}\cite[Thm. 4.8]{Henriques08}\cite{NSW11},
  the non-abelian cohomology
  $H^1(-;\, \mathrm{String}(G))$
  classifying principal 2-bundles \cite{NikolausWaldorf11}
  with structure 2-group $\mathrm{String}(G)$
  is, equivalently \cite{BaezStevenson09},
  \vspace{-3mm}
  \begin{equation}
    \label{IsomorphismBetweenPrincipal2BundlesAndMapsToBG}
    \overset{
      \mathclap{
      \raisebox{3pt}{
        \tiny
        \color{darkblue}
        \bf
        \begin{tabular}{c}
          classification of
          \\
          String-bundles
        \end{tabular}
      }
      }
    }{
    H^1
    \big(
      -;
      \mathrm{String}(G)
    \big)
    }
    \;\simeq\;
    H
    \big(
      -;
      B \mathrm{String}(G)
    \big)
  \end{equation}

  \vspace{-1mm}
\noindent  non-abelian cohomology in the general sense of Def. \ref{NonAbelianCohomology}
  with coefficients in the classifying space
  $B \mathrm{String}(G)$.
\end{example}

\begin{example}[Non-abelian gerbes]
  \label{NonAbelianGerbes}
  For $G$ a well-behaved topological group,
  a non-abelian \emph{$G$-gerbe} \cite{Giraud71}\cite{Breen09}
  is equivalently \cite[\S 4.4]{NSS12a} a fiber 2-bundle
  %with typical 2-fiber
  %of homotopy type of the classifying space $B G$ \eqref{ClassifyingSpace},
  associated to principal 2-bundles
  with a certain topological structure 2-group $\mathrm{Aut}(\mathbf{B}G)$
  (the automorphism 2-group of the moduli stack of $G$,
  see Rem. \ref{StructuredNonAbelianCohomology}).
  Hence, as in Example \ref{NonAbelianCohomologyInDegree2},
  $G$-gerbes are classified
  by non-abelian cohomology with
  coefficients in $B \mathrm{Aut}(\mathbf{B} G)$ \cite[Cor 4.51]{NSS12a}:
  \vspace{-2mm}
  $$
    \overset{
      \mathclap{
      \raisebox{3pt}{
        \tiny
        \color{darkblue}
        \bf
        \begin{tabular}{c}
          classification of
          \\
          non-abelian gerbes
        \end{tabular}
      }
      }
    }{
      G\mathrm{Gerbes}(X)_{\!/\sim}
    }
    \;\;\simeq\;\;
    H^1
    \big(
      X;
      \,
      \mathrm{Aut}(\mathbf{B} G)
    \big)
    \;\;\simeq\;\;
    H
    \big(
      X;
      \,
      B \mathrm{Aut}(\mathbf{B} G)
    \big)
    \,.
  $$
\end{example}

\begin{example}[Non-abelian cohomology in unbounded degree]
  \label{NonAbelianCohomologyInUnboundedDegree}
  For any $\infty$-group $\mathcal{G}$
  (see \cite[\S 2.2]{NSS12a}\cite[\S 3.5]{NSS12b}),
  the non-abelian cohomology
  $H^1\big(-;\, \mathcal{G}\big)$
  classifying principal $\infty$-bundles
  \cite{Glenn82}\cite{JardineLuo06}\cite{NSS12a}\cite{NSS12b}
  with structure $\infty$-group $\mathcal{G}$
  is, equivalently \cite{Wendt10}\cite{RobertsStevenson12},
  \vspace{-2mm}
  \begin{equation}
    \label{IsomorphismBetweenPrincipalInfinityBundlesAndMapsToBG}
    \overset{
      \mathclap{
      \raisebox{3pt}{
        \tiny
        \color{darkblue}
        \bf
        \begin{tabular}{c}
          classification of
          \\
          non-abelian $\infty$-gerbes
        \end{tabular}
      }
      }
    }{
    H^1
    (
      -;
      \mathcal{G}
    )
    }
    \;\simeq\;
    H
    (
      -;
      B \mathcal{G}
    )
  \end{equation}

  \vspace{-1mm}
\noindent  non-abelian cohomology in the general sense of Def. \ref{NonAbelianCohomology}
  with coefficients in the classifying space
  $B \mathcal{G}$ (see also \cite{Stevenson12}).
\end{example}

Example \ref{NonAbelianCohomologyInUnboundedDegree}
is, in fact,  universal:

\begin{prop}[Connected homotopy types are higher non-abelian classifying
spaces {\cite{May72}\cite[7.2.2.11]{Lurie09}, \cite[Thm. 2.19]{NSS12a}\cite[Thm. 3.30, Cor. 3.34]{NSS12b}}]
\label{ConnectedHomotopyTypesAreHigherNonAbelianClassifyingSpaces}
Every connected homotopy type
$A \in \HomotopyTypes$ \eqref{ClassicalHomotopyCategory}
is the classifying space of a topological group, namely
of its loop group\footnote{
  A priori, the loop group is an $A_\infty$-group,
  for which classifying spaces are defined as in
  \cite[Rem. 2.23]{NSS12a},
  but each such is weakly equivalent to an actual
  topological group, see \cite[Prop. 3.35]{NSS12b}.
} $\Omega A$
\vspace{-1mm}
\begin{equation}
  \label{LoopingDelooping}
  A \;\simeq\; B (\Omega A)
  \;\;\;\;\;
  \in
  \HomotopyTypes
  \,.
\end{equation}
\end{prop}

This allows to make precise the core nature of non-abelian cohomology:
\begin{remark}[From non-abelian to abelian $\infty$-groups]
  \label{NonAbelianInfinityGroups}
  For $A \simeq B G$ \eqref{LoopingDelooping}, the
  $\infty$-group structure on $G$ is reflected by its
  weak homotopy equivalence
  $G \simeq \Omega B G$ with a based loop space.
  \begin{itemize}
  \vspace{-1mm}
  \item
  There is no commutativity of
  composition of loops in a generic loop space,
  and hence this exhibits $G$ as a \emph{non-abelian} $\infty$-group.

  \vspace{-1.5mm}
  \item  But it may happen that $A$ itself is already
  equivalent to a loop space, which by \eqref{LoopingDelooping}
  means that $A \,\simeq\, B\big( B G\big) \,=:\, B^2 G$
  is a \emph{double delooping}.
  In this case $G \,\simeq\, \Omega \big( \Omega A \big) \,=: \,
  \Omega^2 A$ is an \emph{iterated loop space} \cite{May72},
  specifically a \emph{double loop space};
  hence a \emph{braided $\infty$-group}
  (\cite[\S 1]{GarzonMiranda97}\cite[\S 6]{GarzonMiranda00}\cite[Def. 4.28]{FSS12a}).
  By the Eckmann-Hilton argument \cite[Thm. 1.12]{EckmannHilton61}\cite{SchlankYanovski18},
  this implies a first level of commutativity of the group operation in
  $G$. Indeed, in the special case that such $G$ is also 0-truncated
  \eqref{nTruncationOnHomotopyTypes}, it implies that $G$ is
  an ordinary abelian group.

  \vspace{-1.5mm}
  \item Next, it may happen that $A \,\simeq\, B^3 G$
  is a 3-fold delooping, hence that
  $G  \,\simeq\, \Omega^3 A$ is a 3-fold loop space, hence a
  \emph{sylleptic $\infty$-group}
  (where the terminology follows \cite[\S 5]{DayStreet97}\cite[\S 4]{Crans98}
  see \cite[\S 2.2]{GurskiOsorno12} for relation to our context).
  This is one step ``more abelian'' than a braided $\infty$-group.
\end{itemize}

\newpage

 \begin{itemize}
    \item  In the limiting case that $G$ is an $n$-fold loop space
  for any $n \in \mathbb{N}$, hence an \emph{infinite loop space}
  \cite{May77}\cite{Adams78},
  it is as abelian as possible for an $\infty$-group.
  Such {\it symmetric}
  (in the monoidal $\infty$-category theoretic terminology of \cite{Lurie09TQFT})
  or, we may say,
  \emph{abelian $\infty$-groups} are
  the coefficients of abelian cohomology theories, namely of
  generalized cohomology theories in the sense of Whitehead
  (Example \ref{GeneralizedCohomologyAsNonabelianCohomology}).

  \vspace{-1.5mm}
  \item The fewer deloopings an $\infty$-group $G$ admits, the
  ``more non-abelian'' is the cohomology theory represented by $B G$.
  \end{itemize}
\end{remark}

\hspace{0cm}
\def\arraystretch{1.3}
\begin{tabular}{|r|l||c|l|}
  \hline
  \multicolumn{2}{|c||}{
    {\bf Coefficients}
  }
  &
  $
    H
    (
      X;
      \,
      B G
    )
  $
  &
  {\bf Examples}
  \\
  \hline
  \hline
  non-abelian
  $\infty$-group
  &
  $ G \,\simeq\, \Omega^{\phantom{1}} B^{\phantom{1}} G $
  &
  \multirow{4}{*}{
    {\begin{tabular}{c}
      non-abelian
      \\
      cohomology
    \end{tabular}}
  }
  &
  $\pi^{ n }(-)$ ($n$-Cohomotopy,  Example \ref{CohomotopyTheory})
  \\
  \cline{1-2}
  braided $\infty$-group
    &
  $ G \,\simeq\, \Omega^2 B^2 G $
   &
   &
  $\pi^{3}(-)$ (3-Cohomotopy)
  \\
  \cline{1-2}
  sylleptic $\infty$-group
    &
  $ G \,\simeq\, \Omega^3 B^3 G $
    &
    &
  \\
  \cline{1-2}
  \scalebox{.8}{$\vdots$} $\;\;\;\;\;$
  &
  $ G \,\simeq\, \Omega^n B^n G $
  &
  &
  \\
  \hline
  abelian $\infty$-group
    &
  $ G \,\simeq\, \Omega^\infty B^\infty G $
  &
  abelian cohomology
  &
  $E^n(-)$
  (Whitehead-general. cohom., Ex. \ref{GeneralizedCohomologyAsNonabelianCohomology})
  \\
  \hline
\end{tabular}

\vspace{4mm}
\noindent
The most fundamental connected homotopy types are
the $n$-spheres (all other are obtained by gluing $n$-spheres
to each other):

\begin{example}[Cohomotopy theory]
  \label{CohomotopyTheory}
  The non-abelian cohomology theory (Def. \ref{NonAbelianCohomology})
  with coefficients
  in the homotopy types of $n$-spheres is
  (unstable) {\emph{Cohomotopy} theory}
  \cite{Borsuk36}\cite{Spanier49}\cite{Peterson56}\cite{Taylor09}\cite{KMT12}:
 \vspace{-2mm}
  $$
    \overset{
      \mathclap{
      \raisebox{3pt}{
        \tiny
        \color{darkblue}
        \bf
        Cohomotopy
      }
      }
    }{
      \pi^n(-)
    }
    \;=\;
    H
    (
      -;
      S^n
    )
    \;\simeq\;
    H^1
    (
      -;
      \,
      \Omega S^n
    )
    \;\;\;\;\;\;
    \mbox{for $n \in \mathbb{N}_+$}
    \,.
  $$

  \noindent
  {\bf (i)}
  By Prop. \ref{ConnectedHomotopyTypesAreHigherNonAbelianClassifyingSpaces},
  Cohomotopy theory classifies principal $\infty$-bundles
  (Example \ref{NonAbelianCohomologyInUnboundedDegree})
  with structure $\infty$-group of the homotopy type of the $\infty$-group $\Omega S^n$.

  \noindent
  {\bf (ii)}
  By Remark \ref{NonAbelianInfinityGroups}, Cohomotopy theory is
  a \emph{maximally non-abelian} cohomology theory, in that
  $S^n$ does not admit deloopings, for general $n$
  (it admits a single delooping for $n = 3$ and arbitrary
  deloopings for $n = 0, 1$).

\end{example}

\begin{example}[Bundle gerbes]
  \label{BundleGerbes}
  The classifying space \eqref{ClassifyingSpace}
  of the circle group $\mathrm{U}(1)$ is
  an Eilenberg-MacLane space \eqref{EilenbergMacLaneSpaces}
  \vspace{-2mm}
  $$
    B \mathrm{U}(1)
    \;\simeq\;
    K(\mathbb{Z},2)
    \;\;\;\;\;
    \in
    \;
    \HomotopyTypes
    \,.
  $$

  \vspace{-1mm}
  \noindent
  Since $\mathrm{U}(1)$ is abelian, this space
  carries itself the structure of (the homotopy type of)
  a 2-group, and hence has a higher classifying space
    \vspace{-1mm}
    $$
    B^2 \mathrm{U}(1) := B (B \mathrm{U}(1))
    \;\simeq\;
    K(\mathbb{Z},3)
    \;\;\;\;
    \in
    \;
    \HomotopyTypes \;,
  $$

    \vspace{-1mm}
  \noindent
  in the sense of Example \ref{NonAbelianCohomologyInDegree2},
  which is an Eilenberg-MacLane space in one degree higher.
  The higher principal 2-bundles with topological structure 2-group
  $\mathbf{B} \mathrm{U}(1)$ are
  equivalently \cite[Rem. 4.36]{NSS12a}
  known as \emph{bundle gerbes} \cite{Murray96}\cite{SW07}.
  Therefore,
  Example \ref{NonAbelianCohomologyInUnboundedDegree}
  combined with
  Example \ref{OrdinaryCohomology}
  gives the classification of bundle gerbes
  by ordinary
  integral cohomology  in degree 3:
    \vspace{-1mm}
    $$
    \overset{
      \mathclap{
      \raisebox{3pt}{
        \tiny
        \color{darkblue}
        \bf
        \def\arraystretch{.9}
        \begin{tabular}{c}
          classification of
          \\
          bundle gerbes
        \end{tabular}
      }
      }
    }{
    H^1
    \big(
      -;
      \,
      \mathbf{B} \mathrm{U}(1)
    \big)
    }
    \;\simeq\;
    \overset{
      \raisebox{3pt}{
        \tiny
        \color{darkblue}
        \bf
        \def\arraystretch{.9}
        \begin{tabular}{c}
        \end{tabular}
      }
    }{
    H
    \big(
      -;
      \,
      B^2 \mathrm{U}(1)
    \big)
    }
    \;\simeq\;
    H^3(-;\, \mathbb{Z})\;.
  $$
\end{example}

\begin{example}[Higher bundle gerbes]
  \label{HigherBundleGerbes}
  In fact, Prop. \ref{ConnectedHomotopyTypesAreHigherNonAbelianClassifyingSpaces} implies
  that, for all $n \in \mathbb{N}$,
    \vspace{-2mm}
  \begin{equation}
    \label{KZnAndBnU1}
    B^{n+1} \mathrm{U}(1)
    \;:=\;
    B \big( B^{n} \mathrm{U}(1)\big)
    \;\simeq\;
    K(\mathbb{Z},n+2)
    \;\;\;\;
    \in
    \;
    \HomotopyTypes \;,
  \end{equation}

    \vspace{-1mm}
  \noindent
  in the sense of Example \ref{NonAbelianCohomologyInUnboundedDegree}.
  The higher principal bundles with structure ($n+1$)-group
  $\mathbf{B}^{n} \mathrm{U}(1)$
  \cite{Gajer97}\cite[\S 3.2.3]{FiorenzaSchreiberStasheff10}\cite[\S 2.6]{FSS12b}
  are also known as \emph{higher bundle gerbes}
  (for $n = 2$ see \cite{CMW97}\cite{Stevenson01}).
  On these coefficients,
  Example \ref{NonAbelianCohomologyInUnboundedDegree}
  reduces to the classification of higher bundle gerbes by ordinary
  integral cohomology in higher degree:
    \vspace{-2mm}
  $$
    \overset{
      \mathclap{
      \raisebox{3pt}{
        \tiny
        \color{darkblue}
        \bf
        \def\arraystretch{.9}
        \begin{tabular}{c}
          classification of
          \\
          higher bundle gerbes
        \end{tabular}
      }
      }
    }{
    H^1
    \big(
      -;
      \,
      \mathbf{B}^n \mathrm{U}(1)
    \big)
    }
    \;\simeq\;
    \overset{
      \raisebox{3pt}{
        \tiny
        \color{darkblue}
        \bf
        \def\arraystretch{.9}
        \begin{tabular}{c}
        \end{tabular}
      }
    }{
    H
    \big(
      -;
      \,
      B^{n+1} \mathrm{U}(1)
    \big)
    }
    \;\simeq\;
    H^{n+2}(-;\, \mathbb{Z})\;.
  $$
\end{example}

More generally, the special case of Example \ref{NonAbelianCohomologyInUnboundedDegree}
where the coefficient $\infty$-group happens to be
abelian is ``generalized cohomology'' in the
standard sense of algebraic topology
(including cohomology theories such as K-theory, elliptic cohomology,
stable Cobordism theory, stable Cohomotopy theory, etc.):
\begin{example}[Whitehead-generalized cohomology]
  \label{GeneralizedCohomologyAsNonabelianCohomology}
  For $E$ a generalized cohomology theory in the traditional sense of
  \cite{Whitehead62} (review in \cite{Adams74}\cite{Adams78}\cite{TamakiKono06}),
  Brown's representability theorem (\cite[\S III.6]{Adams74}\cite[\S 3.4]{Kochman96})
  says that there is a
  \emph{spectrum} (``$\Omega$-spectrum'', Example \ref{ModelStructureOnSequentialSpectra})
  of pointed homotopy types
  \vspace{-2mm}
  \begin{equation}
    \label{Spectrum}
    \Big\{
    E_n \in \PointedHomotopyTypes
    \,,
    \xymatrix{
      E_n \ar[r]^-{ \widetilde \sigma_n }_-{\simeq} & \Omega E_{n+1}
    }
    \Big\}_{n \in \mathbb{N}}
  \end{equation}

  \vspace{-1mm}
  \noindent
  such that the generalized $E$-cohomology in degree $n$ is
  equivalently non-abelian cohomology theory in the sense of
  Def. \ref{NonAbelianCohomology} with coefficients in $E_n$:
  \vspace{-2mm}
  \begin{equation}
    \label{GeneralizedECohomologyAsNonabelianCohomology}
    \overset{
      \mathclap{
      \raisebox{3pt}{
        \tiny
        \color{darkblue}
        \bf
        \def\arraystretch{.9}
        \begin{tabular}{c}
          generalized
          \\
          cohomology
        \end{tabular}
      }
      }
    }{
      E^n(-)
    }
    \;\;\simeq\;\;
    H
    (
      -;
      \,
      E_n
    )\;.
  \end{equation}
  Often one is interested in the special case that the representing spectrum
  carries the structure of an
  $E_\infty$-ring (review in \cite{BakerRichter04}\cite{Richter21}), in which case
  $E^\bullet(-)$ is a {\it multiplicative cohomology theory} (e.g. \cite[\S 2.6]{TamakiKono06})
  where, in particular,
  the generalized cohomology groups \eqref{GeneralizedECohomologyAsNonabelianCohomology}
  inherit ordinary ring-structure.
\end{example}

\begin{example}[Topological K-theory]
  \label{ComplexTopologicalKtheory}
  The classifying space \eqref{Spectrum}
  representing complex K-cohomology theory $\mathrm{KU}$
  \cite[\S 2]{AtiyahHirzebruch59} (review in \cite{Atiyah67})
  in degree 0 is
  \cite[\S 1.3]{AtiyahHirzebruch61}:
    \vspace{-2mm}
  \begin{equation}
    \label{ClassifyingSpaceForComplexTopologicalKTheory}
    \mathrm{KU}_0
    \;\simeq\;
    \mathbb{Z} \times B \mathrm{U}
    \,,
  \end{equation}

    \vspace{-2mm}
  \noindent
  where
    \vspace{-1mm}
  \begin{equation}
    \label{ClassifyingSpaceOfInfiniteUnitaryGroup}
    B \mathrm{U} \;:=\;
    \underset{\underset{n}{\longrightarrow}}{\mathrm{lim}}
    \,
    B \mathrm{U}(n)
  \end{equation}

  \vspace{-1mm}
  \noindent
  is the classifying space \eqref{ClassifyingSpace}
  for the infinite unitary group (e.g. \cite{EspinozaUribe14}).
  Hence for the case of complex K-theory, Example \ref{GeneralizedCohomologyAsNonabelianCohomology}
  says that:
    \vspace{-2mm}
  $$
    \overset{
      \mathclap{
      \raisebox{3pt}{
        \tiny
        \color{darkblue}
        \bf
        \def\arraystretch{.9}
        \begin{tabular}{c}
          topological
          \\
          K-theory
        \end{tabular}
      }
      }
    }{
      \mathrm{KU}^0(-)
    }
    \;\simeq\;
    H
    (
      -;
      \,
      \mathbb{Z} \times B \mathrm{U}
    )
    \,.
  $$
\end{example}
\begin{example}[Iterated K-theory]
  \label{IteratedKTheory}
  Given an $E_\infty$-ring spectrum $R$ (Ex. \ref{GeneralizedCohomologyAsNonabelianCohomology}),
  one may form its \emph{algebraic K-theory spectrum}
  $K(R)$ \cite[\S VI]{EKMM97}\cite[\S 9.5]{BGT10}\cite{Lurie14}
  and hence the corresponding
  generalized cohomology theory
  (Example \ref{GeneralizedCohomologyAsNonabelianCohomology}).
  Much like complex topological K-theory
  (Example \ref{ComplexTopologicalKtheory})
  is the K-theory of topological $\mathbb{C}$-module bundles,
  so $K(R)$-cohomology theory is the
  K-theory of $R$-module $\infty$-bundles \cite{Lind13}.
  Specifically, for $R = \mathrm{ku}$ the connective spectrum of
  topological K-theory, its algebraic K-theory $K(\mathrm{ku})$
  \cite{Ausoni09}\cite{AusoniRognes02}\cite{AusoniRognes07}
  has been argued to be the K-theory of certain categorified complex vector
  bundles \cite{BDR03} \cite{BDRR09}.

  Moreover, if $R$ is connective,
  then $K(R)$ itself carries the structure of a
  connective $E_\infty$-ring spectrum
  (by \cite[Thm. 1]{SchwaenzlVogt94}\cite[Thm. 6.1]{EKMM97}),
  so that the
  construction may be iterated
  to yield \emph{iterated algebraic K-theories} \cite{Rognes14}
  $K^{\circ_2}(R) := K(K(R))$, $K^{\circ_3}(R) := K(K(K(R)))$, et cetera.

  For
  $R = \mathrm{ku}$, this generalizes the above
  ``form of elliptic cohomology'' $K(\mathrm{ku})$
  to higher degrees \cite{LindSatiWesterland16}.
  By Example \ref{GeneralizedCohomologyAsNonabelianCohomology}, we
  will regard these
  (connective) iterated algebraic K-theories
  $K^{\circ_n}(\mathrm{ku})$
  of the complex
  topological K-theory spectrum
  as examples of non-abelian cohomology theories (that happen to be abelian):
  \vspace{-1mm}
  $$
    \overset{
      \mathclap{
      \raisebox{3pt}{
        \tiny
        \color{darkblue}
        \bf
        \def\arraystretch{.9}
        \begin{tabular}{c}
          iterated
          K-theory
        \end{tabular}
      }
      }
    }{
      K^{\circ_n}(\mathrm{ku})^0(-)
    }
    \;\;\simeq\;\;
    H
    \big(
      -;
      \,
        K^{\circ_n}(\mathrm{ku})_0
    \big).
  $$
\end{example}

\begin{example}[Stable Cohomotopy]
  \label{StableCohomotopy}
  The generalized cohomology theory
  (Example \ref{GeneralizedCohomologyAsNonabelianCohomology})
  represented
  by the suspension spectra (Example \ref{DerivedStabilizationAdjunction})
  of $n$-spheres is
  called \emph{stable Cohomotopy theory} (e.g. \cite{Stretch81}\cite{Nowak03})
  or \emph{stable framed Cobordism theory}:
  \begin{equation}
    \mathbb{S}^n(-)
    \;\;=\;\;
    H
    \big(
      -;
      \,
      (\Sigma^\infty S^n)_0
    \big)
    \,.
  \end{equation}
\end{example}

\medskip

\noindent {\bf Non-abelian cohomology operations.}

\begin{defn}[Non-abelian cohomology operation]
  \label{NonAbelianCohomologyOperations}
  For $A_1, A_2 \in \HomotopyTypes$ (Example \ref{TheClassicalHomotopyCategory}), we say that a
  natural transformation in non-abelian cohomology
  (Def. \ref{NonAbelianCohomology})
  from $A_1$-cohomology theory to
  $A_2$-cohomology theory \eqref{NonAbelianCohomologyTheory}
  is a (non-abelian) \emph{cohomology operation}
  \vspace{-2mm}
  \begin{equation}
    \label{NaturalTransformationOfNonAbelianCohomologyTheories}
    \phi_\ast
    \;:\;
    \xymatrix{
      H
      (
        -;
        \,
        A_1
      )
      \ar[r]
      &
      H
      (
        -;
        \,
        A_2
      )
    }
    .
  \end{equation}

  \vspace{-2mm}
  \noindent
  By the Yoneda lemma,
  these are in bijective correspondence
  to morphisms of coefficients
    \vspace{-2mm}
     \begin{equation}
    \xymatrix{
      A_1
      \ar[r]^-{ \phi }
      &
      A_2
    }
    \;\;
    \in
    \HomotopyTypes
  \end{equation}

    \vspace{-1mm}
   \noindent
  via the covariant functoriality of the hom-sets \eqref{HomFunctor}:
  \vspace{-2mm}
  \begin{equation}
    \label{InducedCohomologyOperation}
    \phi_\ast
    \;=\;
    H
    (
      -;
      \,
      \phi
    )
    \;:=\;
    \HomotopyTypes
    (
      -;
      \,
      \phi
    )
    \,.
  \end{equation}
\end{defn}

\begin{example}[Cohomology of coefficient spaces parametrizes cohomology operations]
  \label{CohomologyOfCoefficientsIsCohomologyOperations}
  By the Yoneda lemma \eqref{InducedCohomologyOperation}
  in $\HomotopyTypes$ (Example \ref{TheClassicalHomotopyCategory}),
  the set of all cohomology operations (Def. \ref{NonAbelianCohomologyOperations})
  from $A_1$-cohomology theory to $A_2$-cohomology theory
  \eqref{NaturalTransformationOfNonAbelianCohomologyTheories}
  coincides with
  the non-abelian $A_2$-cohomology (Def. \ref{NonAbelianCohomology}) \emph{of}
  the coefficients $A_1$:
    \vspace{-2mm}
  \begin{equation}
    \label{CohomologyOperationParametrizedByCohomologyOfCoeffcients}
    \xymatrix@R=10pt{
    \overset{
      \mathclap{
      \raisebox{3pt}{
        \tiny
        \color{darkblue}
        \bf   \hspace{1.8cm}
        \def\arraystretch{.9}
        \begin{tabular}{c}
          non-abelian $A_2$-cohomology of $A_1$
          \\
          acting as cohomology operations
        \end{tabular}
      }
      }
    }{
    H
    (
      A_1;
      \,
      A_2
    )
    }
    \times
    H
    (
      -;
      A_1
    )
    \ar[rr]^-{ (-) \circ (-) }
    &&
    H
    (
      -;
      \,
      A_2
    )
    }
  \end{equation}

    \vspace{-1mm}
  \noindent
  acting by composition composition in
  $\HomotopyTypes$.
\end{example}

\begin{example}[Cohomology operations in ordinary cohomology]
  In specialization to Example \ref{OrdinaryCohomology}
  the non-abelian cohomology operations according to
  Def. \ref{NonAbelianCohomologyOperations}
  reduce to the classical
  cohomology operations in ordinary cohomology
  \cite{Steenrod72}\cite{MosherTangora08}
  (review in \cite[\S 22.5]{May99}),
  such as Steenrod operations \cite{Steenrod47}\cite{SteenrodEpstein62}
  (review in \cite[\S 2.5]{Kochman96}).
  These operations admit refinements, involving rational/real form data, to differential cohomology operations
  \cite{GS-ops}.
\end{example}

\begin{example}[Cohomology operations in generalized cohomology]

  In specialization to Example \ref{GeneralizedCohomologyAsNonabelianCohomology},
  the non-abelian cohomology operations according to
  Def. \ref{NonAbelianCohomologyOperations}
  on a Whitehead-generalized cohomology theory $E^\bullet(-)$
  regarded as a system of non-abelian cohomology theories
  $\big\{ E^n(-) \big\}$ reduce to the traditional notion of
  {\it unstable} cohomology operations
  on generalized cohomology theories \cite{BJW95},
  such as the Adams operations in K-theory \cite{Adams62}
  (review in \cite[\S 10]{AGP02})
  or the Quillen operations in stable Cobordism theory
  (review in \cite[\S 4,5]{Kochman96}). For differential refinements see \cite{GS-KO}.
\end{example}

\begin{example}[Characteristic classes of principal $\infty$-bundles]
  \label{CharacteristicClassesOfPrincipalBundles}
  For $G$ a topological group,
  the ordinary group cohomology of $G$ (Example \ref{GroupCohomology})
  parametrizes, via Example  \ref{CohomologyOfCoefficientsIsCohomologyOperations},
  the cohomology operations from
  non-abelian cohomology
  classifying $G$-principal bundles
  (Examples \ref{TraditionalNonAbelianCohomology},
  \ref{NonAbelianCohomologyInDegree2}, \ref{NonAbelianCohomologyInUnboundedDegree})
  to ordinary cohomology of the base space (Example \ref{OrdinaryCohomology}):
  \vspace{-2mm}
  \begin{equation}
    \label{PullbackOfUniversalCharacteristicClasses}
    \xymatrix{
      \overset{
        \mathclap{
        \raisebox{3pt}{
          \tiny
          \color{darkblue}
          \bf
          \def\arraystretch{.9}
          \begin{tabular}{c}
            group
            \\
            cohomology
          \end{tabular}
        }
        }
      }{
      H^n_{\mathrm{Grp}}
      \big(
        G;
        \,
        A
      \big)
      }
      \;\times\;
      \overset{
        \mathclap{
        \raisebox{3pt}{
          \tiny
          \color{darkblue}
          \bf
          \def\arraystretch{.9}
          \begin{tabular}{c}
            $G$-principal
            \\
            bundles
          \end{tabular}
        }
        }
      }{
        H^1(-;\, G)
      }
      \ar[rr]^-{
        \mbox{
          \tiny
          \color{darkblue}
          \bf
        \def\arraystretch{.9}
          \begin{tabular}{c}
            characteristic
            \\
            classes
          \end{tabular}
        }
      }_-{
        \mbox{
          \tiny
          \color{darkblue}
          \bf
          \eqref{CohomologyOperationParametrizedByCohomologyOfCoeffcients}
        }
      }
      &&
      \overset{
        \mathclap{
        \raisebox{3pt}{
          \tiny
          \color{darkblue}
          \bf
          \def\arraystretch{.9}
          \begin{tabular}{c}
            ordinary
            \\
            cohomology
          \end{tabular}
        }
        }
      }{
        H^n(-; A)\;.
      }
    }
  \end{equation}

  \vspace{-1mm}
  \noindent
  This is the assignment of \emph{characteristic classes}
  to principal bundles (principal $\infty$-bundles).
  In the case when $A = \mathbb{R}$, this is equivalently
  the \emph{Chern-Weil homomorphism}, by Chern's
  fundamental theorem (see Remark \ref{ChernWeilTheoryAndItsFundamentalTheorem}
  and Theorem \ref{NonAbelianChernDoldSubsumesChernWeil} below).
\end{example}

\begin{example}[Rationalization cohomology operation]
  For fairly general non-abelian coefficients $A$
  (see
  Def. \ref{Rationalization},
  Def. \ref{RationalizationOfCoefficientsInNonabelianCohomology}
  for details),
  their \emph{rationalization}\footnote{
    To make the connection to differential cohomology,
    we consider rationalization over the {\it real} numbers;
    see Remark \ref{RationalHomotopyTheoryOverTheRealNumbers} below.
  }
  $\!\!\xymatrix@C=30pt{A \ar[r]|-{\;\eta^{\mathbb{R}}_A\;} & L_{\mathbb{R}}}\!\! A$
  (Def. \ref{Rationalization}, \ref{Lk} below) induces
  a cohomology operation (Def. \ref{NonAbelianCohomologyOperations})
  from non-abelian $A$-cohomology theory (Def. \ref{NonAbelianCohomology})
  to \emph{non-abelian real cohomology} (Def. \ref{NonAbelianRealCohomology} below):
  \vspace{-2mm}
  \begin{equation}
    \label{RationalizationCohomologyOperation}
    \xymatrix{
      \overset{
        \mathclap{
        \raisebox{3pt}{
          \tiny
          \color{darkblue}
          \bf
          \def\arraystretch{.9}
          \begin{tabular}{c}
            non-abelian
            \\
            cohomology
          \end{tabular}
        }
        }
      }{
      H
      (
        -;
        A
      )
      }
      \ar[rr]^-{   (\eta^{\mathbb{R}}_A)_\ast }_-{
              \raisebox{-3pt}{
            \tiny
            \color{darkblue}
            \bf
            rationalization
          }
          }
         &&
      \overset{
        \mathclap{
        \raisebox{3pt}{
          \tiny
          \color{darkblue}
          \bf
          \def\arraystretch{.9}
          \begin{tabular}{c}
            non-abelian
            \\
            real cohomology
          \end{tabular}
        }
        }
      }{
      H
      \big(
        -;
        L_{\mathbb{R}}A
      \big)
      }
      \,.
    }
  \end{equation}
\end{example}

\begin{remark}[Rationalization as character map]
  Up to composition with an equivalence provided by the
  non-abelian de Rham theorem (Theorem \ref{NonAbelianDeRhamTheorem} below),
  which serves to bring the right hand side of \eqref{RationalizationCohomologyOperation}
  into neat minimal form, this rationalization cohomology
  operation is the \emph{character map in non-abelian cohomology}
  (Def. \ref{NonAbelianChernDoldCharacter} below).
\end{remark}

\begin{example}[Stabilization cohomology operation]
  \label{StabilizationCohomologyOperation}
  For $A \in \HomotopyTypes$, the non-abelian cohomology operation
  (Def. \ref{NonAbelianCohomologyOperations})
  induced \eqref{InducedCohomologyOperation}
  by the unit of the derived stabilization adjunction
  (Example \ref{DerivedStabilizationAdjunction})
  goes from non-abelian $A$-cohomology theory
  (Def. \ref{NonAbelianCohomology})
  to (abelian)
  generalized cohomology theory
  (Example \ref{GeneralizedCohomologyAsNonabelianCohomology})
  represented by the 0th component space of the suspension spectrum of $A$:
   \vspace{-3mm}
  $$
    \xymatrix{
      \overset{
        \mathclap{
        \raisebox{3pt}{
          \tiny
          \color{darkblue}
          \bf
          \def\arraystretch{.9}
          \begin{tabular}{c}
            non-abelian
            \\
            $A$-cohomology
          \end{tabular}
        }
        }
      }{
      H
      \big(
        -;
        \,
        A
      \big)
      }
      \ar[rr]_-{
        \mbox{
          \tiny
          \color{darkblue}
          \bf
          stabilization
        }
      }
      &&
      \overset{
        \mathclap{
        \raisebox{3pt}{
          \tiny
          \color{darkblue}
          \bf
          \def\arraystretch{.9}
          \begin{tabular}{c}
            generalized
            \\
            $\Sigma^\infty A$-cohomology
          \end{tabular}
        }
        }
      }{
      H
      \big(
        -;
        \,
        (\LeftDerived\Sigma^\infty A)_0
      \big)
      }
      \,.
    }
  $$

   \vspace{-2mm}
  \noindent
  Hence a lift through this operation is an
  enhancement of generalized cohomology to non-abelian cohomology.
\end{example}

\begin{example}[Non-abelian enhancement of stable Cohomotopy]
  \label{NonabelianEnhancementOfStableCohomotopy}
  The canonical non-abelian enhancement
  (in the sense of Example \ref{StabilizationCohomologyOperation})
  of stable Cohomotopy (Example \ref{StableCohomotopy})
  is actual Cohomotopy theory (Example \ref{CohomotopyTheory}):
  \vspace{-3mm}
  $$
    \xymatrix{
      \overset{
        \mathclap{
        \raisebox{3pt}{
          \tiny
          \color{darkblue}
          \bf
          \def\arraystretch{.9}
          \begin{tabular}{c}
            Cohomotopy
          \end{tabular}
        }
        }
      }{
        \pi^n(-)
      }
      \ar[rr]_-{
        \mbox{
          \tiny
          \color{darkblue}
          \bf
          stabilization
        }
      }
      &&
      \overset{
        \mathclap{
        \raisebox{3pt}{
          \tiny
          \color{darkblue}
          \bf
          \def\arraystretch{.9}
          \begin{tabular}{c}
            stable
            \\
            Cohomotopy
          \end{tabular}
        }
        }
      }{
        \mathbb{S}^n(-)
      }
      \,.
    }
  $$
\end{example}

\begin{example}[Hurewicz homomorphism and Hopf degree theorem]
  \label{HurewiczHomomorphism}
  By definition of Eilenberg-MacLane spaces \eqref{EilenbergMacLaneSpaces}
  there is, for $n \in \mathbb{N}$, a canonical map
  \vspace{-2mm}
  $$
    \xymatrix{
      S^n
      \ar[r]^-{ e^{(n)} }
      &
      K(\mathbb{Z},n)
      \;\;\;\;\;\;
      \in
      \HomotopyTypes
      \,,
    }
  $$

  \vspace{-2mm}
  \noindent
  which represents the element
  $1 \in \mathbb{Z} \simeq \pi_n\big(K(\mathbb{Z},n)\big)$.
  The non-abelian cohomology operation (Def. \ref{NonAbelianCohomologyOperations})
  induced by this, from degree $n$ Cohomotopy (Example \ref{CohomotopyTheory})
  to degree $n$ ordinary cohomology (Example \ref{OrdinaryCohomology})
     \vspace{-2mm}
  $$
    \xymatrix{
      \pi^n(-)
      \ar[r]^-{ e^{(n)}_\ast }
      &
      H^n(-;\mathbb{Z})
    }
  $$

  \vspace{-2mm}
  \noindent
  is the cohomological version of the
  Hurewicz homomorphism. The \emph{Hopf degree theorem}
  (e.g. \cite[\S IX (5.8)]{Kosinski93})
  is the statement that
  the non-abelian cohomology operation
  $e^{(n)}_\ast$ becomes an isomorphism
  on connected, orientable closed manifolds of dimension $n$.
  These maps, together with their differential refinements,
  are analyzed in more detail via Postnikov towers
  in \cite{GS-Postnikov}.
\end{example}

\medskip

\noindent
{\bf Structured non-abelian cohomology.}
\begin{remark}[Structured non-abelian cohomology]
  \label{StructuredNonAbelianCohomology}
  More generally, it makes sense to consider the analog
  of Def. \ref{NonAbelianCohomology} for the
  homotopy category
  $\mathrm{Ho}(
    \mathbf{H}
  )$ of a model category which is a \emph{homotopy topos}
  \cite{ToenVezzosi05}\cite{Lurie09}\cite{Rezk10}.

\vspace{1mm}
  \noindent {\bf (i)} This yields \emph{structured} non-abelian cohomology
  \cite{Simpson97}\cite{Simpson99}\cite{Toen02}\cite{SSS09}\cite{NSS12a}\cite{NSS12b}\cite{dcct}\cite{FSS19b} \cite{SS20b}:
  \vspace{-2mm}
  $$
    \overset{
      \mathclap{
      \raisebox{3pt}{
        \tiny
        \color{darkblue}
        \bf
        \def\arraystretch{.9}
        \begin{tabular}{c}
          structured
          \\
          non-abelian cohomology
        \end{tabular}
      }
      }
    }{
    H
    \big(
      \mathcal{X};
      \,
      \mathbf{A}
    \big)
    }
    \;\;\;:=\;
    \overset{
      \mathclap{
      \raisebox{3pt}{
        \tiny
        \color{darkblue}
        \bf
        \def\arraystretch{.9}
        \begin{tabular}{c}
          homotopy topos
        \end{tabular}
      }
      }
    }{
    \mathrm{Ho}
    (
      \mathbf{H}
    )
    }
    \big(
      \underset{
        \mathclap{
        \raisebox{3pt}{
          \tiny
          \color{darkblue}
          \bf
          $\infty$-stacks
        }
        }
      }{
        \mathcal{X}
        \,,\,
        \mathbf{A}
      }
    \big)
    \,,
  $$

   \vspace{-1mm}
   \noindent
  including the stacky non-abelian cohomology originally considered
  in \cite{Giraud71}\cite{Breen90} (``gerbes'', see \cite[\S 4.4]{NSS12a}), and, more generally,
  {\it differential-, {\'e}tale-, and equivariant}-
  nonabelian cohomology theories (see \cite[p. 6]{SS20b})
  based on $\infty$-stacks.

\noindent {\bf (ii)} In good cases (cohesive homotopy toposes \cite{dcct}\cite[\S 3.1]{SS20b}),
  the homotopy topos $\mathrm{Ho}(\mathbf{H})$ comes equipped with a
  \emph{shape} operation down to the classical homotopy category
  (Example \ref{TheClassicalHomotopyCategory}):
    \vspace{-3mm}
  \begin{equation}
    \label{FromStructuredNonAbelianCohmologyToPlain}
    \xymatrix@R=-15pt{
      \overset{
        \mathclap{
        \raisebox{3pt}{
          \tiny
          \color{darkblue}
          \bf
          homotopy topos
        }
        }
      }{
        \mathrm{Ho}
        (
          \mathbf{H}
        )
      }
      \ar[rr]^-{ \mathrm{Shp} }
      &&
      \overset{
        \mathclap{
        \raisebox{3pt}{
          \tiny
          \color{darkblue}
          \bf
          \def\arraystretch{.9}
          \begin{tabular}{c}
            classical
            homotopy category
          \end{tabular}
        }
        }
      }{
        \HomotopyTypes
      }
      \\
      \underset{
        \mathclap{
        \raisebox{-3pt}{
          \tiny
          \color{darkblue}
          \bf
          \def\arraystretch{.9}
          \begin{tabular}{c}
            structured
            \\
            non-abelian cohomology
          \end{tabular}
        }
        }
      }{
        H
        \big(
          \mathcal{X};
          \,
          \mathbf{A}
        \big)
      }
      \ar@{|->}[rr]
      &&
      \underset{
        \mathclap{
        \raisebox{-3pt}{
          \tiny
          \color{darkblue}
          \bf
          \def\arraystretch{.9}
          \begin{tabular}{c}
            plain
            \\
            non-abelian cohomology
          \end{tabular}
        }
        }
      }{
        H
        \big(
          \mathrm{Shp}(\mathcal{X});
          \,
          \mathrm{Shp}(\mathbf{A})
        \big)
      }
    }
  \end{equation}

    \vspace{-3mm}
   \noindent
  which takes, for well-behaved group $\infty$-stacks $G$,
  the \emph{classifying stacks}
  $\mathbf{B}G$ of $G$-principal bundles to the traditional
  classifying spaces $B G \;\simeq\; \mathrm{Shp}(\mathbf{B}G)$
  of underlying topological groups \eqref{ClassifyingSpace}.
  This gives a forgetful functor from structured non-abelian
  cohomology to plain non-abelian cohomology in the sense of
  Def. \ref{NonAbelianCohomology}.  A classical example is the
  map from non-abelian {\v C}ech cohomology with coefficients in a well-behaved
  group $G$ to homotopy classes of maps to the classifying space of $G$,
  in which case this comparison map is a bijection
  (Example \ref{TraditionalNonAbelianCohomology}).

 \noindent {\bf (iii)}  All constructions on non-abelian cohomology have their structured
  analogues, for instance non-abelian cohomology operations
  (Def. \ref{NonAbelianCohomologyOperations}) in
  structured cohomology

  \vspace{-.4cm}
  \begin{equation}
    \label{InducedStructuredCohomologyOperation}
    \xymatrix{
      H
      \big(
        \mathcal{X}
        ;
        \,
        \mathbf{A}_1
      \big)
      \ar[r]^-{ \phi_\ast }
      &
      H
      \big(
        \mathcal{X}
        ;
        \,
        \mathbf{A}_2
      \big)
    }
  \end{equation}
  \vspace{-.4cm}

  \noindent
  are induced by postcomposition with morphisms
  $
      \mathbf{A}_1
      \xrightarrow{\; \phi \; }
            \mathbf{A}_2
      $
  of coefficient stacks.
\end{remark}

  Ultimately, one is interested in working
  with structured non-abelian cohomology on the left of \eqref{FromStructuredNonAbelianCohmologyToPlain}.
  However, since this is rich and intricate, it behooves us to
  study its projection into
  plain non-abelian cohomology on the right of \eqref{FromStructuredNonAbelianCohmologyToPlain}.
  This is what we are mainly concerned with here.
  But we provide in \cref{NonabelianDifferentialCohomology}
  a brief discussion of non-abelian differential cohomology
  on smooth $\infty$-stacks.

%%%%%%%%%%%%%%%%%
 \subsection{Twisted non-abelian cohomology.}
%%%%%%%%%%%%%%%%

For $\mathcal{C}$ any category and $B \in \mathcal{C}$
any object, there is the \emph{slice category}
$\mathcal{C}^{/X}$,
whose objects are morphisms in $\mathcal{C}$ to $X$ and whose
morphisms are commuting triangles over $X$ in $\mathcal{C}$.
Basic as this is, hom-sets in the \emph{homotopy category}
$\mathrm{Ho}(\mathbf{C}^{/B})$ (Def. \ref{HomotopyCategory})
of a slice model category $\mathbf{C}^{/B}$ (Example \ref{SliceModelCategory})
are of paramount interest:

\vspace{1mm}
The slicing imposes \emph{twisting} on the corresponding
non-abelian cohomology (Def. \ref{NonAbelianCohomology}),
in that the slicing of the domain space serves as a twist,
the slicing of the coefficient space as a local coefficient bundle,
and the slice morphisms as twisted cocycles.

\vspace{1mm}
\begin{prop}[$\infty$-Actions on homotopy types {\cite{DDK80}\cite[\S 5]{Prezma10}\cite[\S 4]{NSS12a}\cite{Sharma15}\cite[\S 2.2]{SS20b}}]
  \label{GActionsAsFibrations}
  $\,$

\noindent  For any $A \in \HomotopyTypes$ (Ex. \ref{TheClassicalHomotopyCategory})
  and $G$ a topological group,
  homotopy-coherent actions of $G$ on $A$ are
  equivalent to fibrations $\rho$
  with homotopy fiber $A$ (Def. \ref{HomotopyFibers})
  over the classifying space $B G$ \eqref{ClassifyingSpace}
  \vspace{-2mm}
  \begin{equation}
    \label{AutomorphismActionAsHomotopyFiberSequence}
    \xymatrix@R=1.5em{
      A \ar[r] &
      A \!\sslash\! G
      \ar[d]^-{ \rho }
      \\
      &
      B G
      \,.
    }
  \end{equation}

  \vspace{-3mm}
  \noindent
  Here
  $$
    A \!\sslash\! G
    \;\;
    \simeq
    \;\;
    \big(
      A \times E G
    \big)_{
      \!\!/_{\mathrm{diag}G}
    }
  $$

    \vspace{-1mm}
  \noindent
  is the homotopy quotient
  (Borel construction) of the action.
\end{prop}

\begin{defn}[Twisted non-abelian cohomology {\cite[\S 4]{NSS12a}\cite[(10)]{FSS19b}\cite[Rem. 2.94]{SS20b}}]
  \label{NonabelianTwistedCohomology}
 $\,$

  \noindent
  For $X, A \in \HomotopyTypes$ (Def. \ref{TheClassicalHomotopyCategory})
  we say:

  \noindent {\bf (i)}
  A \emph{local coefficient bundle} for twisted $A$-cohomology
  is an $A$-fibration $\rho$ over
  a classifying space
  $B G$ \eqref{ClassifyingSpace}
  as in Prop. \ref{GActionsAsFibrations}:
    \vspace{-2mm}
    \begin{equation}
    \label{LocalCoefficientBundle}
    \raisebox{12pt}{
    \xymatrix@C=11pt@R=-2pt{
      A \ar[rr]
      &&
      A \!\sslash\! G
      \ar[dd]^-{ \rho }
      \\
      &
      \mathclap{
      \mbox{
        \tiny
        \color{darkblue}
        \bf
        \def\arraystretch{.9}
        \begin{tabular}{c}
          local coefficient
          \\
          bundle
        \end{tabular}
      }
      }
      \\
      &&
      B G
      \,.
    }
    }
  \end{equation}

  \vspace{-1mm}
    \noindent {\bf (ii)}
  A \emph{twist} for non-abelian $A$-cohomology theory on $X$
  with local coefficient bundle $\rho$ over $B G$ is a map
  \vspace{-2mm}
  \begin{equation}
    \label{ATwist}
    \xymatrix{
      X
      \ar[r]^-{ \tau }
      &
      B G
    }
    \;\;\;\;\;
    \in
    \;
    \HomotopyTypes\;.
  \end{equation}

  \vspace{-2mm}
  \noindent {\bf (iii)}
  The {\it non-abelian $\tau$-twisted $A$-cohomology} of $X$
  with local coefficients $\rho$ is
  the hom-set from $\tau$ \eqref{ATwist} to $\rho$
   \eqref{AutomorphismActionAsHomotopyFiberSequence}
      \vspace{-2mm}
  \begin{equation}
    \label{TwistedCohomologyHom}
    \overset{
      \mathclap{
      \raisebox{3pt}{
        \tiny
        \color{darkblue}
        \bf
        \def\arraystretch{.9}
        \begin{tabular}{c}
          twisted
          \\
          non-abelian
          \\
          cohomology
        \end{tabular}
      }
      }
    }{
    H^\tau
    (
      X;
      \,
      A
    )
    }
    \;:=\;
    \TwistedHomotopyTypes
    \big(
      \tau
      \,,\,
      \rho
    \big)
    \;
    =
    \;
    \left\{\!\!\!\!\!\!\!\!
    \raisebox{18pt}{
    \xymatrix@C=20pt{
      \;\;\;X\;\;
      \ar@{-->}[rr]^-{
        \overset{
          \mathclap{
          \raisebox{3pt}{
            \tiny
            \color{darkblue}
            \bf
            cocycle
          }
          }
        }{
          c
        }
      }_>>>{\ }="s"
      \ar[dr]_-{
        \mathllap{
          \mbox{
            \tiny
            \color{darkblue}
            \bf
            twist
          }
          \;
        }
        \tau
      }^-{\ }="t"
      &&
      A \!\sslash\! G
      \ar[dl]^-{
        \underset{
          \mathrlap{
          \raisebox{-2pt}{
            \tiny
            \color{darkblue}
            \bf
            \def\arraystretch{.9}
            \begin{tabular}{c}
              local
              \\
              coefficients
            \end{tabular}
          }
          }
        }{
          \rho
        }
      }
      \\
      &
      B G
      \ar@{=>}^-{\simeq} "s"; "t"
    }
    }
    \right\}_{
      \!\!
      \big/
      \!\!\!\!\!\!\!
      \mbox{
        \tiny
        \def\arraystretch{.9}
        \begin{tabular}{c}
          homotopy
          \\
          relative $B G$
        \end{tabular}
      }
    }
  \end{equation}

  \vspace{-1mm}
  \noindent
  in the homotopy category (Def. \ref{HomotopyCategory})
  of the slice model category over $B G$
  (Example \ref{SliceModelCategory})
  of the classical model category on topological spaces
  (Example \ref{ClassicalModelStructureOnTopologicalSpaces}).
\end{defn}

\begin{defn}[Associated coefficient bundle {\cite[\S 4.1]{NSS12a}\cite[Prop. 2.92]{SS20b}}]
  \label{AssociatedBundles}
  Given a local coefficient $A$-fiber bundle $\rho$ \eqref{LocalCoefficientBundle}
  and a twist $\tau$ \eqref{ATwist} on a domain space $X$,
  the corresponding
  \emph{associated $A$-fiber bundle} over $X$ is the
  homotopy pullback (Def. \ref{HomotopyPullback})
  of $\rho$ along $\tau$, sitting in a homotopy pullback
  square \eqref{ConstructionOfHomotopyPullback} of this form:

  \vspace{-2mm}
  \begin{equation}
    \label{AssociatedCoefficientBundle}
    \raisebox{20pt}{
    \xymatrix@C=4em@R=2.5em{
    {\phantom{A}}
    \ar@{}[r]|- {
      \raisebox{3pt}{
          \tiny
          \color{darkblue}
          \bf
          {
          \def\arraystretch{.9}
          \begin{tabular}{c}
            associated
            \\
            $A$-fiber
            bundle
          \end{tabular}}
        }
      }
   &
   E
      \ar[rr]
      \ar[d]_-{
        \mathbb{R}\tau^\ast \rho
      }
      \ar@{}[drr]|-{
        \underset{
          \raisebox{2pt}{
            \tiny
            \color{darkblue}
            \bf
            homotopy pullback
          }
        }{
          \mbox{\tiny(hpb)}
        }
      }
      &&
      A \!\sslash\! G
  \ar[d]^-{ \rho }
  &
      {\phantom{A}}
        \ar@{}[l]|-{\raisebox{3pt}{
          \tiny
          \color{darkblue}
          \bf
          {
          \def\arraystretch{.9}
          \begin{tabular}{c}
            local
            \\
            coefficient
            bundle
          \end{tabular}}
        }
        }
        \\
      &
      X
      \ar[rr]^-{\tau}_-{
        \raisebox{-3pt}{
        \tiny
        \color{darkblue}
        \bf
        twist
        }
      }
      &&
      B G
    }
    }
  \end{equation}

  \vspace{-.3cm}
 \noindent
  We write

  \vspace{-.3cm}
  \begin{equation}
    \label{VerticalHomotopyClassesOfSections}
    \overset{
      \mathclap{
      \raisebox{3pt}{
        \tiny
        \color{darkblue}
        \bf
        \def\arraystretch{.9}
        \begin{tabular}{c}
          sections of
          \\
          associated bundle
        \end{tabular}
      }
      }
    }{
      \Gamma_X(E)_{\!/\sim}
    }
    \;\;:=\;\;
    \mathrm{Ho}
    \Big(
      \TopologicalSpaces_{\mathrm{Qu}}^{/X}
    \Big)
    \big(
      \mathrm{id}_X
      \,,\,
      \mathbb{R}\tau^\ast \rho
    \big)
    \;\;=\;\;
    \left\{\!\!
    \raisebox{15pt}{
    \xymatrix@R=1em@C=2em{
      &&
      E
      \mathrlap{
       \!\!\!\!\!
       \mbox{
          \tiny
          \color{darkblue}
          \bf
          \def\arraystretch{.9}
          \begin{tabular}{c}
            associated
            \\
            bundle
          \end{tabular}
        }
      }
      \ar[d]
      \\
      X
      \ar@{=}[rr]
      \ar@{-->}[urr]^-{
        \overset{
          \mathllap{
          \raisebox{3pt}{
            \tiny
            \color{darkblue}
            \bf
            section
          }
          }
        }{
          \sigma
        }
      }
      &&
      X
    }
    }
    \qquad \quad
    \right\}_{
      \!\!
      \big/
      \!\!\!\!\!\mbox{\tiny
      \def\arraystretch{.9}
      \begin{tabular}{c}
      vertical
      \\
      homotopy
    \end{tabular}}}
  \end{equation}

  \vspace{-1mm}
  \noindent
  for the set of vertical homotopy classes of section of
  the associated bundle,
  hence for the hom-set, from the identity on $X$ to
  the associated bundle projection,
  in the homotopy category (Def. \ref{HomotopyCategory})
  of the slice model category over $X$
  (Example \ref{SliceModelCategory})
  of the classical model category on topological spaces
  (Example \ref{ClassicalModelStructureOnTopologicalSpaces}).
\end{defn}

\begin{prop}[Twisted non-abelian cohomology is sections of associated
coefficient bundle {\cite[Prop. 4.17]{NSS12a}}]
  \label{TwistedNonAbelianCohomologyIsSectionsOfAssociatedBundle}
  Given a local coefficient bundle $\rho$ \eqref{LocalCoefficientBundle}
  and a twist $\tau$ \eqref{ATwist},
  the $\tau$-twisted non-abelian cohomology (Def. \ref{NonabelianTwistedCohomology})
  with local coefficient in $\rho$
  is equivalent to the vertical homotopy classes of sections
  \eqref{VerticalHomotopyClassesOfSections}
  of the associated coefficient bundle $E$ (Def. \ref{AssociatedBundles}):
    \vspace{-2mm}
  \begin{equation}
    \label{TwistedNonAbelianCohomologyAsSections}
    \overset{
      \mathclap{
      \raisebox{3pt}{
        \tiny
        \color{darkblue}
        \bf
        {
        \def\arraystretch{.9}
        \begin{tabular}{c}
          twisted non-abelian
          \\
          cohomology
        \end{tabular}}
      }
      }
    }{
    H^\tau
    (
      X;
      \,
      A
    )
    }
    \;\;\simeq\;\;
    \overset{
      \mathclap{
      \raisebox{3pt}{
        \tiny
        \color{darkblue}
        \bf
        {
        \def\arraystretch{.9}
        \begin{tabular}{c}
          sections of
          \\
          associated bundle
        \end{tabular}}
      }
      }
    }{
      \Gamma_X(E)_{\!\!/\sim}
    }
    \,.
  \end{equation}
\end{prop}
\begin{proof}
  Consider the following sequence of bijections:
  \vspace{-1mm}
  $$
    \begin{aligned}
      H^\tau
      (
        X;
        \,
        A
      )
      & =
      \mathrm{Ho}
      \Big(
        \TopologicalSpaces^{/ B G}_{\mathrm{Qu}}
      \Big)
      \big(
        \tau
        \,,\,
        \rho
      \big)
      \\
      & \simeq
      \mathrm{Ho}
      \Big(
        \TopologicalSpaces^{/ B G}_{\mathrm{Qu}}
      \Big)
      \big(
        \LeftDerived\tau_\ast \mathrm{id}_X
        \,,\,
        \rho
      \big)
      \\
      & \simeq
      \mathrm{Ho}
      \Big(
        \TopologicalSpaces^{/ X}_{\mathrm{Qu}}
      \Big)
      \big(
        \mathrm{id}_X
        \,,\,
        \RightDerived\tau^\ast \rho
      \big)
      \\
      & =
      \Gamma_X(E)_{\!/\sim}
      \,.
    \end{aligned}
  $$

  \vspace{-1mm}
  \noindent
  Here the first line is the definition \eqref{TwistedCohomologyHom}.
  Then the first step is the observation that every
  slice object is the derived
  left base change
  (Ex. \ref{BaseChangeQuillenAdjunction},
   Prop. \ref{DerivedFunctors})
  along itself of the identity on its domain, by \eqref{LeftBaseChange}.
  With this, the second step is the hom-isomorphism
  \eqref{AnAdjunction}
  of the derived base change adjunction
  $\LeftDerived\tau_! \dashv \mathbb{R}\tau^\ast$.
  The last line is \eqref{VerticalHomotopyClassesOfSections}.
\end{proof}

In twisted generalization of Example \ref{OrdinaryCohomology} we have:
\begin{example}[Twisted ordinary cohomology with local coefficients]
  \label{TwistedOrdinaryCohomology}
  Let $n \in \mathbb{N}$,
  let $X \in \HomotopyTypes$ (Ex. \ref{TheClassicalHomotopyCategory})
  be connected
  and consider a traditional \emph{system of local coefficients}
  \cite[\S 3]{Steenrod43} (see also \cite{MQRT77}\cite{ABG10}\cite{GS-Deligne})
      \vspace{-2mm}
  $$
    \xymatrix{
      \Pi_1(X)
      \ar[r]^-{ t }
      &
      \AbelianGroups
      \,,
    }
  $$

  \vspace{-1mm}
  \noindent
  namely, a functor from the fundamental groupoid of $X$ to
  the category of abelian groups. Since the construction
  $A \mapsto K(A,n)$
  of Eilenberg-MacLane spaces \eqref{EilenbergMacLaneSpaces}
  is itself functorial and using the assumption that
  $X$ is connected, this induces
  (see \cite[Def. 3.1]{BFGM03})
  a local coefficient bundle
  \eqref{LocalCoefficientBundle} of the form
     \vspace{-1mm}
  \begin{equation}
    \label{EilenbergMacLaneCoefficientBundle}
    \xymatrix@R=1em{
      K(A,n)
      \ar[r]
      &
      K(A,n) \sslash \pi_1(X)\;.
      \ar[d]^{ \rho_{t} }
      \\
      & B \pi_1(X)
    }
    \end{equation}

     \vspace{-3mm}
  \noindent
  Finally, write
  $\!\!
    \xymatrix@C=12pt{
      X \ar[r]^-{ \tau } & B \pi_1(X)
    }
  \!\!$
  for the classifying map (via Example \ref{TraditionalNonAbelianCohomology})
  of the universal connected cover of
  $X$ (equivalently: for the 1-truncation projection of $X$).
  Then the $\tau$-twisted non-abelian cohomology
  (Def. \ref{NonabelianTwistedCohomology})
  of $X$ with local coefficients in $\rho_t$ \eqref{EilenbergMacLaneCoefficientBundle}
  is equivalently {\it $t$-twisted ordinary cohomology},
  traditionally known as ordinary cohomology
  {\it with local coefficients} $t$:
  \vspace{-2mm}
  $$
    \xymatrix{
      \overset{
        \mathclap{
        \raisebox{3pt}{
          \tiny
          \color{darkblue}
          \bf
          \def\arraystretch{.9}
          \begin{tabular}{c}
            twisted
            \\
            ordinary cohomology
          \end{tabular}
        }
        }
        }
      {
      H^{n+t}
      (
        X;
        \,
        A
      )
      }
      \;\;
      \simeq
      \;\;
      H^\tau
      \big(
        X;
        \,
        K(A,n)
      \big).
    }
  $$

  \vspace{-2mm}
  \noindent
  This is manifest from comparing
  Def. \ref{NonabelianTwistedCohomology}
  with the characterization of
  cohomology with local coefficients
  found in
  \cite[Cor. 1.3]{Hirashima79}\cite[p. 332]{GoerssJardine99}\cite[Lemma 4.2]{BFGM03}.
\end{example}

\begin{example}[Classification of tangential structure]
  \label{ClassificationOfTangentialStructure}
  Let $X$ be a smooth manifold of dimension
  $n$. Its \emph{frame bundle} is an
  $\mathrm{O}(n)$-principal bundle $\mathrm{Fr}(X) \to X$,
  whose class (a diffeomorphism invariant of $X$)

  \vspace{-.4cm}
  \begin{equation}
    \label{ClassOfFrameBundle}
    \xymatrix@R=0pt{
      \mathrm{O}(n)\mathrm{Bundles}(X)_{\!/\sim}
      \ar@{->}[rr]^-{ \simeq }
      &&
      H\big(X;\, B \mathrm{O}(n) \big)
      \\
      \big[
        \mathrm{Fr}(X)
      \big]
      \ar@{}[rr]|-{\leftrightarrow }
      &&
      \big[
        \tau_{\mathrm{fr}}
      \big]
    }
  \end{equation}
  \vspace{-.4cm}

  \noindent
  gives, by Example \ref{TraditionalNonAbelianCohomology},
  the class of a twist $\tau_{\mathrm{fr}}$ \eqref{ATwist} in the non-abelian
  $\mathrm{O}(n)$-cohomology of $X$.

  Now for $B G$ any connected homotopy type (Prop. \ref{ConnectedHomotopyTypesAreHigherNonAbelianClassifyingSpaces})
  and for $\xymatrix@C=12pt{ B G \ar[r]^\rho & B \mathrm{O}(n) }$
  any map (equivalently the delooping of a morphism of $\infty$-groups
  $\xymatrix@C=12pt{G \ar[r] & \mathrm{O}(n)}$), we get a local
  coefficient bundle \eqref{LocalCoefficientBundle}
  with (homotopy-)coset space fiber \cite[Lemma 2.7]{FSS19b}:

  \vspace{-.5cm}
  \begin{equation}
    \label{HomotopyCosetSpaces}
    \xymatrix{
      \mathrm{O}(n) \!\sslash\! G
      \ar[rr]^-{ \mathrm{hofib}(\rho) }
      &&
      B G
      \ar[d]^-{ \rho }
      \\
      &&
      B \mathrm{O}(n)
      \,.
    }
  \end{equation}
  The relative homotopy class of
  a homotopy lift of the frame bundle classifier $\tau_{\mathrm{fr}}$
  \eqref{ClassificationOfTangentialStructure}
  through this map $\rho$

  \vspace{-.4cm}
  \begin{equation}
    \label{TangentialStructures}
    \left[
    \raisebox{20pt}{
    \xymatrix{
      X
      \ar@{-->}[rr]^-{
        \mbox{
          \tiny
          \color{darkblue}
          \bf
          tangential structure
        }
      }_>>>>>{\ }="s"
      \ar[dr]_-{ \tau_{\mathrm{fr}} }^-{\ }="t"
      &&
      B G
      \ar[dl]^-{ \rho }
      \\
      & B \mathrm{O}(n)
      \ar@{=>}^g "s"; "t"
    }
    }
    \right]
    \;\;
      \in
    \;\;
    G\mathrm{TangentialStructures}(X)
  \end{equation}
  is known
  a \emph{topological $G$-structure} or
  \emph{tangential $\rho$-structure} on $X$
  (e.g. \cite[\S 1.4]{Kochman96}\cite[\S 5]{GMTW09}\cite[Def. 4.48]{SS20b}).
  For instance, for $\rho$ a stage in the Whitehead tower of
  $\mathrm{O}(n)$, this is, in turn,
  \emph{Orientation}, \emph{Spin structure}, \emph{String structure},
  \emph{Fivebrane structure} \cite{SSS09}, etc.:

  
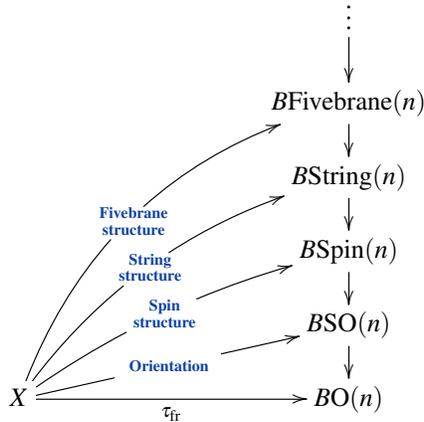
\vspace{-.7cm}
  $$
    \xymatrix@C=40pt@R=12pt{
      &&
      \vdots
      \ar[d]
      \\
      &&
      B \mathrm{Fivebrane}(n)
      \ar[d]
      \\
      &&
      B \mathrm{String}(n)
      \ar[d]
      \\
      &&
      B \mathrm{Spin}(n)
      \ar[d]
      \\
      &&
      B \mathrm{SO}(n)
      \ar[d]
      \\
      X
      \ar[urr]|-{
        \mbox{
          \tiny
          \color{darkblue}
          \bf
          {\begin{tabular}{c}
            Orientation
          \end{tabular}}
        }
      }
      \ar@/^.6pc/[uurr]|-{
        \mbox{
          \tiny
          \color{darkblue}
          \bf
          {
          \def\arraystretch{.9}
          \begin{tabular}{c}
            Spin
            \\
            structure
          \end{tabular}}
        }
      }
      \ar@/^1.2pc/[uuurr]|-{
        \mbox{
          \tiny
          \color{darkblue}
          \bf
          {
          \def\arraystretch{.9}
          \begin{tabular}{c}
            String
            \\
            structure
          \end{tabular}}
        }
      }
      \ar@/^1.9pc/[uuuurr]|-{
        \mbox{
          \tiny
          \color{darkblue}
          \bf
          \def\arraystretch{.9}
          {
          \def\arraystretch{.9}
          \begin{tabular}{c}
            Fivebrane
            \\
            structure
          \end{tabular}}
        }
      }
      \ar[rr]_-{ \tau_{\mathrm{fr}} }
      &&
      B \mathrm{O}(n)
    }
  $$

  By comparison of \eqref{TangentialStructures} with \eqref{TwistedCohomologyHom}
  we see that tangential $G$-structures on $X$
  are classified
  by twisted non-abelian cohomology (Def. \ref{NonabelianTwistedCohomology})
  with coefficients in
  (homotopy-)coset spaces $\mathrm{O}(n)\!\sslash\! G$
  \eqref{HomotopyCosetSpaces}
  and twisted by the class $\tau_{\mathrm{fr}}$ of the frame bundle
  \eqref{ClassOfFrameBundle}:
  \begin{equation}
    \label{TangentialStrucuresAndTwistedCohomology}
    G\mathrm{TangentialStructures}(X)
    \;\;
    \simeq
    \;\;
    H^{ \tau_{\mathrm{fr}} }
    \big(
      X;
      \,
      \mathrm{O}(n) \!\sslash\! G
    \big)
    \,.
  \end{equation}
  \vspace{-.4cm}

  \noindent
  According to Prop. \ref{ConnectedHomotopyTypesAreHigherNonAbelianClassifyingSpaces},
  this example is actually universal
  for $\tau_{\mathrm{fr}}$-twisted non-abelian cohomology.
\end{example}

As a special case of
Example \ref{TwistedOrdinaryCohomology} and
in twisted generalization of Examples \ref{BundleGerbes},
\ref{HigherBundleGerbes} we have:
\begin{example}[Orientifold gerbes]
  Consider the action $\sigma_{\mathrm{U}(1)}$ of $\mathbb{Z}_2$ on the circle
  group $\mathrm{U}(1) \subset \mathbb{C}^\times$ given by complex
  conjugation. This deloops (see \cite[\S 4.4]{FSS15}) to
  an action $\sigma_{B^n \mathrm{U}(1)}$ of $\mathbb{Z}_2$ on
  the classifying spaces $B^n \mathrm{U}(1)$ \eqref{ClassifyingSpace}.
  By Prop. \ref{GActionsAsFibrations}
  there is a corresponding local coefficient bundle
  \vspace{-2mm}
  \begin{equation}
    \label{LocalCoefficientBundleForZ2ActionOn}
    \xymatrix@R=1.5em{
      B^n \mathrm{U}(1)
      \ar[r]
      &
      B^n \mathrm{U}(1) \!\sslash\! \mathbb{Z}_2
      \ar[d]^-{ \sigma_{B^n \mathrm{U}(1)} }
      \\
      &
      B \mathbb{Z}_2
    }
  \end{equation}

    \vspace{-4mm}
  \noindent
  Moreover, consider a smooth manifold $X$, with orientation
  bundle classified by
  $
    \xymatrix@C=12pt{ X \ar[r]^-{\mathrm{or}} & B \mathbb{Z}_2 }
    \,.
  $
  Then the $\mathrm{or}$-twisted cohomology (Def. \ref{NonabelianTwistedCohomology}) of $X$...

  \noindent
  {\bf (i)}
  ...with local coefficients
  in $\sigma_{B^2 \mathrm{U}(1)}$ classifies
  what is equivalently known as Jandl gerbes
  \cite{SSW07}\cite{GSW11} or real gerbes \cite{HMSV19}
  or
  orientifold B-fields;

  \noindent
  {\bf (ii)}
  ...with local coefficients
  in $\sigma_{B^3 \mathrm{U}(1)}$ classifies
  what is equivalently known as
  topological sectors of orientifold C-fields
  \cite[\S 4.4]{FSS15}.

  \noindent More generally, one can consider
  twisted Deligne cohomology \cite{GS-Deligne}
  as well as higher-twisted
  periodic integral- and Deligne-cohomology \cite{GS-HigherDeligne}
  (see also \cref{NonabelianDifferentialCohomology}).
\end{example}

\begin{remark}[The Whitehead principle of non-abelian cohomology]
  \label{WhiteheadPrincipleOfNonAbelianCohomology}
  $\,$

  \noindent
  Let $A \in \HomotopyTypes$ be connected,
  so that $A \simeq B G$ (Prop. \ref{ConnectedHomotopyTypesAreHigherNonAbelianClassifyingSpaces}).

  \noindent {\bf (i)} If $A$ is also $n$-truncated \eqref{nTruncated}, then its
  Postnikov tower (Prop. \ref{PostnikovTower}) says that
  $A$ is the total space of a local coefficient bundle
  \eqref{NonabelianTwistedCohomology}
  of the form

  \vspace{-.7cm}
  $$
    \xymatrix@R=1.5em{
      K( \pi_n(A), n )
      \ar[r]^-{ \mathrm{hfib}(p_n) }
      &
      A
      \ar[d]^-{p^A_n}
      \\
      &
      A(n-1)
      \mathrlap{
        \;
        \simeq
        B \big( G(n-2) \big)
      }
    }
  $$
  \vspace{-.7cm}

  \noindent
  with homotopy fiber an Eilenberg-MacLane space \eqref{EilenbergMacLaneSpaces}.

 \vspace{1mm}
 \noindent {\bf (ii)}
 Accordingly, non-abelian cohomology with coefficients in $A$
  (Def. \ref{NonAbelianCohomology})
  is equivalently the disjoint union, over the space of twists
  $\tau_{{}_n}$ \eqref{ATwist}
  in non-abelian cohohomology with coefficients in
  $A(n-1)$, of $\tau$-twisted non-abelian cohomology
  (Def. \ref{NonabelianTwistedCohomology})
  with coefficients in $K( \pi_{n}(A), n )$:
  \vspace{-2mm}
  \begin{equation}
    \label{UnravellingNonAbelianCohomologyAsTwistedAbelianCohomology}
    \overset{
      \mathclap{
      \raisebox{3pt}{
        \tiny
        \color{darkblue}
        \bf
        \def\arraystretch{.9}
        \begin{tabular}{c}
          non-abelian cohomology
          \\
          in higher degree
        \end{tabular}
      }
      }
    }{
    H
    (
      X;
      \,
      A
    )
    }
    \;\;
    \;
    \simeq
    \;
    \;\;
    \;\;\;\;\;\;\;\;
    \underset{
      \mathclap{
      \underset{
        \mathclap{
        \raisebox{-3pt}{
          \tiny
          \color{darkblue}
          \bf
          \def\arraystretch{.9}
          \begin{tabular}{c}
            twist in
            \\
            non-abelian cohomology
            \\
            of lower degree
          \end{tabular}
        }
        }
      }{
      \tau_{{}_n} \in
      H
      (
        X;
        \,
        A(n-1)
      )
      }
      }
    }{
      \overset{
        \raisebox{3pt}{
          \tiny
          \color{darkblue}
          \bf
        }
      }{
        \bigsqcup
      }
    }
      \;\;\;\;\;\;\;\;\;
      \overset{
        \mathclap{
        \raisebox{3pt}{
          \tiny
          \color{darkblue}
          \bf
          \def\arraystretch{.9}
          \begin{tabular}{c}
            higher twisted
            \\
            ordinary cohomology
          \end{tabular}
        }
        }
      }{
      H^{\tau_{{}_n}}
      \big(
        X;
        \,
        K(\pi_n(A), n)
      \big)
      }
      \,.
  \end{equation}

  \vspace{-2mm}
  \noindent {\bf (iii)}
  But notice that this is just the first step,
  and that iterating this unravelling yields unwieldy formulas:
  \begin{equation}
    \label{UnravellingNonAbelianCohomologyAsTwistedAbelianCohomologyAgain}
    H
    (
      X;
      \,
      A
    )
    \;
    \;
    \;\;
    \simeq
    \;\;\;
    \;
    \;\;
    \;\;
    \underset{
      \mathclap{
      \tau_{{}_{n}} \in
      \;\;\;\;\;\;\;
      \underset{
        \mathclap{
        \tau_{{}_{n-1}}
        \in \;
        H(
          X;
          \,
          A(n-2)
        )
        }
      }{\bigsqcup}
      \;\;\;\;\;\;\;
      H^{\tau_{{}_{n-1}}}
      \big(
        X;
        \,
        K(\pi_{n-1}(A), n-1)
      \big)
      }
    }{\bigsqcup}
      \;\;\;\;
      H^{\tau_{{}_n}}
      \big(
        X;
        \,
        K(\pi_n(A), n)
      \big)
      \,,
    \;\;
    \mbox{then}
    \;\;
    \;
    H
    (
      X;
      \,
      A
    )
    \;\;
    \;
    \;\;
    \simeq
    \;\;
    \;
    \;\;
    \;\;
    \underset{
      \mathclap{
      \tau_{{}_{n}} \in
      \;\;\;\;\;\;\;
      \underset{
        \mathclap{
        \tau_{{}_{n-1}}
        \in
        \;\;\;\;
        \underset{
          \mathclap{
          \tau_{{}_{n-2}}
          \in
          H
          \big(
            X;
            \,
            A(n-3)
          \big)
          }
        }{\bigsqcup}
        \;\;\;\;
        H^{\tau_{{}_{n-2}}}
        \big(
          X;
          \,
          K( \pi_{n-2}(A), n-2 )
        \big)
        }
      }{\bigsqcup}
      \;\;\;\;\;\;\;
      H^{\tau_{{}_{n-1}}}
      \big(
        X;
        \,
        K(\pi_{n-1}(A), n-1)
      \big)
      }
    }{\bigsqcup}
      \;\;\;\;
      H^{\tau_{{}_n}}
      \big(
        X;
        \,
        K(\pi_n(A), n)
      \big)
      \;\;
      \mbox{etc.}
      \,.
  \end{equation}

%  and then
%  \begin{equation}
%    \label{UnravellingNonAbelianCohomologyAsTwistedAbelianCohomologyYetAgainAndAgain}
%    H
%    \big(
%      X;
%      \,
%      A
%    \big)
%    \;\;
%    \;
%    \;\;\;
%    \simeq
%    \;\;\;
%    \;
%    \;\;
%    \;\;
%    \underset{
%      \mathclap{
%      \tau_{{}_{n}} \in
%      \;\;\;\;\;\;\;
%      \underset{
%        \mathclap{
%        \tau_{{}_{n-1}}
%        \in
%        \;\;\;\;
%        \underset{
%          \mathclap{
%          \tau_{{}_{n-2}}
%          \in
%          \;\;\;\;
%            \underset{
%              \mathclap{
%              \tau_{{}_{n-3}}
%              \;
%              \in
%              \;
%              H
%              \big(
%                X;
%                \,
%                A(n-4)
%              \big)
%              }
%            }{\bigsqcup}
%          \;\;\;\;
%          H^{\tau_{{}_{n-3}}}
%          \big(
%            X;
%            \,
%            K(\pi_{n-3}(A), n-3)
%          \big)
%          }
%        }{\bigsqcup}
%        \;\;\;\;
%        H^{\tau_{{}_{n-2}}}
%        \big(
%          X;
%          \,
%          K( \pi_{n-2}(A), n-2 )
%        \big)
%        }
%      }{\bigsqcup}
%      \;\;\;\;\;\;\;
%      H^{\tau_{{}_{n-1}}}
%      \big(
%        X;
%        \,
%        K(\pi_{n-1}(A), n-1)
%      \big)
%      }
%    }{\bigsqcup}
%      \;\;\;\;
%      H^{\tau_{{}_n}}
%      \big(
%        X;
%        \,
%        K(\pi_n(A), n)
%      \big)
%      \,.
%  \end{equation}
%  and so on.

  \vspace{1mm}
  \noindent {\bf (iv)} Thus, non-abelian cohomology in higher degrees
  (Example \ref{NonAbelianCohomologyInUnboundedDegree})
  decomposes as a tower of
  consecutively higher twisted but otherwise ordinary cohomology theories,
  starting with a twist in non-abelian cohomology in degree 1.
  This phenomenon has been called the {\it Whitehead principle
  of non-abelian cohomology} \cite[p. 8]{Toen02}
  and has been interpreted as
  saying that
  ``nonabelian cohomology occurs essentially only in degree 1''
  \cite[p. 1]{Simpson96}.

 \vspace{1mm}
 \noindent {\bf (v)}  But the above formulas
  \eqref{UnravellingNonAbelianCohomologyAsTwistedAbelianCohomology},
  \eqref{UnravellingNonAbelianCohomologyAsTwistedAbelianCohomologyAgain},
  make manifest that there are two perspectives on this phenomenon.
  On the one hand: non-abelian cohomology
  in higher degrees may be \emph{computed}
  by brute force as a sequence of consecutively higher twisted abelian
  cohomologies, with lowest twist starting in degree-1 non-abelian cohomology.
  On the other hand, conversely: intricate such systems of consecutively twisted abelian cohomology theories are neatly \emph{understood}
  as unified by non-abelian cohomology.

  \vspace{1mm}
  \noindent {\bf (vi)}
  Similarly, even though Postnikov towers
  do exist
  (Prop. \ref{PostnikovTower})
  in the classical homotopy category (Example \ref{TheClassicalHomotopyCategory}), the latter is far from
  being
  equivalent to the stable homotopy category \eqref{StabilizationAdjunction}
  ``up to twists in degree 1''.
\end{remark}

In twisted generalization of Example \ref{ComplexTopologicalKtheory},
we have:
\begin{example}[Twisted topological K-theory]
  \label{TwistedKTheory}
  The classifying space $\mathrm{KU}_0 \,\simeq\,
  \mathbb{Z} \times B \mathrm{U}$ \eqref{ClassifyingSpaceForComplexTopologicalKTheory}
  for complex topological K-theory (Example \ref{ComplexTopologicalKtheory})
  is the fiber of a local coefficient bundle \eqref{LocalCoefficientBundle}
  over $K(\mathbb{Z}, 3) \simeq B^3 \mathrm{U}(1)$ \eqref{KZnAndBnU1}:
    \vspace{-2mm}
  \begin{equation}
    \label{LocalCoefficientBundleForTwistedKTheory}
    \raisebox{20pt}{
    \xymatrix@R=1em{
      \mathrm{KU}_0
      \ar[r]
      &
      \mathrm{KU}_0 \!\sslash\! B \mathrm{U}(1)
      \ar[d]
      \\
      &
      B^2 \mathrm{U}(1)
    }
    }
  \end{equation}

    \vspace{-4mm}
  \noindent
  For $\!\xymatrix@C=12pt{X \ar[r]^-{\tau} & B^2 \mathrm{U}(1)}\!$
  a corresponding twist \eqref{ATwist}
  (hence equivalently a bundle gerbe, by Example \ref{BundleGerbes}),
  the corresponding twisted non-abelian cohomology (Def. \ref{NonabelianTwistedCohomology}) is twisted complex topological K-theory \cite{Karoubi68}\cite{DonovanKaroubi70}:
    \vspace{-2mm}
  \begin{equation}
    \label{ReproducingTwistedTopologicalKTheory}
    \overset{
      \mathclap{
      \raisebox{3pt}{
        \tiny
        \color{darkblue}
        \bf
        \def\arraystretch{.9}
        \begin{tabular}{c}
          twisted
          \\
          topological K-theory
        \end{tabular}
      }
      }
    }{
      \mathrm{KU}^\tau(-)
    }
    \;\;
    \simeq
    \;\;
    H^\tau
    \big(
      -;
      \,
      \mathbb{Z}
      \times
      B \mathrm{U}
    \big)
    \,.
  \end{equation}

  \vspace{-1mm}
  \noindent
  This is manifest from comparing \eqref{TwistedCohomologyHom}
  with \cite[(2.6)]{FHT-complex}.
  Alternatively, under Prop. \ref{TwistedNonAbelianCohomologyIsSectionsOfAssociatedBundle},
  this is manifest from comparing the
  equivalent right hand side of
  \eqref{TwistedNonAbelianCohomologyAsSections} with
  \cite[Prop. 2.1]{Rosenberg89} (using \cite[Cor. 4.18]{NSS12a})
  or, more directly, with \cite[\S 3]{AtiyahSegal04}\cite[\S 2.1]{ABG10}.
\end{example}

Generally, in twisted generalization of Example
\ref{GeneralizedCohomologyAsNonabelianCohomology}, we have:

\begin{example}[Local coefficient bundle for twisted Whitehead-generalized cohomology]
  \label{TwistedGeneralizedCohomology}
  Let $R$ be an $E_\infty$-ring spectrum (Ex. \ref{GeneralizedCohomologyAsNonabelianCohomology})
  and write $\GroupOfUnits{R}$ for its
  {\it $\infty$-group of units} \cite[\S 2.3]{Schlichtkrull04}\cite[\S 22.2]{MaySigurdsson04}\cite[\S 3]{ABGHR08}\cite[\S 2]{ABGHR14a}, defined as
  the homotopy pullback (Def. \ref{HomotopyPullback})
  of the component space
  $R_0 \,=\, \Derived\Omega^\infty R$ \eqref{ComponentSpacesOfSpectra}
  fibered over its 0-truncation (i.e. its 1-coskeleton \eqref{nTruncationMorphismsViaCoskeletonUnits})
  to the ordinary group of units of this ordinary ring of connected components:

  \vspace{-.4cm}
  \begin{equation}
    \label{InfinityGroupOfUnits}
    \begin{tikzcd}[row sep=10pt]
      \overset{
        \mathclap{
        \raisebox{3pt}{
          \tiny
          \color{darkblue}
          \bf
          \def\arraystretch{.9}
          \begin{tabular}{c}
            $\infty$-group of units
          \end{tabular}
        }
        }
      }{
        \GroupOfUnits{R}
      }
      \ar[rr]
      \ar[d]
      \ar[
        drr,
        phantom,
        "\mbox{\tiny\rm (hpb)}"
      ]
      &&
      \overset{
        \mathclap{
        \raisebox{3pt}{
          \tiny
          \color{darkblue}
          \bf
          \def\arraystretch{.9}
          \begin{tabular}{c}
            $E_\infty$-ring space
          \end{tabular}
        }
        }
      }{
        R_0
      }
      \ar[
        d,
        "p_0"
      ]
      \\
      \underset{
        \raisebox{-3pt}{
          \tiny
          \color{darkblue}
          \bf
          \def\arraystretch{.9}
          \begin{tabular}{c}
            ordinary group
            \\
            of units
          \end{tabular}
        }
      }{
        \GroupOfUnits{\pi_0(R_0)}
      }
      \ar[
        rr,
        hook
      ]
      &&
      \underset{
        \mathclap{
        \raisebox{-3pt}{
          \tiny
          \color{darkblue}
          \bf
          \def\arraystretch{.9}
          \begin{tabular}{c}
            ordinary ring of
            \\
            connected components
          \end{tabular}
        }
        }
      }{
        \pi_0\big( R_0 \big)
      }
  \end{tikzcd}
\end{equation}
\vspace{-.3cm}

  \noindent
  This makes $\GroupOfUnits{R}$ as an
  $\infty$-group (as in Example \ref{NonAbelianCohomologyInUnboundedDegree})
  with group operation induced from the {\it multiplicative}
  structure on $R_0$.
  The canonical action of
  $\GroupOfUnits{R}$ on
  $R_0$ is given, via Prop. \ref{GActionsAsFibrations},
  by a local coefficient bundle \eqref{LocalCoefficientBundle}
  of this form:

  \vspace{-.3cm}
  \begin{equation}
    \label{LocalCoefficientBundleForSpectra}
    \begin{tikzcd}
      R_0
      \ar[r]
      &
      (
        R_0
      )
        \!\sslash\!
      \GroupOfUnits{R}
      \ar[
        d,
        "{ \rho_R }"
      ]
      \\
      &
      B \GroupOfUnits{R}
      \mathrlap{\,.}
    \end{tikzcd}
  \end{equation}
\end{example}

\begin{prop}[Twisted non-abelian cohomology subsumes twisted generalized cohomology]
  \label{ProofTwistedGeneralizedCohomology}
  For $R$ an $E_\infty$-ring spectrum (Ex. \ref{GeneralizedCohomologyAsNonabelianCohomology}),
  the twisted non-abelian cohomology
  (Def. \ref{NonabelianTwistedCohomology}) with local
  coefficient bundle $\rho_R$ from Example \ref{TwistedGeneralizedCohomology}
  is, equivalently,
  twisted generalized $R$-cohomology in the traditional sense
  (e.g. \cite[\S 22.1]{MaySigurdsson04}):

  \vspace{-2mm}
  \begin{equation}
    \label{ReproducingTwistedGeneralizedCohomology}
    \overset{
      \mathclap{
      \raisebox{3pt}{
        \tiny
        \color{darkblue}
        \bf
        \def\arraystretch{.9}
        \begin{tabular}{c}
          twisted Whitehead-
          \\
          generalized cohomology
        \end{tabular}
      }
      }
    }{
      R^\tau(-)
    }
    \;\;\;\;\;\;\;
    \simeq
    \;\;\;\;\;\;\;
    \overset{
      \mathclap{
      \raisebox{3pt}{
        \tiny
        \color{darkblue}
        \bf
        \def\arraystretch{.9}
        \begin{tabular}{c}
          twisted non-abelian cohomology
          \\
          with local $\rho_R$-coefficients
        \end{tabular}
      }
      }
    }{
    H^\tau
    (
      -;
      \,
      \rho_R
    )
    }
    \;.
  \end{equation}
\end{prop}
\begin{proof}
  Given any twist
  $\!\xymatrix@C=12pt{X \ar[r]^-{ \tau } & B \GroupOfUnits{R}}\!$
  \eqref{NonabelianTwistedCohomology}, write $P \to X$
  for the homotopy pullback (Def. \ref{HomotopyPullback})
  along $\tau$
  of the essentially unique point inclusion:
  \vspace{-2mm}
  \begin{equation}
    \label{CoefficientBundleAssociatedToPrincipalBundle}
    \raisebox{20pt}{
    \xymatrix@R=20pt@C=4em{
      P
      \ar[r]
      \ar[d]
      \ar@{}[dr]|-{\mbox{\tiny\rm(hpb)}}
      &
      \ast
      \ar[d]
      \\
      X
      \ar[r]_-{ \tau }
      &
      B \GroupOfUnits{R}
    }
    }
    \phantom{AAAAAAAAAAAAAAAAAAA}
    \raisebox{20pt}{
    \xymatrix@R=20pt@C=4em{
      \mathllap{
        \scalebox{.9}{$
        \big(
          P \times R_0
        \big)\!\sslash_{\scalebox{.5}{diag}}\! \GroupOfUnits{R}
        $}
        \simeq
        \;
      }
      E
      \ar[r]
      \ar[d]_-{
        \mathbb{R}\tau^\ast \rho_R
      }
      \ar@{}[dr]|-{\mbox{\tiny\rm(hpb)}}
      &
      R_0 \!\sslash\! \GroupOfUnits{R}
      \ar[d]^-{ \rho_R }
      \\
      X
      \ar[r]_-{ \tau }
      &
      B \GroupOfUnits{R}
    }
    }
  \end{equation}
  \vspace{-2mm}

  \noindent
  This $P$ is the $\GroupOfUnits{R}$-principal $\infty$-bundle
  which is classified by $\tau$, \cite[Thm. 3.17]{NSS12a},
  to which the coefficient bundle $E$ \eqref{AssociatedCoefficientBundle} is
  $\GroupOfUnits{R}$-associated \cite[Prop. 4.6]{NSS12a},
  as shown on the right of \eqref{CoefficientBundleAssociatedToPrincipalBundle}.
  Consider then the following sequence of natural bijections:
  \vspace{-1mm}
  \begin{equation}
    \label{TowardsTwistedGeneralizedCohomology}
      \begin{aligned}
        H^\tau
        \big(
          X;
          \,
          R_0
        \big)
        & \simeq
        \Gamma_X(E)
        \\
        & \simeq
        \mathrm{Ho}
        \big(
          \GroupOfUnits{R}\mathrm{Actions}
        \big)
        (
          P;
          \,
          R_0
        )
        \\
        & \simeq
        \mathrm{Ho}
        \big(
          R\mathrm{Modules}
        \big)
        (
          M \tau;
          \,
          R
        )
        \\
        & \simeq
        R^\tau(X)\;.
      \end{aligned}
  \end{equation}

    \vspace{-2mm}
\noindent
  Here  the first step is Prop.
  \ref{TwistedNonAbelianCohomologyIsSectionsOfAssociatedBundle},
  while the second step is \cite[Cor. 4.18]{NSS12a}.
  The  third step is \cite[(2.15)]{ABGHR08}\cite[(3.15)]{ABGHR14a},
  with $M \tau$ denoting the $R$-Thom spectrum of $\tau$
  \cite[Def. 2.6]{ABGHR08}\cite[Def. 3.13]{ABGHR14a}.
  The  last step
  is \cite[\S 2.5]{ABGHR08} \cite[\S 1.4]{ABGHR14a}\cite[\S 2.7]{ABGHR14b}.
  The composite of these natural bijections is the desired
  \eqref{ReproducingTwistedGeneralizedCohomology}.
\end{proof}

\begin{example}[Higher Cohomotopy-twisted K-theory]
  \label{CohomotopyTwistedKTheory}
  For complex topological K-theory $R = \mathrm{KU}$
  (Ex. \ref{ComplexTopologicalKtheory})
  with $\mathrm{KU}_0 \,=\, \mathbb{Z} \times B\mathrm{U}$
  \eqref{ClassifyingSpaceForComplexTopologicalKTheory}
  --
  where the $\mathbb{Z}$-factor encodes the virtual rank of vector bundles
  and the multiplicative operation in the
  ring structure corresponds to tensor product of vector bundles --
  the $\infty$-group of units
  \eqref{InfinityGroupOfUnits} classifies the virtual vector bundles
  of invertible rank in $\{\pm 1\} = \GroupOfUnits{\mathbb{Z}} \subset \mathbb{Z}$:

  \vspace{-.3cm}
  \begin{equation}
    \GroupOfUnits{\mathrm{KU}}
    \;\simeq\;
    \big(
      \{\pm 1\} \times B \mathrm{U}
    \big)_{\otimes}
    \,.
  \end{equation}
  \vspace{-.5cm}

  \noindent
  (Here the subscript just indicates the $\infty$-group structure, now
  with respect to the multiplicative operation corresponding to tensor
  product of virtual vector bundles.)
  Since delooping $B(-)$ shifts up homotopy groups by one, it follows that
  the homotopy groups of $B \GroupOfUnits{\mathrm{KU}}$,
  appearing in \eqref{LocalCoefficientBundleForSpectra},
  are freely generated
  by the powers of the Bott generator $\beta \in \pi_2(\mathrm{KU})$,
  shifted up in degree by one:

  \vspace{-.5cm}
  \begin{equation}
    \label{HomotopyGroupsOfGroupOfUnitsOfKU}
    \pi_\bullet
    \big(
      \mathrm{KU}_0
    \big)
    \;\simeq\;
    \left\{
    \!\!\!
    \begin{array}{ll}
      \mathbb{Z} \simeq \langle\beta^k\rangle
      & \mbox{if $n = 2k$ is even}
      \\
      0 & \mbox{if $n$ is odd}
    \end{array}
    \right.
    {\phantom{AA}}
    \Rightarrow
    {\phantom{AA}}
    \pi_n
    \big(
      \GroupOfUnits{\mathrm{KU}}
    \big)
    \;\simeq\;
    \left\{
    \!\!\!
    \begin{array}{ll}
      \mathbb{Z}_2 & \mbox{if $n = 1$}
      \\
      \mathbb{Z} \simeq \langle \beta^{2k}\rangle
      &
      \mbox{if $n = 2k+3$}
      \\
      0 & \mbox{if $n$ is even.}
    \end{array}
    \right.
  \end{equation}
  \vspace{-.3cm}

  \noindent
  It follows
  that,
  parameterized by any odd-dimensional sphere
  $S^{2k+1}$ for positive $k \in \mathbb{N}_+$,
  there
  are exactly $\mathbb{Z}$ worth of higher twists of
  complex K-theory, up to equivalence, embodied by the local coefficient bundles
  which are the
  homotopy pullback of \eqref{LocalCoefficientBundleForSpectra}
  along the classifying maps of the elements \eqref{HomotopyGroupsOfGroupOfUnitsOfKU}.
  The universal one among these is the pullback along the classifying map
  for the suspended power of the Bott generator itself:

  \vspace{-.3cm}
  \begin{equation}
    \label{LocalCoefficientBundleForCohomotopicallyTwistedKtheory}
    \begin{tikzcd}[column sep=40pt]
      \mathrm{KU}_0
      \ar[r]
      &
      \overset{
        \mathclap{
        \raisebox{3pt}{
          \tiny
          \color{darkblue}
          \bf
          \def\arraystretch{.9}
          \begin{tabular}{c}
            local coefficient bundle for
            \\
            Cohomotopy-twisted K-theory
          \end{tabular}
        }
        }
      }{
        \big( \mathrm{KU}_0\big)
          \!\sslash\!
        \Omega S^{2k+1}
      }
      \ar[rr]
      \ar[
        drr,
        phantom,
        "\mbox{\tiny\rm(hpb)}"
      ]
      \ar[
        d,
        "\rho^{S^{2k+1}}_{\mathrm{KU}}"
      ]
      &&
      \overset{
        \mathclap{
        \raisebox{3pt}{
          \tiny
          \color{darkblue}
          \bf
          \def\arraystretch{.9}
          \begin{tabular}{c}
            universal local coefficient bundle
            \\
            for twisted complex K-theory
          \end{tabular}
        }
        }
      }{
        \big( \mathrm{KU}_0\big)
          \!\sslash\!
        \GroupOfUnits{\mathrm{KU}}
      }
      \ar[
        d,
        "\rho_{\mathrm{KU}}"
      ]
      \\
      &
      S^{2k+1}
      \ar[
        rr,
        "{
          \Sigma(\beta^{2k})
        }"{below},
        "
          \mbox{
            \tiny
            \color{greenii}
            \bf
            shifted power of Bott generator
          }
        "
      ]
      &&
      B \GroupOfUnits{\mathrm{KU}}
    \end{tikzcd}
  \end{equation}
  \vspace{-.3cm}

\noindent  By Def. \ref{NonabelianTwistedCohomology},
  these local coefficient bundles encode higher twists of
  complex K-theory
  (Ex. \ref{ComplexTopologicalKtheory})
  by classes in unstable/non-abelian Cohomotopy
  (Ex. \ref{CohomotopyTheory}) in degree $2k+1$:

  \vspace{-.2cm}
  \begin{equation}
    \label{CohomotopicallyTwistedTopologicalKTheory}
    \overset{
      \mathclap{
     \raisebox{3pt}{
       \tiny
       \color{darkblue}
       \bf
       twist in Cohomotopy
     }
     }
    }{
    [
      \lambda
    ]
    \;\in\;
    \pi^{2k+1}(X)
    }
    {\phantom{AAA}}
    \vdash
    {\phantom{AAA}}
    \overset{
      \mathclap{
      \raisebox{3pt}{
        \tiny
        \color{darkblue}
        \bf
        \def\arraystretch{.9}
        \begin{tabular}{c}
          cohomotopically-twisted
          topological K-theory
        \end{tabular}
      }
      }
    }{
    \mathrm{KU}^{\lambda}(X)
    \;=\;
    H^\lambda
    (
      X;
      \,
      \mathrm{KU}_0
    )
    }
    \,.
  \end{equation}
  \vspace{-.5cm}

  \noindent
  This cohomotopically higher twisted K-theory
  has been considered in \cite[Def. 2.5]{MMS20}.
\end{example}

In twisted generalization of Example \ref{IteratedKTheory},
we have:
\begin{example}[Twisted iterated K-theory]
  \label{TwistedIteratedkTheory}
  Let $r \in \mathbb{N}$, $r \geq 1$.
  By \cite[Prop. 1.5, Def. 1.7]{LindSatiWesterland16}
  and using Prop. \ref{ProofTwistedGeneralizedCohomology},
  there is
  a local
  coefficient bundle \eqref{LocalCoefficientBundle} of the form
  \vspace{-3mm}
  \begin{equation}
    \label{LocalCoefficientBundleForIteratedKTheory}
    \xymatrix@R=1.5em{
      \big(
        K^{2r-2}(\mathrm{ku})
      \big)_0
      \ar[r]
      &
      \Big(
        \big(
          K^{2r-2}(\mathrm{ku})
        \big)_0
      \Big) \!\sslash\! B^{2r-1} \mathrm{U}(1)
      \ar[d]^-{ \rho_{\mathrm{lsw}_{2r-1}} }
      \\
      &
      B^{2r} \mathrm{U}(1)
      \,,
    }
  \end{equation}

  \vspace{-2mm}
  \noindent
  where $K^{2r-2}(\mathrm{ku})_0$ is the 0th space in the
  spectrum \eqref{Spectrum} representing iterated K-theory
  (Ex. \ref{IteratedKTheory}) and
  $B^{2r}\mathrm{U}(1) \,\simeq\, K(\mathbb{Z}, 2r+1)$
  is the classifying space for bundle $(2r-1)$-gerbes
  (Ex. \ref{HigherBundleGerbes}).
  This means that
  for $\!\xymatrix@C=12pt{ X \ar[r]^-\tau & B^{2r}\mathrm{U}(1) }\!$
  a classifying map for such a higher gerbe,
  the $\tau$-twisted non-abelian cohomology (Def. \ref{NonabelianTwistedCohomology})
  with local coefficients in \eqref{LocalCoefficientBundleForIteratedKTheory}
  is equivalently
  (still by Prop. \ref{ProofTwistedGeneralizedCohomology})
  integrally twisted iterated K-theory
  according to \cite{LindSatiWesterland16}:

  \vspace{-.3cm}
  $$
    \overset{
      \mathclap{
      \raisebox{3pt}{
        \tiny
        \color{darkblue}
        \bf
        \def\arraystretch{.9}
        \begin{tabular}{c}
          twisted
          \\
          iterated K-theory
        \end{tabular}
      }
      }
    }{
      \big(
        K^{\circ_{2r-1}}(\mathrm{ku})
      \big)^\tau(-)
    }
    \;\;
      \simeq
    \;\;
    H^{\tau}
    \Big(
      -;
        K^{\circ_{2r-2}}(\mathrm{ku})_0
    \Big)
    \,.
  $$
\end{example}

In twisted generalization of Example \ref{CohomotopyTheory}, we have:
\begin{example}[Twisted Cohomotopy theory {\cite[\S 2.1]{FSS19b}}]
  \label{JTwistedCohomotopyTheory}
  For $n \in \mathbb{N}$, consider the canonical action
  of the orthogonal group $\mathrm{O}(n+1)$ on the homotopy
  type of the $n$-sphere, via the defining
  action on the unit sphere in $\mathbb{R}^{n+1}$, which
  restricts along the canonical inclusion
  $\mathrm{O}(n) \xhookrightarrow{\;} \mathrm{O}(n+1)$ to the
  defining action of $\mathrm{O}(n)$ on the one-point compactification
  $\big(\mathbb{R}^n\big)^{\mathrm{cpt}} \,=\, S^n$.
  By Prop. \ref{GActionsAsFibrations}, this corresponds to
  local coefficient bundles \eqref{LocalCoefficientBundle} for
  twisting Cohomotopy theory (Example \ref{CohomotopyTheory}):
  \vspace{-2mm}
  \begin{equation}
    \label{LocalCoefficientBundleForTwistedCohomotopy}
    \begin{tikzcd}
      S^n
      \ar[r]
      &
      S^n \!\sslash\! \mathrm{O}(n)
      \ar[r]
      \ar[
        dr,
        phantom,
        "\mbox{\tiny\rm(hpb)}"
      ]
      \ar[
        d,
        "{\rho_J}"{left}
      ]
      &
      S^n \!\sslash\! \mathrm{O}(n+1)
      \ar[d]
      \\
      &
      B \mathrm{O}(n)
      \ar[r]
      &
      B \mathrm{O}(n+1)
    \end{tikzcd}
  \end{equation}

  \vspace{-3mm}
  \noindent
  The classifying map
  $B \mathrm{O}(n) \xrightarrow{\; J \; }B\mathrm{Aut}(S^n)$ of
  $\rho_J$ is the unstable \emph{J-homomorphism} (e.g. \cite[\S 4.4]{TamakiKono06}).
  For $X$ a smooth manifold of dimension $d \geq n+1$,
  and equipped with tangential $\mathrm{O}(n+1)$-structure
  (e.g. \cite[Def. 4.48]{SS20b})
  \vspace{-2mm}
  $$
    \xymatrix@R=12pt{
      X
      \ar[dr]_{T X}^-{\ }="t"
      \ar[rr]^-{ \tau }_>>>>>>{\ }="s"
      &&
      B \mathrm{O}(n+1)
      \ar[dl]^-{ B i }
      \\
      & B \mathrm{O}(d)
      \ar@{=>}^\simeq "s"; "t"
    }
  $$

  \vspace{-2mm}
  \noindent
  the $\tau$-twisted non-abelian Cohomology (Def. \ref{NonabelianTwistedCohomology})
  with local coefficients in
  \eqref{LocalCoefficientBundleForTwistedCohomotopy}
  is the \emph{tangentially twisted Cohomotopy theory} of
  \cite{FSS19b}\cite{FSS19c}\cite{SS20a}:

  \vspace{-.4cm}
  $$
    \overset{
      \mathclap{
      \raisebox{3pt}{
        \tiny
        \color{darkblue}
        \bf
        \def\arraystretch{.9}
        \begin{tabular}{c}
          tangentially twisted
          \\
          Cohomotopy
        \end{tabular}
      }
      }
    }{
      \pi^\tau(-)
    }
    \;\;
      :=
    \;\;
    H^\tau
    \big(
       -;
       S^n
    \big)
    \,.
  $$
  This twisted Cohomotopy theory in degree $n  = 4$ encodes,
  in particular, the
  shifted flux quantization condition of the C-field
  \cite[Prop. 3.13]{FSS19b} and the vanishing of the
  residual M5-brane anomaly \cite{SS20a};
  while
  J-twisted Cohomotopy in degree $n = 7$ encodes, in particular,
  level quantization of the Hopf-Wess-Zumino term on the M5-brane
  \cite{FSS19c}.
\end{example}

\medskip

\noindent {\bf Twisted non-abelian cohomology operations.}
In generalization of Def. \ref{NonAbelianCohomologyOperations},
we set:
\begin{defn}[Twisted non-abelian cohomology operation]
  \label{TwistedNonabelianCohomologyOperation}
  Given a transformation of local coefficient bundles
  \eqref{LocalCoefficientBundle} presented
  (under localization \eqref{LocalizationOfAModelCategoryAtWeakEquivalenes}
  to homotopy types \eqref{ClassicalHomotopyCategory})
  as a strictly commuting diagram
      \vspace{-2mm}
  \begin{equation}
    \label{MorphismOfLocalCoefficientBundles}
    \raisebox{20pt}{
    \xymatrix@R=1.5em{
      A_1 \!\sslash\! G_1
      \ar[d]_-{ \rho_1 }
      \ar[rr]^-{ \phi_t }
      &&
      A_2 \!\sslash\! G_2
      \ar[d]^-{ \rho_2 }
      \\
      B G_1
      \ar[rr]^-{ \phi_b }
      &&
      B G_2
    }
    }
    \;\;\;\;\;
    \in
    \TopologicalSpaces_{\mathrm{Qu}}
    \,,
  \end{equation}

  \vspace{-2mm}
  \noindent
  \emph{pasting composition} induces, \footnote{
    We postpone discussing the details
    of forming pasting composites to \cref{TwistedNonabelianCharacterMap},
    where they are provided by
    Def. \ref{RationalizationInTwistedNonAbelianCohomology}.
  }
  for each
  twist $ X \xrightarrow{\;\tau \;}  B G_1$ \eqref{ATwist},
  a map

  \vspace{-2mm}
  \begin{equation}
    \label{ATwistedNonabelianCohomologyOperation}
    \phi_\ast
    \;:\;
    \xymatrix{
      H^{\tau}
      (
        X;
        \,
        A_1
      )
    \ar[rrr]^-{
      \left(
        \phi_t
        \,\circ\,
        (-)
      \right)
      \,\circ\,
      (
        \rho_1
      )_\ast
    }
    &&&
      H^{\phi_b \circ \tau}
      (
        X;
        \,
        A_2
      )
    }
  \end{equation}

  \vspace{-1mm}
  \noindent
  of twisted non-abelian cohomology sets (Def. \ref{NonabelianTwistedCohomology}).
  We call these \emph{twisted non-abelian cohomology operations}.
\end{defn}

\begin{example}[Total non-abelian class of twisted cocycles]
  \label{TotalNonAbelianClassOfTwistedCocycle}
  For any coefficient bundle $\rho$ \eqref{LocalCoefficientBundle}
  there is the
  tautological transformation \eqref{MorphismOfLocalCoefficientBundles}
  to its total space regarded as fibered over the point:
  \vspace{-2mm}
  $$
    \xymatrix@R=12pt{
      A \!\sslash\! G
      \ar[d]_-{ \rho }
      \ar@{=}[rr]
      &&
      A \!\sslash\! G
      \ar[d]
      \\
      B G
      \ar[rr]
      &&
      \ast
      \,.
    }
  $$

  \vspace{-2mm}
\noindent  The induced twisted non-abelian cohomology operation \eqref{ATwistedNonabelianCohomologyOperation}
  goes from twisted cohomology to non-twisted cohomology
  with coefficient in the total space:
  \begin{equation}
    \label{PushforwardOfTwistedCohomologyAlongCoefficientBundle}
    \xymatrix{
      H^\tau
      (
        X;
        \,
        A
      )
      \ar[r]^-{
        \rho_\ast
      }
      &
      H
      \big(
        X;
        \,
        A \!\sslash\! G
      \big).
    }
  \end{equation}
\end{example}

\begin{example}[Hopf cohomology operation in twisted Cohomotopy {\cite[\S 2.3]{FSS19b}}]
  \label{HopfOperationInJTwistedCohomotopy}
  The quaternionic Hopf fibration
  $\!\xymatrix@C=2.5em{ S^7 \ar[r]|-{\; h_{\mathbb{H}} \;} & S^4 }\!$
  is equivariant under the symplectic unitary group
  $\mathrm{Sp}(2) \simeq \mathrm{Spin}(5)$, so that
  after passage to classifying spaces it induces
  a morphism of local coefficient bundles
  \eqref{MorphismOfLocalCoefficientBundles}
  for
  twisted Cohomotopy \eqref{LocalCoefficientBundleForTwistedCohomotopy}
  in degrees 4 and 7:
  \vspace{-2mm}
  \begin{equation}
    \label{EquivariantizedQuaternionicHopfFibration}
    \xymatrix@C=4em{
      S^7 \!\sslash \! \mathrm{Sp}(2)
      \ar[rr]_-{
          h_{\mathbb{H}} \sslash \mathrm{Sp}(2)}^-{
          \raisebox{3pt}{
            \tiny
            \color{darkblue}
            \bf
            \def\arraystretch{.9}
            \begin{tabular}{c}
              Borel-equivariantized
              \\
              quaternionic Hopf fibration
            \end{tabular}
          }
                 }
      \ar[d]_-{ J_7 }
      &&
      S^4 \!\sslash \! \mathrm{Sp}(2)
      \ar[d]^-{ J_4 }
      \\
      B \mathrm{Sp}(2)
      \ar@=[rr]
      &&
      B \mathrm{Sp}(2)\;.
    }
  \end{equation}

  \vspace{-1mm}
  \noindent
  Via \eqref{ATwistedNonabelianCohomologyOperation}, this
  induces for each Spin 8-manifold  $X$ equipped with
  tangential $\mathrm{Sp}(2)$-structure (Example \ref{ClassificationOfTangentialStructure})
  \vspace{-2mm}
  \begin{equation}
    \label{TangentialSp2Structure}
    \xymatrix@R=12pt{
      X
      \ar[dr]_{T X}^-{\ }="t"
      \ar[rr]^-{ \tau }_>>>>>>{\ }="s"
      &&
      B \mathrm{Sp}(2)
      \ar[dl]^-{ B i }
      \\
      & B \mathrm{O}(8)
      \ar@{=>}^\simeq "s"; "t"
    }
  \end{equation}

  \vspace{-2mm}
  \noindent
  a twisted non-abelian cohomology operation (Def. \ref{TwistedNonabelianCohomologyOperation})
  \vspace{-2mm}
  \begin{equation}
    \label{TheHopfTwistedCohomologyOperation}
    \xymatrix{
      \pi^{\tau_{\,7}}(X)
      \ar[rr]^-{
        (h_{\mathbb{H}}\sslash \mathrm{Sp}(2))_\ast
      }
      &&
      \pi^{\tau^4}(X)
    }
  \end{equation}

  \vspace{-2mm}
  \noindent
  in twisted non-abelian Cohomotopy theory
  (Example \ref{JTwistedCohomotopyTheory}).
  Lifting through the twisted non-abelian cohomology
  transformation \eqref{TheHopfTwistedCohomologyOperation}
  encodes vanishing of
  C-field flux up to C-field background charge
  \cite[Prop. 3.14]{FSS19b}.
\end{example}

\begin{example}[Twistorial Cohomotopy {\cite[\S 3.2]{FSS20a}} ]
  \label{TwistorialCohomotopy}
  The equivariantized Hopf morphism \eqref{EquivariantizedQuaternionicHopfFibration}
  of coefficient bundles factors through
  Borel-equivariantizations of the complex Hopf fibration
  $h_{\mathbb{C}}$
  followed by that of the \emph{twistor fibration} $t_{\mathbb{H}}$
  \vspace{-2mm}
  \begin{equation}
    \label{EquivariantizedTwistorFibration}
    \raisebox{20pt}{
    \xymatrix@C=4em@R=2.5em{
      S^7 \!\sslash \! \mathrm{Sp}(2)
      \ar[rr]_-{
          h_{\mathbb{C}} \sslash \mathrm{Sp}(2)
        }^-{
          \raisebox{3pt}{
            \tiny
            \color{darkblue}
            \bf
            \def\arraystretch{.9}
            \begin{tabular}{c}
              Borel-equivariantized
              \\
              complex Hopf fibration
            \end{tabular}
          }
          }
              \ar[d]_-{ J_{S^7} }
      &&
      \mathbb{C}P^3 \!\sslash \! \mathrm{Sp}(2)
      \ar[rr]_-{
          t_{\mathbb{H}} \sslash \mathrm{Sp}(2)
        }^-{
              \raisebox{3pt}{
            \tiny
            \color{darkblue}
            \bf
            \def\arraystretch{.9}
            \begin{tabular}{c}
              Borel-equivariantized
              \\
              twistor fibration
            \end{tabular}
          }
          }
         \ar[d]_-{ J_{\mathbb{CP}^3} }
      &&
      S^4 \!\sslash \! \mathrm{Sp}(2)
      \ar[d]^-{ J_{S^4} }
      \\
      B \mathrm{Sp}(2)
      \ar@=[rr]
      &&
      B \mathrm{Sp}(2)
      \ar@=[rr]
      &&
      B \mathrm{Sp}(2)
    }
    }
  \end{equation}
  The twisted non-abelian cohomology theory
  (Def. \ref{NonabelianTwistedCohomology})
  with local coefficients in the bundle appearing in this
  factorization
  is the \emph{Twistorial Cohomotopy} of \cite{FSS20a}
  \vspace{0mm}
  $$
    \overset{
      \mathclap{
      \raisebox{3pt}{
        \tiny
        \color{darkblue}
        \bf
        \def\arraystretch{.9}
        \begin{tabular}{c}
          Twistorial
          \\
          Cohomotopy
        \end{tabular}
      }
      }
    }{
      \mathcal{T}^\tau(-)
    }
    \;\;:=\;\;
    H^{\tau}
    \big(
      -;
      \mathbb{C}P^3
    \big)
    \,.
  $$

  \vspace{-1mm}
\noindent    Via \eqref{ATwistedNonabelianCohomologyOperation},
  the morphisms \eqref{EquivariantizedTwistorFibration} induce,
  for each spin 8-manifold  $X$ equipped with
  tangential $\mathrm{Sp}(2)$-structure \eqref{TangentialSp2Structure},
  twisted non-abelian cohomology operations (Def. \ref{TwistedNonabelianCohomologyOperation})
 \vspace{-2mm}
  \begin{equation}
    \label{TwistedCohomologyTransformationsDefiningTwistorialCohomotopy}
    \xymatrix@C=4em{
      \overset{
        \mathclap{
        \raisebox{3pt}{
          \tiny
          \color{darkblue}
          \bf
          \def\arraystretch{.9}
          \begin{tabular}{c}
            J-twisted
            \\
            7-Cohomotopy
          \end{tabular}
        }
        }
      }{
      \pi^{\tau^7}(X)
      }
      \ar[rr]^-{
        (h_{\mathbb{C}}\sslash \mathrm{Sp}(2))_\ast
      }
      &&
      \overset{
        \mathclap{
        \raisebox{3pt}{
          \tiny
          \color{darkblue}
          \bf
          \def\arraystretch{.9}
          \begin{tabular}{c}
            Twistorial
            \\
            Cohomotopy
          \end{tabular}
        }
        }
      }{
        \mathcal{T}^{\tau}(X)
      }s
      \ar[rr]^-{
        (t_{\mathbb{H}}\sslash \mathrm{Sp}(2))_\ast
      }
      &&
      \overset{
        \mathclap{
        \raisebox{3pt}{
          \tiny
          \color{darkblue}
          \bf
          \def\arraystretch{.9}
          \begin{tabular}{c}
            J-twisted
            \\
            4-Cohomotopy
          \end{tabular}
        }
        }
      }{
        \pi^{\tau^4}(X)
      }
    }
  \end{equation}

  \vspace{-1mm}
  \noindent
  between J-twisted non-abelian Cohomotopy theory
  (Example \ref{JTwistedCohomotopyTheory}) and
  Twistorial Cohomotopy.

  We turn to the differential refinement of this
  statement in \cref{CohomotopicalChernCharacter} below.
\end{example}

\newpage

%%%%%%%%%%%%%%%%%%%%%%%%%%%%%%%%%%%%%%
\section{Non-abelian de Rham cohomology}
  \label{NonAbelianDeRhamCohomologyTheory}
%%%%%%%%%%%%%%%%%%%%%%%%%%%%%%%%%%%%%%

We formulate
(twisted) non-abelian de Rham cohomology
(Def. \ref{NonabelianDeRhamCohomology},
Def. \ref{TwistedNonabelianDeRhamCohomology})
of differential forms with values in
$L_\infty$-algebras (Example \ref{ChevalleyEilenbergAlgebraOfLInfinityAlgebra})
and prove the (twisted) non-abelian de Rham theorem
(Theorem \ref{NonAbelianDeRhamTheorem},
Theorem \ref{TwistedNonAbelianDeRhamTheorem}),
as a consequence of the fundamental theorem of
dg-algebraic rational homotopy theory,
which we recall (Prop. \ref{FundamentalTheoremOfdgcAlgebraicRationalHomotopyTheory}).

%%%%%%%%%%%%%%%%%%%%%%%%%%%%%%%%%%%%%%%%%%%%%%%%%%%%%%%%
\subsection{Dgc-Algebras and $L_\infty$-algebras}
\label{DGCAlgebrasAndLInfinityAlgebras}
%%%%%%%%%%%%%%%%%%%%%%%%%%%%%%%%%%%%%%%%%%%%%%%%%%%%%%%%

Here we fix notation and conventions for the following
system of categories and functors:
\vspace{-3mm}
\begin{equation}
  \label{dgCategoriesAndAdjunctions}
  \hspace{-5mm}
  \raisebox{40pt}{
  \xymatrix@R=14pt@C=4em{
    \big(
      \overset{
        \raisebox{3pt}{
          \tiny
          \color{darkblue}
          \bf
          Def. \ref{NilpotentLInfinityAlgebras}
        }
      }{
      \LInfinityAlgebrasNil
      }
    \big)^{\mathrm{op}}
    \;
    \ar@{^{(}->}[d]
    \ar@{->}[r]^-{
      \overset{
        \raisebox{3pt}{
          \tiny
          \color{darkblue}
          \bf
          \eqref{NilpotentLInfinityAlgebrasEquivalentToConnectedSullivanModels}
        }
      }{
        \mathrm{CE}
      }
    }_-{ \simeq }
    &
    \overset{
      \raisebox{3pt}{
        \tiny
        \color{darkblue}
        \bf
        Def. \ref{SullivanModels}
      }
    }{
      \SullivanModelsConnected
    }
    \ar@{^{(}->}[d]
    \\
    \overset{
      \raisebox{3pt}{
        \tiny
        \color{darkblue}
        \bf
        Def. \ref{ChevalleyEilenbergAlgebraOfLInfinityAlgebra}
      }
    }{
      \big(
        \LInfinityAlgebras
      \big)^{\mathrm{op}}
    }
    \;
    \ar@{^{(}->}[r]^-{
      \overset{
        \raisebox{3pt}{
          \tiny
          \color{darkblue}
          \bf
          \eqref{CEAlgebraOfLInfinityAlgebra}
        }
      }{
        \mathrm{CE}
      }
    }
    &
    \overset{
      \raisebox{3pt}{
        \tiny \bf
        \color{darkblue}
        Def. \ref{CategoryOfdgAlgebras}
      }
    }{
      \dgcAlgebras{\mathbb{R}}
    }
    \;\;
    \ar@{<-}@<+5pt>[rr]^-{
      \overset{
        \raisebox{3pt}{
          \tiny \bf
          \color{darkblue}
          Def. \ref{FreeDifferentalGradedAlgebras}
        }
      }{
        \mathrm{Sym}
      }
    }
    \ar@<-5pt>[rr]^-{\bot}
    \ar[d]_-{\scalebox{0.6}{$
      \overset{
        \raisebox{3pt}{
          \tiny
          \color{darkblue}
          \bf
          Def. \ref{UnderlyingGradedCommutativeAlgebraOfdgcAlgebra}
        }
      }{
        \mathrm{GrddCmmttvAlgbr}
      }
    $}
    }
    &&
    \;\;
    \overset{
      \raisebox{3pt}{
        \tiny \bf
        \color{darkblue}
        Def \ref{CoochainComplexes}
      }
    }{
      \CochainComplexes
    }
    \ar[d]^-{\scalebox{0.6}{$
      \overset{
        \raisebox{3pt}{
          \tiny
          \color{darkblue}
          \bf
          Def. \ref{UnderlyingGradedVectorSpaceOfCochainComplex}
        }
      }{
        \mathrm{GrddVctrSpc}
      }
    $}
    }
    \\
    &
    \underset{
      \raisebox{-3pt}{
        \tiny \bf
        \color{darkblue}
        Def. \ref{CategoriesOfGradedAlgebras}
      }
    }{
      \GradedAlgebras
    }
    \;\;
    \ar@{<-}@<+5pt>[rr]^-{
      \overset{
        \raisebox{3pt}{
          \tiny \bf
          \color{darkblue}
          Def. \ref{FreeGradedCommutativeAlgebra}
        }
      }{
        \mathrm{Sym}
      }
    }
    \ar@<-5pt>[rr]^-{ \bot }
    &&
    \;\;
    \underset{
      \raisebox{-3pt}{
        \tiny \bf
        \color{darkblue}
        Def. \ref{CategoryOfGradedVectorSpaces}
      }
    }{
      \GradedVectorSpaces
    }
  }
  }
\end{equation}

\begin{remark}[Homotopical grading]
  Our grading conventions, to be detailed in the following,
  are strictly \emph{homotopy theoretic}, in that
  all algebraic data in degree $n$
  always corresponds to homotopy groups in that same degree:

 \vspace{1mm}
  \noindent {\bf (i)}
  Every graded-algebraic object discussed here corresponds,
  under the equivalences of rational homotopy theory
  laid out in \cref{RationalHomotopyTheory} below,
  to a rational space, such that algebraic generators in
  degree $n$ correspond to homotopy groups in
  the same degree $n$.
  Since homotopy groups of spaces are in non-negative degree
  $n \in \mathbb{N}$,
  all dg-algebraic objects discussed,
  both homological as well as cohomological,
  we take to be
  concentrated in non-negative degree.
  This implies that we take linear duality of (co)chain complexes
  (Def. \ref{DegreewiseLinearDual})
  to preserve the degree
  as opposed to changing it by a sign.

\vspace{1mm}
  \noindent {\bf (ii)}
  In particular, our $L_\infty$-algebras
  are in non-negative degree, hence are {\it connective},
  naturally accommodating
  (as in \cite{LadaMarkl95} \cite[\S 2.9]{BFM06})
  the rationalized Whitehead homotopy Lie algebras
  $\pi_\bullet(\Omega X) \otimes_{{}_{\mathbb{Z}}} \mathbb{R}$ of connected spaces $X$,
  with their natural non-negative grading
  induced from that of the homotopy groups of $\Omega X$.
  See Prop. \ref{WhiteheadLInfinityAlgebras} and
  Prop. \ref{RationalHomotopyGroupsInRationalWhiteheadLInfinityAlgebra} below.

\vspace{1mm}
  \noindent {\bf (iii)}
  Accordingly, all Chevalley-Eilenberg (CE) dgc-algebras
  (Ex. \ref{CEAlgebraOfLieAlgebra}) are taken to be in non-negative degree,
  as usual,
  so that their generators in degree $n$ correspond to dual homotopy groups in degree $n$.
  For example, the CE-algebra model for an Eilenberg-MacLane space
  $K(\mathbb{Z},n)$ has a single generator which is in degree $+n$ (Ex. \ref{RationalizationOfEMSpaces}).
\end{remark}

\medskip

\noindent {\bf Graded vector spaces.}
\begin{defn}[Connective graded vector spaces]
  \label{CategoryOfGradedVectorSpaces}
  {\bf (i)} We write
  \vspace{-2mm}
  \begin{equation}
    \label{GradedVectorSpaces}
    \GradedVectorSpaces
    \;\in\;
    \Categories
  \end{equation}

  \vspace{-2mm}
\noindent   for the category whose objects are $\mathbb{N}$-graded
  (i.e. non-negatively $\mathbb{Z}$-graded) vector spaces over the
  real numbers; and we write
  \begin{equation}
    \label{GradedVectorSpacesOfFiniteType}
    \GradedVectorSpacesFin
    \xymatrix{
      \ar@{^{(}->}[r]
      &
    }
    \GradedVectorSpaces
    \;\in\;
    \Categories
  \end{equation}
  for its full subcategory on those objects which are of
  \emph{finite type}, namely degree-wise finite-dimensional.

\noindent {\bf (ii)}  For $V \in \GradedVectorSpaces$ and $k \in \mathbb{N}$,
  we write
 \vspace{-1mm}
  $$
    V^k \;\in\; \VectorSpaces_{\scalebox{.5}{$\mathbb{R}$}}
  $$

   \vspace{-2mm}
\noindent  for the component vector space in degree $k$.
\end{defn}
\begin{example}[The zero-object in graded vector spaces]
  We write
     \vspace{-2mm}
  \begin{equation}
    \label{ZeroObjectInGradedVectorSpaces}
    0
    \;\in\;
    \GradedVectorSpaces
  \end{equation}

     \vspace{-2mm}
\noindent
  for the graded vector space which is the zero vector space
  in each degree. This is both the initial as well as the terminal
  object (hence the zero object) in $\GradedVectorSpaces$.
\end{example}

\begin{example}[Graded linear basis]
  \label{GradedLinearBasis}
  For $n_1, n_2, \cdots, n_k \in \mathbb{N}$ a finite sequence
  of non-negative integers, we write
     \vspace{-1mm}
  $$
    \big\langle
      \alpha_{n_1},
      \alpha_{n_2},
      \cdots,
      \alpha_{n_k}
    \big\rangle
    \;\in\;
    \GradedVectorSpacesFin
  $$

     \vspace{-2mm}
\noindent
  for the graded vector space (Def. \ref{CategoryOfGradedVectorSpaces})
  spanned by elements $\alpha_{n_i}$ in degree $n_i$, respectively.
\end{example}

\begin{defn}[Tensor product of graded vector spaces]
  \label{TensorProductOfGradedVectorSpaces}
  The category of $\GradedVectorSpaces$ (Def. \ref{CategoryOfGradedVectorSpaces})
  becomes a symmetric monoidal category under the graded tensor product
%  $$
%    \xymatrix@R=4pt{
%      \GradedVectorSpaces
%      \times
%      \GradedVectorSpaces
%     \ar[rr]
%      &&
%      \GradedVectorSpaces
%      \\
%      (V, W)
%      \ar@{}[rr]|-{
%        \mapsto
%      }
%      &&
%      V \otimes W
%    }
%  $$
  given by
     \vspace{-1mm}
  $$
    (V \otimes W)^k
    \;:=\;
    \underset{n_1 + n_2 = k}{\bigoplus}
    V^{n_1} \otimes \, W^{n_2}
    \,.
  $$

  \vspace{-2mm}
\noindent  and the symmetric braiding isomorphism given by
     \vspace{-2mm}
  \begin{equation}
    \label{BraidingIsomorphism}
    \xymatrix@R=5pt{
      V \otimes W
      \ar[rr]^-{\scalebox{0.7}{$\sigma^{V,W}$} }_-{ \simeq }
      &&
      W \otimes V
      \\
      \\
      V^{n_1} \otimes W^{n_2}
      \ar@{^{(}->}[uu]
      \ar[rr]^-{\scalebox{0.7}{$ \sigma^{V,W}_{n_1,n_2}$} }_-{ \simeq }
      &&
      W^{n_2}
        \otimes
      V^{n_1}
      \ar@{^{(}->}[uu]
      \\
      (v,w)
      \ar@{}[u]|-{
        \rotatebox{90}{$\in$}
      }
      \ar@{}[rr]|<<<<<<<<{\longmapsto}
      &&
      \!\!\!\!\!\!\!\!\!\!\!\!\!\!\!\!\!
      (-1)^{n_1 n_2} \cdot (w,v)
      \ar@{}[u]|-{
        \rotatebox{90}{$\in$}
      }
    }
  \end{equation}

     \vspace{-2mm}
\noindent
  We denote this by
  \begin{equation}
    \label{SymmetricMonoidalCategoryOfGradedVectorSpaces}
    \big(
      \GradedVectorSpaces
      ,
      \otimes,
      \sigma
    \big)
    \;\in\;
    \mathrm{SymmetricMonoidalCategories}\;.
  \end{equation}
\end{defn}

\begin{defn}[Degreewise linear dual]
  \label{DegreewiseLinearDual}
  For $V \in \GradedVectorSpacesFin$ (Def. \ref{CategoryOfGradedVectorSpaces})
  we write
  \vspace{-1mm}
  $$
    V^\vee \;\in\; \GradedVectorSpacesFin
  $$

  \vspace{-2mm}
  \noindent  for its degree-wise linear dual: \footnote{
    This is in contrast to the intrinsic duality $(-)^\ast$ in
    the monoidal category of
    graded vector spaces in \emph{unbounded} degree
    (not considered here),
    which
    instead goes along with inversion of the degree:
    $(V^\ast)^k \;=\; (V^{-k})^\ast$.
  }
   \vspace{-1mm}
  \begin{equation}
    \label{DegreeiwseDual}
    (V^\vee)^k
      \;:=\;
    (V^k)^\ast
    \,.
  \end{equation}
\end{defn}

\begin{defn}[Degree shift]
  \label{DegreeShift}
  For $V \in \GradedVectorSpaces$ (Def. \ref{CategoryOfGradedVectorSpaces})
  we write
    \vspace{-1mm}
  \begin{equation}
    \label{ShiftedGradedVectorSpace}
    \mathfrak{b}V
    \;\in\;
    \GradedVectorSpaces
  \end{equation}

  \vspace{-1mm}
  \noindent
  for the result of shifting degrees up by $1$:
   \vspace{-3mm}
  $$
    (\mathfrak{b}V)^k
    \;:=\;
    \left\{
    \begin{array}{lcl}
      V^{k-1} & \vert & k \geq 1,
      \\
      0 & \vert & k = 0.
    \end{array}
    \right.
     $$
\end{defn}

\medskip

\noindent {\bf Graded-commutative algebras.}
\begin{defn}[Graded-commutative algebras]
  \label{CategoriesOfGradedAlgebras}
  We write
    \vspace{-2mm}
  \begin{equation}
    \label{CategoryOfGradedAlgebras}
    \GradedAlgebras
    \;:=\;
    \mathrm{CommMonoids}
    \big(
      \GradedVectorSpaces, \otimes, \sigma
    \big)
    \;\in\;
    \Categories
  \end{equation}

    \vspace{-2mm}
\noindent
  for the category whose objects are
  non-negatively $\mathbb{Z}$-graded, graded-commutative
  unital algebras over the real numbers
  (hence commutative unital monoids with respect to the braided tensor product
  of Def. \ref{TensorProductOfGradedVectorSpaces}); and we write
    \vspace{-2mm}
  \begin{equation}
    \GradedAlgebrasFin
    \xymatrix{
      \ar@{^{(}->}[r]
      &
    }
    \GradedAlgebras
    \;\in\;
    \Categories
  \end{equation}

    \vspace{-2mm}
\noindent
  for its full sub-category in those objects which are of \emph{finite type},
  namely degree-wise finite dimensional.
\end{defn}

\begin{defn}[Underlying graded vector space]
  \label{UnderlyingGradedVectorSpaces}
  We write
     \vspace{-2mm}
  \begin{equation}
    \label{UnderlyingGradedVectorSpace}
    \xymatrix{
      \GradedAlgebras
      \ar[rr]^-{\scalebox{0.6}{$
        \mathrm{GrddVctrSpc}$}
      }
      &&
      \GradedVectorSpaces
    }
  \end{equation}

    \vspace{-2mm}
\noindent
  for the functor on graded algebras (Def. \ref{CategoriesOfGradedAlgebras})
  that forgets the algebra structure and
  remembers only the underlying graded vector space
  (Def. \ref{CategoryOfGradedVectorSpaces}).
\end{defn}

\begin{example}[Free graded-commutative algebras]
  \label{FreeGradedCommutativeAlgebra}
  For
  $ V \in \GradedVectorSpaces$
  (Def. \ref{CategoryOfGradedVectorSpaces}),
  we write
    \vspace{-2mm}
  \begin{equation}
    \label{SymmetricGradedAlgebra}
    \mathrm{Sym}(V)
    \;\in\;
    \GradedAlgebras
  \end{equation}

  \vspace{-2mm}
\noindent
  for the graded-commutative algebra (Def. \ref{CategoriesOfGradedAlgebras})
  freely generated by $V$, hence that whose underlying graded vector space
  \eqref{UnderlyingGradedVectorSpace} is
     \vspace{-1mm}
  $$
    \mathrm{GrddVctrSpc}
    \big(
      \mathrm{Sym}(V)
    \big)
    \;=\;
    \mathbb{R}
      \;\oplus\;
    V
      \;\oplus\;
    \big( V \otimes V \big)_{/\mathrm{Sym}(2)}
      \;\oplus\;
    \big( V \otimes V \otimes V\big)_{/\mathrm{Sym}(3)}
      \;\oplus\;
    \cdots
    \,,
  $$

    \vspace{-1mm}
\noindent
  where the symmetric groups $\mathrm{Sym}(n)$ act via the
  braiding \eqref{BraidingIsomorphism}.
\end{example}

\begin{example}[Graded Grassmann algebra]
  \label{GrassmannAlgebra}
  For
  $ V \;\in\; \GradedVectorSpaces$
  (Def. \ref{CategoryOfGradedVectorSpaces}),
  we write
    \vspace{-2mm}
  $$
    \wedge^\bullet V
    \;:=\;
    \mathrm{Sym}
    \big(
      \mathfrak{b}V
    \big)
    \;\in\;
    \GradedAlgebras
  $$

  \vspace{-2mm}
  \noindent
  for the free graded-commutative algebra (Def. \ref{FreeGradedCommutativeAlgebra})
  on $V$ shifted up in degree (Def. \ref{DegreeShift});
  and we call this the \emph{graded Grassmann-algebra} on $V$.
\end{example}

\begin{example}[Graded polynomial algebra]
  \label{GradedPolynomialAlgebra}
  For $n_1, n_2, \cdots, n_k \in \mathbb{N}$ a finite sequence
  of non-negative integers, we write
   \vspace{-1mm}
  $$
    \mathbb{R}
    \big[
      \alpha_{n_1},
      \alpha_{n_2},
      \cdots,
      \alpha_{n_k}
    \big]
    \;:=\;
    \mathrm{Sym}
    \Big(
      \big\langle
        \alpha_{n_1},
        \alpha_{n_2},
        \cdots,
        \alpha_{n_k}
      \big\rangle
    \Big)
    \;\in\;
    \GradedAlgebrasFin
  $$

  \vspace{-2mm}
\noindent  for the free graded-commutative algebras (Def. \ref{FreeGradedCommutativeAlgebra})
  the graded vector space spanned by the $\alpha_{n_i}$ (Def. \ref{GradedLinearBasis}).
\end{example}

\begin{remark}[Incarnations of Grassmann algebras]
  \label{IncarnationsOfGrassmannAlgebras}
  With these notation conventions from
  Examples \ref{FreeGradedCommutativeAlgebra}, \ref{GrassmannAlgebra},
  \ref{GradedPolynomialAlgebra},
  an ordinary Grassmann algebra on $k$ generators is
  equivalently:
     \vspace{-1mm}
  $$
    \wedge^\bullet
    \big(
      \mathbb{R}^k
    \big)
    \;=\;
    \mathrm{Sym}
    \big(
      \mathfrak{b}\mathbb{R}^k
    \big)
    \;=\;
    \mathbb{R}
    \big[
      \theta_1^{(1)},
      \theta_1^{(2)},
      \cdots,
      \theta_1^{(k)}
    \big]
    \,.
  $$
\end{remark}

\medskip

\noindent {\bf Cochain complexes.}
\begin{defn}[Connective cochain complexes]
  \label{CoochainComplexes}
  We write
     \vspace{-2mm}
  $$
    \CochainComplexes
    \;\in\;
    \Categories
  $$

  \vspace{-2mm}
  \noindent
  for the category of cochain complexes
  (i.e. with differential of degree $+1$)
  of real
  vector spaces in non-negative degree.
\end{defn}
\begin{defn}[Underlying graded vector space]
  \label{UnderlyingGradedVectorSpaceOfCochainComplex}
  We write
   \vspace{-2mm}
  \begin{equation}
    \label{ForgetfulFunctorFromCochainComplexesToGradedVectorSpaces}
    \xymatrix{
      \CochainComplexes
      \ar[rr]^-{\scalebox{0.6}{$ \mathrm{GrddVctrSpc}$} }
      &&
      \GradedVectorSpaces
    }
  \end{equation}

   \vspace{-2mm}
\noindent
  for the forgetful functor
  on connective cochain complexes (Def. \ref{CoochainComplexes})
  which forgets the differential and remembers only the
  underlying
  connective graded vector space (Def. \ref{CategoryOfGradedVectorSpaces}).
\end{defn}

\begin{defn}[Tensor product on cochain complexes]
  \label{TensorProductOnCochainComplexes}
  The tensor product and braiding of graded vector spaces
  from Def. \ref{TensorProductOfGradedVectorSpaces}
  lifts, through \eqref{ForgetfulFunctorFromCochainComplexesToGradedVectorSpaces},
  to a tensor product and braiding on
  $\CochainComplexes$ (Def. \ref{CoochainComplexes}), making
  it a symmetric monoidal category:
   \vspace{-2mm}
  \begin{equation}
    \big(
      \CochainComplexes
      ,
      \otimes,
      \sigma
    \big)
    \;\in\;
    \mathrm{SymmetricMonoidalCategories}\;.
  \end{equation}
\end{defn}

\noindent {\bf Differential graded commutative algebras.}
\begin{defn}[Connective differential graded commutative algebras {\cite[V.3.1]{GelfandManin96}}]
  \label{CategoryOfdgAlgebras}
  We write
     \vspace{-1mm}
  $$
    \dgcAlgebras{\mathbb{R}}
    \;:=\;
    \mathrm{CommMonoids}
    \big(
      \CochainComplexes
      , \otimes,
      \sigma
    \big)
    \;\in\;
    \Categories
  $$

   \vspace{-1mm}
\noindent
  for the category whose objects are differential-graded,
  graded-commutative, unital algebras over the real numbers
  concentrated in non-negative degrees
  (hence commutative unital monoids in the symmetric monoidal
  category of Def. \ref{TensorProductOnCochainComplexes}).
\end{defn}

\begin{defn}[Underlying graded-commutative algebra]
  \label{UnderlyingGradedCommutativeAlgebraOfdgcAlgebra}
  We write
     \vspace{-2mm}
  \begin{equation}
    \label{UnderlyingGradedCommutativeAlgebra}
    \xymatrix{
      \dgcAlgebras{\mathbb{R}}
      \ar[rrr]^-{\scalebox{0.6}{$
        \mathrm{GrddCmmttvAlgbr}$}
      }
      &&&
      \GradedAlgebras
    }
  \end{equation}

     \vspace{-1mm}
\noindent
  for the functor on
  dgc-algebras (Def. \ref{CategoryOfdgAlgebras})
  that forgets the differential and
  remembers only the underlying
  graded-commutative algebra
  (Def. \ref{CategoriesOfGradedAlgebras}).
\end{defn}

\begin{defn}[Free differential graded algebras]
  \label{FreeDifferentalGradedAlgebras}
  For $V^\bullet$ in $\CochainComplexes$ (Def. \ref{CoochainComplexes})
  we write
     \vspace{-1mm}
  $$
    \mathrm{Sym}(V^\bullet)
    \;\in\;
    \dgcAlgebras{\mathbb{R}}
  $$

   \vspace{-1mm}
\noindent
  for the free differential graded-commutative algebra on $V^\bullet$,
  (Def. \ref{CategoryOfdgAlgebras}), hence whose underlying
  graded-commutative algebra
  algebra \eqref{UnderlyingGradedCommutativeAlgebra}
  is as in Example \ref{FreeGradedCommutativeAlgebra}.
\end{defn}

\begin{example}[Initial algebra]
  \label{InitialGradedAlgebra}
  The real algebra of real numbers, regarded as
  concentrated in degree-0
  \vspace{-1mm}
  $$
    \mathbb{R}
    \;\in\;
    \GradedAlgebras
    \xymatrix{
      \ar@{^{(}->}[r]
      &
    }
    \dgcAlgebras{\mathbb{R}}
  $$

  \vspace{-1mm}
  \noindent
  is the \emph{initial} object:
  For any other $A \in \GradedAlgebras$ (Def. \ref{CategoryOfGradedAlgebras})
  or $\in \dgcAlgebras{\mathbb{R}}$ (Def. \ref{CategoryOfdgAlgebras})
  there is a unique morphism
  $
    \xymatrix{
      \mathbb{R}
      \; \ar@{^{(}->}[r]^-{ i_{\mathbb{R}} }
      &
      A
    }
  $

  \vspace{-1mm}
  \noindent
  (because our algebras are unital and homomorphims
  need to preserve the
  unit element).
\end{example}

\begin{example}[The terminal algebra]
  \label{TheTerminalAlgebra}
  We write
     \vspace{-2mm}
  \begin{equation}
    \label{TheZeroAlgebra}
    0
    \;\in\;
    \GradedAlgebras
    \xymatrix{
      \ar@{^{(}->}[r]
      &
    }
    \dgcAlgebras{\mathbb{R}}
  \end{equation}

     \vspace{-2mm}
\noindent
  for the unique graded-commutative algebra (Def. \ref{CategoriesOfGradedAlgebras})
  or dgc-algebra (Def. \ref{CategoryOfdgAlgebras})
  whose underlying graded vector space (Def. \ref{UnderlyingGradedVectorSpaces})
  is the zero-vector space \footnote{
    Notice that the algebra $0$ \eqref{TheZeroAlgebra}
    is indeed a unital algebra \eqref{CategoryOfGradedAlgebras}.
  } \eqref{ZeroObjectInGradedVectorSpaces}.
  This is the terminal object \footnote{
    Beware that the corresponding statement in \cite[p. 335]{GelfandManin96}
    is incorrect.
  } in $\GradedAlgebras$: For every
  $A \in \GradedAlgebras$, there is a unique morphism
  $
    \xymatrix{A \ar[r]^-{\exists !} & 0}
    \,.
  $
\end{example}

\begin{example}[Product and co-product algebras]
  \label{ProductAndCoproductAlgebras}
  In the categories $\GradedAlgebras$ (Def. \ref{CategoriesOfGradedAlgebras})
  and $\dgcAlgebras{\mathbb{R}}$ (Def. \ref{CategoryOfdgAlgebras}):

\noindent  {\bf (i)} the coproduct is given by the tensor product (Def. \ref{TensorProductOfGradedVectorSpaces}),

\noindent  {\bf (ii)} the product is given by the direct sum

  \noindent on underlying graded vector spaces (Def. \ref{UnderlyingGradedVectorSpaces}).

\noindent (The first follows by \cite[p. 478, Cor. 1.1.9]{Johnstone02}, while the second holds since \eqref{UnderlyingGradedVectorSpace}
  is a right adjoint.)
\end{example}

\begin{example}[Smooth de Rham complex (e.g. {\cite{BottTu82}})]
  \label{SmoothdeRhamComplex}
  For $X$ be a smooth manifold, its de Rham algebra
  of smooth differential forms is a dgc-algebra in the sense of
  Def. \ref{CategoryOfdgAlgebras}, to be denoted here:
  \vspace{-2mm}
  $$
    \Omega^\bullet_{\mathrm{dR}}(X)
    \;\in\;
    \dgcAlgebras{\mathbb{R}}
    \,.
  $$
\end{example}

\begin{example}[Chevalley-Eilenberg algebras of Lie algebras]
  \label{CEAlgebraOfLieAlgebra}
  For $(\mathfrak{g},[-,-])$ a finite-dimensional real
  Lie algebra, its Chevalley-Eilenberg algebra
  is a dgc-algebra (Def. \ref{CategoryOfdgAlgebras}):
     \vspace{-2mm}
  $$
    \mathrm{CE}(\mathfrak{g})
    \;:=\;
    \big(
      \wedge^\bullet \mathfrak{g}^\ast
      \,,\;
      d\vert_{\wedge^1 \mathfrak{g}^\ast} = [-,-]^\ast
    \big)
    \;\in\;
    \dgcAlgebras{\mathbb{R}}
  $$

  \vspace{-1mm}
  \noindent
  with underlying graded-commutative algebra (Def. \ref{CategoriesOfGradedAlgebras})
  the Grassmann algebra on
  the linear dual space $\mathfrak{g}^\ast$
  (Def. \ref{GrassmannAlgebra}, Remark \ref{IncarnationsOfGrassmannAlgebras}),
  and with differential given on $\wedge^1 \mathfrak{g}^\ast$
  by the linear dual of the Lie bracket.
    More explicitly, for $\{v_a\}_{a = 1}^{\mathrm{dim}_{\mathbb{R}}(\mathfrak{g})}$
   a linear basis for the underlying vector space of
   the Lie algebra
    \vspace{-2mm}
   \begin{equation}
     \label{LinearBasisForLieAlgebra}
     \mathfrak{g}
       \;\simeq\;
     \langle
       v_1, v_2,
       \cdots,
       v_{ \mathrm{dim}(\mathfrak{g}) }
     \rangle
     \,,
   \end{equation}

    \vspace{-2mm}
\noindent
   with Lie brackets
   \begin{equation}
     \label{LieAlgebraStructureConstants}
     [v_a, v_b]
     \;=\;
     f_{a b}^c
     v_c
     \,,
     \;\;\;\;
     \mbox{
       for structure constants
       $f_{a b}^c \,\in\,s \mathbb{R}$
     }
   \end{equation}
  we have
  \begin{equation}
    \label{CEAlgebraOfLieAlbebraInComponents}
    \mathrm{CE}(\mathfrak{g})
    \;\simeq\;
    \mathbb{R}
    \big[
      \theta^{(1)}_1, \theta^{(2)}_1, \cdots
      \theta^{(\mathrm{dim}(\mathfrak{g}))}_1
    \big]
    \!\big/\!
    \big(
    d\,\theta^{(c)}_1
      \,=\,
    f_{a b}{}^c\, \theta^{(b)}_1 \wedge \theta^{(a)}_1
    \big)
    \,.
  \end{equation}
  One observes that the Jacobi identity on $[-,-]$ is equivalent to the condition
  that the differential $d := [-,-]^\ast$ squares to zero,
  so that \eqref{CEAlgebraOfLieAlbebraInComponents} being
  a dgc-algebra is actually equivalent to $(\mathfrak{g}, [-,-])$
  being a Lie algebra.

  This construction is evidently contravariantly functorial and
  constitutes a full subcategory inclusion
  \vspace{-2mm}
  \begin{equation}
    \label{FullSubcategoryInclusionOfLieAlgebrasIntodgcAlgebras}
    \xymatrix{
      \LieAlgebras
      \;
      \ar@{^{(}->}[r]^-{\mathrm{CE}}
      &
      \;
      \big(
        \dgcAlgebras{\mathbb{R}}
      \big)^{\mathrm{op}}
    }\,,
  \end{equation}

  \vspace{-2mm}
  \noindent
  meaning that, in addition, homomorphisms of Lie algebras are
  in natural bijection to dgc-algebra morphisms between their
  CE-algebras.
\end{example}
This observation is the golden route to approaching $L_\infty$-algebras:

\medskip

\noindent {\bf $L_\infty$-algebras.}
\begin{defn}[Chevalley-Eilenberg algebras of $L_\infty$-algebras {\cite[Thm . 2.3]{LadaMarkl95}\cite[Def. 13]{SatiSchreiberStasheff08}\cite[\S 2]{BFM06}}]
  \label{ChevalleyEilenbergAlgebraOfLInfinityAlgebra}
  In direct generalization of \eqref{FullSubcategoryInclusionOfLieAlgebrasIntodgcAlgebras},
  consider
  those $A \in \dgcAlgebras{\mathbb{R}}$ (Def. \ref{CategoryOfdgAlgebras})
  whose underlying graded-commutative algebra \eqref{UnderlyingGradedCommutativeAlgebra}
  is free (Example \ref{FreeGradedCommutativeAlgebra},
  Remark \ref{IncarnationsOfGrassmannAlgebras})
  on the degreewise dual $\mathfrak{b}\mathfrak{g}^\vee$ (Def. \ref{DegreewiseLinearDual})
  of the degree shift $\mathfrak{b}\mathfrak{g}$ (Def. \ref{DegreeShift})
  of some connective finite-type graded vector space (Def. \ref{CategoryOfGradedVectorSpaces})
    \vspace{-1mm}
  \begin{equation}
    \label{GradedVectorSpaceUnderlyingLInfinityAlgebra}
    \mathfrak{g} \;\in\; \GradedVectorSpacesFin
  \end{equation}

    \vspace{-2mm}
 \noindent in that
  \vspace{-1mm}
  \begin{equation}
    \label{SemiFreedgcAlgebra}
    A
    \;:=\;
    \big(
      \wedge^\bullet \mathfrak{g}^\vee
      \,,\,
      d
    \big)
    \;:=\;
    \big(
      \mathrm{Sym}(\mathfrak{b}\mathfrak{g}^\vee)
      \,,\,
      d
    \big)
    \;\in\;
    \dgcAlgebras{\mathbb{R}}
    \,.
  \end{equation}

  \vspace{0mm}
  \noindent
  In this case the differential $d$ restricted to
  $\wedge^1 \mathfrak{g}^\vee$ defines, under linear dualization,
  a sequence of $n$-ary graded-symmetric multilinear maps
  $\{-,\cdots,-\}$ on $\mathfrak{g}$:
  \vspace{-1mm}
  \begin{equation}
    \label{LInfinityBracketsEncodedInCEDifferential}
    \xymatrix@R=1pt@C=6pt{
      d\vert_{\wedge^1 \mathfrak{g}^\vee}(-)
      &
      \ar@{}[r]|-{=}
      &
      &
      \{-\}^\ast
      \ar@{}[r]|-{+}
      &
      \{-,-\}^\ast
      \ar@{}[r]|-{+}
      &
      \{-,-,-\}^\ast
      \ar@{}[r]|-{+}
      &
      \cdots
      \\
      \wedge^1 \mathfrak{g}^\vee
      \ar[rr]^-{ d }
      &
      &
      &
      \wedge^1 \mathfrak{g}^\vee
      \ar@{}[r]|-{\oplus}
      &
      \wedge^2 \mathfrak{g}^\vee
      \ar@{}[r]|-{\oplus}
      &
      \wedge^3 \mathfrak{g}^\vee
      \ar@{}[r]|-{\oplus}
      &
      \cdots
      &
      \ar@{}[r]|-{=}
      &
      \wedge^\bullet \mathfrak{g}^\vee
      \;=\;
      \mathrm{Sym}
      \big(
        \mathfrak{b}\mathfrak{g}^\vee
      \big)
      \,,
    }
  \end{equation}

  \vspace{-1mm}
  \noindent  and the condition $d \circ d = 0$ imposes a sequence
  of compatibility conditions on these brackets,
  generalizing the Jacobi identity in Example \ref{CEAlgebraOfLieAlgebra}.
  The corresponding graded skew-symmetric $n$-ary brackets
  (\cite[(3)]{LadaStasheff93})
  \begin{equation}
    \label{BracketsOnLInfinityAlgebra}
    [a_1, \cdots, a_n]
    \;\;:=\;\;
    (-1)^{
      n +
      \sum_{i \leq n/2}
      \mathrm{deg}(a_i)
    }
    \{a_1, \cdots, a_n\}
  \end{equation}
  subject to these conditions give
  $\mathfrak{g}$ the structure of an \emph{$L_\infty$-algebra}
  (or \emph{strong homotopy Lie algebra}):
  \vspace{-1mm}
  \begin{equation}
    \label{LInfinityAlgebra}
    \Big(
      \mathfrak{g}
      \,,\,
      [-], [-,-], [-,-,-], \cdots
    \Big)
    \;\in\;
    \LInfinityAlgebras
    \,,
  \end{equation}

  \vspace{-1mm}
  \noindent
  which makes $A$ in \eqref{SemiFreedgcAlgebra} its
  Chevalley-Eilenberg algebra:
  \vspace{-1mm}
  \begin{equation}
    \label{CEAlgebraOfLInfinityAlgebra}
    \begin{aligned}
    \mathrm{CE}
    (
      \mathfrak{g}
    )
    & :=\;
    \big(
      \wedge^\bullet \mathfrak{g}^\vee
      \,,\,
      d = \{-\}^\ast + \{-,-\}^\ast + \{-,-,-\}^\ast + \cdots
    \big)
    \\
    \phantom{\mathclap{\vert^{\vert^{\vert}}}}
    & \phantom{:}=\;
    \big(
      \mathrm{Sym}(\mathfrak{b}\mathfrak{g}^\vee)
      \,,\,
      d_{\mathrm{CE}}
    \big)
    \,.
    \end{aligned}
  \end{equation}

  \vspace{-1mm}
  \noindent
  This construction constitutes a full subcategory inclusion
  \vspace{-1mm}
  \begin{equation}
    \label{FullSubcategoryInclusionOfLieAlgebras}
    \xymatrix{
      \LInfinityAlgebras
      \;
      \ar@{^{(}->}[r]^-{\mathrm{CE}}
      &
      \;
      \big(
        \dgcAlgebras{\mathbb{R}}
      \big)^{\mathrm{op}}
      \,.
    }
  \end{equation}

  \vspace{-1mm}
  \noindent
  into the category of dgc-algebras
  of the category of connective finite-type $L_\infty$-algebras,
  with the homotopy-correct morphisms between them
  (known as
  ``weak maps'' \cite[Rem. 5.4]{LadaMarkl95},
  ``sh maps'' \cite[\S 2.11]{Merkulov02}
  or
  ``$L_\infty$-morphisms'' \cite[p. 12]{Kontsevich97}).
\end{defn}

\begin{example}[Differential graded Lie algebras]
  \label{DifferentialGradedLieAlgebras}
  A differential graded Lie algebra is
  an $L_\infty$-algebra \eqref{LInfinityAlgebra}
  whose only possibly non-vanishing brackets
  are the unary bracket $\partial := [-]$ (its differential)
  and the binary bracket $[-,-]$ (its graded Lie bracket).
  In further specialization, a plain Lie algebra
  (Example \ref{CEAlgebraOfLieAlgebra})
  is an $L_\infty$-algebra/dg-Lie algebra concentrated in degree 0:
  \vspace{-1mm}
  \begin{equation}
    \label{DifferentialGradedLieAlgebrasInsideLInfinityAlgebras}
    \xymatrix{
      \LieAlgebras \;
      \ar@{^{(}->}[r]
      &
      \DifferentialGradedLieAlgebras \;
      \ar@{^{(}->}[r]
      &
      \LInfinityAlgebras
      \,.
    }
  \end{equation}
\end{example}

\begin{example}[Line Lie $n$-algebra]
  \label{LineLienPlusOneAlgebras}
  For $n \in \mathbb{N}$, the
  \emph{line Lie $(n+1)$-algebra} is the $L_\infty$-algebra
  (Def. \ref{ChevalleyEilenbergAlgebraOfLInfinityAlgebra})

  \vspace{-.3cm}
  \begin{equation}
    \label{LineLienPlus1Algebra}
    \mathfrak{b}^n \mathbb{R}
    \;\in\;
    \LInfinityAlgebras
  \end{equation}

  \vspace{-1mm}
  \noindent
  whose Chevalley-Eilenberg algebra \eqref{CEAlgebraOfLInfinityAlgebra}
  is the polynomial dgc-algebra (Example \ref{PolynomialAlgebraAndCellAttachment}) on a single
  closed generator in degree $n+1$:
  \vspace{-1mm}
  \begin{equation}
    \label{CEAlgebraOfLineLienAlgebra}
    \mathrm{CE}
    \big(
      \mathfrak{b}^n\mathbb{R}
    \big)
    \;:=\;
    \mathbb{R}
    [
      c_{n+1}
    ]
    \!\big/\!
    (
      d\, c_{n+1} \, = 0
    )\;.
  \end{equation}
  More generally, for $V \,\in\, \VectorSpaces^{\mathrm{fin}}_{\scalebox{.5}{$\mathbb{R}$}}$,
  we have

  \vspace{-.5cm}
  \begin{equation}
    \label{VectorSpaceLienPlus1Algebra}
    \mathfrak{b}^n V
    \;\simeq\;
    \underset{
      \mathrm{dim}(V)
    }{\bigoplus}
    \mathfrak{b}^n \mathbb{R}
    \;\;\;\in\;
    \LInfinityAlgebras
    \,,
    \;\;\;\;\;\;
    \mbox{with}
    \;\;\;\;
    \mathrm{CE}
    \big(
      \mathfrak{b}^n V
    \big)
    \;\simeq\;
    \mathbb{R}
    \big[
      c^{\scalebox{.5}{$(1)$}}_{n+1},
      c^{\scalebox{.5}{$(2)$}}_{n+1},
        \cdots,
      c^{\scalebox{.5}{$(\mathrm{dim}V)$}}_{n+1}
    \big]
    \!\big/\!
   \scalebox{0.7}{$ \left(
      \!\!\!
      \def\arraystretch{1.2}
      \begin{array}{l}
        d\, c^{\scalebox{.5}{$(1)$}}_{n+1}\;\, \, = 0\,,
        \\
        \vdots
        \\
        d\, c^{\scalebox{.5}{$(\mathrm{dim} V )$}}_{n+1} \, = 0
      \end{array}
    \right)
    $}
    \,.
  \end{equation}
\end{example}

\begin{example}[String Lie 2-algebra {\cite[\S 5]{BCSS07}\cite[\S 1.2]{Henriques08}\cite[App.]{FSS12a}}]
  \label{StringLie2Algebra}
  Let $\mathfrak{g} \in \LieAlgebras$ be semisimple
  (such as $\mathfrak{g} = \mathfrak{su}(n+1), \mathfrak{so}(n+3)$,
  for $n \in \mathbb{N}$), hence
  equipped with a non-degenerate, symmetric, $\mathfrak{g}$-invariant
  bilinear form (``Killing form'')
   \vspace{-2mm}
  \begin{equation}
    \label{KillingForm}
    \xymatrix@C=3em{
      \mathfrak{g} \otimes \mathfrak{g}
      \ar[r]^-{ \langle -,-\rangle }
      &
      \mathbb{R}
    }
    \,.
  \end{equation}

  \vspace{-2mm}
  \noindent
  Then the element
  $$
    \mu \;:=\;
    \big\langle
      -,\, [-,-]
    \big\rangle
    \;\;\;
    \in
    \;
    \mathrm{CE}(\mathfrak{g})
  $$
  in the Chevalley-Eilenberg \eqref{CEAlgebraOfLieAlgebra}
  is closed (is a Lie algebra cocycle)
   \vspace{-1mm}
  $$
    d \mu \;=\; 0
    \,.
  $$

  \vspace{-1mm}
  \noindent
  In terms of a linear basis $\{v_a\}$ \eqref{LinearBasisForLieAlgebra}
  with structure constants $\{f_{a b}^c\}$ \
  \eqref{LieAlgebraStructureConstants} and inner product
  $
    k_{a b} \;:=\; \langle v_a, v_b\rangle
  $
  we have, in terms of \eqref{CEAlgebraOfLieAlbebraInComponents}:
    \vspace{-1mm}
  $$
    \mu
      \;:=\;
    f_{a b}{}^{c'} k_{c' c} \,
    \theta^{(c)}_1
      \wedge
    \theta^{(b)}_1
      \wedge
    \theta^{(a)}_1
    \,.
  $$

\noindent    Hence we get an $L_\infty$-algebra (Def. \ref{ChevalleyEilenbergAlgebraOfLInfinityAlgebra})
\vspace{-1mm}
  \begin{equation}
    \label{TheStringLie2Algebra}
    \mathfrak{string}_{\mathfrak{g}}
    \;\;
    \in
    \;
    \LInfinityAlgebras
  \end{equation}

  \vspace{-1mm}
  \noindent
  with the following Chevalley-Eilenberg algebra
  \eqref{CEAlgebraOfLInfinityAlgebra}:
   \vspace{-2mm}
  \begin{equation}
    \label{CEAlgebraOfStringLie2Algebra}
    \mathrm{CE}
    \big(
      \mathfrak{string}_{\mathfrak{g}}
    \big)
    \;:=\;
    \mathbb{R}
    \left[
      \!\!\!\!\!
      \begin{array}{c}
        \{\theta^a_1\},
        \\
        b_2^{\phantom{(a)}}
      \end{array}
      \!\!\!\!\!
    \right]
    \!\big/\!
   \scalebox{0.8}{$ \left(
    \!
    \begin{aligned}
      d\,\theta^{(c)}_1
        & =
        f_{a b}^c
        \,
        \theta^{(b)}_1
          \wedge
        \theta^{(a)}_1
      \\
      d\, b_2 & =
        \underset{
          = \, \mu
        }{
        \underbrace{
        f_{a b}^{c'}
        \,
        k_{c' c}
        \,
        \theta^{(c)}_1
          \wedge
        \theta^{(b)}_1
          \wedge
        \theta^{(a)}_1
        }
        }
    \end{aligned}
        \!
    \right)
    $}
    .
  \end{equation}

  \vspace{-1mm}
  \noindent
  This is known as the \emph{string Lie 2-algebra},
  since it is
  \cite{BCSS07}\cite{Henriques08}
  the $L_\infty$-algebra of the
  String 2-group of Ex. \ref{NonAbelianCohomologyInDegree2}.
\end{example}

\medskip

\noindent {\bf Sullivan models and nilpotent $L_\infty$-algebras.}
\begin{example}[Polynomial dgc-algebras]
  \label{PolynomialAlgebraAndCellAttachment}
  For $A \in \dgcAlgebras{\mathbb{R}}$ (Def. \ref{CategoryOfdgAlgebras}),
  and
    \vspace{-2mm}
\begin{equation}
    \label{CocycleToBeTrivialized}
    \mu \;\in\; A^{n+1} \;\subset\; A\,,
    \phantom{AAAA}
    d\,\mu \;=\; 0
  \end{equation}

    \vspace{-2mm}
\noindent
  a closed element of homogeneous degree $n+1$,
  we write
    \vspace{-2mm}
  \begin{equation}
    \label{PolynomialdgcAlgebraExpression}
    A\!
    \big[
      \alpha_n
    \big]
    \big/
    \!
    \left(
    \!
    \begin{aligned}
      d\, \alpha_n & = \mu
    \end{aligned}
    \!\!\!\!\!\!\!
    \right)
    \;\;
    \in
    \;
    \dgcAlgebras{\mathbb{R}}
  \end{equation}

    \vspace{-2mm}
\noindent
  for the dgc-algebra obtained by adjoining a generator
  $\alpha_n$ of degree $n$ to the underlying graded-commutative algebra
  \eqref{UnderlyingGradedCommutativeAlgebra}
  of $A$
  and extending the differential from $A$ to $A\big[\alpha_n\big]$
  by taking its value on the new
  generator to be $\mu$.
  The polynomial dgc-algebra \eqref{PolynomialdgcAlgebraExpression}
  receives a canonical algebra inclusion of $A$:
  %(the unique $A$-algebra homomorphism):
\vspace{-3mm}
  \begin{equation}
    \label{InjectingAIntoPolynomialdgcAlgebraOverA}
    \xymatrix{
      A
      \;
      \ar@{^{(}->}[r]^-{ i_A }
      &
      \;
    A\!
    [
      \alpha_n
    ]
    \big/
    \!
    (
          d\, \alpha_n  = \mu
   )
    }.
  \end{equation}
\end{example}

\begin{example}[Multivariate polynomial dgc-algebras]
  \label{MultivariatePolynomialdgcAlgebras}
  Let $A \in \dgcAlgebras{\mathbb{R}}$ (Def. \ref{CategoryOfdgAlgebras}),
  $\mu^{(1)} \in A^{n_{1}+1}$, $d \mu^{(1)} = 0$,
  with corresponding polynomial dgc-algebra
  \eqref{PolynomialdgcAlgebraExpression}
  as in
  Example \ref{PolynomialAlgebraAndCellAttachment}.
  Then, for
     \vspace{-1mm}
  $$
    \mu^{(2)}
    \;\in\;
    A
    \big[
      \alpha^{(1)}_{n_1}
    \big]
    \big/
    \big(
      d\,\alpha^{(1)}_{n_1}\;=\; \mu^{(1)}
    \big)
    \,,
    \;\;\;\;\;
    d\,
    \mu^{(2)}
    \;=\; 0
  $$

   \vspace{-1mm}
\noindent
  another closed element of some homogeneous degree
  $n_2 + 1$, in the new algebra \eqref{PolynomialdgcAlgebraExpression}
  we may iterate the construction of
  Example \ref{PolynomialAlgebraAndCellAttachment} to obtain
  the bivariate polynomial dgc-algebra over $A$, to be denoted:
   \vspace{-1mm}
  $$
    A
    \!\!
   \scalebox{0.8}{$ \left[
      \!\!\!
      \begin{array}{c}
        \alpha^{(2)}_{n_2+1}
        \\
        \mathclap{\phantom{\vert^{\vert^{\vert^{\vert^{\vert}}}}}}
        \alpha^{(1)}_{n_1+1},
      \end{array}
      \!\!\!
    \right]
    $}
    \!\Big/\!
   \scalebox{0.8}{$ \left(
    \!\!
    \begin{aligned}
      d\,\alpha^{(2)}_{n_1}
      &
      = \mu^{(2)}_{n_2 + 1}
      \,,
      \\
      d\,\alpha^{(1)}_{n_1}
      &
      = \mu^{(1)}_{n_1 + 1}
    \end{aligned}
    \!\!
    \right)
    $}
    \;\;:=\;\;
    \left(
    A
    \!
    \big[
      \mu^{(1)}_{n_1 + 1}
    \big]
    \!\big/\!
    \big(
      d\,\alpha^{(1)}_{n_1+1}
      \;=\;
      \mu^{(1)}
    \big)
    \right)
    \!\!
    \big[
      \alpha^{(2)}
    \big]
    \!\big/\!
    \big(
      d\,\alpha^{(2)}_{n_2+1}
      \;=\;
      \mu^{(2)}
    \big)
    \,.
  $$

   \vspace{0mm}
\noindent
  Iterating further, we have multivariate polynomial dgc-algebras
  over $A$, to be denoted as follows:
    \vspace{-1mm}
  \begin{equation}
    \label{MultivariatePolynomialdgcAlgebra}
    A
    \!\!
\scalebox{0.8}{$    \left[
      \!\!\!
      \begin{array}{c}
        \alpha^{(k)}_{n_k+1}\,,
        \\
        \mathclap{\phantom{\vert^{\vert^{\vert^{\vert^{\vert}}}}}}
        \vdots
        \\
        \mathclap{\phantom{\vert^{\vert^{\vert^{\vert^{\vert}}}}}}
        \alpha^{(2)}_{n_2+1}
        \\
        \mathclap{\phantom{\vert^{\vert^{\vert^{\vert^{\vert}}}}}}
        \alpha^{(1)}_{n_1+1},
      \end{array}
      \!\!\!
    \right]
    $}
    \!\Big/\!
   \scalebox{0.8}{$ \left(
    \!
    \begin{aligned}
      d\,\alpha^{(k)}_{n_k},
      &
      = \mu^{(k)}
      \\
      & \vdots
      \\
      d\,\alpha^{(2)}_{n_1}
      &
      = \mu^{(2)}
      \,,
      \\
      d\,\alpha^{(1)}_{n_1}
      &
      = \mu^{(1)}
    \end{aligned}
    \!\!
    \right)
    $}
    \;\;
    \in
    \;
    \dgcAlgebras{\mathbb{R}}
  \end{equation}

  \vspace{-1mm}
  \noindent
  with
      \vspace{-1mm}
  $$
    \mu^{r}
    \;\in\;
    A
    \!\!
\scalebox{0.8}{$    \left[
      \!\!\!
      \begin{array}{c}
        \alpha^{(r-1)}_{n_{r-1}+1}\,,
        \\
        \mathclap{\phantom{\vert^{\vert^{\vert^{\vert^{\vert}}}}}}
        \vdots
        \\
        \mathclap{\phantom{\vert^{\vert^{\vert^{\vert^{\vert}}}}}}
        \alpha^{(1)}_{n_1+1},
      \end{array}
      \!\!\!
    \right]
    $}
    \,,
    \phantom{AAA}
    \mbox{
      for $1 \leq r \leq k$.
    }
  $$
  These multivariate polynomial algebras
  \eqref{MultivariatePolynomialdgcAlgebra} receive the canonical inclusion
  \eqref{InjectingAIntoPolynomialdgcAlgebraOverA} of $A$:
     \vspace{-2mm}
  \begin{equation}
    \label{InjectingAIntoMultivariatePolynomialdgcAlgebraOverA}
    \xymatrix{
      A
      \;
      \ar@{^{(}->}[rr]^-{ i_A }
      &&
      \;
    A
    \!\!
\scalebox{0.8}{$    \left[
      \!\!\!
      {\begin{array}{c}
        \alpha^{(k)}_{n_k+1}\,,
        \\
        \mathclap{\phantom{\vert^{\vert^{\vert^{\vert^{\vert}}}}}}
        \vdots
        \\
        \mathclap{\phantom{\vert^{\vert^{\vert^{\vert^{\vert}}}}}}
        \alpha^{(2)}_{n_2+1}
        \\
        \mathclap{\phantom{\vert^{\vert^{\vert^{\vert^{\vert}}}}}}
        \alpha^{(1)}_{n_1+1},
      \end{array}}
      \!\!\!
    \right]
    $}
    \!\Big/\!
   \scalebox{0.8}{$ \left(
    \!
    {\begin{aligned}
      d\,\alpha^{(k)}_{n_k}
      &
      = \mu^{(k)},
      \\
      & \vdots
      \\
      d\,\alpha^{(2)}_{n_1}
      &
      = \mu^{(2)}
      \,,
      \\
      d\,\alpha^{(1)}_{n_1}
      &
      = \mu^{(1)}
    \end{aligned}}
    \!\!
    \right)
    $}
   ,    }
  \end{equation}

    \vspace{-2mm}
\noindent
  these being the composites of the
  stage-wise inclusions \eqref{InjectingAIntoPolynomialdgcAlgebraOverA}.
\end{example}
\begin{defn}[Semifree dgc-Algebras/Sullivan models/FDAs]
  \label{SullivanModels}
  The multivariate polynomial dgc-algebras of
  Example \ref{MultivariatePolynomialdgcAlgebras}
  are sometimes
  called
 {\bf (i)} \emph{semi-free dgc-algebras} over $A$
  (since their underlying graded-commutative algebra
  \eqref{UnderlyingGradedCommutativeAlgebra} is free,
  as in Example \ref{FreeGradedCommutativeAlgebra}),
  but they are traditionally
  known {\bf (ii)}
  in rational homotopy theory as \emph{relative Sullivan models}
  (due to \cite{Sullivan77}, review in \cite[II]{FHT00}\cite{LM13}\cite{FelixHalperin17}),
  or, {\bf (iii)} in supergravity theory
  (following \cite{vanNieuwenhuizen82}\cite{DAuriaFre82}),
  as \emph{FDAs}\footnote{
    Beware that ``FDA'' in the supergravity literature
    is meant to be short-hand for
    ``free differential algebra'', which is misleading,
    because what is really meant are not free dgc-algebras as in
    Example \ref{FreeDifferentalGradedAlgebras}
    (in general) but just ``semi-free'' dcg-algebras,
    only whose underlying graded-commutative algebras
    \eqref{UnderlyingGradedCommutativeAlgebra}
    is required to be free (Example \ref{FreeGradedCommutativeAlgebra}).
  } \cite{CDF91},
  (for translation see \cite{FSS13}\cite{FSS16a}\cite{FSS16b}\cite{HSS18}\cite{BMSS19}\cite{FSS19a}).
  Here we write:
    \vspace{-1mm}
  \begin{equation}
    \label{CategoryOfSullivanModels}
    \xymatrix{
      \SullivanModelsConnected
      \;\;
      \ar@{^{(}->}[r]
      &
      \SullivanModels
      \;\;
      \ar@{^{(}->}[r]
      &
      \dgcAlgebras{\mathbb{R}}
    }
  \end{equation}
  for, from right to left,
  {\bf (a)} the full subcategory of connective
  dgc-algebras (Def. \ref{CategoryOfdgAlgebras})
  on those which are isomorphic to a
  multivariate polynomial dgc-algebra
  over $\mathbb{R}$,
  as in Example \ref{MultivariatePolynomialdgcAlgebras}
  (i.e., the ordering of the generators in \eqref{MultivariatePolynomialdgcAlgebra}
  is not part of the data of a Sullivan model, only the
  resulting dgc-algebra); and
  {\bf (b)} for the further full
  subcategory on those Sullivan model that are generated in
  positive degree  $\geq 1$.
\end{defn}

\begin{example}[Polynomial dgc-algebras as pushouts]
  \label{PolynomialdgcAlgebrasAsPushouts}
For $A \in \dgcAlgebras{\mathbb{R}}$ (Def. \ref{CategoryOfdgAlgebras})
the polynomial dgc-algebras over $A$ (Def. \ref{PolynomialAlgebraAndCellAttachment})
are pushouts in $\dgcAlgebras{\mathbb{R}}$ of the following form:
  \vspace{-2mm}
\begin{equation}
  \label{PolynomialdgcAlgebraAsPushout}
  \raisebox{40pt}{
  \xymatrix@R=1.3em{
    A
    \!
    \big[
      \alpha_{n}
    \big]
    \!\big/\!
    \big(
      d\,\alpha_n \,=\, \mu
    \big)
    \ar@{}[ddrr]|-{
      \mbox{
        \tiny
        (po)
      }
    }
    \ar@{<-^{)}}[dd]_-{
      i_A
    }
    \ar@{<-}[rr]^-{
      \scalebox{.7}{$
        {\begin{array}{ccc}
          \alpha_n &\!\!\!\!\mapsfrom\!\!\!\!& \alpha_n
          \\[-3pt]
          \mu &\!\!\!\mapsfrom\!\!\!& c_{n+1}
        \end{array}}
      $}
    }
    &&
    \mathbb{R}
    \!
    \left[
      \!\!\!
      {\begin{array}{c}
        \alpha_n\,,
        \\
        c_{n+1}
      \end{array}}
      \!\!\!
    \right]
    \!\big/\!
    \left(
      \!
      {\begin{aligned}
        d\;\, \alpha_n\;\;  & = c_{n+1}
        \\[-3pt]
        d\, c_{n+1}  & = 0
      \end{aligned}}
      \!\!
    \right)
    \ar@{<-^{)}}[dd]^-{
      \scalebox{.7}{$
        \arraycolsep=2.2pt
        {\begin{array}{c}
          c_{n+1}
          \\
          \mapsup
          \\
          c_{n+1}
        \end{array}}
      $}
    }
    \\
    \\
    A
    \ar@{<-}[rr]_-{
      \scalebox{.7}{$
      \begin{array}{ccc}
        \mu &\!\!\!\!\mapsfrom\!\!\!\!& c_{n+1}
      \end{array}
      $}
    }
    &&
    \mathbb{R}
    \big[
      c_{n+1}
    \big]
    \big/
    \big(
      d\,c_{n+1} \,=\, 0
    \big)
  }
  }
\end{equation}

  \vspace{-2mm}
\noindent Here on the right we have multivariate polynomial dgc-algebras
(Example \ref{MultivariatePolynomialdgcAlgebras}) over
$\mathbb{R}$ (Example \ref{InitialGradedAlgebra}) as shown.
The horizontal morphisms encode the choice of
$\mu \in A$ \eqref{CocycleToBeTrivialized} and the
left vertical morphism is the canonical inclusion \eqref{InjectingAIntoPolynomialdgcAlgebraOverA}.
\end{example}

\begin{example}[Chevalley-Eilenberg algebras of nilpotent Lie algebras]
 \label{CEAlgebraOfNilpotentLieAlgebra}
  Beware that not every Lie
  algebra $\mathfrak{g}$ has Chevalley-Eilenberg algebra
  (Example \ref{CEAlgebraOfLieAlgebra})
  which satisfies the stratification in the
  Definition \ref{MultivariatePolynomialdgcAlgebras}
  of multivariate polynomial dgc-algebras.

  \vspace{-1mm}
  \item {\bf (i)}  For instance, the Lie algebra $\mathfrak{su}(2)$
  has
  \vspace{-2mm}
  $$
    \mathrm{CE}
    \big(
      \mathfrak{su}(2)
    \big)
    \;=\;
    \mathbb{R}
    \big[
      \theta_1, \theta_2,
      \theta_3
    \big]
    \big/
    \Big(
    d\,\theta_i \,=\,
    \underset{j,k}{\sum}\epsilon_{i j k} \theta_j \wedge \theta_k
    \Big)
  $$

  \vspace{-2mm}
  \noindent  and no ordering of $\{1,2,3\}$ brings this into the
  iterative form
  required in \eqref{MultivariatePolynomialdgcAlgebra}.

  \vspace{-1mm}
    \item {\bf (ii)}  Instead, those Lie algebras whose CE-algebra is of
    the form
  \eqref{MultivariatePolynomialdgcAlgebra} are precisely the
  nilpotent Lie algebras.
\end{example}

In generalization of Example \ref{CEAlgebraOfNilpotentLieAlgebra},
we may say
(by \cite[Thm 2.3]{Berglund15} this matches \cite[Def. 4.2]{Getzler04}):
\begin{defn}[Nilpotent $L_\infty$-algebras]
  \label{NilpotentLInfinityAlgebras}
  An $L_\infty$-algebra \eqref{LInfinityAlgebra}
  is \emph{nilpotent} if its CE-algebra (Def. \ref{ChevalleyEilenbergAlgebraOfLInfinityAlgebra})
  is a multivariate polynomial dgc-algebra
  (Example \ref{MultivariatePolynomialdgcAlgebras}),
  hence is in the sub-category of $\SullivanModels$ \eqref{CategoryOfSullivanModels}:
  \begin{equation}
    \label{CategoryOfNilpotentLInfinityAlgebras}
    \xymatrix{
       L_\infty\Algebras
         ^{\scalebox{.5}{$\geq 0, {\color{blue}\mathrm{nil}}$}}
         _{\scalebox{.5}{$\mathbb{R}, \mathrm{fin}$}}
      \;
      \ar@{^{(}->}[rr]^-{ \mathrm{CE} }
      \ar@{^{(}->}[d]
      \ar@{}[drr]|-{ \mbox{\tiny(pb)} }
      &&
      \big(
        \SullivanModels
      \big)^{\mathrm{op}}
      \ar@{^{(}->}[d]
      \\
      \LInfinityAlgebras
      \;
      \ar@{^{(}->}[rr]^-{ \mathrm{CE} }
      &&
      \dgcAlgebrasOp{\mathbb{R}}
    }
  \end{equation}
  In fact, from \eqref{SemiFreedgcAlgebra} it is clear
  that every connected Sullivan model, hence with generators in
  degrees $\geq 1$, is the Chevalley-Eilenberg algebra of
  a unique nilpotent $L_\infty$-algebra, so that
  the defining inclusion at the top of \eqref{CategoryOfNilpotentLInfinityAlgebras} further restricts
  to an equivalence of homotopy categories:
 \vspace{-3mm}
  \begin{equation}
    \label{NilpotentLInfinityAlgebrasEquivalentToConnectedSullivanModels}
    \xymatrix{
      \LInfinityAlgebrasNil
      \ar[rr]^-{ \mathrm{CE} }_-{ \simeq }
      &&
      \SullivanModelsConnectedOp
      \,.
    }
  \end{equation}
\end{defn}

\noindent {\bf Homotopy theory of connective dgc-Algebras.}
We recall the homotopy theory of connective
differential graded-commutative algebras,
making free use of model category theory
\cite{Quillen67}; for a review see \cite{Hovey99}\cite[A.2]{Lurie09}
and \cref{ModelCategoryTheory}.

\begin{defn}[Homotopical structure on connective cochain complexes]
  \label{ProjectiveHomotopicalStructureOnConnectiveCochainComplexes}
  Consider the following sub-classes of morphisms
  in the category $\CochainComplexes$ (Def. \ref{CoochainComplexes}):

 \noindent  {\bf (i)} $\mathrm{W}$ -- \emph{weak equivalences} are the quasi-isomorphisms;

\noindent   {\bf (ii)} $\mathrm{Fib}$ -- \emph{fibrations} are the degreewise surjections;

 \noindent  {\bf (iii)} $\mathrm{Cof}$ -- \emph{cofibrations} are the
    injections in positive degrees.

  \noindent We call this the \emph{injective} homotopical structure on
   $\CochainComplexes$.
\end{defn}
\begin{prop}[Injective model structure on connective cochain complexes {\cite[p. 6]{Hess07}}]
  \label{ProjectiveModelStructureOnConnectiveCochainComplexes}
  Equipped with the injective homotopical structure of
  Def. \ref{ProjectiveHomotopicalStructureOnConnectiveCochainComplexes}
  the category of $\CochainComplexes$ (Def. \ref{CoochainComplexes})
  becomes a model category (Def. \ref{ModelCategories})
  which is right proper (Def. \ref{ProperModelCategories}).
  We denote this by:
  \vspace{-1mm}
  $$
    \CochainComplexesInj
    \;\in\;
    \mathrm{ModelCategories}
    \,.
  $$
\end{prop}
\begin{proof}
  The proof of the model structure itself is formally dual to the proof
  of the projective model structure on connective
  chain complexes \cite[II.4]{Quillen67}\cite[Thm. 1.5]{GoerssSchemmerhorn06};
  it is spelled out in
  \cite[Thm. 2.4.5]{Dungan10}.
  (Here we are using that for modules
  over a field of characteristic zero, as in our case,
  the condition that kernels of
  epimorphisms be injective is automatic.)
  A proof of right properness with respect to degreewise surjections
  is spelled out in
  \cite[Prop. 24]{Strickland20}.
\end{proof}
\begin{defn}[Homotopical structure on connective dgc-algebras {\cite[\S 4.2]{BousfieldGugenheim76}\cite[\S V.3.4]{GelfandManin96}}]
  \label{HomotopicalStructureOndgcAlgebras}
  Consider the following sub-classes of
  morphisms in the category of
  $\dgcAlgebras{\mathbb{R}}$ (Def. \ref{CategoryOfdgAlgebras}):

\noindent   {\bf (i)}  $\mathrm{W}$ -- \emph{weak equivalences} are the quasi-isomorphisms;

\noindent   {\bf (ii)} $\mathrm{Fib}$ -- \emph{fibrations} are the degreewise surjections;

  \noindent We call this the \emph{projective} homotopical structure on
  $\mathrm{dgcAlgebras}^{\geq 0}_{\mathbb{R}}$.
\end{defn}
\begin{prop}[Projective model structure connective on dgc-algebras]
  \label{ProjectiveModelStructureOndgcAlgebras}
  Equipped with the projective homotopical structure from Def. \ref{HomotopicalStructureOndgcAlgebras},
  the category of $\dgcAlgebras{\mathbb{R}}$ (Def. \ref{CategoryOfdgAlgebras}),
  becomes a model category (Def. \ref{ModelCategories})
  which is right proper (Def. \ref{ProperModelCategories}), in fact this
  is the case over any ground field $k$ of characteristic 0.
  \vspace{-2mm}
  \begin{equation}
    \label{ProjectiveModelStructureOnConnectivedgcAlgebras}
    \dgcAlgebrasProj{k}
    \;\in\;
    \mathrm{ModelCategories}
    \,.
  \end{equation}
\end{prop}
\begin{proof}
 The model structure itself is due to
 {\cite[\S 4.3]{BousfieldGugenheim76},
 the proof is spelled out in \cite[V.3.4]{GelfandManin96}}.
 Right properness follows from the right properness of
 the injective model structure on cochain complexes
 (Prop. \ref{ProjectiveModelStructureOnConnectiveCochainComplexes})
 since the free/forgetful adjunction
 \eqref{FreeForgetfulQuillenAdjunctionBetweenDgcAlgebrasAndCochainComplexes}
 implies that
 underlying pullbacks of dgc-algebras are pullbacks of the
 underlying cochain complexes.
\end{proof}
\begin{prop}[Quillen adjunction between dgc-algebras and cochain complexes]
  \label{QuillenAdjunctionBetweendgcAlgebrasAndCochainComplexes}
  The adjunction \eqref{dgCategoriesAndAdjunctions}
  between $\dgcAlgebras{\mathbb{R}}$ (Def. \ref{CategoryOfdgAlgebras})
  and $\CochainComplexes$ (Def. \ref{CoochainComplexes})
  is a
  Quillen adjunction (Def. \ref{QuillenAdjunction})
  with respect to the model category
  structures
  from Prop. \ref{ProjectiveModelStructureOnConnectiveCochainComplexes}
  and
  Prop. \ref{ProjectiveModelStructureOndgcAlgebras}:
  \vspace{-2mm}
  \begin{equation}
    \label{FreeForgetfulQuillenAdjunctionBetweenDgcAlgebrasAndCochainComplexes}
    \xymatrix{
      \dgcAlgebrasProj{\mathbb{R}}
      \;\;
      \ar@{<-}@<+8pt>[rr]^-{ \scalebox{0.6}{$\mathrm{Sym}  $}}
      \ar@<-8pt>[rr]^-{ \bot_{{}_{\mathrlap{\mathrm{Qu}}}} }_{\mathrm{underlying}}
      &&
      \;\;
      \CochainComplexesInj
    }
    \,.
  \end{equation}
\end{prop}
\begin{proof}
  It is immediate from
  Def. \ref{ProjectiveHomotopicalStructureOnConnectiveCochainComplexes} and
  Def. \ref{HomotopicalStructureOndgcAlgebras}
  that the forgetful right adjoint preserves
  the classes $\mathrm{W}$ and $\mathrm{Fib}$.
\end{proof}

\begin{remark}[All dgc-algebras are projectively fibrant]
  \label{AlldgcAlgebrasAreProjectivelyFibrant}
  Every object
  $A \in \dgcAlgebrasProj{\mathbb{R}}$ \eqref{ProjectiveModelStructureOnConnectivedgcAlgebras}
  is fibrant: By Example \ref{TheTerminalAlgebra} the terminal
  morphism is to the 0-algebra, and this is clearly surjective,
  hence is a fibration, by Def. \ref{HomotopicalStructureOndgcAlgebras}:
  $
    \xymatrix{
      A
      \ar[r]_-{ \in \, \mathrm{Fib} }
      &
      0
      \,.
    }
  $
\end{remark}

\medskip

\noindent {\bf Cofibrant dgc-algebras.}
With all dgc-algebras being fibrant (Rem. \ref{AlldgcAlgebrasAreProjectivelyFibrant}),
the crucial property is cofibrancy.

\begin{lemma}[Generating cofibrations]
  \label{GeneratingCofibrations}
  The following inclusions
  of multivariate polynomial dgc-algebras
  (Example \ref{MultivariatePolynomialdgcAlgebras})
  are cofibrations in $\dgcAlgebrasProj{\mathbb{R}}$ (Def. \ref{ProjectiveModelStructureOndgcAlgebras})
    \vspace{-2mm}
  \begin{equation}
    \label{GeneratingCofibrationsOfdgcAlgebras}
    \xymatrix{
      \mathbb{R}
      \!
      [
        c_{n+1}
      ]
      \big/
      (
        d\,c_{n+1} \,=\, 0
      )
      \;
      \ar@{^{(}->}[rrr]^-{
        c_{n+1}
        \mapsto
        c_{n-1}
      }_-{ \in\,\mathrm{Cof} }
      &&&
      \mathbb{R}
      \left[
        \!\!\!
        {\begin{array}{c}
          \alpha_n\,,
          \\
          c_{n+1}
        \end{array}}
        \!\!\!
      \right]
      \!\big/\!
      \left(
      \!
      {\begin{aligned}
        d\;\,\alpha_n\; & = c_{n+1}\,,
        \\[-3pt]
        d\, c_{n+1} & = 0
      \end{aligned}}
      \!
      \right)
    }
    \phantom{AAA}
    \mbox{
      for $n \in \mathbb{N}$.
    }
\end{equation}
\end{lemma}

\begin{proof}
  Consider the following morphisms of cochain complexes,
  for $n \in \mathbb{N}$:

  \vspace{-.4cm}
  \begin{equation}
    \xymatrix{
      \scalebox{.8}{$
      \left[
        \!\!
        {\begin{array}{c}
          \vdots
          \\
          0
          \\
          \uparrow\mathrlap{\!\scalebox{.7}{$d$}}
          \\
          0
          \\
          \uparrow\mathrlap{\!\scalebox{.7}{$d$}}
          \\
          \langle c_{n+1} \rangle
          \\
          \uparrow\mathrlap{\!\scalebox{.7}{$d$}}
          \\
          0
          \\
          \uparrow\mathrlap{\!\scalebox{.7}{$d$}}
          \\
          0
          \\
          \uparrow\mathrlap{\!\scalebox{.7}{$d$}}
          \\
          \vdots
          \\
          \uparrow\mathrlap{\!\scalebox{.7}{$d$}}
          \\
          0
        \end{array}}
        \!\!
      \right]
      $}
      \;
      \ar@{^{(}->}[rr]^-{ i_n }
      &&
      \;
      \scalebox{.8}{$
      \left[
        \!\!
        {\begin{array}{c}
          \vdots
          \\
          0
          \\
          \uparrow\mathrlap{\!\scalebox{.7}{$d$}}
          \\
          0
          \\
          \uparrow\mathrlap{\!\scalebox{.7}{$d$}}
          \\
          \langle c_{n+1} \rangle
          \\
          \uparrow\mathrlap{\!\scalebox{.7}{$d$}}
          \\
          %RationalizationViaPLDeRhamAdjunction
          \langle \alpha_n \rangle
          \\
          \uparrow\mathrlap{\!\scalebox{.7}{$d$}}
          \\
          0
          \\
          \uparrow\mathrlap{\!\scalebox{.7}{$d$}}
          \\
          \vdots
          \\
          \uparrow\mathrlap{\!\scalebox{.7}{$d$}}
          \\
          0
        \end{array}}
        \!\!
      \right]
      $}
      \mathrlap{
        \;\;\;\;\;\;
        \mbox{
          with
          $d \alpha_n \,=\, c_{n+1}$.
        }
      }
    }
  \end{equation}
  \vspace{-.3cm}

  \noindent
  Since these are injections,
  they are cofibrations in
  $\CochainComplexesInj$
  (Prop. \ref{ProjectiveModelStructureOnConnectiveCochainComplexes}),
  by Def. \ref{ProjectiveHomotopicalStructureOnConnectiveCochainComplexes}.
  Thus also their images
  under $\mathrm{Sym}$ (Def. \ref{FreeDifferentalGradedAlgebras})
  are cofibrations in
  $\dgcAlgebrasProj{\mathbb{R}}$ (Prop. \ref{ProjectiveModelStructureOndgcAlgebras})
  because $\mathrm{Sym}$ is a left Quillen functor,
  by Prop. \ref{QuillenAdjunctionBetweendgcAlgebrasAndCochainComplexes}.
  But $\mathrm{Sym}(i_n)$ manifestly equals \eqref{GeneratingCofibrationsOfdgcAlgebras},
  and so the claim follows.
\end{proof}

\begin{prop}[Relative Sullivan algebras are cofibrations]
  \label{RelativeSullivanModelsAreCofibrations}
  For a multivariate polynomial dgc-algebra
  from Example \ref{MultivariatePolynomialdgcAlgebras}, the canonical
  inclusion \eqref{InjectingAIntoMultivariatePolynomialdgcAlgebraOverA}
  of the base algebra
  is a cofibration in $\dgcAlgebrasProj{\mathbb{R}}$ (Prop. \ref{ProjectiveModelStructureOndgcAlgebras}):

  \vspace{-2mm}
  \begin{equation}
    \label{InjectingAIntoMultivariatePolynomialdgcAlgebraOverA}
    \xymatrix{
      A
      \;
      \ar@{^{(}->}[rr]_-{ i_A }^-{ \in \, \mathrm{Cof} }
      &&
      \;
    A
    \!\!
   \scalebox{0.8}{$ \left[
      \!\!\!
      {\begin{array}{c}
        \alpha^{(k)}_{n_k+1}\,,
        \\
        \mathclap{\phantom{\vert^{\vert^{\vert^{\vert^{\vert}}}}}}
        \vdots
        \\
        \mathclap{\phantom{\vert^{\vert^{\vert^{\vert^{\vert}}}}}}
        \alpha^{(1)}_{n_1+1},
      \end{array}}
      \!\!\!
    \right]
    $}
    \!\Big/\!
    \scalebox{0.8}{$\left(
        {\begin{aligned}
      d\,\alpha^{(k)}_{n_k}
      &
      = \mu^{(k)},
      \\
      & \vdots
      \\
      d\,\alpha^{(1)}_{n_1}
      &
      = \mu^{(1)}
    \end{aligned}}
    \!\!
    \right)
    $}    .
    }
  \end{equation}
  \vspace{-.3cm}

\noindent  In particular, since $\mathbb{R} \in \dgcAlgebras{\mathbb{R}}$
  is the initial object (Example \ref{InitialGradedAlgebra}),
  all multivariate polynomial dgc-algebras over $\mathbb{R}$
  (the Sullivan models, Def. \ref{SullivanModels})
  are cofibrant objects in $\dgcAlgebrasProj{\mathbb{R}}$.
\end{prop}
\begin{proof}
  By Lemma \ref{GeneratingCofibrations},
  the right vertical morphisms in the pushout
  diagram \eqref{PolynomialdgcAlgebraAsPushout} are cofibrations.
  Since the class of cofibrations is preserved under pushout,
  so are hence the left vertical morphisms in
  \eqref{PolynomialdgcAlgebraAsPushout}, which are the
  base algebra inclusions \eqref{InjectingAIntoPolynomialdgcAlgebraOverA}
  of polynomial dgc-algebras.
  The base algebra inclusions into general
  multivariate polynomial dgc-algebras
  are composites of these, and since the class of cofibrations
  is presered under composition, the claim follows.
\end{proof}

\begin{lemma}[Pushout along relative Sullivan algebras preserves
quasi-isomorphisms {\cite[Prop. 6.7 (ii), Lemma 14.2]{FHT00}}]
  \label{PushoutAlongRelativeSullivanModelPreservesQuasiIsomorphism}
  The operation of pushout \eqref{coCartesianSquare}
  along
  the canonical
  inclusion \eqref{InjectingAIntoMultivariatePolynomialdgcAlgebraOverA}
  of a base dgc-algebra
  into a multivariate polynomial dgc-algebra
  (Example \ref{MultivariatePolynomialdgcAlgebras})
  preserves quasi-isomorphisms. In fact, it sends
  quasi-isomorphism between base algebras to
  quasi-isomrophisms of multivariate polynomial dgc-algebras:

  \vspace{-2mm}
  \begin{equation}
    \label{PushoutAlongRelativeSullivanAlgbra}
    \raisebox{46pt}{
    \xymatrix@C=3em{
      A
      \;
      \ar@{^{(}->}[rr]_-{ i_A }^-{ \in \, \mathrm{Cof} }
      \ar[d]_-{ f }^-{ \in \; \mathrm{W} }
      \ar@{}[drr]|>>>>>>>{
        \mbox{
          \hspace{-2cm}
          \tiny \rm(po)
        }
      }
      &&
      \;
    A
    \!\!
  \scalebox{0.8}{$   \left[
      \!\!\!
      {\begin{array}{c}
        \alpha^{(k)}_{n_k+1}\,,
        \\
        \mathclap{\phantom{\vert^{\vert^{\vert^{\vert^{\vert}}}}}}
        \vdots
        \\
        \mathclap{\phantom{\vert^{\vert^{\vert^{\vert^{\vert}}}}}}
        \alpha^{(1)}_{n_1+1},
      \end{array}}
      \!\!\!
    \right]
    $}
    \!\Big/\!
    \scalebox{0.8}{$ \left(
        {\begin{aligned}
      d\,\alpha^{(k)}_{n_k}
      &
      = \mu^{(k)},
      \\
      & \vdots
      \\
      d\,\alpha^{(1)}_{n_1}
      &
      = \mu^{(1)}
    \end{aligned}}
    \!
    \right)
    $}
    \ar@<-34pt>[d]^-{
      (i_A)_\ast f
    }
    \\
      A'
      \;
      \ar@{^{(}->}[rr]_-{ i_A }^-{ \in \, \mathrm{Cof} }
      &&
      \;
    A'
    \!\!
     \scalebox{0.8}{$\left[
      \!\!\!
      {\begin{array}{c}
        \alpha^{(k)}_{n_k+1}\,,
        \\
        \mathclap{\phantom{\vert^{\vert^{\vert^{\vert^{\vert}}}}}}
        \vdots
        \\
        \mathclap{\phantom{\vert^{\vert^{\vert^{\vert^{\vert}}}}}}
        \alpha^{(1)}_{n_1+1},
      \end{array}}
      \!\!\!
    \right]
    $}
    \!\Big/\!
   \scalebox{0.8}{$  \left(
        {\begin{aligned}
      d\,\alpha^{(k)}_{n_k}
      &
      = \mu^{(k)},
      \\
      & \vdots
      \\
      d\,\alpha^{(1)}_{n_1}
      &
      = \mu^{(1)}
    \end{aligned}}
    \!
    \right)
$}
    }
    }
    \phantom{AAAA}
    \Rightarrow
    \phantom{AAAA}
    (i_a)_\ast f
    \;\in\;
    \mathrm{W}
    \,.
  \end{equation}
  \vspace{-.5cm}
\end{lemma}

\begin{lemma}[Weak equivalences of nilpotent $L_\infty$-algebras {\cite[Prop. 14.13]{FHT00}}]
  \label{WeakEquivalencesOfNilpotentLInfinityAlgebras}
  A morphism between Chevalley-Eilenberg algebras (Def. \ref{ChevalleyEilenbergAlgebraOfLInfinityAlgebra})
  of nilpotent $L_\infty$-algebras (Def. \ref{NilpotentLInfinityAlgebras}),
  is a quasi-isomorphism of dgc-algebras
  (hence a weak equivalence according to Def. \ref{HomotopicalStructureOndgcAlgebras})
  precisely if the corresponding
  morphism \eqref{FullSubcategoryInclusionOfLieAlgebrasIntodgcAlgebras}
  of $L_\infty$-algebras is a quasi-isomorphism between the
  chain complexes given by the unary bracket operation $\partial := [-]$
  \eqref{BracketsOnLInfinityAlgebra}:

  \vspace{-3mm}
  $$
    \xymatrix{
      \mathrm{CE}(\mathfrak{g})
      \ar@{<-}[rr]^-{ \mathrm{CE}(\phi) }_-{ \in \; \mathrm{W} }
      &&
      \mathrm{CE}(\mathfrak{h})
    }
    \;\;\;\;\;
    \Leftrightarrow
    \;\;\;\;\;
    \xymatrix{
      \big(
        \mathfrak{g}, [-]_{\mathfrak{g}}
      \big)
      \ar[rr]^-{ \phi }_{ \in \; \mathrm{W} }
      &&
      \big(
        \mathfrak{h}, [-]_{\mathfrak{h}}
      \big).
    }
  $$
\end{lemma}

\begin{remark}[Homotopy theory of nilpotent $L_\infty$-algebras inside all
$L_\infty$-algebras]
  \label{HomotopyTheoryOfLInfinityAlgebras}
  $\,$

  \noindent
  {\bf (i)}
  Prop. \ref{RelativeSullivanModelsAreCofibrations},
  with Remark \ref{AlldgcAlgebrasAreProjectivelyFibrant}
  and
  Def. \ref{ChevalleyEilenbergAlgebraOfLInfinityAlgebra},
  allows to identify
  the homotopy category of
  finite-type nilpotent connective  $L_\infty$-algebras
  (Def. \ref{NilpotentLInfinityAlgebras}), with a
  full subcategory of the homotopy category
  (Def. \ref{HomotopyCategory}) of the opposite
  (Example \ref{OppositeModelStructure}) of dgc-algebras
  (Prop. \ref{ProjectiveModelStructureOndgcAlgebras}):

  \vspace{-4mm}
  \begin{equation}
    \label{HomotopyCategoryOfNilpotentLInfinityAlgebras}
    \xymatrix{
      \mathrm{Ho}
      \big(
        \LInfinityAlgebrasNil
      \big)
      \;
      \ar@{^{(}->}[r]^-{ \mathrm{CE} }
      &
      \mathrm{Ho}
      \big(
        \dgcAlgebrasOpProj{\mathbb{R}}
      \big)
      \,.
    }
  \end{equation}
  \vspace{-.4cm}

  \noindent
  {\bf (ii)} There is also the homotopy theory of
  more general
  $L_\infty$-algebras \cite{Hinich01}\cite{Pridham10}\cite{Vallette14}\cite{Rogers20},
  whose weak equivalences are the quasi-isomorphisms on
  chain complexes formed by the unary bracket $[-]$ \eqref{BracketsOnLInfinityAlgebra}.
  Lemma \ref{WeakEquivalencesOfNilpotentLInfinityAlgebras}
  says that the homotopy theory
  \eqref{HomotopyCategoryOfNilpotentLInfinityAlgebras}
  of finite-type, nilpotent connective $L_\infty$-algebras
  that we are concerned with here
  is fully faithfully embedded into this more general
  $L_\infty$ homotopy theory:

  \vspace{-.2cm}
  $$
    \xymatrix{
      \mathrm{Ho}
      \big(
        \LInfinityAlgebrasNil
      \big)
      \;
      \ar@{^{(}->}[r]
      &
      \mathrm{Ho}
      \big(
        \LInfinityAlgebrasGeneral
      \big)
      \,.
    }
  $$
  \vspace{-.4cm}

\end{remark}

\medskip

\noindent {\bf Minimal Sullivan models}

\begin{defn}[Minimal Sullivan models {\cite[Def. 7.2]{BousfieldGugenheim76}\cite[Def. 1.10]{Hess07}}]
  \label{MinimalSullivanModels}
  A connected (relative) Sullivan model dgc-algebra
  $A \in \SullivanModelsConnected$ (Def. \ref{SullivanModels})
  is called \emph{minimal} if
  it is given by a multivariate polynomial dgc-algebra
  as in \eqref{MultivariatePolynomialdgcAlgebra}
  the degrees $n_i$ of whose generators $\alpha_{n_i}^{(i)}$
  are monotonically increasing
    \vspace{-2mm}
  $$
    i < j
    \;\;\;
    \Rightarrow
    \;\;\;
    n_j \leq n_j\;.
  $$
\end{defn}
\begin{example}[Minimal models of simply connected dgc-algebras {\cite[Prop. 7.4]{BousfieldGugenheim76}}]
  \label{MinimalModelOFSimplyConnecteddgcAlgebras}
  If
  $A \in \SullivanModelsConnected$ (Def. \ref{SullivanModels})
  is trivial in degree 1, then
  it is minimal (Def. \ref{MinimalSullivanModels})
  precisely if the unary bracket $[-]$
  \eqref{LInfinityBracketsEncodedInCEDifferential} of the
  corresponding $L_\infty$-algebra \eqref{NilpotentLInfinityAlgebrasEquivalentToConnectedSullivanModels}
  vanishes:
    \vspace{-1mm}
  $$
    A^1 = 0
    \;\;\;\;
    \Rightarrow
    \;\;\;\;
    \big(
      \mbox{$A$ is minimal}
      \;\;\;
      \Leftrightarrow
      \;\;\;
      [-] = 0
    \big).
  $$
\end{example}

\begin{prop}[Existence of minimal Sullivan models {\cite[Prop. 7.7, 7.8]{BousfieldGugenheim76}\cite[Thm. 14.12]{FHT00}}]
  \label{ExistenceOfMinimalSullivanModels}
  $\,$

  \noindent
  If $A \in \dgcAlgebras{\mathbb{R}}$ is cohomologically connected, in that
  $H^0(A) = \mathbb{R}$, then:

 \noindent {\bf (i)} There exists a minimal Sullivan model
   $A_{\mathrm{min}}$
    (Def. \ref{MinimalSullivanModels})
    with weak equivalence in $\dgcAlgebrasProj{\mathbb{R}}$ \eqref{ProjectiveModelStructureOnConnectivedgcAlgebras} to $A$
    \vspace{-3mm}
    \begin{equation}
      \label{MinimalSullivanModel}
      \xymatrix{
        A_{\mathrm{min}}
        \ar[rr]^-{ p_A^{\mathrm{min}} \in\, \mathrm{W} }
        &&
        A
        \,.
      }
    \end{equation}

   \vspace{-1mm}
   \noindent
   {\bf (ii)} This $A_{\mathrm{min}}$ is unique up to
   isomorphisms of $\dgcAlgebras{\mathbb{R}}$ compatible with
   the weak equivalences in \eqref{MinimalSullivanModel}:
   Any two
   $p_A^{\mathrm{min}}, p_A^{\mathrm{min}'}$ in \eqref{MinimalSullivanModel},
   make a commuting diagram of this form (\cite[Thm. 14.12]{FHT00}):
      \vspace{-2mm}
   \begin{equation}
     \label{FactorizationOfMinimalModelsThroughEachOther}
     \begin{tikzcd}[row sep=-2pt, column sep=30pt]
       A_{\mathrm{min}}
       \ar[
         dr,
         "{
           p_A^{\mathrm{min}}
         }"{above, xshift=3pt}
       ]
       \ar[
         dd,
         "\simeq"{left}
       ]
       \\
       &
       A \;.
       \\
       A_{\mathrm{min}'}
       \ar[
         ur,
         "{
           p_A^{\mathrm{min}'}
         }"{right, xshift=-9pt, yshift=-7pt}
       ]
     \end{tikzcd}
   \end{equation}
\end{prop}

\vspace{-2mm}
More generally:

\vspace{-3mm}
\begin{prop}[Existence of minimal relative Sullivan models {\cite[Thm. 14.12]{FHT00}}]
  \label{ExistenceOfRelativeMinimalSullianModels}
  Let $\!\xymatrix@C=12pt{B \ar[r]^-{\phi} & A}\!$ be a morphism in $\dgcAlgebras{\mathbb{R}}$
  (Def. \ref{CategoryOfdgAlgebras}) such that

  {\bf (a)} $A$ and $B$ are cohomologically connected,
  in that $H^0(A) = \mathbb{R}$ and $H^0(B) = \mathbb{R}$,

  {\bf (b)} $H^1(\phi) : H^1(B) \longrightarrow H^1(A)$ is
  an injection.

  \noindent Then:

  \noindent
  {\bf (i)} There exists a minimal relative Sullivan model
  $B \longhookrightarrow A_{\mathrm{min}_B}$
  (Def. \ref{MinimalSullivanModels})
  equipped with a weak equivalence to $\phi$
  in $\dgcAlgebrasProj{\mathbb{R}}$ (Def. \ref{ProjectiveModelStructureOnConnectivedgcAlgebras}):
  \vspace{-3mm}
  \begin{equation}
  \label{RelativeMinimalSullivanModel}
  \hspace{2cm}
    \xymatrix@R=1pt@C=3em{
      A_{\mathrm{min}_B}
      \ar[rr]^-{ \in\; \mathrm{W} }
      &&
      A
      \\
      &
      \,
      B
      \ar[ur]_-{\phi}
      \ar@{^{(}->}[ul]
    }
  \end{equation}

  \vspace{-2mm}
  \noindent
  {\bf (ii)}
  This $A_{\mathrm{min}_B}$ is unique up to isomorphism in the coslice category
  $\big(\dgcAlgebras{\mathbb{R}}\big)^{B/}$
  compatible with the weak equivalence in \eqref{RelativeMinimalSullivanModel},
  in cosliced generalization of \eqref{FactorizationOfMinimalModelsThroughEachOther}.
\end{prop}

%\newpage

%%%%%%%%%%%%%%%%%%%%%%%%%%%%%%%%%%%%%%%%%
\subsection{$\mathbb{R}$-Rational homotopy theory}
\label{RationalHomotopyTheory}
%%%%%%%%%%%%%%%%%%%%%%%%%%%%%%%%%%%%%%%%%

We recall fundamental facts of dg-algebraic rational homotopy theory
\cite{Sullivan77}\cite{BousfieldGugenheim76}\cite{GriffithMorgan13}
(review in \cite{FHT00}\cite{Hess07}\cite{FOT08} \cite{FelixHalperin17}),
\footnote{
  One may naturally understand rational homotopy theory
  also within the theory of derived algebraic $\infty$-stacks
  \cite{Toen00ChampsAffine}\cite[\S 1]{Lurie11}.
  The real character map in \cref{ChernCharacterInNonabelianCohomology}
  instead expresses rational homotopy theory
  within smooth (differential-geometric) $\infty$-stacks
  in non-abelian generalization of the way it appears in
  differential cohomology theory, see \cref{NonabelianDifferentialCohomology}.
}
with emphasis on its incarnation over the {\it real} numbers
(Rem. \ref{RationalHomotopyTheoryOverTheRealNumbers}) and
streamlined towards the application to non-abelian de Rham theory
below in \cref{NonAbelianDeRhamTheory} and thus to the non-abelian character map in
\cref{ChernCharacterInNonabelianCohomology}.
For the usual technical reasons (Rem. \ref{AssumptionsOnConnctedNilpotentRFiniteHomotopyTypes}),
we focus on the following
class of homotopy types (with little to no restriction
in practice):
\begin{defn}[Connected nilpotent spaces of finite rational type {\cite[9.2]{BousfieldGugenheim76}}]
\label{NilpotentConnectedSpacesOfFiniteRationalType}
Write
\vspace{-2mm}
$$
  \xymatrix{
    \NilpotentConnectedQFiniteHomotopyTypes
    \;
    \ar@{^{(}->}[r]
    &
    \;
    \HomotopyTypes
  }
$$

\vspace{-2mm}
\noindent
for the full subcategory of homotopy types of topological spaces $X$
\eqref{ClassicalHomotopyCategory}
on those which are:

  {\bf (i)} {\it connected}: $\pi_0(X) \simeq \ast$;

  {\bf (ii)} {\it nilpotent}: $\pi_1(X) \,\in\, \mathrm{NilpotentGroups}$,
    and $\pi_{n \geq 2}(X)$ are nilpotent $\pi_1(X)$-modules (e.g. \cite{Hilton82});

  {\bf (iii)} {\it finite rational type}:
  $\mathrm{dim}_{\mathbb{Q}}\big( H^n(X;\, \mathbb{Q}) \big) < \infty\,,\;$
  for all $n \in \mathbb{N}$.
\end{defn}

\begin{remark}[Technical assumptions]
  \label{AssumptionsOnConnctedNilpotentRFiniteHomotopyTypes}
  The connectedness assumption in Def. \ref{NilpotentConnectedSpacesOfFiniteRationalType}
  is a pure convenience; for non-connected spaces all
  of the following applies just by iterating over connected components.
  On the other hand,
  the  nilpotency and $\mathbb{R}$-finiteness condition in
  Def. \ref{NilpotentConnectedSpacesOfFiniteRationalType}
  are strictly necessary for the plain dg-algebraic formulation of
  rational homotopy theory
  (due to \cite{BousfieldGugenheim76}\cite{Sullivan77})
  to satisfy the fundamental theorem
  (Theorem \ref{FundamentalTheoremOfdgcAlgebraicRationalHomotopyTheory} below).
  The generalizations required to drop these assumptions are
  known, but considerably more unwieldy:

  \noindent
  {\bf (i)}
  To drop the nilpotency assumption, all dgc-algebra models need to be
  equipped with the action of the fundamental group (see \cite{FHT15}).

  \noindent
  {\bf (ii)}
  To drop the finite-type assumption
  one needs dgc-coalgebras in place of dgc-algebras, as in
  the original \cite{Quillen69}.

\noindent  Therefore, we expect that the construction of the
  (twisted) non-abelian character map, below in sections \cref{ChernCharacterInNonabelianCohomology}
  and
  \cref{TwistedNonabelianCharacterMap},
  works also without imposing these technical assumptions,
  but a discussion in that generality is beyond the
  scope of the present article.
\end{remark}

\begin{example}[Examples of nilpotent spaces {\cite[\S 3]{Hilton82}\cite[\S 3.1]{MayPonto12}}]
  \label{ExamplesOfNilpotentSpaces}
  Such examples
  % of nilpotent spaces
  (Def. \ref{NilpotentConnectedSpacesOfFiniteRationalType})
  include:
  \begin{itemize}

    \vspace{-.3cm}
    \item[{\bf (i)}]
    every simply connected space $X$, $\pi_1(X) = 1$;

    \vspace{-2.5mm}
    \item[{\bf (ii)}]
    every simple space $X$, i.e. with abelian fundamental group acting
    trivially, such as tori;

    \vspace{-2.5mm}
    \item[{\bf (iii)}]
    hence every connected H-space;

    \vspace{-2.5mm}
    \item[{\bf (iv)}]
    hence every loop space $X \simeq \Omega Y$,
    and hence every $\infty$-group (Prop. \ref{ConnectedHomotopyTypesAreHigherNonAbelianClassifyingSpaces});

    \vspace{-2.5mm}
    \item[{\bf (v)}]
    hence every infinite-loop space,
    i.e., every component space $E_n$ of a spectrum $E$ \eqref{Spectrum};

    \vspace{-2.5mm}
    \item[{\bf (vi)}]
    the classifying spaces $B G$ \eqref{ClassifyingSpace}
    of nilpotent Lie groups $G$;

    \vspace{-2.5mm}
    \item[{\bf (vii)}]
    the mapping spaces $\mathrm{Maps}(X,A)$
    out of manifolds $X$ into nilpotent spaces $A$.

  \end{itemize}
\end{example}

Rational homotopy theory
is concerned with understanding the
following notion:
\begin{defn}[Rationalization
{\cite[p. 133]{BousfieldKan72}\cite[\S 11.1]{BousfieldGugenheim76}\cite[\S 1.4, \S 1.7]{Hess07}}]
  \label{Rationalization}
  $\,$

\noindent {\bf (i)}  A connected nilpotent homotopy type
  $X \in \mathrm{Ho}\big( \TopologicalSpaces_{\mathrm{Qu}}\big)_{\geq 1, \mathrm{nil}}$
  (Def. \ref{NilpotentConnectedSpacesOfFiniteRationalType})
  is called \emph{rational}
  if the following equivalent conditions hold
  \cite[\S V 3.3]{BousfieldKan72}\cite[\S 9.2]{BousfieldGugenheim76}:

  \vspace{-.2cm}
  \begin{itemize}

  \vspace{-.1cm}
  \item
  the higher homotopy groups $\pi_{\bullet \geq 2}(X)$ have the structure of
  $\mathbb{Q}$-vector spaces, and
  the fundamental group $\pi_1(X)$ is
  {\it uniquely divisible} in that each element
  $g$
  has a unique $n$th root $x$, i.e. with $x^n = g$, for all $n \in \mathbb{N}_+$;

  \vspace{-.2cm}
  \item
  the integral homology groups $H_{\bullet \geq 1}(X;\, \mathbb{Z})$
  all carry the structure of $\mathbb{Q}$-vector spaces;

  \end{itemize}
  \vspace{-.2cm}

\noindent {\bf (ii)}
A \emph{rationalization} of $X$ is a map
\vspace{-2mm}
  \begin{equation}
    \label{RationalizationUnit}
    \xymatrix{
      X
      \ar[r]^-{ \eta_X^{\mathbb{Q}} }
      &
      L_{\mathbb{Q}}(X)
      \;\;\;
      \in
      \mathrm{Ho}
      \big(
        \TopologicalSpaces_{\mathrm{Qu}}
      \big)_{\geq 1, \mathrm{nil}}    }
  \end{equation}

  \vspace{-4mm}
  \noindent
  such that:

  {\bf (a)} $L_{\mathbb{Q}}$ is rational in the above sense;

  {\bf (b)} the map $\eta_X^{\mathbb{Q}}$ induces an isomorphism on rational cohomology groups:
  \vspace{-3mm}
    $$
      \xymatrix{
        H^\bullet
        \big(
          L_{\mathbb{Q}} X;\, \mathbb{Q}
        \big)
        \ar[rr]^-{
          H^\bullet
          (
            \eta_X^{\mathbb{Q}}
            ;
            \,
            \mathbb{Q}
          )
        }_-{\simeq}
        &&
        H^\bullet (X;\, \mathbb{Q})
        \,.
      }
    $$
\end{defn}
Rationalization exists essentially uniquely, and defines
a reflective subcategory inclusion

\vspace{-3mm}
\begin{equation}
  \label{RationalizationReflection}
  \xymatrix{
    \overset{
      \raisebox{3pt}{
        \tiny
        \color{darkblue}
        \bf
        \def\arraystretch{.9}
        \begin{tabular}{c}
          connected, nilpotent,
          \\
          {\color{orangeii}rational} homotopy types
        \end{tabular}
      }
    }{
    \mathrm{Ho}
    \big(
      \TopologicalSpaces_{\mathrm{Qu}}
    \big)_{\geq 1, \mathrm{nil}}^{\mathbb{Q}}
    }
    \;\;
    \ar@{<-}@<+7pt>[rr]^-{ L_{\mathbb{Q}} }
    \ar@{^{(}->}@<-7pt>[rr]^-{\bot}
    &&
    \;
    \overset{
      \raisebox{3pt}{
        \tiny
        \color{darkblue}
        \bf
        \def\arraystretch{.9}
        \begin{tabular}{c}
          connected, nilpotent
          \\
          homotopy types
        \end{tabular}
      }
    }{
      \mathrm{Ho}
      \big(
        \TopologicalSpaces_{\mathrm{Qu}}
      \big)_{\geq 1, \mathrm{nil}}
    }
  }
\end{equation}
\vspace{-1mm}
\noindent whose adjunction unit \eqref{AdjunctionUnit} is \eqref{RationalizationUnit}.

\medskip

\noindent {\bf PL de Rham theory.} At the heart
of dg-algebraic rational homotopy theory is the
observation that a variant of the de Rham dg-algebra
of a smooth manifold (Example \ref{SmoothdeRhamComplex})
applies to general topological spaces:
the \emph{PL de Rham complex}\footnote{
  The terminology ``PL'' or ``P.L.'' for this construction
  seems to have been silently introduced in
  \cite{BousfieldGugenheim76}, as shorthand for
  ``piecewise linear'', and has become widely adopted
  (e.g. \cite[\S 9]{GriffithMorgan13}).
  Our
  subscript ``$\mathrm{P}k\mathrm{L}$'' is for ``piecewise $k$-linear'', in this sense.
  But beware that this refers to the
  piecewise-linear structure that a choice of triangulation
  (Ex. \ref{HomotopyTypesOfManifoldsViaTriangulations})
  induces on a topological space; while the actual differential
  forms in the PL de Rham complex are piecewise \emph{polynomial}
  with respect to this piecewise linear structure.
} (Def. \ref{PLdeRhamComplexes}).
This satisfies an appropriate
{\it PL de Rham theorem} (Prop. \ref{PLdeRhamTheorem})
and makes dg-algebras of PL differential forms
detect rational homotopy type (Prop. \ref{FundamentalTheoremOfdgcAlgebraicRationalHomotopyTheory}).
At the same time, over a smooth manifold the PL de Rham complex
is suitably equivalent to the smooth de Rham complex
(Lemma \ref{PLdeRhamComplexOnSmoothManifoldsEquivalentToSmoothDeRhamComplex}).

\begin{defn}[PL de Rham complex and PL de Rham cohomology {\cite[pp. 1-7]{BousfieldGugenheim76}\cite[\S 9.1]{GriffithMorgan13}}]
  \label{PLdeRhamComplexes}
Let $k$ be a field of characteristic zero.

\noindent
{\bf (i)} {\it simplicial dgc-algebra of $k$-polynomial differential
forms on the standard simplices}
(\cite[p. 297]{Sullivan77}\cite[p. 1]{BousfieldGugenheim76}\cite[p. 83]{GriffithMorgan13})
is:
\vspace{-2mm}
\begin{equation}
  \label{pdR}
  \begin{tikzcd}[row sep=-1pt]
    \mathllap{
      \Omega^\bullet_{k\mathrm{pdR}}
      \big(
        \Delta^{(-)}
      \big)
      \;\;:\;\;\;
    }
    \Delta^{\mathrm{op}}
    \ar[rr]
    &&
    \dgcAlgebras{k}
    \\
    {[n]}
    \ar[
      d,
      "f"
    ]
    &\longmapsto&
    k
    \big[
      t^{(0)}_0, \cdots, t^{(n)}_0,
      \,
      \theta^{(0)}_1, \cdots, \theta^{(n)}_1
    \big]
    \mathrlap{
    \Big/
    \!
 \scalebox{0.8}{$   \left(
      \!\!\!\!\!
      \def\arraystretch{1.1}
      \begin{array}{l}
      \;\;\sum_i \; t^{(i)}_0  =  1
      \!,\;
      \\
      \forall_{\!i}\,\;
      d\, t^{(i)}_0 =  \theta^{(i)}_1
      \end{array}
      \!\!\!\!\!\!
    \right)
    $}
    }
    \\[40pt]
    {[m]}
    &\longmapsto&
    k
    \big[
      t^{(0)}_0, \cdots, t^{(m)}_0,
      \,
      \theta^{(0)}_1, \cdots, \theta^{(m)}_1
    \big]
    \mathrlap{
    \Big/
    \!
 \scalebox{0.8}{$      \left(
      \!\!\!\!\!
      \def\arraystretch{1.1}
      \begin{array}{l}
      \;\;\sum_j \; t^{(j)}_0  =  1
      \!,\;
      \\
      \forall_{\!j}\,\;
      d\, t^{(j)}_0 =  \theta^{(j)}_1
      \end{array}
      \!\!\!\!\!\!
    \right)
    $}
    }
    \ar[
      u,
      |->,
      "{  \scalebox{0.8}{$
        \begin{array}{c}
          \underset{
            \mathclap{
              i \, \in \, f^{-1}(\{j\})
            }
          }
          {\sum}\; t^{(i)}_0
          \\
          \mapsup
          \\
          t^{(j)}_0
        \end{array}
        $}
      }"
    ]
  \end{tikzcd}
\end{equation}

\vspace{-1mm}
\noindent {\bf (ii)} For $S \in \SimplicialSets$, its
\emph{PL de Rham complex} is the hom-object
of simplicial objects from $S$ to $\Omega^\bullet_{k\mathrm{pdR}}\big(\Delta^{(-)}\big)$
\eqref{pdR},
hence is the following {\it end}
(e.g. \cite[Def. 6.6.8]{Borceux94})
in $\dgcAlgebras{\mathbb{R}}$:
\vspace{-2mm}
\begin{equation}
  \label{PLdERhamComplex}
  \Omega^\bullet_{\mathrm{P}k\mathrm{LdR}}
  (
    S
  )
  \;:=\;
  \underset{[n] \in \Delta}{\int}
  \,
  \underset{S_n}{\prod}
  \Omega^\bullet_{k\mathrm{pdR}}
  (
    \Delta^n
  )\;.
\end{equation}

\vspace{-2mm}
\noindent
This means that an element
$\omega \,\in\, \Omega^\bullet_{\mathrm{P}k\mathrm{LdR}}(S)$
is a $k$-polynomial differential form
$\omega^{(n)}_{\sigma} \,\in\, \Omega^\bullet_{k\mathrm{pdR}}(\Delta^n)$
\eqref{pdR}
on each $n$-simplex $\sigma \in S_n$
for all $n \in \mathbb{N}$, such that these are compatible under
pullback along all simplex face inclusions $\delta_i$
and along all degenerate simplex projections $\sigma_i$:
\vspace{-3mm}
$$
  \Omega^\bullet_{\mathrm{P}k\mathrm{LdR}}(S)
  \;\;
  =
  \;\;
  \left\{
  \raisebox{68pt}{
  \xymatrix@C=3em{
    \vdots
    \ar@{..}@<+21pt>[d]
    \ar@{..}@<+14pt>[d]
    \ar@{..}@<+7pt>[d]
    \ar@{..}@<+0pt>[d]
    \ar@{..}@<-7pt>[d]
    \ar@{..}@<-14pt>[d]
    \ar@{..}@<-21pt>[d]
    &&
    \vdots
    \ar@{..}@<+21pt>[d]
    \ar@{..}@<+14pt>[d]
    \ar@{..}@<+7pt>[d]
    \ar@{..}@<+0pt>[d]
    \ar@{..}@<-7pt>[d]
    \ar@{..}@<-14pt>[d]
    \ar@{..}@<-21pt>[d]
    \\
    S_2
    \ar@{->}@<-18pt>[d]|-{
      \mathclap{\phantom{\vert^{\vert}_{\vert}}}
      \scalebox{.6}{$
        \delta^\ast_0
      $}
    }
    \ar@{<-}@<-9pt>[d]|-{
      \mathclap{\phantom{\vert^{\vert}_{\vert}}}
      \scalebox{.6}{$
        \sigma^\ast_0
      $}
    }
    \ar@{->}@<+0pt>[d]|-{
      \mathclap{\phantom{\vert^{\vert}_{\vert}}}
      \scalebox{.6}{$
        \delta^\ast_1
      $}
    }
    \ar@{<-}@<+9pt>[d]|-{
      \mathclap{\phantom{\vert^{\vert}_{\vert}}}
      \scalebox{.6}{$
        \sigma^\ast_1
      $}
    }
    \ar@{->}@<+18pt>[d]|-{
      \mathclap{\phantom{\vert^{\vert}_{\vert}}}
      \scalebox{.6}{$
        \delta^\ast_2
      $}
    }
    \ar@{-->}[rr]^{
      \omega^{(2)}_{(-)}
    }
    &&
    \Omega^\bullet_{k\mathrm{pdR}}(\Delta^2)
    \ar@{->}@<-18pt>[d]|-{
      \mathclap{\phantom{\vert^{\vert}_{\vert}}}
      \scalebox{.6}{$
        \delta^\ast_0
      $}
    }
    \ar@{<-}@<-9pt>[d]|-{
      \mathclap{\phantom{\vert^{\vert}_{\vert}}}
      \scalebox{.6}{$
        \sigma^\ast_0
      $}
    }
    \ar@{->}@<+0pt>[d]|-{
      \mathclap{\phantom{\vert^{\vert}_{\vert}}}
      \scalebox{.6}{$
        \delta^\ast_1
      $}
    }
    \ar@{<-}@<+9pt>[d]|-{
      \mathclap{\phantom{\vert^{\vert}_{\vert}}}
      \scalebox{.6}{$
        \sigma^\ast_1
      $}
    }
    \ar@{->}@<+18pt>[d]|-{
      \mathclap{\phantom{\vert^{\vert}_{\vert}}}
      \scalebox{.6}{$
        \delta^\ast_2
      $}
    }
    \\
    S_1
    \ar@{-->}[rr]^{
      \omega^{(1)}_{(-)}
    }
    \ar@{->}@<-9pt>[d]|-{
      \mathclap{\phantom{\vert^{\vert}_{\vert}}}
      \scalebox{.6}{$
        \delta^\ast_0
      $}
    }
    \ar@{<-}@<+0pt>[d]|-{
      \mathclap{\phantom{\vert^{\vert}_{\vert}}}
      \scalebox{.6}{$
        \sigma^\ast_0
      $}
    }
    \ar@{->}@<+9pt>[d]|-{
      \mathclap{\phantom{\vert^{\vert}_{\vert}}}
      \scalebox{.6}{$
        \delta^\ast_1
      $}
    }
    &&
    \Omega^\bullet_{k\mathrm{pdR}}(\Delta^1)
    \ar@{->}@<-9pt>[d]|-{
      \mathclap{\phantom{\vert^{\vert}_{\vert}}}
      \scalebox{.6}{$
        \delta^\ast_0
      $}
    }
    \ar@{<-}@<+0pt>[d]|-{
      \mathclap{\phantom{\vert^{\vert}_{\vert}}}
      \scalebox{.6}{$
        \sigma^\ast_0
      $}
    }
    \ar@{->}@<+9pt>[d]|-{
      \mathclap{\phantom{\vert^{\vert}_{\vert}}}
      \scalebox{.6}{$
        \delta^\ast_1
      $}
    }
    \\
    S_0
    \ar@{-->}[rr]^{
      \omega^{(0)}_{(-)}
    }
    &&
    \Omega^\bullet_{k\mathrm{pdR}}(\Delta^0)
  }
  }
  \right\}
  \,.
$$

\vspace{-1mm}
\noindent {\bf (iii)}  For $X \in \TopologicalSpaces$,
its \emph{PL de Rham complex} is that of its singular simplicial set,
according to \eqref{PLdERhamComplex}:
\vspace{-1mm}
\begin{equation}
  \label{PLdeRhamComplexOfTopologicalSpace}
  \Omega^\bullet_{\mathrm{P}k\mathrm{LdR}}
  (
    X
  )
  \;:=\;
  \Omega^\bullet_{\mathrm{P}k\mathrm{LdR}}
  \big(
    \mathrm{Sing}(X)
  \big)
  \,.
\end{equation}

\vspace{-2mm}
\noindent
By pullback of differential forms, this extends to a
functor of the form
\vspace{-2mm}
\begin{equation}
  \label{PLdeRhamComplexFunctor}
  \Omega^\bullet_{\mathrm{P}k\mathrm{LdR}}
  \;:\;
  \xymatrix{
    \SimplicialSets
    \ar[rr]^-{ }
    &&
    \dgcAlgebrasOp{\mathbb{R}}
  }.
\end{equation}

\vspace{-3mm}
\noindent {\bf (iv)} We write
\begin{equation}
  \label{PLdeRhamCohomology}
  H^\bullet_{\mathrm{P}k\mathrm{LdR}}(-)
  \;:=\;
  H \Omega^\bullet_{\mathrm{P}k\mathrm{LdR}}(-)
\end{equation}
for \emph{PL de Rham cohomology}, the cochain cohomology
of the PL de Rham complex.
\end{defn}

\begin{prop}[PL de Rham theorem {\cite[Thm. 2.2]{BousfieldGugenheim76}\cite[Thm. 9.1]{GriffithMorgan13}}]
  \label{PLdeRhamTheorem}
  The evident operation of integrating differential
  forms over simplices induces a quasi-isomorphism
  \vspace{-2mm}
  $$
    \xymatrix@C=3em{
      \Omega^\bullet_{\mathrm{P}k\mathrm{LdR}}(-)
      \ar[r]^-{ \in \, \mathrm{qIso} }
      &
      \; C^\bullet(-;\, k)
    }
  $$

\vspace{-2mm}
\noindent
  from the PL de Rham complex (Def. \ref{PLdeRhamComplexes})
  to the cochain complex of ordinary singular cohomology
  with coefficients in $k$. Hence on cochain cohomology
  this induces an isomorphism
  \vspace{-2mm}
  $$
    \xymatrix{
      H^\bullet_{\mathrm{P}k\mathrm{LdR}}(-)
      \ar[r]^-{ \simeq }
      &
    \;  H^\bullet(-;\, k)
    }
  $$

\vspace{-2mm}
\noindent
  between PL de Rham cohomology \eqref{PLdeRhamCohomology}
  and ordinary cohomology with coefficients in $k$.
\end{prop}

\begin{example}[PL de Rham complex of the interval]
  \label{PLdeRhamComplexOfTheInterval}
  The PL de Rham complex (Def. \ref{PLdeRhamComplexes}) of the
  1-simplex,
  hence the polyonial differential forms \eqref{pdR} on $\Delta^1$,
  is isomorphic to the multivariate polynomial dgc-algebra
  (Ex. \ref{MultivariatePolynomialdgcAlgebras}) of the form
  \begin{equation}
    \label{PLdeRhamComplexOfInterval}
    \Omega^\bullet_{\mathrm{PL}k\mathrm{dR}}(\Delta^1)
    \;\simeq\;
    \Omega^\bullet_{k\mathrm{pdR}}(\Delta^1)
    \;\;\simeq\;\;
    k
    \big[
      t_0,
      \theta_1
    \big]
    \big/
    \big(
      d\, t_0 \,=\, \theta_1
    \big)
    \,.
  \end{equation}
  For $A \in \dgcAlgebrasProj{k}$ (Prop. \ref{ProjectiveModelStructureOnConnectiveChainComplexes}), its tensor product with \eqref{PLdeRhamComplexOfInterval}

  \vspace{-.4cm}
  $$
    A
      \otimes_{k}
    \Omega^\bullet_{\mathrm{PL}k\mathrm{dR}}(\Delta^1)
    \;\;
      \simeq
    \;\;
    A
    \big[
      t_0, \,
      \theta_1
    \big]
    \big/
    \big(
      d\, t_0 \,=\, \theta_1
    \big)
  $$
  \vspace{-.4cm}

  \noindent
  is a path space object for $A$ (Def. \ref{PathSpaceObject}), in that
  it fits into the following diagram

  \vspace{-.4cm}
  \begin{equation}
    \label{CanonicalPathSpaceObjectFordgcAlgebras}
    \begin{tikzcd}
      A
      \ar[
        rr,
        "a \,\mapsto\, a"{above},
        "\in \mathrm{W}"{below}
      ]
      \ar[
      rrrrr,
      rounded corners,
      to path={
           -- ([yshift=-10pt]\tikztostart.south)
           --node[below]{
               \scalebox{.8}{$
                 \Delta_A
               $}
             } ([yshift=-8pt]\tikztotarget.south)
           -- (\tikztotarget.south)}
      ]
      &&
      A
      \big[
        t_0,
        \,
        \theta_1
      \big]
      \ar[
        rr,
        "{
          a
          \,\mapsto\,
          \big(
            a(t_0 = 0, \theta_1 = 0)
            +
            a(t_0 = 1, \theta_1 = 0)
          \big)
        }"{above},
        "\in \mathrm{Fib}"{below}
      ]
      &[+76pt] &
      A \oplus A
      \ar[r, phantom, "\simeq"]
      &[-10pt]
      A \times A
      \,,
    \end{tikzcd}
  \end{equation}
  \vspace{-.4cm}

  \noindent
  where the morphism on the right is given by evaluation of polynomials
  as shown, and where the equivalence on the right is by Ex. \ref{ProductAndCoproductAlgebras}.

  Here the morphism on the right is a degreewise surjection
  (a pre-image for $(a_0, a_1) \in A \oplus A$ is
  $(1 - t_0) \wedge a_0   + t_0 \wedge a_1 \,\in\, A[t_0, \theta_1]$),
  hence a fibration according to Def. \ref{HomotopicalStructureOndgcAlgebras};
  while the morphism on the left is a quasi-isomorphism, in fact a chain homotopy
  equivalence (with homotopy operator $t_0 \cdot \frac{\partial}{\partial \theta_1}$), hence
  a weak equivalence according to Def. \ref{HomotopicalStructureOndgcAlgebras}.
\end{example}

The following type of argument will be greatly expanded on in \cref{NonAbelianDeRhamTheory}:
\begin{lemma}[Homotopical formulation of ordinary cohomology]
  \label{HomotopicalFormulationOfCochainCohomology}
  $\,$

  \noindent
  {\bf (i)}
  The cochain cohomology of any $A \,\in\, \dgcAlgebras{\mathbb{R}}$
  (Def. \ref{CategoryOfdgAlgebras})
  in positive degree
  is naturally and $\mathbb{R}$-linearly identified,
  \begin{equation}
    \label{CochainCohomologyOfDGCAlgebrasAsHomSetsInHomotopyCategory}
    H^{\bullet + 1}(A)
    \;\;
      \simeq
    \;\;
    \mathrm{Ho}
    \Big(
      \dgcAlgebrasProj{\mathbb{R}}
    \big)
    \big(
      \mathrm{CE}(\mathfrak{b}^\bullet\mathbb{R})
      ,\,
      A
    \big)
    \,,
  \end{equation}
  with the hom-sets out of the
  CE-algebras \eqref{CEAlgebraOfLineLienAlgebra} of
  the line Lie $(\bullet+1)$-algebras (Ex. \ref{LineLienPlusOneAlgebras})
  in the homotopy category (Def. \ref{HomotopyCategory})
  of the projective model category (Prop. \ref{ProjectiveModelStructureOndgcAlgebras}).

  \noindent
  {\bf (ii)}
  For $X \,\in\, \TopologicalSpaces \xrightarrow{\mathrm{Sing}} \SimplicialSets$,
  its real cohomology in positive degree
  is naturally identified with these hom-sets into its
  PL de Rham complex (Def. \ref{PLdeRhamComplexes})
  \begin{equation}
    \label{RealCohomologyAsHomSetsInHomotopyCategory}
    H^{\bullet + 1}(X;\, \mathbb{R})
    \;\;
    \simeq
    \;\;
    \mathrm{Ho}
    \Big(
      \dgcAlgebrasProj{R}
    \big)
    \big(
      \mathrm{CE}(\mathfrak{b}^\bullet\mathbb{R})
      ,\,
      \Omega^\bullet_{\mathrm{P}\mathbb{R}\mathrm{LdR}}(X)
    \big)
    \,.
  \end{equation}

  \noindent
  {\bf (iii)}
  If the above $X$ is equipped with a base-point $\ast \xrightarrow{x} X$,
  then the real cohomology of $X$ in positive degrees
  is equivalently computed by the homotopy classes
  of morphisms of augmented dgc-algebras, hence with respect to the slice model structure
  (Ex. \ref{SliceModelCategory}) over the initial dgc-algebra $\mathbb{R}$ (Ex. \ref{InitialGradedAlgebra}),
  as follows:
  \begin{equation}
    \label{RealCohomologyAsHomSetsInSlicedHomotopyCategory}
    H^{\bullet + 1}(X;\, \mathbb{R})
    \;\;
    \simeq
    \;\;
    \mathrm{Ho}
    \Big(
      \dgcAlgebrasProj{R}^{/\mathbb{R}}
    \big)
    \left(
      \mathrm{CE}(\mathfrak{b}^\bullet\mathbb{R})
      ,\,
      \Omega^\bullet_{\mathrm{P}\mathbb{R}\mathrm{LdR}}(X)
      \!\!\!
      \underset{
          \Omega^\bullet_{\mathrm{P}\mathbb{R}\mathrm{LdR}}(\ast)
      }{\times}
      \!\!\!
      \mathbb{R}
    \right)
    \,.
  \end{equation}
\end{lemma}
\begin{proof}
  Observing that

  {(a)}
  $\mathrm{CE}(\mathfrak{b}^n\mathbb{R})$ is cofibrant, by Prop. \ref{RelativeSullivanModelsAreCofibrations};

  {(b)} $A$ is fibrant, by Rem. \ref{AlldgcAlgebrasAreProjectivelyFibrant};

  {(c)} $A[t_0, \theta_1]/(d\, t_0 \,=\, \theta_1)$ is a path space object for $A$,
  by Expl. \ref{PLdeRhamComplexOfTheInterval};

  \noindent
  we may identify, by Prop. \ref{RightHomotopyReflectsHomotopyClasses},
  the morphisms in the homotopy category
  with equivalence classes of dgc-algebra morphisms
  $
    \mathrm{CE}(\mathfrak{b}^n \mathbb{R})
    \xrightarrow{\;c\;}
    A
  $
  under the corresponding equivalence relation of right homotopy
  (Def. \ref{RightHomotopy}):

  \vspace{-.4cm}
  $$
    c \,\sim_r\, c'
    \;\;\;\;\;\;\;\;\;\;\;\;
    \Leftrightarrow
    \;\;\;\;\;\;\;\;\;\;\;\;
    \;\;
    \begin{tikzcd}[row sep=small]
      &&
      A
      \\
      \mathrm{CE}(\mathfrak{b}^n\mathbb{R})
      \ar[
        rr,
        dashed,
        "\exists"
      ]
      \ar[
        urr,
        "c"{above}
      ]
      \ar[
        drr,
        "c'"{below}
      ]
      &&
      \frac{
        A[t_0, \theta_1]
      }
      {
        (d\, t_0 \,=\, \theta_1)
      }
      \ar[
        u,
        "{
          (-)_{\vert t_0 = 0, \theta_1 = 0}
        }"{right}
      ]
      \ar[
        d,
        "{
          (-)_{\vert t_0 = 1, \theta_1 = 0}
        }"{right}
      ]
      \\
      &&
      A
      \mathrlap{\,.}
    \end{tikzcd}
  $$
  \vspace{-.3cm}

  Since $\mathrm{CE}(\mathfrak{b}^n \mathbb{R})$ is free on a single
  generator $\theta_n$ in degree $n$, subject only to the differential relation
  $d \, \theta_n \,=\, 0$,
  dgc-algebra homomorphisms
  $\mathrm{CE}(\mathfrak{b}^n \mathbb{R}) \xrightarrow{\;} A$ are
  in bijection to closed degree-$n$ elements of $A$
  (see also Ex. \ref{OrdinaryClosedFormsAreFlatLineLInfinityAlgebraValuedForms}).
  Hence, under this identification it remains to see that
  existence of coboundaries is equivalent to existence of right homotopies:

  \vspace{-2mm}
  $$
    \underset{ h \in A }{\exists}
    \;
    c' \,=\, c + d h
    \;\;\;\;\;\;\;\;\;\;\;
    \Leftrightarrow
    \;\;\;\;\;\;\;\;\;\;\;
    \underset{
      \eta
      \,\in\,
      \frac{
        A[t_0, \theta_1]
      }{
        (d\, t_0 \,=\, \theta_1)
      }
    }{\exists}
    \;
    \left\{
    \begin{array}{l}
      d \eta \, = 0 \,,
      \\
      \eta(t_0 = 0, \theta_1 = 0) = c \,,
      \\
      \eta(t_0 = 1, \theta_1 = 0) = c'
      \,.
    \end{array}
    \right.
  $$
  Indeed: If $h$ is given as on the left, then
  $
    \eta
      \,\coloneqq\,
    t_0 \wedge c'
      \,+\,
    (1- t_0)  \wedge c
      \,+\,
    \theta_1 \wedge ( c - c')
  $
  is as required on the right;
  while if any $\eta$ is given as on the right, then
  $
    h
    \;\coloneqq\;
    \int_0^1 \eta \, d t_0
  $
  is as required on the left (by Stokes, as in Lem. \ref{FiberwiseStokesTheorem}).
  This proves the first statement, whence the second follows via
  Prop. \ref{PLdeRhamTheorem}.

  To deduce from this the third statement, observe that:

  \begin{itemize}

  \vspace{-.2cm}
  \item[(a)]
  for $A \xrightarrow{\epsilon_A} \mathbb{R}$ an augmented dgc-algebra,
  the canonical path object \eqref{CanonicalPathSpaceObjectFordgcAlgebras} for $A$
  yields a path space object in the slice over $\mathbb{R}$
  by equipping it with the induced augmentation
   $
     \frac{A [t_0, \theta_1]}{(d\, t_0 = \theta_1)}
     \xrightarrow{ a \,\mapsto\, a(t_0 = 0, \theta_1 = 0) }
     A
     \xrightarrow{\epsilon_A}
     \mathbb{R}
     \,,
   $

  \vspace{-.2cm}
  \item[(b)]
   the projection
   \begin{equation}
     \label{QuasiIsomorphismBetweenAugmentedPLdRComplexAndOriginalPLdRham}
     \Omega^\bullet_{\mathrm{P}\mathbb{R}\mathrm{LdR}}(X)
     \!\!\
     \underset{
       \Omega^\bullet_{\mathrm{P}\mathbb{R}\mathrm{LdR}}(\ast)
     }{\times}
     \!\!\!
     \mathbb{R}
     \xrightarrow{\;\;\; \in \mathrm{W} \;\;\;}
     \Omega^\bullet_{\mathrm{P}\mathbb{R}\mathrm{LdR}}(X)
   \end{equation}
   is a quasi-isomorphism,
   by right-properness (Def. \ref{ProperModelCategories})
   of the model structure
   (Prop. \ref{ProjectiveModelStructureOndgcAlgebras}),
   since this is the pullback of the quasi-isomorphism
   $\mathbb{R} \xrightarrow{\in \mathrm{W}} \Omega^\bullet_{\mathrm{P}\mathbb{R}\mathrm{LdR}}$
   (by Prop. \ref{PLdeRhamTheorem}) along the morphism
   $\Omega^\bullet_{\mathrm{P}\mathbb{R}\mathrm{LdR}}(\ast \xrightarrow{\;} X)$,
   which is a projective fibration by the fact that $\ast \to X$ is an injection
   and hence a cofibration (Ex. \ref{ClassicalModelStructureOnSimplicialSets})
   and that $\Omega^\bullet_{\mathrm{P}\mathbb{R}\mathrm{LdR}}$ is a left Quillen functor
   (Prop. \ref{QuillenAdjunctionBetweendgcAlgebrasAndSimplicialSets})
   to the opposite model structure (Ex. \ref{OppositeModelStructure}).

   \end{itemize}

\vspace{-1mm}
   \noindent
   Therefore, since the generators $\mathrm{CE}(\mathfrak{b}^n \mathbb{R})$
   are in positive degree and hence unaffected by the augmentation slicing,
   the right homotopy classes in the slice \eqref{RealCohomologyAsHomSetsInSlicedHomotopyCategory}
   are computed as in case (2) above
   and hence yield the cochain cohomology of
   $\Omega^\bullet_{\mathrm{P}\mathbb{R}\mathrm{LdR}}(X)
      \times_{\Omega^\bullet_{\mathrm{P}\mathbb{R}\mathrm{LdR}}(\ast) }
     \mathbb{R}
   $,
   which by \eqref{QuasiIsomorphismBetweenAugmentedPLdRComplexAndOriginalPLdRham}
   equals the real cohomology of $X$.
\end{proof}

In fact, before passing to cochain cohomology,
the PL de Rham complex captures the full rational homotopy type.
This is the Fundamental Theorem
which we recall as Prop. \ref{FundamentalTheoremOfdgcAlgebraicRationalHomotopyTheory}:

\begin{lemma}[Extension lemma for polynomial differential forms {\cite[Lemma 9.4]{GriffithMorgan13}}]
  \label{ExtensionLemmaForPolynomialDifferentialForms}
  For $n \in \mathbb{N}$, the operation of
  pullback of piecewise polynomial differential forms
  (Def. \ref{PLdeRhamComplexOnSmoothManifoldsEquivalentToSmoothDeRhamComplex})
  along the boundary inclusion of the $n$-simplex
  $
    \xymatrix@C=12pt{
      \partial \Delta^n
      \ar[r]^-{ i_n }
      &
      \Delta^n
    }
  $
  is an epimorphism:

  \vspace{-4mm}
  $$
  \begin{tikzcd}
      \Omega^\bullet_{\mathrm{P}k\mathrm{LdR}}(\Delta^n)
      \ar[
        rr,
        ->>,
        "{ i_n^\ast }"
      ]
      &&
      \Omega^\bullet_{\mathrm{P}k\mathrm{LdR}}(\partial \Delta^n)
      \,.
  \end{tikzcd}
  $$
\end{lemma}

\begin{prop}[PL de Rham Quillen adjunction {\cite[\S 8]{BousfieldGugenheim76}}]
  \label{QuillenAdjunctionBetweendgcAlgebrasAndSimplicialSets}
  For all ground fields $k$ of characteristic zero,
  the PL de Rham complex functor (Def. \ref{PLdeRhamComplexes})
  is the left adjoint in a Quillen adjunction (Def. \ref{QuillenAdjunction})
  \vspace{-4mm}
  \begin{equation}
    \label{QuillenAdjunctBetweendgcAlgsAndSimplicialSets}
    \xymatrix{
      \dgcAlgebrasOpProj{k}
      \;\;
      \ar@{<-}@<+7pt>[rr]^-{
        \Omega^\bullet_{\mathrm{P}k\mathrm{LdR}}
      }
      \ar@<-7pt>[rr]^-{\bot_{\mathrlap{\mathrm{Qu}}}}_-{
        \Bexp_{\mathrm{P}k\mathrm{L}}
      }
      &&
      \;\;
      \SimplicialSets_{\mathrm{Qu}}
    }
  \end{equation}

  \vspace{-3mm}
  \noindent
  between
  the opposite (Def. \ref{OppositeModelStructure}) of
  the model category of dgc-algebras
  (Prop. \ref{ProjectiveModelStructureOndgcAlgebras})
  and the classical model structure on simplicial sets
  (Prop. \ref{ClassicalModelStructureOnSimplicialSets});
  where the right adjoint sends a dgc-algebra $A$ to
  \vspace{-1mm}
  \begin{equation}
    \label{exp}
    \Bexp_{\mathrm{P}k\mathrm{L}}(A)
    \;=\;
    \Big(
    \Delta[n] \;\longmapsto\;
    \dgcAlgebras{k}
    \big(
      \Omega^\bullet_{\mathrm{P}k\mathrm{LdR}}(\Delta^n)
      \,,\,
      A
    \big)
    \Big)
    \;\;\;\;\;
    \in
    \;
    \SimplicialSets\;.
  \end{equation}

\end{prop}
\begin{proof}
  That the right adjoint exists and is given as in
  \eqref{exp} follows by general nerve/realization theory
  \cite{Kan58}, or else by direct inspection.

  For the left adjoint to preserve cofibrations means to take
  injections of simplicial sets to degreewise surjections of
  dgc-algebras. This follows from the extension lemma
  (Lemma \ref{ExtensionLemmaForPolynomialDifferentialForms}).
  Moreover, the left adjoint preserves even all
  weak equivalences, by the PL de Rham theorem (Prop. \ref{PLdeRhamTheorem}).
\end{proof}

\begin{prop}[Fundamental theorem of dgc-algebraic rational homotopy theory]
  \label{FundamentalTheoremOfdgcAlgebraicRationalHomotopyTheory}
  For $k = \mathbb{Q}$,
  the derived adjunction (Prop. \ref{DerivedFunctors})
  \vspace{-2mm}
  \begin{equation}
    \label{PLdeRhamDerivedAdjunction}
    \xymatrix{
      \mathrm{Ho}
      \big(
        \dgcAlgebrasOpProj{\mathbb{Q}}
      \big)
      \;
      \ar@{<-}@<+5pt>[rr]^-{
        \LeftDerived \Omega^\bullet_{\mathrm{P}\mathbb{Q}\mathrm{LdR}}
      }
      \ar@<-5pt>[rr]_-{
        \RightDerived \Bexp_{\mathrm{P}\mathbb{Q}\mathrm{L}}
      }^-{\bot}
      &&
      \;\;
      \mathrm{Ho}
      \big(
        \SimplicialSets_{\mathrm{Qu}}
      \big)
    }
  \end{equation}

  \vspace{-2mm}
  \noindent
  of the Quillen adjunction \eqref{QuillenAdjunctBetweendgcAlgsAndSimplicialSets}
  from Prop. \ref{QuillenAdjunctionBetweendgcAlgebrasAndSimplicialSets}
  is such that:

  \noindent
  {\bf (i)} on connected, nilpotent, $\mathbb{Q}$-finite homotopy types
  (Def. \ref{NilpotentConnectedSpacesOfFiniteRationalType})
  the derived PLdR-adjunction unit \eqref{DerivedAdjunctionUnit}
  is equivalently the unit \eqref{RationalizationUnit}
  of rationalization (Def. \ref{Rationalization}):
  \vspace{-2mm}
  \begin{equation}
    \label{RationalizationViaPLDeRham}
    \hspace{5mm}
    \raisebox{17pt}{
    \xymatrix@R=20pt@C=4em{
      X
      \ar[rr]_-{\mathbb{D}\eta^{\mathrm{P}\mathbb{Q}\mathrm{LdR}}_X}^-{
        \overset{
          \mathclap{
          \raisebox{-1pt}{
            \tiny
            \color{greenii}
            \bf
            \def\arraystretch{.9}
            \begin{tabular}{c}
              derived unit of rational
              \\
              PL de Rham adjunction
            \end{tabular}
          }
          }
        }{
        }
      }
      \ar@{=}[d]
      &&
      \mathbb{R}\Bexp_{\mathrm{P}\mathbb{Q}\mathrm{L}}
      \,\circ\,
      \Omega^\bullet_{\mathrm{P}\mathbb{Q}\mathrm{LdR}}
      (X)
      \ar@{}[d]|-{
        \rotatebox[origin=c]{-90}{$
          \simeq
        $}
      }
      \\
      X
      \ar[rr]^-{
        \eta^{\mathbb{Q}}_X
      }_-{
        \underset{
          \mathclap{
          \raisebox{2pt}{
            \tiny
            \color{greenii}
            \bf
            \def\arraystretch{.9}
            \begin{tabular}{c}
              rationalization unit
            \end{tabular}
          }
          }
        }{
 %         \eta^{\mathbb{R}}
        }
      }
      &&
      L_{\mathbb{Q}} X
    }
    }
    \;
    \in
    \;
    \mathrm{Ho}
    \big(
      \TopologicalSpaces_{\mathrm{Qu}}
    \big)_{\geq 1, \mathrm{nil}}^{\mathrm{fin}_{\mathbb{Q}}}\;.
  \end{equation}

  \vspace{-2mm}
  \noindent {\bf (ii)}
  For
  $X, A$ nilpotent, connected, $\mathbb{Q}$-finite homotopy types
  (Def. \ref{NilpotentConnectedSpacesOfFiniteRationalType}),
  the PL de Rham space functor \eqref{PLdeRhamComplexFunctor}
  from Def. \ref{PLdeRhamComplexes}
  induces natural bijections
  \vspace{-2mm}
  \begin{equation}
    \label{FundamentalRationalizationEquivalence}
    \hspace{-4mm}
    \xymatrix{
    \mathrm{Ho}
    \!\left(
      \TopologicalSpaces_{\mathrm{Qu}}
    \right)
    \!\!
    \big(
      X
      ,
      L_{\mathbb{Q}} A
    \big)
    \ar[rr]^-{ \simeq }_-{ \Omega^\bullet_{\mathrm{P}\mathbb{Q}\mathrm{LdR}} }
    &&
    \mathrm{Ho}
    \!\!
    \left(\!\!
      \dgcAlgebrasProj{\mathbb{Q}}
    \right)
    \!\!
    \Big(
      \Omega^\bullet_{\mathrm{P}\mathbb{Q}\mathrm{LdR}}(A)
      \,,\,
      \Omega^\bullet_{\mathrm{P}\mathbb{Q}\mathrm{LdR}}(X)
    \!\Big).
    }
  \end{equation}
\end{prop}
\begin{proof}
  {\bf (i)} This is \cite[Thm 11.2]{BousfieldGugenheim76}.

  \noindent
  {\bf (ii)}
  This follows via \cite[Thm 9.4(i)]{BousfieldGugenheim76},
  which says that the derived adjunction
  \eqref{PLdeRhamDerivedAdjunction} restricts on
  connected, nilpotent, $\mathbb{Q}$-finite
  (Def. \ref{NilpotentConnectedSpacesOfFiniteRationalType})
  rational homotopy types (Def. \ref{Rationalization}) to an
  equivalence of homotopy categories:
  \vspace{-1mm}
  \begin{equation}
    \label{BousfieldGugenheimTheorem}
    \xymatrix{
      \mathrm{Ho}
      \big(
        \dgcAlgebrasOpProj{\mathbb{Q}}
      \big)^{\geq 1}_{\mathrm{fin}}
     \;\; \ar@{<-}@<+6pt>[rr]^-{
        \LeftDerived \Omega^\bullet_{\mathrm{P}\mathbb{Q}\mathrm{LdR}}
      }
      \ar@<-6pt>[rr]_-{
        \RightDerived \Bexp_{\mathrm{P}\mathbb{Q}\mathrm{L}}
      }^-{\simeq}
      &&
     \;\; \mathrm{Ho}
      \big(
        \SimplicialSets_{\mathrm{Qu}}
      \big)_{\geq 1, \mathrm{nil}}^{\mathbb{Q}, \mathrm{fin}_{\mathbb{Q}}}
    }.
  \end{equation}

\vspace{-2mm}
\noindent
In detail, this is
witnessed by the following sequence of natural bijections of hom-sets:
\vspace{-0mm}
  \begin{equation}
    \label{TowardsTheFundamentalTheoremOfdgcAlgebraicRationalHomotopyTheory}
    \begin{aligned}
    &
    \mathrm{Ho}
    \left(
      \TopologicalSpaces_{\mathrm{Qu}}
    \right)
    \big(
      X
      \,,\,
      L_{\mathbb{Q}} A
    \big)
    \\
    & \simeq
    \mathrm{Ho}
    \left(
      \SimplicialSets_{\mathrm{Qu}}
    \right)
    \big(
      \mathrm{Sing}(X)
      \,,\,
      L_{\mathbb{Q}}
      \mathrm{Sing}(A)
    \big)
    \\
    & \simeq
    \mathrm{Ho}
    \left(
      \SimplicialSets_{\mathrm{Qu}}
    \right)
    \big(
      L_{\mathbb{Q}} \mathrm{Sing}(X)
      \,,\,
      L_{\mathbb{Q}}
      \mathrm{Sing}(A)
    \big)
    \\
    & \simeq
    \mathrm{Ho}
    \left(
      \SimplicialSets_{\mathrm{Qu}}
    \right)
    \big(
      \mathbb{R}\Bexp_{\mathrm{P}\mathbb{Q}\mathrm{L}}
        \,\circ\,
      \Omega^\bullet_{\mathrm{P}\mathbb{Q}\mathrm{LdR}}(X)
      \,,\,
      \mathbb{R}\Bexp_{\mathrm{PL}}
        \,\circ\,
      \Omega^\bullet_{\mathrm{P}\mathbb{Q}\mathrm{LdR}}(A)
    \big)
    \\
    & \simeq
    \mathrm{Ho}
    \left(
      \dgcAlgebrasOpProj{\mathbb{Q}}
    \right)
    \big(
      \Omega^\bullet_{\mathrm{P}\mathbb{Q}\mathrm{LdR}}
      \,\circ\,
      \mathbb{R}\Bexp_{\mathrm{P}\mathbb{Q}\mathrm{L}}
        \,\circ\,
      \Omega^\bullet_{\mathrm{P}\mathbb{Q}\mathrm{LdR}}(X)
      \,,\,
      \Omega^\bullet_{\mathrm{P}\mathbb{Q}\mathrm{LdR}}
      \,\circ\,
      \mathbb{R}\Bexp_{\mathrm{P}\mathbb{Q}\mathrm{L}}
        \,\circ\,
      \Omega^\bullet_{\mathrm{P}\mathbb{Q}\mathrm{LdR}}(A)
    \big)
    \\
    & \simeq
    \mathrm{Ho}
    \left(
      \dgcAlgebrasOpProj{\mathbb{Q}}
    \right)
    \big(
      \Omega^\bullet_{\mathrm{P}\mathbb{Q}\mathrm{LdR}}(X)
      \,,\,
      \Omega^\bullet_{\mathrm{P}\mathbb{Q}\mathrm{LdR}}(A)
    \big)
    \\
    &
    \simeq
    \mathrm{Ho}
    \left(
      \dgcAlgebrasProj{\mathbb{Q}}
    \right)
    \Big(
      \Omega^\bullet_{\mathrm{P}\mathbb{Q}\mathrm{LdR}}(A)
      \,,\,
      \Omega^\bullet_{\mathrm{P}\mathbb{Q}\mathrm{LdR}}(X)
    \Big)
    \,.
    \end{aligned}
  \end{equation}

  \vspace{-1mm}
  \noindent
  Here
  the first step is \eqref{ClassicalHomotopyCategory};
  the second step uses that rationalization
  is a reflection \eqref{RationalizationReflection};
  the third step uses \eqref{RationalizationViaPLDeRham};
  the fourth is the equivalence \eqref{BousfieldGugenheimTheorem}
    along $\LeftDerived\Omega^\bullet_{\mathrm{P}\mathbb{Q}\mathrm{LdR}}$
    (using, with Example \ref{DerivedFunctorsByCoFibrantReplacement},
    that every simplicial set is already cofibrant \eqref{EverySimplicialSetIsCofibrantInClassicalModelStructure},
    Example \ref{ClassicalModelStructureOnSimplicialSets});
  the fifth step is the statement from \eqref{BousfieldGugenheimTheorem}
   that $\mathbb{R}\Bexp_{\mathrm{P}\mathbb{Q}\mathrm{L}}$ is the inverse equivalence.
  The last step is just the definition of the opposite
    of a category.
  The composite of the bijections \eqref{TowardsTheFundamentalTheoremOfdgcAlgebraicRationalHomotopyTheory}
  is the desired bijection \eqref{FundamentalRationalizationEquivalence}.
\hfill \end{proof}

In view of \eqref{RationalizationViaPLDeRham} the following
notation is convenient, keeping in mind that $L_k$ is a localization
in the sense of localization of spaces only for $k = \mathbb{Q}$:
\begin{defn}[Rationalization over $\mathbb{R}$]
  \label{Lk}
  For $k $ a field of characteristic zero, we write
  $
    L_k
    \;\coloneqq\;
    \RightDerived
    \Bexp_{\mathrm{P}k\mathrm{L}}
    \,\circ\,
    \LeftDerived \Omega^\bullet_{\mathrm{P}k\mathrm{LdR}}
  $
  for the monad given by the derived functors (Prop. \ref{DerivedFunctors})
  of the $k$-PL de Rham Quillen adjunction (Prop. \ref{QuillenAdjunctionBetweendgcAlgebrasAndSimplicialSets}).
  Our focus here is on the case over the real numbers:
  \vspace{-2mm}
  \begin{equation}
    \label{RealificationMonad}
    (-)
    \xrightarrow{
      \;
      \eta^{\mathbb{R}}
      \,\coloneqq\,
      \mathbb{D}\eta^{\mathrm{PL}\mathbb{R}{dR}}
      \;
    }
    L_{\mathbb{R}}(-)
    \;\coloneqq\;
    \RightDerived
    \Bexp_{\mathrm{P}{\mathbb{R}}\mathrm{L}}
    \,\circ\,
    \LeftDerived \Omega^\bullet_{\mathrm{P}{\mathbb{R}}\mathrm{LdR}}(0) \;.
  \end{equation}
\end{defn}
We may refer to $L_{\mathbb{R}}$ as {\it rationalization over $\mathbb{R}$}.
Because, while the derived PLdR-adjunction (Prop. \ref{QuillenAdjunctionBetweendgcAlgebrasAndSimplicialSets})
is a localization of homotopy types only over $k = \mathbb{Q}$,
(Prop. \ref{FundamentalTheoremOfdgcAlgebraicRationalHomotopyTheory},
Rem. \ref{FailureOfRationalizationOverRToBeIdempotent}),
for general $k$ it is the suitable {\it change of scalars} of $\mathbb{Q}$-localization:

\begin{lemma}[Derived change of scalars {\cite[Lem. 11.6]{BousfieldGugenheim76}}]
  \label{DerivedChangeOfScalars}
  For $k$ a field of characteristic zero, the
  {\it extension/restriction of scalars}-adjunction along
  $\mathbb{Q} \xhookrightarrow{\;} k$
  is a Quillen adjunction (Def. \ref{QuillenAdjunction})
  between
  the corresponding projective model categories of dgc-algebras
  (from Prop. \ref{ProjectiveModelStructureOnConnectiveChainComplexes}):
 \vspace{-4mm}
  $$
    \begin{tikzcd}
      \dgcAlgebrasProj{{\color{blue}k}}
      \ar[
        rr,
        shift right=8pt,
        "\mathrm{res}_{\mathbb{Q}}"{below}
      ]
      \ar[
        rr,
        phantom,
        "\scalebox{.8}{$\bot_{\mathrlap{\mathrm{Qu}}}$}"
      ]
      &&
      \dgcAlgebrasProj{\color{blue}\mathbb{Q}}
      \,.
      \ar[
        ll,
        shift right=8pt,
        "(-)\otimes_{{}_{\mathbb{Q}}} k"{above}
      ]
    \end{tikzcd}
  $$
\end{lemma}
\begin{proof}
  Since restriction of scalars $\mathrm{res}_{\mathbb{Q}}$
  is the identity on the underlying sets of a dgc-algebra,
  it manifestly preserves all fibrations
  (since the are the surjections of underlying sets )
  and all weak equivalences
  (since these are the bijections on underlying cochain cohomology groups).
\end{proof}

\begin{prop}[P$k$LdR-Adjunction factors through rationalization]
  \label{kPLdRAdjunctionFactorsThroughRationalization}
  The following
  diagram of derived functors
  (Prop. \ref{DerivedFunctors} --
  with the left derived functors from Prop. \ref{QuillenAdjunctionBetweendgcAlgebrasAndSimplicialSets}
  and the right
  derived functor from Lem. \ref{DerivedChangeOfScalars})
  commutes up to natural isomorphism:
\begin{equation}
  \label{DerivedExtensionOfScalarsCompatibleWithOmegaPLdR}
  \begin{tikzcd}[row sep=small]
    \mathrm{Ho}
    \Big(
      \big(
        \mathrm{dgcAlg}^{\geq 0 }_{\color{blue}\mathbb{Q}}
      \big)^{\mathrm{op}}_{\mathrm{proj}}
    \Big)
    \ar[
      dd,
      "{
        \RightDerived
        \big(
          (-) \otimes_{\mathbb{Q}} \mathbb{R}
        \big)
      }"{left}
    ]
    &&
    \;\;
    \mathrm{Ho}
    \big(
      \SimplicialSets_{\mathrm{Qu}}
    \big)^{\mathrm{fin}_{\mathbb{Q}}}
    \ar[
      ll,
      "
        \LeftDerived
        \Omega^\bullet_{\mathrm{P}{\color{blue}\mathbb{Q}}\mathrm{LdR}}
      "{above}
    ]
    \ar[
      dd,-,
      shift right=1pt
    ]
    \ar[
      dd,-,
      shift left=1pt
    ]
    \\
    \\
    \mathrm{Ho}
    \Big(
    \big(
      \mathrm{dgcAlg}^{\geq 0}_{\color{blue}\mathbb{R}}
    \big)^{\mathrm{op}}_{\mathrm{proj}}
    \Big)
    &&
    \;\;
    \mathrm{Ho}
    \big(
      \SimplicialSets_{\mathrm{Qu}}
    \big)^{\mathrm{fin}_{\mathbb{Q}}}
    \ar[
      ll,
      "
        \LeftDerived
        \Omega^\bullet_{\mathrm{P}{\color{blue}k}\mathrm{LdR}}
      "{below}
    ]
  \end{tikzcd}
\end{equation}
\end{prop}
\begin{proof}
Via the formula for derived functors in terms of (co)fibrant replacement
(Ex. \ref{DerivedFunctorsByCoFibrantReplacement}) and using that
Sullivan models are cofibrant in $\dgcAlgebrasProj{k}$ (Prop. \ref{RelativeSullivanModelsAreCofibrations}), hence fibrant in
$\dgcAlgebrasOpProj{k}$, this follows by \cite[Lem. 11.7]{BousfieldGugenheim76}.
\end{proof}

\begin{remark}[Rational homotopy theory over the real numbers]
  \label{RationalHomotopyTheoryOverTheRealNumbers}
  Below in Prop. \ref{RealificationIsRationalizationFollowedByExtensionOfScalars}
  we recast Prop. \ref{kPLdRAdjunctionFactorsThroughRationalization}
  as the statement that the real character map on non-abelian cohomology
  factors through the rational character map via extension of scalars.

  This fact motivates and justifies the focus on
  {\it rational homotopy theory over the real numbers}
  (as in \cite{DGMS75} \cite{GriffithMorgan13})
  in all of the following.
  Rational homotopy theory over the real numbers is the version that
  connects to differential geometry (e.g. \cite{FOT08}),
  since the \emph{smooth} de Rham complex is not defined
  over $\mathbb{Q}$ but over $\mathbb{R}$
  (see Lemma
  \ref{PLdeRhamComplexOnSmoothManifoldsEquivalentToSmoothDeRhamComplex}).
  The original account \cite{BousfieldGugenheim76}
  of rational homotopy theory is, for the most part,
  formulated over an arbitrary field $k$ of characteristic zero;
  and \cite[Lem. 11.7]{BousfieldGugenheim76} (Prop. \ref{kPLdRAdjunctionFactorsThroughRationalization})
  makes explicit that the choice of this base field does not
  change the form of the classical theorems.
  For example, the ``real-ified'' homotopy groups of a space $X$
  \vspace{-1mm}
  $$
    \pi_\bullet(X) \otimes_{\scalebox{.5}{$\mathbb{Z}$}} \mathbb{R}
    \;\simeq\;
    \big(
      \pi_\bullet(X) \otimes_{\scalebox{.5}{$\mathbb{Z}$}} \mathbb{Q}
    \big) \otimes_{\scalebox{.5}{$\mathbb{Q}$}} \mathbb{R}
  $$

  \vspace{-1mm}
  \noindent
  form a real vector space with real dimension equal to the
  rational dimension of the corresponding rationalized homotopy groups

  \vspace{-.5cm}
  $$
    \mathrm{dim}_{\mathbb{Q}}
    \big(
      \pi_\bullet(X) \otimes_{\scalebox{.5}{$\mathbb{Z}$}} \mathbb{Q}
    \big)
    \;\; = \;\;
    \mathrm{dim}_{\mathbb{R}}
    \big(
      \pi_\bullet(X) \otimes_{\scalebox{.5}{$\mathbb{Z}$}} \mathbb{R}
    \big)
    \,,
  $$
  \vspace{-.4cm}

  \noindent
  and hence
  the rational Whitehead $L_\infty$-algebras (Prop. \ref{WhiteheadLInfinityAlgebras} below) have the same set
  of generators and their Chevalley-Eilenberg algebras
  (Def. \ref{ChevalleyEilenbergAlgebraOfLInfinityAlgebra})
  have the same structure constants, irrespective of
  whether they come as algebras over $\mathbb{Q}$
  or over $\mathbb{R}$.
  \footnote{
    While rational homotopy theory has the same form
    over all ground fields of characteristic zero,
    there is, of course, a difference between rational homotopy equivalences
    over different ground fields:
    Two minimal Sullivan models
    (Prop. \ref{ExistenceOfMinimalSullivanModels})
    over the real numbers may be isomorphic as real dgc-algebras
    but not as rational dgc-algebras.
    This happens when the isomorphism
    is given by an irrational linear transformation
    between the generators.
    For example, for any $b \in \mathbb{R}$, $b \geq 0$
    there is a dgc-algebra isomorphism over the real numbers

    \vspace{-.4cm}
    $$
      \xymatrix{
      \mathbb{R}
      \left[
        \!\!\!
        {\begin{array}{c}
          \omega_3
          \\
          \alpha_2, \beta_2
        \end{array}}
        \!\!\!
      \right]
      \!\big/\!
      \left(
        \!
        {\begin{aligned}
          d\, \omega_3 & =
            \alpha_2 \wedge \alpha_2
            +
            \beta_2 \wedge \beta_2
          \\[-2pt]
          d \, \alpha_2 & = 0\,,\; d\, \beta_2 \,=\, 0
        \end{aligned}}
        \!
      \right)
      \ar[rr]^-{
        \scalebox{.7}{$
          \arraycolsep=1.4pt
          {\begin{array}{rcr}
            \omega_3 & \mapsto & \omega_3
            \\[-3pt]
            \alpha_2 & \mapsto & \alpha_2
            \\[-2pt]
            \beta_2 & \mapsto & {\color{blue} \sqrt{b}}\, \beta_2
          \end{array}}
        $}
      }_-{ \simeq }
      &&
      \mathbb{R}
      \left[
        \!\!\!
        {\begin{array}{c}
          \omega_3
          \\
          \alpha_2, \beta_2
        \end{array}}
        \!\!\!
      \right]
      \!\big/\!
      \left(
        \!
        {\begin{aligned}
          d\, \omega_3
            & =
          \alpha_2 \wedge \alpha_2
            +
          {\color{blue} b} \,
          \beta_2 \wedge \beta_2
          \\[-2pt]
          d \, \alpha_2 & = 0\,,\; d\, \beta_2 \,=\, 0
        \end{aligned}}
        \!
      \right)
      \,.
      }
    $$
    \vspace{-.3cm}

    \noindent
    But over the rational numbers this exists only when the
    square root of $b$ is rational.

    Notice that the comparison between the homomotopy types over,
    in this order,
    the integers, the rational numbers and then the real numbers
    is provided by the character map
    (Def. \ref{DifferentialNonabelianCharacterMap} below);
    and the theory which
    embodies the distinction between working over these
    coefficients is that of homotopy fiber products of the
    character map, which is the theory of
    non-abelian differential cohomology
    (Def. \ref{DifferentialNonAbelianCohomology} below),
    where for instance the homotopy fiber
    $\xymatrix@C=10pt{
      \mathbb{R}/\mathbb{Q}
      \ar[r]
      &
      B \mathbb{Q}
      \ar[r]
      &
      B \mathbb{R} }$
    is being detected
    (e.g. \cite{GS-tAHSS}\cite{GS-RR}).
  }
  Therefore, we regard the case $k \,=\, \mathbb{R}$ as our default and
  abbreviate the PL de Rham Quillen adjunction
  (Prop. \ref{QuillenAdjunctionBetweendgcAlgebrasAndSimplicialSets})
  in this case by:
  \begin{equation}
    \label{RealPLdeRhamQuillenAdjunction}
    \Omega^\bullet_{\mathrm{PLdR}}
    \,\dashv_{{}_{\mathrm{Qu}}}\,
    \Bexp_{\mathrm{PL}}
    \;\;\;\;\;\;\;
    \coloneqq
    \;\;\;\;\;\;\;
    \Omega^\bullet_{\mathrm{P}{\scalebox{.6}{\color{blue}$\mathbb{R}$}}\mathrm{LdR}}
    \,\dashv_{{}_{\mathrm{Qu}}}\,
    \Bexp_{\mathrm{P}\scalebox{.6}{\color{blue}$\mathbb{R}$}\mathrm{L}}
    \,.
  \end{equation}
\end{remark}

\medskip

\noindent {\bf PS de Rham theory.} The point of
using piecewise \emph{polynomial} differential forms in
the PL de Rham complex (Def. \ref{PLdeRhamComplexes})
is that these, but not the piecewise smooth differential forms,
can be defined over the field $\mathbb{Q}$ of rational numbers.
But since
we may and do use the real numbers as the rational ground field
(Remark \ref{RationalHomotopyTheoryOverTheRealNumbers}),
it is expedient to also consider piecewise smooth de Rham complexes:

\begin{defn}[PS de Rham complex {\cite[p. 91]{GriffithMorgan13}}]
  \label{PSDeRhamComplex}
  For $n \in \mathbb{N}$, we write, in variation of \eqref{pdR},
  \vspace{-1mm}
  $$
    \Omega^\bullet_{\rm dR}
    \big(
      \mathbb{R}^n \times \Delta^{(-)}
    \big)
    \;:\;
    \xymatrix{
      \Delta^{\mathrm{op}}
      \ar[r]
      &
      \dgcAlgebras{\mathbb{R}}
    }
  $$

  \vspace{-2mm}
  \noindent
  for the simplicial dgc-algebra of smooth differential forms
  on the product manifold of
  $n$-dimensional Cartesian space with the standard simplices
  (i.e., of smooth differential forms
  on an ambient Cartesian space
  (Example \ref{SmoothdeRhamComplex}),
  restricted to the simplex and identified there if they agree on some open neighbourhood).
  As in Def. \ref{PLdeRhamComplexes},
  this induces for each $S \in \SimplicialSets$
  the corresponding \emph{piecewise smooth de Rham complexes}
  \begin{equation}
    \label{PiecewiseSmoothDeRhamComplexes}
    \Omega^\bullet_{\mathrm{PSdR}}(\mathbb{R}^n \times S)
    \;\;
    :=
    \;\;
    \underset{
      [k] \in \Delta
    }{\int}
      \underset{S_n}{\prod}
      \,
      \Omega^\bullet_{\mathrm{dR}}(\mathbb{R}^n \times \Delta^n)
  \end{equation}
  and by pullback of differential forms these extend to functors
  \vspace{-2mm}
  \begin{equation}
    \label{PiecewiseSmoothDeRhamComplexFunctor}
    \xymatrix{
      \SimplicialSets
      \ar[rrr]^-{
        \Omega^\bullet_{\mathrm{PSdR}}(\mathbb{R}^n \times (-))
      }
      &&&
      \dgcAlgebrasOp{\mathbb{R}}
      \,.
    }
  \end{equation}
\end{defn}

\begin{prop}[Fundamental theorem for piecewise smooth de Rham complexes]
  \label{FundamentalTheoremForPiecewiseSmoothDeRhamComplex}
  For all $n \in \mathbb{N}$
  the piecewise smooth de Rham complex functors
  (Def. \ref{PSDeRhamComplex})
  participate in a Quillen adjunction analogous to
  the PL de Rham adjunction (Prop. \ref{QuillenAdjunctionBetweendgcAlgebrasAndSimplicialSets})
  over the real numbers

  \vspace{-.4cm}
  \begin{equation}
    \label{PSdRQuillenAdjunctBetweendgcAlgsAndSimplicialSets}
    \xymatrix@C=4em{
      \dgcAlgebrasOpProj{\mathbb{R}}
      \;\;
      \ar@{<-}@<+7pt>[rr]^-{
        \Omega^\bullet_{\mathrm{PSdR}}(\mathbb{R}^n \times (-))
      }
      \ar@<-7pt>[rr]^-{\bot_{\mathrlap{\mathrm{Qu}}}}_-{
        \Bexp_{\mathrm{PS},n}
      }
      &&
      \;\;
      \SimplicialSets_{\mathrm{Qu}}
    }
  \end{equation}

  \vspace{-2mm}
  \noindent with right adjoint given as in \eqref{exp}:
  \vspace{-2mm}
  \begin{equation}
    \label{PSexp}
    \Bexp_{\mathrm{PS},n}(A)
    \;=\;
    \Big(
      \Delta[k]
        \;\longmapsto\;
      \dgcAlgebras{\mathbb{R}}
      \big(
        \Omega^\bullet_{\mathrm{PLdR}}(\mathbb{R}^n \times \Delta^k)
        \,,\,
        A
      \big)
    \Big)
    \;\;\;\;\;
    \in
    \;
    \SimplicialSets
    \,,
  \end{equation}

  \vspace{-2mm}
  \noindent
  \noindent
  whose derived functors (Prop. \ref{DerivedFunctors})
  are naturally equivalent to those of the PL de Rham adjunction
  \eqref{PLdeRhamDerivedAdjunction} over the real numbers:
  \vspace{-2mm}
  \begin{equation}
    \label{NaturalEquivalenceBetweenPLAndPSLeftDerivedFunctors}
    \LeftDerived\Omega^\bullet_{\mathrm{PSdR}}(\mathbb{R}^n \times (-))
    \;\;\simeq\;\;
    \LeftDerived\Omega^\bullet_{\mathrm{PSdR}}(-)
    \;\;\simeq\;\;
    \LeftDerived\Omega^\bullet_{\mathrm{P}\mathbb{R}\mathrm{LdR}}(-)
    \;,
  \end{equation}
  \begin{equation}
    \label{NaturalEquivalenceBetweenPLAndPSRightDerivedFunctors}
    \RightDerived\Bexp_{\mathrm{PS},n}
    \;\;\simeq\;\;
    \RightDerived\Bexp_{\mathrm{PS}}
    \;\;\simeq\;\;
    \RightDerived\Bexp_{\mathrm{P}\mathbb{R}\mathrm{L}}
    \,.
  \end{equation}
\end{prop}

\newpage

\begin{proof}
  {\bf (i)} The proofs of the
  PL de Rham theorem (Prop. \ref{PLdeRhamTheorem})
  as well as of the extension Lemma (Lemma \ref{ExtensionLemmaForPolynomialDifferentialForms})
  apply essentially verbatim also to piecewise-smooth
  differential forms (\cite[Prop. 9.8]{GriffithMorgan13})
  and hence so does the proof of the PL de Rham Quillen adjunction
  in the form given in Prop. \ref{QuillenAdjunctionBetweendgcAlgebrasAndSimplicialSets}.

  \noindent
  {\bf (ii)} We have evident natural transformations
  \vspace{-2mm}
  $$
    \xymatrix@C=3em{
      \Omega^\bullet_{\mathrm{P}\mathbb{R}\mathrm{LdR}}(S)
      \ar[r]_-{ \in\, \mathrm{W} }
      &
      \Omega^\bullet_{\mathrm{PSdR}}(S)
      \ar[r]_-{ \in\, \mathrm{W} }
      &
      \Omega^\bullet_{\mathrm{PSdR}}(\mathbb{R}^n \times S)
    }
    \,,
    \;\;\;\;\;\;\;\;\;\;\;
    \mbox{for $S \,\in\, \SimplicialSets$}
    \,,
  $$

  \vspace{-2mm}
\noindent
  given by inclusion of polynomial differential forms into
  smooth differential forms, and then by pullback of differential forms
  along the projections
  $\!\!\xymatrix@C=12pt{\mathbb{R}^n \times \Delta^k
  \ar[r] & \Delta^k}\!\!$. The corresponding
  component morphisms are quasi-isomorphisms (\cite[Cor. 9.9]{GriffithMorgan13}),
  hence are
  weak equivalences in $\dgcAlgebrasProj{\mathbb{R}}$ (Def. \ref{HomotopicalStructureOndgcAlgebras}).
  Under passage to homotopy categories (Def. \ref{HomotopyCategory})
  and derived functors (Example \ref{DerivedFunctorsByCoFibrantReplacement}),
  these natural weak equivalences become the natural isomorphisms
  \eqref{NaturalEquivalenceBetweenPLAndPSLeftDerivedFunctors}
  (by Prop. \ref{HomotopyCategoryIsLocalization}).
  By essential uniqueness of adjoint functors,
  this implies the natural isomorphisms \eqref{NaturalEquivalenceBetweenPLAndPSRightDerivedFunctors}.
\end{proof}

\medskip

\noindent {\bf Whitehead $L_\infty$-algebras.}

\begin{prop}[Real Whitehead $L_\infty$-algebras]
  \label{WhiteheadLInfinityAlgebras}
  For $X \in \NilpotentConnectedQFiniteHomotopyTypes$
  (Def. \ref{NilpotentConnectedSpacesOfFiniteRationalType}),
  there exists a nilpotent $L_\infty$-algebra
  (Def. \ref{NilpotentLInfinityAlgebras})
   \vspace{-2mm}
  \begin{equation}
    \label{WhiteheadLInfinityAlgebra}
    \mathfrak{l}X
    \;\in\;
    \LInfinityAlgebrasNil
    \,,
  \end{equation}

  \vspace{-1mm}
  \noindent
  unique up to isomorphism,
  whose Chevalley-Eilenberg algebra (Def. \ref{ChevalleyEilenbergAlgebraOfLInfinityAlgebra})
  is the minimal model (Def. \ref{MinimalSullivanModels})
  of the PL de Rham complex of $X$ (Def. \ref{PLdeRhamComplexes}):
  \vspace{-2mm}
  \begin{equation}
    \label{CEAlgebraOfWhiteheadLInfinityAlgebra}
    \mathrm{CE}
    (
      \mathfrak{l}X
    )
    \;:=\;
    \xymatrix{
      \big(
        \Omega^\bullet_{\mathrm{PLdR}}(X)
      \big)_{\mathrm{min}}
      \ar[r]_-{ \in \, \mathrm{W} }^-{
        p^{\mathrm{min}}_X
      }
      &
      \Omega^\bullet_{\mathrm{PLdR}}(X)
    }
    \,.
  \end{equation}
\end{prop}
\begin{proof}
By the PL de Rham theorem (Prop. \ref{PLdeRhamTheorem})
and the assumption that $X$ is connected, it follows that we have
$
  H\Omega^0_{\mathrm{PLdR}}(X)
  \;=\;
  \mathbb{R}
  \,.
$
Therefore Prop. \ref{ExistenceOfMinimalSullivanModels}
applies and says that
$\big( \Omega^\bullet_{\mathrm{PLdR}}(X) \big)_{\mathrm{min}}
 \in \SullivanModelsConnected
$
exists, and is unique up to isomorphism.
With this, the equivalence \eqref{NilpotentLInfinityAlgebrasEquivalentToConnectedSullivanModels}
says that $\mathfrak{l}X$ exists and is unique up to isomorphism.
\end{proof}
\begin{prop}[$\mathbb{R}$-Rationalization as integration of Whitehead $L_\infty$-algebras]
  \label{RationalizationInTermsOfWhiteheadLInfinityAlgebras}
  For $X \in \NilpotentConnectedQFiniteHomotopyTypes$
  (Def. \ref{NilpotentConnectedSpacesOfFiniteRationalType})
  its rationalization over the real numbers (Def. \ref{Lk})
  is equivalently the image under $\Bexp_{\mathrm{PL}}$
  \eqref{RealPLdeRhamQuillenAdjunction} of
  the CE-algebra \eqref{CEAlgebraOfWhiteheadLInfinityAlgebra}
  of its Whitehead $L_\infty$-algebra \eqref{WhiteheadLInfinityAlgebra}:
  \vspace{-.2cm}
  \begin{equation}
    \label{}
    L_{\mathbb{R}}
    (X)
    \;\simeq\;
    \Bexp_{\mathrm{PL}}
    \big(
      \mathrm{CE}(\mathfrak{l}X)
    \big)
    \;\;\;
    \in
    \;
    \mathrm{Ho}
    \big(
      \SimplicialSets_{\mathrm{Qu}}
    \big)
    \,.
  \end{equation}
\end{prop}
\begin{proof}
  By Def. \ref{Lk} and
  the characterization of derived functors (Ex. \ref{DerivedFunctorsByCoFibrantReplacement}),
  $L_{\mathbb{R}}$ is equivalently the image under $\Bexp_{\mathrm{PL}}$
  of any cofibrant replacement of
  $\Omega^\bullet_{\mathrm{PLdR}}(X) \,\in\, \dgcAlgebrasProj{\mathbb{R}}$
  (using that every $X \in \SimplicialSets_{\mathrm{Qu}}$ is already cofibrant \eqref{EverySimplicialSetIsCofibrantInClassicalModelStructure}).
  This is provided by $\mathrm{CE}(\mathfrak{l}X)$,
  according to \eqref{CEAlgebraOfWhiteheadLInfinityAlgebra}
  and by Prop. \ref{RelativeSullivanModelsAreCofibrations}.
\end{proof}

\begin{prop}[Rational homotopy groups in the rational Whitehead $L_\infty$ algebra]
  \label{RationalHomotopyGroupsInRationalWhiteheadLInfinityAlgebra}
  $\,$

  \noindent Let $X \in \NilpotentConnectedQFiniteHomotopyTypes$
  (Def. \ref{NilpotentConnectedSpacesOfFiniteRationalType}).

  \noindent {\bf (i)} If $X$ is simply connected, $\pi_1(X) = 1$
  (Example \ref{ExamplesOfNilpotentSpaces}),
  then there is an isomorphism of graded vector spaces (Def. \ref{CategoryOfGradedVectorSpaces})
  between the graded vector space underlying
  \eqref{GradedVectorSpaceUnderlyingLInfinityAlgebra} the
  Whitehead $L_\infty$-algebra $\mathfrak{l}X$
  (Prop. \ref{WhiteheadLInfinityAlgebras})
  and
  the rationalized homotopy groups of
  the based loop space $\Omega X$:
    \vspace{-3mm}
  $$
    \xymatrix{
      \overset{
        \mathclap{
        \raisebox{6pt}{
          \tiny
          \color{darkblue}
          \bf
          \def\arraystretch{.9}
          \begin{tabular}{c}
            Whitehead
            \\
            $L_\infty$-algebra
          \end{tabular}
        }
        }
      }{
        \mathfrak{l}X
      }
      \;\;\;\;\simeq\;\;\;\;
      \overset{
        \mathclap{
        \raisebox{6pt}{
          \tiny
          \color{darkblue}
          \bf
          \def\arraystretch{.9}
          \begin{tabular}{c}
            rationalized
            \\
            homotopy groups
          \end{tabular}
        }
        }
      }{
\       \pi_{\bullet}(\Omega X) \otimes_{\mathbb{Z}} \mathbb{R}
      }
      \;\;\;\;\;
      \in
      \GradedVectorSpaces
      \,.
    }
  $$

    \vspace{-3mm}
  \noindent
  {\bf (ii)} If $\pi_1(X)$ is not necessarily trivial but abelian,
  then this statement
  still holds with $\mathfrak{l}X$ replaced by its homology
  with respect to the unary differential $[-]$ \eqref{LInfinityBracketsEncodedInCEDifferential}.

  \noindent
  {\bf (iii)} If $\pi_1(X)$ is not abelian, then {(ii)}
  still holds in degrees $\geq 2$.
\end{prop}
\begin{proof}
  Under translation through Prop. \ref{WhiteheadLInfinityAlgebras}
  and Def. \ref{ChevalleyEilenbergAlgebraOfLInfinityAlgebra},
  and using
  $\pi_\bullet(\Omega X) \simeq \pi_{\bullet+1}(X)$,
   claim {\bf (i)} is equivalent to the existence of a dual isomorphism:
     \vspace{-2mm}
  \begin{equation}
    \label{DualIsomorphismFromSullivanGeneratorsToRationalHomotopyGroups}
    \mathrm{CE}
    (
      \mathfrak{l}X
    )_{
      \!\!\big/
      \mathrm{CE}
      (
      \mathfrak{l}X
      )^2
    }
    \;\;\;\simeq\;\;\;
    \mathrm{Hom}_{\mathbb{Z}}
    \big(
      \pi_\bullet(X)
      ,\,
      \mathbb{R}
    \big)
    \;\;\;\;\;
    \in
    \GradedVectorSpaces
    \,,
  \end{equation}

    \vspace{-2mm}
\noindent  where the left hand side denotes the graded vector
  space of indecomposable elements in the Chevalley-Eilenberg algebra
  (the $\alpha_{n_i}^{(i)}$ in \eqref{MultivariatePolynomialdgcAlgebra}).
  In this form, this is the statement of
  \cite[Theorem 11.3 with Def. 6.12]{BousfieldGugenheim76},
  in the special case when,
  with $\pi_1(X) = 1$, the unary differential $[-]$ in
  $\mathfrak{l}X$ vanishes (Example \ref{MinimalModelOFSimplyConnecteddgcAlgebras}).
  The generalizations follow analogously.
\end{proof}

\begin{remark}[Equivalent $L_\infty$-structures on Whitehead products]
  \label{EquivalentLInfinityStructuresOnWhiteheadProducts}
  The original discussion of the Whitehead algebra
  structure on the homotopy groups of a space
  is in terms of differential-graded Lie algebras
  (\cite[Theorem B]{Hilton55}), as are the Quillen models
  of rational homotopy theory \cite{Quillen69}.

  \noindent {\bf (i)} Notice that
  dg-Lie algebras (Example \ref{DifferentialGradedLieAlgebras})
  and $L_\infty$-algebras with minimal CE-algebra (Def. \ref{MinimalSullivanModels}) are two opposite classes
  of $L_\infty$-algebras: The former has
  $k$-ary brackets \eqref{LInfinityBracketsEncodedInCEDifferential} only for $k \leq 2$, the latter only for $k \geq 2$
  (in the simply connected case, by Example \ref{MinimalModelOFSimplyConnecteddgcAlgebras}).
  Yet, quasi-isomorphisms connect algebras in one class to
  those in the other (\cite[p. 28]{KrizMay95}),
  such as to make their homotopy theories equivalent
  (\cite{Pridham10}, see also Rem. \ref{HomotopyTheoryOfLInfinityAlgebras}).
  The transmutation of dg-Lie- into minimal $L_\infty$-algebras
  is described in \cite[Thm. 2.1]{BBMM16};
  that from $L_\infty$- to dg-Lie-algebras in \cite[\S 1.0.2]{FRS13}.

    \noindent {\bf (ii)} The minimal $L_\infty$-algebra structure on $\mathfrak{l}X$
  that we obtained in Prop. \ref{WhiteheadLInfinityAlgebras},
  \ref{RationalHomotopyGroupsInRationalWhiteheadLInfinityAlgebra},
  has the property that its
  $k$-ary brackets are, up to possibly a sign, equal to the
  order-$k$ higher Whitehead products on $X$
  \cite[Prop. 3.1]{BBMM16}.
\end{remark}

\medskip

\noindent
{\bf Examples of rationalizations over the real numbers.}
The following fundamental examples of rationalizations
serve to illustrate the above notation and terminology
and to highlight that rationalization {\it over the real numbers}
(Def. \ref{Lk}), even though it is not a localization
(Rem. \ref{FailureOfRationalizationOverRToBeIdempotent} below),
still acts as real-ification on the homotopy groups of Eilenberg-MacLane
spaces (Ex. \ref{RationalizationOfEMSpaces} below) and, more generally,
of loop spaces (Ex. \ref{WhiteheadLInfinityAlgebraOfLoopSpaces} below).
This is the crucial fact that makes the real character map on non-abelian cohomology
in \cref{ChernCharacterInNonabelianCohomology} reduce to the traditional Chern-Dold characters
on abelian generalized cohomology in \cref{TheChernDoldCharacter}.

\begin{example}[$\mathbb{R}$-Rationalization of $n$-spheres {(e.g. \cite[\S 1.2]{LM13})}]
  \label{RationalizationOfnSpheres}
  The Serre finiteness theorem (see \cite[Thm. 1.1.8]{Ravenel86})
  says that the homotopy groups of $n$-spheres
  for $n \geq 1$ are of the form
  \vspace{-1mm}
  \begin{equation}
    \label{HomotopyGroupsOfSpheres}
    \pi_{n+k}
    \big(
      S^n
    \big)
    \;\simeq\;
    \left\{
    \begin{array}{lcl}
      \mathbb{Z} & \vert & k = 0
      \\
      \mathbb{Z} \oplus \mathrm{fin} &\vert&
      k = 2 m \;\; \mbox{and}\;\; n = 2m -1
      \\
      \mathrm{fin} &\vert& \mathrm{otherwise}
    \end{array}
    \right.
  \end{equation}

  \vspace{-1mm}
  \noindent
  where ``$\mathrm{fin}$'' stands for some finite group.
  Since finite groups are pure torsion, hence
  have trivial rationalization, this means that the rational
  homotopy groups of spheres are:
  \vspace{-1mm}
  $$
    \pi_{n+k}
    \big(
      S^n
    \big)
    \otimes_{\scalebox{.5}{$\mathbb{Z}$}}
    \mathbb{R}
    \;\simeq\;
    \left\{
    \begin{array}{lcl}
      \mathbb{R} & \vert & k = 0
      \\
      \mathbb{R} &\vert&
      k = 2 m \;\,\mbox{and}\;\, n = 2m -1
      \\
      0 &\vert& \mathrm{otherwise}\,.
    \end{array}
    \right.
  $$

  \vspace{-1mm}
  \noindent
  Moreover, the fact that ordinary cohomology is represented
  by Eilenberg-MacLane spaces (Example \ref{OrdinaryCohomology})
  means that
  \begin{equation}
    \label{RealCohomologyOfASphere}
    H^k
    \big(
      S^n;
      \,
      \mathbb{R}
    \big)
    \;\;\simeq\;\;
    \left\{
    \begin{array}{lcl}
      \mathbb{R} &\vert& k \in \{0,n\}
      \\
      0 &\vert& \mbox{otherwise}.
    \end{array}
    \right.
  \end{equation}
  With this, Prop. \ref{RationalHomotopyGroupsInRationalWhiteheadLInfinityAlgebra}
  together with Prop. \ref{PLdeRhamTheorem} implies
  that the Whitehead $L_\infty$-algebras of spheres
  (Prop. \ref{WhiteheadLInfinityAlgebras})
  are as follows:
  \begin{equation}
    \label{SullivanModelForOddDimensionalSpheres}
    \mathrm{CE}
    \big(
      \mathfrak{l}S^n
    \big)
    \;\;\simeq\;\;
    \mathbb{R}
    \big[
      \omega_n
    \big]
    \big/
    \big(
      d\, \omega_n \,=\,0.
    \big)
    \phantom{AAAAAAAAAA}
    \mbox{if $n$ is odd}
  \end{equation}
  and
  \begin{equation}
    \label{SullivanModelForEvenDimensionalSpheres}
    \mathrm{CE}
    \big(
      \mathfrak{l}S^n
    \big)
    \;\;\simeq\;\;
    \mathbb{R}
    \left[
      \!\!\!
      \begin{array}{c}
        \omega_{2n-1}
        \\
        \omega_n
      \end{array}
      \!\!\!
    \right]
    \!\big/\!
    \left(
      \!
      \begin{aligned}
        d\, \omega_{2n-1} & = - \omega_n \wedge \omega_n
        \\[-3pt]
        d\, \omega_n \;\;\;\;\;\; & = 0
      \end{aligned}
      \!
    \right)
    \phantom{AAAA}
    \mbox{if $n > 1$ is even}
  \end{equation}
\end{example}

\begin{example}[$\mathbb{R}$-Rationalization of Eilenberg-MacLane spaces]
  \label{RationalizationOfEMSpaces}
  Since the homotopy types of
  Eilenberg-MacLane-spaces $K(A,n+1) = B^{n+1} A$  (see \eqref{EilenbergMacLaneSpaces})
  are fully characterized
  by their homotopy groups   (for discrete abelian groups $A$, e.g. \cite[\S 6]{AGP02}))
  \vspace{-1mm}
  $$
    \pi_k
    \big(
      B^{n+1} A
    \big)
    \;\simeq\;
    \left\{
    \begin{array}{lcl}
      A &\vert& k = n+1
      \\
      0 &\vert& k \neq n+1
    \end{array}
    \right.
  $$

  \vspace{-1mm}
  \noindent we have, for $n \in \mathbb{N}$:

  \noindent
  {\bf (i)}
  Their Whitehead $L_\infty$-algebra (Prop. \ref{WhiteheadLInfinityAlgebras})
  is,
  by Prop. \ref{RationalHomotopyGroupsInRationalWhiteheadLInfinityAlgebra},
  the direct sum of
  $\mathrm{dim}\big( A \otimes_{{}_{\mathbb{Z}}} \mathbb{R} \big)$
  copies of the line Lie $n$-algebra (Def. \ref{LineLienPlusOneAlgebras}):
  \vspace{-1mm}
  \begin{equation}
    \label{WhiteheadLInfinityAlgebraOfIntegralEMSpace}
    \mathfrak{l}
    \big(
      B^{n+1} A
    \big)
    \;\simeq\;
    \mathfrak{b}^n
    \big(
      A \otimes_{{}_{\mathbb{Z}}} \mathbb{R}
    \big)
    \;\simeq\;
    \underset{
      \mathrm{dim}
      (
        A \otimes_{{}_{\mathbb{Z}}} \mathbb{R}
      )
    }{\bigoplus}
    \mathfrak{b}^n \mathbb{R}
    \;.
  \end{equation}

  \vspace{-2mm}
  \noindent
  {\bf (ii)}
  Their rationalization over $\mathbb{R}$ (Def. \ref{Lk})
  is the Eilenberg-MacLane space on the realification of $A$:
  \begin{equation}
    \label{LineLienAlgebraExponentiation}
    \;\;\;\;\;
    L_{\mathbb{R}}
    \big(
      B^{n+1} A
    \big)
    \;\simeq\;
    B^{n+1}
    \big(
      A \otimes_{{}_{\mathbb{Z}}} \mathbb{R}
    \big)
    \;\;\;
    \in
    \;
    \mathrm{Ho}
    \big(
      \TopologicalSpaces_{\mathrm{Qu}}
    \big).
  \end{equation}
  Observe how this is implied
  via the machinery that we have set up above: For all $k \in \mathbb{N}_+$ we have:

  \vspace{-.4cm}
  $$
    \def\arraystretch{1.8}
    \begin{array}{lll}
      \pi_k
      \Big(
        L_{\mathbb{R}}
        \big(
          B^{n+1} A
        \big)
      \Big)
      &
      \;=\;
      \mathrm{Ho}
      \big(
        \TopologicalSpaces_{\mathrm{Qu}}^{\ast/}
      \big)
      \Big(
        S^k
        ,\,
        L_{\mathbb{R}}
        \big(
          B^{n+1} A
        \big)
      \Big)
      &
      \mbox{by def. of $\pi_k(-)$}
      \\
      & \;=\;
      \mathrm{Ho}
      \big(
        \TopologicalSpaces_{\mathrm{Qu}}^{\ast/}
      \big)
      \Big(
        S^k
        ,\,
        \Bexp_{\mathrm{PL}}
        \big(
          \mathrm{CE}( \mathfrak{b}^n ( A \otimes_{{}_{\mathbb{Z}}} \mathbb{R} ) )
        \big)
      \Big)
      &
      \mbox{by Prop. \ref{RationalizationInTermsOfWhiteheadLInfinityAlgebras} with \eqref{WhiteheadLInfinityAlgebraOfIntegralEMSpace}}
      \\
      & \;\simeq\;
      \mathrm{Ho}
      \Big(
        \dgcAlgebrasProj{\mathbb{R}}^{/\mathbb{R}}
      \Big)
      \Big(
        \mathrm{CE}
        \big(
          \mathfrak{b}^n (A \otimes_{{}_{\mathbb{Z}}} \mathbb{R} )
        \big)
        ,\,
        \Omega^{\bullet}_{\mathrm{PLdR}}
        \big(
          S^k
        \big)
        \!\!\!
        \underset{
          \Omega^\bullet_{\mathrm{PLdR}}(\ast)
        }{ \times }
        \!\!\!
        \mathbb{R}
      \Big)
      &
      \mbox{by Props. \ref{QuillenAdjunctionBetweendgcAlgebrasAndSimplicialSets},
       \ref{InducedQuillenAdjunctionOnPointedObjects},
       \ref{DerivedFunctors}}
      \\
      & \;\simeq\;
      \underset{
        \mathrm{dim}
        (
          A \otimes_{{}_{\mathbb{Z}}} \mathbb{R}
        )
      }{\prod}
      \!\!\!\!\!\!\!
      \mathrm{Ho}
      \big(
        \dgcAlgebrasProj{\mathbb{R}}^{/\mathbb{R}}
      \big)
      \Big(
        \mathrm{CE}
        \big(
          \mathfrak{b}^n \mathbb{R}
        \big)
        ,\,
        \Omega^{\bullet}_{\mathrm{PLdR}}
        \big(
          S^k
        \big)
        \!\!\!
        \underset{
          \Omega^\bullet_{\mathrm{PLdR}}(\ast)
        }{ \times }
        \!\!\!
        \mathbb{R}
      \Big)
      &
      \mbox{by \eqref{VectorSpaceLienPlus1Algebra}}
      \\
      & \;\simeq\;
      \underset{
        \mathrm{dim}
        (
          A \otimes_{{}_{\mathbb{Z}}} \mathbb{R}
        )
      }{\prod}
      \!\!\!\!\!\!\!
      H^{n+1}\big( S^k;\, \mathbb{R} \big)
      &
      \mbox{by Lem. \ref{HomotopicalFormulationOfCochainCohomology}}
      \\
      & \;\simeq\;
      \left\{
      \!\!\!
      \def\arraystretch{1}
      \begin{array}{ccl}
        A \otimes_{{}_{\mathbb{Z}}} \mathbb{R}
        &\vert& k = n + 1
        \\
        0 &\vert& k \neq n + 1
      \end{array}
      \right.
      &
      \mbox{by \eqref{RealCohomologyOfASphere}}
      \,.
    \end{array}
  $$
  The same computation, but with $S^k$ replaced by the point $\ast$
  and without the slicing, shows that
  $\pi_0\left(L_{\mathbb{R}}\big(B^{n+1}A\big)\right) = \ast$.
\end{example}
\begin{remark}[Failure of $\mathbb{R}$-rationalization to be idempotent]
  \label{FailureOfRationalizationOverRToBeIdempotent}
  Example \ref{RationalizationOfEMSpaces} reveals how
  rationalization over $\mathbb{R}$ (or over any field $k$ strictly containg
  the rational numbers, Def. \ref{Lk}) is not idempotent, hence not a
  localization (see also \cite[Rem. 9.7]{BousfieldGugenheim76}):
  Applying \eqref{LineLienAlgebraExponentiation} twice yields

  \vspace{-.4cm}
  \begin{equation}
    \label{DoubleRealificationOfEilenbergMacLaneSpace}
    L_{\mathbb{R}}
    \circ
    L_{\mathbb{R}}
    \big(
      B^{n+1} A
    \big)
    \;\;
    \simeq
    \;\;
    B^{n+1}
    \big(
      A \otimes_{\mathbb{Z}} \mathbb{R} \otimes_{\mathbb{Z}} \mathbb{R}
    \big)
    \,,
  \end{equation}
  \vspace{-.4cm}

  \noindent
  but $\mathbb{R} \otimes_{{}_{\mathbb{Z}}} \mathbb{R}  \,\simeq\!\!\!\!\!\!\!\!/\;\;\, \mathbb{R}$,
  in contrast to $\mathbb{Q} \otimes_{{}_{\mathbb{Z}}} \mathbb{Q} \,\simeq\, \mathbb{Q}$
  (reflecting that $\mathbb{Q}$ is a {\it solid ring} \cite[\S 2.4]{BousfieldKan72Core}, while $\mathbb{R}$ is not).
\end{remark}

\begin{example}[$\mathbb{R}$-Rationalization of complex projective spaces]
  \label{RationalizationOfComplexProjectiveSpace}
  From the defining homotopy fiber sequence
  for complex projective $n$-space $\mathbb{C}P^n$, $n \in \mathbb{N}$
  (e.g. \cite[Ex. 14.22]{BottTu82})

  \vspace{-.3cm}
  \begin{equation}
    \begin{tikzcd}[row sep=4pt]
      S^1
      \ar[
        r,
        hook
      ]
      \ar[
        d,-,
        shift left=1pt
      ]
      \ar[
        d,-,
        shift right=1pt
      ]
      &
      S^{2n+1}
      \ar[
        r,
        ->>
      ]
      \ar[
        d,-,
        shift left=1pt
      ]
      \ar[
        d,-,
        shift right=1pt
      ]
      &
      \mathbb{C}P^n
      \ar[
        d,-,
        shift left=1pt
      ]
      \ar[
        d,-,
        shift right=1pt
      ]
      \\
      \mathrm{U}(1)
      \ar[
        r,
        hook
      ]
      &
      \mathrm{U}(n+1)/\mathrm{U}(n)
      \ar[
        r,
        ->>
      ]
      &
      \mathrm{SU}(n+1)/\mathrm{U}(n)
    \end{tikzcd}
  \end{equation}
  \vspace{-.3cm}

  \noindent
  the corresponding long exact sequence of homotopy groups (e.g. \cite[Thm. 6.1.2]{tomDieck08})
  yields the following homotopy groups:

  \vspace{-.5cm}
  \begin{equation}
    \label{HomotopyGroupsOfComplexProjectiveSpace}
    \pi_k
    \big(
      \mathbb{C}P^n
    \big)
    \;=\;
    \left\{
    \begin{array}{lll}
      \ast & \vert & k \in \{0,1\}
      \\[-2pt]
      \mathbb{Z} &\vert&  k \in \{2, 2 n + 1\}
      \\[-2pt]
      \pi_k\big(S^{2 n + 1}\big) &\vert & k \geq 2n + 1
      \\[-2pt]
      0 &\vert& \mbox{otherwise}
      \,.
    \end{array}
    \right.
    \,,
    {\phantom{AAA}}
    H^k
    \big(
      \mathbb{C}P^n;
      \,
      \mathbb{R}
    \big)
    \;\simeq\;
    \left\{
    \begin{array}{lll}
      \mathbb{R} &\vert& k \in \{0,1, 2, \cdots, 2n\}
      \\
      0 &\vert& \mbox{otherwise.}
    \end{array}
    \right.
  \end{equation}
  \vspace{-.3cm}

  \noindent
  Since these homotopy groups \eqref{HomotopyGroupsOfComplexProjectiveSpace}
  in degrees $\geq 2 n + 1$ are finite,
  hence rationally trivial, by Ex. \ref{RationalizationOfnSpheres},
  it follows,
  with Prop. \ref{RationalHomotopyGroupsInRationalWhiteheadLInfinityAlgebra}
  and Prop. \ref{PLdeRhamTheorem},
  that the Chevalley-Eilenberg algebra of the Whitehead $L_\infty$-algebra of $\mathbb{C}P^n$
  (Prop. \ref{WhiteheadLInfinityAlgebras})
  has exactly one generator $f_{2}$ in degree $2$ and one generator $h_{2k+1}$
  in degree $2k+1$.
  Moreover, since the cohomology groups are
  (e.g. \cite[Ex. 14.22.1]{BottTu82})
  as shown on the right of \eqref{HomotopyGroupsOfComplexProjectiveSpace},
  the first of these must be the closed generator of the cohomology ring,
  and the differential of the latter must exhibit the vanishing of its $(n+1)$st cup product
  in cohomology (see also, e.g., \cite[p. 203]{FHT00}\cite[\S 5.3]{LM13}):

  \vspace{-5.5mm}
  \begin{equation}
    \mathrm{CE}
    \big(
      \mathfrak{l}
      \mathbb{C}P^n
    \big)
    \;=\;
    \mathbb{R}
    \left[
      \!\!\!
      \def\arraystretch{.9}
      \begin{array}{l}
        h_{2n+1},
        \\
        f_2
      \end{array}
      \!\!\!
    \right]
    \big/
    \Big(
    \raisebox{7pt}{$
    \!\!\!
    \def\arraystretch{.9}
    \begin{array}{l}
      d \,h_{2n+1}\,
        =
      \overset{
        \mbox{\tiny \rm $n + 1$ factors}
      }{
      \overbrace{
        f_2 \wedge \cdots \wedge f_2
      }
      }
      \\
    \quad  \; d \, f_2 \;=\; 0
    \end{array}
    \!\!\!
    $}
    \Big)
    \mathrlap{\,.}
  \end{equation}
  \vspace{-.4cm}
\end{example}

\begin{example}[$\mathbb{R}$-Rationalization of loop spaces]
  \label{WhiteheadLInfinityAlgebraOfLoopSpaces}
  The minimal Sullivan model
  (Def. \ref{MinimalSullivanModels})
  of a loop space $A \,\simeq\, \Omega A'$ of $\mathbb{Q}$-finite type
  (which exists by Ex. \ref{ExamplesOfNilpotentSpaces})
  has vanishing differential (e.g. \cite[p. 143]{FHT00}).
  Therefore,
  Prop. \ref{RationalHomotopyGroupsInRationalWhiteheadLInfinityAlgebra}
  implies
  that the rational Whitehead $L_\infty$-algebra
  $\mathfrak{l}A$
  (Prop. \ref{WhiteheadLInfinityAlgebras}) of $A$ is
  the direct sum of line Lie $n$-algebras $\mathfrak{b}^n \mathbb{R}$
  (Example \ref{LineLienPlusOneAlgebras})
  and its Chevalley-Eilenberg algebra (Def. \ref{CEAlgebraOfLieAlgebra})
  is the tensor product of those of the summands:
  \vspace{-1.5mm}
  $$
    \mathfrak{l}A
    \;\simeq\;
    \underset{n \in \mathbb{N}}{\bigoplus}
    \,
    \mathfrak{b}^n
    \big(
      \pi_{n+1}(A) \otimes_{\scalebox{.5}{$\mathbb{Z}$}} \mathbb{R}
    \big)
    \;\;\;\;
    \in
    \;
    \LInfinityAlgebrasNil
    \,,
    \;\;\;\;\;\;\;\;\;
    \mathrm{CE}
    \big(
      \mathfrak{l}A
    \big)
    \;\simeq\;
    \underset{
      n \in \mathbb{N}
    }{\bigotimes}
    \,
    \mathrm{CE}
    \Big(
      \mathfrak{b}^n
      \big(
        \pi_{n+1}(A) \otimes_{\scalebox{.5}{$\mathbb{Z}$}} \mathbb{R}
      \big)
   \Big)
   \;\;\;\;\;
   \in
   \dgcAlgebras{\mathbb{R}}
   \,.
  $$

  \vspace{-2mm}
\noindent  Noticing that the tensor product of dgc-algebras is the
  coproduct in the category of $\dgcAlgebras{\mathbb{R}}$
  (Ex. \ref{ProductAndCoproductAlgebras}),
  and hence the Cartesian product in the opposite category,
  the right adjoint functor $\Bexp_{\mathrm{PL}}$ \eqref{QuillenAdjunctBetweendgcAlgsAndSimplicialSets}
  preserves this,
  so that
  \vspace{-1.5mm}
  $$
    \Bexp_{\mathrm{PL}}
    \,\circ\,
    \mathrm{CE}
    \big(
      \mathfrak{l}A
    \big)
    \;\;
    \simeq
    \;\;
    \underset{
      n \in \mathbb{N}
    }{ \prod }
    \Big(
      \Bexp_{\mathrm{PL}}
        \,\circ\,
      \mathrm{CE}
      \big(
        \mathfrak{b}^n
        (
          \pi_{n+1}(A) \otimes_{\scalebox{.5}{$\mathbb{Z}$}} \mathbb{R}
        )
      \big)
    \Big)
    \;\;\;\;\;
    \in
    \mathrm{Ho}
    \big(
      \TopologicalSpaces_{\mathrm{Qu}}
    \big)
    \,.
  $$

  \vspace{-1mm}
\noindent  Finally, by Ex. \ref{RationalizationOfEMSpaces}, this means that
  the rationalization over $\mathbb{R}$ (Def. \ref{Lk}) of a loop space
  is a Cartesian product of Eilenberg-MacLane spaces:
  \begin{equation}
    \label{RationalizationOverTheRealNumbersOfLoopSpace}
    A \,\simeq\, \Omega A'
    \;\;\;\;\;\;\;\;\;\;
    \Rightarrow
    \;\;\;\;\;\;\;\;\;\;
    L_{\mathbb{R}}(A)
    \;\;\simeq\;\;
    \underset{ n \in \mathbb{N} }{\prod}
    B^{n+1}
    \big(
      \pi_{n+1}(A)
        \otimes_{{}_{\mathbb{Z}}}
      \mathbb{R}
    \big)
    \,.
  \end{equation}
\end{example}

\begin{example}[$\mathbb{R}$-Rationalization of spectra]
  \label{RealRationalizationOfSpectra}
  For $E$ a spectrum (Ex. \ref{GeneralizedCohomologyAsNonabelianCohomology}),
  Ex. \ref{WhiteheadLInfinityAlgebraOfLoopSpaces}
  says
  that its degree-wise
  rationalization (Def. \ref{Rationalization})
  and $\mathbb{R}$-rationalization (Def. \ref{Lk})
  are both direct sums of the same form
  of rational Eilenberg-MacLane spectra:
  \begin{equation}
  \label{RationalizedSpectra}
  \begin{aligned}
    L_{\color{blue}\mathbb{Q}}(E)
    &
    \;\simeq\;
    \underset{ k \in \mathbb{Z} }{\oplus}
    H\big( \pi_\bullet(E) \otimes_{{}_{\mathbb{Z}}} {\color{blue}\mathbb{Q}} \big)
    \\
    L_{\color{blue}\mathbb{R}}(E)
    &
    \;\simeq\;
    \underset{ k \in \mathbb{Z} }{\oplus}
    H\big( \pi_\bullet(E) \otimes_{{}_{\mathbb{Z}}} {\color{blue}\mathbb{R}} \big)
    \,,
  \end{aligned}
  \end{equation}
  But rationalization
  of spectra is also known
  (review in \cite[Ex. 8.12]{Lawson20}\cite[Ex. 1.7 (4)]{Bauer14})
  to be
  given by forming the smash product of spectra (e.g. \cite{EKMM97})
  with the rational Eilenberg-MacLane spectrum:
  \begin{equation}
    \label{RationalizationOfSpectraIsSmashing}
    L_{\mathbb{Q}}(E_n)
    \;\simeq\;
    \big(
      E \wedge H \mathbb{Q}
    \big)_n
    \,.
  \end{equation}
  Observing with\eqref{RationalizedSpectra} that

  \vspace{-.5cm}
  \begin{equation}
    \label{RealRationalizationOfSpectraIsRationalizationSmashedOverHQWithHR}
    L_{\mathbb{R}}(E)
    \;\simeq\;
    \big(
      L_{\mathbb{Q}}(E)
    \big)
      \wedge_{{}_{H \mathbb{Q}}}
    H \mathbb{R}
    \,.
  \end{equation}
  this implies that the componentwise
  $\mathbb{R}$-rationalization (Def. \ref{Lk})
  of spectra is analogously given by the smash product
  with the real Eilenberg-MacLane:

  \vspace{-.4cm}
  \begin{equation}
    \label{RealificationOfSpectraIsSmashWithHR}
    \begin{array}{lll}
      L_{\mathbb{R}}(E)
      &
      \;\simeq\;
      \big(
        L_{\mathbb{Q}}(E)
      \big)
        \wedge_{{}_{H \mathbb{Q}}}
      H \mathbb{R}
      &
      \mbox{by \eqref{RationalizationOfSpectraIsSmashing}}
      \\
      & \;\simeq\;
      E
        \wedge
      H \mathbb{Q}
        \wedge_{{}_{H \mathbb{Q}}}
      H \mathbb{R}
      &
      \mbox{by \eqref{RealRationalizationOfSpectraIsRationalizationSmashedOverHQWithHR}}
      \\
      & \;\simeq\;
      E
        \wedge
      H \mathbb{R}
      \,.
    \end{array}
  \end{equation}
  \vspace{-.3cm}

  \noindent
  It is in this smashing form that $\mathbb{R}$-rationalization of spectra
  traditionally appears in the construction of
  differential generalized cohomology theories, see
  Ex. \ref{DifferentialGeneralizedCohomology} below.
\end{example}

\medskip

\noindent {\bf Non-abelian real cohomology.} Using these
$\mathbb{R}$-rational models, we obtain the first key concept
formation towards the character map:

\begin{defn}[Non-abelian real cohomology]
  \label{NonAbelianRealCohomology}
  Let $X, A \in \TopologicalSpaces$.
  Then the \emph{non-abelian real cohomology}
  of $X$ with coefficients in $A$
  is the non-abelian cohomology
  (Def. \ref{NonAbelianCohomology})
   of $X$ with
  coefficients in the $\mathbb{R}$-rationalization $L_{\mathbb{R}} A$
  (Def. \ref{Lk})
    \vspace{-2mm}
    \begin{equation}
    \label{NonAbelianRealCohomologyByHomotopyClassesOfMaps}
    H
    \big(
      X;
      \,
      L_{\mathbb{R}}A
    \big)
    \;:=\;
    \mathrm{Ho}
    \big(
      \TopologicalSpaces_{\mathrm{Qu}}
    \big)
    \big(
      X
      \,,\,
      L_{\mathbb{R}}A
    \big)
    \,.
  \end{equation}
\end{defn}
\begin{example}[Non-abelian real cohomology subsumes ordinary real cohomology]
  For $n \in \mathbb{N}$,
  non-abelian real cohomology (Def. \ref{NonAbelianRealCohomology})
  with coefficients in the
  $\mathbb{R}$-rationalized classifying space
  $
    L_{\mathbb{R}}
    \big(
      B^{n+1} \mathbb{Z}
    \big)
    \;\simeq\;
    B^{n+1} \mathbb{R}
  $
  (by Ex. \ref{RationalizationOfEMSpaces})
  is naturally equivalent
  (by Ex. \ref{OrdinaryCohomology})
  to ordinary real cohomology in degree $n$:
  \vspace{-1mm}
  $$
    H
    \big(
      X;
      \,
      L_{\mathbb{R}}
      B^{n+1} \mathbb{Z}
    \big)
    \;\simeq\;
    H
    \big(
      X;
      \,
      B^{n+1} \mathbb{R}
    \big)
    \;\simeq\;
    H^{n+1}
    (
      X;
      \,
      \mathbb{R}
    )\;.
  $$
\end{example}

More generally:
\begin{prop}[Non-abelian real cohomology with coefficients in loop spaces]
  \label{NonAbelianRealCohomologyWithCoefficientsInLoopSpaces}
  $\,$

  \noindent
  Let $A \in \NilpotentConnectedQFiniteHomotopyTypes$
  (Def. \ref{NilpotentConnectedSpacesOfFiniteRationalType})
  such that it admits loop space structure,
  hence such that there exists $A'$
  with
  \vspace{-1mm}
  $$
    A
    \;\simeq\;
    \Omega A'
    \;\;\;
    \in
    \;
    \mathrm{Ho}
    \big(
      \TopologicalSpaces_{\mathrm{Qu}}
    \big)
    \,.
  $$

  \vspace{-1mm}
  \noindent
  Then the non-abelian real cohomology (Def. \ref{NonAbelianRealCohomology})
  with coefficients in $L_{\mathbb{R}}A$ is
  naturally equivalent to ordinary real cohomology
  with coefficients in the rationalized homotopy groups of $A$:
  \vspace{-1mm}
  \begin{equation}
    \label{NonAbelianRealCohomologyWithCoefficientsInLoopSpace}
    H
    \big(
      X;
      \,
      L_{\mathbb{R}}A
    \big)
    \;\;
    \simeq
    \;\;
    \underset{n \in \mathbb{N}}{\bigoplus}
    \,
    H^{ n+1  }
    \big(
      X;
      \,
      \pi_{n+1}(A) \otimes_{\scalebox{.5}{$\mathbb{Z}$}}
      \mathbb{R}
    \big).
  \end{equation}
\end{prop}
\begin{proof}
  This follows from the following sequence of natural bijections:
  $$
    \begin{array}{lll}
      H
      \big(
        X;
        \,
        L_{\mathbb{R}} A
      \big)
      &
      \simeq
      H
      \Big(
        X;
        \,
        \underset{
          n \in \mathbb{N}
        }{\prod}
        \,
        B^{n+1}
        \big(
          \pi_{n+1}(A) \otimes_{\scalebox{.5}{$\mathbb{Z}$}} \mathbb{R}
        \big)
      \Big)
      &
      \mbox{by Ex. \ref{WhiteheadLInfinityAlgebraOfLoopSpaces}}
      \\
      & \simeq
      \underset{
        n \in \mathbb{N}
      }{\prod}
      \,
      H
      \Big(
        X;
        \,
        B^{n+1}
        \big(
          \pi_{n+1}(A) \otimes_{\scalebox{.5}{$\mathbb{Z}$}} \mathbb{R}
        \big)
      \Big)
      &
      \mbox{by Def. \ref{NonAbelianCohomology} \& Prop. \ref{HomFunctorOfAHomotopyCategoryRespectsCoProducts}}
      \\
      & =
      \underset{
        n \in \mathbb{N}
      }{\prod}
      \,
      H^{n+1}
      \big(
        X;
        \,
        \pi_{n+1}(A) \otimes_{\scalebox{.5}{$\mathbb{Z}$}} \mathbb{R}
      \big)
      &
      \mbox{by Ex. \ref{OrdinaryCohomology}}
      \\
      & =
      \underset{
        n \in \mathbb{N}
      }{\bigoplus}
      \,
      H^{n+1}
      \big(
        X;
        \,
        \pi_{n+1}(A) \otimes_{\scalebox{.5}{$\mathbb{Z}$}} \mathbb{R}
      \big)
      &
      \mbox{by finite type}
      \,.
    \end{array}
  $$
  \vspace{-.8cm}

\end{proof}

\medskip

\noindent {\bf Relative rational Whitehead $L_\infty$-algebras.}
In parameterized generalization of Prop. \ref{WhiteheadLInfinityAlgebras}, we have:

\begin{prop}[Relative real Whitehead $L_\infty$-algebras]
  \label{WhiteheadLInfinityAlgebrasRelative}
  For $A, B, F \in \NilpotentConnectedQFiniteHomotopyTypes$
  (Def. \ref{NilpotentConnectedSpacesOfFiniteRationalType})
  and $p$ a Serre fibration (Example \ref{ClassicalModelStructureOnTopologicalSpaces})
  from $A$ to $B$ with fiber $F$
  \vspace{-2mm}
  $$
    \xymatrix@R=1em@C=3em{
      F \ar[r]^-{ \mathrm{fib}(p) }
      &
      A
      \ar[d]_-{p}^-{ \in\; \mathrm{Fib} }
      \\
      & B
    }
  $$
  \vspace{-2mm}

  \noindent there exists a nilpotent $L_\infty$-algebra
  (Def. \ref{NilpotentLInfinityAlgebras})
   \vspace{-1mm}
  \begin{equation}
    \label{WhiteheadLInfinityAlgebra}
    \mathfrak{l}_{\scalebox{.7}{$B$}}A
    \;\in\;
    \LInfinityAlgebrasNil
    \,,
  \end{equation}

  \vspace{-1mm}
  \noindent
  unique up to isomorphism,
  whose Chevalley-Eilenberg algebra (Def. \ref{ChevalleyEilenbergAlgebraOfLInfinityAlgebra})
  is the relative minimal model
  (Def. \ref{MinimalSullivanModels}, Prop. \ref{ExistenceOfRelativeMinimalSullianModels})
  of the PL de Rham complex of $p$ (Def. \ref{PLdeRhamComplexes}),
  relative to $\mathrm{CE}(\mathfrak{l}B)$ (from Prop. \ref{WhiteheadLInfinityAlgebras}):
  \vspace{-3mm}
  \begin{equation}
    \label{CEAlgebraOfWhiteheadLInfinityAlgebraRelative}
    \mathrm{CE}
    (
      \mathfrak{l}_{\scalebox{.5}{$B$}} A
    )
    \;:=\;
    \xymatrix@C=5em@R=1.4em{
      \big(
        \Omega^\bullet_{\mathrm{PLdR}}(A)
      \big)_{\mathrm{min}_{\mathrm{CE}(\mathfrak{l}B)}}
      \ar[rr]_-{ \in \, \mathrm{W} }^-{
        p^{\mathrm{min}_B}_A
      }
      &&
      \Omega^\bullet_{\mathrm{PLdR}}(A)
      \\
      &
      \mathrm{CE}(\mathfrak{l}B)
      \ar@{_{(}->}[ul]^-{
        \mathllap{
          \mbox{
            \tiny
            \color{darkblue}
            \bf
            relative minimal model
          }
        }
        \mathrm{CE}(\mathfrak{l}p)
      }
      \ar[r]_-{ \in\;\mathrm{W} }^-{
        p^{\mathrm{min}}_{B}
      }
      &
      \Omega^\bullet_{\mathrm{PLdR}}(B)\,.
      \ar[u]_-{
        \Omega^\bullet_{\mathrm{PLdR}}(p)
      }
    }
  \end{equation}
\end{prop}
\begin{proof}
By the PL de Rham theorem (Prop. \ref{PLdeRhamTheorem})
and the assumption that $A$ and $B$ are connected, it follows that we have
$
  H\Omega^0_{\mathrm{PLdR}}(A)
  \;=\;
  \mathbb{R}
$
and
$
  H\Omega^0_{\mathrm{PLdR}}(B)
  \;=\;
  \mathbb{R}
  \,.
$
Moreover, by the assumption that $p$ is a Serre fibration
with connected fiber, it follows that
$H^1(\Omega^\bullet_{\mathrm{PLdR}}(p))$ is injective
(e.g. \cite[p. 196]{FHT00}).
Therefore, Prop. \ref{ExistenceOfRelativeMinimalSullianModels}
applies and says that
$\big( \Omega^\bullet_{\mathrm{PLdR}}(A) \big)_{\mathrm{min}_{B}}
 \in \SullivanModelsConnected
$
exists, and is unique up to isomorphism.
With this, the equivalence \eqref{NilpotentLInfinityAlgebrasEquivalentToConnectedSullivanModels}
says that $\mathfrak{l}_{\scalebox{.5}{$B$}}A$ exists and is unique up to isomorphism.
\end{proof}

In parameterized generalization of Prop. \ref{RationalizationInTermsOfWhiteheadLInfinityAlgebras}
we have:
\begin{prop}[Relative $\mathbb{R}$-rationalization as integration of
 relative Whitehead $L_\infty$-algebras]
 \label{RelativeRationalizationAsIntegralOfRelativeWhiteheadLinfinityAlgebras}
 For a Serre fibration $A \xrightarrow{\;p\;} B$ as in
 Prop. \ref{WhiteheadLInfinityAlgebrasRelative},
 its rationalization over the real numbers (Def. \ref{Lk})
 is equivalently the image under $\Bexp$ \eqref{QuillenAdjunctBetweendgcAlgsAndSimplicialSets}
 of the image under forming CE-algebras \eqref{CEAlgebraOfWhiteheadLInfinityAlgebra}
 of its relative Whitehead $L_\infty$-algebra \eqref{WhiteheadLInfinityAlgebra}:
 \vspace{-2mm}
 $$
   L_{\mathbb{R}}
   \Bigg(
     \!\!\!
     \raisebox{5pt}{
     \begin{tikzcd}[row sep=7pt]
       A
       \ar[
         d,
         "\, p"
       ]
       \\
       B
     \end{tikzcd}
     }
     \!\!\!
   \Bigg)
   \;\;
     \simeq
   \;\;
     \raisebox{5pt}{
     \begin{tikzcd}[row sep=9pt]
       \Bexp_{\mathrm{PL}} \mathrm{CE}( \mathfrak{l}_{{}_B} A )
       \ar[
         d,
         "{
           \; \Bexp_{\mathrm{PL}} \mathrm{CE}( \mathfrak{l}p )
         }"
       ]
       \\
       \Bexp_{\mathrm{PL}} \mathrm{CE}( \mathfrak{l} B )
     \end{tikzcd}
     }
 $$
\end{prop}
\begin{proof}
  As in Prop. \ref{RationalizationInTermsOfWhiteheadLInfinityAlgebras},
  now appealing to
  Prop. \ref{WhiteheadLInfinityAlgebrasRelative} for the
  (co)fibrant replacement.
\end{proof}

\begin{lemma}[Minimal relative Sullivan models preserve homotopy fibers
  {\cite[\S 15 (a)]{FHT00}\cite[Thm. 5.1]{FHT15}}]
  \label{MinimalRelativeSullivanModelsPreserveHomotopyFibers}
  Consider $F, A, B \in$
  $\NilpotentConnectedQFiniteHomotopyTypes$
  (Def. \ref{NilpotentConnectedSpacesOfFiniteRationalType})
  and let $p$ be a Serre fibration from $A$ to $B$
  (Ex. \ref{ClassicalModelStructureOnTopologicalSpaces})
  such that the homology groups
  $H_\bullet(F,\mathbb{R})$
  of the fiber are nilpotent as $\pi_1(B)$-modules
  (for instance in that $B$ is simply-connected or that the fibration
  is principal).
  Then the cofiber of the minimal relative Sullivan model for $p$
  \eqref{CEAlgebraOfWhiteheadLInfinityAlgebraRelative}
  is
  the minimal Sullivan model \eqref{CEAlgebraOfWhiteheadLInfinityAlgebra}
  for the homotopy fiber $F$ (Def. \ref{HomotopyFibers}):
    \vspace{-2mm}
  \begin{equation}
    \raisebox{20pt}{
    \xymatrix@R=1.5em@C=3em{
      F
      \ar[rr]^-{
        \mathrm{fib}(p)
      }
      &&
      A
      \ar[d]_-{
        p
      }^-{ \in\, \mathrm{Fib} }
      &&
      \mathrm{CE}(\mathfrak{l}F)
      \ar@{<-}[rr]^-{
         \mathrm{cofib}
         (
           \mathrm{CE}(\mathfrak{l}p)
         )
      }
      &
      &
      \mathrm{CE}
      \big(
        \mathfrak{l}_{\scalebox{.5}{$B$}} A
      \big)
      \ar@{<-^{)}}[d]^-{
        \mathrm{CE}(\mathfrak{l}p)
      }
      \\
      &&
      B
      && &&
      \mathrm{CE}(\mathfrak{l}B)
    }
    }
  \end{equation}
\end{lemma}

\noindent {\bf Twisted non-abelian real cohomology.}

\begin{prop}[$\mathbb{R}$-Rationalization of local coefficients { -- the {\it fiber lemma} \cite[\S II]{BousfieldKan72}}]
  \label{RationalizationOfLocalCoefficients}
  Let
  \vspace{-2mm}
  $$
    \xymatrix@R=1.2em{
      A \ar[r] & A \!\sslash\! G
      \ar[d]^-{ \rho }
      \\
      & B G
    }
  $$

  \vspace{-3mm}
  \noindent  be a local coefficient bundle
  (Def. \ref{NonabelianTwistedCohomology})
  such that all spaces are connected, nilpotent and of $\mathbb{Q}$-finite type:
  $A,\, BG,\, A \!\sslash\! G \,\in\, \NilpotentConnectedQFiniteHomotopyTypes$
  (Def. \ref{NilpotentConnectedSpacesOfFiniteRationalType}),
  and such that the action of $\pi_1(BG)$ on $H_{\bullet}(A, \mathbb{R})$
  is nilpotent (for instance in that $BG$ is simply connected). Then:
  $\mathbb{R}$-Rationalization (Def. \ref{Lk})
  preserves the homotopy fiber:
  \vspace{-4mm}
  \begin{equation}
  \label{RationalizationOfLocalCoefficientBundle}
  \begin{tikzcd}[row sep={between origins,30pt}, column sep={between origins, 61pt}]
    &&
    \Bexp_{\mathrm{PL}} \mathrm{CE}(\mathfrak{l} A)
    \ar[
      rrr,
      "
        \Bexp_{\mathrm{Pl}}
        \big(
          \mathrm{cof}(\mathrm{CE}(\mathfrak{l}p))
        \big)
      "
    ]
    & & &
    \Bexp_{\mathrm{PL}} \mathrm{CE}\big(\mathfrak{l}_{{}_{B G}} A \!\sslash\! G  \big)
    \ar[
      ddd,
      "{
        \Bexp_{\mathrm{PL}} \mathrm{CE}( \mathfrak{l}p )
      }"
    ]
    \\
    &
    L_{\mathbb{R}} A
    \ar[
      rrr,
      "{
        \mathrm{hofib}
        ( L_{\mathbb{R}}p )
      }"
    ]
    \ar[
      ur,-,
      shift left=1pt
    ]
    \ar[
      ur,-,
      shift right=1pt
    ]
    &&&
    \big(
      L_{\mathbb{R}}
    \big)
    \!\sslash\!
    \big(
      L_{\mathbb{R}}G
    \big)
    \ar[
      ddd,
      "{
        L_{\mathbb{R}}(\rho)
      }"
    ]
    \ar[
      ur,-,
      shift left=1pt
    ]
    \ar[
      ur,-,
      shift right=1pt
    ]
    \\
    A
    \ar[
      rrr,
      "{
        \mathrm{hofib}(\rho)
      }"
    ]
    \ar[
      ur,
      "{
        \Derived \eta^{\mathbb{R}}_{A}
      }"{above, xshift=-5pt}
    ]
    &&&
    A
      \!\sslash\!
    G
    \ar[
      ddd,
      "\rho"
    ]
    \ar[
      ur,
      "{
        \Derived \eta^{\mathbb{R}}_{ A \!\sslash\! G }
      }"{below, xshift=5pt}
    ]
    \\
    &&&
    &&
    \Bexp_{\mathrm{PL}} \mathrm{CE}(\mathfrak{l} B G)
    \\
    &&&
    &
    L_{\mathbb{R}}(B G)
    \ar[
      ur,-,
      shift left=1pt
    ]
    \ar[
      ur,-,
      shift right=1pt
    ]
    \\
    &&&
    B G
    \ar[
      ur,
      "{
        \Derived \eta^{\mathbb{R}}_{B G}
      }"{below, xshift=5pt}
    ]
  \end{tikzcd}
  \;\;
  \in
  \mathrm{Ho}
  \big(
    \SimplicialSets_{\mathrm{Qu}}
  \big)
  \end{equation}
\end{prop}
\begin{proof}
  By
  Prop. \ref{WhiteheadLInfinityAlgebrasRelative},
  Prop. \ref{RelativeRationalizationAsIntegralOfRelativeWhiteheadLinfinityAlgebras}
  and since $\Bexp_{\mathrm{PL}}$ preserves fibrations, being a right adjoint,
  the homotopy fiber of $L_{\mathbb{R}}(p)$ is the image under
  $\Bexp_{\mathrm{PL}}$ of the cofiber of $\mathrm{CE}(\mathfrak{l}p)$.
  That this is a claimed is
  the content of Lemma \ref{MinimalRelativeSullivanModelsPreserveHomotopyFibers}.
\end{proof}

Due to Prop. \ref{RationalizationOfLocalCoefficients},
it makes sense to say, in generalization of Def. \ref{NonAbelianRealCohomology}:

\begin{defn}[Twisted non-abelian real cohomology]
  \label{TwistedNonAbelianRealCohomology}
  Let $X \in \TopologicalSpaces$
  and let $\xymatrix@C=12pt{A \!\sslash G \ar[r]^-{ \rho } & B G}$
  be a local coefficient bundle
  (Prop. \ref{GActionsAsFibrations}, Def. \ref{NonabelianTwistedCohomology})
  in $\NilpotentConnectedQFiniteHomotopyTypes$ (Def. \ref{NilpotentConnectedSpacesOfFiniteRationalType}).
    Then the \emph{twisted non-abelian real cohomology}
  of $X$ with local coefficients $\rho$ is the
  twisted non-abelian $L_{\mathbb{R}}A$-cohomology (Def. \ref{NonabelianTwistedCohomology})
  of $X$ with coefficients in the rationalized
  local coefficient bundle $L_{\mathbb{R}}(\rho)$ from
  Prop. \ref{RationalizationOfLocalCoefficients}:
      \vspace{-2mm}
  $$
    H^\tau
    \big(
      X;
      \,
      L_{\mathbb{R}} A
    \big)
    \;\;:=\;\;
    \RationallyTwistedHomotopyTypes
    \big(
      \tau
      \,,\,
      L_{\mathbb{R}}(\rho)
    \big)
    \,.
  $$
\end{defn}

%%%%%%%%%%%%%%%%%%%%%%%%%%%%%%%%%%%%%%%%%%%%%%%%%%
\subsection{Non-abelian de Rham theorem}
\label{NonAbelianDeRhamTheory}
%%%%%%%%%%%%%%%%%%%%%%%%%%%%%%%%%%%%%%%%%%%%%%%%%%%

We establish non-abelian de Rham theory
for differential forms with values in (nilpotent)
$L_\infty$-algebras, following
\cite{SatiSchreiberStasheff08} \cite{FiorenzaSchreiberStasheff10}.
The main result is the
non-abelian de Rham theorem, Theorem \ref{NonAbelianDeRhamTheorem},
and its generalization to the
twisted non-abelian de Rham theorem,
Theorem \ref{TwistedNonAbelianDeRhamTheorem}.

\medskip

\noindent {\bf $L_\infty$-Algebra valued differential forms.}

\begin{defn}[Flat $L_\infty$-algebra valued differential forms {\cite[\S 6.5]{SatiSchreiberStasheff08}\cite[\S 4.1]{FiorenzaSchreiberStasheff10}}]
  \label{FlatLInfinityAlgebraValuedDifferentialForms}
  $\,$

  \noindent {\bf (i)}   For $X \in \mathrm{SmoothManifold}$
  and $\mathfrak{g} \in \LInfinityAlgebras$ (Def. \ref{ChevalleyEilenbergAlgebraOfLInfinityAlgebra}),
  a \emph{flat $\mathfrak{g}$-valued differential form on $X$}
  is a morphism of dgc-algebras (Def. \ref{CategoryOfdgAlgebras})
  \vspace{-1mm}
  \begin{equation}
    \label{FlatLInfinityAlgebraValuedDifferentialForm}
    \xymatrix{
      \Omega^\bullet_{\mathrm{dR}}(X)
      \ar@{<-}[r]^-{ A }
      &
      \mathrm{CE}(\mathfrak{g})
    }
    \;\;\;\;
    \in
    \;
    \dgcAlgebras{\mathbb{R}}
  \end{equation}

  \vspace{-1mm}
  \noindent
  to the smooth de Rham dgc-algebra of $X$ (Example \ref{SmoothdeRhamComplex})
  from the
  Chevalley-Eilenberg dgc-algebra of $\mathfrak{g}$
  (Def. \ref{ChevalleyEilenbergAlgebraOfLInfinityAlgebra}).

  \noindent {\bf (ii)}  We write
  \vspace{-1mm}
  \begin{equation}
    \label{SetOfFlatLInfinityAlgebraValuedForms}
    \Omega_{\mathrm{dR}}
    (
      X;
      \,
      \mathfrak{g}
    )_{\mathrm{flat}}
    \;:=\;
    \dgcAlgebras{\mathbb{R}}
    \big(
      \mathrm{CE}(\mathfrak{g})
      \,,\,
      \Omega^\bullet_{\mathrm{dR}}(X)
    \big)
  \end{equation}
  for the set of all flat $\mathfrak{g}$-valued forms on $X$.
\end{defn}

\begin{example}[Flat Lie algebra valued differential forms]
  \label{FlatLieAlgebraValuedDifferentialForms}
  Let $\mathfrak{g} \in \LieAlgebras$ be a Lie algebra
  \eqref{DifferentialGradedLieAlgebrasInsideLInfinityAlgebras}
  with Lie bracket $[-,-]$.
  Then a flat $\mathfrak{g}$-valued differential form
  in the sense of Def. \ref{FlatLInfinityAlgebraValuedDifferentialForms}
  is the traditional concept: a $\mathfrak{g}$-valued 1-form
  satisfying the Maurer-Cartan equation:
  \vspace{-2mm}
  \begin{equation}
    \label{FlatLieAlgebrasValuedFormsInComponents}
    \Omega^\bullet_{\mathrm{dR}}
    (
      X;
      \,
      \mathfrak{g}
    )_{\mathrm{flat}}
    \;\simeq\;
    \Big\{
    \left.
      A \in \Omega^1_{\mathrm{dR}}(X)\otimes \mathfrak{g}
      \;\right\vert\;
      d A + [A\wedge A] \;=\; 0
    \; \Big\}
    \,.
  \end{equation}

    \vspace{-2mm}
\noindent  One way to see this is to appeal to the classical
  fact that the Chevalley-Eilenberg algebra of a
  finite-dimensional Lie algebra (Example \ref{CEAlgebraOfLieAlgebra}) is
  isomorphic to the dgc-algebra of left invariant differential
  forms on the corresponding Lie group $G$, which is generated
  from the Maurer-Cartan form
  $\theta \in \Omega^1_{\mathrm{dR}}(G) \otimes \mathfrak{g}$
  satisfying $\theta_{\vert T_e G} = \mathrm{id}_{\mathfrak{g}}$
  and $d \theta = [\theta \wedge \theta]$.
    More explicitly, for $\{v_a\}$ a linear basis for $\mathfrak{g}$ \eqref{LinearBasisForLieAlgebra} with structure constants
  $\{f_{a b}^c\}$ \eqref{LieAlgebraStructureConstants},
  we see from \eqref{CEAlgebraOfLieAlbebraInComponents}
  that a dgc-algebra homomorphims \eqref{FlatLInfinityAlgebraValuedDifferentialForm}
  has the following components (second line)
  and constraints (third line):
      \vspace{-2mm}
  \begin{equation}
    \label{ComponentsAndConstraintsOfFlatLieAlgebraValuedForms}
    \hspace{-2cm}
    \xymatrix@R=4pt{
      \Omega^\bullet_{\mathrm{dR}}(X)
      \ar@{<-}[rr]_-{A}^-{
        \overset{
          \mathclap{
          \raisebox{3pt}{
            \tiny
            \color{darkblue}
            \bf
            \def\arraystretch{.9}
            \begin{tabular}{c}
            flat Lie algebra valued
            differential form
            \end{tabular}
          }
          }
        }{
%          A
        }
      }
      &&
      \mathbb{R}
      \big[
        \{\theta^{(a)}_1\}
      \big]
      \!\big/\!
      \left(
        d\, \theta^{(c)}_1
        =
        f_{a b}^c \,
        \theta^{(b)}_1
          \wedge
        \theta^{(a)}_1
      \right)
      \mathrlap{
        \;
        \simeq
        \;
        \mathrm{CE}(\mathfrak{g})\;.
      }
      \\
      A^{(c)}
      \ar@{<-|}[rr]^-{
        \mbox{
          \tiny
          \color{darkblue}
          \bf
          components
        }
      }
      \ar@{|->}[dd]_-{ d }
      &&
      \theta^{(c)}_1
      \ar@{|->}[dd]^-{ d }
      \\
      \\
      d A^{(c)}
      \ar@{=}[r]^-{
        \mbox{
          \tiny
          \color{darkblue}
          \bf
          constraints
        }
      }
      &
      f_{a b}^c \;
      A^{(b)}
        \wedge
      A^{(a)}
      \ar@{<-|}[r]
      &
      f_{a b}^c \,
      \theta^{(a)}_1 \wedge \theta^{(b)}_1
    }
  \end{equation}
\end{example}

\begin{example}[Ordinary closed forms are flat line $L_\infty$-algebra valued forms]
  \label{OrdinaryClosedFormsAreFlatLineLInfinityAlgebraValuedForms}
  For $n \in \mathbb{N}$, consider $\mathfrak{g} = \mathfrak{b}^n \mathbb{R}$
  the line Lie $(n+1)$-algebra (Example \ref{LineLienPlusOneAlgebras}). Then
  the corresponding flat $\mathfrak{g}$-valued differential
  forms (Def. \ref{FlatLInfinityAlgebraValuedDifferentialForms})
  are in natural bijection to ordinary closed $(n+1)$-forms:
      \vspace{-1mm}
  \begin{equation}
    \label{FlatLineLinenPlus1AlgebraValuedFormsAreClosedForms}
    \Omega_{\mathrm{dR}}
    (
      X;
      \,
      \mathfrak{b}^n \mathbb{R}
    )_{\mathrm{flat}}
    \;\simeq\;
    \Omega^{n+1}_{\mathrm{dR}}(X)_{\mathrm{closed}}
    \,.
  \end{equation}

    \vspace{0mm}
\noindent
 That is, by \eqref{CEAlgebraOfLineLienAlgebra}, we see
  that the elements on the left
  of \eqref{FlatLineLinenPlus1AlgebraValuedFormsAreClosedForms}
  have the following component
  (second line) subject to the follows constraint (third line):

  \vspace{-5mm}
  \begin{equation}
    \label{FlatLineLInfinityAlgebraValuedFormsInComponents}
    \hspace{-2cm}
    \xymatrix@R=2pt@C=4em{
      \Omega^\bullet_{\mathrm{dR}}(X)
      \ar[rr]^-{
        \mbox{
          \tiny
          \color{darkblue}
          \bf
          \def\arraystretch{.9}
          \begin{tabular}{c}
            flat
            \\
            line Lie $(n+1)$-algebra-valued
            \\
            differential form
          \end{tabular}
        }
      }
      &&
      \mathbb{R}
      [
        c_{n+1}
      ]
      \!\big/\!
      (
        d\, c_{n+1} = 0
      )
      \mathrlap{
        \;\simeq\;
        \mathrm{CE}
        (
          \mathfrak{b}^{n}\mathbb{R}
        )\;.
      }
      \\
      C_{n+1}
      \ar@{<-}[rr]^-{
        \mbox{
          \tiny
          \color{darkblue}
          \bf
          component
        }
      }
      \ar@{|->}[dd]_-{ d }
      &&
      c_{n+1}
      \ar@{|->}[dd]^-{d}
      \\
      \\
      d C_{n+1}
      \ar@{=}[r]^-{
        \mbox{
          \tiny
          \color{darkblue}
          \bf
          constraint
        }
      }
      &
      0
      \ar@{<-|}[r]
      &
      0
    }
  \end{equation}
\end{example}

\begin{example}[Flat String Lie 2-algebra valued differential forms]
  \label{FlatStringLie2AlgebraValuedDifferentialForms}
  Flat $L_\infty$-algebras valued forms (Def. \ref{FlatLInfinityAlgebraValuedDifferentialForms})
  with values in a String Lie 2-algebra $\mathfrak{string}_{\mathfrak{g}}$
  (Example \ref{StringLie2Algebra})
  are pairs consisting of a flat $\mathfrak{g}$-valued 1-form $A_1$
  (Example \ref{FlatLieAlgebraValuedDifferentialForms})
  and a coboundary 2-form $B_2$ for its Chern-Simons form
  $\mathrm{CS}(A) := c \big\langle A \wedge [A \wedge A] \big\rangle$:
      \vspace{-1mm}
  $$
    \Omega_{\mathrm{dR}}
    \big(
      X;
      \,
      \mathfrak{string}_{\mathfrak{g}}
    \big)_{\mathrm{flat}}
    \;\simeq\;
    \left\{\!\!\!
      \left.
      \begin{array}{c}
        B_2,
        \\
        A_1
      \end{array}
      \in
      \Omega^\bullet_{\mathrm{dR}}(X)
      \;\right\vert\;
      \begin{aligned}
        d\, B_2 & =  \tfrac{1}{c} \mathrm{CS}(A)\,,
        \\
        d\, A_1 & = - [A_1 \wedge A_1]
      \end{aligned}
    \right\}
    \,.
  $$

  \vspace{-1mm}
  \noindent
  Namely, from \eqref{CEAlgebraOfStringLie2Algebra} we see that in degree 1 the components of
  and constraints on such a differential form datum are
  exactly as in \eqref{ComponentsAndConstraintsOfFlatLieAlgebraValuedForms},
  while in degree 2 they are as follows:
      \vspace{-2mm}
  \begin{equation}
  \hspace{-3.7cm}
    \xymatrix@R=4pt{
      \Omega^\bullet_{\mathrm{dR}}(X)
      \ar@{<-}[rr]^-{
        \mbox{
          \tiny
          \color{darkblue}
          \bf
          flat String Lie 2-algebra valued form
        }
      }
      &&
      \mathbb{R}
      \left[
        \!\!\!\!
        {\begin{array}{c}
          b_2,
          \\
          \{\theta_1^{(a)}\}
        \end{array}}
        \!\!\!\!
      \right]
      \!\big/\!
   \scalebox{0.8}{$   \left(
      {\begin{aligned}
        d\, b_2\;
          & =
        \mu_{a b c}\, \theta_1^{(c)} \wedge \theta_1^{(b)} \wedge \theta_1^{(a)}
        \\[-2pt]
        d \theta_1^{(c)}
          & =
        f^c_{a b} \, \theta_1^{(b)} \wedge \theta_1^{(a)}
      \end{aligned}}
      \right)
      $}
      \mathrlap{
        \;\simeq\;
        \mathrm{CE}
        \big(
          \mathfrak{string}_{\mathfrak{g}}
        \big).
      }
      \\
      B_2
      \ar@{<-|}[rr]^-{
        \mbox{
          \tiny
          \color{darkblue}
          \bf
          component in degree 2
        }
      }
      \ar@{|->}[dd]_-{ d }
      &&
      b_2
      \ar@{|->}[dd]^-{ d }
      \\
      \\
      d B_2
      \ar@{=}[r]^-{
        \mbox{
          \tiny
          \color{darkblue}
          \bf
          constraint
        }
      }
      &
      \mu_{a b c}
      \,
      A_1^{(c)}
        \wedge
      A_1^{(b)}
        \wedge
      A_1^{(a)}
      \ar@{<-|}[r]
      &
      \mu_{a b c}
      \,
      \theta_1^{(c)}
        \wedge
      \theta_1^{(b)}
        \wedge
      \theta_1^{(a)}
    }
  \end{equation}
\end{example}

\vspace{0mm}
\begin{example}[Flat sphere-valued differential forms]
  \label{FlatSphereValuedDifferentialForms}
  Flat $L_\infty$-algebras valued forms (Def. \ref{FlatLInfinityAlgebraValuedDifferentialForms})
  with values in the rational Whitehead $L_\infty$-algebra
  (Prop. \ref{WhiteheadLInfinityAlgebras}) of a sphere
  (Ex. \ref{RationalizationOfnSpheres})
  of positive even dimension $2k$ are pairs consisting of a closed
  differential $2k$-form and a $(4k-1)$-form whose differential
  equals minus the wedge square of the $2k$-form:
  $$
    \Omega_{\mathrm{dR}}
    \big(
      -;
      \mathfrak{l}S^{2k}
    \big)
    \;\simeq\;
    \left\{
      \!
      \begin{array}{l}
        G_{4k-1},
        \\[-2pt]
        G_{2k}
      \end{array}
      \;
      \in
      \:
      \Omega^\bullet_{\mathrm{dR}}(X)
      \,\left\vert\,
      \begin{aligned}
        d\, G_{4k-1} & = - G_{2k} \wedge G_{2k},
        \\[-2pt]
        d \, G_{2k\phantom{-1}} & = 0
      \end{aligned}
      \right.
      \!
    \right\}
    \
    \,.
  $$

  \vspace{-2mm}
\noindent  Namely, from \eqref{SullivanModelForEvenDimensionalSpheres}
  one sees that the components of and the constraints on
  an $\mathfrak{l}S^{2k}$-valued form are as follows:
  \vspace{-2mm}
  \begin{equation}
  \hspace{-2cm}
    \xymatrix@R=12pt{
      \Omega^\bullet_{\mathrm{dR}}
      \big(
        X
      \big)
      \ar@{<-}[rr]^-{
        \mathclap{
        \mbox{
          \tiny
          \color{darkblue}
          \bf
          flat $\mathfrak{l}S^{2k}$-valued form
        }
        }
      }
      &&
      \mathbb{R}
      \!
      \left[
      \!\!\!
      {\begin{array}{l}
        \omega_{4k-1},
        \\[-2pt]
        \omega_{2k}
      \end{array}}
      \!\!\!
      \right]
      \!\big/\!
\scalebox{0.8}{$      \left(
      {\begin{aligned}
        d\, \omega_{4k-1} & = -\omega_{2k} \wedge \omega_{2k},
        \\[-2pt]
        d\, \omega_{2k\phantom{-1}} & = 0
      \end{aligned}}
      \!
      \right)
      $}
      \mathrlap{
        \;
        =
        \;
        \mathrm{CE}
        \big(
          \mathfrak{l}S^{2k}
        \big)
      }
      \\
      G_{2k}
      \ar[d]_-{ d }
      \ar@{<-|}[rr]^-{
        \mbox{
          \tiny
          \color{darkblue}
          \bf
          component in degree $2k$
        }
      }
      &&
      \omega_{2k}
      \ar@{|->}[d]^-{ d }
      \\
      d\, G_{2k}
      \ar@{=}[r]^-{
        \mbox{
          \tiny
          \color{darkblue}
          \bf
          constraint
        }
      }
      &
      0
      \ar@{<-|}[r]
      &
      0
      \\
      G_{4k-1}
      \ar@{|->}[d]_-{ d }
      \ar@{<-|}[rr]^-{
        \mbox{
          \tiny
          \color{darkblue}
          \bf
          component in degree $4k-1$
        }
      }
      &&
      \omega_{4k-1}
      \ar@{|->}[d]
      \\
      d\, G_{4k-1}
      \ar@{=}[r]^-{
        \mbox{
          \tiny
          \color{darkblue}
          \bf
          constraint
        }
      }
      &
      - G_{2k} \wedge G_{2k}
      \ar@{<-|}[r]
      &
      -\omega_{4k} \wedge \omega_{4k}
    }
  \end{equation}

  \vspace{-1.5mm}
\noindent  For $2k = 4$ this is the structure
  of the equations of
  motion of the C-field in 11-dimensional supergravity
  (modulo the Hodge self-duality constraint $G_7 = \star G_4$)
  \cite[\S 2.5]{Sati13}.
\end{example}

\newpage

\begin{example}[PL de Rham right adjoint via $L_\infty$-algebra valued forms]
  \label{PLDeRhamRightAdjointViaLInfinityAlgebraValuedForms}
  For $n \in \mathbb{N}$,
  the right adjoint functor
  in the PS de Rham adjunction \eqref{PSdRQuillenAdjunctBetweendgcAlgsAndSimplicialSets}
  sends the Chevalley-Eilenberg algebra
  (Def. \ref{ChevalleyEilenbergAlgebraOfLInfinityAlgebra})
  of any
  $\mathfrak{g} \in \LInfinityAlgebrasNil$
  (Def. \ref{NilpotentLInfinityAlgebras})
  to a simplicial set of flat $\mathfrak{g}$-valued
  differential forms (Def. \ref{FlatLInfinityAlgebraValuedDifferentialForms}):
  $$
    \flatBexp_{\mathrm{PL}}(\mathfrak{g})(\mathbb{R}^n)
    \;:=\;
    \Bexp_{\mathrm{PS},n}
    \big(
      \mathrm{CE}(\mathfrak{g})
    \big)
    \;:\;
    [k]
    \;\longmapsto\;
    \Omega_{\mathrm{dR}}
    \big(
      \mathbb{R}^n \times \Delta^k
      ;
      \,
      \mathfrak{g}
    \big)_{\mathrm{flat}}
    \;\;\;\;
    \in
    \;
    \SimplicialSets
  $$
  (by direct comparison of \eqref{PSexp} with \eqref{SetOfFlatLInfinityAlgebraValuedForms}).
  Regarded as a simplicial presheaf over
  $\CartesianSpaces$ (Def. \ref{SimplicialPresheavesOverCartesianSpaces}),
  this construction is the moduli $\infty$-stack
  of flat $L_\infty$-algebra valued differential forms
  (see \cref{NonabelianDifferentialCohomology} below).
\end{example}

\medskip
\noindent {\bf Non-abelian de Rham cohomology.}

\begin{defn}[Coboundaries between flat $L_\infty$-algebra valued forms]
  \label{CoboundariesBetweenFlatLInfinityAlgebraValuedForms}
  Let $X \in \SmoothManifolds$ and (from Def. \ref{ChevalleyEilenbergAlgebraOfLInfinityAlgebra})
  $\mathfrak{g} \in \LInfinityAlgebras$.
  For
      \vspace{-1mm}
  $$
    A^{(0)}\!,\; A^{(1)}
    \;\in\;
    \Omega_{\mathrm{dR}}
    (
      X;
      \,
      \mathfrak{g}
    )_{\mathrm{flat}}
  $$

  \vspace{-1mm}
  \noindent
  a pair of flat $\mathfrak{g}$-valued differential forms on
  $X$ (Def. \ref{FlatLInfinityAlgebraValuedDifferentialForms}),
  we say that a \emph{coboundary} between them
  is a flat $\mathfrak{g}$-valued differential form on the
  cylinder manifold over $X$
  (its Cartesian product manifold with the real line):
      \vspace{-2mm}
  \begin{equation}
    \label{FlatFormOnCylinder}
    \widetilde A
    \;\in\;
    \Omega
    (
      X \times \mathbb{R};
      \,
      \mathfrak{g}
    )_{\mathrm{flat}}
  \end{equation}

      \vspace{-3mm}
  \noindent
  such that its restrictions along
  \vspace{-2mm}
  $$
    \xymatrix{
      X
      \simeq
      X \times \{0\}
     \; \ar@{^{(}->}[r]^-{ i^{\, X}_0 }
      &
      \; X \times \mathbb{R}
     \; \ar@{<-^{)}}[r]^-{ i^{\, X}_1 }
      &
     \; X \times \{1\}
      \simeq
      X
    }
  $$

  \vspace{-3mm}
  \noindent
  are equal to $A^{(0)}$ and to $A^{(1)}$, respectively:
  \vspace{-1mm}
  \begin{equation}
    \label{RestrictionOfFormsOnCylinderToBoundary}
    (i^{\, X}_0)^\ast \widetilde A \;=\; A^{(0)}
        \phantom{AAAA}
    \mbox{and}
    \phantom{AAAA}
    (i^{\, X}_1)^\ast \widetilde A \;=\; A^{(1)}
    \,.
  \end{equation}

  \vspace{-2mm}
  \noindent
  If such a coboundary exists, we say that $A^{(0)}$
  and $A^{(1)}$ are \emph{cohomologous}, to be denoted
      \vspace{-2mm}
  $$
    A^{(0)}
    \;\sim\;
    A^{(1)}
    \,.
  $$
\end{defn}

\begin{defn}[Non-abelian de Rham cohomology]
  \label{NonabelianDeRhamCohomology}
  Let $X \in \SmoothManifolds$
  and $\mathfrak{g} \in \LInfinityAlgebras$
  (Def. \ref{ChevalleyEilenbergAlgebraOfLInfinityAlgebra}).
  Then the \emph{non-abelian de Rham cohomology}
  of $X$ with coefficients in $\mathfrak{g}$ is the set
  \vspace{-1mm}
  \begin{equation}
    \label{NonabelianDeRhamCohomologyAsQuotientSet}
    H_{\mathrm{dR}}
    (
      X;
      \,
      \mathfrak{g}
    )
    \;:=\;
    \big(
    \Omega_{\mathrm{dR}}
    (
      X;
      \,
      \mathfrak{g}
    )_{\mathrm{flat}}
    \big)_{\!\!/\sim}
  \end{equation}

  \vspace{-2mm}
  \noindent
  of equivalence classes with respect to the
  coboundary relation from Def. \ref{CoboundariesBetweenFlatLInfinityAlgebraValuedForms}
  on the set of flat $\mathfrak{g}$-valued differential forms
  on $X$
  (Def. \ref{FlatLInfinityAlgebraValuedDifferentialForms}).
\end{defn}

We recall the following basic facts (e.g. \cite[Rem 3.1]{GomiTerashima00}):
\begin{lemma}[Fiberwise Stokes theorem and Projection formula]
  \label{FiberwiseStokesTheorem}
  Let $X$ be a smooth manifold
  and let $F$ be a compact smooth manifold with corners,
  e.g. $F = \Delta^k$ a standard $k$-simplex,
  which for $k = 1$ is the interval $F = [0,1]$.

  Then fiberwise integration over $F$ of differential forms
  on the Cartesian product manifold $X \times F$
  \vspace{-2mm}
  $$
    \xymatrix{
      \Omega^\bullet_{\mathrm{dR}}
      (
        X
          \times
        F
      )
      \ar[rr]^-{ \int_{F} }
      &&
      \Omega^{\bullet- \mathrm{dim}(F) }_{\mathrm{dR}}
      (
        X
      )
    }
    \;\;\;\;\;\mbox{e.g.}\;\;\;\;\;
    \xymatrix{
      \Omega^\bullet_{\mathrm{dR}}
      (
        X
          \times
        \mathbb{R}
      )
      \ar[rr]^-{ \int_{[0,1]} }
      &&
      \Omega^{\bullet-1}_{\mathrm{dR}}
      (
        X
      )
    }
  $$

  \vspace{-2mm}
  \noindent
  satisfies, for all
  $\alpha \in \Omega_{\mathrm{dR}}^\bullet(X \times F)$
  and
  $\beta \in \Omega^\bullet_{\mathrm{dR}}(X)$:

  \noindent
  {\bf (i)} The fiberwise Stokes formula:
    \vspace{-2mm}
  \begin{equation}
    \label{StokesFormula}
    \int_F d \alpha
    \;\;=\;\;
    (-1)^{\mathrm{dim}(F)}
    \,
    d \int_F \alpha
    +
    \int_{\partial F} \alpha
    \;\;\;\;\;\;\;
    \mathrm{e.g.}
    \;\;\;\;\;\;\;
    d \int_{[0,1]} \alpha
    \;\;=\;\;
    (i^X_1)^\ast\alpha
    \;-\;
    (i^X_0)^\ast \alpha
    \;-\;
    \int_{[0,1]} d \alpha
    \,,
  \end{equation}

  \vspace{-2mm}
  \noindent
  where
  $$
    \xymatrix{
      X
      \;\simeq\;
      X \times \{0\}
      \; \ar@{^{(}->}[r]^-{ i^X_0 }
      &
      X \times \mathbb{R}
      \ar@{<-^{)}}[r]^-{ i^X_1 }
      &
   \;   X \times \{1\}
      \;\simeq\;
      X
    }
  $$

  \vspace{-2mm}
  \noindent
  are the boundary inclusions.

  \noindent
  {\bf (ii)} The projection formula
  \begin{equation}
    \label{ProjectionFormula}
    \int_F
      \big( \mathrm{pr}_X^\ast \beta \big) \wedge \alpha
    \;=\;
    (-1)^{\mathrm{dim}(F) \, \mathrm{deg}(\beta) }
    \beta \wedge \int_{F} \alpha\,,
    \;\;\;\;\;
    \mbox{e.g.}
    \;\;\;\;\;
    \int_{[0,1]}
      \big( \mathrm{pr}_X^\ast \beta \big) \wedge \alpha
    \;=\;
    (-1)^{\mathrm{deg}(\beta)}
    \beta \wedge \int_{[0,1]} \alpha\,,
  \end{equation}

  \vspace{-1mm}
  \noindent
  where

  \vspace{-1mm}
  \noindent
  $$
    \xymatrix{
      X \times F
      \ar[r]^-{ \mathrm{pr}_X }
      &
      X
    }
  $$
  is projection onto the first factor.
\end{lemma}

\begin{prop}[Non-abelian de Rham cohomology subsumes ordinary de Rham cohomology]
  \label{NonAbelianDeRhamCohomologySubsumesOrdinaryDeRhamCohomology}
  For any $n \in \mathbb{N}$, let
  $\mathfrak{g} = \mathfrak{b}^n \mathbb{R}$
  be the line Lie $(n+1)$-algebra
  (Example \ref{LineLienPlusOneAlgebras}).
  Then the non-abelian de Rham cohomology
  with coefficients in $\mathfrak{g}$
  (Def. \ref{NonabelianDeRhamCohomology})
  is naturally equivalent
  to ordinary de Rham cohomology in degree $n+1$:
  \vspace{-2mm}
  \begin{equation}
    \label{NonAbelianDeRhamCohomologyWithCoefficientsInLineLienAlgebra}
    H_{\mathrm{dR}}
    (
      -;
      \,
      \mathfrak{b}^n \mathbb{R}
    )
    \;\;\simeq\;\;
    H^{n+1}_{\mathrm{dR}}
    (
      -
    )
    \,.
  \end{equation}
\end{prop}
\begin{proof}
  From Example \ref{OrdinaryClosedFormsAreFlatLineLInfinityAlgebraValuedForms},
   we know that the canonical cocycle sets are
   in natural bijection
       \vspace{-1mm}
   $$
     \Omega_{\mathrm{dR}}
     (
       X;
       \,
       \mathfrak{b}^n \mathbb{R}
     )_{\mathrm{flat}}
     \;\simeq\;
     \Omega^{n+1}_{\mathrm{dR}}
     (
       X
     )_{\mathrm{closed}}
     \,.
   $$

  \vspace{-2mm}
  \noindent
  Therefore, it just remains to see that the
  coboundary relations in both cases coincide.
  By the explicit nature \eqref{FlatLineLInfinityAlgebraValuedFormsInComponents}
  of the above natural bijection and by the
  Definition \ref{CoboundariesBetweenFlatLInfinityAlgebraValuedForms}
  of non-abelian coboundaries,
  we hence need to see that a pair of closed forms
    \vspace{-2mm}
  $$
    C^{(0)}_{n+1}
    \!,\;
    C^{(1)}_{n+1}
    \;\;\;
    \in
    \Omega^{n+1}_{\mathrm{dR}}(X)_{\mathrm{closed}}
  $$

    \vspace{-2mm}
  \noindent
  has a de Rham coboundary, i.e.,
      \vspace{-2mm}
  \begin{equation}
    \label{ExistenceOfDeRhamCoboundary}
    \exists
    \;
    h_n \;\in\; \Omega^n_{\mathrm{dR}}(X)
    \,,
    \phantom{A}
    \mbox{
      such that}
      \phantom{A}
        C^{(0)}_{n+1}
        +
        d h_n
        \;=\;
        C^{(1)}_{n+1}
      \;,
      \end{equation}

      \vspace{-2mm}
  \noindent
  precisely if the pair extends to a closed $(n+1)$-form
  on the cylinder over $X$, as in \eqref{FlatFormOnCylinder}
  \eqref{RestrictionOfFormsOnCylinderToBoundary}:
  \vspace{-2mm}
  \begin{equation}
    \label{ExistenceOfExtensionOverCylinderManifold}
    \exists
    \;
    {\widetilde C}_{n+1}
    \;\in\;
    \Omega^{n+1}_{\mathrm{dR}}
    (
      X \times \mathbb{R}
    )_{\mathrm{closed}}
    \,,
    \phantom{A}
    \mbox{
      such that
      }
     \phantom{A}
     \big(i^{\, X}_0\big)^{\!\! \ast}
       \widetilde C_{n+1}
       \;=\;
       C_{n+1}^{(0)}
     \phantom{A}
       \mbox{and}
     \phantom{A}
     \big(i^{\, X}_1\big)^{\!\! \ast}
     \widetilde C_{n+1} \;=\; C_{n+1}^{(1)}
    \;.
  \end{equation}

  \vspace{-1mm}
  \noindent
  That \eqref{ExistenceOfDeRhamCoboundary}
  $\Leftrightarrow$
  \eqref{ExistenceOfExtensionOverCylinderManifold}
  is a standard argument:
  Let $t$ denote the canonical coordinate function on
  $\mathbb{R}$.
  In one direction, given $h_n$ as in \eqref{ExistenceOfDeRhamCoboundary},
  the choice
  \vspace{-3mm}
  $$
    \widetilde C_{n+1} \;:=\;
    (1-t)\, \mathrm{pr}_X^\ast\big(C^{(0)}_{n+1}\big)
    +
    t \, \mathrm{pr}_X^\ast\big( C^{(1)}_{n+1} \big)
    +
    d t \wedge
    \mathrm{pr}_X^\ast\big( h_n\big)
  $$

      \vspace{-2mm}
\noindent
  clearly satisfies \eqref{ExistenceOfExtensionOverCylinderManifold}.
  In the other direction, given $\widetilde C_{n+1}$ as in
  \eqref{ExistenceOfExtensionOverCylinderManifold}, the choice
    \vspace{-2mm}
  $$
    h_n
    \;:=\;
    \int_{[0,1]}
    \widetilde C_{n+1}
  $$

  \vspace{-1mm}
  \noindent
  satisfies \eqref{ExistenceOfDeRhamCoboundary},
  by the fiberwise Stokes theorem (Lemma \ref{FiberwiseStokesTheorem}).
\end{proof}

\medskip

\noindent {\bf The non-abelian de Rham theorem.}

\begin{theorem}[Non-abelian de Rham theorem]
  \label{NonAbelianDeRhamTheorem}
  Let $X \in \HomotopyTypes$
  and $A \in \NilpotentConnectedQFiniteHomotopyTypes$
  (Def. \ref{NilpotentConnectedSpacesOfFiniteRationalType}),
  and let $X$ admit the structure of
  a smooth manifold.
  Then the non-abelian de Rham cohomology (Def. \ref{NonabelianDeRhamCohomology})
  of $X$ with
  coefficients in the real Whitehead $L_\infty$-algebra
  $\mathfrak{l}A$ (Prop. \ref{WhiteheadLInfinityAlgebras})
  is in natural bijection with the non-abelian
  real cohomology (Def. \ref{NonAbelianRealCohomology})
  of $X$ with coefficients in $L_{\mathbb{R}}A$
  (Def. \ref{Lk}):
  \vspace{-1mm}
  \begin{equation}
    \label{EquivalenceBetweenNonabelianRealAndNonaebelianDeRhamCohomology}
    H
    \big(
      X;
      \,
      L_{\mathbb{R}}A
    \big)
    \;\simeq\;
    H_{\mathrm{dR}}
    (
      X;
      \,
      \mathfrak{l}A
    )\;.
  \end{equation}
\end{theorem}
\begin{proof}
This follows by the following sequence of bijection:
\begin{equation}
  \label{TowardsTheNonAbelianDeRhamTheorem}
  \begin{array}{lll}
    H
    \big(
      X;
      \,
      L_{\mathbb{R}}A
    \big)
    & =
    \HomotopyTypes
    \big(
      X
      \,,\,
      L_{\mathbb{R}}A
    \big)
    &
    \mbox{by Def. \ref{NonAbelianRealCohomology}}
    \\
    & \simeq
    \mathrm{Ho}
    \big(
      \dgcAlgebrasProj{\mathbb{R}}
    \big)
    \big(
      \Omega^\bullet_{\mathrm{PLdR}}(A)
      \,,\,
      \Omega^\bullet_{\mathrm{PLdR}}(X)
    \big)
    &
    \mbox{by Def. \ref{Lk} \& Prop. \ref{QuillenAdjunctionBetweendgcAlgebrasAndSimplicialSets} }
    \\
    & \simeq
    \mathrm{Ho}
    \big(
      \dgcAlgebrasProj{\mathbb{R}}
    \big)
    \big(
      \mathrm{CE}(\mathfrak{l}A)
      \,,\,
      \Omega^\bullet_{\mathrm{dR}}(X)
    \big)
    &
    \mbox{by Prop. \ref{WhiteheadLInfinityAlgebras} \& Lem. \ref{PLdeRhamComplexOnSmoothManifoldsEquivalentToSmoothDeRhamComplex}}
    \\
    & \simeq
    H_{\mathrm{dR}}
    (
      X;
      \,
      \mathfrak{l}A
    )
    &
    \mbox{by Lem. \ref{NonAbelianDeRhamCohomologyViaThedgcHomotopyCategory}}.
  \end{array}
\end{equation}
The two lemmas invoked are proved next.
\end{proof}

\begin{lemma}[De Rham complex over cylinder of manifold is path space object]
  \label{DeRhamComplexOverCylinderManifoldIsPathSpaceObject}
  For $X \,\in\,\SmoothManifolds$,
  consider the following morphisms of
  dgc-algebras (Def. \ref{CategoryOfdgAlgebras})
  \vspace{-2mm}
  \begin{equation}
    \label{PathSpaceObjectFordeRhamComplex}
    \xymatrix{
      \Omega^\bullet_{\mathrm{dR}}(X)
      \ar[rr]^{ (\mathrm{pr}_X)^\ast }
      &&
      \Omega^\bullet_{\mathrm{dR}}\big(X \times \mathbb{R} \big)
      \ar[rr]^-{ ( i_0^\ast, \, i_1^\ast )   }
      &&
      \Omega^\bullet_{\mathrm{dR}}(X)
      \oplus
      \Omega^\bullet_{\mathrm{dR}}(X)
    }
  \end{equation}

  \vspace{-2mm}
  \noindent
  (from the de Rham complex of $X$ (Example \ref{SmoothdeRhamComplex})
  to that of its cylinder manifold $X \times \mathbb{R}$,
  to its Cartesian product with itself, by Example \ref{ProductAndCoproductAlgebras}),
  given by pullback of differential forms along these smooth functions:
  \vspace{-2mm}
  $$
    \xymatrix{
      X
      \ar@{<-}[rr]^-{ \mathrm{pr}_X }
      &&
      X \times \mathbb{R}
      \ar@{<-^{)}}[rr]^-{ (i_0, \, i_1) }
      &&
      \; \big(X \times \{0\}\big)
        \sqcup
      \big(X \times \{1\}\big)
      \;\simeq\;
      X \sqcup  X
      \,.
    }
  $$

  \vspace{-3mm}
  \noindent
  This is a path space object (Def. \ref{PathSpaceObject})
  for $\Omega^\bullet_{\mathrm{dR}}(X)$
  in $\dgcAlgebrasProj{\mathbb{R}}$ (Prop. \ref{ProjectiveModelStructureOndgcAlgebras}).
\end{lemma}
\begin{proof}
  {\bf (i)} It is clear by construction that the composite morphism
  is the diagonal.

\noindent  {\bf (ii)}
  That $(\mathrm{pr}_X)^\ast$ is a weak equivalence, hence a
  quasi-isomorphism, follows from the de Rham theorem, using
  that ordinary cohomology is homotopy invariant:
  $H^\bullet(X \times \mathbb{R}; \mathbb{R})
    \simeq H^\bullet(X; \mathbb{R})$.

  \noindent {\bf (iii)} That $( i_0^\ast, i_1^\ast )$ is a fibration, namely
  degreewise surjective, is seen from the fact that any pair of
  forms on the boundaries $X \times \{0\}$, $X \times \{1\}$
  may be smoothly interpolated to zero along any small
  enough positive parameter length,
  and then glued to a form on $X \times \mathbb{R}$.
\end{proof}

\begin{lemma}[Non-abelian de Rham cohomology via the dgc-homotopy category]
  \label{NonAbelianDeRhamCohomologyViaThedgcHomotopyCategory}
  Let $X \in \SmoothManifolds$ and
  $\mathfrak{g} \in \LInfinityAlgebrasNil$ (Def. \ref{NilpotentLInfinityAlgebras}).
  Then the non-abelian de Rham cohomology of $X$ with
  coefficients in $\mathfrak{g}$ (Def. \ref{NonabelianDeRhamCohomology})
  is in natural bijection with the hom-set in the
  homotopy category of $\dgcAlgebrasProj{\mathbb{R}}$ (Prop. \ref{ProjectiveModelStructureOndgcAlgebras})
  from $\mathrm{CE}(\mathfrak{g})$ (Def. \ref{ChevalleyEilenbergAlgebraOfLInfinityAlgebra})
  to $\Omega^\bullet_{\mathrm{dR}}(X)$ (Example \ref{SmoothdeRhamComplex}):
  \vspace{-2mm}
  \begin{equation}
    \label{BijectionBetweenNonabelianDeRhamCohomologyAndHomSetsIndgcHomotopyCategory}
    H_{\mathrm{dR}}
    \big(
      X;
      \,
      \mathfrak{g}
    \big)
    \;\simeq\;
    \mathrm{Ho}
    \big(
      \dgcAlgebrasProj{\mathbb{R}}
    \big)
    \big(
      \mathrm{CE}(\mathfrak{g})
      \,,\,
      \Omega^\bullet_{\mathrm{dR}}(X)
    \big).
  \end{equation}
\end{lemma}
\begin{proof}
  Consider a pair of dgc-algebra homomorphisms
  \vspace{-2mm}
  \begin{equation}
    \label{APairOfdgcHomomorphisms}
    A^{(0)}\!,\; A^{(1)}
    \;\in\;
    \dgcAlgebras{\mathbb{R}}
    \big(
      \mathrm{CE}(\mathfrak{g})
      \,,\,
      \Omega^\bullet_{\mathrm{dR}}(X)
    \big)
  \end{equation}

        \vspace{-2mm}
\noindent
  hence of flat $\mathfrak{g}$-valued differential forms,
  according to Def. \ref{FlatLInfinityAlgebraValuedDifferentialForms}.
  Observe that:

  \begin{itemize}
  \vspace{-.3cm}
  \item [\bf (i)] $\mathrm{CE}(\mathfrak{g})$ is cofibrant
  in
  $\dgcAlgebrasProj{\mathbb{R}}$ \eqref{ProjectiveModelStructureOnConnectivedgcAlgebras}.
  (by Prop. \ref{RelativeSullivanModelsAreCofibrations},
  and since $\mathfrak{g}$ is assumed to be nilpotent \eqref{CategoryOfNilpotentLInfinityAlgebras});

  \vspace{-.1cm}
  \item [\bf (ii)]
  $\Omega^\bullet_{\mathrm{dR}}(X)$ is fibrant
  in
  $\dgcAlgebrasProj{\mathbb{R}}$ \eqref{ProjectiveModelStructureOnConnectivedgcAlgebras}.
  (by Remark \ref{AlldgcAlgebrasAreProjectivelyFibrant});

  \vspace{-.1cm}
  \item[\bf (iii)]
  A right homotopy (Def. \ref{RightHomotopy})
  between the pair \eqref{APairOfdgcHomomorphisms} of morphisms,
  with respect to the
  path space object
  $\Omega^\bullet_{\mathrm{dR}}(X \times \mathbb{R})$
  from Lemma \ref{DeRhamComplexOverCylinderManifoldIsPathSpaceObject},
  namely a morphism $\widetilde A$ making the following diagram
  commute

  \vspace{-3mm}
  \begin{equation}
    \label{NonAbelianCoboundaryAsRightHomotopyOfdgcAlgebras}
    \raisebox{40pt}{
    \xymatrix@R=1.5em{
      \Omega^\bullet_{\mathrm{dR}}(X)
      \ar@{<-}[d]_-{
        i_0^\ast
      }
      \ar@{<-}[drr]^-{ A^{(0)} }
      \\
      \Omega^\bullet_{\mathrm{dR}}(X \times \mathbb{R})
      \ar[d]_-{
        i_1^\ast
      }
      \ar@{<-}[rr]|-{ \;\widetilde A\; }
      &&
      \mathrm{CE}(\mathfrak{g})
      \\
      \Omega^\bullet_{\mathrm{dR}}(X)
      \ar@{<-}[urr]_-{ A^{(1)} }
    }
    }
  \end{equation}

    \vspace{-2mm}
\noindent
  is manifestly the same as a coboundary
  $\widetilde A$ between the corresponding
  flat $\mathfrak{g}$-valued forms according to Def. \ref{CoboundariesBetweenFlatLInfinityAlgebraValuedForms}.
  \end{itemize}

  \vspace{-2mm}

  \noindent
  Therefore, Prop. \ref{RightHomotopyReflectsHomotopyClasses}
  says
  that the quotient set \eqref{NonabelianDeRhamCohomologyAsQuotientSet}
  defining the non-abelian de Rham cohomology
  is in natural bijection to the hom-set in the homotopy category.
\end{proof}

\begin{lemma}[PL de Rham complex on smooth manifold is equivalent to smooth de Rham complex]
  \label{PLdeRhamComplexOnSmoothManifoldsEquivalentToSmoothDeRhamComplex}
  Let $X$ be a smooth manifold. Then

  \noindent {\bf (i)} There exists a zig-zag
  of weak equivalences (Def. \ref{HomotopicalStructureOndgcAlgebras})
  in $\dgcAlgebrasProj{\mathbb{R}}$ \eqref{ProjectiveModelStructureOnConnectivedgcAlgebras}
  between the smooth de Rham complex of $X$ (Example \ref{SmoothdeRhamComplex})
  and the PL de Rham complex of its underlying topological space
  (Def. \ref{PLdeRhamComplexes}).

    \noindent {\bf (ii)}  In particular, both are isomorphic in the homotopy category:
    \vspace{-2mm}
  $$
    \mbox{\rm $X$ smooth manifold}
    \;\;\;\;\;
    \Rightarrow
    \;\;\;\;\;
    \Omega^\bullet_{\mathrm{dR}}(X)
    \;\simeq\;
    \Omega^\bullet_{\mathrm{PLdR}}(X)
    \;\;\;
    \in
    \mathrm{Ho}
    \Big(
      \dgcAlgebrasProj{\mathbb{R}}
    \Big).
  $$
\end{lemma}
\begin{proof}
  Let $\Omega^\bullet_{\mathrm{PSdR}}(-)$
  (for ``piecewise smooth'') be defined as
  the PL de Rham complex in Def. \ref{PLdeRhamComplexes},
  but with smooth differential forms on each simplex.
  Notice that this
  comes with the canonical natural inclusion
      \vspace{-3mm}
  $$
    \xymatrix{
      \Omega^\bullet_{\mathrm{PLdR}}(-)
      \;
      \ar@{^{(}->}[r]^-{ i_{\mathrm{poly}} }
      &
      \;
      \Omega^\bullet_{\mathrm{PSdR}}(-)
    }\,.
  $$

  \vspace{-2mm}
  \noindent Let then $\mathrm{Tr}(X) \in \SimplicialSets$
  be any smooth triangulation of $X$
  (Ex. \ref{HomotopyTypesOfManifoldsViaTriangulations}).
  This means that we have a homeomorphism
  out of its geometric realization
  \eqref{QuillenAdjunctionBetweenTopologicalSpacesAndSimplicialSets}
  to $X$

  \vspace{-.4cm}
  \begin{equation}
    \label{HomeomorphismFromRealizationOfTriangulation}
    \xymatrix@C=3em{
      \left\vert\mathrm{Tr}(X)\right\vert
      \ar[r]^-{p}_-{ \mathrm{homeo} }
      &
      X
      \,,
    }
  \end{equation}
  \vspace{-.5cm}

  \noindent
  which restricts on the interior of each simplex to a
  diffeomorphism onto its image; and that we have an inclusion
  \vspace{-2mm}
  \begin{equation}
    \label{InclusionOfTriangulationIntoSingularSimplicialSet}
   \xymatrix@C=4em{
     \mathrm{Tr}(X)
     \;
     \ar@{^{(}->}[r]^-{ \eta_{\mathrm{Tr}(X)} }_-{ \in\, \mathrm{W} }
     &
     \mathrm{Sing}
     \big(
       \left\vert
         \mathrm{Tr}(X)
       \right\vert
     \big)
     \ar[r]^-{
       \mathrm{Sing}(p)
     }_-{ \in\, \mathrm{Iso} }
     &
     \mathrm{Sing}(X)
     \,,
   }
  \end{equation}

  \vspace{-2mm}
  \noindent  which is a weak equivalence (by Example \ref{SimplicialSetsEquivalentToSingularSimplicialSetsOfTheirRealization}).
  In summary, this gives us the following zig-zag
  of dgc-algebra homomorphisms:
  $$
    \xymatrix@C=-2pt{
      &
      \Omega^\bullet_{\mathrm{PLdR}}
      \big(
        \mathrm{Tr}(X)
      \big)
      \ar[dr]^-{ i_{\mathrm{poly}} }
      &&
      \Omega^\bullet_{\mathrm{dR}}
      (
        X
      )
      \ar[dl]_-{ p^\ast }
      \\
      \mathllap{
        \Omega^\bullet_{\mathrm{PLdR}}(X)
        \;=\;
      }
      \Omega^\bullet_{\mathrm{PLdR}}
      \big(
        \mathrm{Sing}(X)
      \big)
      \ar[ur]^-{ (\eta_S)^\ast }
      &&
      \Omega^\bullet_{\mathrm{PSdR}}
      \big(
        \mathrm{Tr}(X)
      \big)
    }
  $$
  Here the two morphisms on the right are
  quasi-isomorphisms by
  {\cite[Cor. 9.9]{GriffithMorgan13}}
  (as in Prop. \ref{FundamentalTheoremForPiecewiseSmoothDeRhamComplex}).
  The morphism on the left is a quasi-isomorphism because
  $i$ is a weak homotopy equivalence \eqref{InclusionOfTriangulationIntoSingularSimplicialSet}
  and since $\Omega^\bullet_{\mathrm{PLdR}}$ preserves
  weak equivalences (by Lem. \ref{KenBrownLemma}),
  since it is a Quillen left adjoint
  (by Prop. \ref{QuillenAdjunctionBetweendgcAlgebrasAndSimplicialSets})
  and since every simplicial set is cofibrant (Ex. \ref{ClassicalModelStructureOnSimplicialSets}).
\end{proof}

\newpage

%%%%%%%%%%%%%%%%%%%%%%%%%%%%%%%%%%%%%%%%%%%%%%%%%%%%%%%
%\subsubsection{Twisted non-abelian de Rham theory}
%%%%%%%%%%%%%%%%%%%%%%%%%%%%%%%%%%%%%%%%%%%%%%%%%%%%%%%

\noindent {\bf Flat twisted $L_\infty$-algebra valued differential forms.}
We generalize the above discussion to include twistings.

\begin{defn}[Local $L_\infty$-algebraic coefficients]
  \label{LocalLInfinityAlgebraicCoefficients}
  We say that a \emph{local $L_\infty$-algebraic coefficient bundle}
  is a fibration
      \vspace{-3mm}
  \begin{equation}
    \xymatrix@R=8pt{
      \mathfrak{g}
      \ar[r]
      &
      \widehat{\mathfrak{b}}
      \ar[d]^-{\mathfrak{p}}
      \\
      & \mathfrak{b}
    }
    \label{Eq-bhat}
 \end{equation}

      \vspace{-2mm}
\noindent
  in $\LInfinityAlgebras$
  (Def. \ref{ChevalleyEilenbergAlgebraOfLInfinityAlgebra}), hence
  a morphism
  such that under passage to Chevalley-Eilenberg algebras \eqref{FullSubcategoryInclusionOfLieAlgebras} we have
  a cofibration
        \vspace{-2mm}
  \begin{equation}
    \label{CEAlgebraOfLInfinityAlgebraicCoefficientBundle}
    \raisebox{20pt}{
    \xymatrix@R=12pt@C=3em{
      \mathrm{CE}
      (
        \mathfrak{g}
      )
      \ar@{<-}[rr]^-{
        \mathrm{cofib}
        (
          \mathrm{CE}(\mathfrak{p})
        )
      }
      &&
      \mathrm{CE}
      \big(\,
        \widehat{\mathfrak{b}}
      \,\big)
      \ar@{<-}[d]_-{
        \mathrm{CE}(\mathfrak{p})
      }^-{
        \in\;\mathrm{Cof}
      }
      \\
      &&
      \mathrm{CE}(\mathfrak{b})
    }
    }
  \end{equation}

      \vspace{-2mm}
\noindent
  in $\dgcAlgebrasProj{\mathbb{R}}$ (Prop. \ref{ProjectiveModelStructureOndgcAlgebras}).
\end{defn}

In generalization of Def. \ref{FlatLInfinityAlgebraValuedDifferentialForms}, we say:

\begin{defn}[Flat twisted $L_\infty$-algebra valued differential forms]
  \label{FlatTwistedLInfinityAlgebraValuedDifferentialForms}
  $\,$

  \noindent
  {\bf (i)} Let $X \in \SmoothManifolds$ and $\widehat{\mathfrak{b}}$ \eqref{Eq-bhat}
%  $
%    \xymatrix@R=6pt{
%      \mathfrak{g} \ar[r]
%      &
%      \widehat{\mathfrak{b}}
%      \ar[d]^-{ \mathfrak{p} }
%      \\
%      &
%      \mathfrak{b}
%    }
%  $
  a local $L_\infty$-algebraic coefficient bundle
  (Def. \ref{LocalLInfinityAlgebraicCoefficients}).
  For
    \vspace{-1mm}
  \begin{equation}
    \label{TwistForTwistedDifferentialForms}
    \tau_{\mathrm{dR}}
    \;\in\;
    \Omega_{\mathrm{dR}}
    (
      X;
      \,
      \mathfrak{b}
    )_{\mathrm{flat}}
  \end{equation}

    \vspace{-1mm}
  \noindent
  a flat $\mathfrak{b}$-valued differential form on $X$
  (Def. \ref{FlatLInfinityAlgebraValuedDifferentialForms}),
  we say that
  a \emph{flat $\tau$-twisted $\mathfrak{g}$-valued differential form}
  on $X$ is a morphism of dgc-algebras (Def. \ref{CategoryOfdgAlgebras})
  in the slice over $\mathrm{CE}(\mathfrak{b})$
      \vspace{-3mm}
\noindent
  \begin{equation}
    \label{MorphismInSliceRepresentingFlatTwistedForms}
    \xymatrix@C=6em@R=1.5em{
      \Omega_{\mathrm{dR}}^\bullet(X)
      \ar@{<-}[dr]_-{\!\!
   %     \mathllap{
          \mbox{
            \tiny
            \color{darkblue}
            \bf
            twist
                      }
               }^-{
        \tau_{\mathrm{dR}}
      }
      \ar@{<-}[rr]^-{
          \mathclap{
          \raisebox{3pt}{
            \tiny
            \color{darkblue}
            \bf
            \def\arraystretch{.9}
            \begin{tabular}{c}
              flat $\tau_{\mathrm{dR}}$-twisted
              \\
              $\mathfrak{g}$-valued differential form
            \end{tabular}
          }
         }
      }_-{
        A
      }
      &&
      \mathrm{CE}
      \big(\,
        \widehat{\mathfrak{b}}
      \,\big)
      \ar@{<-^{)}}[dl]_>>>>>>>>>{
        \mathrm{CE}(\mathfrak{p})
        }^{
       % \mathrlap{
          \mbox{
            \tiny
            \color{darkblue}
            \bf
            \def\arraystretch{.9}
            \begin{tabular}{c}
              local
              \\
              $L_\infty$-algebraic
              \\
              coefficients
            \end{tabular}
  %        }
        }
      }
      \\
      &
      \mathrm{CE}(\mathfrak{b})
    }
  \end{equation}

      \vspace{-2mm}
  \noindent
  {\bf (ii)} We write
  $$
    \Omega^{\tau_{\mathrm{dR}}}_{\mathrm{dR}}
    (
      X;
      \,
      \mathfrak{g}
    )_{\mathrm{flat}}
    \;\;
    :=\;\;
    \big(
      \dgcAlgebras{\mathbb{R}}
    \big)_{\!\!/ \mathrm{CE}(\mathfrak{b}) }
    (
      \tau_{\mathrm{dR}}
      \,,\,
      \mathfrak{p}
    )
  $$
  for the set of all flat $\tau_{\mathrm{dR}}$-twisted
  $\mathfrak{g}$-valued differential forms on $X$.
\end{defn}

\begin{remark}[Underlying flat forms of flat twisted forms]
  \label{UnderlyingFlatFormsOfFlatTwistedForms}
  Let $X \in \SmoothManifolds$, let
  $
    \xymatrix@C=10pt{
      \mathfrak{g}
      \ar[r]
      &
      \widehat{\mathfrak{b}}
      \ar[r]^-{ \mathfrak{p} }
      &
      \mathfrak{b}
    }
  $
  be a local $L_\infty$-algebraic coefficient bundle
  (Def. \ref{LocalLInfinityAlgebraicCoefficients}),
  and let $\tau_{\mathrm{dR}} \,\in\,
    \Omega_{\mathrm{dR}}\big(X;\, \mathfrak{b} \big)
  $.
  Then there is a canonical forgetful natural transformation
  \vspace{-2mm}
  \begin{equation}
    \label{ForgetfulMapFromTwistedFlatFormsToFlatForms}
    \xymatrix{
      \Omega^{\tau_{\mathrm{dR}}}
      (
        X;
        \,
        \mathfrak{g}
      )_{\mathrm{flat}}
      \ar[r]
      &
      \Omega
      \big(
        X;
        \,
        \widehat{\mathfrak{b}} \;
      \big)_{\mathrm{flat}}
    }
  \end{equation}

      \vspace{-2mm}
\noindent
  from flat $\tau_{\mathrm{dR}}$-twisted $\mathfrak{g}$-valued
  differential forms (Def. \ref{FlatTwistedLInfinityAlgebraValuedDifferentialForms})
  to flat $\widehat{\mathfrak{b}}$-valued differential forms
  (Def. \ref{FlatLInfinityAlgebraValuedDifferentialForms}),
  given by remembering only the top morphism in
  \eqref{MorphismInSliceRepresentingFlatTwistedForms}.
\end{remark}

\begin{example}[$L_\infty$-coefficient bundle for $H_3$-twisted differential
forms {\cite[\S 4]{FSS16a}\cite[\S 4]{FSS16b}\cite[Lem. 2.31]{BMSS19}}]
  \label{RecoveringH3TwistedDifferentialForms}
  Consider the local $L_\infty$-algebraic coefficient bundle
  (Def. \ref{LocalLInfinityAlgebraicCoefficients})
  given by the following multivariate polynomial dgc-algebras
  (Def. \ref{MultivariatePolynomialdgcAlgebras}):
      \vspace{-3mm}
  $$
  \hspace{-1.5cm}
    \xymatrix@R=1.1em{
      \mathllap{
        \mathrm{CE}
        \big(
          \mathfrak{l}
          \mathrm{ku}_1
        \big)
        \;=\;
      }
      \mathbb{R}
      \!\!
      \left[
      \!\!
      {\begin{array}{c}
        \vdots
        \\[-3pt]
        f_5,
        \\[-3pt]
        f_3,
        \\[-3pt]
        f_1,
      \end{array}}
      \!\!
      \right]
      \!\Big/\!
      \left(
      {\begin{aligned}
        \vdots
        \\[-3pt]
        d\, f_5 & = 0
        \\[-3pt]
        d\, f_3 & = 0
        \\[-3pt]
        d\, f_1 & = 0
      \end{aligned}}
      \!\!\!\!\!\!\!
      \right)
      \ar@{<-}[rr]^-{
        \omega_{2k+1} \!\mapsfrom \, \omega_{2k+1}
      }
      &&
      \mathbb{R}
      \!\!
      \left[
      \!\!
      {\begin{array}{c}
        \vdots
        \\[-3pt]
        f_5,
        \\[-3pt]
        f_3,
        \\[-3pt]
        f_1,
        \\[-3pt]
        h_3
      \end{array}}
      \!\!
      \right]
      \!\Big/\!
      \left(
            {\begin{aligned}
        \vdots
        \\[-3pt]
        d\, f_5 & = h_3 \wedge f_3,
        \\[-3pt]
        d\, f_3 & = h_3 \wedge f_1,
        \\[-3pt]
        d\, f_1 & = 0,
        \\[-3pt]
        d\, h_3 & = 0
      \end{aligned}}
      \!\!\!\!\!\!\!\!
      \right)
      \mathrlap{
        \;
        =
        \mathrm{CE}
        \Big(
          \mathfrak{l}
          \big(
            \mathrm{ku}_1 \!\sslash\! B \mathrm{U}(1)
          \big)
        \Big)
      }
      \ar@{<-}[dd]^-{
        \scalebox{.7}{$
          \begin{array}{c}
            h_3
            \\
            \mapsup
            \\
            h_3
          \end{array}
        $}
      }
      \\
      \\
      &&
      \mathbb{R}
      \big[
        h_3
      \big]
      \big(
        d\,h_3 \;= 0
      \big)
      \mathrlap{
        \; = \;
        \mathrm{CE}
        \big(
          \mathfrak{b}^2 \mathbb{R}
        \big)
      }
    }
  $$

      \vspace{-2mm}
\noindent
  Here the rational model of the classifying space
  $\mathrm{ku}_1$ for complex topological K-theory in degree 1
  and for its twisted version is as in
  \cite[\S 4]{FSS16a}\cite[\S 4]{FSS16b}\cite[Lem. 2.31]{BMSS19}.
  In this case:

  \noindent
  {\bf (i)} A twist \eqref{TwistForTwistedDifferentialForms}
  is equivalently an ordinary closed 3-form form
  (by Example \ref{OrdinaryClosedFormsAreFlatLineLInfinityAlgebraValuedForms}):
      \vspace{-1mm}
  \begin{equation}
    \label{TwistDeRhamForlKU1}
    H_3 \;\in\;
    \Omega_{\mathrm{dR}}
    \big(
      X;
      \,
      \mathfrak{b}^2\mathbb{R}
    \big)_{\mathrm{flat}}
    \;\simeq\;
    \Omega^3_{\mathrm{dR}}(X)_{\mathrm{closed}}
    \,.
  \end{equation}

      \vspace{-1mm}
  \noindent
  {\bf (ii)} The flat $\tau_{\mathrm{dR}} \sim H_3$-twisted
  $\mathfrak{l}\mathrm{ku}_1$-valued
  differential forms according to Def. \ref{FlatTwistedLInfinityAlgebraValuedDifferentialForms}
  are equivalently sequences of odd-degree differential forms
  $F_{2k+1} \in \Omega^{2k+1}_{\mathrm{dR}}(X)$ satisfying the
  $H_3$-twisted de Rham closure condition (see \cite[(23)]{RohmWitten86}\cite{GS-RR}):
  \vspace{-2mm}
   \begin{equation}
    \label{SetOfH3TwistedlKU1ValuedForms}
    \Omega^{\tau_{\mathrm{dR}}}
    \big(
      X;
      \,
      \mathfrak{l}\mathrm{ku}_1
    \big)_{\mathrm{flat}}
    \;\;
    \simeq
    \;\;
    \bigg\{
    F_{2\bullet + 1} \in \Omega^{2\bullet + 1}_{\mathrm{dR}}
    \,\Big\vert\,\;
    d\, \underset{k}{\sum} F_{2k + 1}
    \;=\;
    H_3 \wedge \underset{k}{\sum} F_{2k - 1}
  %  \right.
    \bigg\}
  \end{equation}

  \vspace{-2mm}
\noindent
  (where we set $F_{2k - 1} := 0$ if $2k-1 < 0$, for convenience of notation).
\end{example}

In direct generalization of Example \ref{RecoveringH3TwistedDifferentialForms}, we have:

\begin{example}[$L_\infty$-coefficient bundle for higher twisted differential forms {\cite[Def. 2.14]{FSS18}}]
  \label{RecoveringHigherTwistedDifferentialForms}
  For $r \in \mathbb{N}$, $r \geq 1$,
  consider the local $L_\infty$-algebraic coefficient bundle
  (Def. \ref{LocalLInfinityAlgebraicCoefficients})
  given by the following multivariate polynomial dgc-algebras
  (Def. \ref{MultivariatePolynomialdgcAlgebras}):
  \vspace{-2mm}
  \begin{equation}
  \label{LInfinityCoefficientBundleForHigherTwistedForms}
  \hspace{-1.5cm}
  \raisebox{30pt}{
    \xymatrix@R=1.1em@C=3em{
      \mathrm{CE}
      \Big(
        \underset{k \in \mathbb{N}}{\oplus} \mathfrak{b}^{2rk}\mathbb{R}
      \Big)
      \ar@{=}[d]
      &&
      \mathrm{CE}
      \Big(
        \Big(
          \underset{k \in \mathbb{N}}{\oplus} \mathfrak{b}^{2rk}\mathbb{R}
        \Big)
          \!\sslash\!
          B^{2r-1} \mathrm{U}(1)
      \Big)
      \ar@{=}[d]
      \\
      \mathbb{R}
      \!\!
      \left[
      \!\!
      {\begin{array}{c}
        \vdots
        \\[-3pt]
        f_{4r+1},
        \\[-3pt]
        f_{2r+1},
        \\[-3pt]
        f_1,
      \end{array}}
      \!\!
      \right]
      \!\Big/\!
      \left(
      {\begin{aligned}
        \vdots
        \\[-3pt]
        d\, f_{4r+1} & = 0
        \\[-3pt]
        d\, f_{2r+1} & = 0
        \\[-3pt]
        d\, f_1\;\; & = 0
      \end{aligned}}
      \!\!\!\!\!\!\!
      \right)
      \ar@{<-}[rr]^-{
        f_{2rk+1} \, \mapsfrom  \, f_{2rk+1}
      }
      &&
      \mathbb{R}
      \!\!
      \left[
      \!\!
      {\begin{array}{c}
        \vdots
        \\[-3pt]
        f_{4 r + 1},
        \\[-3pt]
        f_{2r + 1},
        \\[-3pt]
        f_1,
        \\
        h_{2r + 1}
      \end{array}}
      \!\!
      \right]
      \!\Big/\!
      \left(
            {\begin{aligned}
        \vdots
        \\[-3pt]
        d\, f_{4r+1} & = h_{2r+1} \wedge f_{2r+1},
        \\[-3pt]
        d\, f_{2r+1} & = h_{2r+1} \wedge f_1,
        \\[-3pt]
        d\, f_1 \;\; & = 0,
        \\[-3pt]
        d\, h_{2r+1} & = 0
      \end{aligned}}
      \!\!\!\!\!\!\!\!
      \right)
      \ar@{<-}[dd]^-{
        \scalebox{.7}{$
          \begin{array}{c}
            h_{2r+1}
            \\
            \mapsup
            \\
            h_{2r+1}
          \end{array}
        $}
      }
      \\
      \\
      &&
      \mathbb{R}
      \big[
        h_{2r+1}
      \big]
      \big(
        d\,h_{2r+1} \;= 0
      \big)
      \\
      &&
      \mathrm{CE}
      \big(
        \mathfrak{b}^{2r} \mathbb{R}
      \big)
      \ar@{=}[u]
    }
    }
  \end{equation}

  \vspace{-2mm}
  \noindent
  In this case:

  \noindent
  {\bf (i)} A twist \eqref{TwistForTwistedDifferentialForms}
  is equivalently an ordinary closed $(2r+1)$-form form
  (by Example \ref{OrdinaryClosedFormsAreFlatLineLInfinityAlgebraValuedForms}):
  \vspace{-2mm}
  \begin{equation}
    \label{TwistDeRhamForlK2rKU1}
    H_{2r+1} \;\in\;
    \Omega_{\mathrm{dR}}
    \big(
      X;
      \,
      \mathfrak{b}^{2r}\mathbb{R}
    \big)_{\mathrm{flat}}
    \;\simeq\;
    \Omega^{2r+1}_{\mathrm{dR}}(X)_{\mathrm{closed}}
    \,.
  \end{equation}

  \vspace{-2mm}
  \noindent
  {\bf (ii)} The flat $\tau_{\mathrm{dR}} \sim H_{2r+1}$-twisted
  $\underset{k \in \mathbb{N}}{\oplus} \mathfrak{b}^{2rk}\mathbb{R}$-valued
  differential forms according to Def. \ref{FlatTwistedLInfinityAlgebraValuedDifferentialForms}
  are equivalently sequences of differential forms
  $F_{2r\bullet +1} \in \Omega^{2k\bullet +1}_{\mathrm{dR}}(X)$ satisfying the
  $H_{(2r+1)}$-twisted de Rham closure condition
  \eqref{HigherTwistedDeRhamCohomology}:
   \vspace{-2mm}
  \begin{equation}
    \label{SetOfHigherTwistedlKKU1ValuedForms}
    \Omega^{\tau_{\mathrm{dR}}}
    \big(
      X;
      \,
      \underset{k \in \mathbb{N}}{\oplus} \mathfrak{b}^{2rk}\mathbb{R}
    \big)_{\mathrm{flat}}
    \;\;
    \simeq
    \;\;
    \bigg\{
    F_{2r\bullet + 1} \in \Omega^{2r\bullet + 1}_{\mathrm{dR}}
    \,\Big\vert\,\;
    d\, \underset{k}{\sum} F_{2rk + 1}
    \;=\;
    H_{2r+1} \wedge \underset{k}{\sum} F_{2rk - 1}
    \bigg\}
  \end{equation}

\vspace{-2mm}
\noindent  (where we set $F_{2rk - 1} := 0$ if $2rk-1 < 0$, for convenience of notation).
\end{example}

In twisted generalization of Example \ref{FlatSphereValuedDifferentialForms},
we have the following:

\begin{example}[Flat twisted differential forms with values in
Whitehead $L_\infty$-algebras of spheres and twistor space]
  \label{FlatTwistedDifferentialFormsWithValuesInTwistorSpace}
  The $L_\infty$-algebraic local coefficient bundles
  (Def. \ref{LocalLInfinityAlgebraicCoefficients})
  given as the relative Whitehead $L_\infty$-algebras
  (Prop. \ref{WhiteheadLInfinityAlgebrasRelative})
  of the local coefficient bundles \eqref{EquivariantizedTwistorFibration}
  for twisted and twistorial Cohomotopy
  (Example \ref{TwistorialCohomotopy})
  are as shown on the right of the following diagram
  \cite[Lem. 3.19]{FSS19b}\cite[Thm. 2.14]{FSS20a}:
  \vspace{-2mm}
  $$
      \xymatrix@C=-13pt@R=2.5em{
      &&
      \mathrm{CE}
      \big(
        \mathfrak{l}_{\scalebox{.5}{$B \mathrm{Sp}(2)$}}
        ( \mathbb{C}P^3 \!\sslash\! \mathrm{Sp}(2))
      \big)
      \ar@{<-}[d]^-{ (t_{\mathbb{H}} \sslash \mathrm{Sp}(2) )^\ast }
      \ar@{}[rr]|-{=}
      &&
      {
        \mathrm{CE}(\mathfrak{l}B \mathrm{Sp}(2))
        \!
        \left[
        \!\!\!
        \begin{array}{l}
          h_3,
          \\[-3pt]
          f_2,
          \\[-3pt]
          \omega_7,
          \\[-3pt]
          \omega_4
        \end{array}
        \!\!\!
        \right]
        \!\big/\!
        \left(
        \begin{aligned}
          d \, h_3 & = \omega_4 - \tfrac{1}{4}p_1 - f_2 \wedge f_2
          \\[-3pt]
          d \, f_2 & = 0
          \\[-3pt]
          d\, \omega_7
          & =
          - \omega_4 \wedge \omega_4
          + (\tfrac{1}{4}p_1)^2 - \rchi_8
          \\[-3pt]
          d\, \omega_4 & = 0
        \end{aligned}
       \!\! \right)
      }
      \ar@{<-^{)}}[d]
      \\
      \Omega^\bullet_{\mathrm{dR}}(X)
      \ar@{<--}[rr]^-{
        (G_4, 2G_7)
      }
      \ar@{<-}[dr]_-{ \tau_{\mathrm{dR}} }
      \ar@{<--}[urr]^-{
        (G_4, G_7, F_2, H_3) \;\;
      }
      %\ar@{<--}[uurr]|>>>>>>>>>>>>>>>>>>>>>>{
      %  (G_4, G_7, F_2, H_3, H_1)
      %}
      &&
      \mathrm{CE}
      \big(
        \mathfrak{l}_{\scalebox{.5}{$B \mathrm{Sp}(2)$}}
        ( S^4 \!\sslash\! \mathrm{Sp}(2))
      \big)
      \ar@{<-}[dl]
      \ar@{}[rr]|-{=}
      &
      {\phantom{AAA}}
      &
      {
        \mathrm{CE}(\mathfrak{l}B \mathrm{Sp}(2))
        \!
        \left[
        \!\!\!
        \begin{array}{l}
          \omega_7,
          \\[-3pt]
          \omega_4
        \end{array}
        \!\!\!
        \right]
        \!\big/\!
        \left(
        \begin{aligned}
          d\, \omega_7
          & =
          - \omega_4 \wedge \omega_4
          + (\tfrac{1}{4}p_1)^2 - \rchi_8
          \\[-3pt]
          d\, \omega_4 & = 0
        \end{aligned}
        \right)
      }
      \ar@{<-^{)}}[dl]
      \\
      &
      \mathrm{CE}
      \big(
        \mathfrak{l}
        B \mathrm{Sp}(2)
      \big)
      \ar@{}[r]|-{=}
      &
      \!\!\!\!\!\!\!\!\!\!\!\!\!\!\!
      \mathrlap{
      \mathbb{R}
      \left[
      \!\!\!
      {\begin{array}{l}
        \phantom{\tfrac{1}{2}}\rchi_8,
        \\[-3pt]
        \tfrac{1}{2}p_1
      \end{array}}
      \!\!\!
      \right]
      \!\big/\!
      \left(
        {\begin{aligned}
          d\, \phantom{\tfrac{1}{2}}\rchi_8 & = 0
          \\[-3pt]
          d\, \tfrac{1}{2}p_1 & = 0
        \end{aligned}}
      \right)
      }
      &
      {\phantom{\vert^{\vert^{\vert^{\vert^{\vert^{\vert^{\vert^{\vert^{\vert^{\vert^{\vert^{\vert^{\vert^{\vert^{\vert^{\vert^{\vert^{\vert^{\vert^{\vert^{\vert^{\vert^{\vert^{\vert^{\vert^{\vert^{\vert^{\vert}}}}}}}}}}}}}}}}}}}}}}}}}}}}}
    }
  $$

  \vspace{-2mm}
  \noindent
  Therefore, given a smooth 8-dimensional spin-manifold $X$
  equipped with tangential $\mathrm{Sp}(2)$-structure
  $\tau$ \eqref{TangentialSp2Structure},
  the flat $\tau_{\mathrm{dR}}$-twisted
  $\mathfrak{l}S^4$-
  and
  $\mathfrak{l}\mathbb{C}P^3$-valued differential forms
  (Def. \ref{FlatTwistedLInfinityAlgebraValuedDifferentialForms})
  are of the following form
  \cite[Prop. 3.20]{FSS19b}\cite[Prop. 3.9]{FSS20a}:

 \vspace{-3mm}
  \begin{equation}
    \label{TwistorialCohomotopyFormData}
    \begin{aligned}
    \Omega_{\mathrm{dR}}^{\tau_{\mathrm{dR}}}
    \big(
      X;
      \,
      \mathfrak{l}S^4
    \big)
    &
    =\;
    \left\{
    \!\!\!
    \begin{array}{l}
      2G_7,
      \\[-3pt]
      \phantom{2}G_4
    \end{array}
    \!\!\!
    \in
    \Omega^\bullet_{\mathrm{dR}}(X)
    \,\left\vert\,
    \begin{aligned}
      d\, 2G_7 & =
        -
        \big(
          G_4
          -
          \tfrac{1}{4}p_1(\nabla)
        \big)
        \wedge
        \big(
          G_4
          +
          \tfrac{1}{4}p_1(\nabla)
        \big)
        -
        \rchi_8(\nabla),
      \\[-3pt]
      d\, \phantom{2}G_4 & = 0
    \end{aligned}
    \right.
    \right\}
    \\
        \Omega_{\mathrm{dR}}^{\widetilde \tau_{\mathrm{dR}}}
    \big(
      X;
      \,
      \mathfrak{l}\mathbb{C}P^3
    \big)
    &
    =\;
    \left\{
    \!\!\!
    \begin{array}{l}
      \phantom{2}H_3,
      \\[-3pt]
      \phantom{2}F_2,
      \\
      2G_7,
      \\
      \phantom{2}G_4
    \end{array}
    \!\!\!
    \in
    \Omega^\bullet_{\mathrm{dR}}(X)
    \,\left\vert\,
    \begin{aligned}
      d\, H_3 & = G_4 - \tfrac{1}{4}p_1(\nabla) - F_2 \wedge F_2,
      \\[-3pt]
      d\, F_2 & = 0,
      \\[-3pt]
      d\, 2G_7 & =
        -
        \big(
          G_4
          -
          \tfrac{1}{4}
        \big)
        \wedge
        \big(
          G_4
          +
          \tfrac{1}{4}
        \big)
        -
        \rchi_8(\nabla),
      \\[-3pt]
      d\, \phantom{2}G_4 & = 0,
    \end{aligned}
    \right.
    \right\}
    \end{aligned}
  \end{equation}
  Notice:

  \noindent
  {\bf (a)}
  Here we are using (Example \ref{DeRhamRepresentativeOfTangentialSp2Twist})
  that the de Rham image $\tau_{\mathrm{dR}}$
  of the rationalization $L_{\mathbb{R}}\tau$ of the twist $\tau$ is given
  by evaluating characteristic forms (Def. \ref{CharacteristicForms})
  on any $\mathrm{Sp}(2)$-connection $\nabla$.

  \noindent
  {\bf (b)}
  In the second equation of \eqref{TwistorialCohomotopyFormData}
  we are using the above minimal model
  for $\mathbb{C}P^3 \!\sslash\! \mathrm{Sp}(2) $
  relative to
  $S^4 \!\sslash\! \mathrm{Sp}(2)$
  (instead of relative to $B \mathrm{Sp}(2)$).
\end{example}

\medskip
\noindent {\bf Twisted non-abelian de Rham cohomology.}
In generalization of Def. \ref{CoboundariesBetweenFlatLInfinityAlgebraValuedForms}, we set:

\begin{defn}[Coboundaries between flat twisted $L_\infty$-algebraic forms]
  \label{CoboundariesBetweenFlatTwistedLInfinityAlgebraValuedForms}
  Let $X \in \SmoothManifolds$, let
  $\!\!
    \xymatrix@C=9pt{
      \mathfrak{g} \ar[r]
      &
      \widehat{\mathfrak{b}}
      \ar[r]^-{ \mathfrak{p} }
      &
      \mathfrak{b}
    }
  $
  be a local $L_\infty$-algebraic coefficient bundle
  (Def. \ref{LocalLInfinityAlgebraicCoefficients}),
  and let $\tau_{\mathrm{dR}} \,\in\,
    \Omega_{\mathrm{dR}}(X;\, \mathfrak{b})
  $.
    Then for
        \vspace{-2mm}
  $$
    A^{(0)}\!,\, A^{(1)}
    \;\in\;
    \Omega^{\tau_{\mathrm{dR}}}_{\mathrm{dR}}
    (
      X; \mathfrak{g}
    )
  $$

  \vspace{-2mm}
  \noindent
  a pair of flat $\tau_{\mathrm{dR}}$-twisted
  $\mathfrak{g}$-valued differential forms on $X$
  (Def. \ref{FlatTwistedLInfinityAlgebraValuedDifferentialForms})
  a \emph{coboundary} between them is a coboundary
  \vspace{-2mm}
  \begin{equation}
    \label{CoboundaryBetweenFlatTwistedForms}
    \widetilde A
    \;\in\;
    \Omega_{\mathrm{dR}}
    \big(
      X \times \mathbb{R};
      \,
      \widehat{\mathfrak{b}}
    \, \big)
  \end{equation}

      \vspace{-1mm}
\noindent
  in the sense of Def. \ref{CoboundariesBetweenFlatLInfinityAlgebraValuedForms}
  between the underlying
  flat $\widehat{\mathfrak{b}}$-valued forms (via Remark \ref{UnderlyingFlatFormsOfFlatTwistedForms}),
  such that the underling $\mathfrak{b}$-valued form
  of $H$ equals the pullback of the twist $\tau_{\mathrm{dR}}$
  along $\xymatrix@C=14pt{X \times \mathbb{R} \ar[r]^-{\mathrm{pr}_X} & X }$

  \vspace{-3mm}
 \begin{equation}
    \label{ConditionOnCoboundaryBetweenFlatTwistedForms}
    \mathfrak{p}_\ast(H)
    \;=\;
    \mathrm{pr}_X^\ast(\tau_{\mathrm{dR}})
    \,.
  \end{equation}

      \vspace{-1mm}
  \noindent
  If such a coboundary exists, we say that $A^{(0)}$
  and $A^{(1)}$ are \emph{cohomologous}, to be denoted
      \vspace{-2mm}
  $$
    A^{(0)}
    \;\sim\;
    A^{(1)}
    \,.
  $$
\end{defn}

In generalization of Def. \ref{NonabelianDeRhamCohomology}, we set:

\begin{defn}[Twisted non-abelian de Rham cohomology]
  \label{TwistedNonabelianDeRhamCohomology}
  Let $X \in \SmoothManifolds$, let
  $\!\!
    \xymatrix@C=12pt{
      \mathfrak{g} \ar[r]
      &
      \widehat{\mathfrak{b}}
      \ar[r]^-{ \mathfrak{p} }
      &
      \mathfrak{b}
    \!\!}
  $
  be a local $L_\infty$-algebraic coefficient bundle
  (Def. \ref{LocalLInfinityAlgebraicCoefficients}),
  and let $\tau_{\mathrm{dR}} \,\in\,
    \Omega_{\mathrm{dR}}\big(X;\, \mathfrak{b} \big)
  $.
  Then the \emph{$\tau_{\mathrm{dR}}$ -twisted non-abelian de Rham cohomology}
  of $X$ with coefficients in $\mathfrak{g}$ is the set
  \vspace{-1mm}
  \begin{equation}
    \label{TwistedNonabelianDeRhamCohomologyAsQuotientSet}
    H^{\tau_{\mathrm{dR}}}_{\mathrm{dR}}
    (
      X;
      \,
      \mathfrak{g}
    )
    \;:=\;
    \big(
    \Omega^{\tau_{\mathrm{dR}}}_{\mathrm{dR}}
    (
      X;
      \,
      \mathfrak{g}
    )_{\mathrm{flat}}
    \big)_{\!\!/\sim}
  \end{equation}

  \vspace{-1mm}
  \noindent
  of equivalence classes with respect to the
  coboundary relation from Def. \ref{CoboundariesBetweenFlatTwistedLInfinityAlgebraValuedForms}
  on the set of flat $\tau_{\mathrm{dR}}$-twisted
  $\mathfrak{g}$-valued differential forms
  on $X$
  (Def. \ref{FlatTwistedLInfinityAlgebraValuedDifferentialForms}).
\end{defn}
\begin{remark}[Independence of the representative of the twist]
  \label{IndependenceOfTheRepresentativeOfTheTwist}
  The twisted non-abelian de Rham theorem below (Thm. \ref{TwistedNonAbelianDeRhamTheorem})
  makes manifest that
  the twisted non-abelian de Rham cohomology in
  Def. \ref{TwistedNonabelianDeRhamCohomology} depends
  on the twist $\tau_{\mathrm{dR}}$ only through its
  class
  $
    [\tau_{\mathrm{dR}}]
    \,\in\,
    H_{\mathrm{dR}}(X;\,\mathfrak{b})
  $
  in (un-twisted) non-abelian de Rham cohomology (Def. \ref{NonabelianDeRhamCohomology}).
\end{remark}

\noindent
{\bf The example of traditional twisted de Rham cohomology.}
Twisted de Rham cohomology is traditionally familiar
in the form of degree-3 twisted cohomology of even/odd degree differential
forms
\cite[\S III, Appendix]{RohmWitten86}\cite[\S 9.3]{BCMMS02}\cite[\S 3]{MathaiStevenson03}\cite[\S 2]{FHT-complex}\cite[Prop. 3.7]{Teleman04}\cite[\S I.4]{Cavalcanti05}\cite{tcu}\cite{MathaiWu}\cite{GS-HigherDeligne}
(which is the target of the twisted Chern character in
degree-3 twisted K-theory, see Prop. \ref{TwistedChernCharacterInTwistedTopologicalKTheory}).

\medskip
We discuss now how this archetypical example (Def. \ref{Degree3TwistedAbelianDeRhamCohomology})
and its higher-degree generalization (Def. \ref{HigherTwistedAbelianDeRhamCohomology})
are subsumed by our general Def. \ref{TwistedNonabelianDeRhamCohomology}.

\begin{defn}[Degree-3 twisted abelian de Rham cohomology]
  \label{Degree3TwistedAbelianDeRhamCohomology}
  For $X \in \SmoothManifolds$, and
  $H_3 \in \Omega^3_{\mathrm{dR}}(X)_{\mathrm{closed}}$
  a closed differential 3-form,
  the \emph{$H_3$-twisted de Rham cohomology} of $X$
  is the
  cochain cohomology
  \footnote{The notation ``$H_3$'' for the twist (and of ``$H_{2r+1}$'' for the higher twists later) originates in the physics literature and has made
  it as a convention in differential geometry as well.
  Not to be confused with a third homology group, of course}
  \vspace{-2mm}
  \begin{equation}
    \label{H3TwistedDeRhamCohomology}
    H^{\bullet+H_3}_{\mathrm{dR}}(X)
    \;:=\;
    \frac{
      \mathrm{ker}^\bullet\big( d - H_3 \wedge (-)\big)
    }{
      \mathrm{im}^\bullet\big( d - H_3 \wedge (-)\big)
    }
  \end{equation}

  \vspace{-3mm}
  \noindent
  of the following 2-periodic cochain complex:
  \vspace{-2mm}
  \begin{equation}
    \label{ChainComplexForH3TwisteddeRhamCohomology}
    \xymatrix{
      \cdots
      \ar[r]
      &
      \underset{k}{\bigoplus}
      \,
      \Omega^{(n-1)+2k}_{\mathrm{dR}}(X)
      \ar[rr]^-{
        (d - H_3 \wedge (-))
      }
      &&
      \underset{k}{\bigoplus}
      \,
      \Omega^{n + 2k}_{\mathrm{dR}}(X)
      \ar[rr]^-{
        (d - H_3 \wedge (-))
      }
      &&
      \underset{k}{\bigoplus}
      \,
      \Omega^{(n+1) + 2k}_{\mathrm{dR}}(X)
      \ar[r]
      &
      \cdots
      \,.
    }
  \end{equation}
\end{defn}

We show that this is a special case of twisted
non-abelian de Rham cohomology according to Def. \ref{TwistedNonabelianDeRhamCohomology}:

\begin{prop}[Twisted non-abelian de Rham cohomology subsumes $H_3$-twisted abelian de Rham cohomology]
  \label{TwistedNonabelianDeRhamCohomologySubsumesH3TwisteddeRhamCohomology}
  Given a twisting 3-form as in \eqref{TwistDeRhamForlKU1}

  \vspace{-5mm}
  $$
    \xymatrix@R=6pt{
      \tau_{\mathrm{dR}}
      \ar@{<->}[r]
      \ar@{}[d]|-{
        \rotatebox[origin=c]{-90}{$\in$}
      }
      &
      H_3
      \ar@{}[d]|-{
        \rotatebox[origin=c]{-90}{$\in$}
      }
      \\
      \Omega
      \big(
        X;
        \,
        \mathfrak{b}^2 \mathbb{R}
      \big)_{\mathrm{flat}}
      \ar@{}[r]|-{\simeq}
      &
      \Omega^3
      (
        X
      )_{\mathrm{closed}}
    }
  $$

  \vspace{-2mm}
  \noindent
  the $\tau_{\mathrm{dR}}$-twisted non-abelian de Rham cohomology
  (Def. \ref{TwistedNonabelianDeRhamCohomology}) of
  flat twisted $\mathfrak{l}\mathrm{ku}_1$-valued differential forms
  (Example \ref{RecoveringH3TwistedDifferentialForms})
  is naturally equivalent to
  $H_3$-twisted abelian de Rham cohomology (Def. \ref{Degree3TwistedAbelianDeRhamCohomology})
  in odd degree\footnote{The discussion for even degrees is
  directly analogous and we omit it for brevity.}

  \vspace{-2mm}
  $$
  \overset{
  \mathclap{
      \raisebox{3pt}{
        \tiny
        \color{darkblue}
        \bf
        \def\arraystretch{.9}
        \begin{tabular}{c}
          $\mathfrak{b}^2 \mathbb{R}$-twisted
          $\mathfrak{l}\mathrm{ku}_1$-valued
          \\
          non-abelian de Rham cohomology
          \\
        \end{tabular}
      }
      }
      }{
      H^{\tau_{\mathrm{dR}}}_{\mathrm{dR}}
    (
      X;
      \,
      \mathfrak{l} \mathrm{ku}_1
    )
    }
            \;\;\;\;\simeq\;\;\;\;
        \overset{
      \mathclap{
       \raisebox{3pt}{
        \tiny
        \color{darkblue}
        \bf
        \def\arraystretch{.9}
        \begin{tabular}{c}
          traditional $H_3$-twisted
          \\
          de Rham cohomology
        \end{tabular}
        }
        }
        }{
           H^{1 + H_3}_{\mathrm{dR}}(X)
      }
              $$
\end{prop}

\begin{proof}
By \eqref{SetOfH3TwistedlKU1ValuedForms} in
Example \ref{RecoveringH3TwistedDifferentialForms} the
cocycle sets on both sides are in natural bijection. Hence it
is sufficient to see that the coboundary relations
on the cocycle sets coincide, under this identification.
In one direction, consider a coboundary in the sense of
twisted non-abelian de Rham cohomology
(Def. \ref{CoboundariesBetweenFlatTwistedLInfinityAlgebraValuedForms})
with coefficients as in Example \ref{RecoveringH3TwistedDifferentialForms}:
\vspace{-1mm}
$$
  \widetilde F_{2\bullet+1}
  \;\in\;
  \Omega_{\mathrm{dR}}
  \big(
    X \times \mathbb{R};
    \,
    \mathfrak{l}
    \mathrm{ku}_1
  \big).
$$
\vspace{-2mm}
  We claim that

  \vspace{-2mm}
  \begin{equation}
  \label{FromNonabalianCoboundaryToH3TwistedCoboundary}
  h_{2\bullet} \;:=\; \underset{[0,1]}{\int} \widetilde F_{2\bullet+1}
\end{equation}

\vspace{-2mm}
  \noindent
satisfies the coboundary condition \eqref{H3TwistedDeRhamCohomology}:
\vspace{-2mm}
\begin{equation}
  \label{CoboundaryInHrTwistedDeRhamCohomology}
  \big(
    d - H_3 \wedge
  \big)
  \underset{k}{\sum}
  h_{2k}
  \;=\;
  \underset{k}{\sum}
  \big(
    F^{(1)}_{2k + 1}
    -
    F^{(0)}_{2k + 1}
  \big)
  \,.
\end{equation}

\vspace{-2mm}
\noindent To see this, we may compute as follows:
\vspace{-3mm}
$$
  \begin{aligned}
    d \,
      \underset{k}{\sum}
      h_{2k}
      & = \;
      \underset{k}{\sum}
      \Bigg(
      F_{2k+1}^{(1)}
      -
      F_{2k+1}^{(0)}
      -
      \underset{[0,1]}{\int} d \widetilde F_{2k+1}
      \Bigg)
      \\[-2pt]
      &
      =
      \underset{k}{\sum}
      \Bigg(
      F_{2k+1}^{(1)}
      -
      F_{2k+1}^{(0)}
      -
      \underset{[0,1]}{\int}
      \big(
        \mathrm{pr}_X^\ast H_3
      \big)
      \wedge
      \widetilde F_{2k - 1}
      \Bigg)
      \\[-2pt]
      &
      =
      \underset{k}{\sum}
      \Bigg(
      F_{2k+1}^{(1)}
      -
      F_{2k+1}^{(0)}
      +
      H_3
      \wedge
      \underset{[0,1]}{\int}
         \widetilde F_{2k - 1}
      \Bigg)
      \\[-2pt]
      & =
      \underset{k}{\sum}
      \left(
      F_{2k+1}^{(1)}
      -
      F_{2k+1}^{(0)}
      +
      H_3
      \wedge
      h_{2 k - 2}
      \right)
      \,,
  \end{aligned}
$$

\vspace{-3mm}
  \noindent
where the first step is the fiberwise Stokes formula \eqref{StokesFormula}
together with the defining restrictions \eqref{RestrictionOfFormsOnCylinderToBoundary}
of $\widetilde F_{2\bullet+1}$;
the second step is the
cocycle condition \eqref{SetOfH3TwistedlKU1ValuedForms} on $\widetilde F_{2\bullet + 1}$ using the constraint \eqref{ConditionOnCoboundaryBetweenFlatTwistedForms};
the third step is the projection formula \eqref{ProjectionFormula};
and the last step uses again the definition \eqref{FromNonabalianCoboundaryToH3TwistedCoboundary}.

Conversely, given $h_{2 \bullet}$ satisfying
\eqref{CoboundaryInHrTwistedDeRhamCohomology}, we claim that
\vspace{-2mm}
\begin{equation}
  \label{FromH3TwistedCobound}
  \widetilde F_{2\bullet + 1}
  \;;=\;
  (1-t) \, \mathrm{pr}^\ast_1 \big( F_{2\bullet + 1}^{(0)} \big)
  +
  t \, \mathrm{pr}^\ast_1 \big( F_{2\bullet + 1}^{(1)} \big)
  +
  d t
  \wedge
  \mathrm{pr}_X^\ast(h_{2\bullet})
\end{equation}

\vspace{-1mm}
  \noindent
is a coboundary of twisted non-abelian cocycles, in the sense of
Def. \ref{CoboundariesBetweenFlatTwistedLInfinityAlgebraValuedForms}:
It is immediate that \eqref{FromH3TwistedCobound}
has the required restrictions \eqref{RestrictionOfFormsOnCylinderToBoundary}.
We check by direct computation
that it satisfies the required differential equation:
\vspace{-1mm}
  $$
  \begin{aligned}
    d\,
    \underset{k}{\sum}
    \widetilde F_{2k + 1}
    & =
    \underset{k}{\sum}\Big(
    - dt \wedge  \mathrm{pr}_X^\ast \big( F_{2k + 1}^{(0)} \big)
    +
    (1-t) \,
    \mathrm{pr}_X^\ast\big( H_3 \big) \wedge \mathrm{pr}_X^\ast \big( F_{2k - 1}^{(0)} \big)
    \\[-2pt]
    &
    \phantom{= \underset{k}{\sum}\Big( }\;
    +
    d t \wedge
    \mathrm{pr}_X^\ast \big( F_{2k + 1}^{(1)} \big)
    + t\,
    \mathrm{pr}_X^\ast\big( H_3 \big)
    \wedge \mathrm{pr}_X^\ast \big(  F_{2k - 1}^{(1)} \big)
    \\[-2pt]
    &
    \phantom{= \underset{k}{\sum}\Big( }\;
    -
    d t
    \wedge
    \mathrm{pr}_X^\ast
    \big(\,
      \underset{
        \mathclap{
        =
        F_{2k + 1}^{(1)} - F_{2k + 1}^{(0)}
        +
        H_3 \wedge \, h_{2k}
        }
      }{
      \underbrace{
        d\, h_{2k}
      }
      }
   \, \big)
    \;\;
    \Big)
    \\[-2pt]
    & =
    \underset{k}{\sum}
    \left(
    \mathrm{pr}_X^\ast
    (
      H_3
    )
      \wedge
    \widetilde F_{2k - 1}
    \right)
    .
  \end{aligned}
$$

\vspace*{-1.8\baselineskip}
\end{proof}

In generalization of Def. \ref{Degree3TwistedAbelianDeRhamCohomology},
there are twisted abelian Rham complexes
with twist any odd-degree closed form
\cite[\S 3]{Teleman04}\cite{Sa-Higher}\cite{MathaiWu}\cite{tcu}\cite{GS-HigherDeligne}
(these serve as the targets
for the Chern character \cite{MMS20} on higher-twisted
ordinary K-theory
\cite{Teleman04}\cite{Gomez}\cite{DP13}\cite{Pennig},
see Example \ref{ChernCharacterOnHigherCohomotopicallyTwistedKTheory} below;
and for the LSW-character on
twisted higher K-theories \cite[\S 2.1]{LindSatiWesterland16}, see Prop. \ref{TwistedChernCharacterInTwistedHigherKTheory} below):

\begin{defn}[Higher twisted abelian de Rham cohomology]
  \label{HigherTwistedAbelianDeRhamCohomology}
  For $X \in \SmoothManifolds$,
  $r \in \mathbb{N}$
  %, $r \geq 1$,
  and
  $H_{2r+1} \in \Omega^{2r+1}_{\mathrm{dR}}(X)_{\mathrm{closed}}$
  a closed differential $(2r+1)$-form,
  the \emph{$H_{2r+1}$-twisted de Rham cohomology} of $X$
  is the
  cochain cohomology
 \vspace{-1mm}
  \begin{equation}
    \label{HigherTwistedDeRhamCohomology}
    \Omega^{\bullet+H_{2r+1}}_{\mathrm{dR}}(X)
    \;:=\;
    \frac{
      \mathrm{ker}^\bullet\big( d - H_{2r+1} \wedge (-)\big)
    }{
      \mathrm{im}^\bullet\big( d - H_{2r+1} \wedge (-)\big)
    }
  \end{equation}
  of the following $2r$-periodic cochain complex:
   \vspace{-2mm}
  $$
    \xymatrix@C=1.8em{
      \cdots
      \ar[r]
      &
      \underset{k}{\bigoplus}
      \,
      \Omega^{(n-1)+2 r k}_{\mathrm{dR}}(X)
      \ar[rrr]^-{
        (d - H_{2r+1} \wedge (-))
      }
      &&&
      \underset{k}{\bigoplus}
      \,
      \Omega^{n + 2 r k}_{\mathrm{dR}}(X)
      \ar[rrr]^-{
        (d - H_{2 r + 1} \wedge (-))
      }
      &&&
      \underset{k}{\bigoplus}
      \,
      \Omega^{(n+1) + 2 r k}_{\mathrm{dR}}(X)
      \ar[r]
      &
      \cdots
      \,.
    }
  $$
\end{defn}

In direct generalization of Prop. \ref{TwistedNonabelianDeRhamCohomologySubsumesH3TwisteddeRhamCohomology},
we find:
\begin{prop}[Twisted non-abelian de Rham cohomology subsumes higher twisted abelian de Rham cohomology]
  \label{TwistedNonabelianDeRhamCohomologySubsumesHigherTwisteddeRhamCohomology}
  For $r \in \mathbb{N}$, $r \geq 1$,
  consider a twisting $(2r+1)$-form as in \eqref{TwistDeRhamForlK2rKU1}

  \vspace{-.3cm}
  $$
    \xymatrix@R=6pt{
      \tau_{\mathrm{dR}}
      \ar@{<->}[r]
      \ar@{}[d]|-{
        \rotatebox[origin=c]{-90}{$\in$}
      }
      &
      H_{2r+1}
      \ar@{}[d]|-{
        \rotatebox[origin=c]{-90}{$\in$}
      }
      \\
      \Omega
      \big(
        X;
        \,
        \mathfrak{b}^{2r} \mathbb{R}
      \big)_{\mathrm{flat}}
      \ar@{}[r]|-{\simeq}
      &
      \Omega^{2r+1}
      \big(
        X
      \big)_{\mathrm{closed}}
    }
  $$

  \vspace{-2mm}
  \noindent
  The $\tau_{\mathrm{dR}}$-twisted non-abelian de Rham cohomology
  (Def. \ref{TwistedNonabelianDeRhamCohomology}) of
  flat twisted $\mathfrak{l}K^{2r-2}(\mathrm{ku})_1$-valued differential forms
  (Example \ref{RecoveringHigherTwistedDifferentialForms})
  is naturally equivalent to
  $H_{2r+1}$-twisted abelian de Rham cohomology
  (Def. \ref{HigherTwistedAbelianDeRhamCohomology})
  in degree\footnote{The discussion for other degrees is
  directly analogous, and we omit it for brevity.} $1 \; \mathrm{mod} \; 2r$.
  \vspace{-3mm}
  $$
    \overset{
      \mathclap{
      \raisebox{3pt}{
        \tiny
        \color{darkblue}
        \bf
        \def\arraystretch{.9}
        \begin{tabular}{c}
          twisted
          \\
          non-abelian de Rham cohomology
          \\
        \end{tabular}
      }
      }
    }{
    H^{\tau_{\mathrm{dR}}}_{\mathrm{dR}}
    \Big(
      X;
      \,
      \underset{k \in \mathbb{N}}{\oplus} \mathfrak{b}^{2rk}\mathbb{R}
    \Big)
    }
    \;\;\;\;\simeq\;\;\;\;
    \overset{
      \mathclap{
      \raisebox{3pt}{
        \tiny
        \color{darkblue}
        \bf
        \def\arraystretch{.9}
        \begin{tabular}{c}
          higher $H_{2r+1}$-twisted
          \\
          de Rham cohomology
        \end{tabular}
      }
      }
    }{
      H^{1 + H_{2r+1}}_{\mathrm{dR}}(X)
    }
    \,.
  $$
\end{prop}
\begin{proof}
  By Example \ref{RecoveringHigherTwistedDifferentialForms},
  the cocycle sets on both sides are in natural bijection.
  Hence it remains to see that the coboundary relations
  correspond to each other, under this identification.
  This proceeds verbatim, up to degree shifts, as in the proof
  of Prop. \ref{TwistedNonabelianDeRhamCohomologySubsumesH3TwisteddeRhamCohomology}
  (which is the special case of $r = 1$).
\end{proof}

\vspace{-1mm}
\noindent
\begin{example}[Degree-1 twisted non-abelian de Rham cohomology]
  \label{Degree1TwistedDeRhamCohomology}
  Def. \ref{HigherTwistedAbelianDeRhamCohomology} subsumes also
  the case of a twist in (``lower'') degree 1, for $k = 0$.
  By classical theory of sheaf cohomology for local systems
  (see e.g. \cite[Prop. 2.3]{ChenYang18} following
  \cite[\S II 5.1.1 ]{Voisin02})
  the {\it degree-1 twisted de Rham cohomology}
  in the sense of Def. \ref{HigherTwistedAbelianDeRhamCohomology}
  is equivalently
  classical sheaf cohomology with coefficients in the flat local
  sections of a trivial line bundle with flat connection.
  Beware that for more general local systems of
  lines (or even of vector spaces)
  some authors still speak
  of ``twisted de Rham cohomology'' (e.g. \cite[\S 2.1]{ChenYang18}),
  though the twist itself is then
  no longer in real/de Rham cohomology,
  whence this more general case is, in our terminology,
  no longer an example of Def. \ref{HigherTwistedAbelianDeRhamCohomology},
  but is an example of (torsion-)twisted differential cohomology \cite{GS-Deligne}.
\end{example}

\begin{example}[Cohomology operation in (higher-) twisted de Rham cohomology]
  \label{OperationInTwistedOrdinaryCohomology}
Degree-3 twisted de Rham cohomology (Def. \ref{Degree3TwistedAbelianDeRhamCohomology})
supports the following twisted cohomology operations
(Def. \ref{TwistedNonabelianCohomologyOperation}):

\noindent
{\bf (i)} {\it wedge product with $H_3$}:
\vspace{-1mm}
$$
  \xymatrix@R=-4pt{
    H^{\bullet + H_3}_{\mathrm{dR}}(X)
    \ar[r]
    &
    H^{\bullet + 3 + H_3}_{\mathrm{dR}}(X)
    \\
    \underset{k}{\sum}F_k
    \ar@{}[r]|-\longmapsto
    &
    \underset{k}{\sum}F_k \wedge H_3
  }
$$

\vspace{1mm}
\noindent
{\bf (ii)} {\it wedge square}:
\vspace{-1mm}
$$
  \xymatrix@R=-2pt{
    \underset{r}{\bigoplus}
    \,
    H^{2r + H_3}_{\mathrm{dR}}(X)
    \ar[r]
    &
    \underset{r}{\bigoplus}
    \,
    H^{2r + 2H_3}_{\mathrm{dR}}(X)
    \\
    \underset{k}{\sum}F_k
    \ar@{}[r]|-\longmapsto
    &
    \Big(
      \underset{k}{\sum}
      F_k
    \Big)
    \wedge
    \Big(
      \underset{k}{\sum}
      F_k
    \Big)
  }
$$

\vspace{-1mm}
\noindent {\bf (iii)} {\it compositions of these}:
\vspace{-1mm}
$$
  \xymatrix@R=-2pt{
    \underset{r}{\bigoplus}
    \,
    H^{2r + H_3}_{\mathrm{dR}}(X)
    \ar[r]
    &
    \underset{r}{\bigoplus}
    \,
    H^{2r + 1 + 2 H_3}_{\mathrm{dR}}(X)
    \\
    \underset{k}{\sum}F_k
    \ar@{}[r]|-\longmapsto
    &
    \Big(
      \underset{k}{\sum}
      F_k
    \Big)
    \wedge
    \Big(
      \underset{k}{\sum}
      F_k
    \Big)
    \wedge H_3
  }
$$
It is noteworthy that
terms of the form {\bf (iii)} arise
in type IIA string theory, together with
terms of the form $I_8 \cup [H_3]$ \eqref{TheI8Polynomial},
see \cite{GS-RR}.

This evidently generalizes to higher twisted
de Rham cohomology (Def. \ref{HigherTwistedAbelianDeRhamCohomology})
and higher twisted real cohomology
in the sense of \cite{GS-HigherDeligne}, with $H_3$ replaced by $H_{2r+1}$ for $r\in \mathbb{N}$.
\end{example}

\medskip
\noindent
{\bf Homotopical formulation of twisted non-abelian de Rham cohomology.}
In preparation of the twisted non-abelian de Rham theorem (Thm. \ref{TwistedNonAbelianDeRhamTheorem})
we give a homotopy-theoretic reformulation of twisted non-abelian de Rham cohomology
(Def. \ref{TwistedNonabelianDeRhamCohomology}):

\medskip

\begin{lemma}[Pullback to de Rham complex over cylinder of manifold is relative
path space object]
  \label{PullbackTodeRhamComplexOverCylinderOfManifoldIsRelativePathSpaceObject}
  $\,$

\noindent  Let $X \,\in\,\SmoothManifolds$,
  let $\mathfrak{b} \in \LInfinityAlgebras$ (Example \ref{CEAlgebraOfLieAlgebra}) with Chevalley-Eilenberg
  algebra
    $\mathrm{CE}(\mathfrak{b}) \in$  $\dgcAlgebras{\mathbb{R}}$ \eqref{CEAlgebraOfLInfinityAlgebra}, and let
  $
    \raisebox{0pt}{
    \xymatrix{
\big\{      \Omega_{\mathrm{dR}}^\bullet(X)
      \ar@{<-}[r]^-{ \tau_{\mathrm{dR}}^\ast }
      &
      \mathrm{CE}(\mathfrak{b})
      \big\}
    }
    }
    \;
    \in
    \;
    \dgcAlgebrasProj{\mathbb{R}}^{ \mathrm{CE}(\mathfrak{b})/ }
  $
  be a morphism of dgc-algebras to the
  de Rham complex of $X$ (Example \ref{SmoothdeRhamComplex}),
  regarded as an object in
  the coslice model category (Example \ref{SliceModelCategory})
  of $\dgcAlgebrasProj{\mathbb{R}}$ (Prop. \ref{ProjectiveModelStructureOndgcAlgebras})
  under $\mathrm{CE}(\mathfrak{b})$.
  Then a path space object (Def. \ref{PathSpaceObject})
  for $\tau^\ast_{\mathrm{dR}}$ is given by this diagram:
    \vspace{-2mm}
  $$
    \xymatrix@R=8pt@C=4em{
      \Omega_{\mathrm{dR}}^\bullet(X)
      \ar[rr]^-{
        \mathrm{pr}_X^\ast
     \in \, \mathrm{W} }
      \ar@{<-}[ddrr]_-{ \tau^\ast_{\mathrm{dR}} }
      &&
      \Omega_{\mathrm{dR}}^\bullet(X \times \mathbb{R})
      \ar[rr]^-{
        (i_0^\ast\,,\, i_1^\ast)
     \in \, \mathrm{Fib} }
      \ar@{<-}[dd]|-{
        \mathrm{pr}_X^\ast
        \,\circ\,
        \tau^\ast_{\mathrm{dR}}
      }
      &&
      \Omega_{\mathrm{dR}}^\bullet(X)
      \,\oplus\,
      \Omega_{\mathrm{dR}}^\bullet(X)
      \ar@{<-}[ddll]^-{\;\;
        (
          \tau^\ast_{\mathrm{dR}}
          \,,\,
          \tau^\ast_{\mathrm{dR}}
        )
      }
      \\
      \\
      &&
      \mathrm{CE}(\mathfrak{b})
      \,,
    }
  $$

    \vspace{-2mm}
\noindent
  where the top morphisms are from \eqref{PathSpaceObjectFordeRhamComplex}.
\end{lemma}
\begin{proof}
  It is clear that the diagram commutes, by construction.
  Moreover, the top morphisms are a weak equivalence
  followed by a fibration in
  $\dgcAlgebrasProj{\mathbb{R}}$, by Lemma \ref{DeRhamComplexOverCylinderManifoldIsPathSpaceObject}.
  Therefore, by the nature of the coslice model structure
  (Example \ref{SliceModelCategory}) the total diagram
  constitutes a factorization of the diagonal on
  $\tau^\ast_{\mathrm{dR}}$ through a weak equivalence
  followed by a fibration, as required \eqref{PathSpaceFactorization}.
  (To see that the composite really is still the diagonal
  morphism in the coslice, observe that Cartesian products
  in any coslice category are reflected in the underlying category.)
  It only remains to observe that $\tau^\ast_{\mathrm{dR}}$
  is actually a fibrant object in the coslice model category.
  But the terminal object in the coslice is clearly the
  unique morphism from $\mathrm{CE}(\mathfrak{b})$ to the zero-algebra
  (Example \ref{TheTerminalAlgebra}), so that in fact
  every object in the coslice is still fibrant
    \vspace{-3mm}
  \begin{equation}
    \label{EveryObjectInCosliceOfProjectiveModelStrucOndgcAlgebrasIsCofibrant}
    \raisebox{20pt}{
    \xymatrix@R=.1em@C=4em{
      \Omega^\bullet_{\mathrm{dR}}(X)
      \ar[rr]^-{ \in \, \mathrm{Fib} }
      \ar@{<-}[dr]_-{ \tau^\ast_{\mathrm{dR}} }
      &&
      0
      \ar@{<-}[dl]
      \\
      &
      \mathrm{CE}(\mathfrak{b})
    }
    }
  \end{equation}

  \vspace{-2mm}
  \noindent
  as in Remark \ref{AlldgcAlgebrasAreProjectivelyFibrant}.
\end{proof}

\begin{prop}[Twisted non-abelian de Rham cohomology via the coslice dgc-homotopy category]
  \label{TwistedNonAbelianDeRhamCohomologyViaThedgcHomotopyCategory}
  Consider $X \,\in\,\SmoothManifolds$,
  let
  \vspace{-2mm}
  $$
    \raisebox{10pt}{
    \xymatrix@R=8pt{
      \mathfrak{g}
      \ar[r]
      &
      \widehat{\mathfrak{b}}
      \ar[d]^-{ \mathfrak{p} }
      \\
      & \mathfrak{b}
    }
    }
    \;\;\;\;
    \in
    \;
    \LInfinityAlgebrasNil
  $$

  \vspace{-1mm}
  \noindent
  be an $L_\infty$-algebraic local coefficient bundle
  (Def. \ref{LocalLInfinityAlgebraicCoefficients}) of
  nilpotent $L_\infty$-algebras
  (Def. \ref{NilpotentLInfinityAlgebras})
  with Chevalley-Eilenberg
  algebra
  $\mathrm{CE}(\widehat{\mathfrak{b}}),\,\mathrm{CE}(\mathfrak{b}) \in \dgcAlgebras{\mathbb{R}}$ \eqref{CEAlgebraOfLInfinityAlgebra}, and let
    \vspace{-2mm}
  \begin{equation}
    \label{dgcTwistingMorphism}
    \raisebox{0pt}{
    \xymatrix{
      \Omega_{\mathrm{dR}}^\bullet(X)
      \ar@{<-}[r]^-{ \tau_{\mathrm{dR}}^\ast }
      &
      \mathrm{CE}(\mathfrak{b})
    }
    }
    \;\;\;\;\;\;
    \in
    \;
    \dgcAlgebrasProj{\mathbb{R}}^{ \mathrm{CE}(\mathfrak{b})/ }
  \end{equation}

  \vspace{-2mm}
  \noindent
  be a morphism of dgc-algebras to the
  de Rham complex of $X$ (Example \ref{SmoothdeRhamComplex}),
  hence a flat $\mathfrak{b}$-valued differential form
  (Def. \ref{FlatLInfinityAlgebraValuedDifferentialForms})
    \vspace{-2mm}
  $$
    \tau_{\mathrm{dR}}
    \;\in\;
    \Omega_{\mathrm{dR}}
    (
      X;
      \,
      \mathfrak{b}
    )
    \,,
  $$

  \vspace{-1mm}
  \noindent
  equivalently regarded as an object in
  the coslice model category (Example \ref{SliceModelCategory})
  of $\dgcAlgebrasProj{\mathbb{R}}$ (Prop. \ref{ProjectiveModelStructureOndgcAlgebras})
  under $\mathrm{CE}(\mathfrak{b})$.
  Then the $\tau_{\mathrm{dR}}$-twisted non-abelian de Rham cohomology
  of $X$ with
  coefficients in $\mathfrak{g}$ (Def. \ref{TwistedNonabelianDeRhamCohomology})
  is in natural bijection with the hom-set in the
  homotopy category (Def. \ref{HomotopyCategory})
  of the coslice model category
  $\dgcAlgebrasProj{\mathbb{R}}^{\mathrm{CE}(\mathfrak{b})}$
  (Example \ref{SliceModelCategory})
  of the projective model structure on dgc-algebras
  (Prop. \ref{ProjectiveModelStructureOndgcAlgebras})
  from $\mathrm{CE}(\mathfrak{p})$ \eqref{CEAlgebraOfLInfinityAlgebraicCoefficientBundle}
  to $\tau^\ast_{\mathrm{dR}}$ \eqref{dgcTwistingMorphism}:
  \vspace{-1mm}
  \begin{equation}
    \label{BijectionBetweenTwistedNonabelianDeRhamCohomologyAndHomSetsIndgcHomotopyCategory}
    H^{\tau_{\mathrm{dR}}}_{\mathrm{dR}}
    \big(
      X;
      \,
      \mathfrak{g}
    \big)
    \;\simeq\;
    \mathrm{Ho}
    \Big(
      \dgcAlgebrasProj{\mathbb{R}}^{\mathrm{CE}(\mathfrak{b})/}
    \Big)
    \big(
      \mathrm{CE}(\mathfrak{p})
      \,,\,
      \tau^\ast_{\mathrm{dR}}
    \big).
  \end{equation}
\end{prop}
\begin{proof}
  Consider a pair of dgc-algebra homomorphisms in the coslice
  \begin{equation}
    \label{PairOfdgcAlgebraHomomorphismsInCoslice}
    \raisebox{20pt}{
    \xymatrix@R=1em{
      \Omega^\bullet_{\mathrm{dR}}
      \big(
        X
      \big)
      \ar@{<-}@/^.7pc/[rr]|-{ \;A^{(0)}\; }
      \ar@{<-}@/_.7pc/[rr]|-{ \;A^{(1)}\; }
      \ar@{<-}[dr]_-{
        \tau^\ast_{\mathrm{dR}}
      }
      &&
      \mathrm{CE}(\widetilde{\mathfrak{b}})
      \ar@{<-}[dl]^-{
        \mathrm{CE}(\mathfrak{p})
      }
      \\
      &
      \mathrm{CE}(\mathfrak{b})
    }
    }
    \;\;\;\;\;\;\;
    \in
    \dgcAlgebrasProj{\mathbb{R}}^{\mathrm{CE}(\mathfrak{b})/}
    \big(
      \mathrm{CE}(\mathfrak{p})
      \,,\,
      \tau^\ast_{\mathrm{dR}}
    \big)
    \,,
  \end{equation}

  \vspace{-2mm}
  \noindent
  hence of flat $\tau_{\mathrm{dR}}$-twisted
  $\mathfrak{g}$-valued differential forms,
  according to Def. \ref{FlatTwistedLInfinityAlgebraValuedDifferentialForms}.
  Observe that:
  \begin{itemize}
    \vspace{-.1cm}
    \item[\bf (i)] $\mathrm{CE}(\mathfrak{p})$
    is cofibrant in $\dgcAlgebrasProj{\mathbb{R}}^{\mathrm{CE}(\mathfrak{p})/}$,
    since:

    \vspace{-2mm}
    \noindent
    {\bf (a)} the initial object in the coslice is
    $\xymatrix@C=12pt{
      \mathrm{CE}(\mathfrak{b})
        \ar@{<-}[r]^-{\mathrm{id}}
        &
        \mathrm{CE}(\mathfrak{b})
      }$,

      \vspace{-1mm}
      \noindent
      {\bf (b)} the unique morphism from this object to
      $\mathrm{CE}(\mathfrak{p})$ is

        \vspace{-4mm}
      \begin{equation}
        \label{ACosliceCofibration}
        \xymatrix@R=6pt@C=3em{
          \mathrm{CE}(\mathfrak{b})
          \ar[rr]^-{
            \mathrm{CE}( \mathfrak{p} )
        \in\, \mathrm{Cof} }
          \ar@{<-}[dr]_-{ \mathrm{id} }
          &&
          \mathrm{CE}\big(\,\widehat{\mathfrak{b}}\,\big)
          \ar@{<-}[dl]^-{ \mathrm{CE}(\mathfrak{p}) }
          \\
          &
          \mathrm{CE}(\mathfrak{b})
        }
      \end{equation}

      \vspace{-3mm}
      {\bf (c)}
      $\mathrm{CE}(\mathfrak{p})$ is a cofibration
      in $\dgcAlgebrasProj{\mathbb{R}}$, by  \eqref{CEAlgebraOfLInfinityAlgebraicCoefficientBundle},
      so that the diagram
      \eqref{ACosliceCofibration}
      is a cofibration in the coslice model category,
      by Example \ref{SliceModelCategory}.

    \vspace{-2mm}
    \item[\bf (ii)]
    $ \mathrm{pr}_X^\ast \,\circ\, \tau^\ast_{\mathrm{dR}}$ is fibrant in
    $\dgcAlgebrasProj{\mathbb{R}}^{\mathrm{CE}(\mathfrak{b})/}$, by
    \eqref{EveryObjectInCosliceOfProjectiveModelStrucOndgcAlgebrasIsCofibrant};

    \vspace{-.1cm}
    \item[\bf (iii)]
    A right homotopy (Def. \ref{RightHomotopy})
    between the pair \eqref{PairOfdgcAlgebraHomomorphismsInCoslice}
    of coslice morphisms, with respect to the path space object
    from Lemma \ref{PullbackTodeRhamComplexOverCylinderOfManifoldIsRelativePathSpaceObject},
    namely a $\widetilde A$ that makes the following diagram commute
    \vspace{-2mm}
    \begin{equation}
    \label{TwistedNonAbelianCoboundaryAsRightHomotopyOfdgcAlgebras}
    \raisebox{40pt}{
    \xymatrix{
      \Omega^\bullet_{\mathrm{dR}}(X)
      \ar@{<-}[d]_-{
        i_0^\ast
      }
      \ar@{<-}[drr]|-{ \;A^{(0)}\; }
      \\
      \Omega^\bullet_{\mathrm{dR}}(X \times \mathbb{R})
      \ar@{<-}[dr]^{ \scalebox{0.6}{$
        \mathrm{pr}_X^\ast
        ,
        \tau^\ast_{\mathrm{dR}}
     $}
      }|>>>>{ \phantom{AA} }
      \ar[d]_-{
        i_1^\ast
      }
      \ar@{<-}[rr]|-{ \;\widetilde A\; }
      &&
      \mathrm{CE}(\, \widehat{\mathfrak{b}}\,)
      \ar@{<-}[dl]^-{ \mathrm{CE}(\mathfrak{p}) }
      \\
      \Omega^\bullet_{\mathrm{dR}}(X)
      \ar@{<-}[urr]|>>>>>>>>>>>>>>{ \;A^{(1)}\; }
      &
      \mathrm{CE}(\,\widehat{\mathfrak{b}}\,)
    }
    }
  \end{equation}

  \vspace{-3mm}
  \noindent
  is manifestly the same as a coboundary $\widetilde A$
  between the corresponding flat twisted $\mathfrak{g}$-valued
  forms according to Def. \ref{CoboundariesBetweenFlatTwistedLInfinityAlgebraValuedForms}:

  \vspace{-1mm}
  {\bf (a)} The top part of \eqref{TwistedNonAbelianCoboundaryAsRightHomotopyOfdgcAlgebras}
  is, just as in \eqref{NonAbelianCoboundaryAsRightHomotopyOfdgcAlgebras},
  the flat twisted $\widehat{\mathfrak{g}}$-valued form
  on the cylinder over $X$
  that is required by \eqref{CoboundaryBetweenFlatTwistedForms};

  \vspace{-1mm}
  {\bf (b)} the bottom part of \eqref{TwistedNonAbelianCoboundaryAsRightHomotopyOfdgcAlgebras}
  is the condition \eqref{ConditionOnCoboundaryBetweenFlatTwistedForms}
  on the extension of the twist to the cylinder over $X$.

  \end{itemize}

  \vspace{-2mm}
  \noindent
  Therefore, Prop. \ref{RightHomotopyReflectsHomotopyClasses}
  says
  that the quotient set \eqref{TwistedNonabelianDeRhamCohomologyAsQuotientSet}
  defining the twisted non-abelian de Rham cohomology
  is in natural bijection to the hom-set in the coslice homotopy category.
\end{proof}

\medskip
\noindent {\bf The twisted non-abelian de Rham theorem.} With this in hand we
may finally prove the main result in this section, generalizing the
non-abelian de Rham theorem (Thm. \ref{NonAbelianDeRhamTheorem}) to
the twisted case:

\begin{theorem}[Twisted non-abelian de Rham theorem]
  \label{TwistedNonAbelianDeRhamTheorem}
  Let $X \,\in\, \HomotopyTypes$
  equipped with the structure of
  a smooth manifold, and let
  \vspace{1mm}
  \begin{equation}
    \label{LocalCoefficientBundleInTwistedNonabeliandeRhamTheorem}
   \hspace{-2cm}
   \begin{tikzcd}
      A
      \ar[rr]
      &
      \ar[
        d,
        phantom,
        "
          \mbox{
            \tiny
            \color{darkblue}
            \bf
            local coefficient bundle
          }
        "
      ]
      &
      A \!\sslash\! G
      \ar[
        d,
        "\rho"
      ]
      \\
      &{}& B G
    \end{tikzcd}
    {\phantom{AAAA}}
    \in
    \;
    \NilpotentConnectedQFiniteHomotopyTypes
  \end{equation}

  \vspace{-3mm}
  \noindent
  be a local coefficient bundle \eqref{LocalCoefficientBundle}
  of connected $\mathbb{Q}$-finite nilpotent homotopy types
  (Def. \ref{NilpotentConnectedSpacesOfFiniteRationalType})
  such that the action of $\pi_1(B G) = \pi_0(G)$ on the
  real homology groups of $A$ is nilpotent.
  Consider, via Prop. \ref{RationalizationOfLocalCoefficients},
  the rationalized
  coefficient bundle $L_{\mathbb{R}}(\rho)$
  with corresponding
  $L_\infty$-algebraic coefficient
  bundle $\mathfrak{l}\rho$ (Def. \ref{LocalLInfinityAlgebraicCoefficients})
  of the relative real Whitehead $L_\infty$-algebra
  (Prop. \ref{WhiteheadLInfinityAlgebrasRelative}):
  \vspace{-2mm}
  $$
      \xymatrix@R=1.5em{
      L_{\mathbb{R}}A
      \ar[rr]
      &
      \ar@{}[d]|-{
        \mbox{
          \tiny
          \color{darkblue}
          \bf
          \def\arraystretch{.9}
          \begin{tabular}{c}
            rationalized
            \\
            local coefficient bundle
          \end{tabular}
        }
      }
      &
      \big(L_{\mathbb{R}}A\big)
      \!\sslash\!
      \big(L_{\mathbb{R}}G\big)
      \ar[d]^-{ L_{\mathbb{R}}(\rho) }
      \\
      &&
      L_{\mathbb{R}}B G
    }
    \phantom{AAAAAAA}
    \xymatrix@R=1.5em{
      \mathfrak{a}
      \ar[rr]
      &
      \ar@{}[d]|-{
        \mbox{
          \tiny
          \color{darkblue}
          \bf
          \def\arraystretch{.9}
          \begin{tabular}{c}
            $L_\infty$-algebraic coefficient bundle
            \\
            of Whitehead $L_\infty$-algebras
          \end{tabular}
        }
        \;\;\;\;\;\;\;\;\;
      }
      &
      \widehat{\mathfrak{b}} \;.
      \ar[d]^{ \mathfrak{p} }
      \\
      && \mathfrak{b}
    }
  $$

  \vspace{-1mm}
  \noindent
  Then, for
  $$
    \xymatrix{
      X \ar[r]^-\tau & L_{\mathbb{R}} B G
    }
    \;\;\;
    \in
    \HomotopyTypes
  $$

  \vspace{-1mm}
  \noindent
  a twist,
  the $\tau$-twisted non-abelian
  real cohomology (Def. \ref{TwistedNonAbelianRealCohomology})
  of $X$ with local coefficients in $L_{\mathbb{R}}(\rho)$
  (Prop. \ref{RationalizationOfLocalCoefficients})
  is in natural bijection with
  the $\tau_{\mathrm{dR}}$-twisted non-abelian de Rham cohomology
  (Def. \ref{TwistedNonabelianDeRhamCohomology})
  of $X$ with local coefficients in $\mathfrak{p}$,

  \vspace{-3mm}
  \begin{equation}
    \label{EquivalenceBetweenTwistedNonabelianRealAndNonaebelianDeRhamCohomology}
    \overset{
      \mathclap{
      \raisebox{3pt}{
        \tiny
        \color{darkblue}
        \bf
        \def\arraystretch{.9}
        \begin{tabular}{c}
          $\tau$-twisted
          non-abelian
          \\
          real cohomology
        \end{tabular}
      }
      }
    }{
    H^{\tau}
    \big(
      X;
      \,
      L_{\mathbb{R}}A
    \big)
    }
    \;\;\;\;\simeq\;\;\;\;
    \overset{
      \mathclap{
      \raisebox{3pt}{
        \tiny
        \color{darkblue}
        \bf
        \def\arraystretch{.9}
        \begin{tabular}{c}
          $\tau_{\mathrm{dR}}$-twisted
          non-abelian
          \\
          de Rham cohomology
        \end{tabular}
      }
      }
    }{
    H^{\tau_{\mathrm{dR}}}_{\mathrm{dR}}
    (
      X;
      \,
      \mathfrak{a}
    )
    }
    \,,
  \end{equation}

  \vspace{-1mm}
  \noindent
  where the twists are related by the plain
  non-abelian de Rham theorem (Theorem \ref{NonAbelianDeRhamTheorem}):
  \vspace{-2mm}
  \begin{equation}
    \label{RelationOfTwistsInTwistedNonAbelianDeRhamTheorem}
    \xymatrix@R=5pt{
      \;[\tau]\;
      \ar@{<->}[r]
      \ar@{}[d]|-{
        \rotatebox[origin=c]{-90}{$\in$}
      }
      &
      \;[\tau_{\mathrm{dR}}]\;
      \ar@{}[d]|-{
        \rotatebox[origin=c]{-90}{$\in$}
      }
      \\
      H
      \big(
        X;
        \,
        L_{\mathbb{R}}B G
      \big)
      \ar@{}[r]|-{\simeq}
      &
      H_{\mathrm{dR}}
      \big(
        X;
        \,
        \mathfrak{b}
      \big)
    }
  \end{equation}
\end{theorem}
\begin{proof}
  This is established by the following sequence of natural bijections
  of hom-sets
  (where on the right we are illustrating the structure of their elements):
  \vspace{-4mm}
  $$
  \def\arraystretch{3.8}
  \begin{array}{llll}
    H^\tau
    \big(
      X;
      \,
      L_{\mathbb{R}} A
    \big)
    &
    \;=\;
    \mathrm{Ho}
    \big(
        \SimplicialSets_{\mathrm{Qu}}^{/L_{\mathbb{R}}B G}
    \big)
    \big(
      \tau
      \,,\,
      L_{\mathbb{R}}(\rho)
    \big)
    &&
    =
    \left\{
    \begin{tikzcd}[column sep={between origins, 33}, row sep=10]
      X
      \ar[dr, "\tau"{left, pos=.6, xshift=-3pt}]
      \ar[
        rr,
        dashed
      ]
      &&
      L_{\mathbb{R}}
      \big(
        A \!\sslash\! G
      \big)
      \ar[
        dl,
        "L_{\mathbb{R}}(\rho)"{pos=.2}
      ]
      \\
      &
      L_\mathbb{R} B G
    \end{tikzcd}
    \right\}
    \\
    &
    \;=\;
    \mathrm{Ho}
    \big(
      \SimplicialSets_{\mathrm{Qu}}^{/\Bexp_{\mathrm{PL}} \mathrm{CE}(\mathfrak{b})}
    \big)
    \big(
      \tau
      \,,\,
      \Bexp_{\mathrm{PL}} \mathrm{CE}(\mathfrak{p})
    \big)
    &&
    =
    \left\{
    \begin{tikzcd}[column sep={between origins, 33}, row sep=10]
      X
      \ar[dr, "\tau"{left, pos=.6, xshift=-3pt}]
      \ar[
        rr,
        dashed
      ]
      &&
      \Bexp_{\mathrm{PL}}
      \mathrm{CE}(\mathfrak{a})
      \ar[
        dl,
        "{
          \scalebox{.7}{$
            \Bexp_{\mathrm{PL}}(\mathrm{CE}(\mathfrak{p}))
          $}
        }"{pos=.1}
      ]
      \\
      &
      \Bexp_{\mathrm{PL}}\mathrm{CE}(\mathfrak{b})
    \end{tikzcd}
    \right\}
    \\
    &
    \;\simeq\;
    \mathrm{Ho}
    \Big(
      \big(
        \dgcAlgebrasOpProj{\mathbb{R}}
      \big)^{ /\mathrm{CE}(\mathfrak{g}) }
    \Big)
    \big(
      \widetilde{\tau},
      \,
      \mathrm{CE}(\mathfrak{p})
    \big)
    &&
    =
    \left\{
    \begin{tikzcd}[column sep={between origins, 33}, row sep=10]
      \Omega^\bullet_{\mathrm{PLdR}}(X)
      &&
      \mathrm{CE}(\widehat{\mathfrak{b}})
      \ar[
        ll,
        dashed
      ]
      \\
      &
      \mathrm{CE}(\mathfrak{b})
      \ar[
        ul,
        "{
          \widetilde{\tau}
        }"{left, pos=.1, xshift=-4pt}
      ]
      \ar[
        ur,
        "{
          \scalebox{.7}{$
            \mathrm{CE}(\mathfrak{p})
          $}
        }"{right, pos=.2, xshift=3pt}
      ]
    \end{tikzcd}
    \right\}
    \\
    &
    \;\simeq\;
    \mathrm{Ho}
    \Big(
      \big(
        \dgcAlgebrasOpProj{\mathbb{R}}
      \big)^{ /\mathrm{CE}(\mathfrak{g}) }
    \Big)
    \big(
      \tau_{\mathrm{dR}},
      \,
      \mathrm{CE}(\mathfrak{p})
    \big)
    &&
    =
    \left\{
    \begin{tikzcd}[column sep={between origins, 33}, row sep=10]
      \Omega^\bullet_{\mathrm{dR}}(X)
      &&
      \mathrm{CE}(\widehat{\mathfrak{b}})
      \ar[
        ll,
        dashed
      ]
      \\
      &
      \mathrm{CE}(\mathfrak{b})
      \ar[
        ul,
        "{
          \tau^\ast_{\mathrm{dR}}
        }"{left, pos=.1, xshift=-3pt}
      ]
      \ar[
        ur,
        "{
          \scalebox{.7}{$
            \mathrm{CE}(\mathfrak{p})
          $}
        }"{right, pos=.1, xshift=3pt}
      ]
    \end{tikzcd}
    \right\}
    \\[-30pt]
    &
    \;\simeq\;
    H^{\tau_{\mathrm{dR}}}
    \big(
      X;
      \,
      \mathfrak{l}A
    \big)
    \,.
  \end{array}
$$

  \vspace{-2mm}
  \noindent
  Here the first line is the definition
  of twisted non-abelian real cohomology (Def. \ref{TwistedNonAbelianRealCohomology}),
  while the second line inserts the definition of
  $L_{\mathbb{R}}$ (Def. \ref{Lk}),
  with $\mathrm{CE}(\mathfrak{l}\rho)$
  serving as the required \eqref{RightDerivedFunctorByFibrantReplacement}
  fibrant resolution \eqref{CoFibrantObjectsInSliceModelStructure}
  of $\Omega^\bullet_{\mathrm{PLdR}}(\rho)$.

  The key step is the third line, which uses the hom-isomorphism
  \eqref{AnAdjunction} of the derived adjunction \eqref{DerivedAdjunction}
  of the sliced Quillen adjunction (Ex. \ref{SlicedQuillenAdjunction})
  of the PLdR-adjunction (Prop. \ref{QuillenAdjunctionBetweendgcAlgebrasAndSimplicialSets}),
  using the form \eqref{LeftSliceAdjointProducesAdjuncts} of its left adjoint
  with the observation that this is already derived \eqref{LeftDerivedFunctorByCofibrantReplacement}
  since $\tau$ is necessarily cofibrant,
  by
  \eqref{EverySimplicialSetIsCofibrantInClassicalModelStructure}
  and
  \eqref{CoFibrantObjectsInSliceModelStructure}.

  The fourth step is composition with the slice morphism
  exhibiting \eqref{RelationOfTwistsInTwistedNonAbelianDeRhamTheorem}
  $$
    \begin{tikzcd}[column sep={between origins, 36}, row sep=8]
      \Omega^\bullet_{\mathrm{dR}}(X)
      &&
      \Omega^\bullet_{\mathrm{PLdR}}(X)
      \ar[
        ll, dashed, "\in \mathrm{W}"{above}
      ]
      \\
      &
      \mathrm{CE}(\mathfrak{b})
      \mathrlap{\,,}
      \ar[
        ur,
        "\widetilde \tau"{right, pos=.1, xshift=4pt}
      ]
      \ar[
        ul,
        "\tau^\ast_{\mathrm{dR}}"{left, pos=.1, xshift=-4pt}
      ]
    \end{tikzcd}
  $$
  which is an isomorphism in the homotopy category
  by
  Lemma \ref{PLdeRhamComplexOnSmoothManifoldsEquivalentToSmoothDeRhamComplex}
  (as, in the untwisted case,
  in the last step of \eqref{TowardsTheNonAbelianDeRhamTheorem}).
  The last step is Prop. \ref{TwistedNonAbelianDeRhamCohomologyViaThedgcHomotopyCategory},
\end{proof}

\newpage

%%%%%%%%%%%%%%%%%%%%%%%%%%%%%%%%%%%%%%%%%%%%%%%%%%%
\section{The (differential) non-abelian character map}
  \label{ChernCharacterInNonabelianCohomology}
%%%%%%%%%%%%%%%%%%%%%%%%%%%%%%%%%%%%%%%%%%%%%%%%%%

We introduce the character map in non-abelian cohomology
(Def. \ref{NonAbelianChernDoldCharacter})
and then discuss how it specializes to:

\cref{TheChernDoldCharacter} -- the Chern-Dold character on generalized cohomology;

\cref{TheChernWeilHomomorphism} -- the Chern-Weil homomorphism on degree-1 non-abelian cohomology;

\cref{NonabelianDifferentialCohomology} --
 the Cheeger-Simons differential characters on degree-1 non-abelian cohomology.

\medskip

\begin{defn}[Rationalization and realification in non-abelian cohomology]
  \label{RationalizationOfCoefficientsInNonabelianCohomology}
  Let $A \in \NilpotentConnectedQFiniteHomotopyTypes$
  (Def. \ref{NilpotentConnectedSpacesOfFiniteRationalType}).

  \noindent
  {\bf (i)} We write

  \vspace{-.4cm}
  \begin{equation}
    \label{RationalizationOfCoefficients}
    (\eta_A^{\color{blue}\mathbb{Q}})_\ast
    \;:\;
    \begin{tikzcd}
      \overset{
        \mathclap{
        \raisebox{3pt}{
          \tiny
          \color{darkblue}
          \bf
          \def\arraystretch{.9}
          \begin{tabular}{c}
            non-abelian
            \\
            cohomology
          \end{tabular}
        }
        }
      }{
        H
        (
          -;
          \,
          A
        )
      }
      \ar[
        rr,
        "{
          H(-;\, \eta^{\scalebox{.7}{$\color{blue}\mathbb{Q}$}}_A)
          \,=\,
          H(-;\, \mathbb{D}\eta^{\mathrm{P}{\color{blue}\mathbb{Q}}\mathrm{L}}_A)
        }",
        "{
          \mbox{
            \tiny
            \color{greenii}
            \bf
            rationalization
          }
        }"{below}
      ]
      &[+50pt]&
      \overset{
        \mathclap{
        \raisebox{3pt}{
          \tiny
          \color{darkblue}
          \bf
          \def\arraystretch{.9}
          \begin{tabular}{c}
            non-abelian
            \\
            rational cohomology
          \end{tabular}
        }
        }
      }{
        H
        \big(
          -;
          \,
          L_{\color{blue}\mathbb{Q}}A
        \big)
      }
    \end{tikzcd}
  \end{equation}

  \vspace{-1mm}
  \noindent
  for the cohomology operation (Def. \ref{NonAbelianCohomologyOperations})
  from non-abelian $A$-cohomology
  (Def. \ref{NonAbelianCohomology})
  to non-abelian rational cohomology (Def. \ref{NonAbelianRealCohomology}),
  which is induced \eqref{InducedCohomologyOperation}
  by the rationalization map $\eta^{\mathbb{Q}}_A$
  (Def. \ref{Rationalization}),
  or equivalently, via the Fundamental Theorem (Prop. \ref{FundamentalTheoremOfdgcAlgebraicRationalHomotopyTheory}),
  by the derived unit of the rational PL de Rham adjunction.

  \noindent
  {\bf {(ii)}}
  Analogously, we write
  \begin{equation}
    \label{RealificationOfCoefficients}
    (\eta_A^{\color{blue}\mathbb{R}})_\ast
    \;:\;
    \begin{tikzcd}
      \overset{
        \mathclap{
        \raisebox{3pt}{
          \tiny
          \color{darkblue}
          \bf
          \def\arraystretch{.9}
          \begin{tabular}{c}
            non-abelian
            \\
            cohomology
          \end{tabular}
        }
        }
      }{
        H
        (
          -;
          \,
          A
        )
      }
      \ar[
        rr,
        "{
          H(-;\, \mathbb{D}\eta^{\mathrm{P}{\color{blue}\mathbb{R}}\mathrm{L}}_A)
        }",
        "{
          \mbox{
            \tiny
            \color{greenii}
            \bf
            real-ification
          }
        }"{below}
      ]
      &[+50pt]&
      \overset{
        \mathclap{
        \raisebox{3pt}{
          \tiny
          \color{darkblue}
          \bf
          \def\arraystretch{.9}
          \begin{tabular}{c}
            non-abelian
            \\
            real cohomology
          \end{tabular}
        }
        }
      }{
        H
        \big(
          -;
          \,
          L_{\color{blue}\mathbb{R}}A
        \big)
      }
    \end{tikzcd}
  \end{equation}
  for the cohomology operation to non-abelian real cohomology
  that is induced by the derived PL de Rham adjunction unit over the real numbers
  into $L_{\mathbb{R}}$ \eqref{RealificationMonad}.

  \noindent
  {\bf (iii)}
  For, moreover, $X \in \NilpotentConnectedQFiniteHomotopyTypes$
  (Def. \ref{NilpotentConnectedSpacesOfFiniteRationalType}),
  we consider the
  cohomology operation shown by the dashed arrow here:
  \begin{equation}
    \label{ExtensionOfScalarsFromRationalToRealCohomology}
  \end{equation}

  \vspace{-1.2cm}
  $$
    \begin{tikzcd}
      H
      \big(
        X;\,
        L_{\color{blue}\mathbb{Q}} A
      \big)
      \ar[
        d,
        "\widetilde{(-)}"{left},
        "\rotatebox{90}{$\sim$}"{right,xshift=0pt}
      ]
      \ar[
        rr,
        dashed,
        "(-) \otimes_{{}_{\mathbb{Q}}} \mathbb{R}"
      ]
      &&
      H
      \big(
        X;\,
        L_{\color{blue}\mathbb{R}} A
      \big)
      \\
      \mathrm{Ho}
      \big(
        \dgcAlgebrasOpProj{\color{blue}\mathbb{Q}}
      \big)
      \big(
        \LeftDerived
        \Omega^\bullet_{\mathrm{P}{\color{blue}\mathbb{Q}}\mathrm{LdR}}(X),
        \,
        \LeftDerived
        \Omega^\bullet_{\mathrm{P}{\color{blue}\mathbb{Q}}\mathrm{LdR}}(A)
      \big)
      \ar[
        rr,
        "{
          \RightDerived
          \big(
            (-) \otimes_{{}_{\mathbb{Q}}} \mathbb{R}
          \big)
        }"
      ]
      &&
      \mathrm{Ho}
      \big(
        \dgcAlgebrasOpProj{\color{blue}\mathbb{R}}
      \big)
      \big(
        \LeftDerived
        \Omega^\bullet_{\mathrm{P}{\color{blue}\mathbb{R}}\mathrm{LdR}}(X),
        \,
        \LeftDerived
        \Omega^\bullet_{\mathrm{P}{\color{blue}\mathbb{R}}\mathrm{LdR}}(A)
      \big)
      \,,
      \ar[
        u,
        "\widetilde{(-)}"{right},
        "\rotatebox{90}{$\sim$}"{left,xshift=0pt}
      ]
    \end{tikzcd}
  $$
  hence the composition of:

  \noindent
  (i) the hom-isomorphisms $\widetilde{(-)}$ \eqref{AnAdjunction} of the
  derived
  \eqref{DerivedAdjunction}
  PL de Rham Quillen adjunction
  (Prop. \ref{QuillenAdjunctionBetweendgcAlgebrasAndSimplicialSets})
  over the rational and over the real numbers,
  respectively;

  \noindent
  (ii) the corresponding hom-component of the
  right derived extension-of-scalars functor
  from Lem. \ref{DerivedChangeOfScalars}
  (the operation of ``tensoring a space with $\mathbb{R}$''
  from \cite[Footn. 5]{DGMS75});
\end{defn}

While real-ification \eqref{RealificationOfCoefficients},
in constrast to rationalization \eqref{RationalizationOfCoefficients},
is not directly induced by a {localization of spaces},
it is equivalent to rationalization followed by
derived extension of scalars:
\begin{prop}[Realification is rationalization followed by extension of scalars]
\label{RealificationIsRationalizationFollowedByExtensionOfScalars}
 The operation of real-ification \eqref{RealificationOfCoefficients}
 factors through rationalization \eqref{RationalizationOfCoefficientsInNonabelianCohomology}
 via extension of scalars \eqref{ExtensionOfScalarsFromRationalToRealCohomology}
 in that the following diagram commutes:

 \vspace{-.2cm}
 \begin{equation}
   \label{RealificationFactoringThroughRationalization}
   \begin{tikzcd}[row sep=0pt, column sep=40pt]
     &
     H(X;\, L_{\mathbb{Q}} A)
     \ar[
       dd,
       "(-) \otimes_{{}_{\mathbb{Q}}} \mathbb{R}"
     ]
     \\
     H(X;\, A)
     \ar[
       dr,
       "{
         (\eta^{\mathbb{R}}_A)_\ast
       }"{below}
     ]
     \ar[
       ur,
       "{
         (\eta^{\mathbb{Q}}_A)_\ast
       }"
     ]
     \\
     &
     H(X;\, L_{\mathbb{R}}A)
     \,.
   \end{tikzcd}
 \end{equation}
 \vspace{-.5cm}

\end{prop}
\begin{proof}
Consider the following diagram:
\begin{equation}
  \label{DiagramFactoringRealificationThroughRationalization}
  \begin{tikzcd}
    \mathllap{
      \mbox{
        \tiny
        \color{darkblue}
        \bf
        \def\arraystretch{.9}
        \begin{tabular}{c}
          nonabelian
          \\
          cohomology
        \end{tabular}
      }
      \!\!
    }
    H(X; \, A)
    \ar[
      rr,
      "{
        H(X;\, \mathbb{D}\eta^{\scalebox{.7}{\color{blue}$\mathbb{Q}$}}_{\mathbb{A}})
      }"{description},
      "{
        \mbox{
          \tiny
          \color{greenii}
          \bf
          rational character map
        }
      }"{above, yshift=4pt}
    ]
    \ar[
      drr,
      "{
        \LeftDerived\Omega^\bullet_{\mathrm{P}{\color{blue}\mathbb{Q}}\mathrm{LdR}}
      }"{sloped, below}
    ]
    \ar[
      rrrr,
      rounded corners,
      to path={
           -- ([yshift=+12pt]\tikztostart.north)
           --node[above]{
               \scalebox{.8}{$
                 \overset{
                   \raisebox{4pt}{
                     \scalebox{.8}{
                     \color{greenii}
                     \bf
                     real character map
                     }
                   }
                 }{
                   H
                   \big(
                     X;\,
                     \mathbb{D}\eta^{\scalebox{.66}{\color{blue}$\mathbb{R}$}}_A
                   \big)
                 }
               $}
             } ([yshift=+12pt]\tikztotarget.north)
           -- (\tikztotarget.north)}
    ]
    \ar[
      rrrrd,
      rounded corners,
      to path={
           -- ([yshift=-54pt]\tikztostart.south)
           --node[below]{
               \scalebox{.8}{$
                 \LeftDerived \Omega^\bullet_{\mathrm{P}{\color{blue}\mathbb{R}}\mathrm{LdR}}
               $}
             } ([yshift=-12pt]\tikztotarget.south)
           -- (\tikztotarget.south)}
    ]
    &&
    \overset{
      \mathrlap{
      \mbox{
        \tiny
        \color{darkblue}
        \bf
        nonabelian rational cohomology
      }
      }
    }{
      H(X; L_{\color{blue}\mathbb{Q}} A)
    }
    \ar[
      rr,
      "{
        (-) \otimes_{{}_{\mathbb{Q}}} \mathbb{R}
      }"
    ]
    \ar[
      d,
      "\rotatebox{90}{$\sim$}"{right},
      "\widetilde{(-)}"{left}
    ]
    &&
    H(X; L_{\color{blue}\mathbb{R}} A)
    \mathrlap{
      \!\!\!\!\!
      \mbox{
        \tiny
        \color{darkblue}
        \bf
        \def\arraystretch{.9}
        \begin{tabular}{c}
          nonabelian
          \\
          real cohomology
        \end{tabular}
      }
    }
    \\
    &&
    H
    \big(
      \LeftDerived\Omega^\bullet_{\mathrm{P}{\color{blue}\mathbb{Q}}\mathrm{LdR}}(X),
      \,
      \LeftDerived\Omega^\bullet_{\mathrm{P}{\color{blue}\mathbb{Q}}\mathrm{LdR}}(A)
    \big)
    \ar[
      rr
    ]
    \ar[
      rr,
      "{
        \mathbb{R}
        \big(
          (-) \otimes_{{}_{\mathbb{Q}}} \mathbb{R}
        \big)
      }"{below},
      "{
        \mbox{
          \tiny
          \color{greenii}
          \bf
          \def\arraystretch{.9}
          \begin{tabular}{c}
            extension
            \\
            of scalars
          \end{tabular}
        }
      }"
    ]
    &&
    H
    \big(
      \LeftDerived\Omega^\bullet_{\mathrm{P}{\color{blue}\mathbb{R}}\mathrm{LdR}}(X),
      \,
      \LeftDerived\Omega^\bullet_{\mathrm{P}{\color{blue}\mathbb{R}}\mathrm{LdR}}(A)
    \big)
    \ar[
      u,
      "\rotatebox{90}{$\sim$}"{right},
      "\widetilde{(-)}"{left}
    ]
  \end{tikzcd}
\end{equation}
Here the triangle on the left as well as the outer rectangle commute by
general properties of adjunctions
(the naturality of the hom-isomorphism \eqref{AnAdjunction}
 combined with the definition
\eqref{AdjunctionUnit} of the adjunction unit).
The square on the right commutes by
definition \eqref{ExtensionOfScalarsFromRationalToRealCohomology},
and the bottom part commutes by Prop. \ref{kPLdRAdjunctionFactorsThroughRationalization}.
Together this implies that the top rectangle commutes, which is
the statement to be shown.
\end{proof}

\begin{defn}[Non-abelian character map]
  \label{NonAbelianChernDoldCharacter}
  Let $X, A \in \NilpotentConnectedQFiniteHomotopyTypes$
  (Def. \ref{NilpotentConnectedSpacesOfFiniteRationalType})
  such that $X$ admits the structure of a smooth manifold.
  Then we say that the \emph{non-abelian character map}
  in non-abelian $A$-cohomology
  (Def. \ref{NonAbelianCohomology})
  is the cohomology operation (Def. \ref{RationalizationOfCoefficientsInNonabelianCohomology})
  \vspace{-3mm}
  \begin{equation}
    \label{NonAbelianChernDoldAsCompositeOfRationalizationWithdeRham}
    \mathllap{
      \mbox{
        \tiny
        \color{darkblue}
        \bf
        \def\arraystretch{.9}
        \begin{tabular}{c}
          non-abelian
          \\
          character map
        \end{tabular}
      }
      \;\;
    }
    \mathrm{ch}_A
    \;:\;
    \xymatrix{
      \overset{
        \mathclap{
        \raisebox{3pt}{
          \tiny
          \color{darkblue}
          \bf
          \def\arraystretch{.9}
          \begin{tabular}{c}
            non-abelian
            \\
            cohomology
          \end{tabular}
        }
        }
      }{
      H
      (
        X;
        \,
        A
      )
      }
      \ar[rr]^-{
        (\eta^\mathbb{R}_A)_\ast
      }_-{
        \mbox{
          \tiny
          \color{greenii}
          \bf
          $\mathbb{R}$-rationalization
        }
      }
      &&
      \overset{
        \mathclap{
        \raisebox{3pt}{
          \tiny
          \color{darkblue}
          \bf
          \def\arraystretch{.9}
          \begin{tabular}{c}
            non-abelian
            \\
            real cohomology
          \end{tabular}
        }
        }
      }{
      H
      \big(
        X;
        \,
        L_{\mathbb{R}}A
      \big)
      }
      \ar[rr]^-{ \simeq }_-{
        \mathclap{
        \mbox{
          \tiny
          \color{greenii}
          \bf
          \def\arraystretch{.9}
          \begin{tabular}{c}
            non-abelian
            \\
            de Rham theorem
          \end{tabular}
        }
        }
      }
      &&
   \;  \overset{
        \mathclap{
        \raisebox{3pt}{
          \tiny
          \color{darkblue}
          \bf
          \def\arraystretch{.9}
          \begin{tabular}{c}
            non-abelian
            \\
            de Rham cohomology
          \end{tabular}
        }
        }
      }{
      H_{\mathrm{dR}}
      (
        X;
        \,
        \mathfrak{l}A
      )
      }
    }
  \end{equation}

  \vspace{-3mm}
  \noindent
  from non-abelian $A$-cohomology (Def. \ref{NonAbelianCohomology})
  to non-abelian de Rham cohomology (Def. \ref{NonabelianDeRhamCohomology})
  with coefficients
  in the rational Whitehead $L_\infty$-algebra $\mathfrak{l}A$ of $A$
  (Prop \ref{WhiteheadLInfinityAlgebras}), which is the composite of

  {\bf (i)} the operation  \eqref{RealificationOfCoefficients}
  of real rationalization of coefficients
  (Def. \ref{RationalizationOfCoefficientsInNonabelianCohomology}),

  {\bf (ii)} the equivalence \eqref{EquivalenceBetweenNonabelianRealAndNonaebelianDeRhamCohomology}
  of the non-abelian de Rham theorem
  (Theorem \ref{NonAbelianDeRhamTheorem}).
\end{defn}

Unwinding the definitions
and theorems that go into Def. \ref{NonAbelianChernDoldCharacter},
shows that the non-abelian character
map on a non-abelian cohomology theory with classifying space
a (connected, nilpotent and $\mathbb{Q}$-finite) homotopy type $A$
assigns
flat non-abelian differential form data (Def. \ref{FlatLInfinityAlgebraValuedDifferentialForms})
satisfying the differential relations
of the CE-algebra of the Whitehead $L_\infty$-algebra of $A$
(Prop. \ref{WhiteheadLInfinityAlgebras}):

\begin{example}[Non-abelian character on Cohomotopy theory]
\label{NonAbelianCharacterOnCohomotopyTheory}
The non-abelian character (Def. \ref{NonAbelianChernDoldCharacter})
of

\noindent
{\bf (i)}
a class $[c] \in \pi^n(X) = H^1\big(X; \Omega S^n \big)$ in Cohomotopy (Ex. \ref{CohomotopyTheory})
is (by Ex. \ref{RationalizationOfnSpheres}, Ex. \ref{FlatSphereValuedDifferentialForms})
of this form:

\vspace{-.3cm}
$$
  \mathrm{ch}_{S^n}(c)
  \;=\;
  \left\{
  \def\arraystretch{2}
  \begin{array}{lll}
    \big[
      G_{n}
      \,\in\,
      \Omega^n(X)
      \,\vert\,
      d\, G_n \;= 0\;
    \big]
    & \mbox{if}
    &
    \mbox{$n = 2k+ 1$ is odd}
    \\
    \left[
      \!\!\!
      \left.
      \def\arraystretch{1}
      \begin{array}{l}
        G_{2n-1}
        \,\in\,
        \Omega_{\mathrm{dR}}^{2n-1}(X)\,,
        \\
        G_{n}
        \,\in\,
        \Omega_{\mathrm{dR}}^n(X)
      \end{array}
      \!\right\vert\!
      \def\arraystretch{1}
      \begin{array}{l}
        d\, G_{2n-1} \;=\; - G_n \wedge G_n\,,
        \\
        d\, G_n \;\;\;\;= 0\;
      \end{array}
      \!\!\!
    \right]
    &
    \mbox{if}
    &
    \mbox{$n = 2k$ is even }
  \end{array}
  \right.
$$
\vspace{-.3cm}

\noindent
{\bf (ii)}
a class $[c] \in H^1\big(X; \Omega \mathbb{C}P^n \big)$
in the non-abelian cohomology theory represented by
complex projective $n$-space is (by Ex. \ref{RationalizationOfComplexProjectiveSpace})
of this form:

\vspace{-.5cm}
$$
  \mathrm{ch}_{\mathbb{C}P^n}(c)
  \;=\;
  \Bigg[
    \!\!\!
    \left.
    \begin{array}{l}
      H_{2n+1} \,\in\, \Omega^{2n+1}_{\mathrm{dR}}(X)\,,
      \\[-2pt]
      \;\;F_2\; \,\in\, \Omega^2_{\mathrm{dR}}(X)
    \end{array}
    \,\right\vert\,
    \begin{array}{l}
      d \,H_{2k+1}\,
        =
      \overset{
        \mathclap{
          \mbox{\tiny \rm $n+1$ factors}
        }
      }{
        \overbrace{
          F_2 \wedge \cdots \wedge F_2\,,
        }
      }
      \\
      d \,F_2\, = 0
    \end{array}
    \!\!\!
  \Bigg]
$$
\vspace{-.4cm}

\end{example}

We come back to these new and deeply non-abelian examples in
\cref{CohomotopicalChernCharacter} below. First we now turn attention to verifying that the
non-abelian character map of Def. \ref{NonAbelianChernDoldCharacter}
correctly subsumes more classical structures of differential topology.

%%%%%%%%%%%%%%%%%%%%%%%%%%%%%%%%%%%%%%%%%%%%%%%%%%%%%%%%%%%%%%%%%%%
\subsection{Chern-Dold character}
 \label{TheChernDoldCharacter}
%%%%%%%%%%%%%%%%%%%%%%%%%%%%%%%%%%%%%%%%%%%%%%%%%%%%%%%%%%%%%%%%%%%%%

We prove (Theorem \ref{NonAbelianChernCharacterSubsumesChernDoldCharacter})
that the non-abelian character map
reproduces the Chern-Dold character on generalized cohomology theories
(recalled as Def. \ref{ChernDoldCharacter}) and in particular
the Chern character on topological K-theory (Example \ref{ChernCharacterInKTheory}).

\begin{remark}[Chern-Dold character over the real numbers]
  In view of Prop. \ref{RealificationIsRationalizationFollowedByExtensionOfScalars}
  and Ex. \ref{RealRationalizationOfSpectra}
  we may and will regard Dold's equivalence (Prop. \ref{DoldEquivalence})
  and the Chern-Dold character (Def. \ref{ChernDoldCharacter}) over the real numbers
  instead of over the rational numbers.
  This does not affect the information contained in the character
  but serves to allow, over smooth manifolds,
  for composition with the de Rham isomorphism.
\end{remark}

\begin{prop}[Dold's equivalence {\cite[Cor. 4]{Dold65}\cite[Thm. 3.18]{Hilton71}\cite[\S II.3.17]{Rudyak98}}]
  \label{DoldEquivalence}
  Let $E$ be a generalized cohomology theory
  (Example \ref{GeneralizedCohomologyAsNonabelianCohomology}).
  Then its $\mathbb{R}$-rationalization $E_{\mathbb{R}}$ is equivalent to
  ordinary cohomology with coefficients in the
  rationalized stable homotopy groups of $E$:
  \vspace{-2mm}
  $$
    \xymatrix{
    E^n_{\mathbb{Q}}(X)
    \ar[rr]^-{ \mathrm{do}_E }_-{\simeq}
    &&
    \underset{
      k \in \mathbb{Z}
    }{\bigoplus}
    H^{n + k}
    \big(
      X;
      \,
      \pi_k(E) \otimes_{\scalebox{.5}{$\mathbb{Z}$}} \mathbb{Q}
    \big).
    }
  $$
\end{prop}
\begin{remark}[Rational stable homotopy theory]
  In modern stable homotopy theory, Dold's equivalence
  (Prop. \ref{DoldEquivalence}) is a direct consequence
  of the fundamental theorem \cite[Thm. 5.1.6]{SchwedeShipley01}
  that rational spectra
  are direct sums of Eilenberg-MacLane spectra
  with coefficients in the rationalized stable homotopy groups
  \cite[Prop. 2.17]{BMSS19}.
\end{remark}
But we may explicitly re-derive Dold's equivalence
using the unstable rational homotopy theory from \cref{NonAbelianDeRhamCohomologyTheory}
and passing to rationalization over the real numbers.

\begin{prop}[Dold's equivalence via non-abelian real cohomology]
  \label{DoldEquivalenceViaNonAbelianRealCohomology}
  Let $E$ be a generalized cohomology theory (Example \ref{GeneralizedCohomologyAsNonabelianCohomology})
  and let $n \in \mathbb{N}$ such that the $n$th coefficient space
  \eqref{Spectrum} is of rational finite
  homotopy type (Def. \ref{NilpotentConnectedSpacesOfFiniteRationalType})
  $
    E_n
    \;\in\;
    \QFiniteHomotopyTypes
    \,.
  $
  Then there is a natural equivalence between
  the non-abelian real cohomology (Def. \ref{NonAbelianRealCohomology})
  with coefficients in $E_n$
  and ordinary
  cohomology with coefficients in the $\mathbb{R}$-rationalized homotopy groups of $E$:
  \vspace{0mm}
  \begin{equation}
    \label{DoldEquivalenceInNonabelianRealCohomology}
    H
    \big(
      -;
      \,
      L_{\mathbb{R}}E_n
    \big)
    \;\;
    \simeq
    \;\;
    \underset{
      k \in \mathbb{N}
    }{\bigoplus}
    \,
    H^{n+k}
    \big(
      -;
      \,
      \pi_{k}(E) \otimes_{\scalebox{.5}{$\mathbb{Z}$}} \mathbb{R}
    \big)
    \,.
  \end{equation}
\end{prop}

\vspace{-4mm}
$\,$
\begin{proof}
  Since $E_n$ is an infinite-loop space, it
  is nilpotent (Example \ref{ExamplesOfNilpotentSpaces}).
  We may assume without restriction that it is also connected, for
  otherwise we apply the following argument to each connected component
  (Remark \ref{AssumptionsOnConnctedNilpotentRFiniteHomotopyTypes}).
  Hence $E_n \in \NilpotentConnectedQFiniteHomotopyTypes$
  (Def. \ref{NilpotentConnectedSpacesOfFiniteRationalType})
  and the discussion in \cref{RationalHomotopyTheory} applies.
Again, since $E_n$ is a loop space \eqref{Spectrum},
  Prop. \ref{NonAbelianRealCohomologyWithCoefficientsInLoopSpaces}
  gives
  $
    H
    (
      -;
      \,
      L_{\mathbb{R}}E_n
    )
    \simeq
    \underset{
      k \in \mathbb{N}
    }{\oplus}
    H^{k}
    (
      -;
      \,
      \pi_{k}(E_n) \otimes_{\scalebox{.5}{$\mathbb{Z}$}} \mathbb{R}
    )
    \,.
  $
  The claim follows from the definition of stable homotopy groups
  as $\pi_{k-n}(E) \,=\, \pi_k(E_n)$ for $k,n \geq 0$,
  (as $E$ is an ``$\Omega$-spectrum'' \eqref{Spectrum}).
\end{proof}

\begin{defn}[Real Chern-Dold character {\cite{Buchstaber70}\cite[p. 50]{Hilton71}}]
  \label{ChernDoldCharacter}
  Let $E$ be a generalized cohomology theory (Example \ref{GeneralizedCohomologyAsNonabelianCohomology}).
  The real \emph{Chern-Dold character} in $E$-cohomology theory is the
  cohomology operation to ordinary cohomology
  which is the composite of
  rationalization
  in $E$-cohomology
  with Dold's equivalence (Prop. \ref{DoldEquivalence}):
\vspace{-2mm}
  \begin{equation}
    \label{ChernDoldCharacterAsComposite}
    \overset{
      \mathclap{
      \raisebox{3pt}{
        \tiny
        \color{darkblue}
        \bf
        \def\arraystretch{.9}
        \begin{tabular}{c}
          Chern-Dold
          \\
          character
        \end{tabular}
      }
      }
    }{
      \mathrm{ch}_E
    }
    \;:\;
    \xymatrix@C=3em@R=2em{
    E^\bullet(-)
    \ar[rr]^-{
      \overset{
        \mathclap{
        \raisebox{3pt}{
          \tiny
          \color{darkblue}
          \bf
          \def\arraystretch{.9}
          \begin{tabular}{c}
            $\mathbb{R}$-rationalization
            \\
            in $E$-cohomoloy
          \end{tabular}
        }
        }
      }{
      }
    }
    \ar[d]_-{
                \mbox{
           \tiny
           \color{darkblue}
           \bf
           \eqref{GeneralizedECohomologyAsNonabelianCohomology}
         }
       }^-{
       \scalebox{.8}{{$\simeq$}}
    }
    &&
    E^\bullet_{\mathbb{R}}(-)
    \ar[rr]_-{ \simeq }^-{
      \overset{
        \mathclap{
        \raisebox{3pt}{
          \tiny
          \color{darkblue}
          \bf
          \def\arraystretch{.9}
          \begin{tabular}{c}
            Dold's equivalence
          \end{tabular}
        }
        }
      }{
        \mathrm{do}_E
      }
    }
    \ar[d]^-{\!\!
                \mbox{
           \tiny
           \color{darkblue}
           \bf
           \eqref{GeneralizedECohomologyAsNonabelianCohomology}
         }
       }_-{
       \scalebox{.8}{{$\simeq$}}
    }    &&
    \underset{k}{\bigoplus}
    \,
    H^{\bullet + k}
    \big(
      -;
      \pi_k(E) \otimes_{\scalebox{.5}{$\mathbb{Z}$}} \mathbb{R}
    \big)
    \\
    H
    (
      -;
      E_\bullet
    )
    \ar[rr]^-{
     \left(
        \eta^{\mathbb{R}}_{E_\bullet}
      \right)_\ast
    }_-{
      \mbox{
        \tiny
        \color{darkblue}
        \bf
        \eqref{RealificationOfCoefficients}
      }
    }
    &&
    H
    (
      -;
      L_{\mathbb{R}} E_\bullet
    )
    \ar[urr]^-{
      \scalebox{.8}{\rotatebox[origin=c]{19}{$\simeq$}}
    }_-{
      \mbox{
        \tiny
        \color{darkblue}
        \bf
        \eqref{DoldEquivalenceInNonabelianRealCohomology}
      }
    }
    }
    \,.
  \end{equation}

  \vspace{-2mm}
  \noindent  Here the bottom part
  serves to make the nature of the
  top maps fully explicit,
  using Example \ref{GeneralizedCohomologyAsNonabelianCohomology},
  Def. \ref{RationalizationOfCoefficientsInNonabelianCohomology}
  and Prop. \ref{DoldEquivalenceViaNonAbelianRealCohomology}.
\end{defn}

\begin{remark}[Rationalization in the Chern-Dold character]
  That the first map in the Dold-Chern character \eqref{ChernDoldCharacterAsComposite}
  is the rationalization localization
  (here shown exended to the real numbers)
  is stated
  somewhat indirectly in the original definition \cite{Buchstaber70}
  (the concept of rationalization was fully formulated later in \cite{BousfieldKan72}).
  The role of rationalization in the Chern-Dold character
  is made fully explicit in \cite[\S 2.1]{LindSatiWesterland16}.
  The same rationalization construction
  of the generalized Chern character, but without attribution to \cite{Buchstaber70} or \cite{Dold65},
  is considered in
  \cite[\S 4.8]{HopkinsSinger05} (see also \cite[p. 17]{BunkeNikolaus14}).
\end{remark}

We now come to the main result in this section:

\begin{theorem}[Non-abelian character subsumes Chern-Dold character]
  \label{NonAbelianChernCharacterSubsumesChernDoldCharacter}
  Let $E$ be a generalized cohomology theory (Example \ref{GeneralizedCohomologyAsNonabelianCohomology})
  and let $n \in \mathbb{N}$ such that the $n$th coefficient space
  \eqref{Spectrum} is of rational finite
  homotopy type (Def. \ref{NilpotentConnectedSpacesOfFiniteRationalType}).
  Let moreover
  $X$ be a smooth manifold.

  \vspace{0mm}
  \noindent  Then
  the non-abelian character (Def. \ref{NonAbelianChernDoldCharacter})
  coincides with the Chern-Dold character (Def. \ref{ChernDoldCharacter})
  on $E$-cohomology in degree $n$, in that the following
  diagram commutes:
   \vspace{-1mm}
  \begin{equation}
    \label{NonabelianChernCharacterRelatedToChernDoldCharacter}
    \raisebox{20pt}{
    \xymatrix@R=14pt{
      H
      \big(
        X;
        E_n
      \big)
      \ar[rr]_-{ \mathrm{ch}_{E_n} }
      \ar@{<-}[d]_-{
          \mbox{
            \tiny
            \color{darkblue}
            \bf
            \eqref{GeneralizedECohomologyAsNonabelianCohomology}
          }
          }^-{\simeq}
          &&
      H_{\mathrm{dR}}
      \big(
        \mathfrak{l}
        E_n
      \big)
      \ar[d]_-{\simeq}^-{\!\!\!
               \mathrlap{
          \;
          \mbox{
            \tiny
            \color{darkblue}
            \bf
            \eqref{EquivalenceBetweenNonabelianRealAndNonaebelianDeRhamCohomology}
            \eqref{DoldEquivalenceInNonabelianRealCohomology}
          }
        }
      }
      \\
      E^n(X)
      \ar[rr]_-{
        (\mathrm{ch}_E)^n
      }
      &&
      \underset{ k }{\bigoplus}
      \,
      H^{n+k}
      \big(
        X;
        \,
        \pi_{k}(E)
          \otimes_{\scalebox{.5}{$\mathbb{Z}$}}
        \mathbb{R}
      \big).
    }
    }
  \end{equation}

  \vspace{-2mm}
  \noindent
  Here the equivalence on the left is from Example \ref{GeneralizedCohomologyAsNonabelianCohomology},
  while the equivalence on the right is the inverse
  non-abelian de Rham theorem (Theorem \ref{NonAbelianDeRhamTheorem})
  composed with that from Prop. \ref{DoldEquivalenceViaNonAbelianRealCohomology}.
\end{theorem}
\begin{proof}
  Since $E_n$ is an infinite-loop space, it
  is necessarily nilpotent (Example \ref{ExamplesOfNilpotentSpaces}).
  We may assume without restriction that it is also connected, for
  otherwise we apply the following argument to each connected component
  (Remark \ref{AssumptionsOnConnctedNilpotentRFiniteHomotopyTypes}).
  Hence $E_n \in \NilpotentConnectedQFiniteHomotopyTypes$
  (Def. \ref{NilpotentConnectedSpacesOfFiniteRationalType})
  and the discussion in \cref{RationalHomotopyTheory}
  and \cref{NonAbelianDeRhamTheory} applies:

  The non-abelian de Rham isomorphism \eqref{EquivalenceBetweenNonabelianRealAndNonaebelianDeRhamCohomology}
  in the definition \eqref{NonAbelianChernDoldAsCompositeOfRationalizationWithdeRham} of the non-abelian character cancels against its inverse on the
  right of \eqref{NonabelianChernCharacterRelatedToChernDoldCharacter}.
  Commutativity of the remaining diagram
   \vspace{-2mm}
  $$
    \xymatrix@R=14pt{
      H
      \big(
        X;
        \,
        E_n
      \big)
      \ar[rr]_-{
        (\eta^{\mathbb{R}}_{E_n})_\ast
      }
      \ar@{<-}[d]_-{
        \mathllap{
          \mbox{
            \tiny
            \color{darkblue}
            \bf
            \eqref{GeneralizedECohomologyAsNonabelianCohomology}
          }}
  }^-{
  \simeq}
            &&
      H
      \big(
        X;
        \,
        L_{\mathbb{R}}E_n
      \big)
      \ar[d]_-{\simeq}^-{\!\!
        \mathrlap{
          \;
          \mbox{
            \tiny
            \color{darkblue}
            \bf
            \eqref{DoldEquivalenceInNonabelianRealCohomology}
          }
        }
      }
      \\
      E^n(X)
      \ar[rr]_-{ \mathrm{ch}_{E_n}}
      &&
      \underset{ k }{\bigoplus}
      \,
      H^{n+k}
      \big(
        X;
        \,
        \pi_{k}(E)
          \otimes_{\scalebox{.5}{$\mathbb{Z}$}}
        \mathbb{R}
      \big)
    }
  $$

\vspace{-2mm}
\noindent  is the very definition of the Chern-Dold character
  (Def. \ref{ChernDoldCharacter}).
\end{proof}

%\begin{remark}[Generality of the non-abelian character on spectra]
%  $\,$
%
%  \noindent
%  {\bf i)}
%  The Dold-Chern character (Def. \ref{ChernDoldCharacter})
%  applies to any
%  generalized cohomology theory, given by any spectrum.
%  Of course it may be trivial if the given spectrum has
%  no rational homotopy.
%
%  \noindent
%  {\bf ii)}
%  The restriction to classifying spaces
%  of $\mathbb{R}$-finite homotopy type in
%  Theorem \ref{NonAbelianChernCharacterSubsumesChernDoldCharacter}
%  is only so as to make the non-abelian de Rham theorem
%  apply with the assumptions as stated there (Theorem %\ref{NonAbelianDeRhamTheorem}).
%\end{remark}

\begin{example}[de Rham homomorphism in ordinary cohomology]
  \label{deRhamHomomorphism}
  On ordinary integral cohomology (Example \ref{OrdinaryCohomology}),
  the non-abelian character
  (Def. \ref{NonAbelianChernDoldCharacter}) reduces to
  extension of scalars from the integers to the real numbers,
  followed by the de Rham isomorphism,
  in that the following diagram commutes:
  \vspace{-2mm}
  $$
    \xymatrix@C=4em{
      H
      \big(
        -;
        \,
        B^{n+1}\mathbb{Z}
      \big)
      \ar[rr]_-{\mathrm{ch}_{B^{n+1}\mathbb{Z}}}^-{
        \overset{
          \mathclap{
          \raisebox{3pt}{
            \tiny
            \color{darkblue}
            \bf
            \def\arraystretch{.9}
            \begin{tabular}{c}
              non-abelian character
              \\
              on ordinary cohomology
            \end{tabular}
          }
          }
        }{
  %        \mathrm{ch}_{B^{n+1}\mathbb{Z}}
        }
      }
      \ar@{<-}[d]_-{
                  \mbox{
            \tiny
            \color{darkblue}
            \bf
            \eqref{OrdinaryCohomologyAsNonAbelianCohomology}
          }}^-{
        \simeq
      }
      &&
      H_{\mathrm{dR}}
      \big(
        -;
        \,
        \mathfrak{l}B^{n+1}\mathbb{Z}
      \big)
      \ar[d]_-{\simeq}^-{\!
                \mbox{
            \tiny
            \color{darkblue}
            \bf
            \eqref{WhiteheadLInfinityAlgebraOfIntegralEMSpace}
            \eqref{NonAbelianDeRhamCohomologyWithCoefficientsInLineLienAlgebra}
          }
        }
      \\
      H^{n+1}
      (
        -;
        \mathbb{Z}
      )
      \ar[r]_-{
        \mathclap{
        \mbox{
          \tiny
          \color{darkblue}
          \bf
          \def\arraystretch{.9}
          \begin{tabular}{c}
            extension
            \\
            of scalars
          \end{tabular}
        }
        }
      }
      &
      H^{n+1}
      (
        -;
        \mathbb{R}
      )
      \ar[r]^-{\simeq}_-{
        \mathclap{
        \mbox{
          \tiny
          \color{darkblue}
          \bf
          \def\arraystretch{.9}
          \begin{tabular}{c}
            ordinary
            \\
            de Rham isomorphism
          \end{tabular}
        }
        }
      }
      &
      H^{n+1}_{\mathrm{dR}}
      (
        -
      )
    }
  $$
\end{example}
%This follows as
%
%\noindent {\bf (a)} the special case of the relation to the
%Chern-Dold character (Theorem %\ref{NonAbelianChernCharacterSubsumesChernDoldCharacter} below)
%with ordinary cohomology regarded as the generalized cohomology
%represented by Eilenberg-MacLane spectra;
%
%or alternatively
%
%\noindent {\bf (b)} as in the proof of the relation to the
%Chern-Weil homomorphism (Theorem \ref{NonAbelianChernDoldSubsumesChernWeil} %below), with $G := B^n \mathbb{Z} \simeq B^{n-1} \mathrm{U}(1)$
%as in Example \ref{HigherBundleGerbes}.

\begin{example}[Chern character on complex K-theory]
  \label{ChernCharacterInKTheory}
  The spectrum \eqref{Spectrum}
  representing complex K-theory
  has 0th component space
  $\mathrm{KU}_0 \simeq \mathbb{Z} \times B \mathrm{U}$
  \eqref{ClassifyingSpaceForComplexTopologicalKTheory}.
  Here the connected components $B \mathrm{U}$,
  the classifying space of the infinite unitary group
  \eqref{ClassifyingSpaceOfInfiniteUnitaryGroup},
  are clearly
  of finite rational type (since their rational cohomology is the
  ring of universal Chern classes, e.g. \cite[Thm. 2.3.1]{Kochman96}).
  Therefore, Theorem \ref{NonAbelianChernCharacterSubsumesChernDoldCharacter}
  applies and says that the
  non-abelian character map (Def. \ref{NonAbelianChernDoldCharacter})
  for coefficients in $\mathbb{Z} \times B \mathrm{U}$
  reduces to the Chern-Dold character on complex K-theory.
  This, in turn, is equivalent (by \cite[Thm. 5.8]{Hilton71})
  to the original Chern character $\mathrm{ch}$
  on complex K-theory \cite[\S 12.1]{Hirzebruch56}\cite[\S 9.1]{BorelHirzebruch58}\cite[\S 1.10]{AtiyahHirzebruch61}
  (review in \cite[\S V]{Hilton71}):
  \vspace{-3mm}
%  $$
%    \overset{
%      \mathclap{
%      \;\;\;\;\;\;
%      \raisebox{3pt}{
%        \tiny
%        \color{darkblue}
%        \bf
%        \begin{tabular}{c}
%          Chern character on
%          \\
%          complex K-theory
%        \end{tabular}
%      }
%      }
%    }{
%      \mathrm{ch}
%    }
%    \;\;\;\;\;\;\;\; \simeq \;\;\;\;
%    \mathrm{ch}_{ \mathbb{Z} \times B \mathrm{U} }
%    \,.
%  $$
\begin{equation}
  \label{OrdinaryChernCharacterReproduced}
  \xymatrix@R=-2pt@C=5em{
    \overset{
      \mathclap{
      \raisebox{3pt}{
        \tiny
        \color{darkblue}
        \bf
        \def\arraystretch{.9}
        \begin{tabular}{c}
          non-abelian cohomology
          \\
          with $\mathbb{Z} \times B \mathrm{U}$-coefficients
        \end{tabular}
      }
      }
    }{
      H(
        X;
        \,
        \mathbb{Z} \times B \mathrm{U}
      )
    }
    \ar[rr]_-{
      \mathrm{ch}_{\mathbb{Z} \times B \mathrm{U}}
    }^-{
      \mbox{
        \tiny
        \color{greenii}
        \bf
        \def\arraystretch{.9}
        \begin{tabular}{c}
          non-abelian character map
          \\
          with $\mathbb{Z} \times B \mathrm{U}$-coefficients
        \end{tabular}
      }
    }
    \ar@{=}[ddd]
    &&
    \overset{
      \mathclap{
      \raisebox{3pt}{
        \tiny
        \color{darkblue}
        \bf
        \def\arraystretch{.9}
        \begin{tabular}{c}
          non-abelian de Rham homology
          \\
          with $\mathfrak{l}(\mathbb{Z} \times B \mathrm{U})$-coefficients
        \end{tabular}
      }
      }
    }{
    H_{\mathrm{dR}}
      \big(
        X;
        \,
        \mathfrak{l}
        (
        \mathbb{Z} \times B \mathrm{U}
        )
      \big)
    }
    \ar@{=}[ddd]
    \\
    \\
    \\
    \mathrm{KU}^0(X)
    \ar[rr]^-{ \mathrm{ch} }_-{
      \mbox{
        \tiny
        \color{greenii}
        \bf
        traditional Chern character
      }
    }
    &&
    \underset{k \in \mathbb{Z}}{\bigoplus}
    \,
    H^{2k}_{\mathrm{dR}}(X)
    \\
    [\nabla]
    \ar@{}[rr]|-{ \longmapsto }
    &&
    \mathclap{
    \big[
    \mathrm{tr}
    \circ
    \mathrm{exp}
    \left(
      \frac
        {i \, F_\nabla}
        { 2\pi }
    \right)
    \big]
    }
    \,.
  }
\end{equation}
\vspace{-.4cm}

On the bottom we are showing the classical component-formula,
which to the K-theory class of a complex vector bundle with any choice
of connection $\nabla$ assigns the
de Rham cohomology class of the
trace of the exponential series
of its curvature 2-form \eqref{CurvatureDifferentialForm}.
(More on this {\it Chern-Weil formalism} in \cref{TheChernWeilHomomorphism}).
\end{example}

\begin{example}[Pontrjagin character on real K-theory]
  \label{PontrjaginCharacterInKSO}
  The \emph{Pontrjagin character} $\mathrm{ph}$ on
  real topological K-theory
  (see \cite[\S 9.4]{GHV73}\cite{IaokaKuwana99}\cite{Igusa08}\cite[\S 2.1]{GS-KO})
  is defined
  to be the composite
  \vspace{-2mm}
  $$
  \begin{tikzcd}
    \mathrm{KSpin}^\bullet(-)
    \ar[r]
    &
    \mathrm{KSO}^\bullet(-)
    \ar[
      r
    ]
    \ar[
      rrr,
      rounded corners,
      to path={
           -- ([yshift=+8pt]\tikztostart.north)
           --node[above]{
               \scalebox{.7}{
                 $\widetilde{\mathrm{ph}}^\bullet$
               }
             } ([yshift=+8pt]\tikztotarget.north)
           -- (\tikztotarget.north)}
    ]
    &
    \mathrm{KO}^\bullet(-)
    \ar[
      d,
      "\mathrm{cplx}"
    ]
    \ar[
      rr,
      "\mathrm{ph}^\bullet"
    ]
    &&
    \underset{k}{\bigoplus}
    \,
    H^{\bullet + 4k}
    \big(
      -
      ;\,
      \mathbb{R}
    \big)
    \ar[
      d
    ]
    \\
    &
    &
    \mathrm{KU}^\bullet(-)
    \ar[
      rr,
      "\mathrm{ch}^\bullet"
    ]
    &&
    \underset{k}{\bigoplus}
    \,
    H^{\bullet + 2k}
    \big(
      -
      ;\,
      \mathbb{R}
    \big)
  \end{tikzcd}
  $$

  \vspace{-2mm}
  \noindent  of the complexification map (on representing virtual vector bundles)
  with the Chern character $\mathrm{ch}$ on complex K-theory (Example \ref{ChernCharacterInKTheory}).

  \noindent {\bf (i)} By naturality of the complexification map and since the
  complex Chern character is a Chern-Dold character
  (by \cite[Thm. 5.8]{Hilton71}), so is the Pontrjagin character, as well as its restriction
  $\widetilde {\mathrm{ph}}$ to oriented real K-theory $\mathrm{KSO}$
  and further to ph on KO-theory and to $\mathrm{Spin}$ K-theory, etc.

 \noindent {\bf (ii)}  The connected components $B \mathrm{O}$
  of the classifying space $\mathrm{K\mathrm{O}}_0$
  for real topological K-theory
  are of finite $\mathbb{R}$-type
  (since the real cohomology is the ring of universal Pontrjagin classes).
  Therefore, Theorem \ref{NonAbelianChernCharacterSubsumesChernDoldCharacter}
  applies and
  says that the
  non-abelian Chern character (Def. \ref{NonAbelianChernDoldCharacter})
  for coefficients in $\mathbb{Z} \times B \mathrm{SO}$
  coincides with the Pontrjagin character $\widetilde{\mathrm{ph}}$ in $\mathrm{KSO}$-theory:
  \vspace{-2mm}
  $$
    \overset{
      \mathclap{
      \raisebox{3pt}{
        \tiny
        \color{darkblue}
        \bf
        \def\arraystretch{.9}
        \begin{tabular}{c}
          Pontrjagin character
          \\
          on oriented real K-theory
        \end{tabular}
      }
      }
    }{
      \widetilde{\mathrm{ph}}
    }
    \;\;\;\;\;\; \;\;\;\simeq \;\;\;\;
    \mathrm{ch}_{\mathbb{Z} \times B \mathrm{SO}}
    \,.
  $$

 \noindent {\bf (iii)}   By Remark \ref{AssumptionsOnConnctedNilpotentRFiniteHomotopyTypes}, the construction
  extends to the Pontrjagin character $\mathrm{ph}$ on KO-theory.

 \noindent {\bf (iv)}  The same applies to the further restriction of the
  Pontrjagin character  to $\mathrm{KSpin}$; see \cite{LiDuan}\cite{Thomas}
  for some subtleties involved and  \cite[\S 7]{KSpin} for  interpretation and applications.
\end{example}

\begin{example}[Chern-Dold character on Topological Modular Forms]
\label{CharacterInTMF}
The connective ring spectrum $\mathrm{tmf}$
of \emph{topological modular forms}
\cite[\S 9]{Hopkins94}\cite[\S 4]{Hopkins02} (see \cite{DFHH14})
is, essentially by design, such that
under rationalization it yields the graded ring of rational modular forms
(e.g \cite[p. 2]{DouglasHenriques11}):
\vspace{-3mm}
\begin{equation}
  \label{RationalizationOfTopologicalModularForms}
  \xymatrix@C=3em{
    \overset{
      \mathclap{
      \raisebox{3pt}{
        \tiny\color{darkblue}
        \bf
        \def\arraystretch{.9}
        \begin{tabular}{c}`
          topological
          \\
          modular forms
        \end{tabular}
      }
      }
    }{
      \pi_\bullet(\mathrm{tmf})
    }
    \;
    \ar[rr]^-{
      (-) \otimes_{\scalebox{.5}{$\mathbb{Z}$}} \mathbb{R}
    }
    &&
    \;\;\;
    \overset{
      \mathclap{
      \raisebox{3pt}{
        \tiny
        \color{darkblue}
        \bf
        \def\arraystretch{.9}
        \begin{tabular}{c}
          rational
          \\
          modular forms
        \end{tabular}
      }
      }
    }{
      \mathrm{mf}^{\, \mathbb{R}}_\bullet
    }
    \;\simeq\;
    \mathbb{R}
    \big[\,
      \overset{
        \scalebox{.5}{$
          \mathrm{deg} = 8
        $}
      }{
        \overbrace{
          c_4
        }
      }
      ,
      \overset{
        \scalebox{.5}{$
          \mathrm{deg = 12}
        $}
      }{
      \overbrace{
        c_6
      }
      }
    \,\big].
  }
\end{equation}

\vspace{-2mm}
\noindent
It follows that the Chern-Dold character (Def. \ref{ChernDoldCharacter})
on $\mathrm{tmf}$ takes values in real cohomology with coefficients
in modular forms
\vspace{-2mm}
\begin{equation}
  \label{ChernDoldCharacterIntmf}
  \xymatrix@C=4.5em{
    \mathrm{tmf}^\bullet(-)
    \ar[rr]_-{
        \mathrm{ch}^\bullet_{\mathrm{tmf}}
      }^-{
        \raisebox{3pt}{
          \tiny
          \color{darkblue}
          \bf
          \def\arraystretch{.9}
          \begin{tabular}{c}
            Chern-Dold character
            \\
            on topological modular forms
          \end{tabular}
        }
      }
    &&
    H^{\bullet}
    \big(
      -;
      \mathrm{mf}^{\, \mathbb{R}}_{\bullet}
    \big).
  }
\end{equation}
(This is often considered over the rational numbers,
sometimes over the complex numbers \cite[Fig. 1]{Berwick-Evans13};
we may just as well stay over the real numbers, by Remark
\ref{RationalHomotopyTheoryOverTheRealNumbers}, to retain
contact to the de Rham theorem.)

By Theorem \ref{NonAbelianChernCharacterSubsumesChernDoldCharacter},
this is an instance of the non-abelian character map:
\vspace{-2mm}
$$
  \overset{
    \mathclap{
    \raisebox{3pt}{
      \tiny
      \color{darkblue}
      \bf
      \def\arraystretch{.9}
      \begin{tabular}{c}
        Chern-Dold character on
        \\
        topological nodular forms
      \end{tabular}
    }
    }
  }{
    \mathrm{ch}^{\bullet}_{\mathrm{tmf}}
  }
  \;\;\;\;\;\;\;\;\;\;
  \simeq
  \;\;\;\;
  \mathrm{ch}_{ \mathrm{tmf}_{\bullet} }\;.
$$
\end{example}
\begin{example}[The Hurewicz/Boardman homomorphism on topological modular forms]
  \label{TheBoardmanHomomorphismIntmf}
  The spectrum $\mathrm{tmf}$ (Example \ref{CharacterInTMF})
  carries the structure of an $E_\infty$-ring spectrum (Ex. \ref{GeneralizedCohomologyAsNonabelianCohomology})
  and hence receives an essentially unique homomorphism of
  ring spectra from the sphere spectrum:
  \vspace{-2mm}
  $$
    \xymatrix{
      \Sigma^\infty
      S^0
      \;=\;
      \mathbb{S}
      \ar[rr]^-{ e_{\mathrm{tmf}} }
      &&
      \mathrm{tmf}
      \,.
    }
  $$

  \vspace{-2mm}
  \noindent
  This is also known as the \emph{Hurewicz homomorphism}
  or rather the \emph{Boardman homomorphism}
  (e.g. \cite[\S II.7]{Adams74}\cite[\S 4.3]{Kochman96})
  for $\mathrm{tmf}$.
  The Boardman homomorphism on $\mathrm{tmf}$
  happens to be a stable weak
  equivalence in degrees $\leq 6$,
  in that it is an isomorphism on stable homotopy groups in these
  degrees  \cite[Prop. 4.6]{Hopkins02}\cite[\S 13]{DFHH14}:
    \vspace{-2mm}
  $$
    \xymatrix{
      \pi^s_{\bullet \leq 6}
      \;=\;
      \pi_{\bullet \leq 6}(\mathbb{S})
      \ar[rr]^-{
        \pi_{\bullet\leq 6}(e_{\mathrm{tmf}})
      }_-{ \simeq }
      &&
      \pi_{\bullet \leq 6}(\mathrm{tmf})
      \,.
    }
  $$

  \vspace{-2mm}
  \noindent
  Hence (by Prop. \ref{HomotopyClassesOfMapsOutOfnManifolds})
  when $X^{10}$ is a manifold of dimension
  $\mathrm{dim}(X) \leq 6 + 4 = 10$,
  then the Boardman homomorphism identifies the
  stable Cohomotopy (Example \ref{StableCohomotopy})
  of $X^{10}$
  in degree 4 with $\mathrm{tmf}^{\, 4}\big(X^{10}\big)$:
  \vspace{-2mm}
  \begin{equation}
    \label{BoardmanHomomorphismintmfOn10fold}
    \xymatrix@R=1em@C=3.5em{
      \overset{
        \mathclap{
        \raisebox{3pt}{
          \tiny
          \color{darkblue}
          \bf
          {
          \def\arraystretch{.9}
          \begin{tabular}{c}
            stable
            \\
            4-Cohomotopy
          \end{tabular}}
        }
        }
      }{
      \pi^4_s
      \big(
        X^{10}
      \big)
      }
      \ar[dr]_{ \mathrm{ch}_{\mathbb{S}^4} }
      \;=\;
      \mathbb{S}^4
      (
        X^{10}
      )
      \ar[rr]^-{
        \overset{
          \raisebox{3pt}{
            \tiny
            \color{darkblue}
            \bf
            Boardman homomorphism
          }
        }{
          e^4_{\mathrm{tmf}}
        }
      }_-{ \simeq }
      &&
      \overset{
        \raisebox{3pt}{
          \tiny
          \color{darkblue}
          \bf
          {
          \def\arraystretch{.9}
          \begin{tabular}{c}
            $\mathrm{tmf}$-cohomology
            \\
            in degree 4
          \end{tabular}}
        }
      }{
        \mathrm{tmf}^{\,4}
        (
          X^{10}
        )\;.
      }
      \ar[dl]^{ \mathrm{ch}_{\mathrm{tmf}^{\, 4}} }
      \\
      &
      H^4_{\mathrm{dR}}
      (
        X^{10}
      )
    }
  \end{equation}

  \vspace{-2mm}
  \noindent
  In this situation, the character map from Example \ref{CharacterInTMF}
  extracts exactly the datum of a real 4-class.
\end{example}

\begin{remark}[Clarifying the role of $\mathrm{tmf}$ in string theory]
  \label{ClarifyingTheRoleOfTmfInStringTheory}
  Ever since the famous computation of \cite{Witten87}
  (following \cite{SchellekensWarner86} \cite{SchellekensWarner87})
  showed that the partition function of a 2d super-conformal field theory
  lands in modular forms,
  and since the theorem of \cite{AHS01}\cite{AHR10} showed
  that, mathematically, this statement lifts through (what we call above) the
  $\mathrm{tmf}$-Chern-Dold character \eqref{ChernDoldCharacterIntmf},
  there have been proposals about a possible role of
  $\mathrm{tmf}$-cohomology theory in controlling
  elusive aspects of string theory
  (see
  \cite{KS3}\cite{tcu}\cite{DouglasHenriques11}\cite{StolzTeichner11}\cite{Sa-tert}\cite{GJF18}\cite{GPPV18}\cite{Sa-6d}).
  While good progress has been made, it might be fair to say that the situation has remained inconclusive.

\vspace{1mm}
  \noindent
  {\bf (i) Non-abelian enhancement of $\mathrm{tmf}^4(X^{10})$.}
  But with the non-abelian generalization
  (Def. \ref{NonAbelianChernDoldCharacter})
  of the Chern-Dold character
  in hand, we may ask for a
  non-abelian enhancement
  (Example \ref{StabilizationCohomologyOperation})
  of $\mathrm{tmf}$-theory on string background spacetimes.
  By Example \ref{TheBoardmanHomomorphismIntmf}, this is,
  in degree 4, equivalent to asking for a non-abelian enhancement
  of stable Cohomotopy theory
  (Example \ref{NonabelianEnhancementOfStableCohomotopy}).
  This exists (not uniquely but) canonically:
  given by actual Cohomotopy theory
  (Example \ref{CohomotopyTheory}).
  We work out the non-abelian character map on twisted 4-Cohomotopy
  in
  Example \ref{CharacterMapOnJTwistedCohomotopyAndTwistorialCohomotopy}
  below.
  The concluding Prop. \ref{ChargeQuantizationInJTwistedCohomotopy}
  shows that this
  %non-abelian lift of $\mathrm{tmf}$-theory
  does capture crucial non-linear phenomena of non-perturbative string theory.

\vspace{1mm}
  \noindent
  {\bf (ii) Non-Torsion classes in $\mathrm{tmf}^\bullet$.}
  Part of the statement \eqref{BoardmanHomomorphismintmfOn10fold} is that
  the higher non-torsion generators \eqref{RationalizationOfTopologicalModularForms}
  of $\pi_\bullet(\mathrm{tmf})$
  (hence the actual or ``non-topological'' modular forms)
  do not contribute to $\mathrm{tmf}^4$ on 10-manifolds:
  These start to contribute only on manifolds of dimensionl $ 4 + \mathrm{deg}(c_4) = 12$,
  where, in
  string theory language,
  one computes not fluxes of fields but their (Green-Schwarz-)anomaly densities.
  Indeed, the
  original computation of what came to be known as the ``Witten genus''
  interprets it as the generating function for just these anomalies
  \cite{SchellekensWarner86}\cite{SchellekensWarner87}\cite{LNSW88}\cite{Sati08Anomalies}.
  While the character map \eqref{ChernDoldCharacterIntmf} still applies
  in these higher dimensions,
  the non-abelian enhancement by Cohomotopy is restricted
  exactly to dimension 10,
  and is what makes the character pick up just those non-linear relations,
  discussed in \cref{CohomotopicalChernCharacter},
  that are expected to cancel the anomalies \cite{FSS19b}\cite{FSS20a}\cite{SS20c}.

\vspace{1mm}
  \noindent
  {\bf  (iii) Torsion classes in $\mathrm{tmf}^\bullet$.}
  Indeed, the deep motivation behind topological modular forms
  is the suggestion that these capture mathematical aspects of
  2d supersymmetric field theories even in their non-rational torsion
  elements --  and the beauty of \eqref{BoardmanHomomorphismintmfOn10fold} is
  to show that in the relevant degrees and dimensions these aspects
  are equivalently seen in Cohomotopy.
  Concretely, a famous conjecture orginating with
  \cite{StolzTeichner11}\cite{DouglasHenriques11}
  and cast in more pronounced form in
  \cite[\S 5]{GJF18}
  says that the elements of $\mathrm{tmf}^\bullet(X)$
  correspond bijectively to, roughly,
  the deformation classes of 2-dimensional supersymmetric field
  theories with target space $X$.
  Specifically the torsion elements in
  $\pi_3(\mathrm{tmf}) \,\simeq\, \pi_3(\mathbb{S}) \,\simeq\,\mathbb{Z}/24$,
  have, conjecturally, been identified,
  with certain supersymmetric $\mathrm{SU}(2)$-WZW models
  \cite[p. 17]{GJTW19}\cite{GJF19}\cite{JT20}, whose ``meaning'', however,
  has remained somewhat elusive.
  But under
  the equivalence \eqref{BoardmanHomomorphismintmfOn10fold} with Cohomotopy,
  these same elements could be understood in \cite{SS21} in their role
  in non-perturbative string theory.
\end{remark}

\begin{example}[Chern-Dold character on integral Morava K-theory]
  \label{PontrjaginCharacterInMorava}
  We highlight that a particularly interesting example of the Chern-Dold character, which is not widely known,
  is that on integral Morava K-theory,
  whose codomain in real cohomology has a rich coefficient system.
Morava K-theories $K(n)$ \cite{JohnsonWilson75}
(reviewed in \cite{Wurgler89}\cite[\S IX.7]{Rudyak98})
form a sequence of spectra labeled by
chromatic level $n \in \mathbb{N}$ and by a prime $p$
(notationally left implicit). Their coefficient ring
is pure torsion, and hence vanishes upon rationalization.
%one chromatic level $n$ at a time
%and one prime $p$ at a time.
%However, there is  a version which is integral (or $p$-adic),
However, there is an integral version $\widetilde K(n)$,
highlighted in \cite{KS1}\cite{tcu}\cite{Buh}\cite{SW}\cite{GS-AHSS},
%and which is amenable to rationalization (in the $p$-adic sense)
which has an integral $p$-adic coefficient ring:
%Indeed, $K(n)$ is the mod $p$ reduction of an integral (or $p$-adic) lift %$\widetilde{K}(n)$ with coefficient ring
\vspace{-1mm}
\begin{equation}
  \label{CoefficientRingOfIntegralMoravaKTheory}
  \widetilde{K}(n)_*
  \;=\;
  \mathbb{Z}_p[v_n, v_n^{-1}]
  \,,
  \phantom{AAAAA}
  \mbox{
    with $\mathrm{deg}(v_n) = 2(p^n-1)$.
  }
\end{equation}

\vspace{-1mm}
\noindent
This theory more closely resembles complex K-theory than is the case for
$K(n)$;
%the mod $p$ versions;
in fact, for $n=1$, it coincides with the $p$-completion of complex K-theory.
%\cite{KS1}\cite{tcu}\cite{Buh}\cite{SW}.
%The character for this theory is implicit in differential refinements %\cite{GS-AHSS}) but we make it
%more explicit here.
%
%
%Rationally, Morava K-theory  is isomorphic to $v_n$-periodic `rational' %cohomology
% \vspace{-1mm}
% $$
%\widetilde{K}(n)^*(X) \otimes \mathbb{Q} \cong H^*\big(X; %\widetilde{K}(n)_* \otimes \mathbb{Q}\big).
% $$
%
%  \vspace{-1mm}
%\noindent
%We make this precise as follows.
%Rationalizing the $p$-adic integers $\mathbb{Z}_p$ gives

Therefore, the Chern-Dold character (Def. \ref{ChernDoldCharacter})
on integral Morava K-theory \cite[p. 53]{GS-AHSS}
is of the form
\vspace{-2mm}
\begin{equation}
  \label{MoravaCharacter}
  \mathrm{ch}_{\rm Mor}
  \;:\;
  \xymatrix{
    \widetilde{K}(n)(-)
    \ar[rr]
    &&
    H^*\big(-;\, \mathbb{Q}_p[v_n, v_n^{-1}] \otimes_{\scalebox{.5}{$\mathbb{Q}$}} \mathbb{R} \big)
    \,,
  }
\end{equation}

\vspace{-2mm}
\noindent where we used \eqref{CoefficientRingOfIntegralMoravaKTheory} in
\eqref{ChernDoldCharacterAsComposite} together with the fact that
the rationalization of the $p$-adic integers is the
rational (here: real, by Remark \ref{RationalHomotopyTheoryOverTheRealNumbers})
$p$-adic numbers\footnote{Note, parenthetically, that the classical Chern character
$\mathrm{ch}$ itself can be extended to
cohomology theories with values in graded $\mathbb{Q}$-algebras;
see, e.g., \cite{Maak}.}
$
  \mathbb{Z}_p \otimes_{\scalebox{.5}{$\mathbb{Z}$}} \mathbb{R}
   \;\simeq\;
  \mathbb{Q}_p \otimes_{\scalebox{.5}{$\mathbb{Q}$}} \mathbb{R}
$.

Now $\mathbb{Q}_p$ is not finite-dimensional over $\mathbb{Q}$,
whence $\mathbb{Q}_p \otimes \mathbb{R}$ is not finite-dimensional
over $\mathbb{R}$, so that the classifying space for
integral Morava K-theory is not of $\mathbb{R}$-finite type
(Def. \ref{NilpotentConnectedSpacesOfFiniteRationalType}).
Therefore, our \emph{proof} of the non-abelian de Rham theorem (Theorem \ref{NonAbelianDeRhamTheorem}),
being based on the fundamental theorem of dgc-algebraic
rational homotopy theory (Prop. \ref{FundamentalTheoremOfdgcAlgebraicRationalHomotopyTheory}),
does not immediately apply to integral Morava K-theory
coefficients; and hence the non-abelian character
on integral Morava K-theory
with de Rham codomain,
in the form defined in Def. \ref{NonAbelianChernDoldCharacter},
is not established here. While this is
a purely technical issue, as discussed in Remark \ref{AssumptionsOnConnctedNilpotentRFiniteHomotopyTypes},
further discussion is beyond the scope of the present article.

%    Morava K-theory plays an important role in stratifying the sphere %spectrum
%    into layers by localizing at a hierarchy of chromatic levels
%    (\cite{MahowaldRavenel87} review in \cite{Ravenel86}\cite{Lurie10}).
%    At the lowest
%    level $n=0$, $K(0) \cong H\mathbb{Q}$ so that
%    rationalization can be viewed as localization $L_{K(0)}$ with respect %to  $K(0)$
%     \vspace{-2mm}
%    $$
%      \pi_*\big(L_{K(0)} E\big) = \pi_*(E) \otimes \mathbb{Q}\;.
%    $$
\end{example}

%%%%%%%%%%%%%%%%%%%%%%%%%%%%%%%%%%%%%%%%%%%%%%%%%%%%
\subsection{Chern-Weil homomorphism}
 \label{TheChernWeilHomomorphism}
%%%%%%%%%%%%%%%%%%%%%%%%%%%%%%%%%%%%%%%%%%%%%%%%%%%%

We prove (Theorem \ref{NonAbelianChernDoldSubsumesChernWeil}) that the non-abelian character subsumes the Chern-Weil homomorphism
(recalled as Prop. \ref{ChernWeilHomomorphism},
review in \cite[\S III]{Chern51}\cite[\S XII]{KobayashiNomizu63}\cite[\S 2]{ChernSimons74}\cite[\S C]{MilnorStasheff74}\cite[\S 2.1]{FiorenzaSchreiberStasheff10}) in degree-1 non-abelian cohomology.

\medskip

\noindent {\bf Chern-Weil theory.}
  For definiteness, we recall the statements
  of Chern-Weil theory that we need to prove Theorem \ref{NonAbelianChernDoldSubsumesChernWeil} below.

 \begin{remark}[Attributions in Chern-Weil theory]
  \label{ChernWeilTheoryAndItsFundamentalTheorem}

  \noindent {\bf (i)} What came to be known as the \emph{Chern-Weil homomorphism}
  (recalled as Def. \ref{ChernWeilHomomorphism} below)
  seems to be first publicly described by H. Cartan (in May 1950),
  in his prominent {\it S{\'e}minaire} \cite[\S 7]{Cartan50a},
  published as \cite{Cartan50b}.
  Later that year at the ICM (in Aug.-Sep. 1950),
  Chern discusses this construction in a talk \cite[(10)]{Chern50},
  including a brief reference to unpublished work by Weil
  (which remained unpublished until appearance in Weil's
  collected works \cite{Weil49}) for the proof that the
  construction is independent of the choice of connection
  (which is stated with an announcement of a
  proof in \cite[\S 7]{Cartan50a}).

  \noindent {\bf (ii)}   The new result of Chern's talk
  was the observation \cite[(15)]{Chern50}
  -- later called the {\it fundamental theorem} in
  \cite[\S III.6]{Chern51}, recalled as Prop. \ref{FundamentalTheoremOfChernWeilTheory} below --
  that this differential-geometric
  construction coincides with the topological construction
  of real characteristic classes (Example \ref{CharacteristicClassesOfPrincipalBundles}).
  This crucially uses the identification \cite[(11)]{Chern50}
  of the real cohomology of classifying space $B G$
  with invariant polynomials, later expanded on by Bott \cite[p. 239]{Bott73}.
  (Various subsequent authors, e.g. \cite[(1.14)]{Freed02},
  suggest to prove Chern's equation (15)
  by making sense of a connection on the universal $G$-bundle
  -- which is possible though notoriously subtle, e.g. \cite{Mostow79}
  --
  but the proof in \cite{Chern50} simply observes that
  for any fixed bound $\leq d$ on the dimension of
  the domain space,
  the classifying space
  for $G$-principal bundles may
  be truncated to a $d+1$-dimensional sub-complex
  (as follows by the cellular approximation theorem
  \cite[p. 404]{Spanier66}),
  this carrying a smooth $G$-principal bundle with ordinary connection,
  which is universal for $G$-principal bundles over $\leq d$-manifolds.
  This argument was later worked out in \cite{NarasimhanRamanan61}\cite{NarasimhanRamanan63}\cite{Schlafly80}).

  \noindent {\bf (iii)}  It is this fundamental theorem \cite[(15)]{Chern50}\cite[\S III.6]{Chern51} which
  allows to  identify the Chern-Weil homomorphism as an instance of
  the non-abelian character,
  in Theorem \ref{NonAbelianChernDoldSubsumesChernWeil} below.
\end{remark}

\begin{notation}[Principal bundles with connection]
  \label{PrincipalBundlesWithConnection}
  For $G \,\in\, \mathrm{LieGroups}$ and
  $X \,\in\, \SmoothManifolds$, we write
  \vspace{-2mm}
  \begin{equation}
    \label{ForgetfulMapFromGBundlesWithConnection}
    \xymatrix{
      G\mathrm{Connections}(X)_{\!/\sim}
      \ar@{->>}[r]
      &
      G \mathrm{Bundles}(X)_{\!/\sim}
    }
  \end{equation}

  \vspace{-2mm}
  \noindent  for the forgetful map from the set of
  isomorphism classes, over $X$,
  of $G$-principal bundles
  equipped with principal connections
  (review in \cite[\S 9]{Nakahara03}\cite[\S 1]{RudolphSchmidt17})
  to the underlying bundles without connection, .
\end{notation}
The function \eqref{ForgetfulMapFromGBundlesWithConnection}
is surjective and admits sections, corresponding to a choice
of the class of a principal connection on any class of
$G$-principal bundles.

\begin{defn}[Invariant polynomials {\cite{Weil49}\cite[\S 7]{Cartan50a}}]
  \label{InvariantPolynomials}
  For $\mathfrak{g} \in \LieAlgebras$, we write
  \vspace{-2mm}
  $$
    \mathrm{inv}^\bullet(\mathfrak{g})
    \;:=\;
    \mathrm{Sym}
    \big(
      \mathfrak{b}^2 \mathfrak{g}^\ast
    \big)^G
    \;\;
    \in
    \;
    \GradedAlgebras
  $$

  \vspace{-2mm}
  \noindent
  for the graded sub-algebra \eqref{CategoryOfGradedAlgebras}
  on those elements in the symmetric algebra \eqref{SymmetricGradedAlgebra}
  of the linear dual of $\mathfrak{g}$ shifted up (Def. \ref{DegreeShift})
  into degree 2, which are invariant under the adjoint action of
  $G$ on $\mathfrak{g}^\ast$.
\end{defn}

\begin{defn}[Characteristic forms {\cite[\S 7]{Cartan50a}\cite[(10)]{Chern50}}]
  \label{CharacteristicForms}
 Let $G$ be a finite-dimensional Lie group with Lie algebra
 $\mathfrak{g}$, and let
 $\!\!\xymatrix@C=11pt{P \ar[r]^p &  X}\!\!$ be $G$-principal bundle
 with connection $\nabla$
 (Def. \ref{PrincipalBundlesWithConnection}).
  Then for $\omega \in \mathrm{inv}^{2n}(\mathfrak{g})$ an
 invariant polynomial (Def. \ref{InvariantPolynomials}),
 its evaluation on the curvature 2-form
 $F_\nabla \in \Omega^2(P)\otimes \mathfrak{g}$
 of the connection
 yields a differential form
  \vspace{-2mm}
 $$
   \omega(F_\nabla)
   \;\in\;
   \xymatrix{
     \Omega_{\mathrm{dR}}^{2n}(X)
     \ar[r]^-{ p^\ast }
     &
     \Omega_{\mathrm{dR}}^{2n}(P)
   }
 $$

 \vspace{-1mm}
 \noindent
 which, by the second condition on an Ehresmann connection,
 is \emph{basic}, namely in the image of the
 pullback operation along the bundle projection $p$, as shown.
 Regarded as a differential form on $X$, this
 is called the \emph{characteristic form}
 corresponding to $\omega$.
\end{defn}

\begin{lemma}[Characteristic de Rham classes of characteristic forms {\cite{Weil49}\cite[p. 401]{Chern50}\cite[\S III.4]{Chern51}}]
 \label{CharacteristicClassesOfCharacteristicForms}
 The class in de Rham cohomology
  \vspace{-2mm}
  $$
   \big[
     \omega(F_\nabla)
   \big]
   \;\;
   \in
   \;
   H^{2n}_{\mathrm{dR}}
   (
     X
   )
 $$

 \vspace{-1mm}
 \noindent
 of a characteristic form
 in Def. \ref{CharacteristicForms}
 is independent of the choice of connection $\nabla$
 and depends only  on the isomorphism class of
 the principal bundle $P$.
\end{lemma}

\begin{defn}[Chern-Weil homomorphism {\cite[\S 7]{Cartan50a}\cite[(10)]{Chern50}}]
  \label{ChernWeilHomomorphism}
  Let $G$ be a finite-dimensional Lie group, with classifying space denoted
  $B G$. The \emph{Chern-Weil homomorphism} is the
  composite map
   \vspace{-2mm}
\noindent
  \begin{equation}
    \label{ClassicalChernWeilConstruction}
    \xymatrix@R=6pt{
              \mathclap{
        \raisebox{0pt}{
          \tiny
          \color{darkblue}
          \bf
          \def\arraystretch{.9}
          \begin{tabular}{c}
            Chern-Weil
            \\
            homomorphism
          \end{tabular}
        }
        }
     \phantom{AAAA}
      {
        \mathrm{cw}_G
      }
      \ar@{}[r]|-{:}
      &
      G\mathrm{Bundles}(X)_{\!\!/\sim}
      \ar[r]
      &
      G\mathrm{Connections}(X)_{\!\!/\sim}
      \ar[rr]
      &&
      \mathrm{Hom}
      \big(
        \mathrm{inv}^\bullet(\mathfrak{g})
        ,\,
        H_{\mathrm{dR}}^\bullet
        (
          X
        )
      \big)
      \\
      &
      \overset{
        \mathclap{
        \raisebox{3pt}{
          \tiny
          \color{darkblue}
          \bf
          principal bundle
        }
        }
      }{
        [P]
      }
      \ar@{|->}[r]
      &
      \overset{
        \mathclap{
        \raisebox{3pt}{
          \tiny
          \color{darkblue}
          \bf
          with connection
        }
        }
      }{
        [P,\nabla]
      }
      \ar@{|->}[rr]
      &&
           \big(\phantom{AA}
        \overset{
          \mathclap{
          \raisebox{8pt}{
            \tiny
            \color{darkblue}
            \bf
            \def\arraystretch{.9}
            \begin{tabular}{c}
              invariant
              \\
              polynomial
            \end{tabular}
          }
          }
        }{
          \omega
        }
        \;\;\;\mapsto\;\;\;
        \overset{
          \mathclap{
          \raisebox{4pt}{
            \tiny
            \color{darkblue}
            \bf
            \def\arraystretch{.9}
            \begin{tabular}{c}
              de Rham class of
              \\
              characteristic
              form
            \end{tabular}
          }
          }
        }{
          \big[\omega(F_\nabla)\big]
        }
      \phantom{A}\big)
      \,,
    }
  \end{equation}

\vspace{-2mm}
\noindent  where the first map is any section of \eqref{ForgetfulMapFromGBundlesWithConnection},
  given by choosing any connection on a given principal bundle;
  and the second map is the construction of characteristic forms
  according to Def. \ref{CharacteristicForms}.
  (The $\mathrm{Hom}$ on the right is that in $\GradedAlgebras$.)
  By Lemma \ref{CharacteristicClassesOfCharacteristicForms}
  the second map is well-defined (and its composition
  with the first turns out to be independent of the choices made,
  by Prop. \ref{FundamentalTheoremOfChernWeilTheory} below).
\end{defn}
That this construction is useful, in that it produces interesting
real characteristic classes of $G$-principal bundles
(Example \ref{CharacteristicClassesOfPrincipalBundles}), is the following statement:

\begin{prop}[Abstract Chern-Weil homomorphism {\cite[(11)]{Chern50}\cite[\S III.5]{Chern51}\cite[p. 239]{Bott73}}]
  \label{AbstractChernWeilHomomorphism}
  Let $G$ be a finite-dimensional, compact
  Lie group, with Lie algebra denoted $\mathfrak{g}$.
  Then the real cohomology algebra of its classifying space
  $B G$ is isomorphic to the algebra of invariant polynomials
  (Def. \ref{InvariantPolynomials}):
   \vspace{-1mm}
  \begin{equation}
    \label{RealCohomologyOfClassifyingSpaceIsInvariantPolynomials}
    \mathrm{inv}^\bullet(\mathfrak{g})
    \;\simeq\;
    H^\bullet
    (
      B G;
      \,
      \mathbb{R}
    )
    \;\;\;
    \in
    \;
    \GradedAlgebras
    \,.
  \end{equation}
\end{prop}

We can also obtain the following:

\begin{prop}[Fundamental theorem of Chern-Weil theory {\cite[(15)]{Chern50}\cite[\S III.6]{Chern51}} (Rem. \ref{ChernWeilTheoryAndItsFundamentalTheorem})]
\label{FundamentalTheoremOfChernWeilTheory}
Let $G$ be a finite-dimensional compact Lie group.
Then the Chern-Weil homomorphism (Def. \ref{ChernWeilHomomorphism})
coincides with the operation of pullback of universal
characteristic classes along the classifying maps of $G$-bundles
(Example \ref{CharacteristicClassesOfPrincipalBundles}),
in that the following diagram commutes:
\vspace{-2mm}
\begin{equation}
  \label{NonabelianChernDoldRelatedToTraditionalChernWeil}
  \raisebox{20pt}{
  \xymatrix@R=20pt@C=5em{
    H
    (
      X;
      \,
      B G
    )
    \ar[rr]_-{
      c \;\mapsto\; c^\ast(-)
    }^-{
        \mathclap{
        \raisebox{0pt}{
          \tiny
          \color{darkblue}
          \bf
          \def\arraystretch{.9}
          \begin{tabular}{c}
            pullback of
            \\
            universal characteristic classes
            \\
            along classifying map \eqref{PullbackOfUniversalCharacteristicClasses}
          \end{tabular}
        }
        }
      }
    \ar@{<-}[d]_-{
          \mbox{
          \tiny
          \color{darkblue}
          \bf
          \eqref{IsomorphismBetweenPrincipalBundlesAndMapsToBG}
      }}^-{\simeq}
    &&
    \mathrm{Hom}
    \big(
      H^\bullet
      (
        B G;
        \,
        \mathbb{R}
      )
      \,
      ,
      \,
      H^\bullet
      (
        X;
        \,
        \mathbb{R}
      )
    \big)
    \ar[d]_-{
     \simeq}^-{
        \mbox{
          \tiny
          \color{darkblue}
          \bf
          \eqref{RealCohomologyOfClassifyingSpaceIsInvariantPolynomials}
        }
      }
     \\
    G \mathrm{Bundles}(X)_{\!/\sim}
    \ar[rr]^-{\mathrm{cw}_G}_-{
        \mathclap{
        \raisebox{-3pt}{
          \tiny
          \color{darkblue}
          \bf
          Chern-Weil homomorphism \eqref{ClassicalChernWeilConstruction}
        }
        }
        }
    &&
      \mathrm{Hom}
      \big(
        \mathrm{inv}^\bullet(\mathfrak{g})
        ,\,
        H_{\mathrm{dR}}^\bullet
        (
          X
        )
      \big)
    }
  }
\end{equation}

\vspace{-1mm}
\noindent
Here the isomorphism on the left is from Example \ref{TraditionalNonAbelianCohomology},
while that from the right is from Prop. \ref{AbstractChernWeilHomomorphism}
and using the de Rham theorem.
\end{prop}

\begin{example}[Characteristic forms of classical Lie groups (e.g. {\cite[Ex. 11.5-7]{Nakahara03}})]
\label{PontrjaginForms}
Let $G = \mathrm{SU}(n)$ be the special unitary group, for $n \in \mathbb{N}$.
Then the fundamental Chern-Weil theorem Prop. \ref{FundamentalTheoremOfChernWeilTheory}
identifies,
for any connection $\nabla$ \eqref{ClassicalChernWeilConstruction},
on a given $\mathrm{SU}(n)$-principal bundle, with
associated curvature differential form
\vspace{-2mm}
\begin{equation}
  \label{CurvatureDifferentialForm}
  \begin{tikzcd}[row sep=-2pt, column sep=0pt]
    \Omega^2\big(X;\, \mathfrak{u}(n)\big)
      \ar[
        rr,
        hook
      ]
      &&
    \Omega^2\big(X;\, \mathrm{Mat}_{n \times n}(\mathbb{C})\big)
    \\
    F_\nabla
    &\longmapsto&
    \big( (F_\nabla)_{a b}\big)_{ \mathrlap{1 \leq a,b \leq n} }
  \end{tikzcd}
\end{equation}

\vspace{-2.5mm}
\noindent the following (de Rham cohomology classes of)
classical characteristic forms (Def. \ref{CharacteristicForms}):

\noindent
{\bf (i) Chern forms.}
The real-cohomology images of the first couple {\it Chern classes}
$c_i \,\in\, H^{2i}\big( B \mathrm{U}(n);\, \mathbb{Z} \big)$
are identified with the de Rham cohomology classes
$
  \mathrm{cw}_{{}_{\mathrm{SU}(n)}}(c_i)
  \;=\;
  \big[ c_i(\nabla) \big]
  \;\;
  \in
  \;
  H^{2i}_{\mathrm{dR}}\big(X\big)
$
of the polynomials in the curvature differential form
\eqref{CurvatureDifferentialForm}
which are the homogeneous components of the following
{\it total Chern form} \footnote{The
standard normalization factor $i / 2\pi$ appearing here
results from identifying $\mathrm{U}(1)$ with $\mathbb{R}/h\mathbb{Z}$
for the choice $h = 2\pi$.}

\vspace{-2mm}
$$
  c(\nabla)
  \;:=\;
  \underset{k \in \mathbb{N}}{\sum}
  \;
  \underset{
    \mathrm{deg} = 2k
  }{
    \underbrace{
      c_k(\nabla)
    }
  }
  \;:=\;
  \mathrm{det}
  \left(
    1 + \tfrac{i}{2 \pi } F_\nabla
  \right)
  \,.
$$
\vspace{-2mm}

\noindent
{\bf (ii) Pontrjagin forms.}
When the structure group $G$ is reduced along the canonical inclusion
$\mathrm{SO}(n) \xhookrightarrow{\iota} \mathrm{SU}(n)$ of the
special orthogonal group, then
the real images of the first couple {\it Pontrjagin classes}
$p_k \,\in\, H^{4k}\big( B \mathrm{SO}(n);\, \mathbb{Z} \big)$
are identified with the de Rham cohomology classes
of the corresponding Chern forms \eqref{FormulaForChernForms},
up to a signs:
$$
  \mathrm{cw}_{\mathrm{SO}(n)}(p_k)
  \;=\;
  \big[ p_k(\nabla) \big]
  \;=\;
  (-1)^k
  \big[ c_{2k}(\iota_\ast \nabla) \big]
  \;\;
  \in
  \;
  H^{4k}_{\mathrm{dR}}\big(X\big)
  \,.
$$
One finds
\begin{equation}
\label{FormulaForChernForms}
\begin{aligned}
  p_1(\nabla)
  &
  \;:=\;
  - \tfrac{1}{8 \pi^2} \; \mathrm{tr}\big( F_\nabla \wedge F_\nabla \big)
  \,,
  \\
  p_2(\nabla)
  &
  \;=\;
  \tfrac{1}{128 \pi^4}
  \Big(
    \mathrm{tr}
    \left(
      F_\nabla \wedge F_\nabla
    \right)
    \wedge
    \mathrm{tr}
    \left(
      F_\nabla \wedge F_\nabla
    \right)
    -
    2 \cdot \mathrm{tr}
    \left(
       F_\nabla
         \wedge
       F_\nabla
         \wedge
       F_\nabla
         \wedge
       F_\nabla
    \right)
  \Big)
  \,.
\end{aligned}
\end{equation}
The following rational combination
of these forms plays a central role
in \cref{CohomotopicalChernCharacter}:

\vspace{-2mm}
\begin{equation}
  \label{TheI8Polynomial}
  I_8(\nabla)
  \;\coloneqq\;
  \tfrac{1}{48}
  \left(
    p_2(\nabla)
    -
    \tfrac{1}{4}
    p_1(\nabla)
    \wedge
    p_1(\nabla)
  \right)
  \,.
\end{equation}
\vspace{-4mm}

\noindent
{\bf (iii) Euler form.}
If, moreover, $n = 2k$ is even, then
the real image of the {\it Euler class}
$\rchi_n \,\in\, H^{n}\big( B \mathrm{SO}(n);\, \mathbb{Z} \big)$
is identified with the de Rham cohomology class
$$
  \mathrm{cw}_{\mathrm{SO}(n)}(\rchi_n)
  \;=\;
  \big[ \rchi_n(\nabla) \big]
  \;\;
  \in
  \;
  H^{n}_{\mathrm{dR}}\big(X\big)
$$
of the {\it Pfaffian} wedge-product polynomial of the matrix \eqref{CurvatureDifferentialForm}:
\vspace{-2mm}
\begin{equation}
  \label{FormulaForEulerForms}
  \rchi_{2k}(\nabla)
  \;\coloneqq\;
  \tfrac
    {(-1)^{n/2}}
    {(4\pi)^{n/2} \cdot (n/2)!}
  \underset
    {\sigma \in \mathrm{Sym}(n)}
    {\sum}
  \mathrm{sgn}(\sigma)
  \cdot
  (F_\nabla)_{\sigma(1)\sigma(2)}
  \wedge
  (F_\nabla)_{\sigma(3)\sigma(4)}
  \wedge
  \cdots
  \wedge
  (F_\nabla)_{\sigma(n-1)\sigma(n)}
  \,.
\end{equation}
\end{example}

\medskip

\noindent {\bf Chern-Weil homomorphism as a special case of the non-abelian character.}

\begin{lemma}[Sullivan model of classifying space]
  \label{SullivanModelOfClassifyingSpace}
  Let $G$ be a finite-dimensional, compact and simply-connected
  Lie group, with Lie algebra denoted $\mathfrak{g}$.
  Then the minimal Suillvan model (Def. \ref{MinimalSullivanModels})
  of its classifying space $B G$ is the graded algebra of
  invariant polynomials (Def. \ref{InvariantPolynomials}),
  regarded as a dgc-algebra with vanishing differential:
     \vspace{-2mm}
  \begin{equation}
    \label{SullivanModelOfClassifyingSpaceIsInvariantPolynomials}
    \big(
      \mathrm{inv}(\mathfrak{g}),
      \,
      d = 0
    \big)
    \;\simeq\;
    \mathrm{CE}
    (
      \mathfrak{l} B G
    )
    \;\;\;\;
    \in
    \;
    \dgcAlgebras{\mathbb{R}} \;.
  \end{equation}
\end{lemma}

\begin{proof}
  According to \cite[Example 2.42]{FOT08}, we have
   \vspace{-2mm}
  \begin{equation}
    \label{SullivanModelOfClassifyingSpaceIsCohomologyRing}
    \mathrm{CE}
    (
      \mathfrak{l}BG
    )
    \;\simeq\;
    \big(
      H^\bullet
      (
        B G;
        \,
        \mathbb{R}
      ),
      \,
      d = 0
    \big).
  \end{equation}

   \vspace{-2mm}
\noindent
  The composition of
  \eqref{SullivanModelOfClassifyingSpaceIsCohomologyRing}
  with the isomorphism
  \eqref{RealCohomologyOfClassifyingSpaceIsInvariantPolynomials}
  from Prop. \ref{AbstractChernWeilHomomorphism}
  yields the desired \eqref{SullivanModelOfClassifyingSpaceIsInvariantPolynomials}.
\end{proof}

\begin{lemma}[Non-abelian de Rham cohomology with coefficients in a classifying space]
  \label{NonAbelianDeRhamCohomologyWithCoefficientsInAClassifyingSpace}
  Let $G$ be a finite-dimensional, compact and simply-connected
  Lie group, with Lie algebra denoted $\mathfrak{g}$.
  Then the non-abelian de Rham cohomology (Def. \ref{NonabelianDeRhamCohomology})
  with coefficients in the rational Whitehead $L_\infty$-algebra
  $\mathfrak{l} B G$
  (Prop. \ref{WhiteheadLInfinityAlgebras})
  of the classifying space is canonically identified with the
  codomain of the classical Chern-Weil construction
  \eqref{ClassicalChernWeilConstruction}:
   \vspace{-2mm}
  \begin{equation}
    \label{NonabelianDeRhamRelatedToTraditionalChernWeilCodomain}
    \overset{
      \mathclap{
      \raisebox{3pt}{
        \tiny
        \color{darkblue}
        \bf
        \def\arraystretch{.9}
        \begin{tabular}{c}
          nonabelian
          \\
          de Rham cohomology
        \end{tabular}
      }
      }
    }{
      H_{\mathrm{dR}}
      \big(
        X;
        \,
        \mathfrak{l}BG
      \big)
    }
    \;\;
      \simeq
    \;\;
    \overset{
      \mathclap{
      \raisebox{3pt}{
        \tiny
        \color{darkblue}
        \bf
        \def\arraystretch{.9}
        \begin{tabular}{c}
          traditional codomain of
          \\
          Chern-Weil construction
        \end{tabular}
      }
      }
    }{
      \mathrm{Hom}
      \big(
        \mathrm{inv}^\bullet(\mathfrak{g})
        ,\,
        H_{\mathrm{dR}}^\bullet
        (
          X
        )
      \big).
    }
  \end{equation}
\end{lemma}

\newpage

\begin{proof}
  Consider the following sequence of natural bijections:
  \vspace{-1mm}
  $$
    \begin{aligned}
      H_{\mathrm{dR}}
      \big(
        X;
        \,
        \mathfrak{l}BG
      \big)
      &
      :=
      \dgcAlgebras{\mathbb{R}}
      \big(
        \mathrm{CE}
        \big(
          \mathfrak{l} BG
        \big)
        \,,\,
        \Omega^\bullet_{\mathrm{dR}}(X)
      \big)_{\!\!/\sim}
      \\
      & \simeq
      \dgcAlgebras{\mathbb{R}}
      \Big(
        \big(
          \mathrm{inv}^\bullet(\mathfrak{g}),
          \,
          d = 0
        \big)
        \,,\,
        \Omega^\bullet_{\mathrm{dR}}(X)
      \Big)_{\!\!/\sim}
      \\
      & \simeq
      \GradedAlgebras
      \Big(
        \mathrm{inv}^\bullet(\mathfrak{g})
        \,,\,
        \Omega^\bullet_{\mathrm{dR}}(X)_{\mathrm{closed}}
      \Big)_{\!\!/\sim}
      \\
      & \simeq
      \GradedAlgebras
      \Big(
        \mathrm{inv}^\bullet(\mathfrak{g})
        \,,\,
        \big(
          \Omega^\bullet_{\mathrm{dR}}(X)_{\mathrm{closed}}
        \big)_{\!\!/\sim}
      \Big)
      \\
      & \simeq
      \GradedAlgebras
      \big(
        \mathrm{inv}^\bullet(\mathfrak{g})
        \,,\,
        H^\bullet_{\mathrm{dR}}(X)
      \big)
      \\
      & =:
      \mathrm{Hom}
      \big(
        \mathrm{inv}^\bullet(\mathfrak{g})
        \,,\,
        H^\bullet_{\mathrm{dR}}(X)
      \big).
    \end{aligned}
  $$

  \vspace{-1mm}
  \noindent
  Here the first line is the definition (Def. \ref{NonabelianDeRhamCohomology}).
  After that, the first step is Lemma \ref{NonAbelianDeRhamCohomologyWithCoefficientsInAClassifyingSpace}.
  The second step unwinds what it means to hom out of
  a dgc-algebra with vanishing differential
  (which is generator-wise as in Example \ref{OrdinaryClosedFormsAreFlatLineLInfinityAlgebraValuedForms}),
  while the third and fourth steps unwind what this means for
  the coboundary relations
  (which is generator-wise as in Prop. \ref{NonAbelianDeRhamCohomologySubsumesOrdinaryDeRhamCohomology}).
  The last line just matches the result to the abbreviated notation
  used in \eqref{ClassicalChernWeilConstruction}.
\end{proof}

\begin{theorem}[Non-abelian character map subsumes Chern-Weil homomorphism]
\label{NonAbelianChernDoldSubsumesChernWeil}
Let $G$ be a finite-dimensional compact, connected and simply-connected Lie group,
with Lie algebra $\mathfrak{g}$. Let $X \in \NilpotentConnectedQFiniteHomotopyTypes$
(Def. \ref{NilpotentConnectedSpacesOfFiniteRationalType}) be equipped
with the structure of a smooth manifold. Then the non-abelian character
$\mathrm{ch}_{B G}$ (Def \ref{NonAbelianChernDoldCharacter})
on non-abelian cohomology (Def. \ref{NonAbelianCohomology}) of $X$
with coefficients in
%the classifying space
$B G$ coincides with the Chern-Weil homomorphism
$\mathrm{cw}_G$ (Def. \ref{ChernWeilHomomorphism}) with coefficients in $G$,
in that the following diagram (of cohomology sets) commutes:
 \vspace{-2mm}
\begin{equation}
  \label{NonabelianChernDoldRelatedToTraditionalChernWeil}
  \raisebox{20pt}{
  \xymatrix@R=20pt@C=4em{
    H
    (
      X;
      \,
      B G
    )
    \ar[rr]_-{ \mathrm{ch}_{B G} }^-{
        \mathclap{
        \raisebox{0pt}{
          \tiny
          \color{darkblue}
          \bf
          \def\arraystretch{.9}
          \begin{tabular}{c}
            non-abelian character
          \end{tabular}
        }
        }
      }
    \ar@{<-}[d]_-{
          \mbox{
          \tiny
          \color{darkblue}
          \bf
          \eqref{IsomorphismBetweenPrincipalBundlesAndMapsToBG}
                   }}^-{\simeq}
    &&
    H_{\mathrm{dR}}
    (
      X;
      \,
      \mathfrak{l}B G
    )
    \ar[d]_-{
     \simeq}^-{
                \mbox{
          \tiny
          \color{darkblue}
          \bf
          \eqref{NonabelianDeRhamRelatedToTraditionalChernWeilCodomain}
        }
      }
       \\
    G \mathrm{Bundles}(X)_{\!/\sim}
    \ar[rr]^-{\mathrm{cw}_G}_-{
        \mathclap{
        \raisebox{-3pt}{
          \tiny
          \color{darkblue}
          \bf
          Chern-Weil homomorphism
        }
        }
        }
    &&
      \mathrm{Hom}
      \big(
        \mathrm{inv}^\bullet(\mathfrak{g})
        ,\,
        H_{\mathrm{dR}}^\bullet
        (
          X
        )
      \big)
    }
  }
\end{equation}

\vspace{-1mm}
\noindent
Here the isomorphism on the left is from Example \ref{TraditionalNonAbelianCohomology},
while that
on the right is from Lemma \ref{NonAbelianDeRhamCohomologyWithCoefficientsInAClassifyingSpace}.
\end{theorem}
\begin{proof}
  First, notice that $B G$ is simply connected
  (hence nilpotent), by the assumption
  that $G$ is connected, and that it is of finite rational type
  by Prop. \ref{AbstractChernWeilHomomorphism}. Hence, with Def. \ref{NilpotentConnectedSpacesOfFiniteRationalType},
   \vspace{-2mm}
  \begin{equation}
    \label{ClassifyinSpaceInConnectedNilpotentRFiniteHomotopyTypes}
    B G
    \;\in\;
    \NilpotentConnectedQFiniteHomotopyTypes
    \,.
  \end{equation}

     \vspace{-1mm}
\noindent
  Now, by Definition \ref{NonAbelianChernDoldCharacter},
  the non-abelian character
  map on the top of \eqref{NonabelianChernDoldRelatedToTraditionalChernWeil}
  \vspace{-3mm}
  $$
    \mathrm{ch}_{BG}
    \;:\;
    \xymatrix@C=3em{
      H
      (
        X;
        \,
        B G
      )
      \ar[r]^-{ (\eta^{\mathbb{R}}_{B G}) }
      &
      H
      \big(
        X;
        \,
        L_{\mathbb{R}} B G
      \big)
      \ar[r]^-{ \simeq }
      &
      H_{\mathrm{dR}}
      \big(
        X;
        \,
        L_{\mathbb{R}} B G
      \big)
    }
  $$

  \vspace{-3mm}
  \noindent
  sends a classifying map
     \vspace{-3mm}
  $$
    \xymatrix@C=15pt{X \ar[r]^-{ c } & B G}
    \;\in\;
    H(X;\, B G)
    \;=\;
    \HomotopyTypes
    (
      X
      \,,\,
      B G
    )
  $$

  \vspace{-1mm}
  \noindent
  first to its composite with the rationalization
  map (Def. \ref{Rationalization}). By the
  fundamental theorem  (Theorem \ref{FundamentalTheoremOfdgcAlgebraicRationalHomotopyTheory} {(i)},
  using \eqref{ClassifyinSpaceInConnectedNilpotentRFiniteHomotopyTypes}),
  this is given by the derived adjunction unit $\mathbb{D}\eta_{B G}$
  of
  $\RightDerived\Bexp_{\mathrm{PL}} \,\dashv\, \Omega^\bullet_{\mathrm{PLdR}}$
  \eqref{PLdeRhamDerivedAdjunction}:
     \vspace{-2mm}
  $$
    \xymatrix{
      X
      \ar[r]^-{ c }
      &
      B G
      \ar[rr]^-{
        \LeftDerived_{\mathbb{R}} B G
        \;\simeq\;
        \mathbb{D}\eta_{B G}
      }
      &&
      \RightDerived\Bexp_{\mathrm{PL}}
        \,\circ\,
      \Omega^\bullet_{\mathrm{PLdR}}
      (
        B G
      )
    }
    \;\;\;\;\;
    \in
    \;
    \HomotopyTypes
    \big(
      X
      \,,\,
      L_{\mathbb{R}}B G
    \big)
    \;=\;
    H
    \big(
      X;
      \,
      L_{\mathbb{R}}BG
    \big)
    \,.
  $$

  \vspace{-1mm}
  \noindent
  Moreover, by part (ii) of the fundamental theorem,
  the adjunct of the morphism
   $  \mathbb{D}\eta_{B G} \,\circ\, c $
  under \eqref{PLdeRhamDerivedAdjunction} is

  \vspace{-3mm}
  $$
    \xymatrix{
      \Omega^\bullet_{\mathrm{PLdR}}(X)
      \ar@{<-}[r]^-{\;\;\; c^\ast }
      &
      \Omega^\bullet_{\mathrm{PLdR}}
      (
        B G
      )
      \;\;\;\;\;
    \in
    \;
    \mathrm{Ho}
    \big(
      \dgcAlgebrasProj{\mathbb{R}}
    \big)
    }
  $$

  \vspace{-2mm}
  \noindent
  (using that $\Omega^\bullet_{\mathrm{PLdR}}(\mathbb{D}\eta^{\mathbb{R}})$
  is an equivalence, by reflectivity of rationalization \eqref{RationalizationReflection}).
  Hence it is the pullback operation
  of rational cocycles on $B G$
  along the classifying map $c$.
  Sending this further along the isomorphism  to the
  bottom right in \eqref{NonabelianChernDoldRelatedToTraditionalChernWeil}
  (via Theorem \ref{NonAbelianDeRhamTheorem} and Lemma \ref{NonAbelianDeRhamCohomologyWithCoefficientsInAClassifyingSpace})
  gives, by \eqref{TowardsTheNonAbelianDeRhamTheorem}:
     \vspace{-2mm}
  \begin{equation}
    \label{TowardsSeeingThatNonabelianChernDoldIsChernWeil}
    \mathrm{ch}_{BG}
    \;:\;
    c
    \;\;\mapsto\;\;\;\;
    \xymatrix{
      \Omega^\bullet_{\mathrm{dR}}(X)
      \ar@{<-}[r]^-{\;\; c^\ast }
      &
      \Omega^\bullet_{\mathrm{PLdR}}
      (
        B G
      )
      \ar@{<-}[r]^-{ \simeq }
      &
      \mathrm{inv}^\bullet(\mathfrak{g})
    }
    \;\;\;\;\;
    \in
    \;
    \mathrm{Ho}
    \big(
      \dgcAlgebrasProj{\mathbb{R}}
    \big)
    \,.
  \end{equation}

  \vspace{-1mm}
  \noindent
  In conclusion, we have found that the commutativity of \eqref{NonabelianChernDoldRelatedToTraditionalChernWeil}
  is equivalent to the statement that the characteristic forms
  obtained by the Chern-Weil construction \eqref{ClassicalChernWeilConstruction}
  represent the pullback \eqref{TowardsSeeingThatNonabelianChernDoldIsChernWeil}
  of the universal real characteristic classes
  on $B G$ along the classifying map $c$ of the underlying
  principal bundle (Example \ref{CharacteristicClassesOfPrincipalBundles}). This is the case by the fundamental theorem
  of Chern-Weil theory, Prop. \ref{FundamentalTheoremOfChernWeilTheory}.
\end{proof}

\begin{example}[de Rham representative of tangential $\mathrm{Sp}(2)$-twist]
  \label{DeRhamRepresentativeOfTangentialSp2Twist}
  For $X$ a smooth 8-dimensional spin manifold
  equipped with tangential $\mathrm{Sp}(2)$-structure
  $\tau$ \eqref{TangentialSp2Structure},
  Thm. \ref{NonAbelianChernDoldSubsumesChernWeil}
  says that there exists a smooth $\mathrm{Sp}(2)$-principal bundle
  on $X$ equipped with an Ehresmann connection $\nabla$
  such that the $\mathbb{R}$-rationalization (Def. \ref{Lk})
  of the twist $\tau$ corresponds,
  under the non-abelian de Rham theorem
  (Theorem \ref{NonAbelianDeRhamTheorem}) to a
  flat $\mathfrak{l} B \mathrm{Sp}(2)$-valued differential form
  whose components are the characteristic forms of the
  $\mathrm{Sp}(2)$-principal connection $\nabla$:

  \vspace{-3mm}
  $$
  \hspace{-2.5cm}
    \begin{array}{ccccc}
    H\big(
      X;
      \,
      B \mathrm{Sp}(2)
    \big)
    &
    \mathclap{
    \xrightarrow{\;
      (\eta^{\mathbb{R}}_{B G})_\ast\;}
    }
    &
    H\big(
      X;
      \,
      L_{\mathbb{R}}
      B \mathrm{Sp}(2)
    \big)
    &\simeq&
    \!\!\!\!\!\!\!\!\!\!\!
    \!\!\!\!\!\!\!\!\!\!\!
    \!\!\!\!\!\!\!\!\!\!\!
    \!\!\!\!\!\!\!\!\!\!\!
    \!\!\!\!\!\!\!\!\!\!\!
    \!\!\!\!\!\!\!\!\!\!\!
    H_{\mathrm{dR}}
    \big(
      X;
      \,
      \mathfrak{l}
      B \mathrm{Sp}(2)
    \big)
    \\
    \tau
    &\longmapsto&
    L_{\mathbb{R}}\tau
    &\longleftrightarrow&
    {\xymatrix@R=-2pt@C=1.8em{
      \Omega^\bullet_{\mathrm{dR}}(X)
      \ar@{<-}[rr]^-{ \tau_{\mathrm{dR}} }
      &&
      \mathbb{R}
   \scalebox{0.8}{$   \left[
      \!\!\!
      {\begin{array}{l}
        {\phantom{\tfrac{1}{2}}}\rchi_8,
        \\[-2pt]
        \tfrac{1}{2}p_1
      \end{array}}
      \!\!\!
      \right]
      $}
      \!\big/\!
      \scalebox{0.8}{$
      \left(
      {\begin{aligned}
        d\, \tfrac{1}{2}p_1 & = 0
        \\[-2pt]
        d\, {\phantom{\tfrac{1}{2}}}\rchi_8 & = 0
      \end{aligned}}
      \right)
      $}
      \mathrlap{
        \;
        =
        \;
        \mathrm{CE}
        \big(
          \mathfrak{l} B \mathrm{Sp}(2)
        \big)
      }
      \\
      \tfrac{1}{2}p_1(\nabla)
      \ar@{<-|}[rr]
      &&
      \tfrac{1}{2}p_1
      \\
      {\phantom{\tfrac{1}{2}}}\rchi_8(\nabla)
      \ar@{<-|}[rr]
      &&
      {\phantom{\tfrac{1}{2}}}\rchi_8
    }}
    \end{array}
  $$

  \vspace{-2mm}
  \noindent
  Here on the right we are using \cite[Thm . 8.1]{CV98},
  see \cite[Lemma 2.12]{FSS20a} to identify generating universal
  characteristic classes on $B \mathrm{Sp}(2)$:
  $\tfrac{1}{2}p_1$ is the first Pontrjagin class
  (see Ex. \ref{PontrjaginForms})
  and
  $
    \rchi_8
    =
    \left(
      \tfrac{1}{2}
      p_2  - \big(\tfrac{1}{2}p_1\big)^2
    \right)
  $ is the Euler 8-class, which here on $B \mathrm{Sp}(2)$
  happens to be proportional to $I_8$ \eqref{TheI8Polynomial},
  see \cite[Prop. 3.7]{FSS19b}.
\end{example}

%%%%%%%%%%%%%%%%%%%%%%%%%%%%%%%%%%%%%%%%%%%%%%%%%%%%%%%%%
\subsection{Cheeger-Simons homomorphism}
 \label{NonabelianDifferentialCohomology}
%%%%%%%%%%%%%%%%%%%%%%%%%%%%%%%%%%%%%%%%%%%%%%%%%%%%%%%%%

We show (Theorem \ref{SecondaryDifferentialNonAbelianCharacterSubsumesCheegerSimonsHomomorphism}) that the non-abelian character map
induces \emph{secondary} non-abelian cohomology operations
(Def. \ref{SecondaryNonabelianCohomologyOperations})
which subsume the Cheeger-Simons homomorphism,
recalled around \eqref{CheegerSimonsHomomorphism} below,
with values in ordinary differential cohomology, recalled
around \eqref{DifferentialCohomologyDiagram} below.
We follow \cite{FiorenzaSchreiberStasheff10} \cite{SSS09}\cite{dcct}
where the Cheeger-Simons homomorphism, generalized to
higher principal bundles, is called the
{\it $\infty$-Chern-Weil homomorphism}.
Underlying this is a differential enhancement of the
non-abelian character map (Def. \ref{DifferentialNonabelianCharacterMap}),
and an induced notion of
differential non-abelian cohomology (Def. \ref{DifferentialNonAbelianCohomology})
on smooth $\infty$-stacks (recalled as Def. \ref{SmoothInfinityStacks}).

\medskip

\noindent {\bf The differential non-abelian character map.}
We introduce (in Def. \ref{DifferentialNonabelianCharacterMap} below)
the differential refinement of the non-abelian character map;
given as before by rationalization, but now followed not by a map to non-abelian
de Rham cohomology, but to its refinement by the full cocycle space of
flat non-abelian differential forms (Def. \ref{SmoothCocycleSpaceOfNonAbelianDeRhamCoefficients} below). It
is this refinement of the codomain of the character map that allows
it to be fibered over the smooth space (Ex. \ref{SmoothSpaces})
of actual flat non-abelian differential forms (instead of just their non-abelian de Rham classes),
thus producing differential non-abelian cohomology
(Def. \ref{DifferentialNonAbelianCohomology} below).

  \begin{defn}[Moduli $\infty$-stack of flat $L_\infty$-algebra valued forms {\cite[4.4.14.2]{dcct}}]
    \label{SmoothCocycleSpaceOfNonAbelianDeRhamCoefficients}
    Let $A \in \SimplicialSets$ be of
    connected, nilpotent, $\mathbb{R}$-finite homotopy type
    (Def. \ref{NilpotentConnectedSpacesOfFiniteRationalType}).
    By means of the system of sets (Def. \ref{FlatLInfinityAlgebraValuedDifferentialForms})
   \vspace{-2mm}
    $$
      X
      \;\;
      \longmapsto
      \;\;
      \Omega_{\mathrm{dR}}
      \big(
        X;
        \,
        \mathfrak{l}A
      \big)
      \;\;\;
      \in
      \;
      \mathrm{Sets}
    $$

    \vspace{-2mm}
    \noindent
    of flat non-abelian differential forms
    with coefficient in the Whitehead $L_\infty$-algebra
    $\mathfrak{l}A$ of $A$ (Prop. \ref{WhiteheadLInfinityAlgebras}),
    which are
    contravariantly assigned to smooth manifolds $X$,
    we consider in
    $\Stacks$ (Def. \ref{SmoothInfinityStacks}):

    \noindent
    {\bf (i)}
    the
    {\it smooth space (Ex. \ref{SmoothSpaces}) of flat $\mathfrak{l}A$-valued differential forms}
    \vspace{-2mm}
    \begin{equation}
      \label{PresheafOfFlatlAValuedDifferentialForms}
      \Omega_{\mathrm{dR}}
      \big(
        -;
        \mathfrak{l}A
      \big)_{\mathrm{flat}}
      \;:=\;
        \bigg(
          \mathbb{R}^n
          \mapsto
          \Big(
            \Delta[k]
            \mapsto
            \Omega_{\mathrm{dR}}
            \big(
             \mathbb{R}^n; \, \mathfrak{l}A
            \big)_{\mathrm{flat}}
          \Big)
        \bigg)
        \,,
    \end{equation}

    \vspace{-2mm}
    \noindent
    regarded as a simplicially constant simplicial presheaf
    \eqref{SimplicialPresheavesOnCartesianSpaces};

    \noindent
    {\bf (ii)}
    the {\it smooth $\infty$-stack of
    flat $\mathfrak{l}A$-valued differential forms}
    (Example \ref{PLDeRhamRightAdjointViaLInfinityAlgebraValuedForms})
    \vspace{-2mm}
    \begin{equation}
      \label{flatdeRhamCoefficients}
      \flatBexp(\mathfrak{l}A)
        \;:=\;
        \bigg(
          \mathbb{R}^n
          \mapsto
          \Big(
            \Delta[k]
            \mapsto
            \Omega_{\mathrm{dR}}
            \big(
             \mathbb{R}^n \times \Delta^k;\, \mathfrak{l}A
            \big)_{\mathrm{flat}}
          \Big)
        \bigg)
    \end{equation}

    \vspace{-2mm}
    \noindent
    which to any Cartesian space
    assigns the simplicial set that in degree $k$ is
    the set of flat $\mathfrak{l}A$-valued differential forms
     on the product manifold of the Cartesian space
     with the standard smooth $k$-simplex $\Delta^k \subset \mathbb{R}^k$;

    \noindent
    {\bf (iii)}
    the canonical inclusion
     \vspace{-3mm}
    \begin{equation}
     \label{InclusionOfFlatNonAbelianDifferentialForms}
        \hspace{-4mm}
        \xymatrix@R=7pt@C=1em{
        \overset{
          \mathclap{
          \raisebox{3pt}{
            \tiny
            \color{darkblue}
            \bf
            \def\arraystretch{.9}
            \begin{tabular}{c}
              smooth space of
              \\
              flat $\mathfrak{l}A$-valued forms
            \end{tabular}
          }
          }
        }{
        \Omega
        (
          -;
          \mathfrak{l}A
        )_{\mathrm{flat}}
        }
        \ar[r]^-{
          \mbox{
            \tiny
            \color{darkblue}
            \bf
            atlas
          }
        }
        \ar@{=}[d]
        &
        \overset{
          \mathclap{
          \raisebox{3pt}{
            \tiny
            \color{darkblue}
            \bf
            \def\arraystretch{.9}
            \begin{tabular}{c}
              smooth $\infty$-stack of
              \\
              flat $\mathfrak{l}A$-valued forms
            \end{tabular}
          }
          }
        }{
        \flatBexp(\mathfrak{l}A)
        }
        \ar@{=}[d]
        \\
        \bigg(
          \mathbb{R}^n
          \mapsto
          \Big(
            \Delta[k]
            \mapsto
            \Omega_{\mathrm{dR}}
            \big(
             \mathbb{R}^n; \, \mathfrak{l}A
            \big)_{\mathrm{flat}}
          \Big)
        \!\!\bigg)
        \;
        \ar@{^{(}->}[r]
        &
        \;
        \bigg(
          \mathbb{R}^n
          \mapsto
          \Big(
            \Delta[k]
            \mapsto
            \Omega_{\mathrm{dR}}
            \big(
             \mathbb{R}^n \times \Delta^k; \, \mathfrak{l}A
            \big)_{\mathrm{flat}}
          \Big)
        \!\!\bigg)
      }
    \end{equation}

    \vspace{-2mm}
    \noindent
    exhibiting
    $\Omega(-; \mathfrak{l}A)$ \eqref{PresheafOfFlatlAValuedDifferentialForms}
    as the presheaf of 0-simplices in the simplicial presheaf
    $\flatBexp(\mathfrak{l}A)$ \eqref{flatdeRhamCoefficients}
    (more abstractly: this is the canonical \emph{atlas} of the
    smooth moduli $\infty$-stack, see \cite[Prop. 2.70]{SS20b}).
\end{defn}

\begin{lemma}[Moduli $\infty$-stack of flat forms
  is equivalent to discrete rational $\infty$-stack]
  \label{ModuliInfinityStackOfFlatlAValuedFormsEquivalentToDiscreteRationalization}
  For $A \in \NilpotentConnectedQFiniteHomotopyTypes$
  (Def. \ref{NilpotentConnectedSpacesOfFiniteRationalType}),
  the evident inclusion
  (by inclusion of polynomial forms into smooth differential forms
  followed by pullback along $\mathrm{pr}_{\Delta^k}$)
  \vspace{-2mm}
  \begin{equation}
    \label{WeakEquivalenceFromDiscreteRationalizationToModuliOfFlatForms}
    \xymatrix@R=9pt{
       \mathllap{
         \mathrm{Disc}
         \big(
           L_{\mathbb{R}}A
         \big)
         \;
         \simeq
         \;\,
       }
        \mathrm{Disc}
        \,\circ\,
        \RightDerived \Bexp_{\mathrm{PL}}
        \,\circ\,
        \mathrm{CE}
        \big(
          \mathfrak{l}A
        \big)
        \ar[r]^-{ \in \, \mathrm{W} }
        \ar@{=}[d]
        &
        \flatBexp \big( \mathfrak{l} A\big)
        \ar@{=}[d]
        \\
        \bigg(\!
          \mathbb{R}^n
          \mapsto
          \Big(
            \Delta[k]
            \mapsto
            \Omega_{\mathrm{PLdR}}
            \big(
              \Delta^k;\, \mathfrak{l}A
            \big)_{\mathrm{flat}}
          \Big)
        \!\!\! \bigg)
        \;
        \ar@{^{(}->}[r]
        &
        \;
        \bigg(\!
          \mathbb{R}^n
          \mapsto
          \Big(
            \Delta[k]
            \mapsto
            \Omega_{\mathrm{dR}}
            \big(
             \mathbb{R}^n \times \Delta^k;\, \mathfrak{l}A
            \big)_{\mathrm{flat}}
          \Big)
        \!\!\! \bigg)
    }
  \end{equation}

  \vspace{-2mm}
\noindent   of the image under $\mathrm{Disc}$ \eqref{ConstantPresheafOverCartesianSpaces}
  of the dg-algebraic model \eqref{RationalizationViaPLDeRham}
  for the rationalization of $A$ (Def. \ref{Rationalization}),
  given by the fundamental theorem (Prop. \ref{FundamentalTheoremOfdgcAlgebraicRationalHomotopyTheory}),
  into the moduli $\infty$-stack of flat
  $\mathfrak{l}A$-valued differential forms (Def. \ref{SmoothCocycleSpaceOfNonAbelianDeRhamCoefficients})
  is an equivalence in $\Stacks$ (Def. \ref{SmoothInfinityStacks}).
\end{lemma}
\begin{proof}
  By Prop. \ref{FundamentalTheoremForPiecewiseSmoothDeRhamComplex},
  the inclusion is for each $\mathbb{R}^n$ a weak equivalence
  \eqref{NaturalEquivalenceBetweenPLAndPSLeftDerivedFunctors}
  in $\SimplicialSets_{\mathrm{Qu}}$ (Example \ref{ClassicalModelStructureOnSimplicialSets}),
  hence is a weak equivalence already in the global projective
  model structure on simplicial presheaves,
  and therefore also in the local projective model structure
  (Example \ref{ModelStructureOnSimplicialPresheavesOverCartesianSpaces}).
\end{proof}

\begin{lemma}[Moduli $\infty$-stack of closed differential forms is shifted de Rham complex]
  \label{ModuliInfinityStackOfClosedFormsIsShiftedDeRhamComplex}
 For $n \in \mathbb{N}$,

  \hspace{-.9cm}
  \begin{tabular}{ll}
  \begin{minipage}[l]{8cm}
   \noindent {\bf (i)}  we have an equivalence
    in $\Stacks$ (Def. \ref{SmoothInfinityStacks})
    from the moduli $\infty$-stack
    $\flatBexp\big( \mathfrak{b}^n\mathbb{R}\big)$
    of flat
    differential forms (Def. \ref{SmoothCocycleSpaceOfNonAbelianDeRhamCoefficients})
    with values in the line Lie $(n+1)$-algebra $\mathfrak{b}^n \mathbb{R}$
    (Example \ref{LineLienPlusOneAlgebras})
    to the image under the Dold-Kan construction
    (Def. \ref{DoldKanConstruction}) of the
    smooth de Rham complex
    $\Omega^\bullet_{\mathrm{dR}}(-)$ (Example \ref{SmoothdeRhamComplex}).

    \vspace{1mm}
\noindent {\bf (ii)}  This is    naturally regarded as a presheaf on $\CartesianSpaces$
    \eqref{CartesianSpaces}
    with values in connective chain complexes (Example \ref{ProjectiveModelStructureOnConnectiveChainComplexes})
    (i.e., with de Rham differential lowering the chain degree)
    shifted up in degree by $n+1$  and
    then homologically truncated in degree 0, as
    shown on the right.
  \end{minipage}
  &
  \;\;\;\;\;
  {\small
  $
    \xymatrix{
      \flatBexp
      \big(
        \mathfrak{b}^n \mathbb{R}
      \big)
      \ar[r]^-{\simeq}
      &
      \mathrm{DK}
      \scalebox{0.8}{$
      \left(\!\!\!
      {\begin{array}{c}
        \vdots
        \\[-3pt]
        \downarrow
        \\[-3pt]
        0
        \\[-3pt]
        \downarrow
        \\[-3pt]
        0
        \\[-3pt]
        \downarrow
        \\[-3pt]
        \Omega^{0}_{\mathrm{dR}}(-)
        \\[-3pt]
        \downarrow\mathrlap{\scalebox{.7}{$d$}}
        \\[-3pt]
        \Omega^{1}_{\mathrm{dR}}(-)
        \\[-3pt]
        \downarrow\mathrlap{\scalebox{.7}{$d$}}
        \\[-3pt]
        \vdots
        \\[-3pt]
        \downarrow\mathrlap{\scalebox{.7}{$d$}}
        \\[-3pt]
        \Omega^{n+1}_{\mathrm{dR}}(-)_{\mathrm{clsd}}
      \end{array}}
     \!\!\! \right)
     $}
    }
    \;
    \in \Stacks
  $
  }
  \end{tabular}
\end{lemma}
\begin{proof}
  This follows by an enhancement of the proof of Prop. \ref{NonAbelianDeRhamCohomologySubsumesOrdinaryDeRhamCohomology}.
  First observe,
  with Example \ref{OrdinaryClosedFormsAreFlatLineLInfinityAlgebraValuedForms},
  that the simplicial presheaf
  \vspace{-1mm}
  \begin{equation}
    \label{ModuliStackOfFlatNFormsAsPresheaf}
    \flatBexp\big( \mathfrak{b}^n \mathbb{R}\big)(-)
    \;\;=\;\;
    \Big(
      \Delta[k]
      \,\mapsto\,
      \Omega^{n+1}_{\mathrm{dR}}
      \big(
        (-) \times \Delta^k
      \big)_{\mathrm{clsd}}
    \Big)
  \end{equation}

   \vspace{-1mm}
\noindent  naturally carries the structure of a presheaf of
  simplicial abelian groups, given by addition of differential forms.
  Therefore, by the Dold-Kan Quillen equivalence (Prop \ref{DoldKanQuillenEquivalence}), it is sufficient to
  prove that we have a quasi-isomorphism of presheaves of chain complexes
  from the corresponding normalized chain complex \eqref{NormalizedChainComplex} of \eqref{ModuliStackOfFlatNFormsAsPresheaf}
  to the shifted and truncated de Rham complex itself:
   \vspace{-2mm}
   \begin{equation}
    \label{FiberIntegrationOfDifferentialFormsAsAChainMap}
    \hspace{-5mm}
    \xymatrix{
      N
      \Big(
        \Delta[k]
        \,\mapsto\,
        \Omega^{n+1}_{\mathrm{dR}}
        \big(
          (-) \times \Delta^k
        \big)_{\mathrm{clsd}}
      \Big)
      \ar[r]_-{ \simeq }^-{
        \int_{\Delta^\bullet}
      }
      &
      \Big(
        \cdots \to 0 \to 0
        \to
        \Omega^0_{\mathrm{dR}}(-)
        \overset{d}{\to}
        \Omega^1_{\mathrm{dR}}(-)
        \overset{d}{\to}
        \cdots
        \overset{d}{\to}
        \Omega^{n+1}_{\mathrm{dR}}(-)_{\mathrm{clsd}}
      \Big)
      \,.
    }
  \end{equation}

   \vspace{-2mm}
\noindent
  We claim that such is given by fiber integration of differential
  forms over the simplices $\Delta^k$:

  First, to see that fiber integration does constitute a chain map,
  we compute for
  $\omega \in \Omega^\bullet_{\mathrm{dR}}\big( (-) \times \Delta^k  \big)_{\mathrm{clsd}}$
  on the left of \eqref{FiberIntegrationOfDifferentialFormsAsAChainMap}:

  \vspace{-2mm}
  \begin{equation}
    \label{FiberIntegrationOverSimplicesIsChainMap}
    \int_{\Delta^k} \partial \omega
    \;\;=\;\;
    (-1)^k
    \int_{\partial \Delta^k} \omega
    \;\;=\;\;
    d \int_{\Delta^k} \omega
    \,,
  \end{equation}
  where the first step is the definition of the
  differential in the normalized chain complex
  \eqref{NormalizedChainComplex}
  and the second step is the
  fiberwise Stokes formula \eqref{StokesFormula}.

  Finally, to see that $\int_{\Delta^\bullet}$ is a quasi-isomorphism,
  notice that the chain homology groups on both sides are

  \vspace{-.3cm}
  $$
    H_k(-)
    \;=\;
    \left\{
    \begin{array}{lcl}
      \mathbb{R} & \vert & k = n + 1
      \\
      0 &\vert& \mbox{otherwise}
    \end{array}
    \right.
      $$
  \vspace{-.3cm}

  \noindent
  over each Cartesian space:
  For the left hand side this follows
  via the weak equivalence \eqref{NaturalEquivalenceBetweenPLAndPSLeftDerivedFunctors}
  from the fundamental theorem (Prop. \ref{FundamentalTheoremOfdgcAlgebraicRationalHomotopyTheory})
  via Example \ref{RationalizationOfEMSpaces},
  while for the right hand side this follows from the
  Poincar{\'e} lemma.

  Hence it is sufficient to see that fiber integration
  over $\Delta^{n+1}$
  is an isomorphism on the $(n+1)$st chain homology groups.
  But a generator of this group on the left is clearly given by
  the pullback
  $\mathrm{pr}_{\Delta^{n+1}}^\ast \omega$
  of any $\omega \in \Omega^{n+1}_{\mathrm{dR}}(\Delta^{n+1})$
  of unit weight and supported in the interior of the simplex.
  That this is sent under $\int_{\Delta^{n+1}}$ to a generator
  $\pm 1 \in \mathbb{R} \simeq \Omega^0_{\mathrm{dR}}(-)_{\mathrm{clsd}}$
  on the right follows by the projection formula \eqref{ProjectionFormula}.
\end{proof}

\begin{remark}[Moduli of closed forms via stable Dold-Kan correspondence]
  \label{ModuliOfClosedFormsViaStableDoldKanCorrespondence}
  Expressed in terms of the stable Dold-Kan construction
  $\mathrm{DK}_{\mathrm{st}}$
  (Prop. \ref{StableDoldKanCorrespondence})
  via the derived stabilization adjunction
  (Example \ref{DerivedStabilizationAdjunction}),
  Lemma \ref{ModuliInfinityStackOfClosedFormsIsShiftedDeRhamComplex}
  says, equivalently, that:

  \vspace{-2mm}
  \begin{equation}
    \label{TheModuliOfClosedFormsViaStableDoldKanCorrespondence}
    \flatBexp\big( \mathfrak{b}^n \mathbb{R} \big)
    \;\;\simeq\;\;
    \mathbb{R}
    \Omega^\infty
    \Big(
    \mathrm{DK}_{\mathrm{st}}
    \big(
      \;
      \Omega^\bullet_{\mathrm{dR}}(-)
        \otimes_{\scalebox{.5}{$\mathbb{R}$}}
      \mathfrak{b}^{n+1}\mathbb{R}
      \;
    \big)
    \Big)
    \;\;\;\;\;
    \in
    \;
    \Stacks
    \,,
  \end{equation}

  \vspace{-1mm}
\noindent  where now
  $\Omega^\bullet_{\mathrm{dR}}(-) \in
    \mathrm{PSh}\big( \CartesianSpaces\,,\,\mathrm{ChainComplexes}_{\scalebox{.5}{$\mathbb{R}$}}\big)$
  is in non-positive degrees, with $\Omega^0_{\mathrm{dR}}(-)$ in degree 0,
  and where $\mathfrak{b}^{n+1}\mathbb{R}$ (Def. \ref{DegreeShift})
  is concentrated on $\mathbb{R}$ in degree $n + 1$.
\end{remark}

\begin{defn}[Differential non-abelian character map {\cite[\S 4]{FSS15b}}]
    \label{DifferentialNonabelianCharacterMap}
%    $\,$
%
%    \noindent
    Given $A \in \NilpotentConnectedQFiniteHomotopyTypes$
    (Def. \ref{NilpotentConnectedSpacesOfFiniteRationalType}),
    the {\it differential non-abelian character map}
    in $A$-cohomology theory,
    to be denoted  $\mathbf{ch}_A$,
    is the morphism in
        $\Stacks$
    \eqref{HomotopyCategoryOfSimplicialPresheaves}
    from $\mathrm{Disc}(A)$ \eqref{ConstantPresheafOverCartesianSpaces}
    to the moduli $\infty$-stack of flat $\mathfrak{l}A$-valued forms
    $\flatBexp(\mathfrak{l}A)$ \eqref{flatdeRhamCoefficients}
    given by the composite
    \vspace{-3mm}
    \begin{equation}
      \label{DifferentialNonabelianCharacterMapPresented}
       \hspace{-7mm}
    \mbox{
      \tiny
      \color{darkblue}
      \bf
      \def\arraystretch{.9}
      \begin{tabular}{c}
        coefficient space as
        \\
        geometrically discrete
        \\
        moduli $\infty$-stack
      \end{tabular}
    }
  \!\!\!\!\!\!\!\!
  \begin{tikzcd}
    \mathrm{Disc}(A)
    \ar[
      r,
      "{
        \mathrm{Disc}(\eta_A^{\mathrm{PSdR}})
      }"
    ]
    \ar[
      rr,
      rounded corners,
      to path={
           -- ([yshift=-10pt]\tikztostart.south)
           --node[below]{
               \scalebox{.7}{
                 $
                   \underset{
                     \mathclap{
                       \mbox{
                         \color{greenii}
                         {\bf
                         $\mathbb{R}$-rationalization}
                         \eqref{RealificationMonad}
                       }
                     }
                   }{
                     \mathrm{Disc}
                     \big(
                       \Derived \eta^{\mathrm{PSdR}}_A
                     \big)
                   }
                 $
               }
             } ([yshift=-10pt]\tikztotarget.south)
           -- (\tikztotarget.south)},
    ]
    \ar[
      rrr,
      rounded corners,
      to path={
           -- ([yshift=+10pt]\tikztostart.north)
           --node[above]{
               \scalebox{.7}{
                 $
                   \overset{
                     \mathclap{
                       \mbox{
                         \color{greenii}
                         \bf
                         differential non-abelian character map
                       }
                     }
                   }{
                     \mathbf{ch}_A
                   }
                 $
               }
             } ([yshift=+10pt]\tikztotarget.north)
           -- (\tikztotarget.north)},
    ]
    &[13pt]
    \mathrm{Disc}
      \circ
    \Bexp_{\mathrm{PS}}
      \circ
    \Omega^\bullet_{\mathrm{PLdR}}
    (A)
    \ar[
      r,
      "{
        \mathrm{Disc}
          \circ
        \Bexp_{\mathrm{PL}}
        (p^{\mathrm{min}})
      }"{above},
      "{
        \mbox{
          \tiny
          \color{darkblue}
          \eqref{CEAlgebraOfWhiteheadLInfinityAlgebra}
        }
      }"{below, yshift=-2pt}
    ]
    &[32pt]
    \mathrm{Disc}
      \circ
    \Bexp_{\mathrm{PL}}
      \circ
    \mathrm{CE}
    ( \mathfrak{l}A )
    \ar[
      r,
      "\in \mathrm{W}"{above},
      "
        \mbox{
          \tiny
          \color{darkblue}
          \bf
          \eqref{WeakEquivalenceFromDiscreteRationalizationToModuliOfFlatForms}
        }
      "{below, yshift=-2pt}
    ]
    &[-10pt]
    \underset{
      \mathclap{
      \mbox{
        \tiny
        \color{darkblue}
        \bf
        \def\arraystretch{.9}
        \begin{tabular}{c}
          moduli $\infty$-stack of
          \\
          flat $\mathfrak{l}A$-valued
          \\
          differential forms
        \end{tabular}
      }
      }
    }{
      \flatBexp(\mathfrak{l}A)
    }
  \end{tikzcd}
    \end{equation}

 \vspace{-3mm}
 \noindent
 of

 \noindent {\bf (a)} the image under $\mathrm{Disc}$ \eqref{ConstantPresheafOverCartesianSpaces}
    of the
    derived adjunction unit $\mathbb{D}\eta^{\mathrm{PLdR}}_A$
    \eqref{DerivedAdjunctionUnit}
    of the PS de Rham adjunction \eqref{PSdRQuillenAdjunctBetweendgcAlgsAndSimplicialSets},
    specifically with (co-)fibrant replacement
    $p^{\mathrm{min}}$ being the
    minimal Sullivan model replacement \eqref{MinimalSullivanModel};
    (recalling that $\Bexp_{\mathrm{PL}}$ is a contravariant functor),
   with

\noindent {\bf (b)} the weak equivalence from Lemma \ref{ModuliInfinityStackOfFlatlAValuedFormsEquivalentToDiscreteRationalization}.
  \end{defn}

\begin{remark}[Differential non-abelian character map is independent of choices]
  The differential non-abelian character map (Def. \ref{DifferentialNonabelianCharacterMap})
  is independent, up to equivalence,
  of the choice of comparison morphism $p^{\mathrm{min}}$ to
  a minimal model for the coefficients, since, by
  \eqref{FactorizationOfMinimalModelsThroughEachOther}
  in Prop. \ref{ExistenceOfMinimalSullivanModels},
  any two choices factor through each other by an isomorphism of dgc-algebras.

  It is this uniqueness which makes minimal models provide
  canonical form coefficients for non-abelian differential cohomology,
  see also
  the second item of Ex. \ref{DifferentialGeneralizedCohomology} below.
\end{remark}

\medskip

\noindent {\bf Differential non-abelian cohomology.}

\begin{defn}[Differential non-abelian cohomology {\cite[\S 4]{FSS15b}}]
  \label{DifferentialNonAbelianCohomology}
  For $A \in \NilpotentConnectedQFiniteHomotopyTypes$ (Def. \ref{NilpotentConnectedSpacesOfFiniteRationalType})
  we say that:

  \noindent
  {\bf (i)}
  the {\it moduli $\infty$-stack of $\Omega A$-connections}
  is the object
  $
    A_{\mathrm{diff}}
      \;\in\
    \Stacks
  $
  in the homotopy category of smooth $\infty$-stacks
  (Def. \ref{SmoothInfinityStacks}),
  which is
  given by the homotopy pullback (Def. \ref{HomotopyPullback})
  of the smooth space of flat non-abelian differential forms
  $\Omega_{\mathrm{dR}}(-; \mathfrak{l}A)_{\mathrm{flat}}$
  \eqref{InclusionOfFlatNonAbelianDifferentialForms}
  along the differential non-abelian character map
  $\mathbf{ch}_A$ (Def. \ref{DifferentialNonabelianCharacterMap}):
   \vspace{-2mm}
  \begin{equation}
    \label{PullbackForModuliInfinityStackOfConnections}
    \raisebox{0pt}{
    \xymatrix@C=5em@R=1.2em{
      \overset{
        \mathllap{
        \raisebox{3pt}{
          \tiny
          \color{darkblue}
          \bf
          \def\arraystretch{.9}
          \begin{tabular}{c}
            moduli $\infty$-stack
            \\
            of $\Omega A$-connections
          \end{tabular}
        }
        }
      }{
        A_{\mathrm{diff}}
      }
      \ar[d]^-{c_A}_-{
        \mathllap{
          \mbox{
            \tiny
            \color{darkblue}
            \bf
            {
            \def\arraystretch{.9}
            \begin{tabular}{c}
              universal
              characteristic class
              \\
              in non-abelian $A$-cohomology
            \end{tabular}}
          }
          \!\!\!\!
        }
  %      c_A
      }
      \ar[rr]_-{F_A}^-{
        \overset{
          \raisebox{3pt}{
            \tiny
            \color{darkblue}
            \bf
            {
            \def\arraystretch{.9}
            \begin{tabular}{c}
              $\mathfrak{l}A$-valued
              \\
              curvature forms
            \end{tabular}}
          }
        }{
%          F_A
        }
      }
      \ar@{}[drr]|-{
        \mbox{
          \tiny
          (hpb)
        }
      }
      &&
      \overset{
        \mathrlap{
        \raisebox{3pt}{
          \tiny
          \color{greenii}
          \bf
          \def\arraystretch{.9}
          \begin{tabular}{c}
            smooth space of
            \\
            flat $\mathfrak{l}A$-valued forms
          \end{tabular}
        }
        }
      }{
        \Omega_{\mathrm{dR}}
        (
          -;
          \,
          \mathfrak{l}A
        )_{\mathrm{flat}}
      }
      \ar[d]^-{
        \mbox{
          \tiny
          \color{greenii}
          \bf
          atlas
        }
      }
      \\
      \mathrm{Disc}(A)
      \ar[rr]^-{ \mathbf{ch}_A}_-{
        \underset{
          \mathclap{
          \raisebox{-3pt}{
            \tiny
            \color{greenii}
            \bf
            {
            \def\arraystretch{.9}
            \begin{tabular}{c}
              differential non-abelian
              \\
              character map
            \end{tabular}}
          }
          }
        }{
%          \mathbf{ch}_A
        }
      }
      &&
      \underset{
        \mathrlap{
        \raisebox{-3pt}{
          \tiny
          \color{greenii}
          \bf
          \def\arraystretch{.9}
          \begin{tabular}{c}
            moduli $\infty$-stack of
            \\
            flat $\mathfrak{l}A$-valued forms
          \end{tabular}
        }
        }
      }{
        \flatBexp(\mathfrak{l}A)
      }
    }
    }
    \;\;\;\;\;\;\;\;\;\;\;\;
    \in
    \;\;\;
    \Stacks
    \,;
  \end{equation}

 \vspace{-2mm}
  \noindent
  {\bf (ii)}
  the \emph{differential non-abelian cohomology}
  of a smooth $\infty$-stack $\mathcal{X} \,\in\,\Stacks$
  \eqref{HomotopyCategoryOfSimplicialPresheaves}
  with coefficients in $A$
  is the structured non-abelian cohomology
  (Remark \ref{StructuredNonAbelianCohomology})
  with coefficients in the
  moduli $\infty$-stack $A_{\mathrm{diff}}$ of $\Omega A$-connections \eqref{PullbackForModuliInfinityStackOfConnections},
  hence the hom-set
  in the homotopy category of $\infty$-stacks
  (Def. \ref{SmoothInfinityStacks})
  from $\mathcal{X}$ to $A_{\mathrm{diff}}$
 \vspace{-1mm}
    \begin{equation}
    \label{DifferentialNonabelianCohomology}
    \widehat H
    \big(
      \mathcal{X};
      \,
      A
    \big)
    \;\;
      :=
    \;\;
    H
    \big(
      \mathcal{X}
      ;\,
      A_{\mathrm{diff}}
    \big)
    \;\;
      :=
    \;\;
    \Stacks
    \big(
      \mathcal{X}
      \,
      ,
      \,
      A_{\mathrm{diff}}
    \big).
  \end{equation}

  \noindent
  {\bf (iii)} We call the non-abelian cohomology operations
  induced from the maps in \eqref{PullbackForModuliInfinityStackOfConnections} as follows
  (see \eqref{SystemsOfCohomologyOperationsOnDifferentialCohomology}):

  \vspace{-.5cm}
  \begin{align}
    \label{CharacteristicClassOnDifferentialCohomology}
    \mbox{ {\bf (a)} {\it characteristic class}: \phantom{aa} }
    &
    \xymatrix{
      \widehat H
      \big(
        \mathcal{X};
        \,
        A
      \big)
      \ar[rr]^-{ (c_A)_\ast }
      &&
      H
      \big(
        \mathrm{Shp}(\mathcal{X});
        \,
        A
      \big)
    }
    &
    \mathllap{
      \mbox{(Def. \ref{NonAbelianCohomology})}
    }
    \\
    \label{CurvatureOnDifferentialCohomology}
    \mbox{ {\bf (b)} {\it curvature}: \phantom{characterisa} }
    &
    \xymatrix{
      \widehat H
      \big(
        \mathcal{X};
        \,
        A
      \big)
      \ar[rr]^-{ (F_A)_\ast }
      &&
      \Omega_{\mathrm{dR}}
      \big(
        \mathcal{X};
        \,
        \mathfrak{l}A
      \big)_{\mathrm{flat}}
    }
    &
    \mathllap{
      \mbox{(Def. \ref{FlatLInfinityAlgebraValuedDifferentialForms})}
    }
    \\
    \label{DifferentialCharacterOnDifferentialCohomology}
    \mbox{ {\bf (c)} {\it differential character}:\phantom{a}\, }
    &
    \xymatrix{
      \widehat H
      \big(
        \mathcal{X};
        \,
        A
      \big)
      \ar[rr]^-{ (\mathbf{ch}_A \circ c_A)_\ast }
      &&
      H_{\mathrm{dR}}
      \big(
        \mathcal{X};
        \,
        \mathfrak{l}A
      \big)
    }
    &
    \mathllap{
      \mbox{
        (Def. \ref{NonabelianDeRhamCohomology})}
    }
  \end{align}

\end{defn}

In differential enhancement of Example \ref{GeneralizedCohomologyAsNonabelianCohomology},
we have the following.

\medskip

\noindent {\bf Differential generalized cohomology.}

\begin{example}[Differential Whitehead-generalized cohomology]
  \label{DifferentialGeneralizedCohomology}
  Let $E^\bullet$ be a generalized cohomology theory
  (Example \ref{GeneralizedCohomologyAsNonabelianCohomology})
  with representing spectrum $E$ \eqref{Spectrum}
  which is connective and whose
  component spaces $E_n$ are of finite $\mathbb{R}$-type,
  so that their connected components are, by
  Example \ref{ExamplesOfNilpotentSpaces},  in
  $\NilpotentConnectedQFiniteHomotopyTypes$ (Def. \ref{NilpotentConnectedSpacesOfFiniteRationalType}).

\vspace{-1mm}
\item {\bf (i)}   Then differential non-abelian cohomology,
  in the sense of Def. \ref{DifferentialNonAbelianCohomology},
  with coefficients in the component spaces $E_\bullet$,
  coincides with
  canonical differential generalized $E$-cohomology
  in the traditional sense of \cite[\S 4.1]{HopkinsSinger05}\cite[Def. 4.53]{Bunke12}\cite[\S 2.2]{BunkeGepner13}\cite[\S 4.4]{BNV13}:

  \vspace{-.6cm}
  \begin{equation}
   \label{RecoveringDifferentialGeneralizedCohomology}
    \overset{
      \mathclap{
      \raisebox{3pt}{
        \tiny
        \color{darkblue}
        \bf
        \def\arraystretch{.9}
        \begin{tabular}{c}
          generalized
          \\
          differential cohomology
        \end{tabular}
      }
      }
    }{
      \widehat E^n(-)
    }
    \;\;\;\;\;\;\simeq\;\;
    \widehat H
    (
      -;
      E_n
    )
    \,.
  \end{equation}

\vspace{-1mm}
  \item {\bf (ii)}
  Here ``canonical'', in the sense of \cite[Def. 4.46]{Bunke12},
  refers to choosing the curvature differential form coefficients
  to be $\pi_\bullet(E) \otimes \mathbb{R}$
  (instead of some chain complex quasi-isomorphic to this).
  By  Example \ref{WhiteheadLInfinityAlgebraOfLoopSpaces},
  this choice corresponds in
  our Def. \ref{DifferentialNonAbelianCohomology}
  to the \emph{minimality} (Def. \ref{MinimalSullivanModels})
  of the minimal Sullivan model
  $\mathrm{CE}(\mathfrak{l}E_n)$ for $E_n$
  (Prop. \ref{WhiteheadLInfinityAlgebras})
  that controls the flat
  $L_\infty$-algebra valued differential forms
  $\Omega_{\mathrm{dR}}
  (-;\, \mathfrak{l}E_n)_{\mathrm{flat}}$ (Def. \ref{FlatLInfinityAlgebraValuedDifferentialForms})
  in the top right of \eqref{RecoveringOrdinaryDifferentialCohomology}.

\vspace{-1mm}
  \item {\bf (iii)} Hence for canonical/minimal curvature coefficients,
  we have from
  Ex. \ref{WhiteheadLInfinityAlgebraOfLoopSpaces},
  Lem. \ref{ModuliInfinityStackOfClosedFormsIsShiftedDeRhamComplex}
  and
  Rem. \ref{ModuliOfClosedFormsViaStableDoldKanCorrespondence}
   that
   \vspace{-2mm}
  \begin{align}
    \flatBexp\big( \mathfrak{l}E_n\big)
    \;\;
    &  \simeq
    \;\;
    \mathbb{R}\Omega^\infty
    \Big(
      \mathrm{DK}_{\mathrm{st}}
      \big(
        \Omega^\bullet_{\mathrm{dR}}(-)
        \otimes_{\scalebox{.5}{$\mathbb{Z}$}}
        \pi_\bullet(E_n)
      \big)
    \Big)
    \;\;\;\;\;
    \in
    \;
    \Stacks
    \\
%  \end{equation}
%  and
%  \begin{equation}
    \label{AtlasFromTruncation}
    \Omega_{\mathrm{dR}}
    \big(
      -;
      \mathfrak{l}E_n
    \big)_{\mathrm{flat}}
    \;\;
    &  \simeq
    \;\;
    \mathbb{R}\Omega^\infty
    \Big(
      \mathrm{DK}_{\mathrm{st}}
      \big(
        \Omega^\bullet_{\mathrm{dR}}(-)
        \otimes_{\scalebox{.5}{$\mathbb{Z}$}}
        \pi_\bullet(E_n)
      \big)_{\leq 0}
    \Big)
    \;\;\;\;\;
    \in
    \Stacks\;.
  \end{align}

\vspace{-2mm}
 \item {\bf (iv)} With this, the equivalence \ref{RecoveringDifferentialGeneralizedCohomology}
  follows by
  Ex. \ref{RealRationalizationOfSpectra} and
  observing that the defining
  homotopy pullback diagram \eqref{PullbackForModuliInfinityStackOfConnections}
  for differential non-abelian cohomology with coefficients
  in $A := E_n$ \eqref{ComponentSpacesOfSpectra}
  is the image under $\mathbb{R}\Omega^\infty$
  \eqref{StabilizationAdjunction} of the
  defining homotopy pullback diagram for canonical
  differential $E$-cohomology according to
  \cite[(4.12)]{HopkinsSinger05} \cite[Def. 4.51]{Bunke12}\cite[(24)]{BNV13},
  and using that right adjoints preserve homotopy pullbacks:
  \vspace{-2mm}
  \begin{equation}
    \label{RelatingModuliStacksOfConnectionsToDifferentialFunctionSpectra}
    \underset{
      \mathclap{
      \raisebox{-9pt}{
        \tiny
        \color{darkblue}
        \bf
       \def\arraystretch{.9}
        \begin{tabular}{c}
          moduli $\infty$-stack
          \\
          of $\Omega E_0$-connections
        \end{tabular}
      }
      }
    }{
    \raisebox{24pt}{
    \xymatrix@C=15pt{
      (E_0)_{\mathrm{diff}}
      \ar[d]_-{ c_{E_0} }
      \ar[rr]^-{ F_{E_0} }
      \ar@{}[drr]|-{
        \mbox{
          \tiny
          \rm
          (hpb)
        }
      }
      &&
      \Omega_{\mathrm{dR}}
      (
        -;
        \mathfrak{l}E_0
      )_{\mathrm{flat}}
      \ar[d]^-{
        \mathclap{\phantom{\vert}}
        \mathrm{atlas}
      }
      \\
      \mathrm{Disc}
      (
       E_0
      )
      \ar[rr]_-{ \mathbf{ch}_{E_0} }
      &&
      \flat
      \Bexp
      (
        \mathfrak{l}E_0
      )
    }
    }
    }
    \;
      \simeq
    \;\;\;
    \mathbb{R}
    \Omega^\infty
    \underset{
      \mathclap{
      \raisebox{-6pt}{
        \tiny
        \color{darkblue}
        \bf
        \def\arraystretch{.9}
        \begin{tabular}{c}
          ``differential function spectrum''
          \\
          of differential generalized $E$-cohomology
        \end{tabular}
      }
      }
    }{
    \left(
    \!\!\!\!
    \raisebox{24pt}{
    \xymatrix@C=15pt{
      \mathrm{Diff}
      (
        E, \mathrm{can}
      )
      \ar[d]
      \ar[rr]
      \ar@{}[drr]|-{
        \mbox{
          \tiny
          \rm
          (hpb)
        }
      }
      &&
      \big(
        \Omega^\bullet_{\mathrm{dR}}(-)
        \otimes_{\scalebox{.5}{$\mathbb{Z}$}}
        \pi_\bullet(E)
      \big)_{\leq 0}
      \ar[d]
      \\
      \mathrm{Disc}
      (
        E
      )
      \ar[rr]_-{
        H\mathbb{R}\wedge(-)
      }
      &&
      \Omega^\bullet_{\mathrm{dR}}(-)
        \otimes_{\scalebox{.5}{$\mathbb{Z}$}}
      \pi_\bullet(E)
    }
    }
    \!\!\!\!
    \right)
    }
  \end{equation}

  \vspace{-2mm}
\noindent  The same applies to $(E_n)_{\mathrm{diff}}$,
  by replacing $E$ with $\LeftDerived\Sigma^n E$ \eqref{StabilizationAdjunction}
  on the right of \eqref{RelatingModuliStacksOfConnectionsToDifferentialFunctionSpectra}.
\end{example}

\begin{remark}[The canonical atlas for the moduli stack of connections]
  The operation $(-)_{\leq 0}$ in \eqref{AtlasFromTruncation}
  is the naive truncation functor on the category of chain complexes

  \vspace{-.4cm}
  $$
    \xymatrix@R=1pt{
      \mathrm{ChainComplexes}_{\scalebox{.5}{$\mathbb{Z}$}}
      \ar[rr]^-{ (-)_{\leq 0} }
      &&
      \mathrm{ChainComplexes}^{\scalebox{.5}{$\leq 0$}}_{\scalebox{.5}{$\mathbb{Z}$}}
      \\
      \big(
        \cdots
        \overset{\partial_1}{\longrightarrow}
        V_1
        \overset{\partial_0}{\longrightarrow}
        V_0
        \overset{\partial_{-1}}{\longrightarrow}
        V_{-1}
        \overset{\partial_{-2}}{\longrightarrow}
        V_{-1}
        \to
        \cdots
      \big)
      \ar@{}[rr]|-{ \longmapsto }
      &&
      \big(
        V_0
        \overset{\partial_{-1}}{\longrightarrow}
        V_{-1}
        \overset{\partial_{-2}}{\longrightarrow}
        V_{-1}
        \to
        \cdots
      \big)
      \,.
    }
  $$
  \vspace{-.2cm}

  \noindent
  In contrast to the homological truncation
  involved in $\Omega^\infty$ \eqref{HomologicalTruncationFromBelow},
  this naive truncation is not homotopy-invariant and does
  not have a derived functor.
  Instead, as seen from \eqref{AtlasFromTruncation}
  and \eqref{InclusionOfFlatNonAbelianDifferentialForms},
  once regarded in differential non-abelian cohomology, this
  operation serves to construct the canonical
  \emph{atlas} \cite[Prop. 2.70]{SS20b}
  of the moduli $\infty$-stack of flat
  $\mathfrak{l}E_n$-valued differential forms.
  Via the defining homotopy pullback \eqref{PullbackForModuliInfinityStackOfConnections},
  \eqref{RelatingModuliStacksOfConnectionsToDifferentialFunctionSpectra}
  this becomes hallmark of differential cohomology:
  Differential cohomology is the universal solution to
  lifting the values of the character map from cohomology
  classes to cochain representatives, namely to curvature forms.
\end{remark}

In differential enhancement of
Example \ref{ComplexTopologicalKtheory} and
Example \ref{ChernCharacterInKTheory},
we have:
\begin{example}[Differential complex K-theory]
  \label{CurvatureInDifferentialComplexKTheory}
  With the coefficient space
  $A \,:=\, \mathrm{KU}_0 \,=\, \mathbb{Z} \times B \mathrm{U}$
  \eqref{ClassifyingSpaceForComplexTopologicalKTheory}
  for topological complex K-theory (Example \ref{ComplexTopologicalKtheory}),
  the corresponding differential non-abelian cohomology theory
  (Def. \ref{DifferentialNonAbelianCohomology})
  is, by Example \ref{DifferentialGeneralizedCohomology},
  differential K-theory, whose diagram
  \eqref{SystemsOfCohomologyOperationsOnDifferentialCohomology}
  of cohomology operations is of this form
  (see \cite{HopkinsSinger05}\cite{BunkeSchick09}\cite{BunkeSchick12}\cite{GS-AHSS})
  \vspace{-2mm}
  \begin{equation}
    \label{CohomologyOperationsDiagramForDifferentialKTheory}
    \xymatrix@R=11pt{
      \mathllap{
        \widehat H
        \big(
          \mathcal{X};
          \,
          \mathrm{KU}_0
        \big)
        \;
        \simeq
        \,
      }
      \widehat{\mathrm{KU}}^0(\mathcal{X})
      \ar[rr]^-{ F_{\mathrm{KU}_0} }
      \ar[d]_{ c_{\mathrm{KU}_0} }
      &&
      \Big\{
        \left.
          \big\{
            F_{2k}
            \,\in\,
            \Omega^{2k}_{\mathrm{dR}}(\mathcal{X})
          \big\}_{k \in \mathbb{N}}
        \,\right\vert\,
        d\, F_{2k} = 0
      \Big\}
      \ar[d]
      \\
      \mathrm{KU}^0(\mathcal{X})
      \ar[rr]^-{ \mathrm{ch} }
      &&
      \underset{k \in \mathbb{N}}{\bigoplus}
      H^{2k}_{\mathrm{dR}}
      \big(
        \mathcal{X}
      \big)
      \,,
    }
  \end{equation}

  \vspace{-2mm}
  \noindent  where the bottom map is the ordinary Chern character
  from Example \ref{ChernCharacterInKTheory}, and the
  curvature differential forms are identified as in
  Example \ref{RecoveringH3TwistedDifferentialForms}.
\end{example}
\begin{remark}[Differential K-theory via equivalence classes of principal connections]
\label{DifferentialKTheoryViaEquivalenceClassesOfPrincipalConnections}
In our context of non-abelian cohomology it is
worth highlighting the well-known fact that
differential K-theory classes (Ex. \ref{CurvatureInDifferentialComplexKTheory})
may equivalently be expressed
(\cite{Karoubi87}\cite{Lott94}\cite{SimonsSullivan08Structured}\cite[\S 6]{BNV13},
brief review in \cite[\S 4.1]{BunkeSchick12})
in terms of equivalence classes
of {\it vector bundles with connection}, hence equipped with
principal connections (Nota. \ref{PrincipalBundlesWithConnection})
on the underlying $\mathrm{U}(n)$-principal bundles.
\end{remark}

\medskip

\noindent {\bf Examples of differential non-abelian cohomology.}
In differential enhancement of Example \ref{TraditionalNonAbelianCohomology},
we have:
\begin{prop}[Differential non-abelian cohomology of principal connections]
  \label{DifferentialCohomologyOfPrincipalConnection}
  Let $G$ be a compact Lie group
  with classifying space $B G$ \eqref{ClassifyingSpace}.
  Then there is a natural map over smooth manifolds $X$,
  shown dashed in \eqref{FromGConnectionsToDifferentialNonabelianCohomology},
  from equivalence classes of $G$-principal connections
  (Notation \ref{PrincipalBundlesWithConnection})
  to differential non-abelian cohomology with coefficients in
  $B G$ (Def. \ref{DifferentialNonAbelianCohomology})
  which covers the classification of $G$-principal bundles
  by plain non-abelian cohomology with coefficients in
  $B G$ (Example \ref{TraditionalNonAbelianCohomology}),
  in that the following diagram commutes:
  \vspace{-4mm}
  \begin{equation}
    \label{FromGConnectionsToDifferentialNonabelianCohomology}
    \xymatrix@R=1em{
      G\mathrm{Connections}(X)_{/\sim}
      \ar@{-->}[rr]
      \ar[d]_-{
        \mbox{
          \tiny
          \color{darkblue}
          \bf
          \def\arraystretch{.9}
          \begin{tabular}{c}
            forget
            \\
            connection
          \end{tabular}
        }
        \!\!
      }
      &&
      \overset{
        \mathclap{
        \raisebox{3pt}{
          \tiny
          \color{darkblue}
          \bf
          \def\arraystretch{.9}
          \begin{tabular}{c}
            differential
            \\
            non-abelian cohomology
          \end{tabular}
        }
        }
      }{
        \widehat H
        (
          X;
          B G
        )
      }
      \ar[d]^-{ c_{B G} }
      \\
      G\mathrm{Bundles}(X)_{/\sim}
      \ar@{->}[rr]_-{ \simeq }
      &&
      \underset{
        \mathclap{
        \raisebox{-3pt}{
          \tiny
          \color{darkblue}
          \bf
          \begin{tabular}{c}
            non-abelian
            cohomology
          \end{tabular}
        }
        }
      }{
        H
        (
          X;
          \,
          B G
        )
      }
    }
  \end{equation}
\end{prop}
\begin{proof}
  By Lemma \ref{NonAbelianDeRhamCohomologyWithCoefficientsInAClassifyingSpace},
  the differential form coefficient in the given case is
    \vspace{-2mm}
  $$
    \Omega_{\mathrm{dR}}
    (
      -;
      \mathfrak{l}BG
    )_{\mathrm{flat}}
    \;\;\simeq\;\;
    \mathrm{Hom}_{\mathbb{R}}
    \Big(
      \mathrm{inv}^\bullet(\mathfrak{g})
      \,,\,
      \Omega^\bullet_{\mathrm{dR}}(-)_{\mathrm{clsd}}
    \Big)
    \,.
  $$

  \vspace{-2mm}
  \noindent
  Therefore, with Example \ref{RationalizationOfEMSpaces}, we find that
    \vspace{-2mm}
  $$
    \Big(
      \Delta[k]
      \,\mapsto\,
      \mathrm{Hom}_{\mathbb{R}}
      \big(
        \mathrm{inv}^\bullet(\mathfrak{g})
        \,,\,
        \Omega^\bullet_{\mathrm{dR}}(\Delta^k)_{\mathrm{clsd}}
      \big)
    \Big)
    \;\;
    \simeq
    \;\;
    \underset{k}{\prod}
    \,
    K
    \big(
      \mathrm{inv}^n(\mathfrak{g}),
      n
    \big)
    \;\;\;\;
    \in
    \mathrm{Ho}
    \big(
      \SimplicialSets_{\mathrm{Qu}}
    \big)
  $$

  \vspace{-2mm}
  \noindent
  is a product of Eilenberg-MacLane spaces \eqref{EilenbergMacLaneSpaces}
  for real coefficient groups spanned by the
  invariant polynomials, and so the defining
  homotopy pullback \eqref{PullbackForModuliInfinityStackOfConnections}
  is here of the following form:
  \vspace{-2mm}
  \begin{equation}
    \label{HomotopyPullbackDefinitionForBGdiff}
    \xymatrix@R=1.5em@C=3em{
      B G_{\mathrm{diff}}
      \ar[rr]
      \ar[d]
      \ar@{}[drr]|-{\mbox{\tiny(hpb)}}
      &&
      \mathrm{Hom}_{\mathbb{R}}
      \big(
        \mathrm{inv}^\bullet(\mathfrak{g})
        \,,\,
        \Omega^\bullet_{\mathrm{dR}}(-)_{\mathrm{clsd}}
      \big)
      \ar[d]
      \\
      \mathrm{Disc}(B G)
      \ar[rr]_-{
        (c_k)_{k \in \mathbb{N}}
      }
      &&
      \mathrm{Disc}
      \Big(\;
        \underset{
          k
          \in
          \mathbb{N}
        }{\prod}
        K
        \big(
          \mathrm{inv}^n(\mathfrak{g}),
          n
         \big)
      \Big)
      \,,
    }
  \end{equation}

  \vspace{-2mm}
  \noindent
  where the bottom map classifies the real characteristic
  classes of $B G$ via Example \ref{OrdinaryCohomology}.
  It follows by Example \ref{HomotopyPullbackViaTriples}
  that maps into $B G_{\mathrm{diff}}$ are equivalence classes
  of triples
    \vspace{-1mm}
  \begin{equation}
    \label{DataTripleInDifferentialNonabelianCohomology}
    \widehat H
    (
      X;\,
      B G
    )
    \;\;
     \simeq
    \;\;
    \left\{
      \big(
        f
        ,
        \phi
        ,
        (\alpha_k)
      \big)
      \,\left\vert\,
    \raisebox{25pt}{
    \xymatrix@R=1.5em@C=5em{
      X
      \ar@{-->}[r]^-{  (\alpha_k) }_->>>{\ }="s"
      \ar@{-->}[d]_-{f}^->>>{\ }="t"
      &
      \mathrm{Hom}_{\mathbb{R}}
      \big(
        \mathrm{inv}^\bullet(\mathfrak{g})
        \,,\,
        \Omega^\bullet_{\mathrm{dR}}(-)_{\mathrm{clsd}}
      \big)
      \ar[d]^-{}
      \\
      B G
      \ar[r]_-{}
      &
      \mathrm{Disc}
      \Big(\;
        \underset{
          k
          \in
          \mathbb{N}
        }{\prod}
        K
        \big(
          \mathrm{inv}^n(\mathfrak{g}),
          n
         \big)
      \Big)
      \ar@{==>}^\phi "s"; "t"
    }
    }
    \right.
    \right\}
  \end{equation}

    \vspace{-1mm}
    \noindent
  consisting of
  {\bf (a)} a classifying map $f$ for a $G$-principal bundle
  (Example \ref{TraditionalNonAbelianCohomology}),
  {\bf (b)} a set of closed differential forms $\alpha$ labeled by
  the invariant polynomials,
  and {\bf (c)} a set of coboundaries $\phi$ in real cohomology
  between these differential
  forms and the pullbacks $f^\ast c_k$.

  Now, given a $G$-connection $\nabla$
  on a $G$-principal bundle $f^\ast E G$ over $X$, we
  obtain such a triple by
  {\bf (a)} taking $f$ to be
  the classifying map of the underlying $G$-principal bundle,
  {\bf (b)} taking $\alpha_k := \omega_k(F_\nabla)$
  to be the characteristic forms (Def. \ref{CharacteristicForms})
  of the connection, and {\bf (c)} taking $\phi$
  to be given by the relative Chern-Simons forms
  \cite{ChernSimons74} between the given connection and the
  pullback along $f$ of the universal connection
  (see Remark \ref{ChernWeilTheoryAndItsFundamentalTheorem}).
  This construction is an invariant of the
  isomorphism class of the connection
  (see \cite[p. 28]{HopkinsSinger05}) and hence
  defines the desired map \eqref{FromGConnectionsToDifferentialNonabelianCohomology}:
    \vspace{-2mm}
  \begin{equation}
    \label{GConnectionsToTriples}
    \xymatrix@R=-3pt{
      G\mathrm{Connections}(X)_{/\sim}
      \ar[rr]
      &&
      \widehat H
      (
        X;
        \,
        B G
      )
      \\
      \big[
        f^\ast E G, \nabla
      \big]
      \ar@{}[rr]|-{
        \longmapsto
      }
      &&
      \big[
        f,
        \,
        \big(
          \mathrm{cs}_k
          (
            \nabla, f^\ast \nabla_{\mathrm{univ}}
          )
        \big),
        \,
        \big(
          \omega_k(F_\nabla)
        \big)
      \big]
    }
  \end{equation}

  \vspace*{-1.5\baselineskip}
\end{proof}

\begin{remark}[The role of principal connections in non-abelian differential cohomology]
  \label{PrincipalConnectionsInNonabelianDifferentialCohomology}
  $\,$

  \noindent
  {\bf (i)}
  It seems unlikely that the map
  \eqref{FromGConnectionsToDifferentialNonabelianCohomology} in
  Prop. \ref{DifferentialCohomologyOfPrincipalConnection}
  would not be a bijection, but we do not have a proof that it is,
  in general. A notable case where it is known to be a bijection is
  the abelian case of the circle group $G = \mathrm{U}(1)$;
  this case is Prop. \ref{DifferentialNonabelianCohomologySubsumesOrdinaryDifferentialCohomology}
  below.

  \noindent
  {\bf (ii)}
  However,
  Thm. \ref{SecondaryDifferentialNonAbelianCharacterSubsumesCheegerSimonsHomomorphism}
  below shows that the image of
  $G$-principal connections in differential non-abelian cohomology
  $\widehat{H}\big(X;  B G\big)$, under this map
  \eqref{FromGConnectionsToDifferentialNonabelianCohomology},
  supports the construction of
  all the secondary characteristic classes of $G$-principal bundles,
  hence retains all the relevant information
  extractable from $G$-principal connections.

  \noindent
  {\bf (iii)}
  On the other hand, for each Lie group $G$ with
  Lie algebra
  denoted $\mathfrak{g}$, there exists a smooth stack (Prop. \ref{SmoothInfinityStacks})

  \vspace{-.4cm}
  \begin{equation}
    \label{BGConnAsQuotientStackOfSmoothFormsByGaugeTransformations}
    \mathbf{B}G_{\mathrm{conn}}
    \;\simeq\;
    \Omega^1(-;\mathfrak{g}) \!\sslash\! G
    \;\;\;
    \in
    \;
    \Stacks
  \end{equation}
  \vspace{-.5cm}

  \noindent
  which is the {\it moduli stack} of smooth $G$-principal connections
  (\cite[Def. 3.2.4]{FiorenzaSchreiberStasheff10}\cite{FreedHopkins13},
  exposition in \cite[\S 2.4]{FSS13a}) in that it not only
  makes
  the
  analogue of the
  map \eqref{FromGConnectionsToDifferentialNonabelianCohomology}
  provably a bijection

  \vspace{-.4cm}
  $$
    \begin{tikzcd}
      G\mathrm{Connections}(X)_{/\sim}
      \ar[
        r,
        "\sim"{above}
      ]
      &
      H\big(X;, \mathbf{B}G_{\mathrm{conn}} \big)
    \end{tikzcd}
    \;\;\;
    \in
    \;
    \Sets
  $$
  \vspace{-.5cm}

  \noindent
  but even such that the
  full mapping space \eqref{MappingSpaceBetweenSmoothStacks}
  into it is equivalent
  (\cite[Prop. 3.2.5]{FiorenzaSchreiberStasheff10})
  to the groupoid
  (via Ex. \ref{SimplicialNervesOfGroupoids})
  of
  gauge transformations between
  $G$-principal connections:

  \vspace{-.4cm}
  $$
    \begin{tikzcd}
      G\mathrm{Connections}(X)
      \ar[
        r,
        "\sim"{above}
      ]
      &
      \mathrm{Maps}\big(X;, \mathbf{B}G_{\mathrm{conn}} \big)
    \end{tikzcd}
    \;\;\;
    \in
    \;
    \mathrm{Ho}
    \big(
      \SimplicialSets_{\mathrm{Qu}}
    \big)
    \,.
  $$
  \vspace{-.5cm}

  \noindent
  {\bf (iv)}
  But, while
  $\mathbf{B}G_{\mathrm{conn}}$
  can explicitly be defined as in \eqref{BGConnAsQuotientStackOfSmoothFormsByGaugeTransformations},
  it seems to lack (unless $G$ is abelian, see Prop. \ref{DifferentialNonabelianCohomologySubsumesOrdinaryDifferentialCohomology})
  a more general abstract characterization
  of the kind that defines $B G_{\mathrm{diff}}$
  in \eqref{HomotopyPullbackDefinitionForBGdiff},
  via the systematic Def. \ref{DifferentialNonAbelianCohomology}.
  In particular, it is the construction principle of
  $B G_{\mathrm{diff}}$ -- but apparently not that of
  $\mathbf{B}G_{\mathrm{conn}}$ -- which properly generalizes
  from ordinary non-abelian Lie groups to higher non-abelian groups
  \cite[\S 4.3]{FiorenzaSchreiberStasheff10} such as the
  String 2-group (Ex. \ref{NonAbelianCohomologyInDegree2}),
  again for the fact that $BG_{\mathrm{diff}}$
  canonically supports the secondary characteristic classes: see \cite[\S 3-4]{FSS12a}.
\end{remark}

In differential enhancement of Example \ref{CohomotopyTheory},
we have:
\begin{example}[Differential Cohomotopy {\cite{FSS15b}}]
  \label{DifferentialCohomotopy}
  The canonical differential enhancement of (unstable)
  Cohomotopy theory (Example \ref{CohomotopyTheory})
  in degree $n$
  is differential non-abelian cohomology (Def. \ref{DifferentialNonAbelianCohomology})
  with coefficients in $S^n$:
  \vspace{-2mm}
  $$
    \overset{
      \mathclap{
      \raisebox{3pt}{
        \tiny
        \color{darkblue}
        \bf
        \def\arraystretch{.9}
        \begin{tabular}{c}
          differential
          \\
          Cohomotopy
        \end{tabular}
      }
      }
    }{
      \widehat \pi^{\, n}(-)
    }
    \;:=\;
    \widehat
    H
    \big(
      -;
      S^n
    \big)
    \,.
  $$

  \vspace{-1mm}
  \noindent {\bf (i)}
  By Example \ref{FlatSphereValuedDifferentialForms},
  a cocycle
  $
    \widehat C_3
      \;\in\;
    \widehat
    \pi^{\, 4}(X)
  $
  in differential 4-Cohomotopy has
  as curvature \eqref{PullbackForModuliInfinityStackOfConnections} a pair
  consisting of a differential 4-form $G_4$ and a differential
  7-form $G_7$, satisfying the
  {\it Cohomotopical Bianchi identity}
  shown here:
  \vspace{-4mm}
  \begin{equation}
    \label{CurvatureInWidehatpi4}
    \raisebox{33pt}{
    \xymatrix@R=-15pt{
      \overset{
        \mathclap{
        \raisebox{3pt}{
          \tiny
          \color{darkblue}
          \bf
          \def\arraystretch{.9}
          \begin{tabular}{c}
            differential
            \\
            4-Cohomotopy
          \end{tabular}
        }
        }
      }{
        \widehat \pi^{\, 4}(X)
      }
      \ar[rr]^-{
        \overset{
          \mathclap{
          \raisebox{3pt}{
            \tiny
            \color{darkblue}
            \bf
            cohomotopical curvature
          }
          }
        }{
          F_{S^4}
        }
      }
      &&
      \Omega
      \big(
        X;
        \,
        \mathfrak{l}S^4
      \big)_{\mathrm{flat}}
      \\
      \underset{
        \mathclap{
        \raisebox{-3pt}{
          \tiny
          \color{greenii}
          \bf
          \def\arraystretch{.9}
          \begin{tabular}{c}
            cohomotopically
            \\
            charge-quantized
            \\
            $C_3$-field
          \end{tabular}
        }
        }
      }{
        \widehat C_3
      }
      \ar@{}[rr]|-{\longmapsto}
      &&
      \left\{
        \!\!\!
        {\begin{array}{c}
          G_7(\widehat C_3),
          \\
          G_4(\widehat C_3)
        \end{array}}
        \!\!\!
        \in
        \Omega^\bullet_{\mathrm{dR}}(X)
        \,\left\vert\,
        {\begin{aligned}
          d\, G_7(\widehat C_3)
            & =
            - G_4(\widehat C_3) \wedge G_4(\widehat C_3)
          \\
          d\, G_4(\widehat C_3) & = 0
        \end{aligned}}
        \right.
      \right\}.
    }
    }
  \end{equation}

  \vspace{-2mm}
\noindent  Such differential form data is exactly what characterizes
  the flux densities of the $C_3$-field in 11-dimensional supergravity
  (up to the self-duality constraint $G_7 = \star G_4$).
  By comparison with Dirac's charge quantization
  \eqref{DiracChargeQuantization}, we thus see that
  a natural candidate for charge quantization of
  the supergravity $C_3$-field is (nonabelian/unstable)
  4-Cohomotopy theory $\pi^4$
  \cite[\S 2.5]{Sati13}\cite[\S 2]{FSS16a}\cite[\S 3]{BMSS19}
  (review in \cite[\S 7]{FSS19a})
  or rather:
  differential 4-Cohomotopy theory $\widehat \pi^{\, 4}$
  \cite[p. 9]{FSS15b}\cite[\S 3.1]{GS-Postnikov}.

  \vspace{0mm}
  \noindent {\bf (ii)}
  The consequence of this Cohomotopical charge quantization
  is readily seen from the
  Hurewicz operation on Cohomotopy theory (Example \ref{HurewiczHomomorphism}):
  The de Rham class of the 4-flux density is constrained to be
  integral, hence to be in the image of the de Rham homomorphism
  (Example \ref{deRhamHomomorphism})
  and its cup square is forced to vanish
  \vspace{-2mm}
  \begin{equation}
    \label{ChargeQuantizationInPlain4Cohomotopy}
    \big[
      G_4(\widehat C_3)
    \big]
    \;\in\;
    \xymatrix{
      H^4
      \big(
        X;
        \,
       \mathbb{Z}
      \big)
      \ar[r]
      &
      H^4_{\mathrm{dR}}
      \big(
        X
      \big)
      \,,
    }
    \phantom{AAAA}
    \big[
      G_4(\widehat C_3)
    \big]
    \cup
    \big[
      G_4(\widehat C_3)
    \big]
    \;=\;
    0
    \,.
  \end{equation}

  \vspace{-2mm}
  \noindent
  This innocent-looking
  but {\it non-linear}
  cup-square relation
  is the source of
  the ``quadratic functions in M-theory'' \cite{HopkinsSinger05},
  revealed here
  as originating from a deep phenomenon in
  unstable, hence ``non-abelian'', homotopy theory,
  revolving around Hopf maps and Massey products
  \cite[\S 4.4]{KS3}\cite{SS21} (see \cite{GS-Massey} for differential refinement).

\vspace{0mm}
\noindent {\bf (iii)}
  Passing from 11-dimensional supergravity to
  M-theory, the curvature data in \eqref{CurvatureInWidehatpi4}
  is expected
  (see \cite[Table 1]{FSS19b})
  to be subjected to more refined topological constraints,
  forcing the class of $G_4$ to be integral
  up to a fractional shift by the first Pontrjagin class
  of the tangent bundle, and deforming its cup square to
  a quadratic function with non-trivial ``background charge''.
  We see, in Prop. \ref{ChargeQuantizationInJTwistedCohomotopy} below,
  that these more subtle M-theoretic
  constraints on the $C_3$-field flux densities
  are, once more, imposed by charge quantization in --
  hence lifting through the non-abelian character map of --
  the corresponding
  \emph{twisted} non-abelian cohomology theory,
  namely: \emph{J-twisted} 4-Cohomotopy \cite{FSS19b}\cite{FSS20a}
  (Example \ref{CharacterMapOnJTwistedCohomotopyAndTwistorialCohomotopy}
  below).
\end{example}

\noindent {\bf Ordinary differential cohomology.}
The \emph{ordinary differential cohomology}
$\widehat H^\bullet(X)$ \cite{SimonsSullivan08}
of a smooth manifold $X$ combines
ordinary integral cohomology classes (Example \ref{OrdinaryCohomology})
with closed differential forms that
represent the same class in real cohomology, in that it makes
a diagram of the following form commute:
\vspace{-3mm}
\begin{equation}
  \label{DifferentialCohomologyDiagram}
  \raisebox{20pt}{
  \xymatrix@R=1em@C=3em{
    \overset{
      \mathllap{
      \raisebox{3pt}{
        \tiny
        \color{darkblue}
        \bf
        \def\arraystretch{.9}
        \begin{tabular}{c}
          ordinary
          \\
          differential cohomology
        \end{tabular}
      }
      }
    }{
      \widehat H^\bullet(X)
    }
    \ar[d]_-{
      \mbox{
        \tiny
        \color{greenii}
        \bf
        \def\arraystretch{.9}
        \begin{tabular}{c}
          underlying
          \\
          integral class
        \end{tabular}
      }
      \!\!\!
    }
    \ar[rr]^-{
      \mbox{
        \tiny
        \color{greenii}
        \bf
        curvature
      }
    }
    &&
    \Omega^\bullet_{\mathrm{dR}}(X)_{\mathrm{clsd}}
    \ar[d]^-{
      \!\!\!
      \mbox{
        \tiny
        \color{greenii}
        \bf
        \def\arraystretch{.9}
        \begin{tabular}{c}
          via
          \\
          de Rham theorem
        \end{tabular}
      }
    }
    \\
    H^\bullet(X;\,\mathbb{Z})
    \ar[rr]^-{
      \mbox{
        \tiny
        \color{greenii}
        \bf
        rationalization
      }
    }
    &&
    H^\bullet(X;\, \mathbb{R})
  }
  }
\end{equation}

\vspace{-2mm}
\noindent
In fact, differential cohomology is universal with this property,
but not at the coarse level of cohomology sets shown above
(where the universal property is shallow) but at the fine level of
complexes of sheaves of coefficients
(i.e. of moduli $\infty$-stacks),
as made precise in Prop.
\ref{DifferentialNonabelianCohomologySubsumesOrdinaryDifferentialCohomology} below
(see Rem. \ref{FromHomotopyPullbackOfOrdinaryDiffCohomologyToCohomology}).

\medskip
In degree 2, ordinary differential cohomology classifies
ordinary $\mathrm{U}(1)$-principal bundles
(equivalently: complex line bundles) with connection \cite[\S II]{Brylinski93},
and the curvature map in \eqref{DifferentialCohomologyDiagram}
assigns their traditional curvature 2-form.
In degree 3, ordinary differential cohomology classifies
bundle gerbes with connection \cite{Murray96}\cite{SW07} with their
curvature 3-form. In general degree, it classifies
higher bundle gerbes with connection \cite{Gajer97},
or equivalently higher $\mathrm{U}(1)$-principal bundles with connection
\cite[2.6]{FSS12b}.

\medskip
One construction of ordinary differential cohomology over smooth manifolds
is given in \cite[\S 1]{CheegerSimons85}, now known as
\emph{Cheeger-Simons characters}. An earlier construction
over schemes,
now known as \emph{Deligne cohomology} (Example \ref{OrdinaryDifferentialGeometry}), due
independently to \cite[\S 2.2]{Deligne71}\cite[\S 3.1.7]{MazurMessing74}\cite[\S III.1]{ArtinMazur77}
and brought to seminal application in \cite{Beilinson85}
(review in \cite{EV88}),
is readily adapted to smooth manifolds
\cite[\S I.5]{Brylinski93}\cite{Gajer97}.
The advantage of Deligne cohomology over Cheeger-Simons characters
is that is immediately
generalizes from smooth manifolds to smooth $\infty$-stacks,
\cite[\S 3.2.3]{FiorenzaSchreiberStasheff10}\cite[\S 2.5]{FSS12b},
such as to orbifolds \cite{SS20a} and to
moduli $\infty$-stacks of higher principal connections
where it yields higher Chern-Simons functionals
\cite{SSS09}\cite{FSS12a}\cite{FSS13a}\cite{FSS15}, as well
as allowing for twists in a systematic manner \cite{GS-Deligne}\cite{GS-HigherDeligne}.

\medskip

In differential enhancement of Example \ref{HigherBundleGerbes}, we have:
\begin{example}[Ordinary differential cohomology on smooth
$\infty$-stacks {\cite[\S 3.2.3]{FiorenzaSchreiberStasheff10}\cite[\S 2.5]{FSS12b}}]
  \label{OrdinaryDifferentialGeometry}
  Let $n \in \mathbb{N}$.

  \noindent
  {\bf (i)}
  The smooth \emph{Deligne-Beilinson complex}
  in degree $n+1$ is the presheaf of connective chain complexes
  (Example \ref{ProjectiveModelStructureOnConnectiveChainComplexes})
  over $\CartesianSpaces$ \eqref{CartesianSpaces}
  given by the truncated and shifted smooth de Rham complex
  (Example \ref{SmoothdeRhamComplex}) with a copy of
  the integers included in degree $n + 1$
  (as integer valued 0-forms, hence as smooth real-valued functions
  constant on an integer):
\vspace{-3mm}
  \begin{equation}
    \label{DeligneComplex}
    \mathrm{DB}^{n+1}_\bullet
    \;:=\;
    \Big(
       \xymatrix{
         \cdots
         \ar[r]
         &
         0
         \ar[r]
         &
         0
         \ar[r]
         &
         \mathbb{Z}
         \;
         \ar@{^{(}->}[r]
         &
         \;
         \Omega^0_{\mathrm{dR}}(-)
         \ar[r]^-{d}
         &
         \Omega^1_{\mathrm{dR}}(-)
         \ar[r]^-{d}
         &
         \cdots
         \ar[r]^-{d}
         &
         \Omega^n_{\mathrm{dR}}(-)
       }
    \big)
    %\;\;\;
    %\in
    %\mathrm{PSh}
    %\Big(
    %  \CartesianSpaces
    %  \,,\,
    %  \ZChainComplexes
    %\big)
    \,.
  \end{equation}

\vspace{-2mm}
  \noindent
  {\bf (ii)}
  The de Rham differential in degree 0 gives a morphism of
  presheaves of complexes
   \vspace{-2mm}
  \begin{equation}
    \label{DeRhamDifferentialOnDeligneComplex}
    \xymatrix{
      \mathrm{DB}^{n+1}_\bullet
      \ar[rr]^-{
        (0,0,\cdots,0,d)
      }
      &&
      \Omega^{n+1}_{\mathrm{dR}}(-)_{\mathrm{clsd}}
    }
  \end{equation}

  \vspace{-3mm}
\noindent
 from the Deligne-Beilinson complex \eqref{DeligneComplex}
  to the presheaf of closed $(n+1)$-forms, regarded
  as a presheaf of chain complexes  in degree 0.

  \noindent
  {\bf (iii)}
  \emph{Ordinary differential cohomology} is sheaf hypercohomology with coefficients in the Deligne complex. This means that
  if we look at the Deligne-Beilinson complex
  \eqref{DeligneComplex}
   as a smooth $\infty$-stack (Def. \ref{SmoothInfinityStacks})
  by
  first applying the Dold-Kan construction
  from Example \ref{DoldKanConstruction} and then $\infty$-stackifying the resulting simplicial presheaf, then ordinary differential cohomology
  is stacky non-abelian cohomology (Remark \ref{StructuredNonAbelianCohomology})
  with coefficients in the Deligne-Beilinson complex:
  \vspace{-2mm}
  \begin{equation}
    \label{OrdinaryDifferentialCohomologyAsSheafHypercohomology}
    \overset{
      \mathclap{
      \raisebox{3pt}{
        \tiny
        \color{darkblue}
        \bf
        \def\arraystretch{.9}
        \begin{tabular}{c}
          ordinary
          \\
          differential cohomology
        \end{tabular}
      }
      }
    }{
    \widehat H^{n+1}
    \big(
      \mathcal{X}
    \big)
    }
    \;\;:=\;\;
    \Stacks
    \Big(
      \mathcal{X}
      \;\;,\;\;
      \overset{
        \mathclap{
        \raisebox{3pt}{
          \tiny
          \color{darkblue}
          \bf
          \def\arraystretch{.9}
          \begin{tabular}{c}
            Dold-Kan
            \\
            correspondence
          \end{tabular}
        }
        \;\;
        }
      }{
        L^{\mathrm{loc}}
        \,\circ\,
        \mathrm{DK}
      }
      \big(
        \overset{
          \mathclap{
          \;\;
          \raisebox{6pt}{
          \tiny
          \color{darkblue}
          \bf
          \def\arraystretch{.9}
          \begin{tabular}{c}
            Deligne-Beilinson
            \\
            complex
          \end{tabular}
          }
          }
        }{
          \mathrm{DB}^{n+1}_\bullet
        }
      \big)
    \Big)
    \,.
  \end{equation}

\vspace{-3mm}
  \noindent
  {\bf (iv)} The {\it curvature map} on ordinary differential cohomology
  is the cohomology operation induced by \eqref{DeRhamDifferentialOnDeligneComplex}:
  \vspace{-2mm}
  \begin{equation}
    \label{CurvatureMapOnOrdinaryDifferentialCohomology}
    \hspace{-4mm}
   {\small
    \xymatrix@C=4.5em@R=.1em{
    \overset{
      \mathclap{
      \raisebox{3pt}{
        \tiny
        \color{darkblue}
        \bf
        \def\arraystretch{.9}
        \begin{tabular}{c}
          ordinary
          \\
          differential cohomology
        \end{tabular}
      }
      }
    }{
    \widehat H^{n+1}
    \big(
      \mathcal{X}
    \big)
    }
    \ar@{=}[d]
    \ar[rr]_-{ F }^-{
      \mbox{
        \tiny
        \color{darkblue}
        \bf
        curvature
      }
    }
    &&
    \Omega^{n+1}_{\mathrm{dR}}
    \big(
      \mathcal{X}
    \big)_{\mathrm{clsd}}
    \ar@{=}[d]
    \\
    \Stacks
    \Big(\!\!
      \mathcal{X}
      ,
      L_{\mathrm{loc}}
      \circ
      \mathrm{DK}
      \big(
        \mathrm{DB}^{n+1}_\bullet
      \big)
   \! \Big)
    \ar[rr]^-{
\scalebox{0.6}{$\Stacks
    \left(
      \mathcal{X}
      ,
      L_{\mathrm{loc}}
      \circ
      \mathrm{DK}
      (
        d
      )
    \right)
    $}
    }
    &&
    \Stacks
    \Big(\!\!
      \mathcal{X}
      ,
      L_{\mathrm{loc}}
      \circ
      \mathrm{DK}
      \big(
        \Omega^{n+1}_{\mathrm{dR}}(-)_{\mathrm{clsd}}
      \big)
    \!\Big).
    }
    }
  \end{equation}

\end{example}

\begin{prop}[Differential non-abelian cohomology subsumes differential ordinary cohomology {\cite[Prop. 3.2.26]{FiorenzaSchreiberStasheff10}}]
  \label{DifferentialNonabelianCohomologySubsumesOrdinaryDifferentialCohomology}

  Let $n \in \mathbb{N}$ and consider
   $A \,=\, B^n \mathrm{U}(1) \simeq K(\mathbb{Z}, n+1)$
  (Example \ref{HigherBundleGerbes}). Then:

  \noindent
  {\bf (i)}
  Differential non-abelian $A$-cohomology (Def. \ref{DifferentialNonAbelianCohomology})
  coincides with ordinary differential cohomology
  (Def. \ref{OrdinaryDifferentialGeometry}):

  \vspace{-2mm}
  \begin{equation}
    \label{RecoveringOrdinaryDifferentialCohomology}
    \overset{
      \mathclap{
      \raisebox{3pt}{
        \tiny
        \color{darkblue}
        \bf
        \def\arraystretch{.9}
        \begin{tabular}{c}
          ordinary
          \\
          differential cohomology
        \end{tabular}
      }
      }
    }{
    \widehat H^{n+1}
    \big(
      \mathcal{X}
    \big)
    }
    \;\;
    \simeq
    \;\;
    \widehat H
    \big(
      \mathcal{X};
      \,
      B^n \mathrm{U}(1)
    \big)
    \,.
  \end{equation}

  \noindent
  {\bf (ii)}
  The abstract curvature map in
  differential $A$-cohomology
  \eqref{PullbackForModuliInfinityStackOfConnections}
  reproduces the ordinary curvature map \eqref{CurvatureMapOnOrdinaryDifferentialCohomology}.
\end{prop}
\begin{proof}
 % {\bf (i)}
  In order to compute the defining
  homotopy pullback \eqref{PullbackForModuliInfinityStackOfConnections},
  we use the Dold-Kan correspondence (Prop. \ref{DoldKanCorrespondence})
  to obtain a convenient presentation
  of the differential character map along which to pull back:

  {\bf (a)}
  Since the Dold-Kan construction $\mathrm{DK}$ (Def. \ref{DoldKanConstruction})
  realizes homotopy groups from homology groups \eqref{DoldKanHomotopyGroups},
  and since Eilenberg-MacLane spaces
  are characterized by their homotopy groups \eqref{EilenbergMacLaneSpaces},
we  have the vertical identifications on the left of the following diagram:

  \vspace{-.4cm}
  \begin{equation}
    \label{PresentingTheDifferentialCharacterForOrdinaryCohomology}
    {\footnotesize
    \raisebox{40pt}{
    \xymatrix@C=3em@R=1.5em{
      \mathrm{Disc}
      \big(
        B^{n+1}\mathbb{Z}
      \big)
      \ar@/^1.3pc/[rrrr]^-{ \mathbf{ch}_{B^n \mathrm{U}(1)} }
      \ar[rr]_-{ \eta^{\mathbb{R}}_{B^{n+1}\mathbb{Z}} }
      \ar@{=}[d]
      &&
      \mathrm{Disc}
      \big(
        B^{n+1}\mathbb{R}
      \big)
      \ar@{=}[d]
      \ar[rr]_-{\simeq}
      &&
      \flatBexp( \mathfrak{b}^n \mathbb{R} )
      \ar[d]^-{\simeq}_-{\int_{\Delta^\bullet}}
      \\
      \mathllap{\mathrm{DK}}
      \left(
      {\begin{array}{c}
        \mathbb{Z}
        \\[-3pt]
        \downarrow
        \\[-3pt]
        0
        \\[-3pt]
        \downarrow
        \\[-3pt]
        \vdots
        \\[-3pt]
        \downarrow
        \\[-3pt]
        0
      \end{array}}
      \right)\;
      \ar@{^{(}->}[rr]
      &&
      \;
      \!\!\!\!
      \mathrm{DK}
      \left(
      {\begin{array}{c}
        \mathbb{R}
        \\[-3pt]
        \downarrow
        \\[-3pt]
        0
        \\[-3pt]
        \downarrow
        \\[-3pt]
        \vdots
        \\[-3pt]
        \downarrow
        \\[-3pt]
        0
      \end{array}}
      \right)
      \;
      \ar@{^{(}->}[rr]
      &&
      \;
      \mathrm{DK}
      \left(\!\!
      {\begin{array}{c}
        \Omega^{0}_{\mathrm{dR}}(-)
        \\[-3pt]
        \downarrow\mathrlap{\scalebox{.7}{$d$}}
        \\[-3pt]
        \Omega^{1}_{\mathrm{dR}}(-)
        \\[-3pt]
        \downarrow\mathrlap{\scalebox{.7}{$d$}}
        \\[-3pt]
        \vdots
        \\[-3pt]
        \downarrow\mathrlap{\scalebox{.7}{$d$}}
        \\[-3pt]
        \Omega^{n+1}_{\mathrm{dR}}(-)_{\mathrm{clsd}}
      \end{array}}
     \!\! \right)
    }
    }
    }
  \end{equation}

    \vspace{-2mm}
\noindent  Under this identification, it is clear that the rationalization map
  $\eta^{\mathbb{R}}_{B^{n+1} \mathbb{Z}}$
  (Def. \ref{Rationalization})
  is presented by the
  canonical inclusion of the integers into the real numbers,
  as on the bottom left of \eqref{PresentingTheDifferentialCharacterForOrdinaryCohomology}.

  Moreover, the right vertical equivalence in \eqref{PresentingTheDifferentialCharacterForOrdinaryCohomology}
  is that from Lemma \ref{ModuliInfinityStackOfClosedFormsIsShiftedDeRhamComplex}.

  {\bf (b)} Since the differential character \eqref{DifferentialNonabelianCharacterMapPresented}
  in the present case evidently comes from a morphism of
  (presheaves of) simplicial abelian groups,
  with group structure given by
  addition of ordinary differential forms (Example \ref{OrdinaryClosedFormsAreFlatLineLInfinityAlgebraValuedForms}),
  we may, using the Dold-Kan correspondence (Prop. \ref{DoldKanCorrespondence}),
  analyze the remainder of the diagram
  on normalized chain complexes $N(-)$ \eqref{DoldKanEquivalenceOfCategories}.

  Using this, it follows by inspection of the bottom
  map in \eqref{DifferentialNonabelianCharacterMapPresented}
  that the bottom right square in \eqref{PresentingTheDifferentialCharacterForOrdinaryCohomology}
  commutes, with the bottom morphism on the right being the
  canonical inclusion of (presheaves of) chain complexes.

  Now to use this presentation for identifying the
  resulting homotopy fiber product \eqref{PullbackForModuliInfinityStackOfConnections}:

  \vspace{0mm}
  \noindent
  {\bf (i)} Since the $\mathrm{DK}$-construction (Def. \ref{DoldKanConstruction}),
  applied objectwise over $\CartesianSpaces$,
  is a right Quillen functor into the global model structure
  from Example \ref{ModelStructureOnSimplicialPresheavesOverCartesianSpaces},
  and since $\infty$-stackification preserves homotopy pullbacks
  (Lemma \ref{InfinityStackification}), it is now sufficient
  to show, by definition \eqref{OrdinaryDifferentialCohomologyAsSheafHypercohomology},
  that the homotopy pullback (Def. \ref{HomotopyPullback})
  along the bottom map in
  \eqref{PresentingTheDifferentialCharacterForOrdinaryCohomology},
  formed in presheaves of chain complexes is
  the Deligne-Beilinson complex $\mathrm{DB}^{n+1}_\bullet$ \eqref{DeligneComplex}.

  \noindent
  \begin{minipage}[left]{5.3cm}
  For this it is sufficient,
  by \eqref{ConstructionOfHomotopyPullback},
  to find a fibration replacement of
  the bottom map in \eqref{PresentingTheDifferentialCharacterForOrdinaryCohomology}
  whose ordinary fiber product with $\Omega^{n+1}_{\mathrm{dR}}(-)_{\mathrm{clsd}}$
  is the Deligne-Beilinson complex.
  This is provided by a mapping cylinder construction
  (e.g. \cite[\S 1.5.5]{Weibel94}) shown here:
  \end{minipage}

  \vspace{-3cm}

  \begin{equation}
  \label{TowardsRecoveringDeligneComplex}
  {\scriptsize
    \raisebox{40pt}{
    \xymatrix@C=5em@R=2.1em{
      &&
      \mathllap{
        \mathrm{DB}^{n+1}_\bullet
        \;
        =
      }
      \left(
      {\arraycolsep=0pt
      {\begin{array}{c}
        \mathbb{Z}
        \\
        \downarrow\mathrlap{\!\!\raisebox{1pt}{\scalebox{.7}{$i$}}}
        \\
        \Omega^0_{\mathrm{dR}}(-)
        \\
        \downarrow\mathrlap{\!\!\raisebox{1pt}{\scalebox{.7}{$d$}}}
        \\
        \Omega^1_{\mathrm{dR}}(-)
        \\
        \downarrow\mathrlap{\!\!\raisebox{1pt}{\scalebox{.7}{$d$}}}
        \\
        \vdots
        \\
        \downarrow\mathrlap{\!\!\raisebox{1pt}{\scalebox{.7}{$d$}}}
        \\
        \Omega^{n-1}_{\mathrm{dR}}(-)
        \\
        \downarrow\mathrlap{\!\!\raisebox{1pt}{\scalebox{.7}{$d$}}}
        \\
        \Omega^n_{\mathrm{dR}}(-)
      \end{array}}}
      \right)
      \ar[d]_-{ i_1 }
      \ar[rr]^-{
        \left(
        \scalebox{.7}{$
        {\begin{array}{c}
          0
          \\
          0
          \\
          \vdots
          \\
          0
          \\
          d
        \end{array}}
        $}
        \right)
      }
      \ar@{}[drr]|-{
        \mbox{\tiny\rm(pb)}
      }
      &&
      \left(
      {\begin{array}{ccc}
        0
        \\
        \downarrow
        \\
        0
        \\
        \downarrow
        \\
        0
        \\
        \downarrow
        \\
        \vdots
        \\
        \downarrow
        \\
        0
        \\
        \downarrow
        \\
        \Omega^{n+1}_{\mathrm{dR}}(-)_{\mathrlap{\mathrm{clsd}}}
      \end{array}}
      \;\;\;\;
      \right)
      \ar[d]^-{i}
      \\
      \left(
      {\begin{array}{ccc}
        \mathbb{Z}
        \\
        \downarrow
        \\
        0
        \\
        \downarrow
        \\
        0
        \\
        \downarrow
        \\
        \vdots
        \\
        \downarrow
        \\
        0
        \\
        \downarrow
        \\
        0
      \end{array}}
      \right)
      \ar[rr]^-{ n \,\mapsto\, (n,n) }_-{ \in \, \mathrm{W} }
      &&
      \left(
      {\arraycolsep=0pt
      {\begin{array}{ccc}
        \mathbb{Z}
        &
        \oplus
        &
        \,
        \Omega^0_{\mathrm{dR}}(-)
        \\
        \downarrow\mathrlap{\!\!\raisebox{1pt}{\scalebox{.7}{$i$}}}
        & \swarrow\mathrlap{\!\!\!\!\!\raisebox{0pt}{\scalebox{.7}{$-\mathrm{id}$}}}
        &
        \downarrow\mathrlap{\!\!\raisebox{1pt}{\scalebox{.7}{$d$}}}
        \\
        \Omega^0_{\mathrm{dR}}(-)
        &\oplus&
        \,
        \Omega^1_{\mathrm{dR}}(-)
        \\
        \downarrow\mathrlap{\!\!\raisebox{1pt}{\scalebox{.7}{$d$}}}
        &
        \swarrow\mathrlap{\!\!\!\!\!\raisebox{0pt}{\scalebox{.7}{$+\mathrm{id}$}}}
        &
        \downarrow\mathrlap{\!\!\raisebox{1pt}{\scalebox{.7}{$d$}}}
        \\
        \Omega^1_{\mathrm{dR}}(-)
        &\oplus&
        \,
        \Omega^2_{\mathrm{dR}}(-)
        \\
        \downarrow\mathrlap{\!\!\raisebox{1pt}{\scalebox{.7}{$d$}}}
        &
        \swarrow\mathrlap{\!\!\!\!\!\raisebox{0pt}{\scalebox{.7}{$-\mathrm{id}$}}}
        &
        \downarrow\mathrlap{\!\!\raisebox{1pt}{\scalebox{.7}{$d$}}}
        \\
        \vdots & \vdots & \vdots
        \\
        \downarrow\mathrlap{\!\!\raisebox{1pt}{\scalebox{.7}{$d$}}}
        & \swarrow
        &
        \downarrow\mathrlap{\!\!\raisebox{1pt}{\scalebox{.7}{$d$}}}
        \\
        \Omega^{n-1}_{\mathrm{dR}}(-)
        &\oplus&
        \,
        \Omega^n_{\mathrm{dR}}(-)
        \\
        \downarrow\mathrlap{\!\!\raisebox{1pt}{\scalebox{.7}{$d$}}}
        &
        \swarrow
        &
        %\downarrow\mathrlap{\!\!\raisebox{1pt}{\scalebox{.7}{$d$}}}
        \\
        \Omega^n_{\mathrm{dR}}(-)
        &
        &
      \end{array}}}
      \right)
      \ar[rr]^-{
        \left(
        \scalebox{.7}{$
        {\begin{array}{c}
          \mathrm{pr}_2
          \\
          \mathrm{pr}_2
          \\
          \vdots
          \\
          \mathrm{pr}_2
          \\
          d
        \end{array}}
        $}
        \right)
      }_-{ \in \; \mathrm{Fib} }
      &&
      \left(
      {\begin{array}{ccc}
        \Omega^0_{\mathrm{dR}}(-)
        \\
        \downarrow\mathrlap{\!\!\raisebox{1pt}{\scalebox{.7}{$d$}}}
        \\
        \Omega^1_{\mathrm{dR}}(-)
        \\
        \downarrow\mathrlap{\!\!\raisebox{1pt}{\scalebox{.7}{$d$}}}
        \\
        \Omega^2_{\mathrm{dR}}(-)
        \\
        \downarrow\mathrlap{\!\!\raisebox{1pt}{\scalebox{.7}{$d$}}}
        \\
        \vdots
        \\
        \downarrow\mathrlap{\!\!\raisebox{1pt}{\scalebox{.7}{$d$}}}
        \\
        \Omega^n_{\mathrm{dR}}(-)
        \\
        \downarrow\mathrlap{\!\!\raisebox{1pt}{\scalebox{.7}{$d$}}}
        \\
        \Omega^{n+1}_{\mathrm{dR}}(-)_{\mathrlap{\mathrm{clsd}}}
      \end{array}}
      \;\;\;\;
      \right)
    }
    }
    }
  \end{equation}

  \noindent
  By direct inspection, we see in this diagram that:

  \begin{itemize}

  \vspace{-2mm}
  \item
  the total bottom morphism is the total bottom morphism
  from \eqref{PresentingTheDifferentialCharacterForOrdinaryCohomology},
  factored as a weak equivalence (quasi-isomorphism)
  followed by a fibration (positive degreewise surjection);

  \vspace{-2mm}
  \item
  the ordinary pullback of this fibration is
  the Deligne-Beilinson complex
  $\mathrm{DB}^{n+1}_\bullet$ \eqref{DeligneComplex},
  as shown,
  which therefore
  represents the homotopy pullback
  (since all chain complexes are projectively fibrant),
  by Def.  \ref{HomotopyPullback}.

  \vspace{-2mm}
  \item
  the top morphism out of this (homotopy-)pullback
  coincides with the curvature map \eqref{DeRhamDifferentialOnDeligneComplex}
  on the Deligne complex -- which, under the following implication of
  claim (i), implies claim (ii).

  \end{itemize}
  \vspace{-.2cm}

  \noindent
  {\bf (ii)}
  The image of this homotopy pullback \eqref{TowardsRecoveringDeligneComplex}
  under
  $L^{\mathrm{loc}} \circ \mathrm{DK}$
  is still a homotopy pullback
  (because $\mathrm{DK}$ is a right Quillen functor by
  construction \eqref{DK} and using Lem. \ref{InfinityStackification})
  and hence exhibits the Deligne coefficients \eqref{OrdinaryDifferentialCohomologyAsSheafHypercohomology}
  for ordinary differential cohomology as a model for the
  differential $B^{n+1}\mathbb{Z}$-cohomology according to Def. \ref{DifferentialNonAbelianCohomology}:

  \vspace{-.4cm}
  \begin{equation}
    \label{DeligneCoefficientsAsAbstractHomotopyPullback}
    \begin{tikzcd}[column sep={between origins, 90pt}]
      \overset{
        \mathclap{
        \raisebox{3pt}{
          \tiny
          \color{darkblue}
          \bf
          \def\arraystretch{.9}
          \begin{tabular}{c}
            Deligne complex
            as smooth $\infty$-stack
          \end{tabular}
        }
        }
      }{
        L^{\mathrm{loc}}
          \circ
        \mathrm{DK}
        \big(
          \mathrm{DB}^{n+1}_\bullet
        \big)
      }
      \ar[
        rr,
        "F_{B^{n+1}\mathbb{Z}}"
      ]
      \ar[
        d,
        "{ c_{B^{n+1}\mathbb{Z}} }"
      ]
      \ar[
        drr,
        phantom,
        "\mbox{\tiny\rm (hpb)}"
      ]
      &&
      \Omega_{\mathrm{dR}}
      \big(
        -;
        \,
        \mathfrak{l}B^{n+1}\mathbb{Z}
      \big)_{\mathrm{flat}}
      \ar[d]
      \\
      \mathrm{Disc}
      \big(
        B^{n+1}\mathbb{Z}
      \big)
      \ar[
        rr,
        "{ \mathbf{ch}_{B^{n+1}\mathbb{Z}} }"{below}
      ]
      &&
      \flatBexp
      \big(
        \mathfrak{l}B^{n+1}\mathbb{Z}
      \big)
    \end{tikzcd}
  \end{equation}
  This implies claim (i), by the definitions.
\end{proof}
\begin{remark}[The commuting square of ordinary (differential) cohomology groups]
  \label{FromHomotopyPullbackOfOrdinaryDiffCohomologyToCohomology}
  The image of the homotopy-pullback square
  \eqref{DeligneCoefficientsAsAbstractHomotopyPullback}
  under the hom-functor $\Stacks(X,-)$
  out of a smooth manifold $X$ gives the
  commuting square of
  ordinary (differential) cohomology groups shown in \eqref{DifferentialCohomologyDiagram}.
  Since the hom-functor of a homotopy category
  does not preserve homotopy pullbacks, in general
  (only the mapping space functor \eqref{MappingSpaceBetweenSmoothStacks} does),
  the square \eqref{DifferentialCohomologyDiagram} in cohomology
  is not itself a pullback, in general.
\end{remark}

%In differential enhancement of Example %\ref{TraditionalNonAbelianCohomology}
%we have
%
%\begin{example}[Weak principal connections]
%  For $G$ a compact connected Lie group and
%  $B G$ denoting its classifying space, the
%  differential cohomology enhancement
%  produced by Def. \ref{DifferentialNonAbelianCohomology}
%  is
%\end{example}

%\newpage

\medskip
\noindent {\bf Secondary non-abelian cohomology operations.}
We define secondary non-abelian cohomology operations
(Def. \ref{SecondaryNonabelianCohomologyOperations} below)
which generalize the classical notion of secondary characteristic classes
(Theorem \ref{SecondaryDifferentialNonAbelianCharacterSubsumesCheegerSimonsHomomorphism},
see Remark \ref{SecondaryInvariantsOfGConnections} for the terminology)
to higher non-abelian cohomology.
To formulate the concept in this generality, we need a technical condition
(Def. \ref{RelativeFormalMaps}) which happens to be trivially
satisfied in the classical case (Lemma \ref{CharacteristicClassesOfGPrincipalBundlesAreRelativeFormal} below):

\begin{defn}[Absolute minimal model]
  \label{RelativeFormalMaps}
  For $A_1, A_2 \in \NilpotentConnectedQFiniteHomotopyTypes$
  (Def. \ref{NilpotentConnectedSpacesOfFiniteRationalType})
  we say that an \emph{absolute minimal model}
  for a morphism
  $\xymatrix@C=20pt{A_1 \ar[r]|{\;c\;} & A_2}$
  in $\SimplicialSets$
  is a morphism
  $\xymatrix@C=20pt{ \mathfrak{l}A_1 \ar[r]|-{\;\mathfrak{c}\,} & \mathfrak{l} A_2 }$
  between the respective Whitehead $L_\infty$-algebras (Prop. \ref{WhiteheadLInfinityAlgebras})
  which makes the square on the left
  and hence the square on the far right of the following diagram commute:
  \vspace{-3mm}
  \begin{equation}
    \label{TransformationBetweenMinimalDerivedPLdRAdjunctionUnits}
    \hspace{-6mm}
    \raisebox{25pt}{
    \xymatrix@C=2em@R=2em{
      \Omega^\bullet_{\mathrm{dRPL}}(A_1)
      \ar@{<-}[r]^-{ p^{\mathrm{min}}_{A_1} }
      \ar@{<-}[d]
      &
      \mathrm{CE}\big(\mathfrak{l}A_1\big)
      \ar@{<--}[d]^-{ \mathfrak{c} }
      \\
      \Omega^\bullet_{\mathrm{dRPL}}(A_2)
      \ar@{<-}[r]_-{ p^{\mathrm{min}}_{A_2} }
      &
      \mathrm{CE}\big(\mathfrak{l}A_2\big)
      \,,
      \\
      \ar@{}[r]|-{
        \mbox{$\in \; \dgcAlgebras{\mathbb{R}}$}
      }
      &
    }
    }
    \phantom{AA}
    \raisebox{28pt}{
    \xymatrix@C=30pt@R=3em{
      A_1
      \ar[d]_-{c}
      \ar[rr]|-{ \; \eta^{\mathrm{PLdR}}_{A_1} }
      \ar@/^1.3pc/[rrrr]^-{\scalebox{0.7}{$ \mathbb{D} \eta^{\mathrm{PLdR}}_{A_1}$}}
      &&
      \Bexp_{\mathrm{PL}} \,\circ\, \Omega^\bullet_{\mathrm{PLdR}}(A_1)
      \ar[d]^-{\scalebox{0.7}{$
        \Bexp_{\mathrm{PL}} \,\circ\, \Omega^\bullet_{\mathrm{PLdR}}(c)
        $}
      }
      \ar[rr]|-{ \; \scalebox{0.7}{$
        \Bexp_{\mathrm{PL}}(p^{\mathrm{min}}_{A_1})
        $}
        \;
      }
      &&
      \Bexp_{\mathrm{PL}} \,\circ\, \mathrm{CE}(\mathfrak{l}A_1)
      \ar@{-->}[d]^-{
        \mathclap{\phantom{\vert}}
        \scalebox{0.7}{$
        \Bexp_{\mathrm{PL}} \,\circ\, \mathrm{CE}( \mathfrak{c} )
        $}
      }
      \\
      A_2
      \ar[rr]|-{ \; \scalebox{0.7}{$\eta^{\mathrm{PLdR}}_{A_2}$} }
      \ar@/_1.3pc/[rrrr]_-{\scalebox{0.7}{$ \mathbb{D} \eta^{\mathrm{PLdR}}_{A_2}$} }
      &&
      \Bexp_{\mathrm{PL}} \,\circ\, \Omega^\bullet_{\mathrm{PLdR}}(A_2)
      \ar[rr]|-{\;
       \scalebox{0.7}{$\Bexp_{\mathrm{PL}}(p^{\mathrm{min}}_{A_2})$}
      \;}
      &&
      \Bexp_{\mathrm{PL}} \,\circ\, \mathrm{CE}(\mathfrak{l}A_2)
      \,,
      \\
      \ar@{}[rrrr]|-{
        \mbox{$\phantom{AAAAAAAAAAA}
        \in
        \;
        \SimplicialSets
        $}
      }
      &&&&
    }
    }
  \end{equation}

      \vspace{-1mm}
\noindent
  hence a morphism that yields a transformation between exactly those
  derived adjunction units
  $\mathbb{D}\eta^{\mathrm{PLdR}}$ \eqref{DerivedAdjunctionUnit}
  of the PL-de Rham adjunction \eqref{QuillenAdjunctBetweendgcAlgsAndSimplicialSets}
  that are given by \emph{minimal} fibrant replacement.\footnote{
   Notice that the existence of morphisms $\mathfrak{c}$
   making this diagram commute
   is not guaranteed; it is only the existence of the
   \emph{relative} minimal morphism
   $\mathfrak{l}_{\scalebox{.7}{$A_2$}}(c)$
   from Prop. \ref{WhiteheadLInfinityAlgebrasRelative}
   which is guaranteed to make the square \eqref{CEAlgebraOfWhiteheadLInfinityAlgebraRelative}
   commute.
  }
In this case,
the commuting diagram \eqref{TransformationBetweenMinimalDerivedPLdRAdjunctionUnits}
evidently extends to a strict transformation between the
differential non-abelian characters
\eqref{DifferentialNonabelianCharacterMapPresented}
on the $A_i$
(Def. \ref{DifferentialNonabelianCharacterMap}), in that the
following diagram of simplicial presheaves (Def. \ref{SimplicialPresheavesOverCartesianSpaces}) commutes:
\vspace{-2mm}
\begin{equation}
  \label{TransformationOfDifferentialNonaebalianCharacters}
  \raisebox{22pt}{
  \xymatrix@C=3em{
    \mathrm{Disc}(A_1)
    \ar[d]_-{
      \mathrm{Disc}(c)
    }
    \ar[rr]^-{
      \mathbf{ch}_{A_1}
    }
    &&
    \flatBexp(\mathfrak{l}A_1)
    \ar[d]_-{
      \flatBexp(\mathfrak{c})
    }
    \\
    \mathrm{Disc}(A_2)
    \ar[rr]^-{
      \mathbf{ch}_{A_2}
    }
    &&
    \flatBexp(\mathfrak{l}A_2)
  }
  }
  \;\;\;\;\;\;\;\;
  \in
  \;
  \mathrm{PSh}
  \big(
    \CartesianSpaces
    \,,\,
    \SimplicialSets
  \big).
\end{equation}
\end{defn}

In differential enhacement of Def. \ref{NonAbelianCohomologyOperations},
we have:
\begin{defn}[Secondary non-abelian cohomology operation]
  \label{SecondaryNonabelianCohomologyOperations}
  Let
  $\xymatrix@C=12pt{A_1 \ar[r]^{c} & A_2}$
  in $\SimplicialSets$, with induced cohomology operation
  (Def. \ref{NonAbelianCohomologyOperations})
  \vspace{-2mm}
  $$
    \xymatrix{
      H
      (
        -;
        A_1
      )
      \ar[r]^-{ c_\ast }
      &
      H
      (
        -;
        A_2
      )
      \,,
    }
  $$

    \vspace{-1mm}
  \noindent  have an absolute minimal model $\mathfrak{c}$
  (Def \ref{RelativeFormalMaps}). Then the corresponding
  {\it secondary non-abelian cohomology operation}
  is the structured
  cohomology operation (Remark \ref{StructuredNonAbelianCohomology})
    \vspace{-2mm}
  \begin{equation}
    \label{SecondaryOperationInduced}
    \xymatrix@C=3em{
      \widehat
      H
      (
        -;
        \,
        A_1
      )
      \ar[rr]^-{
        (c_{\mathrm{diff}})_\ast
      }_-{
        \mbox{
          \tiny
          \color{darkblue}
          \bf
          {
          \def\arraystretch{.9}
          \begin{tabular}{c}
            secondary
            \\
            non-abelian character
          \end{tabular}}
        }
      }
      &&
      \widehat
      H
      (
        -;
        \,
        A_2
      )
    }
  \end{equation}

    \vspace{-2mm}
\noindent
  on differential non-abelian cohomology
  (Def. \ref{DifferentialNonAbelianCohomology})
  which is induced, as in \eqref{InducedStructuredCohomologyOperation},
  by the dashed morphism $c_{\mathrm{diff}}$ in the following diagram,
  which in turn is induced from $c$ and $\mathfrak{c}$
  \eqref{TransformationOfDifferentialNonaebalianCharacters}
  by the universal property of the
  defining homotopy pullback operation
  \eqref{DifferentialNonabelianCharacterMapPresented}:

  \vspace{-4mm}
  \begin{equation}
    \label{SecondaryCohomologyOperation}
    \xymatrix@R=1.7em{
      \mathllap{
        \mbox{
          \tiny
          \color{darkblue}
          \bf
          \def\arraystretch{.9}
          \begin{tabular}{c}
            secondary/differential
            \\
            cohomology operation
          \end{tabular}
        }
        \;\;\;
      }
      (A_1)_{\mathrm{diff}}
      \ar@{-->}[rr]^-{ c_{\mathrm{diff}} }
      \ar[dd]_<<<<<<<{ c_{A_1} }
      \ar[dr]_-{ F_{A_1} }
      &&
      (A_2)_{\mathrm{diff}}
      \ar[dd]^<<<<<<<{ c_{A_2} }|-{ \phantom{\vert} }
      \ar[dr]^-{ F_{A_2} }
      \\
      &
      \Omega_{\mathrm{dR}}
      \big(
        -;
        \mathfrak{l}A_1
      \big)_{\mathrm{flat}}
      \ar[dd]
      \ar[rr]^<<<<<<<{ \mathfrak{c}_\ast }
      &&
      \Omega_{\mathrm{dR}}
      \big(
        -;
        \mathfrak{l}A_2
      \big)_{\mathrm{flat}}
      \ar[dd]
      \\
      \mathllap{
        \mbox{
          \tiny
          \color{darkblue}
          \bf
          \def\arraystretch{.9}
          \begin{tabular}{c}
            plain/primary
            \\
            cohomology operation
          \end{tabular}
        }
        \;\;\;
      }
      \mathrm{Disc}(A_1)
      \ar[rr]^<<<<<<<<<<<{
        \mathrm{Disc}(c)
      }|-{ \phantom{AA} }
      \ar[dr]_-{
      \mathllap{
        \mbox{
          \tiny
          \color{darkblue}
          \bf
          \def\arraystretch{.9}
          \begin{tabular}{c}
            transformation of
            \\
            differential
            characters
          \end{tabular}
        }
        \;\;\;
      }
        \mathbf{ch}_{A_1}
      }
      &&
      \mathrm{Disc}(A_2)
      \ar[dr]^-{ \mathbf{ch}_{A_2} }
      \\
      &
      \flatBexp(\mathfrak{l}A_1)
      \ar[rr]^-{ \mathfrak{c}_\ast }
      &&
      \flatBexp(\mathfrak{l}A_2)
      \,.
    }
  \end{equation}

  \vspace{-1mm}
\noindent  The left and right squares are the
  homotopy pullback squares defining differential non-abelian
  cohomology (Def. \ref{DifferentialNonAbelianCohomology})
  while the bottom square is the transformation of
  differential non-abelian characters (Def. \ref{DifferentialNonabelianCharacterMap})
  from \eqref{TransformationOfDifferentialNonaebalianCharacters}.
\end{defn}

In differential enhancement of
Examples \ref{HurewiczHomomorphism}, \ref{TheBoardmanHomomorphismIntmf}
we have:
\begin{example}[Secondary non-abelian Hurewicz/Boardman homomorphism to differential K-theory]
  \label{SecondaryNonAbelianBoardmanHomomorphismToKTheory}
  Consider the map

  \vspace{-3mm}
  $$
    \xymatrix{
      S^4
      \ar[rr]^-{ e_{B\mathrm{U}}^4 }
      &&
      B \mathrm{U}
      \;\;\;\;\;\;\;\;
      \in
      \;
      \HomotopyTypes
    }
  $$

    \vspace{-2mm}
\noindent  from the 4-sphere to the classifying space of the infinite
  unitary group \eqref{ClassifyingSpaceOfInfiniteUnitaryGroup}
  which classifies a generator in
  $\pi_4\big( B \mathrm{U} \big) \,\simeq\, \mathbb{Z} $.
  By Example \ref{RationalizationOfnSpheres}
  and
  Examples \ref{WhiteheadLInfinityAlgebraOfLoopSpaces},
  \ref{RecoveringH3TwistedDifferentialForms}
  the corresponding Whitehead $L_\infty$-algebras
  (Prop. \ref{WhiteheadLInfinityAlgebras})
  are as shown here:
  \begin{equation}
    \label{AbsoluteMinimalModelForBeta4ToKU0}
    \hspace{-1.6cm}
    \raisebox{50pt}{
    \xymatrix@C=40pt@R=12pt{
      \mathrm{CE}
      \big(
        \mathfrak{l}S^4
      \big)
      \ar@{<-}[rr]^-{  }
      \ar@{=}[d]
      &&
      \mathrm{CE}
      \big(
        \mathfrak{l} B \mathrm{U}
      \big)
      \mathrlap{
        \;
        \simeq
        \;
        \underset{k \in \mathbb{N}}{\bigotimes}
        \mathrm{CE}
        \big(
          \mathfrak{l}K(\mathbb{Z}, 2k)
        \big)
      }
      \ar@{=}[d]
      \\
      \mathbb{R}
      \!
      \left[
        \!\!\!
        {\begin{array}{c}
          \omega_7,
          \\[-3pt]
          \omega_4
        \end{array}}
        \!\!\!
      \right]
      \!\big/\!
      \left(
      {\begin{aligned}
        d\, \omega_7 & = - \omega_4 \wedge \omega_4
        \\[-3pt]
        d\, \omega_4 & = 0
      \end{aligned}}
      \right)
      \ar@{<-^{)}}[rr]^-{
        \scalebox{.7}{$
        \left.
        \arraycolsep=1.4pt\def\arraystretch{.1}
        {\begin{array}{ccl}
          \omega_4 &\vert& 2k = 4
          \\
          0 &\vert& \mbox{else}
        \end{array}}
        \right\}
        \;\mapsfrom\;
        f_{2k}
        $}
      }
      && \;\;
      \mathbb{R}
      \!
      \left[
      \!\!\!
      {\begin{array}{c}
        \vdots
        \\[-3pt]
        f_4,
        \\[-3pt]
        f_2,
      \end{array}}
      \!\!\!
      \right]
      \,\big/\,
      \left(
      {\begin{aligned}
        \vdots
        \\[-3pt]
        d f_4 & = 0
        \\[-3pt]
        d f_2 & = 0
      \end{aligned}}
      \!\!
      \right)
    }
    }
  \end{equation}
  The morphism shown in \eqref{AbsoluteMinimalModelForBeta4ToKU0}
  evidently restricts to the relative rational Whitehead
  $L_\infty$-algebra inclusion (Prop. \ref{WhiteheadLInfinityAlgebrasRelative})
  on the factor
  $K(\mathbb{R}, 4) \,\subset\, L_{\mathbb{R}} B \mathrm{U}$
  and is zero elsewhere, hence fits
  into the required diagram \eqref{TransformationBetweenMinimalDerivedPLdRAdjunctionUnits}
  exhibiting it as an absolute minimal model
  (Def. \ref{RelativeFormalMaps})
  for $e_{B\mathrm{U}}^4$
  (by the commuting diagram in Prop. \ref{ExistenceOfRelativeMinimalSullianModels}).

\end{example}

\medskip

\noindent {\bf Cheeger-Simons homomorphism.}
Where the construction of the Chern-Weil homomorphism
(Def. \ref{ChernWeilHomomorphism}) invokes connections
on principal bundles
without actually being sensitive to this choice
(by Prop.  \ref{FundamentalTheoremOfChernWeilTheory}),
the \emph{Cheeger-Simons homomorphism}
\cite[\S 2]{CheegerSimons85}\cite[\S 3.3]{HopkinsSinger05}
(based on \cite{ChernSimons74})
is a refinement of the Chern-Weil homomorphism,
now taking values in differential ordinary cohomology
(Example \ref{OrdinaryDifferentialGeometry}),
that does detect connection data (hence ``differential'' data):

\vspace{-3mm}
\begin{equation}
  \label{CheegerSimonsHomomorphism}
  \raisebox{30pt}{
  \xymatrix@C=4em@R=.5em{
    G \mathrm{Connections}(X)_{/\sim}
    \ar[d]_-{
      \mbox{
        \tiny
        \color{greenii}
        \bf
        \def\arraystretch{.9}
        \begin{tabular}{c}
          forget
          \\
          connection
        \end{tabular}
      }
      \!\!\!
    }
    \ar[rr]_-{ \mathrm{cs}_G }^-{
      \mbox{
        \tiny
        \color{darkblue}
        \bf
        \def\arraystretch{.9}
        \begin{tabular}{c}
          Cheeger-Simons
          \\
          homomorphism
        \end{tabular}
      }
    }
    &&
    \mathrm{Hom}_{\mathbb{Z}}
    \Big(
      H^\bullet(B G;\, \mathbb{Z})
      \,,\,
      \overset{
        \mathclap{
        \raisebox{3pt}{
          \tiny
          \color{darkblue}
          \bf
          \def\arraystretch{.9}
          \begin{tabular}{c}
            differential
            \\
            cohomology
          \end{tabular}
        }
        }
      }{
        \widehat H^\bullet(X)
      }
   \; \Big)
    \ar[d]^-{
      \mbox{
        \tiny
        \color{darkblue}
        \bf
        curvature map
      }
    }
    \\
    G \mathrm{Bundles}(X)_{/\sim}
    \ar[rr]^-{ \mathrm{cw}_G }_-{
      \mbox{
        \tiny
        \color{darkblue}
        \bf
        \def\arraystretch{.9}
        \begin{tabular}{c}
          Chern-Weil
          \\
          homomorphism
        \end{tabular}
      }
    }
    &&
    \mathrm{Hom}_{\mathbb{R}}
    \Big(
      \mathrm{inv}^\bullet(\mathfrak{g})
      \,,\,
      \underset{
        \mathclap{
        \raisebox{-6pt}{
          \tiny
          \color{darkblue}
          \bf
          \def\arraystretch{.9}
          \begin{tabular}{c}
            de Rham
            \\
            cohomology
          \end{tabular}
        }
        }
      }{
        H^\bullet_{\mathrm{dR}}(X)
      }
    \Big)
  }
  }
\end{equation}

\noindent
We discuss how the general notion of
secondary non-abelian cohomology operations
(Def. \ref{SecondaryNonabelianCohomologyOperations})
specializes on ordinary principal bundles
to the Cheeger-Simons homomorphism, and hence
generalizes it to higher non-abelian cohomology:

\begin{lemma}[Characteristic classes of $G$-principal bundles have absolute minimal models]
  \label{CharacteristicClassesOfGPrincipalBundlesAreRelativeFormal}
  Let $G$ be a connected compact Lie group
  with classifying space $B G$ \eqref{ClassifyingSpace}.
  For $n \in \mathbb{N}$, let
  $[c] \,\in\, H^{n+1}(B G;\, \mathbb{Z})$
  be an indecomposable universal integral characteristic class for $G$-principal bundles
  (Example \ref{GroupCohomology}). Then every representative
  classifying map
  $\!\!
    \xymatrix@C=12pt{
      B G
      \ar[r]^-{\; c \; }
      &
      B^{n+1} \mathbb{Z}
    }
  \!\!$
  has an absolute minimal model in the sense of Def. \ref{RelativeFormalMaps}.
\end{lemma}
\begin{proof}
  By Lemma \ref{SullivanModelOfClassifyingSpace},
  the minimal Sullivan model for $B G$
  has vanishing differential,
  while the minimal Sullivan model of $B^{n + 1}\mathbb{Z}$
  is a polynomial algebra on a single degree $n+1$ generator (by Example \ref{RationalizationOfEMSpaces}),
  whose inclusion
  is already the relative minimal
  Sullivan model
  $\mathfrak{l}_{\scalebox{.7}{$B^{n+1}\mathbb{Z}$}}(c)$
  (Prop. \ref{WhiteheadLInfinityAlgebrasRelative})
  of $c$. Therefore, setting
  \vspace{-2mm}
  \begin{equation}
    \label{CharacteristicClassOnGBundlesIsRelativeFormal}
    \mathrm{CE}(\mathfrak{c})
      \,:=\,
    \mathrm{CE}
    \big(
      \mathfrak{l}_{\scalebox{.7}{$B^{n+1}\mathbb{Z}$}}(c)
    \big)
    \;:\;
    \xymatrix{
      \mathbb{R}
      [
        c
      ]
      \!\big/\!(d\, c = 0)
      \;
      \ar@{^{(}->}[r]
      &
      \;
      \mathrm{inv}^\bullet(\mathfrak{g})
    }
  \end{equation}

   \vspace{-2mm}
  \noindent  gives the required morphism of minimal models that makes
  makes the square \eqref{TransformationBetweenMinimalDerivedPLdRAdjunctionUnits}
  commute, by \eqref{CEAlgebraOfWhiteheadLInfinityAlgebraRelative}.
\end{proof}

In differential enhancement of Example \ref{CohomologyOfCoefficientsIsCohomologyOperations}
we have:
\begin{defn}[Secondary characteristic classes of differential non-abelian $G$-cohomology]
  \label{SecondaryCharacteristicClassesOfDifferentialNonAbelianGCohomology}
  Let $G$ be a connected compact Lie group
  with classifying space $B G$ \eqref{ClassifyingSpace}.
  By Lemma
  \ref{CharacteristicClassesOfGPrincipalBundlesAreRelativeFormal}),
  the construction of secondary characteristic
  classes
  (Def. \ref{SecondaryNonabelianCohomologyOperations},
  on differential non-abelian $G$-cohomology
  (Example \ref{DifferentialCohomologyOfPrincipalConnection})
  yields a $\mathbb{Z}$-linear map of the form
  \vspace{-4mm}
  $$
  \xymatrix{
    H^\bullet
    \big(
      B G
      ,;
      \mathbb{Z}
    \big)
    \;
    \simeq
    \;
    H
    \big(
      B G;
      \,
      B^{\bullet}\mathbb{Z}
    \big)
    \ar[rr]^-{ (-)_{\mathrm{diff}} }
    &&
    \widehat H
    \big(
      B G_{\mathrm{diff}};
      \,
      B^{\bullet} \mathbb{Z}
    \big)
    \;=\;
    H
    \big(
      B G_{\mathrm{diff}};
      \,
      B^{\bullet} \mathbb{Z}_{\mathrm{diff}}
    \big)
    \,,
  }
$$

    \vspace{-2mm}
  \noindent
  where on the right we have the
  ordinary differential non-abelian cohomology
  (Prop. \ref{DifferentialNonabelianCohomologySubsumesOrdinaryDifferentialCohomology})
  \emph{of} the moduli $\infty$-stack
  $B G_{\mathrm{diff}}$ \eqref{PullbackForModuliInfinityStackOfConnections}.
 Combined with the composition operation in $\Stacks$
\eqref{SmoothInfinityStacks} this gives a map
    \vspace{-2mm}
$$
  \xymatrix@C=20pt{
    \widehat H
    \big(
      X;
      \,
      B G
    \big)
    \!\times\!
    H
    \big(
      B G;
      \,
      B^{\bullet}\mathbb{Z}
    \big)
    \ar[rr]^-{ \mathrm{id} \times (-)_{\mathrm{diff}} }
    &&
    H
    \big(
      X;
      \,
      B G_{\mathrm{diff}}
    \big)
    \!\times\!
    H
    \big(
      B G_{\mathrm{diff}};
      \,
      B^{\bullet}\mathbb{Z}_{\mathrm{diff}}
    \big)
    \ar[r]^-{ \circ }
    &
    H
    \big(
      X;
      \,
      B^{\bullet} \mathbb{Z}_{\mathrm{diff}}
    \big)
    \,=\,
    \widehat H
    \big(
      X;
      \,
      B^{\bullet} \mathbb{Z}
    \big)
  }
$$

    \vspace{-1mm}
\noindent
which is $\mathbb{Z}$-linear in its second argument, and
whose hom-adjunct is

\vspace{-3mm}
\begin{equation}
  \label{SecondaryClassesOnDifferentialNonAbelianGCohomology}
  \xymatrix{
    \widehat
    H
    (
      X;
      \,
      B G
    )
    \ar[rrr]^-{
      \nabla
      \,\mapsto\,
      \left(
        c
        \,\mapsto\,
        c_{\mathrm{diff}}(\nabla)
      \right)
    }
    &&&
    \mathrm{Hom}_{\mathbb{Z}}
    \big(
      H
      (
        B G;
        \,
        B^\bullet \mathbb{Z}
      )
      \,,\,
      \widehat H
      (
        X;
        \,
        B^\bullet \mathbb{Z}
      )
    \big)
    \,.
  }
\end{equation}
\end{defn}

%\newpage

\begin{theorem}[Secondary non-abelian cohomology operations subsume Cheeger-Simons homomorphism]
  \label{SecondaryDifferentialNonAbelianCharacterSubsumesCheegerSimonsHomomorphism}
  Let $G$ be a connected compact Lie group,
  with classifying space denoted $B G$ \eqref{ClassifyingSpace}.
  Then the canonical construction \eqref{SecondaryClassesOnDifferentialNonAbelianGCohomology}
  of secondary characteristic classes
  on differential non-abelian $G$-cohomology
  (Def. \ref{SecondaryCharacteristicClassesOfDifferentialNonAbelianGCohomology})
  coincides with the
  Cheeger-Simons homomorphim \eqref{CheegerSimonsHomomorphism},
  in that the following diagram commutes:
  \vspace{-5mm}
  \begin{equation}
    \label{CheegerSimonsHomomorphismReproduced}
    \raisebox{3pt}{
    \xymatrix@R=.3em@C=4em{
      G\mathrm{Connections}(X)_{/\sim}
      \ar[d]_-{
        \mbox{
          \tiny
          \color{darkblue}
          \bf
          \eqref{FromGConnectionsToDifferentialNonabelianCohomology}
        }
      }
      \ar[rrr]_-{ \mathrm{cs}_G }^-{
        \mbox{
          \tiny
          \color{darkblue}
          \bf
          {
          \def\arraystretch{.9}
          \begin{tabular}{c}
            Cheeger-Simons
            \\
            homomorphism
          \end{tabular}}
        }
      }
      &&&
      \mathrm{Hom}_{\mathbb{Z}}
      \Big(
        H^\bullet(B G;\, \mathbb{Z})
        \,,\,
        \overset{
          \mathclap{
          \raisebox{3pt}{
            \tiny
            \color{darkblue}
            \bf
            \def\arraystretch{.9}
            \begin{tabular}{c}
              differential
              \\
              ordinary
              \\
              cohomology
            \end{tabular}
          }
          }
        }{
          \widehat H^\bullet(X)
        }
     \; \Big)
      \ar@{<-}[d]^-{
        \simeq
        \mathrlap{
          \mbox{
            \tiny
            \color{darkblue}
            \bf
            \eqref{RecoveringOrdinaryDifferentialCohomology}
          }
        }
      }
      \\
      \underset{
        \mathclap{
        \raisebox{3pt}{
          \tiny
          \color{darkblue}
          \bf
          {
          \def\arraystretch{.9}
          \begin{tabular}{c}
            differential non-abelian
            \\
            cohomology
          \end{tabular}}
        }
        }
      }{
      \widehat
      H
      (
        X;
        \,
        B G
      )
      }
      \ar[rrr]^-{
        \nabla
        \,\mapsto\,
        (
          c
          \,\mapsto\,
          c_{\mathrm{diff}}(\nabla)
        )
      }_-{
        \mathclap{
          \raisebox{-3pt}{
            \tiny
            \color{darkblue}
            \bf
            {
            \def\arraystretch{.9}
            \begin{tabular}{c}
              secondary
              \\
              non-abelian cohomology operations
            \end{tabular}}
          }
        }
      }
      &&&
      \mathrm{Hom}_{\mathbb{Z}}
      \left(
        H
        (
          B G;
          \,
          B^\bullet \mathbb{Z}
        )
        \,,\,
        \widehat
        H
        (
          X;
          \,
          B^{\bullet}\mathbb{Z}
        )
      \right)
      \,,
    }
    }
  \end{equation}
  \vspace{-.4cm}

  \noindent
  where on the left we have the map from $G$-connections to
  differential non-abelian $G$-cohomology
  from Prop. \ref{DifferentialCohomologyOfPrincipalConnection},
  and on the right the identification of ordinary differential cohomology
  from Prop. \ref{DifferentialNonabelianCohomologySubsumesOrdinaryDifferentialCohomology}.
\end{theorem}
\begin{proof}
  Let $c \in H\big(B G;\, B^{\bullet}\mathbb{Z}\big)$
  be a characteristic class,
  and let
  $(f^\ast E G, \nabla)$ be
  a $G$-principal bundle equipped with a $G$-connection.
  By Prop. \ref{DifferentialCohomologyOfPrincipalConnection},
  its image in differential non-abelian cohomology is
  given by the first map in the following diagram
  \vspace{-4mm}
  \begin{equation}
    \label{GConnectionsToTriplesToDifferentialCohomology}
    \hspace{-3mm}
    \xymatrix@R=-3pt@C=1.5em{
      G\mathrm{Connections}(X)_{/\sim}
      \ar[r]
      &
      \widehat H
      (
        X;
        \,
        B G
      )
      \ar[r]^-{ (c_{\mathrm{diff}})_\ast }
      &
      \widehat H
      \big(
        X;\, B^{n+1}\mathbb{Z}
      \big)
      \ar[r]^-{\simeq}
      &
      \widehat H^{n+1}(X)
      \\
      \big[
        f^\ast E G, \nabla
      \big]
      \ar@{}[r]|-{
        \longmapsto
      }
      &
      \big[
        f,
        \,
        \big(
          \mathrm{cs}_k
          (
            \nabla, f^\ast \nabla_{\mathrm{univ}}
          )
        \big)
        \,
        \big(
          \omega_k(F_\nabla)
        \big)
      \big]
         \ar@{}[r]|-{
        \longmapsto
      }
      &
      \big[
        f^\ast c,
        \,
          \mathrm{cs}_c
          (
            \nabla, f^\ast \nabla_{\mathrm{univ}}
          ),
        \,
        c(F_\nabla)
      \big]
    }
  \end{equation}
  \vspace{-1mm}

  \noindent  Here the triple of data are the three components
  (Example \ref{HomotopyPullbackViaTriples}) of a map
  into the defining homotopy pullback of differential non-abelian cohomology
  \eqref{DataTripleInDifferentialNonabelianCohomology}.
  Therefore, the secondary operation induced by the
  transformation \eqref{SecondaryCohomologyOperation}
  of these homotopy pullbacks, which in the present case
  is of this form:
  \vspace{-2mm}
\begin{equation}
  \raisebox{0pt}{
  \xymatrix{
     B G_{\mathrm{diff}}
     \ar[rr]_-{
       c_{\mathrm{diff}}
     }^-{
       \mbox{
         \tiny
         \color{darkblue}
         \bf
         \def\arraystretch{.9}
         \begin{tabular}{c}
           secondary
           \\
           characteristic class
         \end{tabular}
       }
     }
     \ar[dr]
     \ar[dd]_<<<<<<{ c_{BG} }
     &&
     B^{n+1} \mathbb{Z}_{\mathrm{diff}}
     \ar[dr]
     \ar[dd]|-{ \phantom{AAA} \atop \phantom{AA} }_<<<<<<{ c_{B^{n+1}\mathbb{Z}} }
     \\
     &
     \scalebox{1}{$
       \Omega_{\mathrm{dR}}(-;\mathfrak{l}B G)_{\mathrm{flat}}
     $}
     \ar[dd]
     \ar[rr]
     %|-{ \phantom{AA} }
     &&
     \scalebox{1}{$
       \Omega_{\mathrm{dR}}(-;\mathfrak{l}B^{n+1}\mathbb{Z})_{\mathrm{flat}}
     $}
     \ar[dd]
     \\
     B G
     \ar[rr]|-{ \phantom{AA} }_>>>>>>>>>>>{ c }^>>>>>>>>>>>>{
       \mbox{
         \tiny
         \color{darkblue}
         \bf
         characteristic class
       }
     }
     \ar[dr]_<<<<<<<<<<{ \mathbf{ch}_{B G} }
     &&
     B^{n+1}\mathbb{Z}
     \ar[dr]^<<<<<<<<<<<<{ \mathbf{ch}_{B^{n+1} \mathbb{Z}} }
     \\
     &
     \scalebox{1}{$
       \flatBexp(\mathfrak{l}B G )
     $}
     \ar[rr]_-{ \mathfrak{c} _\ast }
     &&
     \scalebox{1}{$
       \flatBexp(\mathfrak{l}B^{n+1}\mathbb{Z} )
     $}
     \,,
  }
  }
\end{equation}

\vspace{-2mm}
\noindent
  acts {\bf (a)} on the first component
  in the triple by postcomposition with
  $c$, hence as
  \vspace{-2mm}
  $$
    f \,\mapsto\, f^\ast c \,:=\, c \,\circ\, f
  $$

  \vspace{-2mm}
\noindent  and {\bf (b)} on the other two components by composition
  with $\mathfrak{c}$, which by \eqref{CharacteristicClassOnGBundlesIsRelativeFormal}
  corresponds to projecting out the
  Chern-Simons form and characteristic form corresponding
  to $c$, respectively.
  This is shown as the second map in \eqref{GConnectionsToTriplesToDifferentialCohomology}.
  Hence we are reduced to showing that the total map in
  \eqref{GConnectionsToTriplesToDifferentialCohomology}
  gives the Cheeger-Simons homomorphism.
  This statement is the content
  of \cite[\S 3.3]{HopkinsSinger05}.
\end{proof}

\begin{remark}[Secondary characteristic classes of $G$-connections]
  \label{SecondaryInvariantsOfGConnections}
  The traditional reason for referring to the Cheeger-Simons homomorphism
  \eqref{CheegerSimonsHomomorphismReproduced} as producing
  \emph{secondary} invariants is that Cheeger-Simons classes
  $\mathrm{cs}_{G}(P,\nabla) \in \widehat H(X)$
  may be non-trivial even if the
  underlying characteristic class $\mathrm{cw}_G(P)$
  (the ``primary'' class) vanishes.
  In this case the
  $\mathrm{cs}_{G}(P,\nabla)$ are
  also called \emph{Chern-Simons invariants}.

\vspace{1mm}
  \noindent {\bf (i)} This happens, in particular, when the $G$-connection $\nabla$
  is flat, $F(\nabla) = 0$ (by Def. \ref{CharacteristicForms}).
  Such secondary Chern-Simons invariants exhibit some subtle
  phenomena
  (\cite{Reznikov95}\cite{Reznikov96}\cite{IyerSimpson07}\cite{Esnault09}).

\vspace{1mm}
   \noindent {\bf (ii)} In fact,
  the proof of Theorem \ref{CheegerSimonsHomomorphismReproduced},
  via the triples \eqref{DataTripleInDifferentialNonabelianCohomology}
  of homotopy data, shows that, in this case,
  $\mathrm{cs}_{G}(P,\nabla)$ measures \emph{how} (or ``why'')
  $\mathrm{cw}_G(P)$ vanishes, namely by which class of
  homotopies.

\vspace{1mm}
   \noindent {\bf (iii)} Here we may understand secondary classes
  more abstractly,
  and explicitly related to the non-abelian character map:
  Where a (primary) non-abelian cohomology operation,
  according to Def. \ref{NonAbelianCohomologyOperations},
  is induced by a morphism of coefficient spaces \eqref{InducedCohomologyOperation},
  a secondary non-abelian cohomology operation,
  according to Def. \ref{SecondaryNonabelianCohomologyOperations},
  is induced \eqref{SecondaryOperationInduced}
  by a morphism
  of non-abelian character maps
  \eqref{TransformationOfDifferentialNonaebalianCharacters}
  -- hence by a morphism of morphisms --
  on these coefficient spaces.

\vspace{1mm}
  \noindent {\bf (iv)} Note that classical secondary cohomology operations themselves admit
  differential refinements. For instance, for the case of Massey products as secondary operations for
  the cup product this is worked out in \cite{GS-Massey}.
  While these can also fit into our context on general grounds, we will not demonstrate that explicitly here.
\end{remark}

\newpage

%%%%%%%%%%%%%%%%%%%%%%%%%%%%%%%%%%%%%%%%%%%%%%%%%%%%%%%%%%
\section{The twisted (differential) non-abelian character map}
  \label{TwistedNonabelianCharacterMap}
%%%%%%%%%%%%%%%%%%%%%%%%%%%%%%%%%%%%%%%%%%%%%%%%%%%%%%%%%%

We introduce the character map in twisted non-abelian cohomology
(Def. \ref{TwistedNonAbelianChernDoldCharacter})
and then discuss how it specializes to:

\cref{TwistedChernCharacterOnTopologicalKTheory} -- the twisted Chern character on (higher) K-theory;

\cref{CohomotopicalChernCharacter} -- the twisted character on
Cohomotopy theory.

\medskip

\noindent {\bf Rationalization in twisted non-abelian cohomology.}
In generalization of Def. \ref{RationalizationOfCoefficientsInNonabelianCohomology} we
now define rationalization of local coefficient bundles \eqref{LocalCoefficientBundle}.
This operation is transparent
in the language of $\infty$-category theory (Rem. \ref{InfinityCategoryTheory}),
where it simply amounts to forming the pasting composite with the
homotopy-coherent naturality square of the
$\mathbb{R}$-rationalization unit $\eta^{\mathbb{R}}$ (from Def. \ref{Lk}):

\vspace{-.5cm}
\begin{equation}
  \label{DirectRationalizationOfLocalCoefficients}
  \raisebox{30pt}{
  \xymatrix{
    X
    \ar[dr]_{\tau}^-{\ }="t"
    \ar[rr]^-{
      \overset{
        \mathclap{
        \raisebox{3pt}{
          \tiny
          \color{darkblue}
          \bf
          \def\arraystretch{.9}
          \begin{tabular}{c}
            $\tau$-twisted cocycle with
            \\
            local coefficients $\rho$
          \end{tabular}
        }
        }
      }{
        c
      }
    }_>>>>>>{\ }="s"
    &&
    A \!\sslash\! G
    \ar[dl]^{\rho}
    \\
    &
    B G
    \ar@{=>}^\simeq "s"; "t"
  }
  }
  \;\;\;\;\;\;\;\;
    \overset{
      \mathclap{
      \raisebox{3pt}{
        \tiny
        \color{darkblue}
        \bf
        rationalization
      }
      }
    }{
      \longmapsto
    }
  \;\;\;\;\;\;\;\;
  \raisebox{30pt}{
  \xymatrix{
    X
    \ar[dr]_{\tau}^-{\ }="t"
    \ar[rr]^-{
      \overset{
        \mathclap{
        \raisebox{3pt}{
          \tiny
          \color{darkblue}
          \bf
          \def\arraystretch{.9}
          \begin{tabular}{c}
            $\tau$-twisted cocycle with
            $\mathrlap{\mbox{\tiny\color{darkblue}\bf
              rationalized local coefficients $L_{\mathbb{R}}(\rho)$}}$
            \\
            \phantom{A}
          \end{tabular}
        }
        }
      }{
        c
      }
    }_>>>>>>{\ }="s"
    &&
    A \!\sslash\! G
    \ar[dl]^{\rho}
    \ar[r]^-{ \eta^{\mathbb{R}}_{A \sslash G} }
    &
    L_{\mathbb{R}} \big( A \!\sslash\! G\big)
    \ar[ddl]^-{
      L_{\mathbb{R}}(\rho)
    }_<<<<{\ }="s2"
    \\
    &
    B G
    \ar[dr]_{
      \eta^{\mathbb{R}}_{B G}
    }^-{\ }="t2"
    &&
    \\
    & &
    L_{\mathbb{R}} B G
    \ar@{=>}^\simeq "s"; "t"
    \ar@{=>}_\simeq "s2"; "t2"
  }
  }
\end{equation}

\vspace{-1mm}
\noindent Slightly less directly but equivalently,
this is the composite of
{\bf (a)} derived base change (Ex. \ref{BaseChangeQuillenAdjunction})
along $\eta^{\mathbb{R}}_{BG}$
from the slice over $B G$ to the slice over
$L_{\mathbb{R}} B G$, {\bf (b)} followed by the
composition with its derived naturality square, now regarded
as a morphism in the slice over $L_{\mathbb{R}} B G$:

$\,$
\vspace{-7mm}
$$
\hspace{-2mm}
  \raisebox{30pt}{
  \xymatrix{
    X
    \ar[dr]_{\tau}^-{\ }="t"
    \ar[rr]^-{
      \overset{
        \mathclap{
        \raisebox{3pt}{
          \tiny
          \color{darkblue}
          \bf
          \def\arraystretch{.9}
          \begin{tabular}{c}
            $\tau$-twisted cocycle with
            \\
            local coefficients $\rho$
          \end{tabular}
        }
        }
      }{
        c
      }
    }_>>>>>>{\ }="s"
    &&
    A \!\sslash\! G
    \ar[dl]_-{\rho}
    \\
    &
    B G
    \ar@{=>}^\simeq "s"; "t"
  }
  }
  \;\;\;
  \overset{
    \mathclap{
    \raisebox{3pt}{
      \tiny
      \color{darkblue}
      \bf base change
    }
    }
  }{
    \longmapsto
  }
  \;\;\;
  \raisebox{30pt}{
  \xymatrix{
    X
    \ar[dr]_{\tau}^-{\ }="t"
    \ar[rr]^-{
      \overset{
        \mathclap{
        \raisebox{3pt}{
          \tiny
          \color{darkblue}
          \bf
          \def\arraystretch{.9}
          \begin{tabular}{c}
          \end{tabular}
        }
        }
      }{
        c
      }
    }_>>>>>>{\ }="s"
    &&
    A \!\sslash\! G
    \ar[dl]^{\rho}
    \\
    &
    B G
    \ar[d]_-{ \eta^{\mathbb{R}}_{B G} }
    \\
    &
    L_{\mathbb{R}} B G
    \ar@{=>}^\simeq "s"; "t"
  }
  }
  \;\;\;\;
    \overset{
      \mathclap{
      \raisebox{3pt}{
        \tiny
        \color{darkblue}
        \bf
        \def\arraystretch{.9}
        \begin{tabular}{c}
          composition
          \\
          in slice
        \end{tabular}
      }
      }
    }{
      \longmapsto
    }
  \;\;\;\;
  \raisebox{30pt}{
  \xymatrix{
    X
    \ar[dr]_{\tau}^-{\ }="t"
    \ar[rr]^-{
      \overset{
        \mathclap{
        \raisebox{3pt}{
          \tiny
          \color{darkblue}
          \bf
          \def\arraystretch{.9}
          \begin{tabular}{c}
            $\tau$-twisted cocycle with
            \\
            local coefficients $\rho$
          \end{tabular}
        }
        }
      }{
        c
      }
    }_>>>>>>{\ }="s"
    &&
    A \!\sslash\! G
    \ar[dl]^{\rho}
    \ar[r]^-{
      \eta^{\mathbb{R}}_{ A \sslash G }
    }
    &
    L_{\mathbb{R}}
    (A \!\sslash\! G)
    \ar[ddll]^-{
      L_{\mathbb{R}}(\rho)
    }_<<<<<<{\ }="s2"
    \\
    &
    B G
    \ar[d]_-{ \eta^{\mathbb{R}}_{B G} }^-{\ }="t2"
    \\
    &
    L_{\mathbb{R}} B G
    \ar@{=>}^\simeq "s"; "t"
    \ar@{=>}_\simeq "s2"; "t2"
  }
  }
$$
It is in this second form that the
operation lends itself to formulation
in model category theory (Def. \ref{RationalizationInTwistedNonAbelianCohomology} below).
For that we just need to produce
a rectified (strictly commuting)
model of the $\eta^{\mathbb{R}}$-naturality
square:
\begin{defn}[Rectified rationalization unit on coefficient bundle]
  \label{RectifiedRationalizationUnitOnCoefficientBundle}
  Consider
  a local coefficient bundle \eqref{LocalCoefficientBundle}
  in $\NilpotentConnectedQFiniteHomotopyTypes$
  (Def. \ref{NilpotentConnectedSpacesOfFiniteRationalType})
  with its minimal relative Sullivan model \eqref{CEAlgebraOfWhiteheadLInfinityAlgebraRelative},
  (given by Prop. \ref{WhiteheadLInfinityAlgebrasRelative})

  \vspace{-3mm}
  \begin{equation}
    \label{ForRationalizationNaturalityRelativeSullivanModelForLocalCoefficientBundle}
    \hspace{-2cm}
    \xymatrix@R=2.8em{
      A \ar[rr]
      &
      \ar@{}[d]|-{
        \mbox{
          \tiny
          \color{darkblue}
          \bf
          \def\arraystretch{.9}
          \begin{tabular}{c}
            local coefficient bundle
          \end{tabular}
        }
      }
      &
      A \!\sslash\! G
      \ar[d]^-{ \rho }
      \\
      &&
      B G
      \mathrlap{\,,}
    }
    {\phantom{AAAAAA}}
    \xymatrix@C=3em{
      \Omega^\bullet_{\mathrm{PLdR}}
      \big(
        A \!\sslash\! G
      \big)
      \ar@{<-}[d]_-{
        \Omega^\bullet_{\mathrm{PLdR}}
        (\rho)
      }
      \ar@{<-}[rr]^-{
        p^{\mathrm{min}_{B G}}_{A \sslash G}
      }_-{ \in \; \mathrm{W} }
      &&
      \mathrm{CE}
      \big(
        \mathfrak{l}_{\scalebox{.5}{$B G$}}(A \!\sslash\! G)
      \big)
      \ar@{<-}[d]^-{
        \mathrm{CE}(\mathfrak{l}p)
      }
      \\
      \Omega^\bullet_{\mathrm{PLdR}}
      \big(
        B G
      \big)
      \ar@{<-}[rr]^-{
        p^{\mathrm{min}}_{B G}
      }_-{ \in \; \mathrm{W} }
      &&
      \mathrm{CE}
      \big(
        \mathfrak{l}(B G)
      \big)
    }
  \end{equation}
  \vspace{-3mm}

  \noindent
  Then the composite of the image of \eqref{ForRationalizationNaturalityRelativeSullivanModelForLocalCoefficientBundle}
  under $\Bexp_{\mathrm{PL}}$
  with the $\Omega^\bullet_{\mathrm{PLdR}} \dashv \Bexp_{\mathrm{PL}}$-adjunction unit
  (from Prop. \ref{QuillenAdjunctionBetweendgcAlgebrasAndSimplicialSets}):

  \vspace{-3mm}
  \begin{equation}
    \label{RectifiedRationalizationNaturalitySquare}
    \Derived\eta^{\mathrm{PLdR}}_{\rho}
    \;
      :=
    \;
    \raisebox{23pt}{
    \xymatrix@C=4em@R=15pt{
      A \!\sslash\! G
      \ar[rr]|-{
        \; \scalebox{0.6}{$\eta^{\mathrm{PLdR}}_{A \sslash G}$}\;
      }
      \ar[d]_-{
        \rho
      }
      \ar@/^1.7pc/[rrrr]|-{
        \;
    \scalebox{0.6}{$    \Derived\eta^{\mathrm{PLdR}}_{A \sslash G}
        \;\simeq\;
        \eta^{\mathbb{R}}_{A \sslash G}
        $}
        \;
      }
      &&
      \Bexp_{\mathrm{PL}} \,\circ \,\Omega^\bullet_{\mathrm{PLdR}}
      \big(
        A \!\sslash\! G
      \big)
      \ar@{->}[d]_-{\scalebox{0.6}{$
        \Bexp_{\mathrm{PL}}
        \,\circ\,
        \Omega^\bullet_{\mathrm{PLdR}}
        (\rho)
        $}
      }
      \ar@{->}[rr]|-{
        \;
      \scalebox{0.6}{$
        \Bexp_{\mathrm{PL}}
        \big(
          p^{\mathrm{min}_{BG}}_{A \sslash G}
        \big)
        $}
        \;
      }
      &&
      \Bexp_{\mathrm{PL}}
      \,\circ\,
      \mathrm{CE}
      \big(
        \mathfrak{l}_{\scalebox{.5}{$B G$}}(A \!\sslash\! G)
      \big)
      \ar@{->}[d]_-{\scalebox{0.6}{$
        \Bexp_{\mathrm{PL}}
        \,\circ\,
        \mathrm{CE}(\mathfrak{l}p)
        $}
      }
      \\
      B G
      \ar[rr]|-{ \;
      \scalebox{0.6}{$\eta^{\mathrm{PLdR}}_{B G}$} \; }
      \ar@/_1.7pc/[rrrr]|-{
        \;
       \scalebox{0.6}{$ \mathbb{D}\eta^{\mathrm{PLdR}}_{B G}
        \;\simeq\;
        \eta^{\mathbb{R}}_{B G}
        $}
        \;
      }
      &&
      \Bexp_{\mathrm{PL}}
      \,\circ\,
      \Omega^\bullet_{\mathrm{PLdR}}
      \big(
        B G
      \big)
      \ar@{->}[rr]|-{
        \;
         \scalebox{0.6}{$ \Bexp_{\mathrm{PL}}
          \big(
            p^{\mathrm{min}}_{B G}
          \big)
          $}
        \;
        }
      &&
      \Bexp_{\mathrm{PL}}
      \,\circ\,
      \mathrm{CE}
      \big(
        \mathfrak{l}(B G)
      \big)
    }
    }
  \end{equation}
  \vspace{-3mm}

  \noindent
  is, after passage \eqref{LocalizationOfAModelCategoryAtWeakEquivalenes}
  to the classical homotopy category (Example \ref{TheClassicalHomotopyCategory}),
  the naturality square of the rationalization unit
  on $\rho$ \eqref{RationalizationUnit},
  namely of the derived adjunction unit \eqref{DerivedAdjunctionUnit}
  $\eta^{\mathbb{R}} = \Derived^{\mathrm{PL}\mathbb{R}\mathrm{dR}}$
  (using, with Prop. \ref{RelativeSullivanModelsAreCofibrations},
  that the right part of \eqref{RectifiedRationalizationNaturalitySquare}
  is the image under $\Bexp_{\mathrm{PL}}$ of a fibrant replacement morphism.)
\end{defn}

\newpage

\begin{defn}[Rationalization in twisted non-abelian cohomology]
  \label{RationalizationInTwistedNonAbelianCohomology}
  Given a local coefficient bundle $\rho$
  and its rectified rationalization unit
  $\mathbb{D}\eta^{\mathrm{PLdR}}_\rho$
  (Def. \ref{RectifiedRationalizationUnitOnCoefficientBundle})
  we say that {\it rationalization}
  in twisted non-abelian cohomology with local coefficients $\rho$
  (Def. \ref{NonabelianTwistedCohomology})
  is the twisted non-abelian cohomology operation
  (Def. \ref{TwistedNonabelianCohomologyOperation})
  \vspace{-2mm}
  \begin{equation}
  \label{RationalizationOperationInTwistedNonabelianCohomology}
  \big(
    \eta^{\mathbb{R}}_\rho
  \big)_\ast
  \;:\;
  \xymatrix{
    H^\tau
    (
      X;
      \,
      A
    )
    \ar[rrrr]^-{\scalebox{0.7}{$
      \big(
        \mathbb{D}\eta^{\mathrm{PLdR}}_\rho
        \,\circ\,
        (-)
      \big)
      \;\circ\;
      \LeftDerived
      \big(
        \eta^{\mathbb{R}}_{B G}
      \big)_!
      $}
    }
    &&&&
    H^{L_{\mathbb{R}} \tau}
    \big(
      X;
      \,
      L_{\mathbb{R}}A
    \big)
  }
  \end{equation}

  \vspace{-2mm}
  \noindent   given by the composite of

  \noindent
  {\bf (a)} derived left base change
    $\LeftDerived(\eta^{\mathbb{R}}_{B G})_!$ (Ex. \ref{BaseChangeQuillenAdjunction})
    along the rationalization unit \eqref{RationalizationUnit}
    on the classifying space of twists,

  \noindent
  {\bf (b)} with the rectified
   rationalization unit \eqref{RectifiedRationalizationNaturalitySquare}
   on the coefficient bundle,
   regarded as a morphism in
   the homotopy category \eqref{LocalizationOfAModelCategoryAtWeakEquivalenes} of
   the slice model category
   (Example \ref{SliceModelCategory})
   of $\SimplicialSets_{\mathrm{Qu}}$
   (Example \ref{ClassicalModelStructureOnSimplicialSets})
   over $\Bexp_{\mathrm{PL}} \,\circ\, \mathrm{CE}(\mathfrak{l}B G) )$.
\end{defn}

\begin{remark}[Rationalization of coefficients and/or of twists]
  Def. \ref{RationalizationInTwistedNonAbelianCohomology}
  rationalizes both the coefficients as well as their twist.
  This is of interest because:

  \noindent
  {\bf (a)} the joint rationalization is defined canonically,
  in fact functorically, as highlighted around \eqref{DirectRationalizationOfLocalCoefficients};

  \noindent
  {\bf (b)} the rationalized twisting
    appears in the archetypical examples
    (such as the twisted Chern character on degree-3 twisted K-theory,
    \cref{TwistedChernCharacterOnTopologicalKTheory})
    and
    gives the {\it Bianchi identities}
    on higher form field/flux data relevant in applications
    to physics
    (see \cref{CohomotopicalChernCharacter}).

  One may also consider rationalization of just the coefficients,
  keeping non-rationalized twists;
  but, in general, this
  requires making a choice, namely a choice of dashed morphisms
  in the following transformation diagram of local coefficient bundles:

  \vspace{-3mm}
  \begin{equation}
    \label{RationalizationOfCoefficientsOverFixedTwist}
    \begin{tikzcd}[row sep=small]
      \overset{
        \mathclap{
        \raisebox{3pt}{
          \tiny
          \color{darkblue}
          \bf
          \def\arraystretch{.9}
          \begin{tabular}{c}
            local
            coefficient
            \\
            bundle
          \end{tabular}
        }
        }
      }{
        A
      }
      \ar[
        d,
        "\mathrm{hofib}(\rho)"
      ]
      \ar[
        rr,
        "\eta^{\mathbb{R}}_A"
      ]
      &&
      L_{\mathbb{R}}(A)
      \mathrlap{
        \mbox{
          \tiny
          \color{darkblue}
          \bf
          \def\arraystretch{.9}
          \begin{tabular}{c}
            $\mathbb{R}$-rationalized
            \\
            local coefficients
          \end{tabular}
        }
      }
      \ar[
        d,
        "\mathrm{hofib}(\rho^{\mathbb{R}})"
      ]
      \\
      A
      \!\sslash\!
      G
      \ar[
        d,
        "\rho"
      ]
      \ar[
        rr,
        dashed
      ]
      &&
      \big(
        L_{\mathbb{R}}A
      \big)
      \!\sslash\!
      G
      \ar[
        d,
        dashed,
        "\rho^{\mathbb{R}}"
      ]
      \\
      B G
      \ar[
        rr,-,
        shift left=1pt
      ]
      \ar[
        rr,-,
        shift right=1pt
      ]
      &&
      B G
      \mathrlap{
        \mbox{
          \tiny
          \color{darkblue}
          \bf
          \def\arraystretch{.9}
          \begin{tabular}{c}
            non-rationalized
            \\
            twist
          \end{tabular}
        }
      }
    \end{tikzcd}
    {\phantom{AAA}}
    \in
    \;
    \mathrm{Ho}
    \big(
      \SimplicialSets_{\mathrm{Qu}}
    \big)
    \,.
  \end{equation}
  \vspace{-.4cm}

  \noindent
  The homotopy-commutativity of the bottom square expresses that
  and how {\it rationalization commutes with twisting}.

  Given such a choice, then using the bottom square \eqref{RationalizationOfCoefficientsOverFixedTwist}
  in place of the rationalization unit's naturality square
  on the  right of \eqref{DirectRationalizationOfLocalCoefficients}
  produces a definition,
  directly analogous to Def. \ref{RationalizationInTwistedNonAbelianCohomology},
  of rationalization
  of just the $A$-coefficients in
  twisted $A$-cohomology.
  This is also of interest
  (see for instance the case of twisted $\mathrm{KO}$-theory
  in \cite[Prop. 4]{GS-tKO}), but currently we do not further
  expand on this generalization here.
\end{remark}

\medskip

\noindent {\bf Twisted non-abelian character map.}
In generalization of Def. \ref{NonAbelianChernDoldCharacter},
we set:

\begin{defn}[Twisted non-abelian character map]
  \label{TwistedNonAbelianChernDoldCharacter}
  Let $X \in \NilpotentConnectedQFiniteHomotopyTypes$
  (Def. \ref{NilpotentConnectedSpacesOfFiniteRationalType})
  equipped with the structure of a smooth manifold,
  and
  \vspace{0mm}
  \begin{equation}
    \label{LocalCoefficientBundleForPushforwardToRationalCoefficients}
    \xymatrix@R=1.5em{
      A \ar[rr]
      &
      \ar@{}[d]|-{
        \mbox{
          \tiny
          \color{darkblue}
          \bf
          \def\arraystretch{.9}
          \begin{tabular}{c}
            local coefficient bundle
          \end{tabular}
        }
      }
      &
      A \!\sslash\! G
      \ar[d]^-{ \rho }
      \\
      && B G
    }
  \end{equation}

  \vspace{-2mm}
  \noindent
  be a local coefficient bundle \eqref{LocalCoefficientBundle}
  in $\NilpotentConnectedQFiniteHomotopyTypes$
  (Def. \ref{NilpotentConnectedSpacesOfFiniteRationalType}).
  Then the \emph{twisted non-abelian character map}
  in twisted non-abelian cohomology
  is the twisted cohomology operation
  \vspace{-2mm}
  \begin{equation}
    \label{TwistedNonAbelianChernDoldAsCompositeOfRationalizationWithdeRham}
    \mathllap{
      \mbox{
        \tiny
        \color{darkblue}
        \bf
        \def\arraystretch{.9}
        \begin{tabular}{c}
          twisted
          \\
          non-abelian
          \\
          character map
        \end{tabular}
      }
      \;\;
    }
    \mathrm{ch}_\rho
    \;:\;
    \xymatrix@C=3em{
      \overset{
        \raisebox{3pt}{
          \tiny
          \color{darkblue}
          \bf
          \def\arraystretch{.9}
          \begin{tabular}{c}
            twisted
            non-abelian
            \\
            cohomology
          \end{tabular}
        }
      }{
      H^{\tau}
      (
        X;
        \,
        A
      )
      }
      \ar[rr]^-{
        (\eta^\mathbb{R}_\rho)_\ast
      }_-{
        \mbox{
          \tiny
          \color{greenii}
          \bf
          rationalization
        }
      }
      &&
      \overset{
        \mathclap{
        \raisebox{3pt}{
          \tiny
          \color{darkblue}
          \bf
          \def\arraystretch{.9}
          \begin{tabular}{c}
            twisted
            non-abelian
            \\
            real cohomology
          \end{tabular}
        }
        }
      }{
      H^{L_{\mathbb{R}}\tau}
      \big(
        X;
        \,
        L_{\mathbb{R}}A
      \big)
      }
      \ar[rr]^-{ \simeq }_-{
        \mathclap{
        \mbox{
          \tiny
          \color{greenii}
          \bf
          \def\arraystretch{.9}
          \begin{tabular}{c}
            twisted
            non-abelian
            \\
            de Rham theorem
          \end{tabular}
        }
        }
      }
      &&
      \overset{
        \mathclap{
        \raisebox{3pt}{
          \tiny
          \color{darkblue}
          \bf
          \def\arraystretch{.9}
          \begin{tabular}{c}
            twisted
            non-abelian
            \\
            de Rham cohomology
          \end{tabular}
        }
        }
      }{
      H^{\tau_{\mathrm{dR}}}_{\mathrm{dR}}
      (
        X;
        \,
        \mathfrak{l}A
      )
      }
    }
  \end{equation}

\vspace{-2mm}
\noindent
  from twisted non-abelian $A$-cohomology
  (Def. \ref{NonabelianTwistedCohomology})
  to twisted non-abelian de Rham cohomology
  (Def. \ref{TwistedNonabelianDeRhamCohomology})
  with local coefficients
  in the rational relative Whitehead
  $L_\infty$-algebra $\mathfrak{l}\rho$ of $\rho$
  (Prop. \ref{RationalizationOfLocalCoefficients})
  which is the composite of

  {\bf (i)} the operation  \eqref{RationalizationOperationInTwistedNonabelianCohomology}
  of rationalization of local coefficients
  (Def. \ref{RationalizationInTwistedNonAbelianCohomology}),

  {\bf (ii)} the equivalence \eqref{EquivalenceBetweenTwistedNonabelianRealAndNonaebelianDeRhamCohomology}
  of the twisted non-abelian de Rham theorem
  (Theorem \ref{TwistedNonAbelianDeRhamTheorem}).
\end{defn}

%%%%%%%%%%%%%%%%%%%%%%%%%%%%%%%%%%%%%%%%%%%%%%%%%%%%%%%%%%%%%%
\subsection{Twisted Chern character on higher K-theory}
\label{TwistedChernCharacterOnTopologicalKTheory}
%%%%%%%%%%%%%%%%%%%%%%%%%%%%%%%%%%%%%%%%%%%%%%%%%%%%%%%%%%%%%%%

We discuss (Prop. \ref{TwistedChernCharacterInTwistedTopologicalKTheory})
how the twisted non-abelian character map
reproduces the  twisted Chern character in twisted
topological K-theory  \cite[\S 6.3]{BCMMS02}\cite{MathaiStevenson03}\cite[\S 7]{AtiyahSegal06}
-- see also \cite{TuXu06}\cite[\S 6]{MathaiStevenson06}\cite[\S 2]{FHT-complex}\cite{BGNT08} \cite[\S 4]{GomiTerashima10}\cite[\S 8.3]{Karoubi12}\cite[\S 3.2]{GS-tAHSS}\cite{GS-RR}.
Then we also consider (Prop. \ref{TwistedChernCharacterInTwistedHigherKTheory})
the twisted character
on twisted iterated K-theory \cite[\S 2.2]{LindSatiWesterland16}.

\newpage

\noindent {\bf Character maps on higher-twisted ordinary K-theories.}

\begin{remark}[Twisted Chern character via twisted characteristic forms]
  The twisted Chern character on twisted K-theory was first
  proposed in \cite{BCMMS02}
  via a natural twisted generalization of the
  component-wise construction
  \eqref{OrdinaryChernCharacterReproduced}
  of the ordinary Chern character
  in terms of characteristic curvature forms. Briefly,
given a degree-3 twist $\tau_3 \,\in\, H^3(X;\, \mathbb{Z})$
  on a (compact) smooth manifold,
  then every
  class
  $[(V_1,V_2)] \,\in\, \mathrm{KU}^{\tau_3}(X)$
  in $\tau_3$-twisted $\mathrm{KU}$-theory
  (Ex. \ref{TwistedKTheory}) may be represented by
  a pair $(V_1, V_2)$ of $\tau_3$-twisted complex vector bundles
  (\cite[\S 7.2]{LupercioUribe04}), generally of infinite rank
  (``bundle gerbe modules'' \cite[\S 4]{BCMMS02}).
  Now given a choice of lift of $\tau_3$
  through the characteristic class map \eqref{DifferentialCohomologyDiagram}
  to a Deligne cocycle $[h_0, A_1, B_2, H_3]$
  (Ex. \ref{OrdinaryDifferentialGeometry}) with respect to some
  open cover $p : (\sqcup_i U_i) \twoheadrightarrow X$ (Ex. \ref{GoodOpenCoversAreProjectivelyCofibrantResolutionsOfSmoothManifolds}),
  hence in particular including a choice of ``local $B$-field'' $B_2$,
  then one may further choose $B_2$-twisted connections $\nabla_i$
  \cite{Mackaay03}
  on the twisted vector bundles, inducing curvature 2-forms $F_i$:

  \vspace{-.6cm}
  \begin{equation}
    \label{TwistedCurvature2Forms}
    B_2 \,\in\, \Omega^2(U)
    \,,
    {\phantom{AAA}}
    F_i \;\in\; \Omega^2\big( U;\, \mathrm{End}(V_i) \big)
    \,.
  \end{equation}
  \vspace{-.5cm}

  \noindent
  Now it turns out \cite[Prop. 9.1]{BCMMS02} that

  \begin{itemize}

    \vspace{-2mm}
    \item
    the following trace \eqref{TwistedChernCharacterForm}
    of differences of wedge-product exponentials
    of these 2-forms \eqref{TwistedCurvature2Forms}
    is
    well defined (i.e., the trace exists, which is non-trivial since the
    twisted vector bundles are in general have infinite rank)

    \vspace{-.4cm}
    \begin{equation}
      \label{TwistedChernCharacterForm}
      \exp(B_2) \wedge
      \mathrm{tr}
      \big(
        \exp(F_1)
        -
        \exp(F_2)
      \big)
      \;\;
      =
      \;\;
      p^\ast
      \mathrm{ch}^{B_2}(\nabla_1, \nabla_2)
      \;\;\;\;\;
        \in
      \;
      \Omega^{2\bullet}(U)
    \end{equation}
    \vspace{-.5cm}

    \noindent
    and equals
    the pullback $p^\ast(-)$ to the given cover
    of an even-degree differential form on the base space $X$,

    \vspace{-.2cm}
    \item
    which is closed in the $H_3$-twisted de Rham complex
    \eqref{ChainComplexForH3TwisteddeRhamCohomology}

    \vspace{-.4cm}
    $$
      (d - H_3 \wedge ) \mathrm{ch}^{B_2}(\nabla_1, \nabla_2)
      \;\;
      =
      \;\;
      0
      ,.
    $$
    \vspace{-.6cm}

    \vspace{-0cm}
    \item
    and whose resulting twisted de Rham cohomology class
    (Def. \ref{Degree3TwistedAbelianDeRhamCohomology}) is
    independent of the choices made:

    \vspace{-.4cm}
    \begin{equation}
      \label{TwistedChernCharacterViaCharacteristicForms}
      \mathrm{ch}^{\tau_3}(V_1, V_2)
      \;
      \coloneqq
      \;
      \big[
        \mathrm{ch}^{B_2}(\nabla_1, \nabla_2)
      \big]
      \;\;\;
      \in
      \;
      H^{3 + H_3}_{\mathrm{dR}}(X)
      \,.
    \end{equation}
    \vspace{-.5cm}

  \end{itemize}
  \vspace{-.2cm}

  \noindent
  This class \eqref{TwistedChernCharacterViaCharacteristicForms}
  was proposed \cite[p. 26]{BCMMS02}
  to be
  the {\it twisted Chern character} of the twisted $K$-theory class
  $[(V_1, V_2)]$.
\end{remark}

A more intrinsic characterization of the twisted Chern character was
later found in \cite[\S 2]{FHT-complex}. This is the form in which
one recognizes the twisted Chern character as an example
of the twisted non-abelian character map (Def. \ref{TwistedNonAbelianChernDoldCharacter}),
in twisted enhancement of Example \ref{ChernCharacterInKTheory}:
\begin{prop}[Twisted Chern character in twisted topological K-theory]
 \label{TwistedChernCharacterInTwistedTopologicalKTheory}
 Consider twisted complex topological K-theory
 $\mathrm{KU}^\tau(-)$
 (Example \ref{TwistedKTheory}), for
 degree-3 twists given (via Example \ref{BundleGerbes}) by
 \vspace{-1mm}
 $$
   \tau
   \;\in\;
   H
   \big(
     -;
     \,
     B^2 \mathrm{U}(1)
   \big)
   \;\simeq\;
   H^3(-;, \mathbb{Z})
   \,,
 $$

   \vspace{-1mm}
  \noindent
 and regarded, via \eqref{ReproducingTwistedTopologicalKTheory},
 as twisted non-abelian cohomology
 with local coefficients in
 $\mathbb{Z} \times B \mathrm{U} \sslash B^2 \mathrm{U}(1)$
 \eqref{LocalCoefficientBundleForTwistedKTheory}.
  Then the twisted non-abelian character map
 (Def. \ref{TwistedNonAbelianChernDoldCharacter})
 $\mathrm{ch}^{\tau}_{\mathbb{Z} \times B \mathrm{U}}$
 is equivalent to the traditional twisted Chern character
 $\mathrm{ch}^\tau$ on twisted K-theory
 with values in $H_3$-twisted de Rham cohomology
 (Def. \ref{Degree3TwistedAbelianDeRhamCohomology}):
   \vspace{-1mm}
 $$
   \overset{
     \mathclap{
     \raisebox{3pt}{
       \tiny
       \color{darkblue}
       \bf
       \def\arraystretch{.9}
       \begin{tabular}{c}
         twisted non-abelian
         \\
         character map
       \end{tabular}
     }
     }
   }{
     \mathrm{ch}^\tau_{\mathbb{Z} \times B \mathrm{U}}
   }
   \;\;\;\;
   \simeq
   \;\;\;\;
   \overset{
     \mathclap{
     \raisebox{3pt}{
       \tiny
       \color{darkblue}
       \bf
       \def\arraystretch{.9}
       \begin{tabular}{c}
         twisted
         \\
         Chern character
       \end{tabular}
     }
     }
   }{
   \mathrm{ch}^\tau.
   }
 $$
\end{prop}
\begin{proof}
  That the codomain of the twisted non-abelian character map
  %, in this case,
  $\mathrm{ch}^\tau_{\mathbb{Z} \times B \mathrm{U}}$
  is indeed $H_3$-twisted de Rham cohomology
  is the content of Prop. \ref{TwistedNonabelianDeRhamCohomologySubsumesH3TwisteddeRhamCohomology}.
  With this, and due to the twisted non-abelian de Rham theorem
  (Theorem \ref{TwistedNonAbelianDeRhamTheorem}), it is
  sufficient to see that the general rationalization map
  of local non-abelian coefficients from Def. \ref{RationalizationInTwistedNonAbelianCohomology}
  reproduces the rationalization map underlying the
  twisted Chern character.
  This is manifest from
  comparing the
  rationalization operation \eqref{DirectRationalizationOfLocalCoefficients},
  that is made formally precise by Def. \ref{RationalizationInTwistedNonAbelianCohomology},
  to the description of the twisted Chern character as given in
  \cite[(2.8)-(2.9)]{FHT-complex}.
\end{proof}

\begin{remark}[Twisted Pontrjagin character in twisted KO-theory]
Similarly, an analogous statement holds for the twisted
  Pontrjagin character (as in Example \ref{PontrjaginCharacterInKSO})
  on twisted real K-theory
  \cite[Prop. 2]{GS-tKO}.
\end{remark}

\begin{example}[Twisted Chern character on higher Cohomotopy-twisted K-theory]
  \label{ChernCharacterOnHigherCohomotopicallyTwistedKTheory}
  For $k \in \mathbb{N}_+$,
  consider the
  cohomotopically-twisted complex K-theory from
  Ex. \ref{CohomotopyTwistedKTheory}.

  \noindent
  {\bf (i)}
  For $\lambda \in \pi^{2k+1}(X)$
  a Cohomotopy class (Ex. \ref{CohomotopyTheory})
  regarded now as a twist,
  the corresponding twisted non-abelian character map
  (Def. \ref{TwistedNonAbelianChernDoldCharacter})
  lands,
  by Theorems \ref{NonAbelianDeRhamTheorem}, \ref{TwistedNonabelianDeRhamCohomology}
  and Examples \ref{RationalizationOfEMSpaces}, \ref{RationalizationOfnSpheres},
  in $\lambda_{\mathrm{dR}} =: H_{2k+1}$-twisted de Rham cohomology
  (Example \ref{RecoveringHigherTwistedDifferentialForms},
   Prop. \ref{TwistedNonabelianDeRhamCohomologySubsumesHigherTwisteddeRhamCohomology}):

  \vspace{-.5cm}
  \begin{equation}
    \label{CharacterMapOnCohomotopicallyTwistedKTheory}
    \overset{
      \mathclap{
      \raisebox{3pt}{
        \tiny
        \color{darkblue}
        \bf
        \def\arraystretch{.9}
        \begin{tabular}{c}
          twist in Cohomotopy
        \end{tabular}
      }
      }
    }{
      \lambda
        \,\in\,
      \pi^{2k+1}(X)
    }
    \;\;\;\;
      \vdash
    \;\;\;\;
    \begin{tikzcd}
      \overset{
        \mathclap{
        \raisebox{3pt}{
          \tiny
          \color{darkblue}
          \bf
          higher Cohomotopy-twisted
          K-theory
        }
        }
      }{
      \mathrm{KU}^{\lambda}(X)
      \;=\;
      H^\lambda
      \big(
        X;
        \,
        \mathrm{KU}_0
      \big)
      }
      \ar[
        rrr,
        "{ \mathrm{ch}^{\lambda}_{\mathrm{KU}_0} }"{above},
        "\mbox{
          \tiny
          \color{greenii}
          \bf
          twisted character map
        }"{below,yshift=-2pt}
      ]
      &&&
      \overset{
        \mathclap{
        \raisebox{3pt}{
          \tiny
          \color{darkblue}
          \bf
          \def\arraystretch{.9}
          \begin{tabular}{c}
            higher twisted de Rham cohomology
          \end{tabular}
        }
        }
      }{
        H^{\lambda_{\mathrm{dR}}}_{\mathrm{dR}}
        \big(
          X;
          \mathfrak{l}\mathrm{KU}_0
        \big)
        \;\simeq\;
        H^{ \bullet + H_{2k+1} }(X)
      }
      \,.
    \end{tikzcd}
  \end{equation}

  \noindent {\bf (ii)}
  A map of this form
  has been defined in \cite[\S 2.2]{MMS20}, by direct construction
  on form representatives. However, \cite[Thm. 4.19]{MMS20}
  implies that this component construction coincides
  with the rationalization map (Def. \ref{RationalizationInTwistedNonAbelianCohomology})
  on the local coefficient bundle
  \eqref{LocalCoefficientBundleForCohomotopicallyTwistedKtheory},
  up to application of the de Rham theorem.
  Therefore, the twisted character map \eqref{CharacterMapOnCohomotopicallyTwistedKTheory}
  obtained as a special case of Def. \ref{TwistedNonAbelianChernDoldCharacter},
  reproduces the MMS-Character
  \cite[\S 2.2]{MMS20}
  on higher Cohomotopy-twisted K-theory.

\end{example}

\begin{remark}[Charge quantization of spherical T-duality in M-theory]
  For $k = 3$,
  the character map \eqref{CharacterMapOnCohomotopicallyTwistedKTheory}
  on 7-Cohomotopy-twisted K-theory
  is a candidate for charge quantization \eqref{DiracChargeQuantization}
  of the super-rational M-theory fields
  participating in 3-spherical T-duality
  over 11-dimensional super-spacetime,
  as derived in \cite[Prop. 4.17, Rem. 4.18]{FSS18}
  (review in \cite[(8), (19)]{SS18}).
  However, the 7-Cohomotopy-twisted K-theory character
  has some spurious fields of 2-periodic degree in its image,
  which are not seen in the physics application, where the
  field degrees are 6-periodic \cite[\S 3]{Sa-Higher}\cite[(65)]{FSS18}.
  Another candidate
  for charge-quantization of
  the super-rational M-theory fields
  participating in 3-spherical T-duality,
  possibly more accurately reflecting the physics,
  is the character map on twisted higher K-theory \cite{LindSatiWesterland16},
  which we turn to next (Prop. \ref{TwistedChernCharacterInTwistedHigherKTheory}).
\end{remark}

\medskip

\noindent {\bf Character map on twisted higher K-theory.}

\begin{remark}[Higher twisted de Rham coefficients inside rational twisted iterated K-theory]
  \label{HigherTwistedDeRhamCoefficientsInsideRationalTwistedIteratedKTheory}
  There is a non-trivial
  twisted cohomology operation (Def. \ref{TwistedNonabelianCohomologyOperation})
  from
    {\bf (a)}
  twisted non-abelian de Rham cohomology
  (Def. \ref{TwistedNonabelianDeRhamCohomology})
  with coefficients in the
  relative rational Whitehead $L_\infty$-algebra
  (Prop. \ref{WhiteheadLInfinityAlgebrasRelative})
  of the coefficient bundle \eqref{LocalCoefficientBundleForIteratedKTheory}
  of
  twisted iterated K-theory (Ex. \ref{TwistedIteratedkTheory})
  to
  {\bf (b)}
  higher twisted de Rham cohomology (Def. \ref{HigherTwistedAbelianDeRhamCohomology})
  regarded as twisted non-abelian de Rham cohomology via Prop. \ref{TwistedNonabelianDeRhamCohomologySubsumesHigherTwisteddeRhamCohomology}):
  \vspace{-2mm}
  \begin{equation}
    \label{RationalLSWCharacter}
    \xymatrix{
      H_{\mathrm{dR}}^{\tau_{\mathrm{dR}}}
      \Big(
        -;
        \,
        \mathfrak{l}
        K^{\circ_{2r-2}}(\mathrm{ku})_1
      \Big)
      \ar[rr]^-{ \phi_\ast }
      &&
      H_{\mathrm{dR}}^{\tau_{\mathrm{dR}}}
      \Big(
        -;
        \underset{k \in \mathbb{N}}{\bigoplus}
        \mathfrak{b}^{2rk}\mathbb{R}
      \Big)
      \,,
    }
  \end{equation}

  \vspace{-1mm}
  \noindent
  given, under the twisted non-abelian de Rham theorem
  (Theorem \ref{TwistedNonAbelianDeRhamTheorem})
  by the LSW-character
  from \cite[\S 2.2]{LindSatiWesterland16}
  applied to rational coefficients.
\end{remark}

\begin{prop}[Twisted Chern character in twisted iterated K-theory]
 \label{TwistedChernCharacterInTwistedHigherKTheory}
 For $r \in \mathbb{N}$, $r \geq 1$,
 consider twisted iterated K-theory
 $\big(K^{\circ_{2r-2}}(\mathrm{ku})\big)^\tau$
 (Example \ref{TwistedIteratedkTheory}), for
 degree-$(2r+1)$ twists given (via Example \ref{HigherBundleGerbes}) by
 \vspace{-1mm}
 $$
   \tau
   \;\in\;
   H
   \big(
     -;
     \,
     B^{2r} \mathrm{U}(1)
   \big)
   \;\simeq\;
   H^{2r+1}(-;, \mathbb{Z})
   \,,
 $$

 \vspace{-1mm}
 \noindent
 and regarded, via Example \ref{TwistedIteratedkTheory},
 as twisted non-abelian cohomology
 with local coefficients in
 $\big( K^{\circ_{2r-2}}(\mathrm{ku})\big)_0$.
 Then the twisted non-abelian character map
 (Def. \ref{TwistedNonAbelianChernDoldCharacter})
 $\mathrm{ch}^{\tau}_{K^{\circ_{2r-2}}(\mathrm{ku})_0}$
 composed with the projection operation \eqref{RationalLSWCharacter}
 onto higher twisted de Rham cohomology,
 (Def. \ref{HigherTwistedAbelianDeRhamCohomology})
 from Lemma \ref{HigherTwistedDeRhamCoefficientsInsideRationalTwistedIteratedKTheory},
  is equivalent to the LSW character map
 $\mathrm{ch}_{2r-1}$ \cite[Def. 2.20]{LindSatiWesterland16}
 restricted along the connective inclusion
 \vspace{-1mm}
 $$
   \overset{
     \mathclap{
     \raisebox{3pt}{
       \tiny
       \color{darkblue}
       \bf
       \def\arraystretch{.9}
       \begin{tabular}{c}
         twisted
         \\
         LSW character
       \end{tabular}
     }
     }
   }{
   \mathrm{ch}_{2r-1}^\tau
   }
   \;\;\;
   \simeq
   \;\;\;\;\;\;
   \underset{
     \mathclap{
     \raisebox{-3pt}{
       \tiny
       \color{darkblue}
       \bf
       \def\arraystretch{.9}
       \begin{tabular}{c}
         projection onto
         \\
         higher twisted
         \\
         de Rham cohomology
       \end{tabular}
     }
     }
   }{
     \phi_\ast
   }
   \;\;\;\;\;\circ\;\;
   \overset{
     \mathclap{
     \raisebox{3pt}{
       \tiny
       \color{darkblue}
       \bf
       \def\arraystretch{.9}
       \begin{tabular}{c}
         twisted non-abelian
         \\
         character map
       \end{tabular}
     }
     }
   }{
     \mathrm{ch}^\tau_{K^{\circ_{2r-2}}(\mathrm{ku})_0}
   }.
 $$
\end{prop}
\begin{proof}
  After unwinding the definitions, the statement
  reduces to the commutativity of the square diagram in
  \cite[p. 15]{LindSatiWesterland16}:
  The top morphism there is the plain rationalization map
  (Def. \ref{RationalizationInTwistedNonAbelianCohomology}),
  the right vertical morphism is $\phi_\ast$ from
  Lemma \ref{HigherTwistedDeRhamCoefficientsInsideRationalTwistedIteratedKTheory}
  before passing from real to de Rham cohomology,
  the left morphism is restriction to the connective part and
  the bottom morphism is the LSW character.
\end{proof}

\medskip

%%%%%%%%%%%%%%%%%%%%%%%%%%%%%%%%%%%%%%%%%%%%%%%%%%%%%%%%
\subsection{Twisted differential non-abelian character}
%%%%%%%%%%%%%%%%%%%%%%%%%%%%%%%%%%%%%%%%%%%%%%%%%%%%%%%%

We introduce twisted differential non-abelian cohomology
(Def. \ref{TwistedDifferentialNonAbelianCohomology} below)
and discuss how the corresponding twisted differential
non-abelian character subsumes existing constructions
on twisted differential K-theory
(Examples \ref{TwistedChernCharacterInTwistedDifferentialKTheory}
and \ref{TwistedDifferentialChernCharacterOnTwistedDifferentialKTheory} below).

\medskip

\noindent {\bf Twisted differential non-abelian cohomology.}
From the perspective of structured non-abelian cohomology
(Remark \ref{StructuredNonAbelianCohomology})
that we have developed,
it is now evident how to canonically combine

{\bf (a)} twisted non-abelian cohomology (Def. \ref{NonabelianTwistedCohomology})
with

{\bf (b)} differential non-abelian cohomology
(Def. \ref{DifferentialNonAbelianCohomology})
to get

%{\bf (c)}
\noindent {\it twisted differential} non-abelian cohomology:

\begin{defn}[Differential non-abelian local coefficient bundles]
  \label{DifferentialNonabelianLocalCoefficientBundles}
Let

\vspace{-.3cm}
$$
    \xymatrix@R=1.5em{
      A \ar[rr]
      &
      \ar@{}[d]|-{
        \mbox{
          \tiny
          \color{darkblue}
          \bf
          \def\arraystretch{.9}
          \begin{tabular}{c}
            local coefficient bundle
          \end{tabular}
        }
      }
      &
      A \!\sslash\! G
      \ar[d]^-{ \rho }
      \\
      && B G
    }
  $$

  \vspace{-1mm}
  \noindent
  be a local coefficient bundle \eqref{LocalCoefficientBundle}
  in $\NilpotentConnectedQFiniteHomotopyTypes$
  (Def. \ref{NilpotentConnectedSpacesOfFiniteRationalType}).

\noindent {\bf (i)}   By
Lemma \ref{MinimalRelativeSullivanModelsPreserveHomotopyFibers},
with Def. \ref{RectifiedRationalizationUnitOnCoefficientBundle},
and using that $\Bexp_{\mathrm{PL}}$ preserves fibrations
(Prop. \ref{FundamentalTheoremForPiecewiseSmoothDeRhamComplex}),
this induces a homotopy fibering (Def. \ref{HomotopyFibers})
in $\Stacks$ (Def. \ref{SmoothInfinityStacks})
of differential non-abelian character maps
(Def. \ref{DifferentialNonabelianCharacterMap})
of this form:
\vspace{-3mm}
\begin{equation}
  \label{FibrationOfDifferentialNonabelianCharacterMaps}
  \hspace{-5mm}
  \raisebox{40pt}{
  \xymatrix@C=3.3em{
    \mathrm{Disc}(A)
    \ar[rr]_-{
      \mathbf{ch}_A
    }^-{
      \mathclap{
      \mbox{
        \tiny
        \color{darkblue}
        \bf
        \def\arraystretch{.9}
        \begin{tabular}{c}
          differential non-abelian character map
          \\
          with coefficients in fiber space
        \end{tabular}
      }
      }
    }
    \ar[dr]_-{
      \mathclap{\phantom{\vert}}
      \mathrm{hofib}
      \left(
        \mathrm{Disc}
        (
          \rho
        )
      \right) \phantom{AA}
    }
    &&
    \flat
    \Bexp
    (
      \mathfrak{l}A
    )
    \ar@{<-}[r]^-{\;\mathrm{atlas}\;}
    \ar[dr]|-{
      \mathclap{\phantom{\vert}}
      \mathrm{hofib}
      \left(
        (
          \mathfrak{l}\rho
        )_\ast
      \right)
    }
    &
    \;\;
    \Omega_{\mathrm{dR}}
    (
      -;
      \,
      \mathfrak{l}A
    )_{\mathrm{flat}}
    \ar[dr]^-{\;\;
      \mathclap{\phantom{\vert}}
      \mathrm{hofib}
      \left(
        (
          \mathfrak{l}\rho
        )_\ast
      \right)
    }
    \\
    &
    \mathrm{Disc}
    \big(
      A \!\sslash\! G
    \big)
    \ar[rr]^-{
      \mathbf{ch}^{B G}_{ A \sslash G } \phantom{AAAA}
    }_-{
      \mathclap{
      \mbox{
        \tiny
        \color{darkblue}
        \bf
        \def\arraystretch{.9}
        \begin{tabular}{c}
          twisted differential non-abelian character map
        \end{tabular}
      }
      }
    }
    \ar[dd]|-{
      \mathclap{\phantom{\vert}}
      \mathrm{Disc}(\rho)
    }
    &&
    \flat
    \Bexp\big(
      \mathfrak{l} (A \!\sslash\! G)
    \big)
    \ar@{<-}[r]_-{ \;\mathrm{atlas}\; }
    \ar[dd]^-{
      \mathclap{\phantom{\vert}}
      \left(
        \mathfrak{l}\rho
      \right)_\ast
    }
    &
    \;
    \Omega_{\mathrm{dR}}
    \big(
      -;
      \,
      \mathfrak{l}_{\scalebox{.6}{$B G$}} (A \!\sslash\! G)
    \big)_{\mathrm{flat}}
    \ar[dd]^-{
      \mathclap{\phantom{\vert}}
      \left(
        \mathfrak{l}\rho
      \right)_\ast
    }
    \\
    \\
    &
    \mathrm{Disc}
    (
      B G
    )
    \ar[rr]^-{
      \mathbf{ch}_{ B G }
    }_-{
      \mathclap{
      \mbox{
        \tiny
        \color{darkblue}
        \bf
        \def\arraystretch{.9}
        \begin{tabular}{c}
          differential non-abelian character map
          \\
          with coefficients in space of twists
        \end{tabular}
      }
      }
    }
    &&
    \flat
    \Bexp(
      \mathfrak{l} B G
    )
    \ar@{<-}[r]_-{ \;\mathrm{atlas}\;  }
    &
    \;
    \Omega_{\mathrm{dR}}
    (
      -;
      \,
      \mathfrak{l} B G
    )_{\mathrm{flat}}
  }
  }
\end{equation}

\vspace{-3mm}
\noindent {\bf (ii)}
 Here the {\it twisted differential non-abelian character map}
$\mathbf{ch}^{B G}_{A \sslash G}$
is defined just as in Def. \ref{DifferentialNonabelianCharacterMap},
but with coefficients the relative Whitehead $L_\infty$-algebra
$\mathfrak{l}_{\scalebox{.6}{$B G$}} (A \!\sslash \! G)$
(Prop. \ref{WhiteheadLInfinityAlgebrasRelative}),
as opposed to the absolute Whitehead $L_\infty$-algebra
$\mathfrak{l} (A \!\sslash\! G)$ (Prop. \ref{WhiteheadLInfinityAlgebras}).
\end{defn}

\begin{remark}[Differential local coefficient bundles]
Since homotopy limits commute over each other,
passage to the homotopy fiber products (Def. \ref{HomotopyPullback})
formed from the horizontal stages of \eqref{FibrationOfDifferentialNonabelianCharacterMaps}
yields a homotopy fibering (Def. \ref{HomotopyFibers}) of
moduli $\infty$-stacks of $\infty$-connections \eqref{PullbackForModuliInfinityStackOfConnections}
of this form:

\vspace{-.6cm}
\begin{equation}
  \label{DifferentialLocalCoefficientBundle}
  \raisebox{70pt}{
  \xymatrix@C=12pt@R=7pt{
    &&
    &&&
    \Omega_{\mathrm{dR}}
    \big(
      -
      ;
      \,
      \mathfrak{l}_{\scalebox{.6}{$B G$}}(A \!\sslash\! B G)
    \big)_{\mathrm{flat}}
    \ar[dr]^{\;\;\;\;\;\; \mathrm{atlas}}
    \ar[ddd]_<<<<<<{
        \mathclap{\phantom{\vert}}
        (\mathfrak{l}\rho)_\ast
      }
      |>>>>>>>>>>>>>{ \phantom{AA} \atop {\phantom{AA}  } }
    \\
    \underset{
      \mathclap{
      \raisebox{-6pt}{
        \tiny
        \color{darkblue}
        \bf
        \def\arraystretch{.9}
        \begin{tabular}{c}
          moduli $\infty$-stack of
          \\
          $\Omega A$-connections
        \end{tabular}
      }
      }
    }{
      A_{\mathrm{diff}}
    }
    \ar[rr]^-{
      \mathrm{hofib}
      (
        \rho_{\mathrm{diff}}
      )
    }
    &
    &
    \big(
      A \!\sslash\! G
    \big)_{\mathrm{diff}_{B G}}
    \ar[ddd]^-{   \rho_{\mathrm{diff}}}_-{
      \mathllap{
      \mbox{
        \tiny
        \color{darkblue}
        \bf
        \def\arraystretch{.9}
        \begin{tabular}{c}
          differential non-abelian
          \\
          local coefficient bundle
        \end{tabular}
      }
      }
%      \rho_{\mathrm{diff}}
    }
    \ar[dr]^-{\!\!\!\!
      \;c^{B G}_{A \sslash G}\;
    }
    \ar[urrr]^-{
      \;F^{B G}_{A \sslash G}\;
    }
    &
    & &&
    \flat
    \Bexp
    \big(
      \mathfrak{l}_{\scalebox{.5}{$B G$}} (A \!\sslash\! G)
    \big)
    \ar[ddd]^<<<<<<<<<<{
      \mathclap{\phantom{\vert}}
      (\mathfrak{l}\rho)_\ast
    }
    \\
    &
    &
    {\phantom{AAAAAAAAAAA}}
    &
    \mathrm{Disc}
    \big(
      A \sslash G
    \big)
    \ar[urrr]
      ^<<<<<<<<<<<<<<<<<{
        \mathclap{\phantom{\vert}}
        \;\mathbf{ch}^{B G}_{ A \sslash G }\;
      }
    \ar[ddd]
      |<<<<<<<{
        \mathclap{\phantom{\vert}}
        \mathrm{Disc}(\rho)
      }
    &&
    \\
    && &
    {\phantom{AAAAAAAAAAA}}
    &&
    \Omega_{\mathrm{dR}}
    \big(
      -;
      \mathfrak{l} B G
    \big)_{\mathrm{flat}}
    \ar[dr]^-{\;\;\;\; \mathrm{atlas} }
    \\
    &&
    \underset{
      \mathllap{
      \raisebox{-6pt}{
        \tiny
        \color{darkblue}
        \bf
        \def\arraystretch{.9}
        \begin{tabular}{c}
          moduli $\infty$-stack of
          \\
          $G$-connections
        \end{tabular}
      }
      \;\;\;\;\;
      }
    }{
      B G_{\mathrm{diff}}
    }
    \ar[urrr]
      |<<<<<<<<<<<<<<<<<<<<<<<<<<{ \phantom{AAA} }
      _>>>>>>>>>>>>>>>{
        \mathclap{\phantom{\vert}}
        \;F_{B G}\;
      }
    \ar[dr]^-{
      \mathclap{\phantom{\vert}}
      \;c_{B G}\;
    }
    &&&&
    \flat
    \Bexp
    \big(
      \mathfrak{l} B G
    \big)
    \\
    &&
    &
    \mathrm{Disc}
    \big(
      B G
    \big)
    \ar[urrr]^-{
      \mathclap{\phantom{\vert}}
      \;\mathbf{ch}_{B G}\;
    }
  }
  }
\end{equation}
\end{remark}

\begin{defn}[Twisted differential non-abelian cohomology]
  \label{TwistedDifferentialNonAbelianCohomology}
  Given a differential non-abelian local coefficient
  bundle $\rho_{\mathrm{diff}}$ \eqref{DifferentialLocalCoefficientBundle}
  according to Def. \ref{DifferentialNonabelianLocalCoefficientBundles},
  we say that:

  \noindent
  {\bf (i)} A \emph{differential twist}
  on a $\mathcal{X} \in \Stacks$ (Def. \ref{SmoothInfinityStacks})
  is a cocycle
  $\tau_{\mathrm{diff}}$ in differential non-abelian
  cohomology with coefficients in $B G$ (Def. \ref{DifferentialNonAbelianCohomology})
  \vspace{0mm}
  \begin{equation}
    \label{DifferentialTwist}
    \big[
      \tau_{\mathrm{diff}}
    \big]
    \;\in\;
    \widehat H
    \big(
      \mathcal{X};
      \,
      B G
    \big).
  \end{equation}

    \vspace{-1mm}
  \noindent
  {\bf (ii)} The
  {\it $\tau_{\mathrm{diff}}$-twisted differential non-abelian
  cohomology} with local coefficients in $\rho_{\mathrm{diff}}$
  is the structured (Remark \ref{StructuredNonAbelianCohomology})
  $\tau_{\mathrm{diff}}$-twisted non-abelian cohomology
  (Def. \ref{NonabelianTwistedCohomology}) with coefficients in
  $\rho_\mathrm{diff}$, hence the hom-set in
  the homotopy category (Def. \ref{HomotopyCategory})
  of the slice model structure (Def. \ref{SliceModelCategory})
  of the local projective model structure $\mathrm{SmoothStacks}_\infty$
  on simplicial presheaves
  over $\CartesianSpaces$
  (Example \ref{ModelStructureOnSimplicialPresheavesOverCartesianSpaces})
  from $\tau_{\mathrm{diff}}$ \eqref{DifferentialTwist}
  to $\rho_{\mathrm{diff}}$ \eqref{DifferentialLocalCoefficientBundle}:
    \vspace{-2mm}
  \begin{equation}
    \label{TwistedDifferentialNonabelianCohomology}
    \hspace{-1.8cm}
    \overset{
      \mathclap{
      \raisebox{3pt}{
        \tiny
        \color{darkblue}
        \bf
        \def\arraystretch{.9}
        \begin{tabular}{c}
          twisted differential
          \\
          non-abelian cohomology
        \end{tabular}
      }
      }
    }{
      \widehat H^{\tau_{\mathrm{diff}}}
      \big(
        \mathcal{X};
        \,
        A
      \big)
    }
    \;
    :=
    \;
    \StacksTwisted
    (
      \tau_{\mathrm{diff}}
      \,,\,
      \rho_{\mathrm{diff}}
    )
        =
       \left\{\!\!\!\!\!\!\!\!
    \raisebox{18pt}{
    \xymatrix@C=20pt{
      \;\;\;\mathcal{X}\;\;
      \ar@{-->}[rr]^-{
        \overset{
          \mathclap{
          \raisebox{3pt}{
            \tiny
            \color{darkblue}
            \bf
            \def\arraystretch{.9}
            \begin{tabular}{c}
              differential
              cocycle
            \end{tabular}
          }
          }
        }{
          c_{\mathrm{diff}}
        }
      }_>>>>>{\ }="s"
      \ar[dr]_-{
        \underset{
          \mathllap{
          \raisebox{-3pt}{
          \mbox{
            \tiny
            \color{darkblue}
            \bf
            \def\arraystretch{.9}
            \begin{tabular}{c}
              differential
              \\
              twist
            \end{tabular}
          }
          }
          }
        }{
          \tau_{\mathrm{diff}}
        }
      }^-{\ }="t"
      &&
      (A \!\sslash\! G)_{\mathrm{diff}_{B G}}
      \ar[dl]^-{
        \underset{
          \mathrlap{
          \raisebox{-2pt}{
            \tiny
            \color{darkblue}
            \bf
            \def\arraystretch{.9}
            \begin{tabular}{c}
              differential local
              \\
              coefficients
            \end{tabular}
          }
          }
        }{
          \rho_{\mathrm{diff}}
        }
      }
      \\
      &
      B G_{\mathrm{diff}}
      \ar@{=>}^-{\simeq} "s"; "t"
    }
    }
    \right\}_{
      \mathrlap{
      \!\!
      \big/
      \!\!\!\!\!\!\!\!\!\!\!
      \mbox{
        \tiny
        \def\arraystretch{.9}
        \begin{tabular}{c}
          homotopy
          \\
          relative $B G_{\mathrm{diff}}$
        \end{tabular}
      }
      }
    }
  \end{equation}

  \vspace{-2mm}
  \noindent
  {\bf (iii)} The twisted non-abelian cohomology operations
  induced from the maps in \eqref{DifferentialLocalCoefficientBundle}
  we call
  (see \eqref{SystemsOfCohomologyOperationsOnDifferentialCohomology}):

  \vspace{-.6cm}
  \begin{align}
    \label{CharacteristicClassOnTwistedDifferentialCohomology}
      \mbox{ {\bf (a)} {\it characteristic class}: \phantom{aa} }
    &
    \xymatrix@C=36pt{
      \widehat H^{\tau_{\mathrm{diff}}}
      \big(
        \mathcal{X};
        \,
        A
      \big)
      \ar[rr]^-{
        c^\tau_A
        \,:=\,
        \big(
          c^{B G}_{A \sslash G}
        \big)_\ast
      }
      &&
      H^\tau
      \big(
        \mathrm{Shp}(\mathcal{X});
        \,
        A
      \big)
    }
    &
    \mathllap{
      \mbox{(Def. \ref{NonabelianTwistedCohomology})}
    }
    \\
    \label{CurvatureOnTwistedDifferentialCohomology}
    \mbox{ {\bf (b)} {\it curvature:} \phantom{characterisa} }
    &
    \xymatrix@C=36pt{
      \widehat H^{\tau_{\mathrm{diff}}}
      \big(
        \mathcal{X};
        \,
        A
      \big)
      \ar[rr]^-{
        F_A^{\tau_{\mathrm{dR}}}
        \,:=\,
        \big(
          F^{B G}_{A \sslash G}
        \big)_\ast }
      &&
      \Omega^{\tau_{\mathrm{dR}}}_{\mathrm{dR}}
      \big(
        \mathcal{X};
        \,
        \mathfrak{l}A
      \big)_{\mathrm{flat}}
    }
    &
    \mathllap{
      \mbox{(Def. \ref{FlatTwistedLInfinityAlgebraValuedDifferentialForms})}
    }
     \\
    \label{DifferentialCharacterOnTwistedDifferentialCohomology}
    \mbox{ {\bf (c)} {\it differential character:}\phantom{a}\, }
    &
    \xymatrix@C=40pt{
      \widehat H^{\tau_{\mathrm{diff}}}
      \big(
        \mathcal{X};
        \,
        A
      \big)
      \ar[rr]^-{
        \mathrm{ch}^{\tau}_A
        \,:=\,
        \big(
          \mathbf{ch}^{B G}_{A \sslash G}
          \,\circ\,
          c^{B G}_{A \sslash G}
        \big)_\ast }
      &&
      H^{\tau_{\mathrm{dR}}}_{\mathrm{dR}}
      \big(
        \mathcal{X};
        \,
        \mathfrak{l}A
      \big)
    }
    &
    \mathllap{
      \mbox{
        (Def. \ref{TwistedNonabelianDeRhamCohomology})}
    }
  \end{align}
\end{defn}

\medskip

\noindent
{\bf Twisted differential non-abelian cohomology as non-abelian $\infty$-sheaf hypercohomology.}
While the formulation
of twisted differential non-abelian cohomology
as hom-sets in a slice of $\mathrm{SmoothStacks}_\infty$
(Def. \ref{TwistedDifferentialNonAbelianCohomology})
is natural and useful,
we indicate how
this is equivalently incarnated as a non-abelian sheaf
hypercohomology over $\mathcal{X}$.
This serves to make the connection to existing literature
(in Example \ref{TwistedDifferentialGeneralizedCohomology} below),
but is not otherwise needed for the development here.
We shall be brief, referring to \cite{SS20b} for some technical
background that is beyond the scope of our presentation here.

\begin{prop}[{\'E}tale $\infty$-topos over $\infty$-stacks {\cite[Prop. 3.33, Rem. 3.34]{SS20b}}]
\label{EtaleInfinityToposOverInfinityStack}
For $\mathcal{X} \in \Stacks$ (Def. \ref{SmoothInfinityStacks})
let
\vspace{-2mm}
$$
  \xymatrix{
    \StacksEtaleOverX
    \;
    \ar@{^{(}->}[rr]^-{ \LeftDerived i_{\mathcal{X}} }
    &&
    \;
    \StacksOverX
  }
$$

\vspace{-2mm}
\noindent
be the full subcategory of the homotopy category
(Def. \ref{HomotopyCategory}) of the slice model structure
over $\mathcal{X}$ (Example \ref{SliceModelCategory}) of the
local projective model structure on simplicial presheaves
(Example \ref{ModelStructureOnSimplicialPresheavesOverCartesianSpaces})
on those $\mathcal{E} \to \mathcal{X}$
which are local diffeomorphisms (\cite[Def. 3.26]{SS20b}).

\noindent
{\bf (i)} The inclusion $\LeftDerived i_{\mathcal{X}}$ is a left-exact homotopy co-reflection,
in that it preserves finite homotopy limits and has
a derived right adjoint $\mathbb{R} \mathrm{Et}$
(sending $\infty$-bundles to their $\infty$-sheaves of $\infty$-sections).

\noindent
{\bf (ii)} There is a \emph{global section} functor
$\mathbb{R}\Gamma_{\mathcal{X}}$
from $\StacksEtaleOverX$ to $\HomotopyTypes$
(Example \ref{TheClassicalHomotopyCategory})
which also admits a left exact left adjoint:

\vspace{-6mm}
\begin{equation}
  \label{DerivedGlobalSections}
  \xymatrix@C=3.5em{
    \overset{
      \mathclap{
      \raisebox{3pt}{
        \tiny
        \color{darkblue}
        \bf
        \def\arraystretch{.9}
        \begin{tabular}{c}
          $\infty$-bundles
          over $\mathcal{X}$
        \end{tabular}
      }
      }
    }{
      \StacksOverX
    }
    \;
    \ar@{<-^{)}}@<+6pt>[rr]^-{ \LeftDerived i_{\mathcal{X}} }
    \ar@<-6pt>[rr]_-{
      \underset{}{
        \underset{
          \mathclap{
          \raisebox{-3pt}{
            \tiny
            \color{darkblue}
            \bf
            \def\arraystretch{.9}
            \begin{tabular}{c}
              $\infty$-sheaf
              of local sections
            \end{tabular}
          }
          }
        }{
          \mathbb{R}\mathrm{Et}
        }
      }
    }^-{ \bot }
    &&
    \;\;
    \overset{
      \mathclap{
      \raisebox{3pt}{
        \tiny
        \color{darkblue}
        \bf
        \def\arraystretch{.9}
        \begin{tabular}{c}
          $\infty$-sheaves
          over $\mathcal{X}$
        \end{tabular}
      }
      }
    }{
      \StacksEtaleOverX
    }
    \ar@{<-}@<+6pt>[rr]^-{ \Delta_{{}_{\mathcal{X}}} }
    \ar@<-6pt>[rr]_-{
      \underset{
        \mathclap{
        \raisebox{-3pt}{
          \tiny
          \color{darkblue}
          \bf
          \def\arraystretch{.9}
          \begin{tabular}{c}
            global sections
          \end{tabular}
        }
        }
      }{
        \mathbb{R}\Gamma_{\!\!{}_{\mathcal{X}}}
      }
    }^-{ \bot }
    \;
    &&
  \;  \HomotopyTypes
  }\!.
\end{equation}
\end{prop}

\vspace{-5mm}
\begin{defn}[Non-abelian $\infty$-sheaf hypercohomology over $\infty$-stacks]
  \label{InfinitySheafHypercohomology}
  Given
  $\mathcal{X} \,\in\, \Stacks$ (Def. \ref{SmoothInfinityStacks})
  and
  $\mathcal{A} \,\in\, \StacksEtaleOverX$
  (Prop. \ref{EtaleInfinityToposOverInfinityStack})
  we say that the set of connected components
  of the \emph{derived global sections} \eqref{DerivedGlobalSections}
  of $\mathcal{A}$ over $\mathcal{X}$
  $$
    H
    \big(
      \mathcal{X}
      \,,\,
      \mathcal{A}
    \big)
    \;\;
    :=
    \;\;
    \pi_0
    \big(
      \mathbb{R}\Gamma_{{}_{\!\!\mathcal{X}}}\!(\mathcal{A})
    \big)
  $$
  is the
  {\it non-abelian $\infty$-sheaf hypercohomology} of
  $\mathcal{X}$ with coefficients in $\mathcal{A}$.
\end{defn}

\begin{lemma}[Twisted differential non-abelian cohomology as
non-abelian $\infty$-sheaf hyper-cohomology]
  \label{TwistedDifferentialNonAbelianCohomologyAsSheafHypercohomology}
Given a differential twist $\tau_{\mathrm{diff}}$
\eqref{DifferentialTwist} on some
$\mathcal{X} \in \Stacks$ \eqref{HomotopyCategoryOfSimplicialPresheaves}
consider the object
\vspace{-2mm}
\begin{equation}
  \label{InfinitySheafForDifferentialCohomology}
  \underline{A}_{\tau_{\mathrm{diff}}}
  \;:=\;
  \mathbb{R}\mathrm{LcllCnstnt}_{\mathcal{X}}
  \big(
    \mathbb{R}\tau_{\mathrm{diff}}^\ast
    (A \!\sslash\! G)_{\mathrm{diff}}
  \big)
  \;\;\;\;
  \in
  \;
  \StacksEtaleOverX
\end{equation}

\vspace{-2mm}
\noindent
in the {\'e}tale $\infty$-topos over $\mathcal{X}$
Prop. \ref{EtaleInfinityToposOverInfinityStack}.
The
non-abelian
$\infty$-sheaf hypercohomology (Def. \ref{InfinitySheafHypercohomology})
of $\underline{A}_{\tau_{\mathrm{diff}}}$ over $\mathcal{X}$
coincides with the
$\tau_{\mathrm{diff}}$-twisted differential non-abelian
cohomology of $\mathcal{X}$ (Def. \ref{TwistedDifferentialNonAbelianCohomology}):
\vspace{-2mm}
\begin{equation}
  \overset{
    \mathclap{
    \raisebox{3pt}{
      \tiny
      \color{darkblue}
      \bf
      \def\arraystretch{.9}
      \begin{tabular}{c}
        non-abelian
        \\
        $\infty$-sheaf hypercohomology
      \end{tabular}
    }
    }
  }{
    \pi_0
    \,
    \mathbb{R}\Gamma_{{}_{\!\!\mathcal{X}}}
    \!
    \big(
      \underline{A}_{\tau_{\mathrm{diff}}}
    \big)
  }
    \;\;\;
    \simeq
    \;\;
    \overset{
      \mathclap{
      \raisebox{3pt}{
        \tiny
        \color{darkblue}
        \bf
        \def\arraystretch{.9}
        \begin{tabular}{c}
          twisted differential
          \\
          non-abelian cohomology
        \end{tabular}
      }
      }
    }{
      \widehat H^{\tau_{\mathrm{diff}}}
      \big(
        \mathcal{X}
        \,,\,
        A
      \big).
    }
\end{equation}
\end{lemma}
\begin{proof}
  As in \cite[Remark 3.34]{SS20b}.
\end{proof}

It is useful to decompose this construction
of twisted differential cohomology via
$\infty$-sheaf hypercohomology again as a homotopy pullback
of corresponding $\infty$-sheaves representing plain twisted cohomology
and plain twisted differential forms:

\begin{remark}[Homotopy pullback of $\infty$-sheaves representing twisted differential cohomology]
Given a differential twist $\tau_{\mathrm{diff}}$
\eqref{DifferentialTwist} on some
$\mathcal{X} \in \Stacks$ \eqref{HomotopyCategoryOfSimplicialPresheaves}
with components
$(\tau, \tau_{\mathrm{dR}}, L_{{}_{\mathbb{R}}}\tau)$
(Example \ref{HomotopyPullbackViaTriples}),

\vspace{1mm}
\noindent {\bf (i)} Consider the pullback stacks over $\mathcal{X}$
in the following diagram
\vspace{-1mm}
$$
  \hspace{-.4cm}
  \scalebox{.85}{
  \xymatrix@R=10pt@C=-60pt{
    \mathbb{R}\tau^\ast
    \,
      \mathrm{Disc}
      \big(
        A \!\sslash\! G
      \big)
    \,
    \ar[dddd]
    \ar[rrrr]
    \ar[dddr]
    &&&&
    \mathrm{Disc}
    \big(
      A \!\sslash\! G
    \big)
    \ar[dddr]
    \ar[dddd]^-{ \rho }|>>>>>>>>{ \phantom{AA} \atop \phantom{AA} }
    \\
    \\
    \\
    &
    \mathbb{R}(L_{{}_{\mathbb{R}}}\tau)^\ast
    \flat
    \Bexp
    \big(
      \mathfrak{l}_{\scalebox{.7}{$BG$}} (A \!\sslash\! G)
    \big)
    \ar[dddr]
    \ar[dddd]
    \ar[rrrr]
    &&&&
    \flat
    \Bexp
    \big(
      \mathfrak{l}_{\scalebox{.7}{$BG$}} (A \!\sslash\! G)
    \big)
    \ar[dddr]
    \ar[dddd]|>>>>>>>{ \phantom{AA} \atop \phantom{AA} }
    \\
    \mathcal{X}
    \ar@{=}[dddr]
    \ar[rrrr]_>>>>>>>>>>>>>>{ \;\tau\; }
      |<<<<<<<<<<<<<<<<<<<<{ \phantom{AA} }
      |<<<<<<<<<<<<<<<<<<<<<<<<<<<<<{ \phantom{AAA} }
    &&&&
    \;\;\;
    \mathrm{Disc}
    (
      B G
    )
    \;\;\;
    \ar[dddr]|>>>>>>>>>>{ \phantom{AA} \atop \phantom{AA} }
    \\
    \\
    {\phantom{AAAAAAAAAAAAAAAAAA}}
    &&
    \mathbb{R}\tau_{\mathrm{dR}}^\ast
    \Omega_{\mathrm{dR}}
    \big(
      -;
      \,
      \mathfrak{l}_{\scalebox{.7}{$B G$}}(A \!\sslash\! G)
    \big)_{\mathrm{flat}}
    \ar[rrrr]
    \ar[dddd]
    &&&&
    \Omega_{\mathrm{dR}}
    \big(
      -;
      \,
      \mathfrak{l}_{\scalebox{.7}{$B G$}}(A \!\sslash\! G)
    \big)_{\mathrm{flat}}
    \ar[dddd]
    \\
    &
    \mathcal{X}
    \ar@{=}[dddr]
    \ar[rrrr]_>>>>>>>>>>>>>>>>>>{ \;L_{{}_{\mathbb{R}}}\tau\; }
      |<<<<<<<<<<<<<<<<<<<<{ \phantom{AAA} }
    &&&&
    \;\;
    \flatBexp
    \big(
      \mathfrak{l}B G
    \big)
    \;\;
    \ar[dddr]
    \\
    \\
    \\
    & &
    \mathcal{X}
    \ar[rrrr]^-{
      \;\tau_{\mathrm{dR}}\;
    }
    &
    &
    &
    &
    \Omega_{\mathrm{dR}}
    (
      -;
      \mathfrak{l} B G
    )
    \\
    {\phantom{AAAAAAAAAAAAAAAAAA}}
    &
    {\phantom{AAAAAAAAAAAAAAAAAA}}
    &
    {\phantom{AAAAAAAAAAAAAAAAAA}}
    &
    {\phantom{AAAAAAAAAAAAAAAAAA}}
    &
    {\phantom{AAAAAAAAAAAAAAAAAA}}
    &
    {\phantom{AAAAAAAAAAAAAAAAAA}}
    &
    {\phantom{AAAAAAAAAAAAAAAAAA}}
  }
  }
$$

\vspace{-5mm}
\noindent Here the right hand side is \eqref{FibrationOfDifferentialNonabelianCharacterMaps}
and all front-facing squares are homotopy pullbacks
(Def. \ref{HomotopyPullback}).

\vspace{1mm}
\noindent {\bf (ii)} By commutativity of homotopy limits over each other,
these form a homotopy pullback square as
on the right of the following diagram,
which gives, under
the derived right adjoint $\mathbb{R}\mathrm{LcllCnstnt}$
\eqref{DerivedGlobalSections}
a homotopy pullback diagram
of $\infty$-sheaves of sections as shown on the left:
\vspace{-2mm}
\begin{equation}
  \label{HomotopyPullbackDiagramOfInfinitySheavesForTwistedDifferentialCohomology}
  \hspace{-3mm}
  \raisebox{24pt}{
  \xymatrix@C=10pt@R=18pt{
    \underline{A}_{\tau_{\mathrm{diff}}}
    \ar[d]
    \ar[rr]
    \ar@{}[drr]|-{
      \mbox{\tiny(hpb)}
    }
    &&
    \underline{
      \Omega
      \big(
        -;
        \,
        \mathfrak{l}A
      \big)_{\mathrm{flat}}
    }_{\, \tau_{\mathrm{dR}}}
    \ar[d]
    \\
    \underline{
      A
    }_{\tau}
    \ar[rr]
    &&
    \underline{
      \flatBexp
      \big(
        \mathfrak{l}A
      \big)
    }_{\, L_{{}_{\mathbb{R}}}\!\tau}
  }
  }
  \!\!\!
  :=
  \;
   \mathbb{R}\mathrm{LcllCnstst}
  \left(
  \!\!\!
  \raisebox{24pt}{
  \xymatrix@C=6pt{
    \mathbb{R}\tau_{\mathrm{diff}}^\ast
    (A \!\sslash\! G)_{\mathrm{diff}}
    \ar[d]
    \ar[rr]
    \ar@{}[drr]|-{ \mbox{\tiny\rm(hpb)} }
    &&
    \mathbb{R}\tau_{\mathrm{dR}}^\ast
    \Omega_{\mathrm{dR}}
    \big(
      -;
      \mathfrak{l}_{\scalebox{.7}{$B G$}}(A \!\sslash\! G)
    \big)
    \ar[d]
    \\
    \mathbb{R}\tau^\ast (A \!\sslash\! G)
    \ar[rr]
    &&
    \mathbb{R}(L_{{}_{\mathbb{R}}}\tau)^\ast
    \flat
    \Bexp
    \big(
      \mathfrak{l}( A \!\sslash\! G)
    \big)
  }
  }
  \!\!\!
  \right)
    \in
  \StacksEtaleOverX .
\end{equation}

\vspace{-2mm}
\noindent Here the top left item $\underline{A}_{\tau_{\mathrm{diff}}}$ from \eqref{InfinitySheafForDifferentialCohomology}
is the $\infty$-sheaf
whose global sections give the $\tau_{\mathrm{diff}}$-twisted
differential cohomology, by Lemma \ref{TwistedDifferentialNonAbelianCohomologyAsSheafHypercohomology}.
\end{remark}

In  differential enhancement of Prop. \ref{ProofTwistedGeneralizedCohomology}
and in twisted enhancement of Example \ref{DifferentialGeneralizedCohomology},
we have:

\begin{example}[Twisted differential generalized cohomology]
\label{TwistedDifferentialGeneralizedCohomology}
Let $\mathcal{X} \,=\, X$ be a smooth manifold
(Ex. \ref{SmoothManifoldAsInfinityStack}),
$R$ an $E_\infty$-ring spectrum (Ex. \ref{GeneralizedCohomologyAsNonabelianCohomology}), and let
\vspace{-1mm}
$$
  \xymatrix{
    E_0
    \ar@{}[drr]|-{\mbox{\tiny(hpb)}}
    \ar[rr]
    \ar[d]_-{ \mathbb{R}\tau^\ast \rho_R }
    &&
    (R_0) \!\sslash\! \GroupOfUnits{R}
    \ar[d]^-{ \rho_R }
    \\
    X
    \ar[rr]_-{ \tau }
    &&
    B \GroupOfUnits{R}
  }
$$

\vspace{-1mm}
\noindent
be a twist for twisted generalized $R$-cohomology over $X$
\eqref{CoefficientBundleAssociatedToPrincipalBundle},
as in Lemma \ref{ProofTwistedGeneralizedCohomology}.

\noindent {\bf (i)} Then the corresponding homotopy pullback diagram
\eqref{HomotopyPullbackDiagramOfInfinitySheavesForTwistedDifferentialCohomology},
which
exhibits, by Lemma \ref{TwistedDifferentialNonAbelianCohomologyAsSheafHypercohomology},
twisted differential non-abelian cohomology
(Def. \ref{TwistedDifferentialNonAbelianCohomology})
with coefficients in $E_0$ as $\infty$-sheaf hypercohomology
(Def. \ref{InfinitySheafHypercohomology}),
is the image under $\mathbb{R}\Omega^\infty_X$
of the homotopy pullback diagram of
sheaves of spectra considered in
\cite[Def. 4.11]{BunkeNikolaus14}, shown on the right below,
for canonical/minimal differential refinement as in
Example \ref{DifferentialGeneralizedCohomology}:
\vspace{-2mm}
$$
  \raisebox{27pt}{
  \xymatrix@C=2em{
    \underline{R_0}_{\, \tau_{\mathrm{diff}}}
    \ar[d]
    \ar[rr]
    \ar@{}[drr]|-{
      \mbox{\tiny(hpb)}
    }
    &&
    \underline{
      \Omega
      \big(
        -;
        \,
        \mathfrak{l}R_0
      \big)_{\mathrm{flat}}
    }_{\, \tau_{\mathrm{dR}}}
    \ar[d]
    \\
    \underline{
      R_0
    }_{\, \tau}
    \ar[rr]
    &&
    \underline{
      \flatBexp
      \big(
        \mathfrak{l}R_0
      \big)
    }_{\, L_{{}_{\mathbb{R}}}\!\tau}
  }
  }
  \;\;\;\;
  \simeq
  \;\;\;\;
  \mathbb{R}\Omega_{{}_X}\infty
  \left(
    \!
    \raisebox{26pt}{
    \xymatrix@R=32pt{
      \mathrm{Diff}
      \big(
        E
      \big)
      \ar@{}[drr]|-{ \mbox{\tiny(hpb)} }
      \ar[d]
      \ar[rr]
      &&
      H \mathcal{M}_{\leq 0}
      \ar[d]
      \\
      \mathrm{Disc}(E)
      \ar[rr]
      &&
      H \mathcal{M}
    }
    }
    \!
  \right)
$$

\vspace{-2mm}
\noindent
This is the twisted/parametrized analog of
the relation \eqref{RelatingModuliStacksOfConnectionsToDifferentialFunctionSpectra}.

\noindent {\bf (ii)} Accordingly, the twisted differential generalized $R$-cohomology
according to \cite[Def. 4.13]{BunkeNikolaus14} is subsumed
by twisted differential non-abelian cohomology, via Lemma \ref{TwistedDifferentialNonAbelianCohomologyAsSheafHypercohomology}.
\end{example}

\medskip

In differential enhancement of Prop. \ref{TwistedChernCharacterInTwistedTopologicalKTheory}
and in twisted generalization of Example \ref{CurvatureInDifferentialComplexKTheory},
we have:
\begin{example}[Twisted Chern character in twisted differential K-theory]
  \label{TwistedChernCharacterInTwistedDifferentialKTheory}
  Consider again the local coefficient bundle
  $$
    \xymatrix@R=12pt{
      \mathrm{KU}_0
      \ar[rr]
      &&
      \mathrm{KU}_0 \!\sslash\! B \mathrm{U}(1)
      \ar[d]^-{ \rho }
      \\
      &&
      B^2 \mathrm{U}(1)
    }
  $$
  for complex topological K-theory (Example \ref{TwistedKTheory}).
  By Example \ref{TwistedDifferentialGeneralizedCohomology},
  the twisted differential non-abelian cohomology theory
  (Def. \ref{TwistedDifferentialNonAbelianCohomology})
  induced from these local coefficients is
  twisted differential K-theory, as discussed in
  \cite{CareyMickelssonWang09} for torsion twists (review in \cite[\S 7]{BunkeSchick12}).
  By the diagram \eqref{SystemsOfCohomologyOperationsOnDifferentialCohomology} of cohomology operations
  on twisted differential cohomology,
  one may regard the corresponding twisted curvature map
  \eqref{CurvatureOnTwistedDifferentialCohomology}
  \vspace{-2mm}
  $$
    \xymatrix{
      \widehat K^{\mathrm{\tau}_{\mathrm{diff}}}
      \big(
        \mathcal{X}
      \big)
      \ar[rr]^-{
        \big(
          F^{\tau_{\mathrm{dR}}}_{\mathrm{KU}_0}
        \big)_\ast
      }
      &&
      \Omega^{\tau_{\mathrm{dR}}}_{\mathrm{dR}}
      \big(
        \mathcal{X}
        ;
        \,
        \mathfrak{l}\mathrm{KU}_0
      \big)_{\mathrm{flat}}
    }
  $$
  (with values in flat $\tau_{\mathrm{dR}} \simeq H_3$-twisted differential forms,
  by Example \ref{RecoveringH3TwistedDifferentialForms})
  as an incarnation of the Chern character map on twisted differential
  K-theory.
  Unwinding this abstract construction
  produces the perspective taken in
  \cite[p. 2]{CareyMickelssonWang09}\cite{Park18} for torsion twists,
  and in
  \cite[p. 6]{BunkeNikolaus14} for general twists.
\end{example}
  However, in the spirit of the Cheeger-Simons homomorphism
  \eqref{NonabelianDifferentialCohomology}, any
  lift of a cohomology operation
  (here: rationalization)
  to differential cohomology should be enhanced all the
  way to a secondary cohomology operation
  (Def. \ref{SecondaryNonabelianCohomologyOperations},
  now to be generalized to a
  \emph{twisted} secondary cohomology operation,
  Def. \ref{TwistedSecondaryNonAbelianCohomologyOperations} below)
  whose codomain is itself a (twisted) differential cohomology theory.
  The twisted Chern character enhanced to a secondary
  cohomology operation this way is
  Example \ref{TwistedDifferentialChernCharacterOnTwistedDifferentialKTheory}
  below, following the perspective taken in \cite[\S 3.2]{GS-tAHSS}\cite[\S 2.3]{GS-RR}.

\medskip

\noindent
{\bf Secondary twisted non-abelian cohomology operations.}
We introduce the twisted generalization of secondary non-abelian
cohomology operations (Def. \ref{TwistedSecondaryNonAbelianCohomologyOperations} below).
This requires the
following twisted analog of the technical condition
in Def. \ref{RelativeFormalMaps}:

\begin{defn}[Twisted absolute minimal model]
\label{TwistedAbsoluteMinimalModel}
For
\vspace{-2mm}
$$
  \raisebox{20pt}{
  \xymatrix{
    A_1 \!\sslash\! G_1
    \ar[d]_{ \rho_1 }
    \ar[rr]^-{ c_t }
    &&
    A_2 \!\sslash\! G_2
    \ar[d]_{ \rho_2 }
    \\
    B G_1
    \ar[rr]_-{ c_b }
    &&
    B G_2
  }
  }
  \;\;\;\;\;
  \in
  \;\;
  \SimplicialSets
$$

\vspace{-2mm}
\noindent a transformation \eqref{MorphismOfLocalCoefficientBundles}
between local coefficient bundles \eqref{LocalCoefficientBundle},
and for $\mathfrak{c}_b$ an absolute minimal model
(Def. \ref{RelativeFormalMaps}) of the map $c_b$
between spaces of twists, hence with induced
transformation \eqref{TransformationOfDifferentialNonaebalianCharacters}
\vspace{-2mm}
$$
  \xymatrix@R=4pt@C=2pt{
    \mathrm{Disc}\big(B G_1\big)
    \ar[rrrr]^-{ \;\mathrm{Disc}(c_b)\; }
    \ar[dddr]^-{
      \mathclap{\phantom{\vert}}
      \mathbf{ch}_{B G_1}
    }
    &&&&
    \mathrm{Disc}\big(B G_2\big)
    \ar[dddr]^-{
      \mathclap{\phantom{\vert}}
      \mathbf{ch}_{B G_2}
    }
    \\
    \\
    \\
    &
    \flat
    \Bexp
    \big(
      \mathfrak{l} B G_1
    \big)
    \ar[rrrr]_-{ (\mathfrak{c}_b)_\ast  }
    &&&&
    \flat
    \Bexp
    \big(
      \mathfrak{l} B G_1
    \big)
  }
$$

 \vspace{-2mm}
\noindent
between the differential character maps
(Def. \ref{DifferentialNonabelianCharacterMap}) on the spaces of twists,
we say that a corresponding \emph{twisted absolute minimal model}
is a lift of $\mathfrak{c}_b$ to a morphism
\vspace{-2mm}
\begin{equation}
\label{MorphismForTwistedAbsoluteMinimalModel}
\xymatrix{
  \mathfrak{l}_{\scalebox{.6}{$B G_1$}} (A_1 \!\sslash\! G_1)
  \ar@{-->}[rr]^-{
    \mathfrak{c}_t
  }
  &&
  \mathfrak{l}_{\scalebox{.6}{$B G_1$}} (A_1 \!\sslash\! G_1)
}
\end{equation}

\vspace{-2mm}
\noindent
between the relative rational Whitehead $L_\infty$-algebras
of the local coefficient bundles
(Prop. \ref{WhiteheadLInfinityAlgebrasRelative})
which

\noindent {\bf (i)} yields a transformation
\vspace{-4mm}
$$
  \xymatrix@R=6pt@C=2pt{
    \mathrm{Disc}\big( A_1 \!\sslash\! G_1 \big)
    \ar[rrrr]^-{ \;\mathrm{Disc}(c_t)\; }
    \ar[dddr]^-{
      \mathbf{ch}^{ B G_1 }_{ A_1 \sslash G_1}
    }
    &&&&
    \mathrm{Disc}\big( A_2 \!\sslash\! G_2 \big)
    \ar[dddr]^-{
         \mathbf{ch}^{ B G_2 }_{ A_2 \sslash G_2 }
      \mathclap{\phantom{\vert}}
    }
    \\
    \\
    \\
    &
    \flat
    \Bexp
    \big(
      \mathfrak{l}_{\scalebox{.6}{$B G_1$}} (A_1 \!\sslash\! G_1)
    \big)
    \ar@{-->}[rrrr]_-{ (\mathfrak{c}_t)_\ast  }
    &&&&
    \flat
    \Bexp
    \big(
      \mathfrak{l}_{\scalebox{.6}{$B G_2$}} (A_2 \!\sslash\! G_2)
    \big)
  }
$$

\vspace{-2mm}
\noindent
of the twisted differential characters
\eqref{FibrationOfDifferentialNonabelianCharacterMaps}
(thus being an ``absolute minimal model for $c_t$ relative to $c_b$''),

\noindent
{\bf (ii)} is compatible with the transformation of the
differential characters on the twisting space, in that the following
cube commutes:

\vspace{-6mm}
\begin{equation}
  \label{TransformationBetweenTwistedDifferentialNonAbelianCharacterMaps}
  \hspace{-6mm}
  \scalebox{.9}{
  \raisebox{40pt}{
  \xymatrix@R=14pt@C=-6.5pt{
    \mathrm{Disc}
    \big(
      A_1 \!\sslash\! G_1
    \big)
    \ar[rrrr]|-{ \;\mathrm{Disc}(c_t)\; }
    \ar[dddd]|-{
      \mathclap{\phantom{\vert}}
      \rho_1
    }
    \ar[dddr]|-{ \mathbf{ch}^{B G_1}_{A_1 \sslash G_1} }
    &&&&
    \mathrm{Disc}
    \big(
      A_2 \!\sslash\! G_2
    \big)
    \ar[dddd]
      |-{
        \mathclap{\phantom{\vert}}
        \rho_2
      }
      |>>>>>>>>>{ \phantom{AAA} }
    \ar[dddr]|-{ \mathbf{ch}^{B G_2}_{A_2 \sslash G_2} }
    \\
    \\
    \\
    &
    \flat
    \Bexp
    \big(
      \mathfrak{l}_{\scalebox{.6}{$B G_1$}} (A_1 \!\sslash\! G_1)
    \big)
    \ar[dddd]|-{ \mathfrak{l}_{\scalebox{.6}{$B G_1$}} (A_1 \!\sslash\! G_1) }
    \ar@{<-}[dddr]|>>>>>>>>>>{
      \mathclap{\phantom{\vert}}
      \mathrm{atlas}
    }
    \ar[rrrr]|<<<<<<<<<<<<<<<<{ \;(\mathfrak{c}_t)_\ast\; }
    &&&&
    \flat
    \Bexp
    \big(
      \mathfrak{l}_{\scalebox{.6}{$B G_1$}} (A_2 \!\sslash\! G_2)
    \big)
    \ar@{<-}[dddr]|-{
      \mathclap{\phantom{\vert}}
      \mathrm{atlas}
    }
    \ar[dddd]
      |-{ (\mathfrak{l}_{\scalebox{.6}{$B G_2$}}p_2)_{\ast} }
      |>>>>>>>>>{ \phantom{AA} }
    \\
    \mathrm{Disc}
    \big(
      B G_1
    \big)
    \ar[rrrr]
      |<<<<<<<<<<<<<<<<<<<{ \phantom{AA} }
      |<<<<<<<<<<<<<<<<<<<<<<<<<<<<<<<{ \phantom{AAA} }
      |>>>>>>>>>>>>>>>>>>{ \;\mathrm{Disc}(c_b)\; }
    \ar[dddr]|-{
      \mathclap{\phantom{\vert}}
      \mathbf{ch}_{B G_1}
    }
    &&&&
    \mathrm{Disc}
    \big(
      B G_2
    \big)
    \ar[dddr]
      |<<<<<<<<<<{
      \mathclap{\phantom{\vert}}
      \mathbf{ch}_{B G_2}
      }
      |>>>>>>>>>>>>{ \phantom{AA} }
    \\
    \\
    & &
    \Omega
    \big(
      -;
      \mathfrak{l}_{\scalebox{.6}{$B G_1$}}(A_1 \!\sslash\! G_1)
    \big)_{\mathrm{flat}}
    \ar@{-->}[rrrr]|<<<<<<<<<<{ \;(\mathfrak{c}_t)_\ast\; }
    \ar[dddd]|-{
      \mathclap{\phantom{\vert}}
      (\mathfrak{l}_{\scalebox{.6}{$B G_1$}}p_1)_\ast
    }
    &&&&
    \Omega
    \big(
      -;
      \mathfrak{l}_{\scalebox{.6}{$B G_1$}}(A_1 \!\sslash\! G_1)
    \big)_{\mathrm{flat}}
    \ar[dddd]|-{
      \mathclap{\phantom{\vert}}
      (\mathfrak{l}_{\scalebox{.6}{$B G_2$}}p_2)_\ast
    }
    \\
    &
    \flat
    \Bexp
    \big(
      \mathfrak{l} B G_1
    \big)
    \ar@{<-}[dddr]|-{
      \mathclap{\phantom{\vert}}
      \mathrm{atlas}
    }
    \ar[rrrr]
      |<<<<<<<<<<<<<<<<<<<<<<{ \phantom{AA} }
      |>>>>>>>>>>>>>>>>>{ \;(\mathfrak{c}_b)_\ast\; }
    &&&&
    \flat
    \Bexp
    \big(
      \mathfrak{l} B G_2
    \big)
    \ar@{<-}[dddr]|-{
      \mathclap{\phantom{\vert}}
      \mathrm{atlas}
    }
    \\
    \\
    {\phantom{AAAAAAAAAAA}}
    &
    &&&
    {\phantom{AAAAAAAAAAAA}}
    \\
    & &
    \Omega
    \big(
      -;
      \mathfrak{l}B G_1
    \big)_{\mathrm{flat}}
    \ar[rrrr]|>>>>>>>>>>>>>>>>>>>>>>>>{\; (\mathfrak{c}_b)_\ast }
    &&&&
    \Omega
    \big(
      -;
      \mathfrak{l}B G_2
    \big)_{\mathrm{flat}}
  }
  }
  }
\end{equation}

\end{defn}

At the level of dgc-algebras, the condition
that $\mathfrak{c}_t$ \eqref{MorphismForTwistedAbsoluteMinimalModel}
is a twisted absolute minimal
model for the transformation of local coefficient bundles
means
equivalently that it
makes the following cube commute:

\vspace{-8mm}
\begin{equation}
  \label{ConditionForTwistedAbsoluteMinimalModel}
  \hspace{-4mm}
  \raisebox{50pt}{
  \xymatrix@R=9pt@C=2.6em{
    \Omega^\bullet_{\mathrm{PLdR}}
    \big(
      A_1 \!\sslash\! G_1
    \big)
    \ar@{<-}[rrrr]|-{
      \;
      \Omega^\bullet_{\mathrm{PLdR}}
      (
        c_t
      )
      \;
    }
    \ar@{<-}[dddr]|-{
      \mathclap{\phantom{\vert^{\vert^{\vert}}}}
      \;
      p^{\mathrm{min}_{B G_1}}_{ A_1 \sslash G_1}
      \;
      \mathclap{\phantom{\vert_{\vert_{\vert}}}}
    }
    \ar@{<-}[dddd]|-{
      \mathclap{\phantom{\vert^{\vert^{\vert^{\vert}}}}}
      \Omega^\bullet_{\mathrm{PLdR}}
      (\rho_1)
      \mathclap{\phantom{\vert_{\vert_{\vert_{\vert}}}}}
    }
    &&&&
    \Omega^\bullet_{\mathrm{PLdR}}
    \big(
      A_2 \!\sslash\! G_2
    \big)
    \ar@{<-}[dddr]|-{
      \mathclap{\phantom{\vert^{\vert^{\vert}}}}
      \;
      p^{\mathrm{min}_{B G_1}}_{ A_1 \sslash G_1}
      \;
      \mathclap{\phantom{\vert_{\vert_{\vert}}}}
    }
    \ar@{<-}[dddd]
      |<<<<<<<<<{
      \mathclap{\phantom{\vert^{\vert^{\vert^{\vert}}}}}
      \Omega^\bullet_{\mathrm{PLdR}}
      (\rho_2)
      \mathclap{\phantom{\vert_{\vert_{\vert_{\vert}}}}}
      }
      |>>>>>>>{ \phantom{AA} \atop \phantom{AA} }
    \\
    \\
    \\
    &
    \mathrm{CE}
    \big(
      \mathfrak{l}_{\scalebox{.5}{$B G_1$}}
      (A_1 \!\sslash\! G_1)
    \big)
    \ar@{<--}[rrrr]|<<<<<<<<<<<<<<<<<<<<<{
      \;
      \mathrm{CE}(
        \mathfrak{c}_t
      )
      \;
    }
    \ar[dddd]|-{
      \mathclap{\phantom{\vert^{\vert^{\vert}}}}
      \mathrm{CE}(\mathfrak{l} \rho_1)
      \mathclap{\phantom{\vert_{\vert_{\vert}}}}
    }
    &&&&
    \mathrm{CE}
    \big(
      \mathfrak{l}_{\scalebox{.5}{$B G_2$}}
      (A_2 \!\sslash\! G_2)
    \big)
    \ar[dddd]|-{
      \mathclap{\phantom{\vert^{\vert^{\vert}}}}
      \mathrm{CE}(\mathfrak{l} \rho_2)
      \mathclap{\phantom{\vert_{\vert_{\vert}}}}
    }
    \\
    \Omega^\bullet_{\mathrm{PLdR}}
    (
      B G_1
    )
    \ar@{<-}[dddr]|-{
      p^{\mathrm{min}}_{B G_1}
    }
    \ar@{<-}[rrrr]
      |<<<<<<<<<<<<<<<<<<<<<<<<<<{ \phantom{AAA} }
      |>>>>>>>>>>>>>>>>>>>>{
      \;
      \Omega^\bullet_{\mathrm{PLdR}}
      (
        c_t
      )
      \;
    }
    &&&&
    \Omega^\bullet_{\mathrm{PLdR}}
    (
      B G_2
    )
    \ar@{<-}[dddr]|-{
      p^{\mathrm{min}}_{B G_2}
    }
    \\
    \\
    \\
    &
    \mathrm{CE}
    (
      \mathfrak{l} B G_1
    )
    \ar@{<-}[rrrr]|-{\;
      \mathfrak{c}_b
    \;}
    &&&&
    \mathrm{CE}
    (
      \mathfrak{l} B G_1
    )
  }
  }
\end{equation}

\medskip

In
differential enhancement of Def. \ref{TwistedNonabelianCohomologyOperation}
and in
twisted generalization of Def. \ref{SecondaryNonabelianCohomologyOperations},
we set:
\begin{defn}[Twisted secondary non-abelian cohomology operations]
 \label{TwistedSecondaryNonAbelianCohomologyOperations}
 Let
\vspace{-2mm}
$$
  \raisebox{20pt}{
  \xymatrix{
    A_1 \!\sslash\! G_1
    \ar[d]_{ \rho_1 }
    \ar[rr]^-{ c_t }
    &&
    A_2 \!\sslash\! G_2
    \ar[d]^{ \rho_2 }
    \\
    B G_1
    \ar[rr]^-{ c_b }
    &&
    B G_2
  }
  }
  \;\;\;\;\;
  \in
  \;\;
  \SimplicialSets
$$

\vspace{-2mm}
\noindent
be a transformation \eqref{MorphismOfLocalCoefficientBundles}
between local coefficient bundles \eqref{LocalCoefficientBundle},
together with an absolute minimal model
$\mathfrak{c}_b$ (Def. \ref{RelativeFormalMaps})
for the base map, and a compatible twisted absolute minimal
model $\mathfrak{c}_t$ (Def. \ref{TwistedAbsoluteMinimalModel}) for the total map. Then forming stage-wise homotopy pullbacks
(Def. \ref{HomotopyPullback})
in the required commuting cube \eqref{TransformationBetweenTwistedDifferentialNonAbelianCharacterMaps}
yields a transformation of corresponding differential coefficient
bundles \eqref{DifferentialLocalCoefficientBundle}:
\vspace{-2mm}
\begin{equation}
  \label{TransformationOfDifferentialLocalCoefficientBundles}
  \raisebox{20pt}{
  \xymatrix{
    (A_1 \!\sslash\! G_1)_{\mathrm{diff}}
    \ar[rr]^-{ (c_t)_{\mathrm{diff}} }
    \ar[d]_-{ (\rho_1)_{\mathrm{diff}} }
    &&
    (A_2 \!\sslash\! G_2)_{\mathrm{diff}}
    \ar[d]^-{ (\rho_2)_{\mathrm{diff}} }
    \\
    (B G_1)_{\mathrm{diff}}
    \ar[rr]_-{ (c_b)_{\mathrm{diff}} }
    &&
    (B G_2)_{\mathrm{diff}}
  }
  }
  \;\;\;\;
  \in \mathrm{PSh}
  \big(
    \CartesianSpaces
    \,,\,
    \SimplicialSets
  \big)
  \,.
\end{equation}

\vspace{-2mm}
\noindent
This yields, in turn, a
natural transformation of twisted differential non-abelian cohomology
sets (Def. \ref{TwistedDifferentialNonAbelianCohomology}),
hence a \emph{twisted secondary non-abelian cohomology operation},
by pasting composition,
hence by right derived base change (Ex. \ref{BaseChangeQuillenAdjunction}) along $(\rho_1)_{\mathrm{diff}}$
followed by composition with $(c_t)_{\mathrm{diff}}$ regarded
as a morphism in the slice (Example \ref{SliceModelCategory})
over $(B G_1)_{\mathrm{diff}}$:
\vspace{-3mm}
$$
  \xymatrix@C=5em{
    \widehat H^{\tau_{\mathrm{diff}}}
    \big(
      \mathcal{X}
      ;
      \,
      A_1
    \big)
    \ar[rr]^-{
      \left(
        (c_t)_{\mathrm{diff}}
        \,\circ\,
        (-)
      \right)
      \,\circ\,
      \left(
        (\rho_1)_{\mathrm{diff}}
      \right)_\ast
    }
    &&
    \widehat H^{
      (c_b)_{\mathrm{diff}}
      \,\circ\,
      \tau_{\mathrm{diff}}
    }
    \big(
      \mathcal{X}
      ;
      \,
      A_2
    \big).
  }
$$
\end{defn}

\medskip

In differential enhancement of Prop. \ref{TwistedChernCharacterInTwistedTopologicalKTheory}, we have:

\begin{example}[Twisted differential character on twisted differential K-theory]
  \label{TwistedDifferentialChernCharacterOnTwistedDifferentialKTheory}
  Consider the rationalization
  (Def. \ref{Rationalization})
  over the actual \emph{rational numbers} (see Remark \ref{RationalHomotopyTheoryOverTheRealNumbers})
  of the local
  coefficient bundle \eqref{LocalCoefficientBundleForTwistedKTheory}
  for degree-3 twisted complex topological
  K-theory (Example \ref{TwistedKTheory}).

  \noindent {\bf (i)} This is captured by the diagram
  \vspace{-5mm}
  \begin{equation}
    \label{RationalizationOfLocalCoefficientBundleForKTheory}
    \raisebox{20pt}{
    \xymatrix@C=4em{
      \mathrm{KU}_0 \!\sslash\! B \mathrm{U}(1)
      \ar[d]_-{ \rho }
      \ar[rr]^-{
        \eta^{\mathbb{Q}}_{
          \mathrm{KU}_0 \!\sslash\! B \mathrm{U}(1)
        }
      }
      &&
      L_{{}_{\mathbb{Q}}}
      \big(
        \mathrm{KU}_0 \!\sslash\! B \mathrm{U}(1)
      \big)
      \ar[d]^-{ L_{{}_{\mathbb{R}}} \rho }
      \\
      B^2 \mathrm{U}(1)
      \ar[rr]^-{
         \eta^{\mathbb{Q}}_{
           B^2 \mathrm{U}(1)
         }
      }
      &&
      L_{{}_{\mathbb{Q}}}
      \big(
        B^2 \mathrm{U}(1)
      \big)
    }
    }
  \end{equation}

  \vspace{-2mm}
\noindent
  regarded as a transformation of local coefficient
  bundles from twisted K-theory to twisted
  even-periodic rational cohomology:
  \vspace{-2mm}
  $$
    L_{{}_{\mathbb{Q}}}
    \mathrm{KU}_0
    \;\simeq\;
    \Omega^\infty
    \Big(\;
      \underset{
        =: \,
        H_{\mathrm{per}} \mathbb{Q}
      }{
        \underbrace{
          \underset{k}{\bigoplus}
          \,
          \Sigma^{2k} H \mathbb{Q}
        }
      }
    \; \Big)
    \,.
  $$

  \noindent {\bf (ii)}  Since rationalization is idempotent \eqref{RationalizationReflection},
  which here means that
  $
    L_{{}_{\mathbb{R}}} \,\circ\, L_{{}_{\mathbb{Q}}}
    \;\;
    \simeq
    \;\;
    L_{{}_{\mathbb{R}}}
    \,,
  $
  in this situation an
  absolute minimal model (Def. \ref{RelativeFormalMaps})
  of the base map
  $c_b \,=\, \eta^{\mathbb{R}}_{B^2 \mathrm{U}(1)}$
  and a twisted absolute minimal model
  (Def. \ref{TwistedAbsoluteMinimalModel}) of
  the total map
  $c_t \,=\, \eta^{\mathbb{R}}_{K_0 \sslash B \mathrm{U}(1)}$
  exist and are given, respectively, simply by
  the identity morphisms
 \vspace{-2mm}
  $$
  \mathfrak{c}_b \,:=\, \mathrm{id}_{ \mathfrak{l} B^2 \mathrm{U}(1) }
  \qquad
\mathrm{and}
\qquad
  \mathfrak{c}_t \,:=\,
    \mathrm{id}_{
      \mathfrak{l}_{B^2 \mathrm{U}(1)}
      (K_0 \!\sslash\! B \mathrm{U}(1) )
    }.
  $$

\vspace{-1mm}
   \noindent {\bf (iii)} Therefore, the induced twisted secondary cohomology
  operation Def. \ref{TwistedSecondaryNonAbelianCohomologyOperations}
  exists, and is for each differential twist $\tau_{\mathrm{diff}}$
  a transformation
   \vspace{-1mm}
  \begin{equation}
    \label{TwistedDifferentialChernCharacterAsTwistedSecondaryOperation}
    \xymatrix@C=5em{
      \widehat K^{\tau_{\mathrm{diff}}}
      \big(
        \mathcal{X}
      \big)
      \ar[rr]^-{
        \mathrm{ch}^{\tau_{\mathrm{diff}}}_{\mathrm{diff}}
        \;:=\;
        \left(
          \eta^{\mathbb{R}}_{K_0 \sslash B \mathrm{U}(1)}
        \right)_{\mathrm{diff}}
      }
      &&
      \widehat
      {
        H_{\mathrm{per}}\mathbb{Q}
      }^{\;L_{{}_{\mathbb{Q}}}\!\tau_{\mathrm{diff}}}
      \big(
        \mathcal{X}
      \big)
    }
  \end{equation}

  \vspace{-1mm}
\noindent  from twisted differential K-theory to twisted differential
  periodic rational cohomology theory.

  \noindent {\bf (iv)}  This is the
  twisted differential
  Chern character map on twisted differential complex K-theory
  as conceived in \cite[\S 3.2]{GS-tAHSS}\cite[Prop. 4]{GS-RR}.
  The analogous statement holds for the twisted differential
  Pontrjagin character (as in Example \ref{PontrjaginCharacterInKSO})
  on twisted differential real K-theory
  \cite[Thm. 12]{GS-tKO}.

  \noindent {\bf (v)}  Notice that this construction is close to but more
  structured
  than the plain curvature map on twisted differential K-theory
  (Example \ref{TwistedChernCharacterInTwistedDifferentialKTheory}):
  If we considered the transformation
  of local coefficients as in \eqref{RationalizationOfLocalCoefficientBundleForKTheory}
  but for rationalization $L_{{}_{\mathbb{R}}}$ over the
  real numbers (Remark \ref{RationalHomotopyTheoryOverTheRealNumbers}),
  then the induced
  twisted secondary cohomology operation would be equivalent
  to the twisted curvature map. Instead, \eqref{TwistedDifferentialChernCharacterAsTwistedSecondaryOperation}
  refines the plain curvature map to
  a twisted secondary operation that retains
  information about rational periods.
\end{example}

%\newpage

%%%%%%%%%%%%%%%%%%%%%%%%%%%%%%%%%%%%%%%%%%%%%%%%%%%%%%%%%
\subsection{Twisted character on twisted differential Cohomotopy}
 \label{CohomotopicalChernCharacter}
%%%%%%%%%%%%%%%%%%%%%%%%%%%%%%%%%%%%%%%%%%%%%%%%%%%%%%%%%%

We discuss here
(Example \ref{CharacterMapOnJTwistedCohomotopyAndTwistorialCohomotopy} below)
the twisted non-abelian character map on  J-twisted Cohomotopy
(Example \ref{JTwistedCohomotopyTheory})
in degree 4,
and
on Twistorial Cohomotopy (Example \ref{TwistorialCohomotopy}).
We highlight the induced charge quantization
(Prop. \ref{ChargeQuantizationInJTwistedCohomotopy} below)
and comment on the relevance to high energy physics
(Remark \ref{SummaryAndHypothesisH}).

\medskip
These twisted non-abelian cohomotopical character maps have been introduced
and analyzed in
\cite{FSS19b} and \cite{FSS20a}.
The general theory of non-abelian characters developed here
shows how these cohomotopical characters are
cousins both of
generalized abelian characters such as the Chern character on
twisted higher K-theory (\cref{TwistedChernCharacterOnTopologicalKTheory}),
notably of the character on topological modular forms
(by Example \ref{TheBoardmanHomomorphismIntmf},
and Remark \ref{ClarifyingTheRoleOfTmfInStringTheory})
as well as of non-abelian characters such as the
Chern-Weil homomorphism (\cref{TheChernWeilHomomorphism})
and the
Cheeger-Simons homomorphism (\cref{NonabelianDifferentialCohomology}).

\newpage

\noindent {\bf Cohomotopical character maps.}

\begin{example}[Character map on J-twisted Cohomotopy and on Twistorial Cohomotopy {\cite[Prop. 3.20]{FSS19b} \cite[Prop. 3.9]{FSS20a}}]
\label{CharacterMapOnJTwistedCohomotopyAndTwistorialCohomotopy}
Let $X$ be an 8-dimensional smooth spin manifold
equipped with tangential $\mathrm{Sp}(2)$-structure $\tau$ \eqref{TangentialSp2Structure}.
Then the twisted non-abelian character maps
(Def. \ref{TwistedNonAbelianChernDoldCharacter})
on J-twisted Cohomotopy
(Example \ref{JTwistedCohomotopyTheory})
in degree 4,
and
on Twistorial Cohomotopy (Example \ref{TwistorialCohomotopy})
are of the following form (with $p_1$, $p_2$, $I_8$ from Ex. \ref{PontrjaginForms}):
\vspace{-3mm}
$$
  \hspace{-4mm}
  \xymatrix@R=40pt@C=2pt{
    \mathbb{C}P^3 \!\sslash\! \mathrm{Sp}(2)
    \ar[d]|-{
      \overset{
        \raisebox{3pt}{
          \tiny
          \color{darkblue}
          \bf
          \def\arraystretch{.9}
          \begin{tabular}{c}
            Borel-equivariantized
            \\
            twistor fibration
          \end{tabular}
        }
      }{
        t_{\mathbb{H}} \sslash \mathrm{Sp}(2)
      }
    }
    &&
    \!\!\!\!\!\!
    \overset{
      \mathclap{
      \raisebox{3pt}{
        \tiny
        \color{darkblue}
        \bf
        {
        \def\arraystretch{.9}
        \begin{tabular}{c}
          Twistorial
          \\
          Cohomotopy
        \end{tabular}}
      }
      }
    }{
    \mathcal{T}^\tau(X)
    }
    \ar@{}[r]|-{:=}
    &
    H^\tau
    \big(
      X;
      \,
      \mathbb{C}P^3
    \big)
    \ar[rr]^-{
      \overset{
        \raisebox{3pt}{
          \tiny
          \color{darkblue}
          \bf
          {
          \def\arraystretch{.9}
          \begin{tabular}{c}
          character map on
          \\
          Twistorial Cohomotopy
          \end{tabular}}
        }
      }{
        \mathrm{ch}^\tau_{\mathbb{C}P^3}
      }
    }
    \ar[d]|-{
      \overset{
      \mbox{
        \tiny
        \color{darkblue}
        \bf
        {
        \def\arraystretch{.9}
        \begin{tabular}{c}
          cohomology operation
          \\
          along twistor fibration
        \end{tabular}}
      }
      }{
        (t_{\mathbb{H}})_*
           }
    }
    &{\phantom{AAAA}}&
    H_{\mathrm{dR}}^{\tau_{\mathrm{dR}}}
    \big(
      X;
      \,
      \mathfrak{l}\mathbb{C}P^3
    \big)
    \ar[d]|>>>>>>>>>>{
        (\mathfrak{l}t_{\mathbb{H}})_{\ast}
    }
    \ar@{}[r]|-{=}
    &
    \scalebox{.8}{$
    \left\{
    \!\!\!
    {\begin{array}{l}
      {\phantom{2}}H_3,
      \\
      {\phantom{2}}F_2,
      \\
      2G_7,
      \\
      {\phantom{2}}G_4
    \end{array}}
    \!\!\!
    %\in
    %\Omega^\bullet_{\mathrm{dR}}(X)
    \,\left\vert\;\,
    {\begin{aligned}
      d\, H_3 & = G_4 - \tfrac{1}{4}p_1(\nabla) - F_2 \wedge F_2,
      \\
      d\, F_2 & = 0,
      \\
      d\, 2G_7 & =
        -
        \big(
          G_4
          -
          \tfrac{1}{4}p_1(\nabla)
        \big)
        \wedge
        \big(
          G_4
          +
          \tfrac{1}{4}p_1(\nabla)
        \big)
        -
        \rchi_8(\nabla),
      \\
      d\, {\phantom{2}}G_4 & = 0
    \end{aligned}}
    \right.
    \!
    \right\}_{\!\!\!\big/\sim}
    $}
    \ar@<-75pt>[d]|-{
      \scalebox{.6}{$
      \begin{array}{cccc}
        H_3 & F_2 & 2G_7 & G_4
        \\
        \mapsdown &
        \mapsdown &
        \mapsdown &
        \mapsdown
        \\
        0 & 0 & 2G_7 & G_4
      \end{array}
      $}
    }
    \\
    S^4\!\sslash\! \mathrm{Sp}(2)
    &&
    \underset{
      \mathclap{
      \raisebox{-3pt}{
        \tiny
        \color{darkblue}
        \bf
        {
        \def\arraystretch{.9}
        \begin{tabular}{c}
          J-twisted
          \\
          4-Cohomotopy
        \end{tabular}}
      }
      }
    }{
    \pi^\tau
    \big(
      X
    \big)
    }
    \ar@{}[r]|-{:=}
    &
    H^\tau
    \big(
      X;
      \,
      S^4
    \big)
    \ar[rr]_-{
      \underset{
        \raisebox{-3pt}{
          \tiny
          \color{darkblue}
          \bf
          {
          \def\arraystretch{.9}
          \begin{tabular}{c}
          character map in
          \\
          J-twisted Cohomotopy
          \end{tabular}}
        }
      }{
        \mathrm{ch}^\tau_{S^4}
      }
    }
    &&
    H_{\mathrm{dR}}^{\tau_{\mathrm{dR}}}
    \big(
      X;
      \,
      \mathfrak{l}S^4
    \big)
    \ar@{}[r]|-{=}
    &
    \scalebox{.8}{$
    \left\{
    \!\!\!
    {\begin{array}{l}
      2G_7,
      \\
      {\phantom{2}}G_4
    \end{array}}
    \!\!\!
    %\in
    %\Omega^\bullet_{\mathrm{dR}}(X)
    \,\left\vert\;\,
    {\begin{aligned}
      d\, 2G_7 & =
        -
        \big(
          G_4
          -
          \tfrac{1}{4}p_1(\nabla)
        \big)
        \wedge
        \big(
          G_4
          +
          \tfrac{1}{4}p_1(\nabla)
        \big)
        -
        \rchi_8(\nabla),
      \\
      d\, {\phantom{2}}G_4 & = 0,
    \end{aligned}}
    \right.
    \!
    \right\}_{\!\!\!\big/\sim}
    $}
  }
$$
\end{example}

\vspace{-2mm}
\noindent
Here:
\begin{itemize}
\vspace{-2mm}
\item[{\bf (i)}] The twisted non-abelian de Rham cohomology targets
on the right are as shown, by Example \ref{FlatTwistedDifferentialFormsWithValuesInTwistorSpace}.
(In particular the twisted curvature forms in the first line are
relative to $\mathfrak{l}S^4$.)

\vspace{-3mm}
\item[{\bf (ii)}]
The vertical
twisted non-abelian cohomology operation (Def. \ref{TwistedNonabelianCohomologyOperation})
on the left
is induced from the Borel-equivariantized
twistor fibration \eqref{EquivariantizedTwistorFibration},
and that on the right from its
associated morphism of rational Whitehead $L_{\infty}$-algebras
(Prop. \ref{WhiteheadLInfinityAlgebrasRelative}).
\end{itemize}

\begin{prop}[Charge-quantization in J-twisted Cohomotopy {\cite[Prop. 3.13]{FSS19b}\cite[Cor. 3.11]{FSS20a}}]
\label{ChargeQuantizationInJTwistedCohomotopy}
Consider the twisted non-abelian character maps
(Def. \ref{TwistedNonAbelianChernDoldCharacter})
in J-twisted Cohomotopy and
in Twistorial Cohomotopy from Example \ref{CharacterMapOnJTwistedCohomotopyAndTwistorialCohomotopy}.

\noindent
{\bf (i)}
A necessary condition for a flat $\mathrm{Sp}(2)$-twisted
$\mathfrak{l}S^4$-valued differential form datum
$(G_4, G_7)$
to lift through the J-twisted cohomotopical character map
(i.e. to be in its image) is that the de Rham class of
$G_4$, when shifted by the fourth fraction of the
Pontrjagin form (Ex. \ref{PontrjaginForms}),
is in the image, under the de Rham homomorphism
(Example \ref{deRhamHomomorphism}), of an integral
class:
\vspace{-2.5mm}
\begin{equation}
  \label{IntegralityConditionInJTwistedCohomotopy}
  \big[
    G_4 - \tfrac{1}{4}p_1(\nabla)
  \big]
  \;\in\;
  \xymatrix{
    H^4(X;\, \mathbb{Z})
    \ar[r]
    &
    H^4_{\mathrm{dR}}(X)
    \,.
  }
\end{equation}

\vspace{-1.5mm}
\noindent
{\bf (ii)}
A necessary condition for a flat $\mathrm{Sp}(2)$-twisted
$\mathfrak{l}\mathbb{C}P^3$-valued differential form datum
$(G_4, G_7, F_2, H_3)$ to lift through the character map
in Twistorial Cohomotopy is that the de Rham class of
$G_4$ shifted by the fourth fraction of the Pontrjagin form (Ex. \ref{PontrjaginForms})
is in the image, under the de Rham homomorphism
(Example \ref{deRhamHomomorphism}), of an integral
class, and as such equal to the
$[F_2]$ cup-square:
\vspace{-1mm}
\begin{equation}
  \label{IntegralityConditionInTwistorialCohomotopy}
  \big[
    G_4 - \tfrac{1}{4}p_1(\nabla)
  \big]
  \;=\;
  \big[
    F_2 \wedge F_2
  \big]
  \;\in\;
  H^4_{\mathrm{dR}}(X)\;.
\end{equation}
\end{prop}

\medskip

\noindent {\bf Twisted differential Cohomotopy theory.}

\begin{defn}[Differential twists for twistorial Cohomotopy]
 \label{DifferentialTwistForTwistorialCohomotopy}
Let $X^8$ be an 8-dimensional smooth spin manifold
equipped with tangential $\mathrm{Sp}(2)$-structure $\tau$ \eqref{TangentialSp2Structure}.
By
\eqref{TangentialStrucuresAndTwistedCohomology} in
Example \ref{ClassificationOfTangentialStructure},
by \eqref{PushforwardOfTwistedCohomologyAlongCoefficientBundle} in
Example \ref{TotalNonAbelianClassOfTwistedCocycle},
and by \eqref{IsomorphismBetweenPrincipalBundlesAndMapsToBG}
in Example \ref{TraditionalNonAbelianCohomology}, we have
\vspace{-1mm}
\begin{equation}
  \label{Sp2Twist}
  [
    \tau
  ]
  \;\in\;
  \xymatrix{
  H^{\tau_{\mathrm{fr}}}
  \big(
    X;
    \,
    \mathrm{O}(n)/\mathrm{Sp}(2)
  \big)
  \ar[rr]^-{ (B i)_\ast }
  &&
  H\big(
    X^8;
    \,
    B \mathrm{Sp}(2)
  \big)
  }
  \;\;
    \simeq
  \;\;
  \mathrm{Sp}(2)\mathrm{Bundles}(X)_{\!/\sim}
  \,.
\end{equation}

\vspace{-1mm}
\noindent This gives, in particular,
the class of a smooth principal $\mathrm{Sp}(2)$-bundle
$P \to X$ to which the tangent bundle $T X$ is associated.
With \eqref{ForgetfulMapFromGBundlesWithConnection},
we may choose an $\mathrm{Sp}(2)$-connection $\nabla$
on $P$,
and, by Prop. \ref{DifferentialCohomologyOfPrincipalConnection},
this connection has a class
$[ \tau_{\mathrm{diff}}]$
in differential non-abelian cohomology
(Def. \ref{DifferentialNonAbelianCohomology}) with coefficients in
$B \mathrm{Sp}(2)$:
\vspace{-3mm}
$$
  \xymatrix@R=-2pt{
    H
    \big(
      X^8;
      \,
      B \mathrm{Sp}(2)
    \big)
    \ar@{}[r]|-{\simeq}
    &
    \mathrm{Sp}(2)\mathrm{Bundles}(X^8)_{\!/\sim}
    \ar@{<<-}[rr]
    &&
    \mathrm{Sp}(2)\mathrm{Connections}(X^8)_{\!/\sim}
    \ar@{-->}[rr]
    &&
    \widehat
    H
    \big(
      X^8;
      \,
      B \mathrm{Sp}(2)
    \big)
    \\
    [
      \tau
    ]
    \ar@{<->}[r]
    &
    [
      P
    ]
    \ar@{<-|}[rr]
    &&
    [
      \nabla
    ]
    \ar@{|->}[rr]
    &&
    [
      \tau_{\mathrm{diff}}
    ]
    \,.
  }
$$

\vspace{-2mm}
\noindent Any such $\tau_{\mathrm{diff}}$ serves as a
{\it differential twist} \eqref{DifferentialTwist} {\it for twistorial Cohomotopy} in the following.
\end{defn}

In twisted generalization of Example \ref{DifferentialCohomotopy},
we have:

\begin{example}[Differential twistorial Cohomotopy]
  \label{TwistorDifferentialJTwistedCohomotopy}
  Let $X^8$ be a spin 8-manifold equipped
  with tangential $\mathrm{Sp}(2)$-structure $\tau$ \eqref{TangentialSp2Structure},
  and with a corresponding differential twist
  $\tau_{\mathrm{diff}}$ (Def. \ref{DifferentialTwistForTwistorialCohomotopy}).

\vspace{1mm}
  \noindent
  {\bf (i)}
  Consider the local coefficient bundle \eqref{LocalCoefficientBundleForTwistedCohomotopy},
  $ S^4
      \to
            S^4 \!\sslash\! \mathrm{Sp}(2)
      \overset{ J_{\mathbb{C}P^3} }{\longrightarrow}
           B \mathrm{Sp}(2)
  $,
  for \emph{J-twisted 4-Cohomotopy} (Example \ref{JTwistedCohomotopyTheory})
  pulled back along
  $B \mathrm{Sp}(2)
    \overset{\simeq}{\longrightarrow}
    B \mathrm{Spin}(5)
\to
    B \mathrm{O}(5)
  $.
This induces, via Def. \ref{TwistedDifferentialNonAbelianCohomology},
  a twisted differential non-abelian cohomology theory
  $\widehat {\mathcal{T}}^{\; \tau_{\mathrm{diff}}}(-)$,
  which we call \emph{J-twisted differential 4-Cohomotopy},
  whose value on manifolds
  $\mathcal{X} \,=\, X^8 \times \mathbb{R}^k$
  sits in a cohomology operation diagram \eqref{SystemsOfCohomologyOperationsOnDifferentialCohomology}
  of this form:

  \vspace{-3mm}
  \begin{equation}
    \label{DifferentialJTwistedCohomotopyDiagram}
    \xymatrix@C=5em@R=16pt{
      \overset{
        \mathclap{
        \raisebox{3pt}{
          \tiny
          \color{darkblue}
          \bf
          {
          \def\arraystretch{.9}
          \begin{tabular}{c}
            differential
            \\
            J-twisted
            \\
            4-Cohomotopy
          \end{tabular}}
        }
        }
      }{
        \widehat {\pi}^{\, \tau^4_{\mathrm{diff}}}(\mathcal{X})
      }
      \ar[d]
      \ar[rr]_-{ F^{\tau^4_{\mathrm{dR}}}_{S^4}}^-{
        \overset{
          \mathclap{
          \raisebox{3pt}{
            \tiny
            \color{greenii}
            \bf
            {
            \def\arraystretch{.9}
            \begin{tabular}{c}
              J-twisted
              \\
              cohomotopical
              \\
              curvature
            \end{tabular}}
          }
          }
        }{
%          F^{\tau^4_{\mathrm{dR}}}_{S^4}
        }
      }
      &&
    \overset{
      \mathclap{
      \;\;\;\;\;\;\;\;\;
      \raisebox{3pt}{
        \tiny
        {
          \color{orangeii}
          \bf
          J-twisted cohomotopical Bianchi identities
        }
        (Example \ref{FlatTwistedDifferentialFormsWithValuesInTwistorSpace})
      }
      }
    }
    {
    \scalebox{.8}{$
    \left\{
    \!\!\!
    {\begin{array}{l}
      2G_7,
      \\
      {\phantom{2}}G_4
    \end{array}}
    \!\!\!
    \in
    \Omega^\bullet_{\mathrm{dR}}(\mathcal{X})
    \,\left\vert\;\,
    {\begin{aligned}
      d\, 2G_7 & =
        -
        \big(
          G_4
          -
          \tfrac{1}{4}p_1(\nabla)
        \big)
        \wedge
        \big(
          G_4
          +
          \tfrac{1}{4}p_1(\nabla)
        \big)
        -
        \rchi_8(\nabla),
      \\
      d\, {\phantom{2}}G_4 & = 0
    \end{aligned}}
    \right.
    \!
    \right\}
    $}
    }
    \ar[d]
    \\
    \underset{
      \mathclap{
      \raisebox{-3pt}{
        \tiny
        \def\arraystretch{.9}
        \begin{tabular}{c}
          \color{darkblue}
          \bf
          J-twisted
          \\
          \color{darkblue}
          \bf
          4-Cohomotopy
          \\
          (Example \ref{JTwistedCohomotopyTheory})
        \end{tabular}
      }
      }
    }{
      \pi^{\tau^4}(\mathcal{X})
    }
    \ar[rr]^-{\mathrm{ch}^{\tau^4}_{S^4}}_-{
      \underset{
        \mbox{
          \tiny
          \def\arraystretch{.9}
          \begin{tabular}{c}
            \color{greenii}
            \bf
            character map
            \\
            \color{greenii}
            \bf
            on J-twisted Cohomotopy
            \\
            (Example \ref{CharacterMapOnJTwistedCohomotopyAndTwistorialCohomotopy})
          \end{tabular}
        }
      }{
%        \mathrm{ch}^{\tau^4}_{S^4}
      }
    }
    &&
    \underset{
      \mathclap{
      \raisebox{-3pt}{
        \tiny
        \def\arraystretch{.9}
        \begin{tabular}{c}
          \color{darkblue}
          \bf
          J-twisted
          \\
          \color{darkblue}
          \bf
          de Rham cohomology
          \\
          (Def \ref{TwistedNonabelianDeRhamCohomology})
        \end{tabular}
      }
      }
    }{
    H^{\tau^4_{\mathrm{dR}}}_{\mathrm{dR}}
    \big(
      \mathcal{X};\,
      \mathfrak{l}S^4
    \big)
    \,.
    }
    }
  \end{equation}

\vspace{-2mm}
  \noindent
  {\bf (ii)}
  Consider the local coefficient bundle \eqref{EquivariantizedTwistorFibration}
  $
  \mathbb{C}P^3
    \to
    \mathbb{C}P^3 \!\sslash\! \mathrm{Sp}(2)
      \overset{ J_{\mathbb{C}P^3} }{\longrightarrow}
    B \mathrm{Sp}(2)
  $
  for \emph{twistorial Cohomotopy} (Def. \ref{TwistorialCohomotopy}).
This induces, via Def. \ref{TwistedDifferentialNonAbelianCohomology},
  a twisted differential non-abelian cohomology theory
  $\widehat {\mathcal{T}}^{\; \tau_{\mathrm{diff}}}(-)$,
 which we call \emph{differential twistorial Cohomotopy},
  whose value on manifolds $\mathcal{X} \,=\, X^8 \times \mathbb{R}^k$
  sits in a cohomology operation diagram \eqref{SystemsOfCohomologyOperationsOnDifferentialCohomology}
  of this form:

  \vspace{-.7cm}
  \begin{equation}
    \label{DifferentialTwistorialCohomotopyDiagram}
    \xymatrix@C=5em@R=16pt{
      \overset{
        \mathclap{
        \raisebox{3pt}{
          \tiny
          \color{darkblue}
          \bf
          {
          \def\arraystretch{.9}
          \begin{tabular}{c}
            differential
            \\
            twistorial
            \\
            Cohomotopy
          \end{tabular}}
        }
        }
      }{
        \widehat {\mathcal{T}}^{\;\tau_{\mathrm{diff}}}(\mathcal{X})
      }
      \ar[d]
      \ar[rr]_-{  F^{\tau_{\mathrm{dR}}}_{\mathbb{C}P^3}}^-{
        \overset{
          \mathclap{
          \raisebox{3pt}{
            \tiny
            \color{greenii}
            \bf
            {
            \def\arraystretch{.9}
            \begin{tabular}{c}
              twistorial
              \\
              curvature
            \end{tabular}}
          }
          }
        }{
 %         F^{\tau_{\mathrm{dR}}}_{\mathbb{C}P^3}
        }
      }
      &&
    \overset{
      \mathclap{
      \!\!\!\!\!\!\!\!
      \raisebox{3pt}{
        \tiny
        {
          \color{orangeii}
          \bf
          twistorial Bianchi identities
        }
        (Example \ref{FlatTwistedDifferentialFormsWithValuesInTwistorSpace})
      }
      }
    }
    {
    \scalebox{.8}{$
    \left\{
    \!\!\!
    {\begin{array}{l}
      {\phantom{2}}H_3,
      \\
      {\phantom{2}}F_2,
      \\
      2G_7,
      \\
      {\phantom{2}}G_4
    \end{array}}
    \!\!\!
    \in
    \Omega^\bullet_{\mathrm{dR}}(\mathcal{X})
    \,\left\vert\;\,
    {\begin{aligned}
      d\, H_3 & = G_4 - \tfrac{1}{4}p_1(\nabla) - F_2 \wedge F_2,
      \\
      d\, F_2 & = 0,
      \\
      d\, 2G_7 & =
        -
        \big(
          G_4
          -
          \tfrac{1}{4}p_1(\nabla)
        \big)
        \wedge
        \big(
          G_4
          +
          \tfrac{1}{4}p_1(\nabla)
        \big)
        -
        \rchi_8(\nabla),
      \\
      d\, {\phantom{2}}G_4 & = 0
    \end{aligned}}
    \right.
    \!
    \right\}
    $}
    }
    \ar[d]
    \\
    \underset{
      \mathclap{
      \raisebox{-3pt}{
        \tiny
        \def\arraystretch{.9}
        \begin{tabular}{c}
          \color{darkblue}
          \bf
          twistorial
          \\
          \color{darkblue}
          \bf
          Cohomotopy
          \\
          (Example \ref{TwistorialCohomotopy})
        \end{tabular}
      }
      }
    }{
      \mathcal{T}^{\tau}(\mathcal{X})
    }
    \ar[rr]^-{\mathrm{ch}^\tau_{\mathbb{C}P^3}}_-{
      \underset{
        \mbox{
          \tiny
          \def\arraystretch{.9}
          \begin{tabular}{c}
            \color{greenii}
            \bf
            character map
            \\
            \color{greenii}
            \bf
            on twistorial Cohomotopy
            \\
            (Example \ref{CharacterMapOnJTwistedCohomotopyAndTwistorialCohomotopy})
          \end{tabular}
        }
      }{
%        \mathrm{ch}^\tau_{\mathbb{C}P^3}
      }
    }
    &&
    \underset{
      \mathclap{
      \raisebox{-3pt}{
        \tiny
        \def\arraystretch{.9}
        \begin{tabular}{c}
          \color{darkblue}
          \bf
          twistorial
          \\
          \color{darkblue}
          \bf
          de Rham cohomology
          \\
          (Def \ref{TwistedNonabelianDeRhamCohomology})
        \end{tabular}
      }
      }
    }{
    H^{\tau_{\mathrm{dR}}}_{\mathrm{dR}}
    \big(
      \mathcal{X};\,
      \mathfrak{l}\mathbb{C}P^3
    \big)
    \,.
    }
    }
  \end{equation}

\end{example}

\begin{prop}[Twisted secondary cohomology operation induced by twistor fibration]
  \label{TwistedSecondaryOperationFromTwistorialToJTwistedCohomotopy}
  The defining twisted non-abelian cohomology operation
  \eqref{TwistedCohomologyTransformationsDefiningTwistorialCohomotopy}
  from twistorial Cohomotopy (Example \ref{TwistorialCohomotopy})
  to J-twisted 4-Cohomotopy (Example \ref{JTwistedCohomotopyTheory}),
  induced by the $\mathrm{Sp}(2)$-equivariantized
  twistor fibration $t_{\mathbb{H}} \sslash \mathrm{Sp}(2)$
  \eqref{EquivariantizedTwistorFibration}
  refines to a twisted secondary cohomology operation
  (Def. \ref{TwistedSecondaryNonAbelianCohomologyOperations})
  from differential twistorial Cohomotopy
  to differential J-twisted Cohomotopy (Example \ref{TwistorDifferentialJTwistedCohomotopy}):

  \vspace{-.3cm}
  $$
    \xymatrix@R=15pt@C=9em{
      \mathllap{
        \mbox{
        \tiny
        \color{darkblue}
        \bf
        \def\arraystretch{.9}
        \begin{tabular}{c}
          differential
          \\
          twistorial
          \\
          Cohomotopy
        \end{tabular}
      }
      }
      \widehat {\mathcal{T}}^{\; \tau_{\mathrm{diff}} }
      \big(
        \mathcal{X}
      \big)
      \ar[dd]|-{
        \mathllap{
          \mbox{
            \tiny
            \color{greenii}
            \bf
            \def\arraystretch{.9}
            \begin{tabular}{c}
              twisted secondary
              \\
              cohomology operation
            \end{tabular}
          }
          \!\!\!
        }
        \left(
        \left(
          t_{\mathbb{H}}\sslash \mathrm{Sp}(2)
        \right)_{\mathrm{diff}}
        \right)_\ast
        \mathrlap{
          \!\!\!
          \mbox{
            \tiny
            \color{greenii}
            \bf
            \def\arraystretch{.9}
            \begin{tabular}{c}
              along
              \\
              $\mathrm{Sp}(2)$-equivariantized
              \\
              twistor fibration
            \end{tabular}
          }
        }
      }
      \ar[rr]^-{ c^\tau_{\mathbb{C}P^3} }
      &&
      {\mathcal{T}}^{ \tau }
      \big(
        \mathcal{X}
      \big)
      \ar[dd]|-{
        \left(
          t_{\mathbb{H}}\sslash \mathrm{Sp}(2)
        \right)_\ast
        \mathrlap{
          \!\!\!
          \mbox{
          \tiny
          \color{greenii}
          \bf
          \def\arraystretch{.9}
          \begin{tabular}{c}
            twisted primary
            \\
            cohomology operation
          \end{tabular}
          }
        }
      }
      \\
      \\
      \mathllap{
        \mbox{
          \tiny
          \color{darkblue}
          \bf
          \def\arraystretch{.9}
          \begin{tabular}{c}
            differential
            \\
            J-twisted
            \\
            4-Cohomotopy
          \end{tabular}
        }
        \!\!\!
      }
      \widehat {\pi}^{ \,\tau^4_{\mathrm{diff}} }
      \big(
        \mathcal{X}
      \big)
      \ar[rr]^-{ c^{\tau^4}_{S^4} }
      &&
      {\pi}^{ \tau^4 }
      \big(
        \mathcal{X}
      \big)
    }
  $$
\end{prop}
\begin{proof}
  By Def. \ref{TwistedSecondaryNonAbelianCohomologyOperations}
  we need to show that we have
  a twisted absolute minimal model
  (Def. \ref{TwistedAbsoluteMinimalModel}) for the
  $\mathrm{Sp}(2)$-equivariantized twistor fibration
  \eqref{EquivariantizedTwistorFibration}.
  By \eqref{ConditionForTwistedAbsoluteMinimalModel}
  this means that we can find a morphism

  \vspace{-2mm}
  \begin{equation}
    \label{TwistedAbsoluteModelForEquivariantizedTwistorFibration}
    \xymatrix{
      \mathfrak{l}_{\scalebox{.5}{$B \mathrm{Sp}(2)$}}
      \big(
        \mathbb{C}P^3 \!\sslash\! \mathrm{Sp}(2)
      \big)
      \ar@{-->}[rr]^-{
        \mathfrak{t}_{\mathbb{H}} \sslash \mathfrak{l}\mathrm{Sp}(2)
      }
      &&
      \mathfrak{l}_{\scalebox{.5}{$B \mathrm{Sp}(2)$}}
      \big(
        S^4 \!\sslash\! \mathrm{Sp}(2)
      \big)
    }
  \end{equation}

  \vspace{0mm}
  \noindent
  between the relative Whitehead $L_\infty$-algebras
  (Prop. \ref{WhiteheadLInfinityAlgebrasRelative}) of the two local
  coefficient bundles, which makes the following cube of transformations
  of derived PL-de Rham adjunction units commute:
  \vspace{-1mm}
$$
  \begin{tikzcd}[row sep=35pt, column sep=-64pt]
    \Bexp_{\mathrm{PL}}
      \circ
    \Omega^\bullet_{\mathrm{PLdR}}
    \big(
      \mathbb{C}P^3
        \!\sslash\!
      \mathrm{Sp}(2)
    \big)
    \ar[
      rr,
      "{
        \color{greenii}
        \Bexp_{\mathrm{PL}}
          \circ
        \Omega^\bullet_{\mathrm{PLdR}}
        \big(
          t_{\mathbb{H}} \!\sslash\! \mathrm{Sp}(2)
        \big)
      }"
    ]
    \ar[
      dd,
      "{
        \color{greenii}
        \Bexp_{\mathrm{PL}}
          \circ
        \Omega^\bullet_{\mathrm{PLdR}}
        \big(
          J_{{}_{\mathbb{C}P^3}}
        \big)
      }"{description, pos=.8}
    ]
    \ar[
      dr,
      "{
        \color{greenii}
        \Bexp_{\mathrm{PL}}
        \left(
          p_{\mathbb{C}P^3 \!\sslash\! \mathrm{Sp}(2)}^{\mathrm{min}_{B \mathrm{Sp}(2)}}
        \right)
      }"{description}
    ]
    &
    &[+80pt]
    \Bexp_{\mathrm{PL}}
      \circ
    \Omega^\bullet_{\mathrm{PLdR}}
    \big(
      S^4
        \!\sslash\!
      \mathrm{Sp}(2)
    \big)
    \ar[
      dd
    ]
    \ar[
      dr,
      "{
        \color{greenii}
        \Bexp_{\mathrm{PL}}
        \left(
          p_{S^4 \!\sslash\! \mathrm{Sp}(2)}^{\mathrm{min}_{B \mathrm{Sp}(2)}}
        \right)
      }"{description}
    ]
    \\
    &
    \Bexp_{\mathrm{PL}}
      \circ
    \mathrm{CE}
    \Big(
      \mathfrak{l}_{{}_{B \mathrm{Sp}(2)}}
      \big(
      \mathbb{C}P^3
        \!\sslash\!
      \mathrm{Sp}(2)
      \big)
    \Big)
    \ar[
      rr,
      dashed,
      crossing over,
      "{
        \color{greenii}
        \Bexp_{\mathrm{PL}}
          \circ
        \mathrm{CE}
        (
          t_{\mathbb{H}}
            \!\sslash\!
          \mathfrak{l}\mathrm{Sp}(2)
        )
      }"{above, pos=.36}
    ]
    &
    &[+20pt]
    \Bexp_{\mathrm{PL}}
      \circ
    \mathrm{CE}
    \Big(
      \mathfrak{l}_{{}_{B \mathrm{Sp}(2)}}
      \big(
      S^4
        \!\sslash\!
      \mathrm{Sp}(2)
      \big)
    \Big)
    \\[-10pt]
    \Bexp_{\mathrm{PL}}
      \circ
    \Omega^\bullet_{\mathrm{PLdR}}
    \big(
      \mathrm{Sp}(2)
    \big)
    \ar[
      rr,-,
      shift left=1pt
    ]
    \ar[
      rr,-,
      shift right=1pt
    ]
    \ar[
      dr,
      "{
        \color{greenii}
        \Bexp_{\mathrm{PL}}
        \left(
          p^{\mathrm{min}}_{B\mathrm{Sp}(2)}
        \right)
      }"{description}
    ]
    &&
    \Bexp_{\mathrm{PL}}
      \circ
    \Omega^\bullet_{\mathrm{PLdR}}
    \big(
      \mathrm{Sp}(2)
    \big)
    \ar[
      dr,
      "{
        \color{greenii}
        \Bexp_{\mathrm{PL}}
        \left(
          p^{\mathrm{min}}_{B\mathrm{Sp}(2)}
        \right)
      }"{description}
    ]
    \\
    &
    \Bexp_{\mathrm{PL}}
      \circ
    \mathrm{CE}
    \big(
      \mathfrak{l}\mathrm{Sp}(2)
    \big)
    \ar[
      rr,-,
      shift left=1pt
    ]
    \ar[
      rr,-,
      shift right=1pt
    ]
    \ar[
      from=uu,
      crossing over,
    ]
    &&
    \Bexp_{\mathrm{PL}}
      \circ
    \mathrm{CE}
    \big(
      \mathfrak{l}\mathrm{Sp}(2)
    \big)
    \ar[
      from=uu
    ]
  \end{tikzcd}
$$

\noindent But, from Example \ref{FlatTwistedDifferentialFormsWithValuesInTwistorSpace},
we see that the total object of the relative Whitehead $L_\infty$-algebra
of $\mathbb{C}P^3 \sslash \mathrm{Sp}(2)$, relative to
$\mathfrak{l} B \mathrm{Sp}(2)$, coincides with that relative to
$\mathfrak{l}_{\scalebox{.6}{$B \mathrm{Sp}(2)$}}\big(S^4 \!\sslash\! \mathrm{Sp}(2)\big)$.
Therefore, we may take the twisted absolute minimal model
\eqref{TwistedAbsoluteModelForEquivariantizedTwistorFibration} to be
equal to top arrow in Example \ref{FlatTwistedDifferentialFormsWithValuesInTwistorSpace}.
This makes the front square commute by construction,
and it being a relative minimal model for $t_{\mathbb{H}} \!\sslash\! \mathrm{Sp}(2)$
implies by Prop. \ref{ExistenceOfRelativeMinimalSullianModels}
that there is an essentially unique top left morphism such that the
top square commutes.
\end{proof}

\begin{remark}[Lifting against the twisted differential
twistor fibration]
  \label{LiftingAgainstTheTwistedDifferentialTwistorFibration}
$\,$

\vspace{.7mm}
\hspace{-.9cm}
\begin{tabular}{ll}
  \begin{minipage}[left]{10.5cm}
  In terms of differential moduli $\infty$-stacks \eqref{FibrationOfDifferentialNonabelianCharacterMaps},
  the result of
  Prop. \ref{TwistedSecondaryOperationFromTwistorialToJTwistedCohomotopy}
  with Example \ref{TwistorDifferentialJTwistedCohomotopy} says
  that lifting a twisted differential Cohomotopy cocycle $\widehat C_3$
  with 4-flux density $G_4$
  against the twisted differential
  refinement \eqref{TransformationOfDifferentialLocalCoefficientBundles}
  of the equivariantized twistor fibration \eqref{EquivariantizedTwistorFibration}
  to a differential twistorial Cohomotopy cocycle
  $(\widehat C_3, \widehat C_2, \widehat C_1)$
  involves, on twisted  curvature forms \eqref{CurvatureOnTwistedDifferentialCohomology}
  the appearance of a 2-flux density $F_2$
  and of a 3-form $H_3$ such that
  $      d H_3 \,=\, G_4  - \tfrac{1}{4}p_1(\nabla) - F_2 \wedge F_2$.
  \end{minipage}
  &
  $
    \raisebox{32pt}{
    \xymatrix@R=45pt@C=2em{
      \mathbb{R}^{1,1} \times X^8
      \ar@{-->}[rr]
        _-{
          (
            \widehat C_3,\, \widehat B_2,\, \widehat A_1
          )
        }
        ^-{
          \mathclap{
          \raisebox{3pt}{
            \tiny
            \color{greenii}
            \bf
            \def\arraystretch{.9}
            \begin{tabular}{c}
              lift of C-field through
              \\
              twistor fibration
            \end{tabular}
          }
          }
        }
      \ar@{^{(}->}[d]
      &&
      \big(
        \mathbb{C}P^3 \sslash \mathrm{Sp}(2)
      \big)_{\mathrm{diff}}
      \ar[d]_-{
        (t_{\mathbb{H}}\sslash \mathrm{Sp}(2))_{\mathrm{diff}}
        \mathrlap{\!\!\!
                \mbox{
            \tiny
            \color{greenii}
            \bf
            \def\arraystretch{.9}
            \begin{tabular}{c}
              twisted differential
              \\
              twistor fibration
            \end{tabular}
          }
        }
      }
      \\
      \mathbb{R}^{2,1} \times X^{8}
      \ar[rr]^-{
        \widehat C_3
      }_-{
        \mbox{
          \tiny
          \color{greenii}
          \bf
          C-field
        }
      }
      &&
      \big(
        S^4 \sslash \mathrm{Sp}(2)
      \big)_{\mathrm{diff}}
    }
    }
  $
\end{tabular}
\end{remark}

\begin{remark}[M-theory fields and Hypothesis H]
 \label{SummaryAndHypothesisH}
In summary, we have found:

\noindent
{\bf (i)} A cocycle $\widehat C_3$ in J-twisted differential
4-Cohomotopy (Example \ref{TwistorDifferentialJTwistedCohomotopy})
has as curvature/character forms \eqref{CurvatureOnTwistedDifferentialCohomology}:

{\bf (a)} a closed 4-form $G_4$, hence a 4-flux density,

{\bf (b)} a non-closed 7-form $G_7$,

\noindent
satisfying the following Bianchi identities
(Example \ref{CharacterMapOnJTwistedCohomotopyAndTwistorialCohomotopy})
and integrality conditions (Prop. \ref{ChargeQuantizationInJTwistedCohomotopy}):
\vspace{-2mm}
\begin{equation}
\label{CohomotopyFields}
\hspace{0cm}
  \xymatrix@C=40pt@R=-5pt{
    \overset{
      \mathclap{
      \raisebox{3pt}{
        \tiny
        \color{darkblue}
        \bf
        {
        \def\arraystretch{.9}
        \begin{tabular}{c}
          differential
          \\
          J-twisted
          4-Cohomotopy
        \end{tabular}}
      }
      }
    }{
     \widehat {\pi}^{\; \tau^4_{\mathrm{diff}}}
       (
        \mathcal{X}
      )
    }
    \ar[rr]^{
      \mbox{
        \tiny
        {
        \def\arraystretch{.9}
        \begin{tabular}{c}
          \color{greenii}
          \bf
          curvature
          \\
          (non-abelian character form representative)
        \end{tabular}}
      }
    }
    &&
    \;\;
    \overset{
      \mathclap{
      \raisebox{3pt}{
        \tiny
        \color{darkblue}
        \bf
        \def\arraystretch{.9}
        \begin{tabular}{c}
          flat twisted cohomotopical
          \\
          differential forms
        \end{tabular}
      }
      }
    }{
      \Omega^{\tau_{\mathrm{dR}}}_{\mathrm{dR}}
      \big(
        X;
        \,
        \mathfrak{l}S^4
      \big)_{\mathrm{flat}}
    }
    \\
    \big(
      \widehat C_3
    \big)
    \ar@{}[rr]|-{ \longmapsto }
    &&
 {\footnotesize   \left(
    \!\!
    {\begin{array}{l}
      G_4,
      \\
      2G_7
    \end{array}}
    \,\left\vert\;\,
    {\begin{aligned}
      &
      \mathrlap{
      \overset{
        \raisebox{3pt}{
          \tiny
          \color{orangeii}
          \bf
          shifted C-field flux quantization
        }
      }{
      \big[
        G_4 - \tfrac{1}{4}p_1(\omega)
      \big]
      }
      \in
      H^4(X;\, \mathbb{Z})
      }
      \\
      d\, {\phantom{2}}G_4 & = 0
      \\
      d\, 2G_7 &
      \underset{
        \raisebox{3pt}{
          \tiny
          \color{orangeii}
          \bf
          {
          \def\arraystretch{.9}
          \begin{tabular}{c}
            C-field tadpole cancellation \& M5 Hopf WZ term level quantization
          \end{tabular}}
        }
      }{
        =
        -
        \big(
          G_4
          -
          \tfrac{1}{4}p_1(\omega)
        \big)
        \wedge
        \big(
          G_4
          +
          \tfrac{1}{4}p_1(\omega)
        \big)
        -
        24\, I_8(\omega)
        }
    \end{aligned}}
    \right.
      \!\!  \right)
      \,,
 }
  }
\end{equation}
\vspace{-.3cm}

\noindent where the characteristic forms $p_1$, $p_2 $ and $I_8$ are from
Ex. \ref{PontrjaginForms}.

\noindent
{\bf (ii)} Lifting this cocycle through the
twisted differential twistor fibration
(Prop. \ref{TwistedSecondaryOperationFromTwistorialToJTwistedCohomotopy})
to a cocycle  $\big( \widehat C_3, \widehat B_2, \widehat A_1\big)$
in differential twistorial Cohomotopy
(Example \ref{TwistorDifferentialJTwistedCohomotopy})
involves (Remark \ref{LiftingAgainstTheTwistedDifferentialTwistorFibration})
adjoining to the 4-flux density $G_4$:

{\bf (c)} a closed 2-form curvature $F_2$, hence a 2-flux density,

{\bf (d)} a non-closed 3-form $H_3$,

\noindent
such that these curvature/character forms satisfy the following Bianchi identities (Example \ref{CharacterMapOnJTwistedCohomotopyAndTwistorialCohomotopy})
and integrality conditions (Prop. \ref{ChargeQuantizationInJTwistedCohomotopy}):

\vspace{-.7cm}
\begin{equation}
\label{TwistorialCohomotopyFields}
\hspace{0cm}
  \xymatrix@C=45pt@R=-5pt{
    \overset{
      \mathclap{
      \raisebox{3pt}{
        \tiny
        \color{darkblue}
        \bf
        {
        \def\arraystretch{.9}
        \begin{tabular}{c}
          differential
          \\
          twistorial
          Cohomotopy
        \end{tabular}}
      }
      }
    }{
     \widehat {\mathcal{T}}^{\; \tau_{\mathrm{diff}}}
       (
        \mathcal{X}
      )
    }
    \ar[rr]^{
      \mbox{
        \tiny
        {
        \def\arraystretch{.9}
        \begin{tabular}{c}
          \color{greenii}
          \bf
          curvature
          \\
          (non-abelian character form representative)
        \end{tabular}}
      }
    }
    &&
    \;\;
    \overset{
      \mathclap{
      \raisebox{3pt}{
        \tiny
        \color{darkblue}
        \bf
        \def\arraystretch{.9}
        \begin{tabular}{c}
          flat twistorial
          \\
          differential forms
        \end{tabular}
      }
      }
    }{
      \Omega^{\tau_{\mathrm{dR}}}_{\mathrm{dR}}
      \big(
        X;
        \,
        \mathfrak{l}\mathbb{C}P^3
      \big)_{\mathrm{flat}}
    }
    \\
    \big(
      \widehat C_3
      \,,\,
      \widehat B_2
      \,,\,
      \widehat A_1
    \big)
    \ar@{}[rr]|-{ \longmapsto }
    &&
 {\footnotesize   \left(
    \!\!
    {\begin{array}{l}
      {\phantom{2}}H_3,
      \\
      \phantom{\overset{\raisebox{3pt}{\tiny a}}{H_3}}
      \\
      G_4, F_2,
      \\
      2G_7
    \end{array}}
    \,\left\vert\;\,
    {\begin{aligned}
      d\, H_3 & = G_4 - \tfrac{1}{4}p_1(\omega) - F_2 \wedge F_2,
      \\
      &
      \mathrlap{
      \overset{
        \raisebox{3pt}{
          \tiny
          \color{orangeii}
          \bf
          Ho{\v r}ava-Witten Green-Schwarz mechanism
        }
      }{
      \big[
        G_4 - \tfrac{1}{4}p_1(\omega)
      \big]
      =
      \big[
        F_2 \wedge F_2
      \big]
      }
      \in
      H^4(X;\, \mathbb{Z})
      }
      \\
      d\, {\phantom{2}}G_4 & = 0\,,\;\;\;d\, F_2 \; = 0,
      \\
      d\, 2G_7 &
      \underset{
        \raisebox{3pt}{
          \tiny
          \color{orangeii}
          \bf
          {
          \def\arraystretch{.9}
          \begin{tabular}{c}
            C-field tadpole cancellation \& M5 Hopf WZ term level quantization
          \end{tabular}}
        }
      }{
        =
        -
        \big(
          G_4
          -
          \tfrac{1}{4}p_1(\omega)
        \big)
        \wedge
        \big(
          G_4
          +
          \tfrac{1}{4}p_1(\omega)
        \big)
        -
        24\, I_8(\omega)
        }
    \end{aligned}}
    \right.
      \!\!  \right)
 }
  }
\end{equation}

\noindent
{\bf (iii)}
With these cohomotopical curvature/character forms
interpreted as flux densities, this is the
Bianchi identities and charge quantization expected
in M-theory on the supergravity C-field
($\widehat C_3$), the heterotic B-field $(\widehat B_2)$
and the heterotic $S\big( \mathrm{U}(1)^2\big) \subset E_8$
gauge field $(\widehat A_1)$, with the following prominent features:

\vspace{1mm}
\hspace{-.5cm}
{\bf (a)} {\bf The charge quantization:}

\vspace{-.7cm}
\begin{equation}
\mbox{
\hspace{-.57cm}
\def\arraystretch{2.3}
\begin{tabular}{lll}
 {\bf (1)}
 &
 $\big[  G_4 - \tfrac{1}{4}p_1(\nabla) \big]
 \phantom{= [F_2 \wedge F_2]\;}
 \in H^4(\mathcal{X};\, \mathbb{Z})$
 &
 \begin{minipage}[left]{7.6cm}
 is expected for the C-field in the M-theory bulk
 \\
 (\cite[(1.2)]{Witten97a}\cite[(1.2)]{Witten97b})
 \end{minipage}
 \\
 {\bf (2)}
 &
 $\big[  G_4 - \tfrac{1}{4}p_1(\nabla) \big] = [F_2 \wedge F_2] \in H^4(\mathcal{X};\, \mathbb{Z})$
 &
 \begin{minipage}[left]{7.6cm}
   is expected on heterotic boundaries
   \\
   (\cite[(1.13)]{HoravaWitten96}, review in \cite[\S 1]{FSS20a})
 \end{minipage}
\end{tabular}
}
\end{equation}

\hspace{-.5cm}
{\bf (b)} {\bf The quadratic functions:}

\vspace{-.7cm}
\begin{equation}
\label{QuadraticFunctions}
\hspace{-.57cm}
\mbox{
\begin{tabular}{lll}
\;\; {\bf (1)}
&
  \raisebox{-8pt}{
  $
  \begin{aligned}
  G_4 \;\mapsto\;
  &
  ({\color{orangeii}G_4} - \tfrac{1}{4}p_1(\omega))
    \wedge
  ({\color{orangeii}G_4} + \tfrac{1}{4}p_1(\omega))
  \\
  & + 24 I_8(\omega)
  \end{aligned}
  $
  }
&
 \begin{minipage}[left]{7.6cm}
   constitute the Hopf Wess-Zumino term
   \\
   (\cite[p. 11]{Aharony96}\cite{Intrilligator00}, see \cite{FSS19c}\cite{SS20a})
 \end{minipage}
\\
\;\; {\bf (2)}
&
$F_2 \;\mapsto\;
  {\color{orangeii}F_2}
    \wedge
  {\color{orangeii}F_2}
  \;\;$
&
\begin{minipage}[left]{7.6cm}
  constitute the  2nd Chern class of a
  \\
  $\mathrm{U}(1) \!\!\subset\!\! \mathrm{SU}(2) \!\!\subset\!\! E_8$-bundle
  %\cite{AndersonGrayLukasPalti12}
  \cite{AGLP11}\cite{AGLP12}\cite[(7)]{FSS20a}
\end{minipage}
\end{tabular}
}
\end{equation}

\noindent
{\bf (iv)} These are necessary, not yet sufficient constraints
on cohomotopical lifts. Further constraints follow
by Postnikov tower analysis \cite{GS-Postnikov}
and coincide with further expected conditions
in M-theory (see \cite[Table 1]{FSS19b}).

\medskip

\noindent
All this suggests the {\it Hypothesis H}
\cite{FSS19b}\cite{FSS19c}\cite{SS19a}\cite{SS19b}\cite{BMSS19}\cite{SS20a}\cite{FSS20a}\cite{FSS20TwistedString}\cite{SS21},
following \cite[\S 2.5]{Sati13},
that the elusive cohomology theory which controls M-theory
in analogy to how K-theory controls string theory is:
(twisted, equivariant, differential) non-abelian Cohomotopy theory.
\end{remark}

\medskip

\noindent
{\bf Cohomotopical character into K-theory.}
We may regard the
{\it secondary non-abelian Hurewicz/Boardman homomorphism}
(Example \ref{SecondaryNonAbelianBoardmanHomomorphismToKTheory})
from differential 4-Cohomotopy
(Example \ref{DifferentialCohomotopy})
to differential K-theory (Example \ref{CurvatureInDifferentialComplexKTheory}),
as a non-abelian but K-theory valued character,
lifting the target of the cohomotopical character
(Example \ref{CharacterMapOnJTwistedCohomotopyAndTwistorialCohomotopy})
from
rational cohomology to K-theory (compare \cite[Fig. 1]{BSS19}):

\vspace{-.4cm}
\begin{equation}
  \label{ChargeQuantizationOnKTheory}
  \xymatrix@R=11pt@C=50pt{
  \overset{
    \mathclap{
    \raisebox{3pt}{
      \tiny
      \color{darkblue}
      \bf
      \def\arraystretch{.9}
      \begin{tabular}{c}
        differential
        \\
        Cohomotopy
      \end{tabular}
    }
    }
  }{
    {\widehat{\tau}}^{\, 4}
    \big(
      \mathcal{X}
    \big)
    }
    \ar[rr]^-{ (e_{\mathrm{KU}})^4_{\mathrm{diff}} }_-{
      \mathclap{
      \mbox{
        \tiny
        \def\arraystretch{.9}
        \begin{tabular}{c}
          \color{greenii}
          \bf
          differential non-abelian
          \\
          \color{greenii}
          \bf
          Boardman homomorphism
          \\
          (Example \ref{SecondaryNonAbelianBoardmanHomomorphismToKTheory})
        \end{tabular}
      }
      }
    }
    &&
    \;
    \overset{
      \mathclap{
      \raisebox{3pt}{
        \tiny
        \color{darkblue}
        \bf
        \def\arraystretch{.9}
        \begin{tabular}{c}
          differential
          \\
          K-theory
        \end{tabular}
      }
      }
    }{
      \widehat{\mathrm{KU}}^0
      \big(
        \mathcal{X}
      \big)
    }
    \;
    \ar[rr]^-{
        \mathrm{ch}_{\mathrm{diff}}
    }_-{
        \raisebox{3pt}{
          \tiny
          {
          \def\arraystretch{.9}
          \begin{tabular}{c}
            \color{greenii}
            \bf
            differential
            \\
            \color{greenii}
            \bf
            Chern character
            \\
            (Example \eqref{TwistedDifferentialChernCharacterOnTwistedDifferentialKTheory})
          \end{tabular}}
        }
    }
    &&
    \overset{
      \mathclap{
      \raisebox{3pt}{
        \tiny
        \color{darkblue}
        \bf
        \def\arraystretch{.9}
        \begin{tabular}{c}
          differential
          \\
          rational cohomology
        \end{tabular}
      }
      }
    }{
      \widehat {H_{\mathrm{per}} \mathbb{Q}}^{\, 0}
      \big(
        \mathcal{X}
      \big)
    }
    \\
    \ar@{<~}[rr]|-{
      \mbox{
        \scalebox{.5}{ \color{orangeii} \bf
        \def\arraystretch{.9}
        \begin{tabular}{c}
          charge-quantization
          \\
          in M-theory
        \end{tabular}
        }
      }
    }
    &&
    \ar@{<~}[rr]|-{
      \mbox{
        \scalebox{.5}{ \color{orangeii}  \bf
        \def\arraystretch{.9}
        \begin{tabular}{c}
          charge-quantization
          \\
          in string theory
        \end{tabular}
        }
      }
    }
    &&
  }
\end{equation}

\noindent {\bf (i)} Lifting through this differential
Boardman homomorphism induces
{\it secondary charge quantization conditions}
\emph{on} K-theory,
analogous to \eqref{DiracChargeQuantization} but
invisible even in generalized cohomology,
instead now coming from
non-abelian cohomology theory.

\vspace{1mm}
\noindent {\bf (ii)}  In the plain version \eqref{ChargeQuantizationOnKTheory}
(i.e. disregarding twisting and equivariant enhancement)
the effect of $(e_{\mathrm{KU}})^4_{\mathrm{diff}}$
on curvature forms \eqref{CurvatureOnDifferentialCohomology}
is (by Example \ref{SecondaryNonAbelianBoardmanHomomorphismToKTheory})
to forget the quadratic function \eqref{QuadraticFunctions}
on $G_4$ and to inject what remains as the
4-form curvature component $F_4$ in differential K-theory:

  \vspace{-.4cm}
  \begin{equation}
        \label{DifferentialBoardmanOnCurvatures}
    \raisebox{35pt}{
    \xymatrix@C=6em{
      \overset{
        \mathclap{
        \raisebox{3pt}{
          \tiny
          \color{darkblue}
          \bf
          \def\arraystretch{.9}
          \begin{tabular}{c}
            differential
            \\
            Cohomotopy
          \end{tabular}
        }
        }
      }{
        \widehat \tau^{\, 4}
        \big(
          \mathcal{X}
        \big)
      }
      \ar[d]_-{
        \big(
          F_{S^4}
        \big)_\ast
      }
      \ar[rr]_-{   (e_{\mathrm{KU}})^4_{\mathrm{diff}}}^-{
        \overset{
          \mathclap{
          \raisebox{3pt}{
            \tiny
            \color{greenii}
            \bf
            \def\arraystretch{.9}
            \begin{tabular}{c}
              secondary
              non-abelian
              \\
              Boardman homomorphism
            \end{tabular}
          }
          }
        }{
 %         \beta^4_{\mathrm{diff}}
        }
      }
      &
      \ar@{}[d]|-{
        \mathclap{
        \mbox{
          \tiny
          \color{greenii}
          \bf
          curvatures/flux densities
        }
        }
      }
      &
      \overset{
        \mathclap{
        \raisebox{3pt}{
          \tiny
          \color{darkblue}
          \bf
          \def\arraystretch{.9}
          \begin{tabular}{c}
            differential
            \\
            K-theory
          \end{tabular}
        }
        }
      }{
        \widehat{\mathrm{KU}}^0
        \big(
          \mathcal{X}
        \big)
      }
      \ar[d]^-{
        \big(
          F_{\mathrm{KU}_0}
        \big)_\ast
      }
      \\
      \left\{
        \!\!\!
        {\begin{array}{c}
          2G_7,
          \\[-3pt]
          \phantom{2}G_4
        \end{array}}
        \,\left\vert\;
        {\begin{aligned}
          d\, 2 G_7 & = - G_4 \wedge G_4
          \\[-3pt]
          d\, \phantom{2}G_4 & = 0
        \end{aligned}}
        \right.
      \right\}
      \ar[rr]_-{
        \scalebox{.7}{$
          \arraycolsep=1.4pt%\def\arraystretch{.1}
          {\begin{array}{ccc}
            G_4 &\mapsto& F_4
            \\[-3pt]
            G_7 &\mapsto& 0
          \end{array}}
        $}
      }
      &
      &
      \Big\{
        \!
        \big(
          F_{2k}
        \big)
        \,\left\vert\;
          d\, F_{2k} \,=\, 0
        \right.
      \Big\}
      \,.
    }
    }
  \end{equation}
  \vspace{-.4cm}

\noindent
This is a `cohomotopical enhancement' of the
reduction in \cite{DMW03} of $E_8$ bundles in M-theory to
the K-theory of type IIA string theory,
now characterized by higher Postnikov stages
of the Boardman homomorphism  \cite{GS-Postnikov}.
The remaining RR-flux components in $\{F_{2k}\}$
are also found
in the cohomotopical character, through
cohomological double dimensional reduction formulated
in parametrized homotopy theory:
this is discussed in detail in \cite{BMSS19}.

\vspace{1mm}
\noindent {\bf (iii)}
The twisted generalization of the non-abelian Boardman homomorphism in
\eqref{ChargeQuantizationOnKTheory}
and \eqref{DifferentialBoardmanOnCurvatures}
is more subtle, since the degree-3 twist of K-theory
does not arise from the J-twist of Cohomotopy, but
arises, together with the further RR-flux components,
from $S^1$-equivariantization/double dimensional reduction
of Cohomotopy \cite{FSS16a}\cite{BMSS19},
reproducing the reduction of $E_8$ bundles from M-theory to type IIA
in \cite{MathaiSati}.

%will have twisted entries. The lift of ${\rm ch}_{\rm diff}$ will
%then be replaced by that of the twisted differential Chern character,
%as discussed in detail in \cite{GS-RR}. This is also a cohomotopical %enhancement of the
%reduction of $E_8$ bundles from M-theory to type IIA
%in \cite{MathaiSati}.
%Likewise, the crucially missing
%fluxes are provided in \cite{BMSS19}.

%\vspace{1mm}
%\noindent {\bf (iv)}
%The comparison between what
%$E_8$ bundles in spacetime give, which is essentially only a degree four %integral class, and
%what differential 4-Cohomotopy provides, which is that and additional %subtle torsion constraints,
% is performed via systematic study of  the differential Postnikov  tower %for $S^4$ in \cite{GS-Postnikov}.
%These constraints are consistent with the twisted Cohomotopical
%treatment of prominent M-theory anomalies in \cite{FSS19b}.

\medskip
\medskip
\medskip
\noindent {\bf Outlook: Equivariant enhancement.}
The twisted non-abelian character theory presented
here enhances further to
{\it proper} (i.e. Bredon-style not Borel-style)
\emph{equivariant} non-abelian cohomology on orbi-orientifolds,
by combining it with the techniques developed in
\cite{HSS18}\cite{SS20b} (essentially: parametrizing
the construction here over the orbit category
of the equivariance group).
The resulting
\emph{character map in equivariant non-abelian cohomology}
is discussed in \cite[\S 2,3]{SS20c}.
For example, the equivariant enhancement of the
cohomotopical character into K-theory
\eqref{ChargeQuantizationOnKTheory},
lifting the RR-fields in equivariant K-theory
through the equivariantized enhancement
of the Boardman homomorphism on the left of
\eqref{ChargeQuantizationOnKTheory},
enforces \cite{SS19a}\cite{BSS19} ``tadpole cancellation'' conditions
expected in string theory at orbifold singularities.

\newpage

\appendix

%%%%%%%%%%%%%%%%%%%%%%%%%%%%%%%%%%%%
\section{Model category theory}
\label{ModelCategoryTheory}
%%%%%%%%%%%%%%%%%%%%%%%%%%%%%%%%%%%%

For ease of reference and
to highlight some less widely used aspects needed in the main text,
we record basics of homotopy theory
via model category theory \cite{Quillen67}
(review in \cite{Hovey99}\cite{Hirschhorn02}\cite[A.2]{Lurie09})
and of homotopy topos theory \cite{Rezk10}
via model categories of simplicial presheaves
\cite{Brown73}\cite{Jardine87}\cite{Dugger01}
(review in \cite{Dugger98}\cite[\S A.3.3]{Lurie09}\cite{Jardine15}).

\medskip

\vspace{1mm}
\noindent {\bf Topology.} By
\vspace{-1mm}
\begin{equation}
  \label{ConvenientCategoryOfTopologicalSpaces}
  \TopologicalSpaces
  \;\in\;
  \Categories
\end{equation}

\vspace{-1mm}
\noindent we denote a \emph{convenient} \cite{Steenrod67}
(in particular: cartesian closed)
category of topological spaces such as
compactly-generated spaces \cite{Strickland09}
or $\Delta$-generated spaces \cite{Dugger03},
equivalently known as: numerically-generated spaces \cite{SYH10}
or D-topological spaces \cite[Prop. 2.4]{SS20b}.

\medskip

\noindent {\bf Categories.}
Let $\mathcal{C}$ be a category.

\noindent
{\bf (i)}
For $X,A \in \mathcal{C}$
a pair of objects, we write
\vspace{-2mm}
$$
  \mathcal{C}(X,A) \;:=\; \mathrm{Hom}_{\mathcal{C}}(X,A)
$$

\vspace{-2mm}
\noindent
for the set of {\it morphisms} from $X$ to $A$.

\noindent {\bf (ii)}
For $\mathcal{C}, \mathcal{D}$ two categories, we denote
a pair of {\it adjoint functors} between them by
\vspace{-3mm}
  \begin{equation}
    \label{AnAdjunction}
    \xymatrix{
      \mathcal{D}
      \ar@{<-}@<+7pt>[rr]^-{L}
      \ar@<-7pt>[rr]_-{R}^-{\bot}
      &&
      \mathcal{C}
    }
    \;\;\;\;
    \Leftrightarrow
    \;\;\;\;
    \begin{tikzcd}
      \mathcal{D}(L(-) \,,\, -)
      \ar[
        r, <->,
        "\widetilde{(-)}"{above},
        "\sim"{below, yshift=-1pt}
      ]
      &
      \mathcal{C}(- \,,\, R(-))\,,
    \end{tikzcd}
  \end{equation}
  and the corresponding {\it adjunction unit} and {\it adjunction counit}
  transformations by, respectively:

  \vspace{-.3cm}
  \begin{equation}
    \label{AdjunctionUnit}
    \eta^{RL}_C
    \;\colon\;
    C
    \xrightarrow{
      \widetilde{\mathrm{id}_{L X}}
    }
    R \circ L(C)
    \,,
    {\phantom{AAAAAAAA}}
    \epsilon^{L R}_D
    \;\colon\;
    L \circ R(D)
    \xrightarrow{
      \widetilde{\mathrm{id}_{R D}}
    }
    D
  \end{equation}
  Notice/recall that this means that adjunct morphisms $f \leftrightarrow \widetilde{f}$
  \eqref{AnAdjunction}
  and (co-)units \eqref{AdjunctionUnit} are related as follows:

  \vspace{-.4cm}
  \begin{equation}
  \label{AdjunctsByCompositesWithCoUnits}
  \begin{tikzcd}[column sep={between origins, 33}]
    L(c)
    \ar[rr, "f"]
    &&
    d
    &{\phantom{AAA}}\leftrightarrow{\phantom{AAA}}&
    c
    \ar[rr, "\eta_c"]
    \ar[
      rrrr,
      rounded corners,
      to path={
           -- ([yshift=+12pt]\tikztostart.north)
           --node[above]{
               \scalebox{.7}{$
                 \widetilde{f}
               $}
             } ([yshift=+8.5pt]\tikztotarget.north)
           -- (\tikztotarget.north)}
    ]
    &&
    R \circ L(c)
    \ar[rr, "R(f)"]
    &&
    R(d)
    \mathrlap{\,,}
  \end{tikzcd}
\end{equation}

\vspace{-.5cm}
$$
  \begin{tikzcd}[column sep={between origins, 33}]
    L(c)
    \ar[rr, "L(\widetilde{f})"]
    \ar[
      rrrr,
      rounded corners,
      to path={
           -- ([yshift=+8pt]\tikztostart.north)
           --node[above]{
               \scalebox{.7}{$
                 f
               $}
             } ([yshift=+8.5pt]\tikztotarget.north)
           -- (\tikztotarget.north)}
    ]
    &&
    L \circ R(d)
    \ar[
      rr,
      "\epsilon_d"
    ]
    &&
    d
    &{\phantom{AAA}}\leftrightarrow{\phantom{AAA}}&
    c
    \ar[rr, "\widetilde f"]
    &&
    R(d)
    \mathrlap{\,.}
  \end{tikzcd}
$$
\vspace{-3mm}
\noindent
{\bf (iii)}
A Cartesian square in $\mathcal{C}$ we indicate by
{\it pullback} notation $f^\ast(-)$ and/or by the symbol
``{\scalebox{.8}{(pb)}}'':
  \vspace{-2mm}
\begin{equation}
  \label{CartesianSquare}
  \raisebox{20pt}{
  \xymatrix@C=4em@R=1.2em{
    f^\ast A
    \ar[d]_-{
      f^\ast p
    }
    \ar[d]
    \ar@{}[dr]|-{\mbox{\tiny(pb)}}
    \ar[r]
    &
    A
    \ar[d]^-{ p }
    \\
    B_1
    \ar[r]_-{f}
    &
    B_2
    \,.
  }
  }
\end{equation}

Dually, a co-Cartesian square in $\mathcal{C}$ we indicate by
{\it pushout} notation $f_\ast(-)$ and/or by the symbol
``{\scalebox{.8}{(po)}}'':
  \vspace{-2mm}
\begin{equation}
  \label{coCartesianSquare}
  \raisebox{20pt}{
  \xymatrix@R=1.2em@C=4em{
    A_1
    \ar[r]^-{ f }
    \ar[d]_-{
      q
    }
    \ar@{}[dr]|-{\mbox{\tiny(po)}}
    &
    A_2
    \ar[d]^-{ f_\ast q }
    \\
    B
    \ar[r]
    &
    f_\ast B
    \,.
  }
  }
\end{equation}

\medskip

\noindent {\bf Model categories.}

\vspace{-1mm}
\begin{defn}[Weak equivalences]
  \label{WeakEquivalences}
  A \emph{category with weak equivalences} is a category $\mathcal{C}$
  equipped with a sub-class $\mathrm{W} \subset \mathrm{Mor}(\mathcal{C})$ of its morphisms, to be called the class of \emph{weak equivalences},
  such that

  \begin{itemize}
  \vspace{-.3cm}
  \item[{\bf (i)}] $\mathrm{W}$ contains the class of isomorphisms;

  \vspace{-.2cm}
  \item[{\bf (ii)}] $\mathrm{W}$ satisfies the cancellation property
   (``2-out-of-3''): if in any commuting triangle in $\mathcal{C}$
   \vspace{-2mm}
   \begin{equation}
     \label{TwoOutOfThree}
     \xymatrix@R=5pt{
       & Y
       \ar[dr]^-{ g }
       \\
       X
       \ar[ur]^-{ f }
       \ar[rr]_-{ g \circ f }
       &&
       Z
     }
   \end{equation}

   \vspace{-4mm}
   \noindent
   two morphisms are in $\mathrm{W}$, then so is the third.
   \end{itemize}
\end{defn}

\begin{defn}[Weak factorization system]
  \label{WeakFactorizationSystem}
  A \emph{weak factorization system} in a category $\mathcal{C}$
  is a pair of sub-classes of morphisms
  $\mathrm{Proj, \mathrm{Inj}} \;\subset\; \mathrm{Mor}(\mathcal{C})$
  such that

  \begin{itemize}
    \vspace{-.4cm}
    \item[{\bf (i)}] every morphisms
    $\xymatrix@C=12pt{X \ar[r]^-{f} & Y}$ in $\mathcal{C}$ may be factored through
    a morphism in $\mathrm{Proj}$ followed by one in $\mathrm{Inj}$:
    \vspace{-2mm}
    \begin{equation}
      \label{FactorizationOfMorphismByWeakFactorizationSystem}
      f
      \; :
      \xymatrix@C=3em{
        X
        \ar[r]^-{ \in \; \mathrm{Proj} }
        &
        Z
        \ar[r]^-{ \in \; \mathrm{Inj} }
        &
        Y
      }
    \end{equation}
    \vspace{-.9cm}
    \item[{\bf (ii)}] For every commuting square in $\mathcal{C}$
    with left morphism in $\mathrm{Proj}$ and right morphism in
    $\mathrm{Inj}$, there exists a lift, namely a dashed morphism
    \vspace{-2mm}
    \begin{equation}
      \label{LiftingProperty}
      \xymatrix@R=1.5em@C=4em{
        X
        \ar[r]
        \ar[d]_-{ \in\, \mathrm{Proj} }
        &
        A
        \ar[d]^-{ \in\, \mathrm{Inj} }
        \\
        Y
        \ar@{-->}[ur]|-{ \;\;\exists\;\; }
        \ar[r]
        &
        B
      }
    \end{equation}

    \vspace{-2mm}
    \noindent
    making the resulting triangles commute.

    \vspace{-.2cm}
    \item[{\bf (iii)}] Given $\mathrm{Inj}$ (resp. $\mathrm{Proj}$),
    the class $\mathrm{Proj}$ (resp. $\mathrm{Inj}$)
    is the largest class for which \eqref{LiftingProperty} holds.
  \end{itemize}
\end{defn}

\begin{defn}[Model category {\cite[Def. E.1.2]{Joyal}\cite{Riehl09}}]
  \label{ModelCategories}
  A \emph{model category} is a category $\mathbf{C}$
  that has all small limits and colimits,
  equipped with three sub-classes of its class of morphisms,
  to be denoted

  $\mathrm{W}$ -- {\it weak equivalences}

  $\mathrm{Fib}$ -- {\it fibrations}

  $\mathrm{Cof}$ -- {\it cofibrations}

  \noindent such that

  \noindent  {\bf (i)} The class $\mathrm{W}$ makes $\mathbf{C}$ a category with weak equivalences (Def. \ref{WeakEquivalences});

  \noindent  {\bf (ii)} The pairs $\big( \mathrm{Fib}\,,\, \mathrm{Cof} \cap \mathrm{W} \big)$ and $\big( \mathrm{Fib} \cap \mathrm{W}\,,\, \mathrm{Cof} \big)$
  are weak factorization systems (Def. \ref{WeakFactorizationSystem}).
\end{defn}
\begin{remark}[Minimal assumptions]
  By item {(iii)} in Def. \ref{WeakFactorizationSystem}
  a model category structure is specified already by
  the classes $\mathrm{W}$ and $\mathrm{Fib}$,
  or alternatively by the classes $\mathrm{W}$ and $\mathrm{Cof}$.
    Moreover, it follows from Def. \ref{ModelCategories}
  that also the class $\mathrm{W}$ is stable under retracts
  {\cite[Prop. E.1.3]{Joyal}\cite[Lemma 2.4]{Riehl09}}:
  Given a commuting diagram in the model category $\mathbf{C}$
  of the form on the left here

  \vspace{1mm}
  \begin{equation}
    \label{ClassOfWeakEquivalencesClosedUnderRetracts}
    \raisebox{20pt}{
    \xymatrix@C=3em@R=2em{
      X \ar[r]
      \ar[d]_f
      \ar@/^1pc/@{=}[rr]
      &
      Y
      \ar[d]|-{
        \mathclap{\phantom{\vert^\vert}}
        \in \; \mathrm{W}
        \mathclap{\phantom{\vert_\vert}}
      }
      \ar[r]
      & X
      \ar[d]^f
      \\
      A
       \ar[r]
        \ar@/_1pc/@{=}[rr]
      & B \ar[r] & A
    }
    }
    \;\;\;\;\;\;\;\;
    \Rightarrow
    \;\;\;\;\;\;\;\;
    f \;\in \; \mathrm{W}
  \end{equation}

    \vspace{2mm}
  \noindent
  with the middle morphism a weak equivalence, then also
  $f$ is a weak equivalence.
\end{remark}

\begin{defn}[Proper model category]
  \label{ProperModelCategories}
  A model category $\mathbf{C}$ Def. \ref{ModelCategories}
  is called

  \noindent
  {\bf (i)} \emph{right proper}, if
  pullback along fibrations preserves weak equivalences:
    \vspace{-2mm}
  \begin{equation}
    \label{RightProperness}
    \raisebox{20pt}{
    \xymatrix@R=1.5em@C=4em{
      X
      \ar[d]
      \ar[r]^-{p^\ast f}
      \ar@{}[dr]|-{\mbox{\tiny (pb)}}
      &
      A
      \ar[d]_-{p}^-{ \in \; \mathrm{Fib} }
      \\
      Y
      \ar[r]_-{f \in \; \mathrm{W} }
      &
      B
    }
    }
    \;\;\;\;\;\;\;\;\;\;
    \Rightarrow
    \;\;\;\;\;\;\;\;\;\;
    p^\ast f \;\in\; \mathrm{W}
  \end{equation}

  \vspace{-2mm}
  \noindent
  {\bf (ii)} \emph{left proper}, if
  pushout along cofibrations preserves weak equivalences,
  hence if the opposite model category (Example \ref{OppositeModelStructure})
  is right proper.
\end{defn}

\begin{notation}[Fibrant and cofibrant objects]
  \label{FibrantAndCofibrantObjects}
  Let $\mathbf{C}$ be a model category (Def. \ref{ModelCategories})

  \noindent {\bf (i)} We write $\ast \in \mathbf{C}$ for the terminal object
  and $\varnothing \in \mathbf{C}$ for the initial object.

  \noindent {\bf (ii)} An object $X \in \mathbf{C}$ is called:

  \vspace{-.1cm}

  {\bf (a)} {\it fibrant} if the unique morphism to the terminal
    object is a fibration,
    $\xymatrix@C=25pt{ X \ar[r]^{\in\; \mathrm{Fib}} & \ast }$;

  {\bf (b)} {\it cofibrant} if the unique morphism from the initial object is a cofibration,
  $\xymatrix@C=25pt{ \varnothing \ar[r]^{\in\; \mathrm{Cof}} & X }$.

  \noindent We write
  $\mathbf{C}_{\mathrm{fib}}, \mathbf{C}^{\mathrm{cof}},
  \mathbf{C}^{ \mathrm{cof} }_{\mathrm{fib}}
  \;\subset\; \mathbf{C}$
  for the full subcategories on, respectively,
  fibrant objects, or cofibrant objects or
  objects that are both fibrant and cofibrant.

  \noindent {\bf (iii)} Given an object $X \in \mathbf{C}$

  \begin{itemize}
    \vspace{-.3cm}
    \item[\bf (a)] A {\it fibrant replacement} is a factorization
    \eqref{FactorizationOfMorphismByWeakFactorizationSystem}
        of the terminal morphism as
          \vspace{-2mm}
        \begin{equation}
          \label{FibrantReplacement}
          \xymatrix{
            X
            \ar[rr]^-{ j_X }_-{ \in \; \mathrm{Cof} \, \cap \mathrm{W} }
            &&
            P X
            \ar[rr]^-{ q_X }_-{ \in \; \mathrm{Fib} }
            &&
            \ast \;.
          }
        \end{equation}

           \vspace{-4mm}
    \item[\bf (b)] A {\it cofibrant replacement} is a factorization \eqref{FactorizationOfMorphismByWeakFactorizationSystem}
        of the initial morphism as
           \vspace{-2mm}
        \begin{equation}
          \label{CofibrantReplacement}
          \xymatrix{
            \varnothing
            \ar[rr]^-{ i_X }_-{ \in \; \mathrm{Cof}  }
            &&
            Q X
            \ar[rr]^-{ p_X }_-{ \in \; \mathrm{Fib}\, \cap \mathrm{W} }
            &&
            X\;.
          }
        \end{equation}
  \end{itemize}
\end{notation}

Recall that a continuous function $f$ between topological spaces
$X, Y$
induces homomorphism on all homotopy groups
$\pi_k(f,x) : \pi_k(X,x) \to \pi_k\big(Y,f(x)\big)$
and is called (e.g. \cite[p. 144]{tomDieck08})
\vspace{-2mm}
\begin{equation}
  \label{WeakHomotopyEquivalence}
  \mbox{$f$ is}
  \left\{
  \begin{array}{lcl}
    \mbox{an {\it $n$-equivalence}}
    &\mbox{if}&
    \def\arraystretch{.9}
    \begin{array}{l}
    \mbox{
      $\pi_{\bullet < n}(f,x)$ is an isomorphism and
    }
    \\
    \mbox{
      $\pi_n(f,x)$ is an epimorphism
    }
    \end{array}
    \\
    \mbox{a {\it weak homotopy equivalence}}
    &\mbox{if}&
    \mbox{
      $\pi_\bullet(f,x)$ is an isomorphism
    }
  \end{array}
  \right.
  \end{equation}
for all $x \in X$.

\begin{example}[Classical model structure on topological spaces {\cite[\S II.3]{Quillen67}\cite{Hi15}}]
  \label{ClassicalModelStructureOnTopologicalSpaces}
  The category $\TopologicalSpaces$ \eqref{ConvenientCategoryOfTopologicalSpaces}
  carries a model category structure whose

  {\bf (i)} $\mathrm{W}$ -- weak equivalences are the weak homotopy equivalences \eqref{WeakHomotopyEquivalence};

  {\bf (ii)} $\mathrm{Fib}$ -- fibrations are the Serre fibrations (e.g. \cite[\S 5.5, 6.3]{tomDieck08}).

  \noindent We denote this model category by
  \vspace{-4mm}
  $$
    \TopologicalSpaces_{\mathrm{Qu}}
    \;\in\;
    \mathrm{ModelCategories}\;.
  $$
\end{example}

\begin{example}[Classical model structure on simplicial sets
{\cite[\S II.3]{Quillen67}\cite[\S V.1-2]{GelfandManin96}\cite[\S I.11]{GoerssJardine99}}]
  \label{ClassicalModelStructureOnSimplicialSets}
  The category of $\SimplicialSets$
  of simplicial sets (e.g. \cite{May67}\cite{Curtis71}
  exposition in \cite{Friedman08})
  carries a model category structure whose

  {\bf (i)} $\mathrm{W}$ -- weak equivalences are those
  whose geometric realization is a weak homotopy equivalence;

  {\bf (ii)} $\mathrm{Cof}$ -- cofibrations are the monomorphisms
  (degreewise injections).

  {\bf (iii)} $\mathrm{Fib}$ -- fibrations are the Kan fibrations.

  \noindent We denote this model category by
  \vspace{-4mm}
  $$
    \SimplicialSets_{\mathrm{Qu}}
    \;\in\;
    \mathrm{ModelCategories}\;.
  $$
  \vspace{-.4cm}

  Every simplicial set is cofibrant in the classical model structure
  (Nota. \ref{FibrantAndCofibrantObjects}):

  \vspace{-.4cm}
  \begin{equation}
    \label{EverySimplicialSetIsCofibrantInClassicalModelStructure}
    \big(
      \SimplicialSets_{\mathrm{Qu}}
    \big)_{\mathrm{cof}}
    \;\;
    \;=\;
    \SimplicialSets
    \,.
  \end{equation}
  \vspace{-.5cm}

  \noindent
  while the fibrant simplicial sets are exactly the
  {\it Kan complexes}
  (e.g. \cite[\S I.3]{GoerssJardine99}, exposition in \cite[\S 7]{Friedman08})

  \vspace{-2mm}
  \begin{equation}
    \label{KanComplexes}
    \big(\SimplicialSets_{\mathrm{Qu}}\big)^{\mathrm{fib}}
    \;=\;
    \mathrm{KanComplexes}
    \,.
  \end{equation}
  \vspace{-.4cm}

  \noindent
  which we may think of as {\it $\infty$-groupoids} \cite[\S 1.1.2]{Lurie09}.
\end{example}

\begin{example}[Simplicial nerves of groupoids]
  \label{SimplicialNervesOfGroupoids}

  \noindent
  {\bf (i)}
  Let
  \vspace{-2mm}
  $$
    \mathcal{G}
    \;=\;
    \Big(\!\!
    \begin{tikzcd}
      \mathcal{G}_1 \;{}_{t}\!\!\times_{s} \mathcal{G}_1
      \ar[
        r,
        "{\circ}"
      ]
      &
      \mathcal{G}_1
      \ar[
        r,
        shift left=2pt,
        "{s}"
      ]
      \ar[
        r,
        shift right=2pt,
        "{t}"{below}
      ]
      &
      \mathcal{G}_0
    \end{tikzcd}
    \!\Big)
  $$

  \vspace{-2mm}
\noindent  be a groupoid (exposition in \cite{Weinstein96})
  then its {\it nerve}
  $N(\mathcal{G}) \,\in \, \SimplicialSets^{\mathrm{fib}}$
  (\cite[\S 2]{Segal68})
  is the Kan complex \eqref{KanComplexes}
  whose $k$-cells are the sequences of $k$ composable morphisms in $\mathcal{G}$.
  \begin{equation}
    \label{NerveOfAGroupoid}
    N(\mathcal{G})
    \,:\,
    [k]
      \,\mapsto\,
    \mathcal{G}_1
      \;{}_{t}\!\!\times_s
    \mathcal{G}_1
      \;{}_{t}\!\!\times_s
      \cdots
      \;{}_{t}\!\!\times_s
    \mathcal{G}_1
    \,.
  \end{equation}

  \noindent
  {\bf (ii)}
  For $S \in \Sets$ any set, consider its
  {\it pair groupoid}
  $\mathrm{Pair}(S) := \big( S \times S \rightrightarrows S\big)$
  whose objects are the elements of $S$ and which has
  exactly one morphism $s_0 \xrightarrow{\exists !} s_1$
  between any pair of elements. Its nerve \eqref{NerveOfAGroupoid}
  is  contractible, in that it is
  weakly equivalent
  in the classical model category (Def. \ref{ClassicalModelStructureOnSimplicialSets})
  to the point (the terminal simplicial set, which is constant on the singleton set):
  \begin{equation}
    \label{NerveOfPairGroupoidIsContractible}
    \begin{tikzcd}
    N\big(\mathrm{Pair}(S)\big)
    \ar[
      r,
      "\in \mathrm{W} \cap \mathrm{Fib}"
    ]
    &
    \ast
    \,.
    \end{tikzcd}
  \end{equation}
\end{example}

\begin{example}[Opposite model category {\cite[\S 7.1.8]{Hirschhorn02}}]
  \label{OppositeModelStructure}
  If $\mathbf{C}$ is a model category (Def. \ref{ModelCategories})
  then the opposite underlying category becomes a model
  category $\mathbf{C}^{\mathrm{op}}$ with the same
  weak equivalences (up to reversal) and with fibrations
  (resp. cofibrations) the cofibrations (resp. fibrations)
  of $\mathbf{C}$, up to reversal.
\end{example}

\begin{example}[Slice model categories {\cite[\S 7.6.4]{Hirschhorn02}\cite[Thm. 165.3.6]{MayPonto12}}]
  \label{SliceModelCategory}
  Let $\mathbf{C}$ be a model category (Def. \ref{ModelCategories})

  \noindent
  {\bf (i)}
  For $X \in \mathbf{C}$ any object,
  the slice category $\mathbf{C}^{\!/X}$,
  whose objects are morphisms to $X$ and whose morphisms
  are commuting triangles in $\mathbf{C}$ over $X$
  \vspace{-3mm}
    $$
    \mathbf{C}^{/X}
    \big(
      a
      \,,\,
      b
    \big)
    \;\;\;
    :=
    \;\;\;
    \left\{\!\!\!\!
    \raisebox{8pt}{
    \xymatrix@R=5pt{
      A
      \ar[dr]_-{ a }
      \ar@{-->}[rr]^-{f}
      &&
      B
      \ar[dl]^-{ b }
      \\
      & X
    }
    }
   \!\! \right\}
  $$

  \vspace{0mm}
  \noindent
  becomes itself a
  model category, whose weak equivalence, fibrations and
  cofibrations are those morphims whose underlying morphisms
  $f$
  are such in $\mathbf{C}$. This means in particular that:
  \begin{equation}
    \label{CoFibrantObjectsInSliceModelStructure}
    a
      \,\in\,
    \big(
      \mathbf{C}^{/X}
    \big)^{\mathrm{cof}}
    \;\;\;\Leftrightarrow\;\;\;
    A \,\in\, \mathbf{C}^{\mathrm{cof}}
    {\phantom{AAAAA}}
    \mbox{and}
    {\phantom{AAAAA}}
    b \,\in\,
    \big(
      \mathbf{C}^{/X}
    \big)_{\mathrm{fib}}
    \;\;\;\Leftrightarrow\;\;\;
    b \in \mathbf{C}^{\mathrm{fib}}
    \,.
  \end{equation}

  \noindent
  {\bf (ii)} Dually there is the coslice model category
  $
    \mathbf{C}^{X/}
    \;:=\;
    \big(
    (
      \mathbf{C}^{\mathrm{op}}
    )^{/X}
    \big)^{\!\!\mathrm{op}}
  $,
  being the opposite model category
  (Example \ref{OppositeModelStructure}) of the slice category
  of the opposite of $\mathbf{C}$:
  \vspace{-3mm}
  $$
    \mathbf{C}^{X/}
    (
      a
      \,,\,
      b
    )
    \;\;\;
    :=
    \;\;\;
    \left\{\!\!\!\!
    \raisebox{15pt}{
    \xymatrix@R=5pt{
      & X
      \\
      A
      \ar@{<-}[ur]^-{ a }
      \ar@{-->}[rr]_-{f}
      &&
      B
      \ar@{<-}[ul]_-{ b }
    }
    }
    \!\!\right\}.
  $$
\end{example}

\medskip

\noindent {\bf Homotopy categories.}

\begin{defn}[Cylinder objects and Path space objects {\cite[Def. I.4]{Quillen67}}]
  \label{PathSpaceObject}
  Let $\mathbf{C}$ be a model category (Def. \ref{ModelCategories}).

  \noindent
  {\bf (a)}
  With
  $A \in \mathbf{C}_{\mathrm{fib}}$ a fibrant object
  (Notation \ref{FibrantAndCofibrantObjects}),
  a \emph{path space object} for $A$ is
  a factorization
  of the diagonal morphism $\Delta_A$
  through a weak equivalence followed by a fibration:
  \vspace{-2mm}
  \begin{equation}
    \label{PathSpaceFactorization}
    \xymatrix@C=3em{
      A
      \ar@/_1.2pc/[rrr]_-{ \Delta_A }
      \ar[r]^-{\; \in\, \mathrm{W} \;}
      &
      \mathrm{Paths}(A)
      \ar[rr]^-{ (p_0, p_1) \in\, \mathrm{Fib} \;}
      &&
      A \times A
    }.
  \end{equation}

  \noindent
  {\bf (b)}
  With
  $X \in \mathbf{C}_{\mathrm{cof}}$ a cofibrant object
  (Notation \ref{FibrantAndCofibrantObjects}),
  a \emph{cylinder object} for $X$ is
  a factorization
  of the co-diagonal morphism $\nabla_A$
  through cofibration followed by a weak equivalence:
  \vspace{-3mm}
  \begin{equation}
    \label{CylinderObjectFactorization}
    \xymatrix@C=3em{
      X \sqcup X
      \ar[rr]^-{ (i_0, i_1) \in\, \mathrm{Cof} \;}
      \ar@/_1.2pc/[rrr]_-{ \nabla_X }
      &&
      \mathrm{Cyl}(X)
      \ar[r]^-{\; \in\, \mathrm{W} \;}
      &
      X
    }.
  \end{equation}
\end{defn}
\begin{example}[Standard cylinder object in simplicial sets]
  \label{StandardCylinderObjectInSimplicialSets}
  For $X \in \SimplicialSets_{\mathrm{Q}}$
  (Example \ref{ClassicalModelStructureOnSimplicialSets})
  a cylinder object (Def. \ref{PathSpaceObject})
  is evidently given by Cartesian product
  $X \times \Delta[1]$ with the 1-simplex,
  with $(i_0, i_1)$ being the two endpoint inclusions.
\end{example}

\begin{defn}[Homotopy]
  \label{RightHomotopy}
  Let $\mathbf{C}$ be a model category (Def. \ref{ModelCategories}),
  $X \in \mathbf{C}^{\mathrm{cof}}$
  a cofibrant object,
  $A \in \mathbf{C}_{\mathrm{fib}}$
  a fibrant object (Notation \ref{FibrantAndCofibrantObjects}).
  Then a \emph{homotopy} between a pair of morphisms
  $
    f,g
    \in
    \mathbf{C}(X,A)
  $,
  to be denoted
  \vspace{-4mm}
  $$
    \phi \;:\; f \;\Rightarrow\; g
    \phantom{AAAAA}
    \mbox{or}
    \phantom{AAAAA}
    \xymatrix{
      X
      \ar@/^1.1pc/[rr]^-{ f }_-{\ }="s"
      \ar@/_1.1pc/[rr]_-{ g }^-{\ }="t"
      &&
      A
      \ar@{=>}^-{\phi} "s"; "t"
    }
  $$
  is
  a morphism
  $
    \phi_l \;\in\; \mathbf{C}\big(\mathrm{Cyl}(X), A \big)
  $
  out of a cylinder object for $X$
  {\it or}
  a morphism
  $
    \phi_r \;\in\; \mathbf{C}\big(X, \mathrm{Paths}(A)\big)
  $
  to a path space object for $A$
  $\mathrm{Paths}(A)$  (Def. \ref{PathSpaceObject})
  which make either of these diagrams commute:
  \vspace{-2mm}
  $$
    \raisebox{25pt}{
    \xymatrix@R=1.3em{
      X
      \ar[drr]^-{f}
      \ar[d]_-{ i_0 }
      &&
      \\
      \mathrm{Cyl}(X)
      \ar[rr]|-{\;\phi_l\;}
      &&
      A
      \,,
      \\
      X
      \ar[urr]_-{g}
      \ar[u]^-{ i_1 }
    }
    }
    \phantom{A}
    \phantom{AAAAAAA}
    \raisebox{25pt}{
    \xymatrix@R=1.3em{
      && A
      \\
      X
      \ar[rr]|-{\;\phi_r\;}
      \ar[urr]^-{f}
      \ar[drr]_-{g}
      &&
      \mathrm{Paths}(A)
      \ar[u]_-{ p_0 }
      \ar[d]^-{ p_1 }
      \\
      &&
      A
    }
    }
  $$
\end{defn}

\begin{prop}[Homotopy classes]
  \label{RightHomotopyClasses}
  Let $\mathbf{C}$ be a model category,
  $X \in \mathbf{C}^{\mathrm{cof}}$ and
  $A \in \mathbf{C}_{\mathrm{fib}}$ (Notation \ref{FibrantAndCofibrantObjects}).
  Then homotopy (Def. \ref{RightHomotopy}) is an
  equivalence relation $\sim$ on
  the hom-set $\mathbf{C}(X,A)$.
  We write
     \vspace{-2mm}
  \begin{equation}
    \label{RightHomotopyClasses}
    \mathbf{C}(X,A)_{\!/\sim}
    \;\;
    \in
    \;\;
    \mathrm{Sets}
  \end{equation}

  \vspace{-2mm}
  \noindent
  for the corresponding set of {\it homotopy classes}
  of morphisms from $X$ to $A$.
\end{prop}

\begin{defn}[Homotopy category of a model category]
  \label{HomotopyCategory}
  For $\mathbf{C}$ a model category (Def. \ref{ModelCategories}),

  \noindent {\bf (i)} we write
  \vspace{-2mm}
  \begin{equation}
    \label{HomSetsInHomotopyCategory}
    \mathrm{Ho}
    (
      \mathbf{C}
    )
    \;:=\;
    \big(
      \mathbf{C}^{\mathrm{cof}}_{\mathrm{fib}}
    \big)_{\!/\sim_r}
    \;\;\in\;\;
    \Categories
  \end{equation}

  \vspace{-1mm}
  \noindent
  for the category whose objects are those objects
  of $\mathbf{C}$ that are both fibrant and cofibrant
  (Notation \ref{FibrantAndCofibrantObjects}), and
  whose morphisms are the right homotopy classes
  of morphisms in $\mathbf{C}$ (Def. \ref{RightHomotopyClasses}):
  \vspace{-2mm}
  $$
    X,A \in \mathbf{C}_{\mathrm{fib}}^{\mathrm{cof}}
    \;\;\;\;
      \Rightarrow
    \;\;\;\;
    \mathrm{Ho}
    (
      \mathbf{C}
    )
    (
      X\,,\,A
    )
    \;:=\;
    \mathbf{C}(X,A)_{\!/\sim_r}
  $$

  \vspace{-2mm}
  \noindent
  and composition of morphisms is induced from composition of representatives in $\mathbf{C}$.

 \noindent {\bf (ii)} Given a choice of fibrant replacement $P$
  and of cofibrant replacement $Q$
  for each object of $\mathbf{C}$ (Notation \ref{FibrantAndCofibrantObjects}) we obtain a functor
   \vspace{-2mm}
  \begin{equation}
    \label{LocalizationOfAModelCategoryAtWeakEquivalenes}
    \xymatrix{
      \mathbf{C}
      \ar[r]^-{ \gamma_{\mathbf{C}} }
      &
      \mathrm{Ho}
      (
        \mathbf{C}
      )
      \,,
    }
  \end{equation}

  \vspace{-4mm}
  \noindent
  which {\bf (a)} sends any object $X \in \mathbf{C}$ to $P Q X$
  and sends {\bf (b)} any morphism $\xymatrix@C=12pt{X \ar[r]^-{f} & A}$
  to the right homotopy class \eqref{RightHomotopyClasses}
  of any lift \eqref{LiftingProperty}
  $P Q f$ obtained from any lift $Q f$ in the following diagrams:
    \vspace{-3mm}
  $$
    \raisebox{20pt}{
    \xymatrix@R=1.5em{
      \varnothing
      \ar[rr]
      \ar[d]_{ i_X }
      &&
      Q Y
      \ar[d]^-{ p_Y }
      \\
      Q X
      \ar@{-->}[urr]|-{ \; Q f \; }
      \ar[rr]_-{
        f \,\circ\, p_x
      }
      &&
      Y
    }
    }
    \phantom{AAA}
    \rightsquigarrow
    \phantom{AAA}
    \raisebox{20pt}{
    \xymatrix@R=1.5em{
      Q X
      \ar[rr]^-{
        j_{Q Y} \,\circ\, Q f
      }
      \ar[d]_-{ j_{Q X} }
      &&
      P Q Y
      \ar[d]^-{ q_{Q Y} }
      \\
      P Q X
      \ar[rr]
      \ar@{-->}[urr]|-{ \; P Q f \; }
      &&
      \ast
    }
    }
  $$
\end{defn}
\begin{prop}[Homotopy category is localization]
  \label{HomotopyCategoryIsLocalization}
  For a model category $\mathbf{C}$ (Def. \ref{ModelCategories})
  the functor
  $\xymatrix@C=9pt{
    \mathbf{C}
    \ar[r]^-{ \gamma_{\mathbf{C}} }
    &
    \mathrm{Ho}\big( \mathbf{C} \big)
  }$ \eqref{LocalizationOfAModelCategoryAtWeakEquivalenes}
  from Def. \ref{HomotopyCategory} exhibits the
  homotopy category as the \emph{localization}
  of the model category at its class of weak equivalences:
  $\gamma_{\mathbf{C}}$ sends all weak equivalences
  in $\mathbf{C}$ to
  isomorphisms, and is the universal functor with this property.
\end{prop}

The restriction to fibrant-and-cofibrant objects in
Def. \ref{HomotopyCategory} is convenient for defining
composition of morphisms, but for computing hom-sets in
the homotopy category it is sufficient that the domain object
is cofibrant, and the codomain fibrant:

\begin{prop}[{\cite[\S I.1 Cor. 7]{Quillen67}}]
  \label{RightHomotopyReflectsHomotopyClasses}
  Let $\mathcal{C}$ be a model category (Def. \ref{ModelCategories}).
  For
  $X \in \mathcal{C}^{\mathrm{cof}}$ a cofibrant object
  and $A \in \mathcal{C}_{\mathrm{fib}}$ a fibrant object,
  any choice of fibrant replacement $P X$
  and cofibrant replacement $Q A$
  (Notation \ref{FibrantAndCofibrantObjects}).
  induces a bijection between the
  set of homotopy classes (Def. \ref{RightHomotopy})
  and the hom-set in the homotopy category (Def. \ref{HomotopyCategory})
  between $X$ and $A$:
  \vspace{-2mm}
  $$
    \xymatrix{
      \mathcal{C}(X,A)_{\!/\sim_r}
      \ar[rr]^-{ \mathcal{C}(j_X, p_A) }_-{\simeq}
      &&
      \mathrm{Ho}
      (
        \mathbf{C}
      )
      (X,A)
    }.
  $$
\end{prop}
While the hom-functor of a homotopy category preserves almost no
homotopy (co)limits, we do have:
\begin{prop}[Hom-functor of homotopy category respects (co)products]
  \label{HomFunctorOfAHomotopyCategoryRespectsCoProducts}
  The hom-functor of a homotopy category (Def. \ref{HomotopyCategory})
  respects
  coproducts in the first argument and products in its second argument,
  in that there are natural bijections of the following form:
   $$
     \mathrm{Ho}(\mathbf{C})
     \big(
       \underset{i \in I}{\coprod} X_i
       ,\,
       \underset{j \in J}{\prod} A_j
     \big)
     \;\;\simeq\;\;
     \underset{
       {i \in I}
       \atop
       {j \in J}
     }{\prod}
     \,
     \mathrm{Ho}(\mathbf{C})
     \big(
       X_i
        ,\,
       A_j
     \big)
   $$
\end{prop}
\begin{proof}
Noticing that
coprodcts preserve cofibrancy
and
products preserve fibrancy, evidently,
this follows from Prop. \ref{RightHomotopyReflectsHomotopyClasses}.
\end{proof}

\noindent {\bf Quillen adjunctions.}

\begin{defn}[Quillen adjunction]
  \label{QuillenAdjunction}
  Let $\mathbf{D}, \mathbf{C}$ be model categories
  (Def. \ref{ModelCategories}).
  Then a pair of adjoint functors $(L \dashv R)$ \eqref{AnAdjunction}
  between their underlying categories
  is called a
  \emph{Quillen adjunction}, to be denoted
  \vspace{-3mm}
  \begin{equation}
    \label{AQuillenAdjunction}
    \xymatrix{
      \mathbf{D}
      \ar@{<-}@<+7pt>[rr]^-{L}
      \ar@<-7pt>[rr]_-{R}^-{\bot_{\mathrlap{\mathrm{Qu}}}}
      &&
      \mathbf{C}
    }
  \end{equation}

  \vspace{-3mm}
  \noindent
  if the following equivalent conditions hold:
  \begin{itemize}
    \vspace{-.15cm}
    \item $L$ preserves $\mathrm{Cof}$, and $R$ preserves $\mathrm{Fib}$;
    \vspace{-.15cm}
    \item $L$ preserves $\mathrm{Cof}$ and $\mathrm{Cof} \cap \mathrm{W}$ (``left Quillen functor'');
    \vspace{-.15cm}
    \item $R$ preserves $\mathrm{Fib}$ and $\mathrm{Fib} \cap \mathrm{W}$ (``right Quillen functor'').
  \end{itemize}
\end{defn}

\begin{example}[Base change Quillen adjunction]
  \label{BaseChangeQuillenAdjunction}
  Let $\mathbf{C}$ be a model category (Def. \ref{ModelCategories}),
  $B_1, B_2 \,\in\, \mathbf{C}_{\mathrm{fib}}$ a pair of
  fibrant objects (Notation \ref{FibrantAndCofibrantObjects})
  and
  \vspace{-3mm}
   \begin{equation}
    \label{MorphismForBaseChange}
    \xymatrix{
      B_1 \ar[r]^-{ f } & B_2
    }
    \;\;\;
    \in \mathbf{C}
  \end{equation}

  \vspace{-2mm}
  \noindent
  a morphism. Then we have a Quillen adjunction
  (Def. \ref{QuillenAdjunction})
  \vspace{-2mm}
  \begin{equation}
  \label{BaseChangeQuillenAdjointFunctor}
    \xymatrix{
      \mathbf{C}^{/B_2}
      \ar@<+7pt>@{<-}[rr]^-{ f_! }
      \ar@<-7pt>[rr]_-{ f^\ast }^-{ \bot_{\mathrlap{\mathrm{Qu}}} }
      &&
      \mathbf{C}^{/B_1}
    }
  \end{equation}

    \vspace{-2mm}
  \noindent
  between the slice model categories (Example \ref{SliceModelCategory}),
  where:

  {\bf (i)} The left adjoint functor $f_!$ is given by postcomposition
  in $\mathbf{C}$ with $f$ \eqref{MorphismForBaseChange}:
  \vspace{-3mm}
    \begin{equation}
    \label{LeftBaseChange}
    f_\ast
    \;\;:\;\;
    \raisebox{10pt}{
    \xymatrix@R=4pt{
      X
      \ar[dr]_-{\tau}
      \ar[rr]^-{ c }
      &&
      A
      \ar[dl]^-{ p }
      \\
      & B_1
    }
    }
    \;\;\;\;\;\;\;\;\;
    \longmapsto
    \;\;\;\;\;\;\;\;\;
    \raisebox{10pt}{
    \xymatrix@R=4pt{
      X
      \ar[dr]_-{\tau}
      \ar[rr]^-{ c }
      &&
      A
      \ar[dl]^-{ p }
      \\
      & B_1
      \ar[dd]^-{f}
            \\
            \\
      & B_2
    }
    }
  \end{equation}

  \vspace{-.2cm}
{\bf (ii)} The right adjoint functor $f^\ast$ is given  by pullback
  \eqref{CartesianSquare} along $f$ \eqref{MorphismForBaseChange}.

  \noindent
  That these functors indeed form an adjunction
  $f_! \dashv f^\ast$ follows from the defining universal
  property of the pullback \eqref{CartesianSquare}:
  \begin{equation}
    \label{BaseChangeAdjunctionHomIsomorphism}
    \arraycolsep=1.4pt\def\arraystretch{.1}
    \begin{array}{ccc}
          \mathbf{C}^{/B_2}
      \big(
        f_\ast \tau
        \,,\
        \rho
      \big)
      &\simeq&
      \mathbf{C}^{/B_1}
      \big(
        \tau
        \,,\
        f^\ast \rho
      \big)
      \\
      \\
      \scalebox{.9}{$
      \raisebox{30pt}{
      \xymatrix@R=1pt@C=3em{
        X
        \ar@{-->}[rr]^-{ c }
        \ar[ddddr]_-{ \tau }
        &
        &
        A
        \ar[ddddd]^-{ \rho }
        \\
        \\
        \\
        \\
        &
        B_1
        \ar[dr]_-{ f }
        \\
        &&
        B_2
      }
      }
      $}
      &
      \;\;\;\;\;
      \leftrightarrow
      \;\;\;\;\;
      &
      \scalebox{.9}{$
      \raisebox{30pt}{
      \xymatrix@R=1pt@C=3em{
        X
        \ar@{-->}[r]^-{ \widetilde {c} }
        \ar[ddddr]_-{ \tau }
        &
        f^\ast A
        \ar[dddd]|-{
          \mathclap{\phantom{\vert^{\vert}}}
          f^\ast \rho
          \mathclap{\phantom{\vert_{\vert}}}
        }
        \ar[r]
        \ar@{}[dddddr]|-{ \mbox{\tiny(pb)} }
        &
        A
        \ar[ddddd]^-{ \rho }
        \\
        \\
        \\
        \\
        &
        B_1
        \ar[dr]_-{ f }
        \\
        &&
        B_2
      }
      }
      $}
    \end{array}
  \end{equation}

  \noindent
  That this adjunction is a Quillen adjunction (Def. \ref{QuillenAdjunction})
  follows since $f_!$ \eqref{LeftBaseChange} evidently preserves each
  of $\mathrm{W}$ and  $\mathrm{Cof}$ (even $\mathrm{Fib}$)
  separately,
  by Example \ref{SliceModelCategory}.
\end{example}
\begin{example}[Sliced Quillen adjunction]
  \label{SlicedQuillenAdjunction}
  Given a Quillen adjunction $L \dashv_{{}_{\mathrm{Qu}}} R$ (Def. \ref{QuillenAdjunction})
  and an object of either of the two model categories,
  there is an induced {\it sliced} Quillen adjunctions (Def. \ref{QuillenAdjunction})
  between slice model categories (Ex. \ref{SliceModelCategory})
  as follows:
  \vspace{-2mm}
  \begin{equation}
    \label{FormOfSlicedQuillenAdjunctions}
    \begin{tikzcd}[column sep=small]
      \mathbf{D}
      \ar[
        rr,
        shift right=7pt,
        "{ R }"{below}
      ]
      \ar[
        rr,
        phantom,
        "\scalebox{.7}{$\mathrm{\bot_{\mathrlap{\mathrm{Qu}}}}$}"
      ]
      &&
      \mathbf{C}
      \ar[
        ll,
        shift right=7pt,
        "{ L }"{above}
      ]
    \end{tikzcd}
    {\phantom{AAAA}}
    \Rightarrow
    {\phantom{AAAA}}
    \left(
    \underset{
      c \in \mathbf{C}
    }{\forall}
    \;\;\;
    \begin{tikzcd}[column sep=small]
      \mathbf{D}^{/L(c)}
      \ar[
        rr,
        shift right=7pt,
        "{ R^{/c} }"{below}
      ]
      \ar[
        rr,
        phantom,
        "\scalebox{.7}{$\mathrm{\bot_{\mathrlap{\mathrm{Qu}}}}$}"
      ]
      &&
      \mathbf{C}^{/c}
      \ar[
        ll,
        shift right=7pt,
        "{ L^{/c} }"{above}
      ]
    \end{tikzcd}
    \right)
    \;\;\;
    \mbox{and}
    \;\;\;
    \left(
    \underset{
      d \in \mathbf{D}
    }{\forall}
    \;\;\;
    \begin{tikzcd}[column sep=small]
      \mathbf{D}^{/d}
      \ar[
        rr,
        shift right=7pt,
        "{ R^{/d} }"{below}
      ]
      \ar[
        rr,
        phantom,
        "\scalebox{.7}{$\mathrm{\bot_{\mathrlap{\mathrm{Qu}}}}$}"
      ]
      &&
      \mathbf{C}^{/R(d)}
      \ar[
        ll,
        shift right=7pt,
        "{ L^{/d} }"{above}
      ]
    \end{tikzcd}
    \right)
    \,,
  \end{equation}
  where:

  \begin{itemize}

  \vspace{-.2cm}
  \item[\bf (i)]
  $L^{/c}$ and $R^{/d}$ are given directly by applying $L$ or $R$, respectively,
   to the triangular diagram that defines a morphism in the slice;

   \vspace{-.2cm}
   \item[\bf (ii)]
   $R^{/c}$ and $L^{/d}$ are given by this direct application followed by
   right/left base change \eqref{BaseChangeQuillenAdjunction}
   along the adjunction unit/counit \eqref{AdjunctionUnit}, respectively:

   \vspace{-4mm}
   \begin{equation}
     \label{TheNontrivialAdjointsInTheSliceAdjunctions}
     \begin{tikzcd}[column sep=small]
     R^{/c}
     \;\colon\;
     \mathbf{D}^{/L(c)}
     \ar[
       rr,
       "R"
     ]
     &&
     \mathbf{C}^{/R \circ L(c)}
     \ar[
       rr,
       "(\eta_c)^\ast"
     ]
     &&
     \mathbf{D}^{ /c }
     \mathrlap{\,,}
     \end{tikzcd}
     {\phantom{AAAAAA}}
     \begin{tikzcd}[column sep=small]
       \mathbf{D}^{/d}
       &&
       \mathbf{C}^{/L \circ R(d)}
       \ar[
         ll,
         "(\epsilon_d)_! "{above}
       ]
       &&
       \mathbf{D}^{/R(d)}
       \;\colon\;
       L^{/d}
       \ar[
         ll,
         "L"{above}
       ]
     \end{tikzcd}
   \end{equation}
   \vspace{-.5cm}

   \noindent
   In particular, this means that $L^{/d}$ sends a slicing morphism $\tau$ to its adjunct
   $\widetilde{\tau}$ \eqref{AdjunctsByCompositesWithCoUnits}, in that:
   \begin{equation}
     \label{LeftSliceAdjointProducesAdjuncts}
     L^{/d}
     \left(\!\!\!
       \begin{tikzcd}
         c
         \ar[d, "\tau"]
         \\
         R(d)
       \end{tikzcd}
    \!\!\! \right)
     \;\;
     =
     \;\;
     \left(\!\!\!\!
       \begin{tikzcd}
         L(c)
         \ar[d, "\widetilde{\tau}"]
         \\
         d
       \end{tikzcd}
     \!\!\!\! \right)
     \;\;\;
     \in
     \;
     \mathbf{D}^{/d}
     \,.
   \end{equation}
   \end{itemize}
\end{example}
Aspects of this statement appear in \cite[Prop. 5.2.5.1]{Lurie09}\cite[Prop. 2.5(2)]{Li14}.
Since it is key to the proof of the twisted non-abelian de Rham theorem
(Th. \ref{TwistedNonAbelianDeRhamTheorem}) we spell it out:
\begin{proof}
It is clear that if we have adjunctions as claimed in
\eqref{FormOfSlicedQuillenAdjunctions}, then
$L^{/c}\,$/$\,R^{/d}$
are left/right Quillen functors, respectively, since
these two act as the left/right Quillen functors $L$/$R$ on underlying morphisms
(by item (i) above),
where the classes of slice morphisms are created, by Ex. \ref{SliceModelCategory}.

To see that we have adjunctions as claimed, we may check their hom-isomorphisms
\eqref{AnAdjunction} (for readability we now denote the object being sliced
over by ``$b$'', in both cases, with ``$c$'' and ``$d$'' now being the
variables in the hom-isomorphism):

\noindent
{\bf (1)} For the first case,
consider the following transformations of slice hom-sets:

\vspace{-.7cm}
$$
\def\arraystretch{3}
\arraycolsep=8pt
\begin{array}{ccc}
  \begin{tikzcd}[column sep={between origins, 33}, row sep=4]
    L(c)
    \ar[rr, "f", dashed]
    \ar[dr]
      &&
    d
    \ar[dl, "p"]
    \\
    &
    L(b)
  \end{tikzcd}
  &
  \longmapsto
  &
  \begin{tikzcd}[column sep={between origins, 33}, row sep=4]
    c
    \ar[
      rr,
      "{\eta_c}"
    ]
    \ar[dr]
    \ar[
      rrrr,
      rounded corners,
      to path={
           -- ([yshift=+12pt]\tikztostart.north)
           --node[above]{
               \scalebox{.7}{$
                 \widetilde{f}
               $}
             } ([yshift=+8pt]\tikztotarget.north)
           -- (\tikztotarget.north)}
    ]
    &&
    R \circ L(c)
    \ar[rr, "{R(f)}", dashed]
    \ar[dr]
      &&
    R(d)
    \ar[dl, "R(p)"]
    \\
    &
    b
    \ar[
      rr,
      "\eta_b"{below}
    ]
    &
    &
    R \circ L(b)
  \end{tikzcd}
  \\
  \updownarrow
  &&
  \updownarrow
  \\
  \begin{tikzcd}[column sep={between origins, 33}, row sep=4]
    L(c)
    \ar[dr]
    \ar[rrrr, "L(\widetilde{f})"]
    \ar[
      rrrrrr,
      rounded corners,
      to path={
           -- ([yshift=+8pt]\tikztostart.north)
           --node[above]{
               \scalebox{.7}{$
                 f
               $}
             } ([yshift=+8pt]\tikztotarget.north)
           -- (\tikztotarget.north)}
    ]
    &&
    &&
    L \circ R(d)
    \ar[rr, "{\epsilon_d}"]
    \ar[dl]
    &&
    d
    \ar[dl, "p"]
    \\
    &
    L(b)
    \ar[
      rr,
      "L(\eta_b)"{below}
    ]
    \ar[
      rrrr,
      rounded corners,
      to path={
           -- ([yshift=-8pt]\tikztostart.south)
           --node[below]{
               \scalebox{.7}{$
                 \mathrm{id}
               $}
             } ([yshift=-8pt]\tikztotarget.south)
           -- (\tikztotarget.south)}
    ]
    &
    &
    L \circ R \circ L(b)
    \ar[rr, "\epsilon_{L(b)}"{below}]
    &&
    L(b)
  \end{tikzcd}
  &\rotatebox{180}{$\longmapsto$}&
  \begin{tikzcd}[column sep={between origins, 33}, row sep=4]
    c
    \ar[
      rr,
      dashed
    ]
    \ar[dr]
      \ar[
      rrrr,
      rounded corners,
      to path={
           -- ([yshift=+12pt]\tikztostart.north)
           --node[above]{
               \scalebox{.7}{$
                 \widetilde{f}
               $}
             } ([yshift=+8pt]\tikztotarget.north)
           -- (\tikztotarget.north)}
      ]
    &&
    \eta_b^\ast\big(R(d)\big)
    \ar[rr]
    \ar[dl]
    \ar[
      dr,
      phantom,
      "\mbox{\tiny\rm(pb)}"
    ]
      &&
    R(d)
    \ar[dl, "R(p)"]
    \\
    &
    b
    \ar[
      rr,
      "\eta_b"{below}
    ]
    &
    &
    R \circ L(b)
    \mathrlap{\,.}
  \end{tikzcd}
\end{array}
$$
Here the horizontal transformations are given by
applying the functors and then (post-)composing with (co-)units,
while the left vertical bijection is the formula \eqref{AdjunctsByCompositesWithCoUnits}
for adjuncts and the right vertical bijection is the universal property of the
pullback. Evidently these operations commute in both possible ways, showing that
also the horizontal operations are bijections
(and they are natural by the naturality
of the underlying hom-isomorphism).

{\bf (2)} The second case follows analogously, but more directly
as no pullback is involved here:
$$
\arraycolsep=2pt
\begin{array}{ccccc}
\begin{tikzcd}[column sep={between origins, 33}, row sep=4]
  c
  \ar[rr, dashed, "f"]
  \ar[
    dr,
    "\tau"{left, xshift=-3pt}, pos=.8
  ]
  &&
  R(d)
  \ar[dl]
  \\
  &
  R(b)
\end{tikzcd}
&\mapsto&
\begin{tikzcd}[column sep={between origins, 33}, row sep=4]
  L(c)
  \ar[rr, dashed, "L(f)"]
  \ar[dr]
  \ar[
      rrrr,
      rounded corners,
      to path={
           -- ([yshift=+8pt]\tikztostart.north)
           --node[above]{
               \scalebox{.7}{$
                 \widetilde{f}
               $}
             } ([yshift=+8pt]\tikztotarget.north)
           -- (\tikztotarget.north)}
  ]
  &&
  L \circ R(d)
  \ar[rr, "\epsilon_d"{above}]
  \ar[dl]
  &&
  d
  \ar[dl]
  \\
  &
  L \circ R(b)
  \ar[rr, "\epsilon_b"{below}]
  &&
  b
\end{tikzcd}
&\mapsto&
\begin{tikzcd}[column sep={between origins, 33}, row sep=4]
  c
  \ar[rr, "\eta_c"{above}]
  \ar[dr]
  \ar[
      rrrrrr,
      rounded corners,
      to path={
           -- ([yshift=+12pt]\tikztostart.north)
           --node[above]{
               \scalebox{.7}{$
                 f
               $}
             } ([yshift=+8pt]\tikztotarget.north)
           -- (\tikztotarget.north)}
  ]
  &&
  R \circ L(c)
  \ar[rrrr,"R(\widetilde{f})"]
  \ar[dr]
  &&
  %R \circ L \circ R(d)
  %\ar[rr, "R(\epsilon_d)"{above}]
  %\ar[dl]
  &&
  R(d)
  \ar[dl]
  \\
  &
  R(b)
  \ar[rr, "\eta_{R(b)}"{below}]
    \ar[
      rrrr,
      rounded corners,
      to path={
           -- ([yshift=-8pt]\tikztostart.south)
           --node[below]{
               \scalebox{.7}{$
                 \mathrm{id}
               $}
             } ([yshift=-8pt]\tikztotarget.south)
           -- (\tikztotarget.south)}
    ]
  &&
  R \circ L \circ R(b)
  \ar[rr, "R(\epsilon_b)"{below}]
  &&
  R(b)
\end{tikzcd}
\end{array}
$$
\vspace{-.9cm}

\end{proof}
\begin{example}[Induced Quillen adjunction on pointed objects]
  \label{InducedQuillenAdjunctionOnPointedObjects}
  Given a Quillen adjunction $L \dashv_{{}_{\mathrm{Qu}}} R$ (Def. \ref{QuillenAdjunction})
  there is an induced Quillen adjunction of model categories
  of pointed objects, hence of coslice model structures
  (Ex. \ref{SliceModelCategory}) under the terminal object
  \vspace{-2mm}
  \begin{equation}
    \label{InducedQuillenAdjunctionOnPointedObjectsViaOppositeSlicing}
    \begin{tikzcd}[column sep=small]
      \mathbf{D}
      \ar[
        rr,
        shift right=8pt,
        "R"{below}
      ]
      \ar[
        rr,
        phantom,
        "\scalebox{.7}{$\bot_{\mathrlap{\mathrm{Qu}}}$}"
      ]
      &&
      \mathbf{C}
      \ar[
        ll,
        shift right=8pt,
        "L"{above}
      ]
    \end{tikzcd}
    {\phantom{AAAA}}
    \Rightarrow
    {\phantom{AAAA}}
    \left(
    \begin{tikzcd}[column sep=small]
      \mathbf{D}^{\ast/}
      \ar[
        rr,
        shift right=8pt,
        "R^{\ast/}"{below}
      ]
      \ar[
        rr,
        phantom,
        "\scalebox{.7}{$\bot_{\mathrlap{\mathrm{Qu}}}$}"
      ]
      &&
      \mathbf{C}^{\ast/}
      \ar[
        ll,
        shift right=8pt,
        "L^{\ast/}"{above}
      ]
    \end{tikzcd}
    \right)
    \;\;\;\;
    =
    \;\;\;\;
    \left(
    \begin{tikzcd}[column sep=small]
      \big(
        \mathbf{D}^{\mathrm{op}}_{/\ast}
      \big)^{\mathrm{op}}
      \ar[
        rr,
        shift right=8pt,
        "\big( R^{\mathrm{op}}_{/\ast}\big)^{\mathrm{op}}"{below}
      ]
      \ar[
        rr,
        phantom,
        "\scalebox{.7}{$\bot_{\mathrlap{\mathrm{Qu}}}$}"
      ]
      &&
      \big(
        \mathbf{C}^{\mathrm{op}}_{/R(\ast)}
      \big)^{\mathrm{op}}
      \ar[
        ll,
        shift right=8pt,
        " \big( L^{\mathrm{op}}_{/\ast} \big)^{\mathrm{op}} "{above}
      ]
    \end{tikzcd}
    \right)
    \mathrlap{\,,}
  \end{equation}

  \vspace{-.3cm}
  \noindent
  where

  \vspace{-.1cm}
  \begin{itemize}

  \vspace{-.2cm}
  \item[\bf (i)]
  $R^{\ast/}$ is given directly by applying $R$ to the
  underling triangular diagrams in $\mathbf{D}$;

  \vspace{-.2cm}
  \item[\bf (ii)]
  $L^{\ast/}$ is given by that direct application of $L$ followed
  (using that the right adjoint $R$ preserves the terminal object) by
  pushout along the adjunction counit
  $ L(\ast) \simeq L \circ R(\ast) \xrightarrow{\; \epsilon_{\ast} \;} \ast$:

  \vspace{-.4cm}
  \begin{equation}
    \label{InducedLeftAdjointOnPointedObjects}
    L^{\ast/}
    \;\;\colon\;\;
    \mathbf{C}^{\ast/}
    \xrightarrow{\;\; L \;\;}
    \mathbf{D}^{L(\ast)/}
    \,\simeq\,
    \mathbf{D}^{L(\ast)/}
    \xrightarrow{\;\; (-) \underset{ L\circ R(\ast) }{\sqcup} \ast  \;\;}
    \mathbf{D}^{\ast/}
    \,.
  \end{equation}
  \end{itemize}

  This may be checked directly (e.g. \cite[Prop. 1.3.5]{Hovey99}),
  but it is also a special case of Ex. \ref{SlicedQuillenAdjunction},
  as shown on the right of \eqref{InducedQuillenAdjunctionOnPointedObjectsViaOppositeSlicing},
  observing that pullbacks are pushouts in the opposite category.
\end{example}

\begin{lemma}[Ken Brown's lemma {\cite[Lemma 1.1.12]{Hovey99}}\cite{Brown73}]
  \label{KenBrownLemma}
  Given a Quillen adjunction $L \dashv R$ (Def. \ref{QuillenAdjunction}),

  \noindent {\bf (i)} the right Quillen functor $R$ preserves all weak equivalences between fibrant objects.

  \noindent {\bf (ii)} the left Quillen functor $L$ preserves all weak equivalences
  between cofibrant objects.
\end{lemma}

\begin{prop}[Derived functors]
  \label{DerivedFunctors}
  Given a Quillen adjunction $(L \dashv_{\mathrm{Qu}} R)$
  (Def. \ref{QuillenAdjunction}),
  there are adjoint functors
  $\LeftDerived L \,\dashv\, \RightDerived R$\footnote{
    We avoid the common notation $\mathbb{L}L \dashv \mathbb{R}R$
    for derived functors, since this clashes with the prominent role
    that ``$\mathbb{R}$'' plays as notation for the field of real numbers
    in the main text.
  }
  \eqref{AnAdjunction}
  between the
  homotopy categories (Def. \ref{HomotopyCategory})
  \vspace{-2mm}
  \begin{equation}
    \label{DerivedAdjunction}
    \xymatrix{
      \mathrm{Ho}
      (
        \mathbf{D}
      )
        \ar@{<-}@<+6pt>[rr]^-{ \LeftDerived L }
        \ar@<-6pt>[rr]_-{ \RightDerived R }^-{
          \scalebox{.8}{\scalebox{0.8}{$\bot$}}
        }
      &&
      \;
      \mathrm{Ho}
      (
        \mathbf{C}
      )
    }
  \end{equation}

  \vspace{-3mm}
  \noindent
  whose composites with the localization functors \eqref{LocalizationOfAModelCategoryAtWeakEquivalenes}
  make the following squares commute up to natural isomorphism:
  \vspace{-2mm}
  $$
    \raisebox{20pt}{
    \xymatrix@R=1.2em{
      \mathbf{D}
      \ar[rr]^-{ R }
      \ar[d]_-{ \gamma_{\mathbf{D}} }
      \ar@{}[drr]|-{
        \scalebox{.8}{
          \rotatebox[origin=c]{-45}{$\Downarrow$}
        }_{\simeq}
      }
      &&
      \mathbf{C}
      \ar[d]^-{ \gamma_{\mathbf{C}} }
      \\
      \mathrm{Ho}
      (
        \mathbf{D}
      )
      \ar[rr]_-{ \RightDerived R }
      &&
      \mathrm{Ho}
      (
        \mathbf{C}
      )
    }
    }
    \phantom{AAAAA}
    \raisebox{20pt}{
    \xymatrix@R=1.2em{
      \mathbf{D}
      \ar@{<-}[rr]^-{ L }
      \ar[d]_-{ \gamma_{\mathbf{D}} }
      \ar@{}[drr]|-{
        {_{\simeq}}
        \scalebox{.8}{
          \rotatebox[origin=c]{+45}{$\Downarrow$}
        }
      }
      &&
      \mathbf{C}
      \ar[d]^-{ \gamma_{\mathbf{C}} }
      \\
      \mathrm{Ho}
      (
        \mathbf{D}
      )
      \ar@{<-}[rr]_-{ \LeftDerived L }
      &&
      \mathrm{Ho}
      (
        \mathbf{C}
      )
      \,.
    }
    }
  $$

    \vspace{-1mm}
\noindent
  These are unique up to natural isomorphism, and are called
  the left and right \emph{derived functors} of $L$ and $R$,
  respectively.
\end{prop}

\begin{example}[Derived functors via (co-)fibrant replacement]
  \label{DerivedFunctorsByCoFibrantReplacement}
  It is convenient to leave the localization functors
  $\gamma$ \eqref{LocalizationOfAModelCategoryAtWeakEquivalenes}
  notationally implicit, and understand objects of $\mathbf{C}$
  as objects of $\mathrm{Ho}(\mathbf{C})$, via $\gamma$.  Then:
  \begin{itemize}
  \vspace{-.2cm}
  \item[{\bf (i)}] The value of a left derived functor $\LeftDerived L$
  (Prop. \ref{DerivedFunctors}) on an object $c \in \mathbf{C}$
  is equivalently the value of $L$ on a cofibrant replacement
  $Q c$ \eqref{CofibrantReplacement}:
    \vspace{0mm}
  \begin{equation}
    \label{LeftDerivedFunctorByCofibrantReplacement}
    \LeftDerived L(c) \;\simeq\; L(Q c)
    \;\;\;\;\;
    \in
    \;
    \mathrm{Ho}(\mathbf{D})\;.
  \end{equation}

  \vspace{-.3cm}
  \item[{\bf (ii)}] The value of a right derived functor $\RightDerived R$
  (Prop. \ref{DerivedFunctors}) on an object $d \in \mathbf{D}$
  is equivalently the value of $R$ on a fibrant replacement
  $P d$ \eqref{FibrantReplacement}:
    \vspace{0mm}
  \begin{equation}
    \label{RightDerivedFunctorByFibrantReplacement}
    \RightDerived R(d) \;\simeq\; R(P d)
    \;\;\;\;\;
    \in
    \;
    \mathrm{Ho}(\mathbf{C})\;.
  \end{equation}

  \vspace{-.3cm}
  \item[{\bf (iii)}]
  The derived unit $\mathbb{D}\eta$ \eqref{AdjunctionUnit}
  of the
  derived adjunction \eqref{DerivedAdjunction} is, on any
  cofibrant object $c \in \mathbf{C}^{\mathrm{cof}}$,
  given by
  \vspace{-3mm}
      \begin{equation}
        \label{DerivedAdjunctionUnit}
        \mathbb{D}\eta_c
        \;:\;
        \xymatrix{
          c
          \ar[r]^-{ \eta_c }
          &
          R\big(L(c)\big)
          \ar[rr]^-{
            R
            (
              j_{ L(c) }
            )
          }
          &&
          R\big( P L(c)\big)
        }
        \;\;\;\;\;
        \in
        \;
        \mathrm{Ho}(\mathbf{C})
      \end{equation}
     \vspace{-.8cm}

   \noindent
   where  $\!\!\xymatrix{ L(c) \ar[r]^-{j_{L(c)}} & P L(c) }\!\!$
   is  any fibrant replacement \eqref{FibrantReplacement}.

  \vspace{-.2cm}
  \item[{\bf (iv)}]
  The derived co-unit $\mathbb{D}\epsilon$ \eqref{AdjunctionUnit}
  of the
  derived adjunction \eqref{DerivedAdjunction}, is, on any
  fibrant object $d \in \mathbf{D}_{\mathrm{fib}}$,
  given by
       \vspace{-3mm}
      \begin{equation}
        \label{DerivedAdjunctionCoUnit}
        \xymatrix{
          \mathbb{D}\epsilon_d
          \;:\;
          L\big( Q R(d)\big)
          \ar[rr]^-{
            L
            (
              p_{ R(d) }
            )
          }
          &&
          L\big(R(d)\big)
          \ar[r]^-{ \epsilon_d }
          &
          d
        }
        \;\;\;\;\;
        \in
        \;
        \mathrm{Ho}(\mathbf{D})
      \end{equation}
      \vspace{-.8cm}

  \noindent
  where  $\!\!\xymatrix@C=1.8em{ Q R(d) \ar[r]^-{p_{R(d)}} & R(d) }\!\!$
  is any cofibrant replacement \eqref{CofibrantReplacement}.
  \end{itemize}
\end{example}

\newpage

\noindent {\bf Homotopy fibers and homotopy pullback.}

\begin{defn}[Homotopy fiber]
  \label{HomotopyFibers}
  Let $\mathbf{C}$ be a model category (Def. \ref{ModelCategories}).

  \noindent
  {\bf (i)} For $\xymatrix@C=12pt{ A \ar[r]^-{p} & B }$
  a morphism in $\mathbf{C}$ with
  $B  \in \mathbf{C}_{\mathrm{fib}} \subset \mathbf{C}$
  a fibrant object (Notation \ref{FibrantAndCofibrantObjects}),
  and for $\xymatrix@C=12pt{\ast \ar[r]^-{b} & B}$ a
  morphism from the terminal object (a ``point in $B$''),
  the {\it homotopy fiber} of $p$ over $b$ is the
  image in the homotopy category \eqref{LocalizationOfAModelCategoryAtWeakEquivalenes}
  of the ordinary fiber over $b$,
  i.e. the pullback \eqref{CartesianSquare} along $b$ in $\mathbf{C}$,
  of any fibration $\widetilde p$ weakly equivalent to $p$:

  \vspace{-.3cm}
  \begin{equation}
    \label{ConstructionOfHomotopyFiber}
    \raisebox{20pt}{
    \xymatrix{
      \mathrm{hofib}_b(\rho)
      \ar[r]
      %\ar[d]
      \ar@{}[dr]|-{
        %\mbox{\tiny (pb)}
      }
      &
      A
      \ar[d]^-{ \rho }
      \\
      %\ast
      %\ar[r]_-{ b }
      &
      B
    }
    }
    \;\;\;\;\;
    :=
    \;\;\;\;\;
    \gamma_{\mathbf{C}}
    \left(\!\!\!\!
    \raisebox{20pt}{
    \xymatrix@C=6em{
      \mathrm{fib}_{b}(\widetilde p)
      \ar[r]
      \ar[d]
      \ar@{}[dr]|-{\mbox{\tiny (pb)}}
      &
      \widetilde A
      \ar[d]_-{
        \widetilde \rho
        }^-{
          \scalebox{.7}{$
            \in \, \mathrm{Fib}
          $}
          }
         \ar@{<-}[r]^-{ \in \, \mathrm{W} }
      &
      A
      \ar[dl]^-{ \rho }
      \\
      \ast
      \ar[r]_-{b}
      &
      B
    }
    }
    \!\!
    \right)
    \;\;\;
    \in
    \:
    \mathrm{Ho}
    (
      \mathbf{C}
    )\;.
  \end{equation}
  \vspace{-.2cm}

  \noindent
  This is well-defined in that
  $\mathrm{hofib}_b(p) \in \mathrm{Ho}(\mathbf{C})$
  depends on the choice of fibration replacement $\widetilde p$
  only up to isomorphism in the homotopy category.

  \noindent
  {\bf (ii)} Dually, homotopy co-fibers are homotopy
  fibers in the opposite model category (Def. \ref{OppositeModelStructure}).
\end{defn}

More generally:
\begin{defn}[Homotopy pullback]
  \label{HomotopyPullback}
  Given a model category $\mathbf{C}$ (Def. \ref{ModelCategories})
  and a pair of coincident morphisms
  \vspace{-2mm}
  $$
    \xymatrix@R=1em{
      &
      A
      \ar[d]^-{\rho}
      \\
      X
      \ar[r]_-{ \tau }
      & B
    }
  $$

  \vspace{-3mm}
  \noindent
  between fibrant objects,
  the \emph{homotopy pullback} of $\rho$ along $\tau$
  (or \emph{homotopy fiber product of $\rho$ with $\tau$})
  is the image of $\rho$,
  regarded as
  an object in the homotopy category (Def. \ref{HomotopyCategory})
  of the slice model category (Example \ref{SliceModelCategory})
  under the right derived functor
  (Prop. \ref{DerivedFunctors}) of the right base change functor
  along $\tau$ (Ex. \ref{BaseChangeQuillenAdjunction}):
    \vspace{-2mm}
  \begin{equation}
    \label{ActionOfRightDerivedPullback}
    \mathrm{Ho}
    \big(
      \mathbf{C}^{/B}
    \big)
    \;\;\ni\;\;
    \raisebox{20pt}{
    \xymatrix@R=1.5em{
      A
      \ar[d]^\rho
      \\
      B
    }
    }
    \;\;\;\;
    \overset{
      \raisebox{3pt}{
        \tiny
        \color{greenii}
        \bf
        \def\arraystretch{.9}
        \begin{tabular}{c}
          homotopy
          \\
          pullback
        \end{tabular}
      }
    }{\longmapsto}
    \;\;\;\;
    \raisebox{20pt}{
    \xymatrix@R=1.5em{
      \RightDerived\tau^\ast A
      \ar[d]^-{
        \RightDerived\tau^\ast \rho
      }
      \\
      X
    }
    }
    \;\;\;
    :=
    \;\;\;
    \mathrm{Ho}
    \big(
      \mathbf{C}^{/X}
    \big)
    \,,
  \end{equation}

    \vspace{-2mm}
  \noindent
  By \eqref{BaseChangeAdjunctionHomIsomorphism},
  the derived adjunction counit \eqref{DerivedAdjunctionCoUnit}
  on \eqref{ActionOfRightDerivedPullback}
  gives a commuting square in \eqref{LocalizationOfAModelCategoryAtWeakEquivalenes}
  the homotopy category of $\mathbf{C}$
  \vspace{-1mm}
  \begin{equation}
    \label{ConstructionOfHomotopyPullback}
    \raisebox{20pt}{
    \xymatrix@C=3em{
      \RightDerived\tau^\ast A
      \ar[r]
      \ar[d]_-{
        \RightDerived\tau^\ast \rho
      }
      \ar@{}[dr]|-{
        \mbox{\tiny (hpb)}
      }
      &
      A
      \ar[d]^-{ \rho }
      \\
      X
      \ar[r]_-{ \tau }
      &
      B
    }
    }
    \;\;\;\;\;
    :=
    \;\;\;\;\;
    \gamma_{\mathbf{C}}
    \left(\!\!\!\!
    \raisebox{20pt}{
    \xymatrix@C=6em{
      \tau^\ast \widetilde A
      \ar[r]
      \ar[d]
      \ar@{}[dr]|-{\mbox{\tiny (pb)}}
      &
      \widetilde A
      \ar[d]_-{
        \widetilde \rho
      }^-{
          \scalebox{.7}{$
            \in \, \mathrm{Fib}
          $}
          }
         \ar@{<-}[r]^-{ \in \, \mathrm{W} }
      &
      A
      \ar[dl]^-{ \rho }
      \\
      X
      \ar[r]_-{
        \tau
      }
      &
      B
    }
    }
    \!\!
    \right)
    \;\;\;
    \in
    \:
    \mathrm{Ho}
    (
      \mathbf{C}
    )\;.
  \end{equation}

  \vspace{-1mm}
\noindent  This square in the homotopy category, together with its
  pre-image pullback square in the model category, is the \emph{homotopy pullback square} of $\rho$ along $\tau$.
\end{defn}
\begin{example}[Homotopy fiber is homotopy pullback to the point]
  Homotopy fibers (Def. \ref{HomotopyFibers})
  are the homotopy pullbacks (Def. \ref{HomotopyPullback})
  to the terminal object, by \eqref{RightDerivedFunctorByFibrantReplacement}.
\end{example}

\begin{lemma}[Factorization lemma {\cite[p. 431]{Brown73}}]
  \label{factorizationLemma}
  Let $\mathbf{C}$ be a model category (Def. \ref{ModelCategories})
  and $\xymatrix@C=12pt{A \ar[r]^\rho & B} \;\in \mathbf{C}_{\mathrm{fib}}$
  a morphism between fibrant objects. Then for
  $\mathrm{Paths}(B)$ a path space object for $B$ (Def. \ref{PathSpaceObject})
  the vertical composite in the following diagram
  \begin{equation}
    \xymatrix@C=3em{
      A
      \ar[rr]^-{ \in \; \mathrm{W} }
      \ar[ddrr]_{ \rho }
      & &
      p_1^\ast A
      \ar[d]
      \ar[rr]^-{ \in \, \mathrm{W} }
      \ar@{}[drr]|-{\mbox{\tiny \rm(pb)}}
      %\ar@/_1.6pc/[dd]_<<<<<<<<<<<<<<<{ \in \; \mathrm{Fib} }
      &&
      A
      \ar[d]^-{\rho}
      \\
      &&
      \scalebox{.85}{$
      \mathrm{Paths}(B)
      $}\!\!
      \ar[rr]_-{ \,p_1\, }
      \ar[d]^-{
        \phantom{\mathclap{\vert}}
        p_0
      }
      &&
      B
      \\
      &&
      B
    }
  \end{equation}
  is a fibration, and in fact a fibration resolution of
  $\rho$, in that it factors $\rho$ through a weak equivalence.
\end{lemma}

\begin{example}[Homotopy pullback via triples]
  \label{HomotopyPullbackViaTriples}
  Given a model category $\mathbf{C}$ (Def. \ref{ModelCategories})
  and a pair of coincident morphisms
  \vspace{-2mm}
  $$
    \xymatrix@R=1em{
      &
      A
      \ar[d]^-{\rho}
      \\
      X
      \ar[r]^-{ \tau }
      & B
    }
  $$

  \vspace{-2mm}
  \noindent
  between fibrant objects,
  Lemma \ref{factorizationLemma} says that the
  corresponding homotopy pullback (Def. \ref{HomotopyPullback})
  is computed by the following diagram
  \vspace{-1mm}
  $$
    \raisebox{40pt}{
    \xymatrix@C=3em{
      \RightDerived\tau^\ast A
      \ar[dd]
      \ar[rr]
      \ar@{-->}[dr]|-{
        \;
        \mathclap{\phantom{\vert}}
        \phi
        \;
        }
      &&
      A
      \ar[d]^-{ \rho }
      \\
      &
      \scalebox{.85}{$
        \mathrm{Paths}(B)
      $}
      \ar[d]^-{p_1}
      \ar[r]_-{ p_0 }
      &
      B
      \\
      X
      \ar[r]^-{ \tau }
      &
      B
    }
    }
    \;\;\;\;\;\;\;
    =
    \;\;\;\;
    \raisebox{40pt}{
    \xymatrix@C=3em{
      \tau^\ast
      \big(
        p_0
        \,\circ\,
        p_1^\ast \rho
      \big)
      \ar[dd]
      \ar[r]
      \ar@{}[ddr]|-{\mbox{\tiny(pb)}}
      &
      p_1^\ast A
      \ar[r]
      \ar[d]
      \ar@{}[dr]|-{\mbox{\tiny(pb)}}
      &
      A
      \ar[d]^-{ \rho }
      \\
      &
      \scalebox{.85}{$
        \mathrm{Paths}(B)
      $}
      \ar[d]^-{p_1}
      \ar[r]_-{ p_0 }
      &
      B
      \\
      X
      \ar[r]^-{ \tau }
      &
      B
    }
    }
  $$

    \vspace{-1mm}
\noindent  Here the right hand side exhibits the
  left hand side as a limit cone. This means that the
  homotopy pullback $\RightDerived\tau^\ast A$ is
  universally characterized by the fact that
  morphisms into it are triples
  $(f,g,\phi)$ , consisting of
  a pair of morphisms $f$, $g$ to $A$, $X$, respectively,
  and a right homotopy $\phi$ (Def. \ref{RightHomotopy}) between
  their composites with $\rho$ and $\tau$, respectively:
    \vspace{-2mm}
  \begin{equation}
    \label{MapsIntoHomotopyPullbackAsTriplesOfData}
    \mathbf{C}
    \big(
      -;
      \,
      \RightDerived\tau^\ast A
    \big)
    \;\;\;\simeq\;\;\;
    \left\{
      (f,g,\phi)
      \,\left\vert\,
    \raisebox{18pt}{
    \xymatrix@C=5em{
      \ar@{-->}[r]^-f_>>>>>>{\ }="s"
      \ar@{-->}[d]_-{g}^->>>{\ }="t"
      &
      A
      \ar[d]^-{\rho}
      \\
      X
      \ar[r]^-{\tau}
      &
      B
      \ar@{==>}|{\;\phi\;} "s"; "t"
    }
    }
    \right.
    \right\}.
  \end{equation}
\end{example}

\noindent {\bf Quillen equivalences.}

\begin{lemma}[Conditions characterizing Quillen equivalences]
  \label{ConditionsCharacterizingQuillenEquivalences}
  Given a Quillen adjunction  $L \dashv_{\mathrm{Qu}} R$ (Def. \ref{QuillenAdjunction}),
  the following conditions are equivalent:
  \begin{itemize}
    \vspace{-.3cm}
    \item
      The left derived functor (Prop. \ref{DerivedFunctors}) is an equivalence
      of homotopy categories (Def. \ref{HomotopyCategory})
      $\!\!
        \xymatrix@C=1.5em{
          \mathrm{Ho}(\mathcal{D})
          \ar@{<-}[r]^-{ \LeftDerived L }_-{\simeq}
          &
          \mathrm{Ho}(\mathcal{C})
        \!\!}
      $.
    \vspace{-.3cm}
    \item
      The right derived functor (Prop. \ref{DerivedFunctors})
      is an equivalence
      of homotopy categories (Def. \ref{HomotopyCategory})
      $\!\!
        \xymatrix@C=1.5em{
          \mathrm{Ho}(\mathcal{D})
          \ar[r]^-{ \RightDerived R }_-{\simeq}
          &
          \mathrm{Ho}(\mathcal{C})
        \!\!}
      $.
    \vspace{-.2cm}
    \item Both of the following two conditions hold:
    \begin{enumerate}
    \vspace{-.2cm}
    \item[{\bf (i)}]
      The derived adjunction unit
      $\mathbb{D}\eta$ \eqref{DerivedAdjunctionUnit}
      is a
      natural isomorphism, hence
      \eqref{DerivedAdjunctionUnit} is a weak equivalence in
      $\mathbf{C}$;

    \vspace{-.1cm}
    \item[{\bf (ii)}]
      The derived adjunction counit
      $\mathbb{D}\epsilon$ \eqref{DerivedAdjunctionCoUnit}
      is a natural isomorphism, hence
      \eqref{DerivedAdjunctionCoUnit}
      is a weak equivalence in $\mathbf{D}$.
    \end{enumerate}
    \vspace{-.3cm}
    \item For $c \in \mathbf{C}^{\mathrm{cof}}$
    and $d \in \mathbf{D}_{\mathrm{fib}}$,
    a morphism out of $L(c)$ is a weak equivalence precisely
    if its adjunct into $R(d)$ is:

    \vspace{-.5cm}
    \begin{equation}
      \label{QuillenEquivalenceAdjunctPreservesWeakEquivalences}
      \xymatrix{
        L(c)
        \ar[rr]^-{f}_-{ \in \; \mathrm{W} }
        &&
        d
      }
      \;\;\;\;\;\;\;
      \Leftrightarrow
      \;\;\;\;\;\;\;
      \xymatrix{
        c
        \ar[rr]^-{\widetilde f}_-{ \in \; \mathrm{W} }
        &&
        R(d)\;.
      }
    \end{equation}
   \end{itemize}
\end{lemma}

\begin{defn}[Quillen equivalence]
  \label{QuillenEquivalence}
  If the equivalent conditions from Lemma \ref{ConditionsCharacterizingQuillenEquivalences} are met,
  a Quillen adjunction $L \dashv_{\mathrm{Qu}} R$ (Def. \ref{QuillenAdjunction})
  is called a \emph{Quillen equivalence}, which we
  denote as follows:
  \vspace{-3mm}
  $$
    \xymatrix{
      \mathcal{D}
      \ar@{<-}@<+7pt>[rr]^-{L}
      \ar@<-7pt>[rr]_-{R}^-{ \simeq_{\mathrlap{\mathrm{Qu}}} }
      &&
      \mathcal{C}
      \,.
    }
  $$
\end{defn}

\vspace{-3mm}
Hence:
\begin{prop}[Derived equivalence of homotopy categories]
  \label{DerivedEquivalenceOfHomotopyCategories}
  The derived adjunction (Prop. \ref{DerivedFunctors})
  of a Quillen equivalence
  (Def. \ref{QuillenEquivalence}) is an adjoint equivalence
  of homotopy categories (Def. \ref{HomotopyCategory}):
    \vspace{-2mm}
  \begin{equation}
    \label{DerivedEquivalence}
    \xymatrix{
      \mathrm{Ho}
      (
        \mathbf{D}
      )
        \ar@{<-}@<+7pt>[rr]^-{ \LeftDerived L }
        \ar@<-7pt>[rr]_-{ \RightDerived R }^-{
          \scalebox{.8}{\raisebox{1pt}{$\simeq$}}
        }
      &&
      \;
      \mathrm{Ho}
      (
        \mathbf{C}
      )\;.
    }
  \end{equation}
\end{prop}

\begin{remark}[$\infty$-Category theory]
  \label{InfinityCategoryTheory}
  As each model category (Def. \ref{ModelCategories})
  provides a context of {\it homotopy theory}
  (with its own notion of homotopy-coherent universal constructions such as
  homotopy pullbacks, Def. \ref{HomotopyPullback}, etc.),
  Prop. \ref{DerivedEquivalenceOfHomotopyCategories}
  is a first indication
  that Quillen equivalent (Def. \ref{QuillenEquivalence})
  model categories represent the {\it same} context of homotopy theory,
  for a suitably homotopy-theoretic notion of sameness.
  This suggests that model categories regarded up to Quillen equivalence
  are but
  coordinate presentations of a more
  intrinsic notion of homotopy theories, now known as
  {\it $\infty$-categories}
  \cite{Joyal}\cite{Joyal08}\cite{Joyal08b}\cite{Lurie09}\cite{Cisinski19}\cite{RiehlVerity21}.
\end{remark}

\begin{lemma}[Quillen equivalence when left adjoint creates weak equivalences {\cite[Lemma 3.3]{ErdalIlhan19}}]
  \label{QuillenEquivalenceWhenLeftAdjointCreatesWeakEquivalences}
  Let $L \dashv_{\mathrm{Qu}} R$
  %$$
  %  \xymatrix{
  %    \mathbf{D}
  %    \ar@{<-}@<+7pt>[rr]^-{L}
  %    \ar@<-7pt>[rr]_-{R}^-{\bot_{\mathrlap{\mathrm{Qu}}}}
  %    &&
  %    \mathbf{C}
  %  }
  %$$
  be a Quillen adjunction
  (Def. \ref{QuillenAdjunction}) such that the left
  adjoint functor $L$ creates weak equivalences, in that
  for all morphisms $f$ in $\mathbf{C}$ we have
\vspace{-1mm}
  \begin{equation}
    \label{LCreatedWeakEquivalences}
    f \;\in\; \mathrm{W}_{\mathbf{C}}
    \;\;\;\;\;\;\;\;\;
    \Leftrightarrow
    \;\;\;\;\;\;\;\;\;
    L(f) \;\in\; \mathrm{W}_{\mathbf{D}}
    \,.
  \end{equation}

  \vspace{-1mm}
\noindent  Then $L \dashv_{\mathrm{Qu}} R$ is a Quillen equivalence
  (Def. \ref{QuillenEquivalence}) precisely if the
  adjunction co-unit $\epsilon_d$ is a weak
  equivalence on all fibrant objects $d \in \mathbf{C}_{\mathrm{fib}}$.
\end{lemma}

\begin{proof}
  By Lemma \ref{ConditionsCharacterizingQuillenEquivalences},
  it is sufficient to check that the  {\bf (i)} derived unit
  and {\bf (ii)} derived counit
  of the adjunction are weak equivalences precisely if the
  ordinary counit is a weak equivalence.

  \noindent
  {\bf (ii)} For the derived counit \eqref{DerivedAdjunctionCoUnit}
   \vspace{-2mm}
\noindent
  $$
    \mathbb{D}\epsilon_c
    \;:\;
        \xymatrix{
          L\big( Q R(d)\big)
          \ar[rr]^-{
            L
            (
                p_{ R(d) }
            )
          }_-{ \in \; \mathrm{W} }
          &&
          L\big(R(d)\big)
          \ar[r]^-{ \epsilon_d }_-{  }
          &
          d
        }
  $$

   \vspace{-2mm}
\noindent
  we have that $p_{ R(d) }$ is a weak equivalence
  \eqref{CofibrantReplacement}, and since $L$ preserves
  this, by assumption, so is $L\big(p_{ R(d) }\big)$.
  Therefore $\mathbb{D}\epsilon_d$ is a weak equivalence
  precisely if $\epsilon_d$ is, by 2-out-of-3 \eqref{TwoOutOfThree}.

  \noindent
  {\bf (i)} For the derived unit \eqref{DerivedAdjunctionUnit}
   \vspace{-3mm}
  $$
        \xymatrix{
          c
          \ar[r]^-{ \eta_c }
          &
          R\big(L(c)\big)
          \ar[rr]^-{
            R
            (
              j_{ L(c) }
            )
          }
          &&
          R\big( P L(c)\big)
        }
  $$

   \vspace{-1mm}
\noindent
  consider the composite of its image under $L$
  with the adjunction counit, as shown in the middle row of
  the following diagram:
   \vspace{-3mm}
  $$
    \xymatrix@C=4em{
      L(c)
      \ar[rr]|-{\;
        L
        (
          \eta_c
        )
      \;}
      \ar@/_1.2pc/[rrrrr]_-{
        j_{L(c)}
      \in \; \mathrm{W}}
      \ar@/^1.2pc/[rrrr]^-{
        L
        (
          \mathbb{D}\eta_c
        )
      }
      &&
      L \circ R\big(L(c)\big)
      \ar[rr]|-{\;
        L \circ R
        (
          j_{ L(c) }
        )
      \;}
      &&
      L \circ R
      \big(
        P L(c)
      \big)
      \ar[r]^-{
        \epsilon_{P L(c)}
      }
      &
      P L(c)
    }.
  $$

   \vspace{-2mm}
\noindent
  By the formula \eqref{AdjunctsByCompositesWithCoUnits}
  for adjuncts, this composite equals the adjunct
  of the derived adjunction unit, hence $j_{L(c)}$,
  as shown by the bottom morphism, which is a weak equivalence
  \eqref{FibrantReplacement}.
  Now, since $L$ creates weak equivalences by assumption,
  $L(\mathbb{D}\eta_c)$ is a weak equivalence
  precisely if $\mathbb{D}\eta_c$ is a weak equivalence.
  Therefore it follows, again by 2-out-of-3 \eqref{TwoOutOfThree},
  that this is the case precisely if the adjunction counit
  $\epsilon$ is a weak equivalence on the fibrant object $P L(c)$.
\end{proof}

\begin{prop}[Base change along weak equivalence in right proper model category]
  \label{BaseChangeAlongWeakEquivalenceInRightProperModelCategory}
  Let $\mathbf{C}$ be a right proper model category
  (Def. \ref{ProperModelCategories}).
  Then its base change Quillen adjunction
  (Ex. \ref{BaseChangeQuillenAdjunction})
  along any weak equivalence
   \vspace{-2mm}
  $$
    \xymatrix{
      B_1
      \ar[r]^-{ f }_-{ \in \; \mathrm{W} }
      &
      B_2
    }
    \;\;\;\;\;
    \in
    \;
    \mathbf{C}
  $$

   \vspace{-2mm}
\noindent
  is a Quillen equivalence (Def. \ref{QuillenEquivalence}):
   \vspace{-2mm}
  $$
    \xymatrix{
      \mathbf{C}^{/B_2}
      \ar@<+7pt>@{<-}[rr]^-{ f_! }
      \ar@<-7pt>[rr]_-{ f^\ast }^-{ \simeq_{\mathrlap{\mathrm{Qu}}} }
      &&
     \; \mathbf{C}^{/B_1}
    }.
  $$
\end{prop}

\vspace{-7mm}
\begin{proof}
  Observe that
  $\!\!\xymatrix@C=12pt{ B_2 \ar[r]^{\mathrm{id}}  & B_2 }\!\!$
  is the terminal object of $\mathbf{C}^{/B_2}$, so that
  the fibrant objects of $\mathbf{C}^{/B_2}$ correspond
  to the fibrations in $\mathbf{C}$ over $B_2$.
  Therefore, the condition \eqref{QuillenEquivalenceAdjunctPreservesWeakEquivalences}
  says that for $f_! \dashv f^\ast$ to be a Quillen equivalence
  it is sufficient that in
  \eqref{BaseChangeAdjunctionHomIsomorphism}
  $c$ is a weak equivalence precisely if $\widetilde c$
  is, assuming that $\rho$ is a fibration:
   \vspace{-2mm}
  \begin{equation}
    \label{BaseChangeAdjunctionHomIsomorphismForProofOfEquivalence}
    \begin{array}{ccc}
      \scalebox{.9}{$
      \raisebox{30pt}{
      \xymatrix@R=1pt@C=4em{
        X
        \ar@{-->}[rr]^-{ c }
        \ar[ddddr]_-{ \tau }
        &
        &
        A
        \ar[ddddd]^-{
          \rho
          \scalebox{.6}{$\in \mathrm{Fib}$}
        }
        \\
        \\
        \\
        \\
        &
        B_1
        \ar[dr]_-{
             f  \in \; \mathrm{W}
        }
        \\
        &&
        B_2
      }
      }
      $}
      &
      \;\;\;\;\;
      \leftrightarrow
      \;\;\;\;\;
      &
      \scalebox{.9}{$
      \raisebox{30pt}{
      \xymatrix@R=1pt@C=5em{
        X
        \ar@{-->}[r]^-{ \widetilde {c} }
        \ar[ddddr]_-{ \tau }
        &
        f^\ast A
        \ar[dddd]^-{
          \mathclap{\phantom{\vert^{\vert}}}
          f^\ast \rho
          \mathclap{\phantom{\vert_{\vert}}}
        }
        \ar[r]^-{
              \rho^\ast f  \in \; \mathrm{W}
         }
        \ar@{}[dddddr]|-{ \mbox{\tiny(pb)} }
        &
        A
        \ar[ddddd]^-{
          \rho
          \scalebox{.6}{$\in \mathrm{Fib}$}
        }
        \\
        \\
        \\
        \\
        &
        B_1
        \ar[dr]_-{
          f
            \in \; \mathrm{W}
               }
        \\
        &&
        B_2
      }
      }
      $}
    \end{array}
  \end{equation}
  But under this assumption, right-properness
  implies that $\rho^\ast f$ is a weak equivalence
  \eqref{RightProperness},
  so that the statement follows by 2-out-of-3
  \eqref{TwoOutOfThree}.

  \medskip

  \noindent
  {\it Alternative Proof.} The conclusion also follows with Lemma \ref{QuillenEquivalenceWhenLeftAdjointCreatesWeakEquivalences}:
  The left adjoint functor $L = f_!$ clearly creates weak equivalences
  \eqref{LCreatedWeakEquivalences}
  (by the nature of the slice model structure, Example \ref{SliceModelCategory}),
  so that Lemma
  \ref{QuillenEquivalenceWhenLeftAdjointCreatesWeakEquivalences}
  asserts that we have a Quillen equivalence as soon as the
  ordinary adjunction counit is a weak equivalence
  on all fibrant objects.
  By \eqref{BaseChangeAdjunctionHomIsomorphism},
  the adjunction counit on a fibration $\rho \in \mathrm{Fib}$ is
  the dashed morphism $\rho^\ast f$ in the following diagram on the right:
  \vspace{-2mm}
  \begin{equation}
    \label{BaseChangeAdjunctionHomIsomorphismForProofOfEquivalenceAlternative}
    \begin{array}{ccc}
      \scalebox{.9}{$
      \raisebox{30pt}{
      \xymatrix@R=1pt@C=5em{
        f^\ast A
        \ar@{-->}[r]^-{ \mathrm{id} }
        \ar[ddddr]_-{ f^\ast \rho }
        &
        f^\ast A
        \ar[dddd]^-{
          \mathclap{\phantom{\vert^{\vert}}}
          f^\ast \rho
          \mathclap{\phantom{\vert_{\vert}}}
        }
        \ar[r]^-{
           \rho^\ast f
            \in \; \mathrm{W}
                }
        \ar@{}[dddddr]|-{ \mbox{\tiny(pb)} }
        &
        A
        \ar[ddddd]^-{
          \rho
          \scalebox{.6}{$\in \mathrm{Fib}$}
        }
        \\
        \\
        \\
        \\
        &
        B_1
        \ar[dr]_-{ f
               \in \; \mathrm{W}
           }
        \\
        &&
        B_2
      }
      }
      $}
      &
      \;\;\;\;\;
      \leftrightarrow
      \;\;\;\;\;
      &
      \scalebox{.9}{$
      \raisebox{30pt}{
      \xymatrix@R=1pt@C=4em{
        f^\ast A
        \ar@{-->}[rr]^-{
        \rho^\ast f
            \in \; \mathrm{W}
        }
        \ar[ddddr]_-{ f^\ast \rho }
        &
        &
        A
        \ar[ddddd]^-{
          \rho
          \scalebox{.6}{$\in \mathrm{Fib}$}
        }
        \\
        \\
        \\
        \\
        &
        B_1
        \ar[dr]_-{
          f
            \in \; \mathrm{W}
        }
        \\
        &&
        B_2
      }
      }
      $}
    \end{array}
  \end{equation}

  \vspace{-1mm}
\noindent  And hence this is a weak equivalence, again by right-properness.
\end{proof}

\begin{example}[Quillen equivalence between topological spaces and simplicial sets {\cite{Quillen67}}]
  \label{QuillenEquivalenceBetweenTopologicalSpacesAndSimplicialSets}
 Forming simplicial sets constitutes a
 Quillen equivalence (Def. \ref{QuillenEquivalence})
 \vspace{-2mm}
  \begin{equation}
    \label{QuillenAdjunctionBetweenTopologicalSpacesAndSimplicialSets}
    \xymatrix{
      \TopologicalSpaces_{\mathrm{Qu}}
      \;
      \ar@{<-}@<+7pt>[rr]^-{
        \overset{
          \raisebox{3pt}{
            \tiny
            \color{greenii}
            \bf
            geometric realization
          }
        }{
          \left\vert - \right\vert
        }
      }
      \ar@<-7pt>[rr]_-{
        \underset{
          \raisebox{-3pt}{
            \tiny
            \color{greenii}
            \bf
            singular simplicial complex
          }
        }{
          \mathrm{Sing}
        }
      }^-{
          \simeq_{\mathrlap{\mathrm{Qu}}}
      }
      &&
      \;\;
      \SimplicialSets_{\mathrm{Qu}}
    }
  \end{equation}

  \vspace{-2mm}
  \noindent  between the classical model structure on topological spaces
  (Example \ref{ClassicalModelStructureOnTopologicalSpaces})
  and the classical model structure on simplicial sets
  (Example \ref{ClassicalModelStructureOnSimplicialSets}).
\end{example}

\begin{example}[Classical homotopy category]
  \label{TheClassicalHomotopyCategory}
 % By Prop. \ref{DerivedEquivalenceOfHomotopyCategories} and
  The derived adjunction (Prop. \ref{DerivedFunctors})
  of the $\left\vert-\right\vert \dashv \mathrm{Sing}$-adjunction
  (Example \ref{QuillenEquivalenceBetweenTopologicalSpacesAndSimplicialSets})
  is an equivalence
  between the homotopy categories (Def. \ref{HomotopyCategory})
  of the classical model category of topological spaces
  (Example \ref{ClassicalModelStructureOnTopologicalSpaces})
  and the classical model category of simplicial sets
  (Example \ref{ClassicalModelStructureOnSimplicialSets}):
  \vspace{-2mm}
  \begin{equation}
    \label{ClassicalHomotopyCategory}
    \xymatrix{
      \mathrm{Ho}
      \big(
        \TopologicalSpaces_{\mathrm{Qu}}
      \big)
      \;
      \ar@{<-}@<+7pt>[rr]^-{
        \LeftDerived\left\vert - \right\vert
      }
      \ar@<-7pt>[rr]_-{
        \RightDerived\mathrm{Sing}
      }^-{
        \scalebox{.7}{\raisebox{1pt}{$\simeq$}}
      }
      &&
      \;\;
      \mathrm{Ho}
      \big(
        \SimplicialSets_{\mathrm{Qu}}
      \big).
    }
  \end{equation}

  \vspace{-2mm}
  \noindent  Either of these is the {\it classical homotopy category}.
  We refer to its objects as \emph{homotopy types},
  to be distinguished from the actual topological spaces or simplicial
  sets that represent them.
\end{example}

\begin{example}[Simplicial sets are weakly equivalent to singular simplicial sets of their realization]
  \label{SimplicialSetsEquivalentToSingularSimplicialSetsOfTheirRealization}
  The characterization of Quillen equivalences (Lemma \ref{ConditionsCharacterizingQuillenEquivalences})
  implies, with Example \ref{QuillenEquivalenceBetweenTopologicalSpacesAndSimplicialSets},
  that for each $S \,\in\, \SimplicialSets$ the composite
  \vspace{-3mm}
  $$
    \xymatrix{
      S
      \ar[rr]^-{ \eta_S }
      &&
      \mathrm{Sing}
      (
        \vert
          S
        \vert
      )
      \ar[rr]^-{\scalebox{0.7}{$
        \mathrm{Sing}
        (
          \vert
            j_{
              \left\vert
                S
              \right\vert
            }
          \vert
        )
        $}
      }
      &&
      \mathrm{Sing}
      (
        P
        \left\vert
          S
        \right\vert
      )
    }
  $$

  \vspace{-2mm}
  \noindent is a weak equivalence, where
  $j_{\left\vert S\right\vert}$ is a fibrant replacement
  for $\left\vert S\right\vert$.
  But since all topological spaces are fibrant
  (Example \ref{ClassicalModelStructureOnTopologicalSpaces}),
  it follows that the ordinary unit of the
  adjunction \eqref{QuillenAdjunctionBetweenTopologicalSpacesAndSimplicialSets}
  is already a weak equivalence:
  \vspace{-3mm}
  \begin{equation}
    \label{RealizationSingAdjunctionUnitIsWeakEquivalence}
    \xymatrix{
      S
      \ar[rr]^-{ \eta_S }_-{ \in\,\mathrm{W} }
      &&
      \mathrm{Sing}
      (
        \vert
          S
        \vert
      )
    }.
  \end{equation}
\end{example}

\medskip

\noindent {\bf Cell complexes.}

\begin{prop}[Skeleta and truncation {\cite[\S II.8]{May67}\cite[\S 1.2 (vi)]{DwyerKan84}} ]
  \label{SkeletaAndPostnikovTowers}
  For each $n \in \mathbb{N}$ there is a pair of adjoint functors
    \vspace{-2mm}
    \begin{equation}
    \label{SkeletonAdjunction}
    \xymatrix{
      \SimplicialSets
      \ar@{<-}@<+6pt>[rr]^-{ \mathrm{sk} }
      \ar@<-6pt>[rr]_-{ \mathrm{cosk} }^-{ \bot }
      &&
      \SimplicialSets
      \,,
    }
  \end{equation}

  \vspace{-2mm}
  \noindent
  where $\mathrm{sk}_n(S)$ is the simplicial sub-set generated
  by the simplices in $S$ of dimension $\leq n$ (hence including
  only all their degenerate higher simplices), and where
  \vspace{-2mm}
  $$
    \mathrm{cosk}_n(S)
      \;:\;
    [k]
      \;\mapsto\;
    \SimplicialSets
    \big(
      \mathrm{sk}_n(\Delta[k])
      \,,\,
      S
    \big)
    \,.
  $$
  One says that $S$ is {\it $n$-coskeletal} if the comparison morphism
  $S \xrightarrow{\;} \mathrm{cosk}_n(S)$ is an isomorphism.

  \noindent
  Here $\mathrm{cosk}_{n+1}$ preserves
  {\rm (\cite[p. 141]{DwyerKan84}, for proofs see \cite{Low13}\cite[Lem. 10.12 ]{Deflorin19})}
  fibrant objects
  of the classical model structure (Example \ref{ClassicalModelStructureOnSimplicialSets}),
  hence preserves Kan complexes \eqref{KanComplexes},
  and models \emph{$n$-truncation}, in that:
  \vspace{-2mm}
  $$
    \pi_k
    \big\vert
      \mathrm{cosk}_{n+1}(S)
    \big\vert
    \;=\;
    0
    \phantom{AAAA}
    \mbox{for $k \geq n+1$}
  $$

  \vspace{-1mm}
  \noindent
  and there are natural morphisms

  \vspace{-.3cm}
  \begin{equation}
   \label{nTruncationMorphismsViaCoskeletonUnits}
   \xymatrix{\!
      S \ar[r]^-{ p_n } & \mathrm{cosk}_n(S)
    \!}
  \end{equation}
  \vspace{-.4cm}

  \noindent
  such that
   \vspace{-3mm}
   $$
    \xymatrix{
      \pi_k
      \big\vert
        S
      \big\vert
      \ar[rr]^-{
        \pi_k \left\vert p_n \right\vert
      }_-{ \simeq }
      &&
      \pi_k
      \big\vert
        \mathrm{cosk}_{n+1}(S)
      \big\vert
    }
    \phantom{AAA}
    \mbox{for $k \leq n$ }.
  $$
\end{prop}
\noindent For $A \in \HomotopyTypes$ we write
 \vspace{-2mm}
  \begin{equation}
  \label{nTruncationOnHomotopyTypes}
  A(n)
  \;:=\;
  \left\vert
  \mathrm{cosk}_{n+1}
  \big(
    \mathrm{Sing}(A)
  \big)
  \right\vert
\end{equation}
We say that $A$ is \emph{$n$-truncated}
if it is
equivalent to its $n$-truncation \eqref{nTruncationOnHomotopyTypes}:
\begin{equation}
  \label{nTruncated}
  \mbox{
    $A$ is $n$-truncated
  }
  \;\;\;\;
  \Leftrightarrow
  \;\;\;\;
  A \;\simeq\; A(n)
  \,.
\end{equation}

\begin{example}[Homotopy types of manifolds via triangulations]
  \label{HomotopyTypesOfManifoldsViaTriangulations}
  For $X \in \TopologicalSpaces$ equipped with the
  structure of a smooth $n$-manifold,
  there exists a \emph{triangulation} of $X$
  (e.g. \cite[\S IV.B]{Whitney57}\cite[Thm. 10.6]{Munkres66},
  see also \cite{Manolescu14}),
  namely a
  simplicial set (in fact a simplicial complex)
  which is $n$-sleletal
  (Prop. \ref{SkeletaAndPostnikovTowers})
  \begin{equation}
    \label{ATriangulation}
    \mathrm{Tr}(X) \in \SimplicialSets
    \,,
    \;\;\;\;\;\;
    \mathrm{sk}_n\big(\mathrm{Tr}(X)\big)
      \;=\;
    \mathrm{Tr}(X)
  \end{equation}
  equipped with a
  homeomorphism to $X$
  out of its geometric realization
  \eqref{QuillenAdjunctionBetweenTopologicalSpacesAndSimplicialSets}

  \vspace{-.4cm}
  \begin{equation}
    \label{HomeomorphismFromRealizationOfTriangulation}
    \xymatrix@C=3em{
      \left\vert\mathrm{Tr}(X)\right\vert
      \ar[r]^-{p}_-{ \mathrm{homeo} }
      &
      X
    }
  \end{equation}
  \vspace{-.4cm}

  \noindent
  which restricts in the interior of each simplex to a diffeomorphism onto its image.
  Since the inclusion
  \vspace{-2mm}
  \begin{equation}
    \label{InclusionOfTriangulationIntoSingularSimplicialSet}
   \xymatrix@C=4em{
     \mathrm{Tr}(X)
     \;
     \ar@{^{(}->}[r]^-{ \eta_{\mathrm{Tr}(X)} }_-{ \in\, \mathrm{W} }
     &
     \mathrm{Sing}
     \big(
       \left\vert
         \mathrm{Tr}(X)
       \right\vert
     \big)
     \ar[r]^-{
       \mathrm{Sing}(p)
     }_-{ \in\, \mathrm{Iso} }
     &
     \mathrm{Sing}(X)
     \,,
   }
  \end{equation}

  \vspace{-1mm}
  \noindent  is a weak equivalence (by Example \ref{SimplicialSetsEquivalentToSingularSimplicialSetsOfTheirRealization}),
  the triangulation represents the homotopy type \eqref{ClassicalHomotopyCategory}
  of the manifold.
\end{example}

\begin{prop}[Homotopy classes of maps out of $n$-manifolds]
  \label{HomotopyClassesOfMapsOutOfnManifolds}
  Let $X \in \TopologicalSpaces$ admit the
  structure of an $n$-manifold. Then for any
  $A \in \HomotopyTypes$
  (Example \ref{TheClassicalHomotopyCategory})
  the homotopy classes of maps
  $\xymatrix@C=11pt{X \ar[r] & A}$ are in natural bijection
  to the homotopy classes into the $n$-truncation \eqref{nTruncationOnHomotopyTypes}
  of $A$:
  \begin{equation}
    \label{HomotopyClassesOfMapsOutofNManifoldsEquivalentToMapsIntonminus1Trucation}
    \HomotopyTypes
    \big(
      X
      \,,\,
      A
    \big)
    \;\simeq\;
    \HomotopyTypes
    \big(
      X
      \,,\,
      A(n)
    \big)
  \end{equation}
\end{prop}
\begin{proof}
  Consider the following sequence of natural isomorphisms
  $$
    \begin{aligned}
      &
      \HomotopyTypes
      \big(
        X
        \,,\,
        A
      \big)
      \\
      &
      \simeq
      \mathrm{Ho}
      \big(
        \SimplicialSets_{\mathrm{Qu}}
      \big)
      \big(
        \mathrm{Sing}(X)
        \,,\,
        \mathrm{Sing}(A)
      \big)
      \\
      & \simeq
      \mathrm{Ho}
      \big(
        \SimplicialSets_{\mathrm{Qu}}
      \big)
      \big(
        \mathrm{Tr}(X)
        \,,\,
        \mathrm{Sing}(A)
      \big)
      \\
      & \simeq
      \SimplicialSets
      \Big(
        \mathrm{Tr}(X)
        \,,\,
        \mathrm{Sing}(A)
      \Big)
      \Big/
      \SimplicialSets
      \Big(
        \mathrm{Tr}(X) \times \Delta[1]
        \,,\,
        \mathrm{Sing}(A)
      \Big)
      \\
      & \simeq
      \SimplicialSets
      \Big(
        \mathrm{sk}_{n+1}
        \big(
          \mathrm{Tr}(X)
        \big)
        \,,\,
        \mathrm{Sing}(A)
      \Big)
      \Big/
      \SimplicialSets
      \Big(
        \mathrm{sk}_{n+1}
        \big(
          \mathrm{Tr}(X) \times \Delta[1]
        \big)
        \,,\,
        \mathrm{Sing}(A)
      \Big)
      \\
      & \simeq
      \SimplicialSets
      \Big(
        \mathrm{Tr}(X)
        \,,\,
        \mathrm{cosk}_{n+1}
        \big(
          \mathrm{Sing}(A)
        \big)
      \Big)
      \Big/
      \SimplicialSets
      \Big(
        \mathrm{Tr}(X) \times \Delta[1]
        \,,\,
        \mathrm{cosk}_{n+1}
        \big(
          \mathrm{Sing}(A)
        \big)
      \Big)
      \\
      & \simeq
      \mathrm{Ho}
      \big(
        \SimplicialSets_{\mathrm{Qu}}
      \big)
      \Big(
        \mathrm{Tr}(X)
        \,,\,
        \mathrm{cosk}_{n+1}\big(\mathrm{Sing}(A)\big)
      \Big)
      \\
      & \simeq
      \HomotopyTypes
      \Big(
        \big\vert\mathrm{Tr}(X)\big\vert
        \,,\,
        \big\vert\mathrm{cosk}_{n+1}\big(\mathrm{Sing}(A)\big)\big\vert
      \Big)
      \\
      & \simeq
      \HomotopyTypes
      \big(
        X
        \,,\,
        A(n)
      \big).
    \end{aligned}
  $$
  Here
  the first step is \eqref{TheClassicalHomotopyCategory},
  using, with Example \ref{DerivedFunctorsByCoFibrantReplacement},
  that all topological spaces are fibrant and all simplicial sets cofibrant.
  The second step uses
  \eqref{InclusionOfTriangulationIntoSingularSimplicialSet}.
  The third step
  uses Example \ref{StandardCylinderObjectInSimplicialSets}
  with Prop. \ref{RightHomotopyReflectsHomotopyClasses}
  (observing that $\mathrm{Sing}(A)$ is fibrant as $A$ is
  and $\mathrm{Sing}$ is right Quillen)
  to express the morphisms in the homotopy category
  as equivalence classes of simplicial maps under the relation
  that identifies those pairs of maps that extend to a map
  on the cylinder $\mathrm{Tr}(X) \times \Delta[1]$.
  The fourth step observes that with
  $\mathrm{Tr}(X)$ being $n$-skeletal \eqref{ATriangulation},
  its cylinder is $(n+1)$-skeletal.
  The fifth step is thus the
  $\mathrm{sk}_{n+1} \dashv \mathrm{cosk}_{n+1}$-adjunction
  isomorphism \eqref{SkeletonAdjunction}.
  The sixth step applies again
  Prop. \ref{RightHomotopyReflectsHomotopyClasses},
  using that $\mathrm{cosk}_{n+1}$
  preserves fibrancy (Prop. \ref{SkeletaAndPostnikovTowers}).
  The seventh step is the reverse of the first step, with the
  same argument on (co-)fibrancy.
  The last step uses \eqref{HomeomorphismFromRealizationOfTriangulation}
  in the first argument and \eqref{nTruncationOnHomotopyTypes}
  in the second.
  The composite of these isomorphisms is the desired
  \eqref{HomotopyClassesOfMapsOutofNManifoldsEquivalentToMapsIntonminus1Trucation}.
\end{proof}

\begin{prop}[Postnikov tower {\cite[Cor. 3.7]{GoerssJardine99}}]
  \label{PostnikovTower}
  Let $X \in \HomotopyTypes$ (Example \ref{TheClassicalHomotopyCategory}).
  If $X$ is connected, then
  its sequence of $n$-truncations \eqref{nTruncationOnHomotopyTypes} forms
  a system of maps with homotopy fibers (Def. \ref{HomotopyFibers})
  the Eilenberg-MacLane spaces
  \eqref{EilenbergMacLaneSpaces}
  of the homotopy group in the given degree:

  \vspace{-6mm}
  $$
    \xymatrix@R=1.2em@C=4em{
      &
      \vdots
      \ar[d]
      \\
      K( \pi_3(X), 3 )
      \ar[r]^-{ \mathrm{hfib}(p^X_3) }
      &
      X(3)
      \ar[d]^-{ p^X_3 }
      \\
      K( \pi_2(X), 2 )
      \ar[r]^-{ \mathrm{hfib}(p^X_2) }
      &
      X(2)
      \ar[d]^-{ p^X_2 }
      \\
      K( \pi_1(X), 1 )
      \ar[r]^-{ \mathrm{hfib}(p^X_1) }
      &
      X(1)
      \ar[d]^-{ p^X_1 }
      \\
      &
      X(0)
      \,.
    }
  $$

  \vspace{-1mm}
\noindent   If $X$ is not connected then this applies to each of its
  connected components.
\end{prop}

\newpage

\noindent {\bf Stable model categories.}

\vspace{-1mm}
\begin{example}[Looping/suspension-adjunction]
  \label{LoopingSuspensionAdjunction}
  On the category of pointed topological spaces,
  equipped with the coslice model structure under the point
  (Example \ref{SliceModelCategory}) of the classical model structure
  (Example \ref{ClassicalModelStructureOnTopologicalSpaces}),
  the operation of forming based loop spaces
  $\Omega X \,:=\, \mathrm{Maps}^{\ast/}(S^1, X)$ is the right
  adjoint in a Quillen adjunction (Def. \ref{QuillenAdjunction})
  \vspace{-2mm}
  \begin{equation}
    \label{BasedLoopingSuspensionAdjunction}
    \xymatrix{
      \TopologicalSpaces^{\ast/}_{\mathrm{Qu}}
      \;
      \ar@{<-}@<+7pt>[rr]^-{ \Sigma }
      \ar@<-7pt>[rr]_-{ \Omega }^-{ \bot_{\mathrlap{\mathrm{Qu}}} }
      &&
      \;
      \TopologicalSpaces^{\ast/}_{\mathrm{Qu}}
    }
  \end{equation}

  \vspace{-2mm}
\noindent  whose left adjoint is the reduced suspension operation
  $
    \Sigma X
      \,:=\,
    S^1 \wedge X
      \,:=\,
    \big(
      S^1 \times X
    \big) /
    \big(
      S^1 \times \{\ast_{{}_X}\}
      \;\sqcup\;
      \{\ast_{{}_{S^1}}\} \times X
    \big)
  $.
\end{example}

\begin{example}[Stable model category of sequential spectra
{\cite{BousfieldFriedlander78}\cite[\S X.4]{GoerssJardine99}}]
  \label{ModelStructureOnSequentialSpectra}
  There exists a model category (Def. \ref{ModelCategories})
  $\mathrm{SequentialSpectra}_{\mathrm{BF}}$
  whose objects are sequences
    \vspace{-2mm}
  $$
    E
    \;:=\;
    \Big\{
      E_n
      \;\in\;
      \TopologicalSpaces
      ,
      \,
      \xymatrix@C=12pt{
        \Sigma E_n
        \ar[r]^-{ \sigma_n }
        &
        E_{n+1}
      }
    \big\}_{n \in \mathbb{N}}
  $$

    \vspace{-2mm}
\noindent
  of topological spaces $E_n$ and continuous function $\sigma_n$
  from their suspension $\Sigma E_n$ (Example \ref{LoopingSuspensionAdjunction})
  to the next space in the sequences; and
  whose morphisms $\xymatrix@C=11pt{E \ar[r]^f & F}$ are sequences of
  component maps $\xymatrix@C=12pt{E_n \ar[r]^{f_n} &  F_n}$ that
  commute with the $\sigma$s.
  Moreover:
  \begin{itemize}
  \vspace{-.2cm}
  \item[$\mathrm{W}$ --]  {\it weak equivalences} are the morphisms
    that induce isomorphisms on all \emph{stable homotopy groups}
    $\pi_\bullet(X) := \underset{\underset{n}{\longrightarrow}}{\mathrm{lim}}\pi_{\bullet + k}(X_k)$
    (where the colimit is formed using the $\sigma$'s);

  \vspace{-3mm}
  \item[$\mathrm{Cof}$ --]  {\it cofibrations} are those morphisms
  $\xymatrix@C=11pt{E \ar[r]^f & F}$ such that the maps
    \vspace{-3mm}
  $$
    \xymatrix{
      E_0 \ar[r]^-{f_0}_-{ \in \; \mathrm{Cof} } &  F_0
    }
    \;\;\;\;\;
    \mbox{and}
    \;\;\;\;\;
    \underset{n \in \mathbb{N}}{\forall}
    \;\;
    \xymatrix{
      E_{n+1} \underset{\Sigma E_n}{\sqcup}
      \Sigma F_n
      \;
      \ar[rr]^-{ (f_{n+1}, \sigma^F_{n}) }_-{ \in \; \mathrm{Cof} }
      &&
      \;
      F_{n+1}
    }
  $$

    \vspace{-4mm}
\noindent
  are cofibrations in the classical model structure on
  topological spaces (Example \ref{ClassicalModelStructureOnTopologicalSpaces}).

  \vspace{-.2cm}
  \item[$\mathrm{Fib}$ --] {\it Fibrant objects} are the
  \emph{$\Omega$-spectra}, namely those sequences of spaces
  $\{E_n\}$
  for which the $\Sigma \dashv \Omega$-adjunct \eqref{BasedLoopingSuspensionAdjunction}
  of each $\sigma_n$ is a weak equivalence:
    \vspace{-2mm}
  \begin{equation}
    \label{OmegaSpectrum}
    \Big\{
    E_n \in \TopologicalSpaces^{\ast/}_{\mathrm{Qu}}
    \,,
    \;
    \xymatrix{
      E_n \ar[r]^-{ \tilde \sigma_n }_-{\in\; \mathrm{W}} & \Omega E_{n+1}
    }
    \Big\}_{n \in \mathbb{N}}
  \end{equation}
  \end{itemize}
\end{example}

\begin{example}[Derived stabilization adjunction]
  \label{DerivedStabilizationAdjunction}
  The suspension/looping Quillen adjunction
  on pointed spaces (Example \ref{LoopingSuspensionAdjunction})
  extends to a commuting diagram of Quillen adjunctions
  (Def. \ref{QuillenAdjunction})
  to and on the stable model category of spectra
  (Example \ref{ModelStructureOnSequentialSpectra})
    \vspace{-3mm}
  \begin{equation}
    \label{StabilizationQuillenAdjunction}
    \raisebox{23pt}{
    \xymatrix@R=.7em@C=4em{
      \TopologicalSpaces^{\ast/}_{\mathrm{Qu}}
      \;
      \ar@{<-}@<+6pt>[rr]^-{ \Sigma }
      \ar@<-6pt>[rr]_-{ \Omega }^-{ \bot_{\mathrlap{\mathrm{Qu}}} }
      \ar@{->}@<-9pt>[dd]_-{  \Sigma^\infty}^-{ \dashv_{\mathrm{Qu}} }
      \ar@{<-}@<+9pt>[dd]^-{  \Omega^\infty}
      &&
      \;
      \TopologicalSpaces^{\ast/}_{\mathrm{Qu}}
      \ar@{->}@<-9pt>[dd]_-{  \Sigma^\infty}^-{ \dashv_{\mathrm{Qu}} }
      \ar@{<-}@<+9pt>[dd]^-{  \Omega^\infty}
      \\
      \\
      \mathrm{SequentialSpectra}_{\mathrm{BF}}
      \;
      \ar@{<-}@<+6pt>[rr]^-{ \Sigma }
      \ar@<-6pt>[rr]_-{ \Omega }^-{ \bot_{\mathrlap{\mathrm{Qu}}} }
      &&
      \;
      \mathrm{SequentialSpectra}_{\mathrm{BF}}
      \,.
    }
    }
  \end{equation}

  \vspace{-2mm}
  \noindent
  such that the bottom adjunction is a Quillen equivalence
  (Def. \ref{QuillenEquivalence}), hence such that
  under passage to derived adjunctions (Prop. \ref{DerivedFunctors})
      \vspace{-3mm}
  \begin{equation}
    \label{StabilizationAdjunction}
    \raisebox{40pt}{
    \xymatrix@R=.7em@C=4em{
      \mathrm{Ho}
      \Big(
        \TopologicalSpaces^{\ast/}_{\mathrm{Qu}}
      \big)
      \;
      \ar@{<-}@<+6pt>[rr]^-{
        \LeftDerived\Sigma
      }
      \ar@<-6pt>[rr]_-{
        \RightDerived\Omega
      }^-{
        \bot
      }
      \ar@{->}@<-6pt>[dd]_-{
        \LeftDerived\Sigma^\infty
      }^-{
        \dashv
      }
      \ar@{<-}@<+6pt>[dd]^-{
        \RightDerived\Omega^\infty
      }
      &&
      \;
      \mathrm{Ho}
      \Big(
        \TopologicalSpaces^{\ast/}_{\mathrm{Qu}}
      \big)
      \ar@{->}@<-6pt>[dd]_-{
        \LeftDerived\Sigma^\infty
      }^-{
        \dashv
      }
      \ar@{<-}@<+6pt>[dd]^-{
        \RightDerived\Omega^\infty
      }
      \\
      \\
      \mathrm{Ho}
      \Big(
        \mathrm{SequentialSpectra}_{\mathrm{BF}}
      \big)
      \;
      \ar@{<-}@<+6pt>[rr]^-{
        \LeftDerived\Sigma
      }
      \ar@<-6pt>[rr]_-{
        \RightDerived\Omega
      }^-{
        \simeq
      }
      &&
      \;
      \mathrm{Ho}
      \Big(
        \mathrm{SequentialSpectra}_{\mathrm{BF}}
      \big)
    }
    }
  \end{equation}

      \vspace{-2mm}
  \noindent
  the bottom adjunction is an equivalence, thus exhibiting
  the homotopy category of spectra as being
  \emph{stable} under looping/suspension.

  \noindent
  We say that

  \noindent
  {\bf (i)}
  $\StableHomotopyTypes$ is the \emph{stable homotopy category} of spectra;

  \noindent
  {\bf (ii)}
  the vertical adjunction
  $(\LeftDerived\Sigma^\infty \dashv  \RightDerived\Omega^\infty)$
  is the \emph{stabilization adjunction} between homotopy types
  \eqref{ClassicalHomotopyCategory} and spectra.

  \noindent
  {\bf (iii)} the images of $\Sigma^\infty$ are the \emph{suspension spectra}.

  \noindent
  {\bf (iv)}
  For $E \,\in\, \StableHomotopyTypes$ and $n \in \mathbb{N}$ we
  write (for brevity and in view of \eqref{OmegaSpectrum})
      \vspace{-2mm}
  \begin{equation}
    \label{ComponentSpacesOfSpectra}
    E_n
    \;:=\;
    \RightDerived\Omega^\infty
    \big(
      (
        \LeftDerived\Sigma
      )^{n}
      E
    \big)
    \;\;\;\;
    \in
    \;
    \PointedHomotopyTypes
  \end{equation}

      \vspace{-2mm}
  \noindent
  for the homotopy type of the $n$th component space of
  the spectrum.
\end{example}

\medskip

\noindent {\bf Smooth $\infty$-stacks.}
We briefly highlight some basics of smooth $\infty$-stack theory,
for more details and more exposition see
\cite[\S 2]{FSS12b}\cite{FSS13a}\cite{dcct}\cite[\S 1]{SS20a}.

\begin{defn}[Simplicial presheaves over Cartesian spaces]
 \label{SimplicialPresheavesOverCartesianSpaces}
  We write
  \\
    \noindent
    {\bf (i)}
    \vspace{-2mm}
    \begin{equation}
      \label{CartesianSpaces}
      \CartesianSpaces
      \;\;:=\;\;
      \big\{\!\!
      \xymatrix{
        \mathbb{R}^{n_1}
        \ar[r]^-{ \mathrm{smooth} }
        &
        \mathbb{R}^{n_2}
      }
     \!\! \big\}_{ n_i \in \mathbb{N} }
     \xhookrightarrow{\;\;}
     \SmoothManifolds
    \end{equation}

    \vspace{-1mm}
    \noindent
    for the category whose objects are the
    Cartesian spaces $\mathbb{R}^n$, for $n \in \mathbb{N}$,
    and whose morphisms are the \emph{smooth} functions between these
    (hence the full subcategory of $\SmoothManifolds$
    on the Cartesian spaces).

    \noindent
    {\bf (ii)}
    \begin{equation}
      \label{SimplicialPresheavesOverCartesianSpaces}
      \mathrm{PSh}
      \big(
        \CartesianSpaces,
        \SimplicialSets
      \big)
      \;:=\;
      \mathrm{Functors}
      \big(
        \CartesianSpaces^{\mathrm{op}},
        \SimplicialSets
      \big)
    \end{equation}
    for the category of functors from the opposite of
    $\CartesianSpaces$ \eqref{CartesianSpaces}
    to $\SimplicialSets$ (Ex. \ref{ClassicalModelStructureOnSimplicialSets}).
\end{defn}

\begin{example}[Model structure on simplicial presheaves over Cartesian spaces {\cite{Dugger98}\cite{Dugger01}\cite[\S A]{FiorenzaSchreiberStasheff10}}]
    \label{ModelStructureOnSimplicialPresheavesOverCartesianSpaces}
    The category of simplicial presheaves over
    Cartesian spaces (Prop. \ref{SimplicialPresheavesOverCartesianSpaces})
    carries the following model category structures
    (Def. \ref{ModelCategories}):

    \noindent
    {\bf (i)}
    The global projective model structure
    \vspace{-2mm}
    \begin{equation}
      \label{GlobalSimplicialPresheavesOnCartesianSpaces}
      \mathrm{PSh}
      \big(
        \CartesianSpaces
        ,\,
        \SimplicialSets_{\mathrm{Qu}}
      \big)_{\mathrm{proj}}
      \;\;\;
      \in
      \;
      \mathrm{ModelCategories}
    \end{equation}

    \vspace{-2mm}
    \noindent
    whose

    \begin{itemize}
    \vspace{-4mm}
    \item[$\mathrm{W}$ --] {\it weak equivalences} are the morphisms
    which over each $\mathbb{R}^n$ are
    weak equivalence in $\SimplicialSets_{\mathrm{Qu}}$
    (Example \ref{ClassicalModelStructureOnSimplicialSets}),

    \vspace{-2mm}
    \item[$\mathrm{Fib}$ --]  {\it fibrations} are the morphisms
    which over each $\mathbb{R}^n$ are
    fibrations in $\SimplicialSets_{\mathrm{Qu}}$
    (Example \ref{ClassicalModelStructureOnSimplicialSets}),
    \end{itemize}

    \noindent
    {\bf (ii)}
    The local projective model structure
    \vspace{-2mm}
    \begin{equation}
      \label{SimplicialPresheavesOnCartesianSpaces}
      \mathrm{PSh}
      \big(
        \CartesianSpaces
        ,\,
        \SimplicialSets_{\mathrm{Qu}}
      \big)_{\mathrm{proj} \atop \mathrm{loc}}
      \;\;\;
      \in
      \;
      \mathrm{ModelCategories}
    \end{equation}

    \vspace{-3mm}
    \noindent
    whose:

    \begin{itemize}
    \vspace{-4mm}
    \item[$\mathrm{W}$ --]  {\it weak equivalences} are the morphisms
    whose stalk at $0 \in \mathbb{R}^n$ is a weak
    equivalence in $\SimplicialSets_{\mathrm{Qu}}$
    (Example \ref{ClassicalModelStructureOnSimplicialSets}),
    for all $n \in \mathbb{N}$;

    \vspace{-2mm}
    \item[$\mathrm{Cof}$ --]  {\it cofibrations} are the morphisms
    with the left lifting property \eqref{LiftingProperty}
    against the class of morphisms which over each $\mathbb{R}^n$
    are in $\mathrm{Fib} \cap \mathrm{W}$
    of $\SimplicialSets_{\mathrm{Qu}}$.
    \end{itemize}
\end{example}
\begin{example}[Smooth manifolds as simplicial presheaves]
 \label{SmoothManifoldAsInfinityStack}
  Consider a smooth manifold. $X \in \SmoothManifolds$.

  \noindent
  {(\bf i)}
  The manifold is incarnated as a simplicial presheaf (Def. \ref{ModelStructureOnSimplicialPresheavesOverCartesianSpaces})
  by the rule which assigns to a Cartesian space
  the set of smooth functions
  $\mathbb{R}^n \xrightarrow{ \mathrm{smooth} } X$,
  regarded as a simplicially constant simplicial set.
  This construction constitutes a full subcategory inclusion:
  \vspace{-2mm}
  \begin{equation}
    \label{SmoothManifoldAsASimplicialPresheaf}
    \begin{tikzcd}[row sep=-5pt, column sep=small]
      \SmoothManifolds
      \ar[
        rr,
        hook
      ]
      &&
      \mathrm{PSh}
      \big(
        \CartesianSpaces
        ,\,
        \SimplicialSets
      \big)
      \\
      X
      &\longmapsto&
      \Big(
        \mathbb{R}^n
        \mapsto
        \big(
          [k]
          \mapsto
          \SmoothManifolds
          (\mathbb{R}^n\,,\,X)
        \big)
      \Big)
      \,.
    \end{tikzcd}
  \end{equation}

\vspace{-2mm}
  \noindent
  {\bf (ii)}
  For $p_n \in \mathbb{R}^n$ any point, the {\it stalk}
  of this presheaf is the set of {\it germs of smooth functions}
  from an open neighbourhood of $p_n$ to $X$.
  This set depends, in general, on $n \in \mathbb{N}$,
  but does not depend on the choice of $p_n$.
\end{example}

\begin{example}[Lie groupoids as simplicial presheaves]
  \label{LieGroupoidsAsSimplicialPresheaves}
  Consider a Lie groupoid
  $\mathcal{G} = \big( \mathcal{G}_1 \rightrightarrows \mathcal{G}_0\big)$
  (review in \cite{Mackenzie87}\cite{MoerdijkMrcun03} \cite{Mackenzie05})
  hence a groupoid internal to smooth manifolds
  $\mathcal{G} \in \Groupoids(\SmoothManifolds)$.
  Notice that for each $\mathbb{R}^n \,\in\, \CartesianSpaces$
  there is an induced bare groupoid of smooth functions into the
  component manifolds:
  \vspace{-2mm}
  $$
    \mathbb{R}^n
    \;\;
      \longmapsto
    \;\;
    \mathcal{G}(\mathbb{R}^n)
    \;:=\;
    \big(
      \SmoothManifolds(\mathbb{R}^n ,\, \mathcal{G}_1)
      \rightrightarrows
      \SmoothManifolds(\mathbb{R}^n ,\, \mathcal{G}_0)
    \big)
    \;\;
    \in
    \;
    \Groupoids(\Sets)
    \,.
  $$

  \vspace{-2mm}
  \noindent
  The simplicial nerves (Ex. \ref{SimplicialNervesOfGroupoids}) of these
  mapping groupoids arrange into a simplicial presheaf (Ex. \ref{ModelStructureOnSimplicialPresheavesOverCartesianSpaces})
  and this construction is the inclusion of a full subcategory,
  extending the full inclusion of smooth manifolds \eqref{SmoothManifoldAsASimplicialPresheaf}:
\vspace{-2mm}
  $$
    \begin{tikzcd}[row sep=-5pt]
      \Groupoids
      \big(
        \SmoothManifolds
      \big)
      \ar[
        rr,
        hook
      ]
      &&
      \mathrm{PSh}
      \big(
        \CartesianSpaces
        ,\,
        \SimplicialSets
      \big)
      \\
      \mathcal{G}
      &\longmapsto&
      \big(
        \mathbb{R}^n
        \,\mapsto\,
        N
        \big(
          \mathcal{G}(\mathbb{R}^n)
        \big)
      \big)
    \end{tikzcd}
  $$
\end{example}

\begin{example}[{\v C}ech groupoids of open covers as simplicial presheaves]
  \label{NerveOfOpenCover}
  Let $X$ be a smooth manifold equipped with
  a cover by a set of open subsets
  $
    \big\{
      U_i \xhookrightarrow{\mathrm{open}} X
    \big\}_{i \in I}
  $.

  \noindent
  {\bf {(i)}}
  The {\it {\v C}ech nerve} of the open cover
  is the
  simplicial presheaf (Def. \ref{ModelStructureOnSimplicialPresheavesOverCartesianSpaces})
  \vspace{-2mm}
  \begin{equation}
    \label{CechNerveOfCover}
    N\big(\{U_i\}_{i \in I}\big)
    \;\in\;
    \mathrm{PSh}
    \big(
      \CartesianSpaces
      ,\,
      \SimplicialSets
    \big)
  \end{equation}

\vspace{-2mm}
  \noindent
  whose $k$-cells over any $\mathbb{R}^n$ are the smooth functions
  \vspace{-3mm}
    \begin{equation}
    \label{CechNervePresheafAsFunctor}
    N\big(\{U_i\}_{i \in I}\big)
    \;:\;
    \big(
      \mathbb{R}^n
      ,\,
      [k]
    \big)
    \;\;
      \longmapsto
    \;\;
    \SmoothManifolds
    \Big(
      \mathbb{R}^n
      ,\,
      \big( \sqcup_i U_i \big)^{\times_X^{(k+1)}}
    \Big)
  \end{equation}

  \vspace{-2mm}
  \noindent
  into the $(k+1)$-fold intersections of the patches $U_i$ in $X$:
  \vspace{-2mm}
  \begin{equation}
    \label{FoldIntersectionsOfOpenSubsets}
    \big(\sqcup_i U_i\big)^{\times^{k+1}_X}
      \;\coloneqq\;
    \underset{
      \mbox{\tiny $k+1$ factors}
    }{
    \underbrace{
    \big( \sqcup_i U_i \big)
      \times_X
      \cdots
      \times_X
    \big( \sqcup_i U_i \big)
    }
    }
    \xhookrightarrow{\;\iota_{k + 1}\;}
    X
  \,.
\end{equation}

\vspace{-2mm}
  \noindent
This is, for each $\mathbb{R}^n$, a Kan complex \eqref{KanComplexes}
and as such is the nerve of the Lie groupoid (Ex. \ref{LieGroupoidsAsSimplicialPresheaves})
which is {\it smooth {\v C}ech-groupoid} of the open cover:

\vspace{-.4cm}
\begin{equation}
  N\big(\{U_i\}_{i \in I}\big)
  \;=\;
  N
  \Big(
    \big(
      \sqcup_i U_i
    \big)
    \times_X
    \big(
      \sqcup_i U_i
    \big)
    \rightrightarrows
    \big(
      \sqcup_i U_i
    \big)
  \Big)
  \,.
\end{equation}
\vspace{-.4cm}

\noindent
({\bf ii})
For any point $p_n \in \mathbb{R}^n$, the
{\it stalk} of the {\v C}ech nerve \eqref{CechNerveOfCover}
at $p_n$ is the disjoint union over the germs of smooth functions
$\mathbb{R}^n \xrightarrow{\phi} X$ of the
nerves of the pair groupoids on the subset $I_{\phi(p_n)} \subset X$
of patches $U_i$ that contain $\phi(p_n)$.

\noindent
{\bf (iii)}
Postcomposition with the inclusions \eqref{FoldIntersectionsOfOpenSubsets}
yields a canonical morphism of simplicial presheaves
from the cover's {\v C}ech nerve \eqref{CechNervePresheafAsFunctor}
to the presheaf incarnation \eqref{SmoothManifoldAsASimplicialPresheaf}
of the underlying manifold:
\vspace{-2mm}
\begin{equation}
  \label{ComparisonMapFromCechNerveToManifold}
  \begin{tikzcd}[row sep=-5pt, column sep=small]
    N
    \big(
      \{U_i\}_{i \in I}
    \big)
    \ar[
      rr,
      "\in \mathrm{W}"
    ]
    &&
    X
    \\
    \big(
      \mathbb{R}^n
      \xrightarrow{\phi}
      \big(\sqcup_i U_i\big)^{\times^{k+1}_X}
    \big)
    &\longmapsto&
    \big(
      \mathbb{R}^n
      \xrightarrow{\phi}
      \big(\sqcup_i U_i\big)^{\times^{k+1}_X}
      \xrightarrow{ \iota_{k+1} }
      X
    \big)
    \,.
  \end{tikzcd}
  {\phantom{AAA}}
  \in
  \mathrm{PSh}
  \big(
    \CartesianSpaces
    ,\,
    \SimplicialSets_{\mathrm{Qu}}
  \big)_{ {\mathrm{proj}} \atop {\mathrm{loc}} }
\end{equation}

\vspace{-2mm}
  \noindent
On stalks, this map takes the
nerve of the pair groupoid on the set of factorizations through the patches $U_i$
of
a the germ of a given smooth $\mathbb{R}^n \xrightarrow{\phi} X$
to that germ itself.
Since nerves of pair groupoids are contractible \eqref{NerveOfPairGroupoidIsContractible},
this means that
\eqref{ComparisonMapFromCechNerveToManifold} is a weak equivalence
in the local model structure of Ex. \ref{ModelStructureOnSimplicialPresheavesOverCartesianSpaces}.
\end{example}
This says that ({\v C}ech nerves of) open covers serve as resolutions of
smooth manifolds in the local model structure Ex. \ref{ModelStructureOnSimplicialPresheavesOverCartesianSpaces};
in fact as {\it cofibrant resolutions} if the cover is ``good'':

\begin{prop}[Dugger's cofibrancy recognition {\cite[Cor . 9.4]{Dugger01}}]
  \label{DuggerCofibrancyRecognition}
  A sufficient condition for
  $\mathcal{X}
    \in
  \mathrm{PSh}(\CartesianSpaces, \SimplicialSets)_{ {\mathrm{proj}} \atop {\mathrm{loc}} }$
  (Ex. \ref{ModelStructureOnSimplicialPresheavesOverCartesianSpaces})
  to be cofibrant (Nota. \ref{FibrantAndCofibrantObjects})
  is that in each simplicial degree $k$,
  the component presheaf $X_k$ is

\vspace{1mm}
\noindent  {\bf (i)} a coproduct (as presheaves, using Ex. \ref{SmoothManifoldAsInfinityStack})
    of Cartesian spaces:
    $\mathcal{X}_k \simeq \underset{i_k}{\coprod} \, \mathbb{R}^{n_{i_k}}$;

 \noindent  {\bf (ii)} whose degenerate cells split off as a disjoint summand.
\end{prop}

\begin{example}[Good open covers are projectively cofibrant resolutions of smooth manifolds]
  \label{GoodOpenCoversAreProjectivelyCofibrantResolutionsOfSmoothManifolds}
  Prop. \ref{DuggerCofibrancyRecognition}
  applied to Ex. \ref{NerveOfOpenCover} says that
  the {\v C}ech nerve of an open cover is a cofibrant resolution
  of the underlying manifold if the open cover is
  {\it good}, or rather: {\it differentiably good},
  in that each non-empty intersection of
  a finite number of its patches is
  diffeomorphic to an open ball
  (namely, equivalently: to a Cartesian space):
  \vspace{-2mm}
  \begin{equation}
    \label{GoodOpenCoversAsCofibrantResolutions}
    \mbox{
      $
      \big\{
        U_i \xhookrightarrow{\mathrm{open}} X
      \big\}_{i \in I}
      $
      is good
    }
    {\phantom{AA}}
    \Rightarrow
    {\phantom{AA}}
    \begin{tikzcd}
      \varnothing
      \ar[
        r,
        "{\in \mathrm{Cof}}"{below}
      ]
      &
      N
      \big(
        \{U_i\}_{i \in I}
      \big)
      \ar[
        rr,
        "{\in \mathrm{W} }"{below},
        "{ p_{\{U_i\}_i} }"
      ]
      &&
      X
    \end{tikzcd}
    \;\;
    \in
    \;
    \mathrm{PSh}
    \big(
      \CartesianSpaces
      ,\,
      \SimplicialSets
    \big)_{ {\mathrm{proj}} \atop {\mathrm{loc}} }
    \,.
  \end{equation}

  \vspace{-2mm}
  \noindent
  Notice that every smooth manifold admits a
  differentiably good open cover \cite[Prop. A.1]{FiorenzaSchreiberStasheff10}.
\end{example}

\begin{example}[Hom-complexes of simplicial presheaves]
\label{HomComplexesOfSimplicialPresheaves}
$\,$

\noindent
{\bf (i)}
For
$\mathcal{X} \,\in\, \mathrm{PSh}\big(\CartesianSpaces ,\, \SimplicialSets \big)$
and
(Def. \ref{SimplicialPresheavesOverCartesianSpaces})
$S \,\in\, \SimplicialSets$ (Ex. \ref{ClassicalModelStructureOnTopologicalSpaces})
there is the {\it tensored} simplicial presheaf
$$
  \mathcal{X} \times S
  \;\;
  \in
  \;
  \mathrm{PSh}\big(\CartesianSpaces ,\, \SimplicialSets \big)
$$
given by value-wise Cartesian product of simplicial sets:
\begin{equation}
  \label{TensoringFormulaOfSimplicialPresheavesWithSimplicialSets}
  \mathcal{X} \times S
  \;\;
  :
  \;\;
  \mathbb{R}^n
  \;\;
    \mapsto
  \;\;
  \mathcal{X}(\mathbb{R}^n) \times S
  \,.
\end{equation}

\noindent
{\bf (ii)}
For
$\mathcal{X}, \mathcal{A} \,\in\, \mathrm{PSh}\big(\CartesianSpaces ,\, \SimplicialSets \big)$
the
{\it simplicial hom-complex}
from $\mathcal{X}$ to $\mathcal{A}$
is the simplicial set of morphisms of simplicial presheaves
\begin{equation}
  \label{SimplicialHomComplexOfSimplicialPresheaves}
  \overset{
    \mathclap{
    \raisebox{3pt}{
      \tiny
      \color{darkblue}
      \bf
      simplicial mapping complex
      }
    }
  }{
  \mathrm{Maps}
  \big(
    \mathcal{X}
    ,\,
    \mathcal{A}
  \big)
  }
  \;\;
    \coloneqq
  \;\;
  \mathrm{PSh}
  \big(
    \CartesianSpaces
    ,\,
    \SimplicialSets
  \big)
  \big(
    \mathcal{X} \!\times \Delta[\bullet]
    ,\;
    \mathcal{A}
  \big)
  \;\;\;
  \in
  \;
  \SimplicialSets
  \,.
\end{equation}
into $\mathcal{A}$
out of the tensoring \eqref{TensoringFormulaOfSimplicialPresheavesWithSimplicialSets}
of $\mathcal{X}$ with the simplicial simplices
$\Delta[n] \,\in\, \SimplicialSets$, $n \in \mathbb{N}$.
Its image in the classical homotopy category
(Ex. \ref{TheClassicalHomotopyCategory}) is the {\it mapping space}
\begin{equation}
  \label{MappingSpaceBetweenSmoothStacks}
  \mathrm{Maps}
  \big(
    \mathcal{X}
    ,\,
    \mathcal{A}
  \big)
  \;\;\;
  \in
  \;
  \mathrm{Ho}
  \big(
    \SimplicialSets_{\mathrm{Qu}}
  \big)
  \,.
\end{equation}

\noindent
{\bf (iii)} The (simplicially enriched) Yoneda lemma says that
simplicial hom-complexes \eqref{SimplicialHomComplexOfSimplicialPresheaves}
out of a Cartesian space \eqref{CartesianSpaces}
regarded as a simplicial presheaf via Ex. \ref{SmoothManifoldAsInfinityStack}:
\begin{equation}
  \label{SimplicialYonedaLemma}
  \mathcal{X}
  \big(
    \mathbb{R}^n
  \big)
  \;\;
  \simeq
  \;\;
  \mathrm{Maps}
  \big(
    \mathbb{R}^n
    ,\,
    \mathcal{X}
  \big)
  \,.
\end{equation}
\end{example}

\begin{prop}[Smooth $\infty$-Stacks]
  \label{SmoothInfinityStacks}
  The fibrant objects (Nota. \ref{FibrantAndCofibrantObjects})
  in the local projective model structure
  \eqref{SimplicialPresheavesOnCartesianSpaces},
  are   to be called the {\it smooth $\infty$-stacks}
  (or {\it smooth $\infty$-groupoids})

  \vspace{-3mm}
  \begin{equation}
    \label{SmoothStacksAsFibrantObjects}
    \mathrm{SmoothStacks}_\infty
    \;\coloneqq\;
    \Big(
      \mathrm{PSh}
      \big(
        \CartesianSpaces
        ,\,
        \SimplicialSets_{\mathrm{Qu}}
      \big)_{ {\mathrm{proj}} \atop  {\mathrm{loc}}  }
    \Big)^{\mathrm{fib}}
    \,,
  \end{equation}
  \vspace{-.4cm}

  \noindent
  are precisely those simplicial presheaves
  which:

  \noindent
  {\bf (i)}
   are {\it presheaves of $\infty$-groupoids}
   in that they take values in Kan complexes \eqref{KanComplexes};

  \noindent
  {\bf (ii)}
  respect gluing of patches of good open covers of Cartesian spaces
  (``satisfy descent'') in that for each
  $n \in \mathbb{N}$ and each good open cover
  $\big\{ U_i \xhookrightarrow{\;} \mathbb{R}^n \big\}_{i \in I}$
  (Ex. \ref{GoodOpenCoversAreProjectivelyCofibrantResolutionsOfSmoothManifolds})
  the following map \eqref{DescentMapOfHomComplexes}
  of simplicial hom-complexes \eqref{SimplicialHomComplexOfSimplicialPresheaves}
  --
  induced by precomposition with the comparison morphism
  \eqref{GoodOpenCoversAsCofibrantResolutions}
  from the {\v C}ech nerve
  \eqref{ComparisonMapFromCechNerveToManifold}
  --
  is a weak equivalence of simplicial sets (Ex. \ref{ClassicalModelStructureOnSimplicialSets}):

  \vspace{-2mm}
  \begin{equation}
    \label{DescentMapOfHomComplexes}
    \begin{tikzcd}
      \mathcal{X}(\mathbb{R}^n)
      \simeq
      \mathrm{Maps}
      \big(
        \mathbb{R}^n
        ,\,
        \mathcal{X}
      \big)
      \ar[
        rrr,
        "{
          \mathrm{Maps}
          (
            p_{\{U_i\}_i}
            ,\,
            \mathcal{X}
          )
        }"{above},
        "{\in \mathrm{W}}"{below}
      ]
      &&&
      \mathrm{Maps}
      \Big(
        N
        \big(
          \{U_i\}_{i \in I}
        \big)
        ,\,
        \mathcal{X}
      \Big)
    \end{tikzcd}
    \;\;\;
    \in
    \;
    \SimplicialSets_{\mathrm{Qu}}
    \,.
  \end{equation}
\end{prop}
\begin{proof}
  By the discussion in \cite[\S 5.1]{Dugger01}
  the claimed condition characterizes the fibrant objects
  in the left Bousfield localization of the global projective model category
  \eqref{GlobalSimplicialPresheavesOnCartesianSpaces}
  at the {\v C}ech nerve projections \eqref{ComparisonMapFromCechNerveToManifold}
  By \cite[Prop. 3.4.8]{Dugger98} this left Bousfield localization is
  the local model structure \eqref{SimplicialPresheavesOnCartesianSpaces}.
\end{proof}
\begin{defn}[Homotopy category of smooth $\infty$-stacks]
  \label{SmoothInfinityStacks}
  In view of Prop. \ref{SmoothInfinityStacks},
  we write
  \vspace{-2mm}
  \noindent
  \begin{equation}
    \label{HomotopyCategoryOfSimplicialPresheaves}
    \Stacks
    \;\;
      :=
    \;\;
    \mathrm{Ho}
    \Big(
      \mathrm{PSh}
      \big(
        \CartesianSpaces,
        \SimplicialSets
      \big)_{\mathrm{proj} \atop \mathrm{loc}}
    \Big)
  \end{equation}

  \vspace{-2mm}
  \noindent
  for the homotopy category (Def. \ref{HomotopyCategory})
  of the local projective model category
  of simplicial presheaves over $\CartesianSpaces$
  (Example \ref{ModelStructureOnSimplicialPresheavesOverCartesianSpaces}).
  We say that the objects of $\Stacks$
  \eqref{HomotopyCategoryOfSimplicialPresheaves}
  are \emph{smooth $\infty$-stacks}.
\end{defn}

\begin{example}[Truncated smooth $\infty$-stacks {\cite[Ex. 3.18]{SS20b}}]
  \label{SmoothSpaces}
  $\,$

  \noindent
  {\bf (i)}
  Those smooth $\infty$-stacks (Def. \ref{SmoothInfinityStacks})
  which take values in
  2-coskeletal, hence 1-truncated,
  Kan complexes (Prop. \ref{SkeletaAndPostnikovTowers})
  are 1-groupoid valued, hence are {\it smooth 1-stacks} or
  just smooth {\it stacks} \cite{Jardine01}\cite{Hollander08}.

  \noindent
  {\bf (ii)}
  Those smooth $\infty$-stacks which are 0-truncated
  take values in sets and hence are {\it sheaves} on
  $\CartesianSpaces$. We call these {\it smooth spaces}.
  The {\it concrete sheavees} among these are
  the {\it diffeological spaces}
  (\cite{Souriau80}\cite{Souriau84}\cite{IZ85}, see \cite{BaezHoffnung08}\cite{IZ13}).

  \vspace{-.5cm}
  $$
    \begin{tikzcd}
      \overset{
        \mathclap{
        \raisebox{3pt}{
          \tiny
          \color{darkblue}
          \bf
          \def\arraystretch{.9}
          \begin{tabular}{c}
            diffeological
            spaces
          \end{tabular}
        }
        }
      }{
        \mathrm{PSh}(\CartesianSpaces, \, \Sets)^{\mathrm{fib}}_{\mathrm{conc}}
      }
      \ar[r, hook]
      &
      \overset{
        \raisebox{3pt}{
          \tiny
          \color{darkblue}
          \bf
          \def\arraystretch{.9}
          \begin{tabular}{c}
            smooth sets
            \\
            (smooth spaces)
          \end{tabular}
        }
      }{
        \mathrm{PSh}(\CartesianSpaces, \, \Sets)^{\mathrm{fib}}
      }
      \ar[r, hook]
      &
      \overset{
        \mathclap{
        \raisebox{3pt}{
          \tiny
          \color{darkblue}
          \bf
          \def\arraystretch{.9}
          \begin{tabular}{c}
            smooth
            groupoids
            \\
            (smooth stacks)
          \end{tabular}
        }
        }
      }{
        \mathrm{PSh}(\CartesianSpaces, \SimplicialSets_{\mathrm{cosk}_2})^{\mathrm{fib}}
      }
      \ar[r, hook]
      &
      \overset{
        \mathclap{
        \raisebox{3pt}{
          \tiny
          \color{darkblue}
          \bf
          \def\arraystretch{.9}
          \begin{tabular}{c}
            smooth
            $\infty$-groupoids
            \\
            (smooth $\infty$-stacks)
          \end{tabular}
        }
        }
      }{
        \mathrm{PSh}(\CartesianSpaces, \SimplicialSets)^{\mathrm{fib}}
      }
    \end{tikzcd}
  $$
\end{example}

\begin{lemma}[$\infty$-Stackification preserves finite homotopy limits]
 \label{InfinityStackification}
 The identity functors constitute a Quillen adjunction
 (Def. \ref{QuillenAdjunction})
 between the local and the global projective model categories
 of Example \ref{ModelStructureOnSimplicialPresheavesOverCartesianSpaces}:
 \vspace{-2mm}
 $$
   \xymatrix{
      \mathrm{PSh}
      \big(
        \CartesianSpaces,
        \SimplicialSets
      \big)_{\mathrm{proj} \atop \mathrm{loc} }
\;\;      \ar@{<-}@<+7pt>[rr]^-{ \mathrm{id} }
      \ar@<-7pt>[rr]_-{ \mathrm{id} }^-{ \bot_{\mathrlap{\mathrm{Qu}}} }
      &&
 \;\;     \mathrm{PSh}
      \big(
        \CartesianSpaces,
        \SimplicialSets
      \big)_{\mathrm{proj} }
   }.
 $$

   \vspace{-2mm}
    \noindent
 Moreover, this is such that the derived left adjoint functor
 (Prop. \ref{DerivedFunctors})
   \vspace{-2mm}
 \begin{equation}
   \label{InfinityStackificationFunctor}
   \hspace{-3mm}
   L^{\mathrm{loc}}
   :
   \xymatrix{
      \mathrm{Ho}
      \Big(
      \mathrm{PSh}
      \big(
        \CartesianSpaces,
        \SimplicialSets
      \big)_{\mathrm{proj} }
      \Big)
      \ar[rr]^-{ \LeftDerived\, \mathrm{id} }
      &&
      \Stacks
   }
 \end{equation}

   \vspace{-2mm}
    \noindent
 (the $\infty$-stackification functor)
 preserves homotopy pullbacks (Def. \ref{HomotopyPullback}).
\end{lemma}

  \begin{prop}[Shape Quillen adjunction {\cite[Prop. 4.4.8]{dcct}\cite[Example 3.18]{SS20a}}]
    \label{ShapeQuillenAdjunction}
    We have a Quillen adjunction (Def. \ref{QuillenAdjunction})
      \vspace{-2mm}
    $$
      \xymatrix{
        \mathrm{PSh}
        \big(
          \CartesianSpaces,
          \SimplicialSets
        \big)_{\mathrm{proj} \atop \mathrm{loc}}
        \;\;
        \ar@<+7pt>[rr]^-{ \mathrm{Shp} }
        \ar@{<-}@<-7pt>[rr]_-{ \mathrm{Disc} }^-{
          \bot_{\mathrlap{\mathrm{Qu}}}
        }
        &&
        \;\;
        \SimplicialSets_{\mathrm{Qu}}
      }
    $$

    \vspace{-2mm}
    \noindent
    between the projective local model structure on simplicial
    presheaves over $\CartesianSpaces$ (Example \ref{ModelStructureOnSimplicialPresheavesOverCartesianSpaces})
    and the classical model structure on simplicial sets
    (Example \ref{ClassicalModelStructureOnSimplicialSets}),
    hence a derived adjunction (Prop. \ref{DerivedFunctors})
    between homotopy category of $\infty$-stacks
    (Def. \ref{SmoothInfinityStacks})
    and the classical homotopy category
    (Example \ref{TheClassicalHomotopyCategory})
    \vspace{-2mm}
    $$
      \xymatrix{
        \Stacks
        \;\;
        \ar@<+7pt>[rr]^-{ \LeftDerived\mathrm{Shp} }
        \ar@{<-}@<-7pt>[rr]_-{ \RightDerived\mathrm{Disc} }^-{
          \bot_{\mathrlap{\mathrm{Qu}}}
        }
        &&
    \;\;    \mathrm{Ho}
        \big(
          \SimplicialSets_{\mathrm{Qu}}
        \big)
        \;\;
      }
    $$

    \vspace{-2mm}
\noindent    whose (underived) right adjoint sends a simplicial set to the
    presheaf which is constant on that simplicial set:
      \vspace{-2mm}
    \begin{equation}
      \label{ConstantPresheafOverCartesianSpaces}
      \mathrm{Disc}(S)
      \;:=\;
      \mathrm{const}(S)
      \;:\;
      \big(
        \mathbb{R}^n
        \mapsto
        S
      \big).
    \end{equation}
  \end{prop}

\medskip

\noindent {\bf Homological algebra.}

\begin{example}[Projective model structure on connective chain complexes {\cite[\S II.4 (5.)]{Quillen67}}]
  \label{ProjectiveModelStructureOnConnectiveChainComplexes}
  The category $\ZChainComplexes$
  of connective chain complexes
  of abelian groups
  (i.e. concentrated in non-negative degrees with differential
  of degree -1)
  carries a model category structure  (Def. \ref{ModelCategories})
  whose
  \begin{itemize}
    \vspace{-2mm}
    \item[$\mathrm{W}$ --]
    weak equivalences are the quasi-isomorphisms
      (those inducing isomorphisms on all chain homology groups)

    \vspace{-2mm}
    \item[$\mathrm{Fib}$ --]
    fibrations are the morphisms that are surjections in each positive degree

    \vspace{-2mm}
    \item[$\mathrm{Cof}$ --]  cofibrations are the morphisms with
    degreewise injective kernels.
  \end{itemize}

 \vspace{-2mm}
  \noindent
  We write
  $\ZChainComplexesProj$ for this model category.
\end{example}

More generally:

\vspace{-1mm}
\begin{example}[Projective model structure on presheaves of connective chain complexes {\cite[p. 7]{Jardine03}}]
\label{ProjectiveModelStructureOnPresheavesOfConnectivedChainComplexes}
The category of presheaves of connective chain complexes
over $\CartesianSpaces$ \eqref{CartesianSpaces}
carries the structure of a model category
whose weak equivalences and fibrations are objectwise
those of $\ZChainComplexesProj$ (Ex. \ref{ProjectiveModelStructureOnConnectiveChainComplexes}).
We write
$
  \mathrm{PSh}
  \big(
    \CartesianSpaces
    \,,\,
    \ZChainComplexes
  \big)_{\mathrm{proj}}
$ for this model category.
\end{example}

\begin{prop}[Dold-Kan correspondence {\cite[Thm 1.9]{Dold58}\cite{Kan58}\cite[\S III.2]{GoerssJardine99}\cite[\S 2.1]{SchwedeShipley03}}]
  \label{DoldKanCorrespondence}
  Given $A_\bullet \in \SimplicialAbelianGroups$,
  its \emph{normalized chain complex} is the connective
  chain complex of abelian groups (Example \ref{ProjectiveModelStructureOnConnectiveChainComplexes})
  which in degree $n \in \mathbb{N}$ is the quotient
  of $A_n$ by the degenerate cells
  and whose differential is the alternating sum of the face maps:
    \vspace{-2mm}
  \begin{equation}
    \label{NormalizedChainComplex}
    \vspace{-2mm}
    N(A)_\bullet
    \;:=
    \Big\{
    N(A)_n
    \;:=\;
    A_n / \sigma(A_{n+1})
    \,,
    \;\;
    \partial_n
    \;:=\;
    \sum_{i = 0}^n
    (-1)^i d_i
    \;:
    \xymatrix@C=13pt{
      N(A)_n
      \ar[r]
      &
      N(A)_{n-1}
    }
    \! \Big\}_{n \in \mathbb{N}}
        \in
    \ZChainComplexes
    \,.
  \end{equation}

  \noindent
  {\bf (i)}
  This construction constitutes an adjoint equivalence of categories
    \vspace{-2mm}
  \begin{equation}
    \label{DoldKanEquivalenceOfCategories}
    \xymatrix{
      \ZChainComplexes
     \;\; \ar@{<-}@<+5pt>[rr]^-{ N }
      \ar@<-5pt>[rr]_-{ }^-{ \simeq }
      &&
      \;\;
      \SimplicialAbelianGroups
    }
  \end{equation}

\vspace{-2mm}
  \noindent
  {\bf (ii)}
  such that simplicial homotopy groups
  of $A \in \SimplicialAbelianGroups \to \mathrm{SimplicialSet}$
  are identified with chain homology groups of the normalized
  chain complex (\cite[Cor. III.2.5]{GoerssJardine99}):
    \vspace{-3mm}
  \begin{equation}
    \label{DoldKanHomotopyGroups}
    \pi_\bullet(A)
    \;\simeq\;
    H_\bullet(N A)\;.
  \end{equation}
\end{prop}

\begin{example}[Model structure on simplicial abelian groups {\cite[\S III.2]{Quillen69}\cite[\S 4.1]{SchwedeShipley03}}]
  \label{ModelStructureOnSimplicialAbelianGroups}
  The category  $\SimplicialAbelianGroups$
  carries a model category structure (Def. \ref{ModelCategories})
  whose
  \begin{itemize}
    \vspace{-2mm}
    \item[$\mathrm{W}$ --]  weak equivalences are the morphisms
     which are weak equivaleces as morphisms in
     $\SimplicialSets_{\mathrm{Qu}}$
     (Example \ref{ClassicalModelStructureOnSimplicialSets})

  \vspace{-2mm}
    \item[$\mathrm{Fib}$ -- ] fibrations are the morphisms
     which are fibrations as morphisms in
     $\SimplicialSets_{\mathrm{Qu}}$
     (Example \ref{ClassicalModelStructureOnSimplicialSets})
  \end{itemize}

    \vspace{-2mm}
\noindent  In other words, this is the transferred model structure
  along the free/forgetful adjunction, which
  thus becomes a Quillen adjunction (Def. \ref{QuillenAdjunction}):
    \vspace{-2mm}
  \begin{equation}
    \label{FreeAdjunctionOfSimplicialAbelianGroupsOverSimplicialSets}
    \xymatrix{
      \SimplicialAbelianGroups_{\mathrm{proj}}
     \;\; \ar@{<-}@<+7pt>[rr]^-{ \mathbb{Z}[-] }
      \ar@<-7pt>[rr]_-{  }^-{ \bot_{\mathrlap{\mathrm{Qu}}} }
      &&
    \;\;  \SimplicialSets_{\mathrm{Qu}}
      \,.
    }
  \end{equation}
\end{example}

\begin{prop}[Dold-Kan Quillen equivalence {\cite[\S 4.1]{SchwedeShipley03}\cite[Lemma 1.5]{Jardine03}}]
  \label{DoldKanQuillenEquivalence}
  With respect to the projective model structure
  on connective chain complexes (Example \ref{ProjectiveModelStructureOnConnectiveChainComplexes})
  and the projective model structure on simplicial abelian
  groups (Example \ref{ModelStructureOnSimplicialAbelianGroups})
  the Dold-Kan correspondence
  (Prop. \ref{DoldKanCorrespondence})
  is a Quillen equivalence
  (Def. \ref{QuillenEquivalence}):
    \vspace{-3mm}
  \begin{equation}
    \label{QuillenEquivalenceDoldKan}
    \xymatrix{
      \ZChainComplexesProj
    \;\;  \ar@{<-}@<+7pt>[rr]^-{ N }
      \ar@<-7pt>[rr]_-{ }^-{ \simeq_{\mathrlap{\mathrm{Qu}}} }
      &&
    \;\;  \SimplicialAbelianGroups_{\mathrm{proj}}\,,
    }
  \end{equation}

  \vspace{-2mm}
  \noindent
  where both functors preserve all three classes of morphims,
  $\mathrm{Fib}$, $\mathrm{Cof}$ and $\mathrm{W}$, separately.
\end{prop}

\begin{example}[Dold-Kan construction {\cite[\S 3.2.3]{FiorenzaSchreiberStasheff10}\cite[\S 2.4]{FSS12b}}]
\label{DoldKanConstruction}
{\bf i)} We write $\mathrm{DK}$ for the total right
adjoint in the composite of the free Quillen adjunction
\eqref{FreeAdjunctionOfSimplicialAbelianGroupsOverSimplicialSets}
and the Dold-Kan equivalence \eqref{QuillenEquivalenceDoldKan}:
\vspace{-2mm}
\begin{equation}
  \label{DK}
  \xymatrix{
    \ZChainComplexesProj
    \;\;
      \ar@{<-}@<+7pt>[rr]^-{ N }
      \ar@<-7pt>[rr]_-{ }^-{
        \simeq_{\mathrlap{\mathrm{Qu}}}
      }
    \ar@/_1.4pc/[rrrr]_-{
      \mathrm{DK}
    }
\;\;    &&
    \SimplicialAbelianGroups_{\mathrm{proj}}
      \ar@{<-}@<+7pt>[rr]^-{ \mathbb{Z}[-] }
      \ar@<-7pt>[rr]_-{  }^-{
        \bot_{\mathrlap{\mathrm{Qu}}}
      }
    &&
\;\;    \SimplicialSets_{\mathrm{Qu}}
  }.
\end{equation}

\noindent
{\bf ii)}
This extends to a right Quillen
functor on global projective model categories of presheaves
(Example \ref{ModelStructureOnSimplicialPresheavesOverCartesianSpaces},
Example \ref{ProjectiveModelStructureOnPresheavesOfConnectivedChainComplexes}).
whose right derived functor (Prop. \ref{DerivedFunctors})
$\RightDerived \mathrm{DK}$ composed with the
$\infty$-stackification functor \eqref{InfinityStackificationFunctor}
is thus of the form
\vspace{-2mm}
$$
  \xymatrix@C=3em{
    \mathrm{Ho}
    \Big(
    \mathrm{PSh}
    \big(
      \CartesianSpaces
      \,,\,
      \ZChainComplexes
    \big)_{\mathrm{proj}}
    \Big)
    \ar[rr]^-{
      \overset{
        \mathclap{
        \raisebox{3pt}{
          \tiny
          \color{darkblue}
          \bf
          \def\arraystretch{.9}
          \begin{tabular}{c}
            derived
            \\
            Dold-Kan construction
          \end{tabular}
        }
        }
      }{
        \RightDerived\mathrm{DK}
      }
    }
    \ar[drr]_-{
      \mbox{
        \tiny
        \color{darkblue}
        \bf
        \def\arraystretch{.9}
        \begin{tabular}{c}
          $\infty$-stackified
          \\
          Dold-Kan construction
        \end{tabular}
      }
    }
    &&
    \mathrm{Ho}
    \Big(
    \mathrm{PSh}
    \big(
      \CartesianSpaces
      \,,\,
      \SimplicialSets
    \big)_{\mathrm{proj}}
    \Big)
    \ar[d]^-{
      L^{\mathrm{loc}}
      \mathrlap{
        \mbox{
          \tiny
          \color{darkblue}
          \bf
          $\infty$-stackification
        }
      }
    }
    \\
    &&
    \Stacks
  }
$$
and preserves homotopy pullbacks (by Lemma \ref{InfinityStackification}).
\end{example}

\begin{example}[Projective model structure on unbounded chain complexes {\cite[Thm. 2.3.11]{Hovey99}}]
  \label{ProjectiveModelStructureOnUnboundedChainComplexes}
  The category $\ZChainComplexesUnbounded$
  of unbounded chain complexes
  of abelian groups
  carries a model category structure  (Def. \ref{ModelCategories})
  whose:
  \begin{itemize}
   \vspace{-2mm}
    \item[$\mathrm{W}$ --]  weak equivalences are the quasi-isomorphisms;
    \vspace{-3mm}
    \item[$\mathrm{Fib}$ --]  fibrations are the degreewise
      surjections.
  \end{itemize}

 \vspace{-2mm}
 \noindent  We write
  $\ZChainComplexesUnboundedProj$ for this model category.
\end{example}

\begin{prop}[Stable Dold-Kan construction]
  \label{StableDoldKanCorrespondence}
  The Dold-Kan construction (Def. \ref{DoldKanConstruction})
  lifts along the stabilization adjunction (Example \ref{DerivedStabilizationAdjunction})
  from connective to unbounded
  chain complexes (Example \ref{ProjectiveModelStructureOnUnboundedChainComplexes}),
  such as to make the following
  diagram commute:

  \vspace{-3mm}
  \begin{equation}
    \label{StableDoldKanConstruction}
    \hspace{-3mm}
    \raisebox{30pt}{
    \xymatrix@C=10pt{
      \mathrm{Ho}
      \Big(
        \big(
          \mathrm{ChainComplexes}^{
            \scalebox{.5}{$\geq 0$}
          }_{
            \scalebox{.5}{$\mathbb{Z}$}
          }
        \big)_{\mathrm{proj}}
      \Big)
      \ar@/^1.6pc/[rrrr]^-{
        \overset{
          \mathclap{
          \raisebox{3pt}{
            \tiny
            \color{darkblue}
            \bf
            Dold-Kan correspondence
          }
          }
        }{
          \RightDerived\mathrm{DK}
        }
      }
      \ar[rr]^-{\simeq}
      \ar@{_{(}->}@<-6pt>[d]_-{  }^-{ \dashv }
      \ar@{<-}@<+6pt>[d]^-{
        \RightDerived\Omega^\infty
      }
      &&
      \mathrm{Ho}
      \big(
        \SimplicialAbelianGroups_{\mathrm{proj}}
      \big)
      \ar[rr]
      &&
      \mathrm{Ho}
      \big(
        \SimplicialSets_{\mathrm{Qu}}
      \big)
      \ar@<-6pt>[d]_-{ \LeftDerived\Sigma^\infty }^-{ \dashv }
      \ar@{<-}@<+6pt>[d]^-{ \RightDerived\Omega^\infty }
      \\
      \mathrm{Ho}
      \Big(
        \big(
          \mathrm{ChainComplexes}_{\scalebox{.5}{$\mathbb{Z}$}}
        \big)_{\mathrm{proj}}
      \Big)
      \ar@/_1.6pc/[rrrr]_-{
        \underset{
          \mathclap{
          \raisebox{-3pt}{
            \tiny
            \color{darkblue}
            \bf
            stable Dold-Kan construction
          }
          }
        }{
          \RightDerived\mathrm{DK}_{\mathrm{st}}
        }
      }
      \ar[rr]^-{ \simeq }
      &&
      \mathrm{Ho}
      \big(
        (H\mathbb{Z})\mathrm{ModuleSpectra}
      \big)
      \ar[rr]
      &&
      \StableHomotopyTypes
      \,.
    }
    }
  \end{equation}

  \vspace{-1mm}
\noindent  Here the right adjoint on chain complexes is the
  homological truncation from below:
   \vspace{-2mm}
  \begin{equation}
    \label{HomologicalTruncationFromBelow}
    \RightDerived\Omega^\infty
    \Big(
      \,
      \cdots
      \overset{\partial_2}{\longrightarrow}
      V_2
      \overset{\partial_1}{\longrightarrow}
      V_1
      \overset{\partial_0}{\longrightarrow}
      V_0
      \overset{\partial_{-1}}{\longrightarrow}
      V_{-1}
      \overset{\partial_{-2}}{\longrightarrow}
      \cdots
      \,
    \Big)
    \;\;\;
    =
    \;\;\;
    \Big(
      \,
      \cdots
      \overset{\partial_2}{\longrightarrow}
      V_2
      \overset{\partial_1}{\longrightarrow}
      V_1
      \overset{\partial_0}{\longrightarrow}
      \mathrm{ker}(\partial_{-1})
      \,
    \Big).
  \end{equation}
\end{prop}
\begin{proof}

  \noindent
  {\bf (i)}
  It is clear from inspection that the assignment
  \eqref{HomologicalTruncationFromBelow} is right adjoint to the
  inclusion of connective chain complexes, so that we have a pair of
  adjoint functors

  \vspace{-2mm}
  \begin{equation}
    \label{StabilizationQuillenAdjunctionForChainComplexes}
    \xymatrix{
      \ZChainComplexesUnboundedProj
    \;  \ar@{<-^{)}}@<+7pt>[rr]
      \ar@<-7pt>[rr]_-{ \Omega^\infty }^-{
        \bot_{\mathrlap{\mathrm{Qu}}}
      }
      &&
    \;\;  \ZChainComplexesProj
      \,.
    }
  \end{equation}
  \vspace{-.3cm}

  \noindent
  Moreover, it is immediate that this is a Quillen adjunction
  (Def. \ref{QuillenAdjunction}) between the projective model
  structure on connective chain complexes
  (Example \ref{ProjectiveModelStructureOnConnectiveChainComplexes})
  and that on unbounded chain complexes
  (Example \ref{ProjectiveModelStructureOnUnboundedChainComplexes}):
  $\Omega^\infty$ clearly preserves fibrations
  (using that those between connective chain complexes
  need to be surjective only in positive degrees!)
  and clearly preserves all weak equivalences.
  Finally, since all chain complexes in the projective model structure
  are fibrant, we have that with $\Omega^\infty$
  also $\RightDerived\Omega^\infty$ is given by \eqref{HomologicalTruncationFromBelow},
  via Example \ref{DerivedFunctorsByCoFibrantReplacement}.

\vspace{1mm}
  \noindent
  {\bf (ii)}
  A Quillen adjunction of the form

  \vspace{-.5cm}
  \begin{equation}
    \label{StableDoldKanQuillenEquivalence}
    \hspace{-5mm}
    \xymatrix@C=3.5em{
      \ZChainComplexesUnboundedProj
      \ar@{<-}@<+7pt>[rr]
      \ar@<-7pt>[rr]_>>>>>>>>>{H}^-{ \bot_{\mathrlap{\simeq}} }
      \ar@/_1.8pc/[rrrr]_-{ \mathrm{DK}_{\mathrm{st}} }
      &&
      (H\mathbb{Z})\mathrm{ModuleSpectra}
      \ar@{<-}@<+7pt>[rr]
      \ar@<-7pt>[rr]^-{ \bot_{\mathrlap{\mathrm{Qu}}} }
      &&
      \mathrm{SequentialSpectra}_{\mathrm{BF}}
    }
  \end{equation}
  \vspace{-.4cm}

  \noindent
  is established in \cite[\S B.1]{SchwedeShipley01},
  where

\vspace{-2mm}
\begin{enumerate}[{\bf (a)}]
\setlength\itemsep{-2pt}
  \item the first step is a
  Quillen equivalence (Def. \ref{QuillenEquivalence})
  between the projective model structure on
  unbounded chain complexes (Example \ref{ProjectiveModelStructureOnUnboundedChainComplexes})
  and a model category of module spectra over the
  Eilenberg-MacLane spectrum $H\mathbb{Z}$
  \cite[\S B.1.11]{SchwedeShipley01};

 \item the second step
  is a Quillen adjunction
  \cite[p. 37, item ii)]{SchwedeShipley01}
  to the Bousfield-Friedlander model
  structure
  (Example \ref{ModelStructureOnSequentialSpectra})
  whose right adjoint assigns underlying sequential spectra;
  such that

 \item
  the composite right adjoint
  $\mathrm{DK}_{\mathrm{st}}$
  \eqref{StableDoldKanQuillenEquivalence}
  further composed with
  $\Omega^\infty$ on spectra \eqref{StabilizationQuillenAdjunction}
  equals the composite of $\Omega^\infty$ on chain complexes
  \eqref{StabilizationQuillenAdjunctionForChainComplexes}
  with the unstable Dold-Kan construction \eqref{DK}:

  \vspace{-.3cm}
  $$
    \Omega^\infty
    \,\circ\,
    \mathrm{DK}_{\mathrm{st}}
    \;\;
    \simeq
    \;\;
    \mathrm{DK}
    \,\circ\,
    \Omega^\infty
  $$
  \vspace{-.5cm}

  \noindent
  (by immediate inspection of the construction in \cite[p. 38-39]{SchwedeShipley01}).
\end{enumerate}

\vspace{-2mm}
  \noindent
  {\bf (iii)}
  By uniqueness of adjoints, this implies that the
  Quillen adjunction
  of the stable Dold-Kan construction \eqref{StableDoldKanQuillenEquivalence}
  is intertwined by the Quillen adjunctions involving $\Omega^\infty$
  with the Quillen adjunction of the unstable
  Dold-Kan construction \eqref{DK},
  and hence the commuting diagram of derived functors
  \eqref{StableDoldKanCorrespondence} follows
  (Prop. \ref{DerivedFunctors}).
\end{proof}

\noindent Domenico Fiorenza, {\it Dipartimento di Matematica, La Sapienza Universita di Roma, Piazzale Aldo Moro 2, 00185 Rome, Italy.}
\\
{\tt fiorenza@mat.uniroma1.it}
\\
\\
\noindent  Hisham Sati, {\it Mathematics, Division of Science, New York University Abu Dhabi, UAE.}
\\
{\tt hsati@nyu.edu}
\\
\\
\noindent  Urs Schreiber, {\it Mathematics, Division of Science, New York University Abu Dhabi, UAE;
on leave from Czech Academy of Science, Prague.}
\\
{\tt us13@nyu.edu}

\end{document}